\documentclass[12pt]{amsart}

\pdfoutput=1

\usepackage{amsmath}

\usepackage{amsthm}

\usepackage{mathrsfs}

\usepackage{latexsym}

\usepackage[english]{babel}

 \usepackage{amssymb}

\usepackage{enumitem}

\usepackage[active]{srcltx}

\usepackage{xcolor}
\usepackage[colorlinks=true,linkcolor=blue,citecolor=red,filecolor=yellow,
 linkcolor=purple,urlcolor=green]{hyperref}

\theoremstyle{plain}

\newtheorem{thm}{Theorem}[section] 
\newtheorem{lem}[thm]{Lemma} 
\newtheorem{prop}[thm]{Proposition} 
\newtheorem{cor}[thm]{Corollary} 
\newtheorem{fact}[thm]{Fact}
\newtheorem{convent}[thm]{Convention}
\theoremstyle{definition}
\newtheorem{defn}[thm]{Definition}
\newtheorem{exmp}[thm]{Example}
\newtheorem{rem}[thm]{Remark}
\newtheorem{assump}[thm]{Assumption}

\newtheorem{nota}[thm]{Notation}

\theoremstyle{abstract}
\theoremstyle{remark}

\newtheorem*{note}{Note}

\renewcommand{\leq}{\leqslant}
\renewcommand{\geq}{\geqslant}

\newcommand{\wt}{\widetilde}

\newcommand{\R}{\mathbb{R}}
\newcommand{\Z}{\mathbb{Z}}
\newcommand{\N}{\mathbb{N}}

\newcommand{\ov}{\overline}

\newcommand{\grad}{\hbox{\rm grad}}
\newcommand{\vol}{\hbox{\rm vol}}
\newcommand{\diam}{\hbox{\rm diam}}

\usepackage[applemac]{inputenc}
\DeclareMathAlphabet{\mathpzc}{OT1}{pzc}{m}{it}

\title{Approximating Riemannian manifolds by polyhedra}

\begin{document}{}

\maketitle

\centerline{{\it by} Daniel Meyer {and} Eric Toubiana}
\vskip1mm
\centerline{\tiny Institut Math\'ematique de Jussieu-Paris Rive Gauche, Universit\'e Paris Cit\'e}
\centerline{\tiny Campus des Grands Moulins, B\^atiment Sophie Germain, Bo\^ite Courrier 7012}
\centerline{\tiny 8 place Aur\' elie Nemours, 75205, Paris Cedex 13}
\centerline{\tiny daniel.meyer@imj-prg.fr 
\hskip2mm  \hskip2mm eric.toubiana@imj-prg.fr}

\begin{abstract} {\tiny This is a study on approximating a Riemannian manifold by polyhedra. 
Our scope is understanding Tullio Regge's article \cite{R} in the restricted Riemannian frame. We give a proof of the ``Regge theorem'' along lines close to its original intuition: {\it one can approximate a compact domain of a Riemannian manifold by polyhedra {\em in such a way that}
the integral of the scalar curvature is approximated by
a corresponding polyhedral curvature.}}
\end{abstract}

\vskip1cm
This is a study in two parts on approximating a Riemannian manifold  $(M,g)$ by polyhedra. The third part comprises technical tools involving singularities and transversality, used in the previous parts. 

\vskip1mm
Our scope is understanding Tullio Regge's elegant, anchoring and visionary paper \cite{R} in the restricted mathematical Riemannian frame. One of his deepest ideas is to connect the parallel translation along a loop bounding a disk to Dirac masses of curvature felt over this disk. Those masses are defined in a polyhedron by generalising the curvature proposed in dimension $2$ by A. D. Alexandrov (p. 36, chapter I, §8 in \cite{Alex1}), who investigated the convergence of $2$-polyhedra towards a Riemannian surface (see especially chapter I, chapter III, §6, n° 17, pp. 84-89 in \cite{Alex2}, chapter VII  §5) and settled the question up to the level of curvature.
\vskip1mm
Though we steadily rely on their pioneering work in establishing a good part of our first important theorem below, we diverge from the classical Cheeger-M\"{u}ller-Schrader approach \cite{C-M-S2} by focusing on the above geometrical tool proposed by Regge, in a rather intuitive way, and follow his intuition. 
It leads us to the notions of {\em Regge curvature} and {\em curving}, see below. 
\vskip1mm 
Before going into more details of this introduction, we give the table of contents of this text which might be useful for the reader.

 \tableofcontents{}
 \vspace{-1cm}
\part*{Introduction} 
Let us insist on a postulate which ruled the making of our work: 
 {\it almost nowhere did we seek for optimality.}
 \vskip1mm
 
In part \ref{M-T1}, the central theorem below states a metric convergence of a sequence of piecewise-flat triangulations towards a compact domain $D\subset (M,g)$, while, along a given curve, it also states convergence at the level of parallel translations. Beyond the various meanings of the word {\em polyhedron}, we also make the special use outlined in the

\begin{convent}\label{D.metric133}
{\em Call {\it polyhedron} a linear finite simplicial complex $K$, homeomorphic to a compact connected $n$-dimensional manifold with boundary, {\em endowed} with one of the following geometric features 
\par - either (see definition \ref{D.metric1}) it carries a natural analytic locally Euclidean metric $g_0$ on $K\setminus K_{n-2}$ whose analytic prolongation over the $(n\!-\!2)$-skeleton $K_{n-2}$ shows any $n$-simplex as being isometric to a linear simplex in $\R^n$ (this is a {\em pseudomanifold} in \cite{Ch}, see page 613). Thus, there exists a natural $g_0$-parallel translation along a regular path through $K\setminus K_{n-2}$. Moreover, $(K, g_0)$ is entirely determined by the combinatorial attachments and the collection, over all $n$-simplices, of all the distances between their vertices (see \cite{R}, §4 pp. 564-5). 
\par - or, as the image of such an object $(K,g_0)$ under a continuous (piecewise) smooth {\em simplicial embedding} $T$ (definition \ref{D.Simplicial.Embedding}) in $(M,g)$, it may instead show up $C^1$-discontinuities across the $(n-1)$-faces.}
\end{convent}
\textit{Unless otherwise stated, differentiable, or smooth, means $C^\infty$ and differentiable manifolds are $C^\infty$-manifolds.}
\begin{defn}\label{A.approx} \textit{(Definition {\rm\ref{convpolyhedra}})} Consider a compact domain $D$ in an $n$-dimensional ($n\!\geq\!2$) connected differentiable Riemannian manifold $(M,g)$, 
meaning that $D\subset M$ is an $n$-dimensional differentiable submanifold, possibly   
with boundary, which is  
compact and connected.

{\it A  
strong metric approximation by polyhedra of $(M,g)$ over $D$} requires constants
${\mathfrak C}_1$ and ${\mathfrak C}_2$ - depending on the geometry of the inclusion $D\subset(M,g)$ -
and a family of triplets $(T_\rho,K_\rho,(g_0)_\rho):=(T,K,g_0)\,$ -  indexed by a real $\rho\!>\!0$,
{\it but we drop the subscript $\rho$} - consisting in a polyhedron $(K,g_0)$, and a 
continuous differentiable simplicial embedding (definition \ref{D.Simplicial.Embedding}) $T$ of $K$ into $M\,$ satisfying the following conditions

 $(i)$ $T(K)$ contains $D$ in its interior;

 $(ii)$ for any $n$-simplex $\sigma$ of $K$ with vertices 
$s_0,s_1,\dots,s_n\,$,  
one has 
\begin{gather*}\hskip6mm 
(a)\hskip6mm\forall i,j=0,1,\dots,n
\hskip6mm d_g(T(s_i),T(s_j))=d_{g_0}(s_i,s_j)\,;\\
\hskip6mm (b)\hskip3mm  \forall q\!\in\!\sigma \hskip3mm\forall u\!\in\! T_q\sigma\hskip3mm\vert T^\ast 
g(u,u)\!-\!g_0(u,u)\vert \!\leq \!
{\mathfrak C}_1\,\rho^2\,T^\ast g(u,u)\,,
\end{gather*}
where $T^\ast g$ is on $\sigma$ the pull-back metric of $g$ by the embedding $T_{\mid\sigma}$.

$(iii)$ For any regular curve 
 $\gamma:[0,1]\rightarrow D\setminus T(K_{n-2}),$ denote by ${\mathcal P}_\gamma^g$ 
(respectively by ${\mathcal P}_{T^{-1}\gamma}^{g_0}$) 
the $g$-parallel 
translation 
along $\gamma$ from $\gamma (0)$ to $\gamma (1)$, 
(the corresponding $g_0$-parallel 
translation  along $T^{-1}\gamma$ from $T^{-1}\gamma (0)$ to $T^{-1}\gamma (1)$). 
There exists a dense open set 
${\mathcal V}\subset C^\infty(M,\R^{n-1})$ 
such that, for any $f \in \mathcal {V}$, there exists 
$\rho_f>0$ verifying, for any regular curve $\gamma$ with ${\rm rge}(\gamma)\!\subset\!f^{-1}\{0\}$, for 
any $\rho\in]0,\rho_f]$ and any $u\in T_{\gamma(0)}M$
\vspace{-1mm}
\begin{equation*} 
\hskip4mm (c)\hskip3mm \Vert{\mathcal P}_{\gamma}^g(u)- 
dT ({\mathcal P}_{T^{-1}\gamma}^{g_0}(dT^{-1}(u)))\Vert_g\leq 
{\mathfrak C}_2\,\rho\,\hbox{\rm length}(\gamma)\,\Vert u\Vert_g\  .
\end{equation*}
\end{defn}

\begin{thm} \label{fundametal} \textit{(Theorem {\rm\ref{theorem I}})} Let $D$ be a compact domain in an open $n$-dimensional Riemannian manifold 
$(M,g)$. One can construct a  
{\it strong metric approximation of $(M,g)$ over $D\,$ by polyhedra} $(T,K,g_0)$ (the $(T,K,g_0)$ depend on $\rho\rightarrow0\,$, see the above definition). 
\end{thm}
\vspace{-2mm}
Our construction involves at its core the barycentre map in the Riemannian context, which goes back to Cartan, see \cite{K-NII}, pages 108ff and \cite{B-K}, \cite{K}, but also \cite{B-C-M} and \cite{Sa} for more recent developments.
\par
\begin{rem}\label{lampionSch} A consequence of the above theorem (and of property $(ii)$ stated in definition  \ref{A.approx}) is that, for a a fixed domain, the volume in the Riemannian metric is the limit of the sequence of its volumes computed in the approximating polyhedral metrics. Thus, {\em a Schwarz' lantern (see  \cite{Be1} vol. {\rm 1} or \cite{Mo}) does not approximate the cylinder in this way}. 
This is quite different from the one-dimensional situation where any sequence of polygonal lines, even {\em irregular}, inscribed in the circle converges in length to the circle if its mesh tends to zero, just as well as Archimedes' regular polygonal approximation to the circle does. 
\end{rem}
{\it Nota Bene}.
Several years after theorem \ref{fundametal} had been obtained (this necessitated the work published in \cite{C-M}) but not published, another text (\cite{DVW}) came out quite incidently to the 
knowledge of the authors: its aim is to ensure that, given a set of points in a Riemannian manifold $(M,g)$ and an abstract simplicial complex admitting this set as vertices, the barycentric riemannian coordinates realise a triangulation of $M$. Though the authors do not explicitly state the existence of a polyhedral approximation satisfying conditions $(i)$ and $(ii)$, we found enough statements and ideas related with results of part of our part \ref{M-T1} (their presentation is different from ours). Our task is to write a self-contained proof of the above theorem stating an approximation satisfying all properties of definition \ref{A.approx}, {\em including} $(iii)$.
\par\vskip1mm
{\it In this introduction,} call {\it tile} the elementary piece $\hat\tau$ of a polyhedral approximation {\em to build}.
A tile is constructed as a {\it Riemannian barycentric $n$-dimensional simplex} through Riemannian barycentric coordinates defined by a {\it spread} set of $n+1$ points $p_0,\dots,p_n$ (see section \ref{rbs}). It is shown in theorem \ref{2n} that tiles can be equipped with a natural Euclidean metric $\hat g$ verifying $d_{\hat g}(p_i,p_j)=d_g(p_i,p_j)$. The closeness of the $g$ and $\hat g$ metrics and parallel translations on a tile is studied in theorem \ref{theorem 3}, which is the one-tile version of theorem \ref{fundametal}. 
\vskip2mm
In a triangulation, the {\em mesh} $\rho$ - see the opening of section \ref{6.1} - roughly controls the diameter of the simplices. Starting from an initial given finite triangulation $(T,K)\,$, which is known to exist  (this is J. H. C. Whitehead's theorem), one
produces tiles by refining enough $(T,K)\,$ in a {\it regular} way (the {\em mesh} $\rho$ goes to $0$), then by showing that each $n$-simplex $\tau$ of such a refinement can be replaced by a {\it tile} $\hat\tau\,$, defining in this way a new triangulation $\hat T\,$ which displays the same combinatorics, is piecewise $C^1$-close to $T$ and verifies over each tile the one-tile central theorem: this points to theorems \ref{theorem 3}, \ref{2n} and \ref{hatT}. 
\vskip2mm
A good part of part \ref{M-T1} is based on a careful and, hopefully, faithful analysis of Munkres' book \cite{Mu}, this is mainly seen in sections \ref{2.1}, \ref{6.1}.
\vskip2mm
In section \ref{2.1}, we present the basic notions shaping a linear simplex, like the thickness and openness. In section \ref{2.2} results on the proximity of simplices are given. In section \ref{2.3} we single out the fundamental theorem of I. Schoenberg \cite{Sc}.
\vskip2mm
The local quasi-isometric character (see \cite{B-K}) of the map $\exp_x$ around a point $x$ in $(M,g)$ is developed in section \ref{EuclRiem} and exploited in section \ref{bdvargeod}: useful bounds involving parallel transport along a curve are provided. In section \ref{diffparll} we recall the role of the difference connexion tensor in view of getting a control on the difference of parallel translations associated with two different metrics along a given curve.
\vskip2mm
In section \ref{4}, more Riemannian bounds needed later are produced in a non-classical way (and far from optimality) as an introduction to {\it a central tool} of this work, consisting in Gromov's magnifying glass \cite{GLP}.
\vskip2mm
In establishing the convergence of parallel translations along a given curve, a difficulty shows up, caused by the piecewise character of the metric $\hat g$ constructed over $D\,$. A rip is felt while crossing an interface of two tiles; think to a regular geodesic octahedron on the $2$-sphere, first in the metric of the sphere, then in the Euclidean metric verifying $(ii) (a)$ of definition \ref{A.approx} (see section \ref{6.2}): the unavoidable rip is computable at each vertex. To overcome this difficulty, we produce an {\it a priori} suitable control on the total number of such crossings by a given curve. This leads to a theorem of Ren\' e Thom on singularities (see \cite{Th}, \cite {C-M}) and relies heavily on a universal bound of the {\it local degree of intersection}, which generalizes in our setting the fact that the local degree of intersection between a generic curve and any straight line in the plane is $\leq3\,$.
In order to apply this theorem, we have to produce an object called {\it ``Texture''} in \cite{C-M} requiring a new geometrical construction: the manifold of Riemannian barycentric $k$-simplices included in a convex and small enough ball $\subset (M,g)$, which generalizes the set of straight lines in the Euclidean plane. This construction fills section \ref{6.3} while theorem
\ref{text-11} is the version of the Thom theorem needed here.
Theorem \ref{fundametal} is proved in section \ref{7}, putting together theorems \ref{theorem 3},\,\ref{2n},
\,\ref{hatT},\,\ref{text-11} and propositions  \ref{rip1}, \ref{bds45}. Lemma \ref{simpletext} says that all buildable small Riemannian barycentric $n$-simplices of interest to us slip into the ``leaves'' of the constructed ``textures''.

\vskip1mm
Part \ref{M-T2} is dedicated to results on curvature convergences of polyhedral approximations given by theorem \ref{fundametal} (i. e. theorem \ref{theorem I}). We make the
\begin{assump} \label{assumpM-T2} In the sequel, $T,K,(M,g)$ tacitly denotes an approximation of the type stated in the above theorem \ref{fundametal} which is also a {\em Riemannian barycentric triangulation} satisfying definition \ref{A.approx} (actually theorem \ref{theorem I}, definition \ref{convpolyhedra} and remark \ref{remarquable1}). It is obtained from a regular subdivision of a given finite embedded simplicial complex, along the lines of part \ref{M-T1}. So it also shares the important properties of Thom's local theorem {\rm\ref{text-11}}, even of a refined version, theorem \ref{text-1133}.
\end{assump}
\vskip1mm
In section \ref{holonomy1}, we recall how, in the Riemannian case, along a loop which is trivial in homotopy, the curvature is linked to the holonomy at its base-point $p$. Indeed,
by means of a ``Green-formula'', for any embedded disk having this loop as boundary, one can express this holonomy as a mean of curvature endomorphisms spread over the disk and translated back along privileged paths to $p\,$ (this is later formalized in the definition of the {\em Gauss curving}, see \ref{Gaussmeas}): this fills theorems \ref{theorem A} and \ref{theorem B}, stated and proved in a way that allows later a parallel treatment for polyhedra. The natural object that takes the role of the aforementioned disk is {\em an embedded parametrised square} with all its squaring lines filling it up to the boundary.  It is easily divided into a union of adjacent analogous subsquares.
Section \ref{proximity33} presents the technical proposition \ref{prox1} enabling us to deal with close holonomies in the proof of theorem \ref{theorem G} (see lemma \ref{thG4}).
\vskip2mm
Following Regge (in the steps of A. D. Alexandrov \cite{Alex1}, \cite{Alex2}), the first way to feel polyhedral curvature looks as follows: given a polyhedron $(K,g_0)$, its curvature is revealed by a disturbance acted by {\em an $(n\!-\!2)$-simplex $\xi$ contained in the interior of $K\,$}, called {\em bone}: the parallel translation along a loop turning once around $\xi$ and contained in the open star of $\xi$ (which is the union of the interiors of simplices containing $\xi$) is a Euclidean isometry of the tangent space taken at the base-point of the loop; but {\em this isometry usually differs from} $\rm{Id}$ on a plane orthogonal to the $(n\!-\!2)$-vector subspace parallel to $\xi$, thus is measured by a {\em defect angle}. This picture is the elementary geometrical brick showing curvature at which Regge points in \cite{R}. Sections \ref{polyedres}, \ref{restrholo} present the above material, see in particular lemma \ref{Reggeoperbonebis}.
\vskip2mm
For a given parametrised square, one wraps the above description in the {\em Regge curvature} of definition \ref{Regge2}; in the generic case, one is able to estimate it (see lemmas \ref{geomversion}, \ref{singular1}) after another development. First we sort out embedded parametrised squares by properties they share {\em in the generic case} with respect to intersections with $(n\!-\!2)$-simplices of $T(K)\!\subset \!(M,g)$. Then, section \ref{play} mimicks for polyhedra section \ref{holonomy1}, leading to propositions \ref{Regge6}, \ref{ReggeGB}, \ref{Regge6bis} and \ref{ReggeGBbis} which depict the
{\em Regge curvature} as a geometric action of $K_{n-2}$ on a parametrised square $G$. 
\vskip2mm
Genericity questions on intersections of squares (of squaring lines) with $(n\!-\!2)$-simplices (with $(n\!-\!1)$-simplices) of $T(K)$ are mentioned in the second point of assumption \ref{assumptionX} (section \ref{contacted}) and treated in sections \ref{RBT}, \ref{multiple1088}, \ref{adaptThom}, producing the central theorem \ref{text-1133}, proved in its version \ref{text-1133bis} (see part \ref{annexe1033}). Actually, a ``testable'' (we most often say {\em ``generic''}) square $G$ is defined in definition \ref{testsquare1}. It has two aspects: the first points to transversality properties of its intersections with Riemannian barycentric $(n\!-\!1)$ and $(n\!-\!2)$-simplices {\em in general} (actually with leaves of some textures), while the second refers to a given approximating polyhedron $K$ paired with a fixed $\rho>0$.

One can define an {\em index of intersection} between a bone and a {\em generic} parametrised square: see remark \ref{smallandneat}, definition \ref{indiceGxi}, lemma \ref{indice34}.
\vskip2mm
To manage the questions of convergences, we have to care about dihedral angles and show that the defect angle around a bone $\xi$ is bounded in terms of $\rho^2$ as the mesh $\rho$ tends to $0$, this is done in sections \ref{dihedral33}, \ref{desirable}, see proposition \ref{diedre4}.
\vskip2mm
Finally, after the study of section \ref{play}, the {\em Regge curving}, defined in section \ref{convenable} for a given parametrised square as a sum of antisymetrised Regge curvatures, appears to be the polyhedral counterpart of the {\em Gauss curving}; the above {\em index of intersection} comes in: definition \ref{Regge22} covers the simplest case, leading to definition \ref{ReggeRegge!!} and its rephrasing, definition \ref{Reggemeas}. The closeness of the Regge curving to a sum of elementary {\em Regge curvatures}, when evaluated on a generic parametrised square, is controlled in terms of $\rho^4$ when $\rho\rightarrow0$, see lemma \ref{Regge-2033}. 
\vskip2mm
These constructions are in some sense {\em hybrid}. On one hand, for a given square $G\,$, they are defined in terms of dihedral angles and privileged $g_0$-parallel translations (given with $G$) back to the origin of $G$: such parallel translations are naturally located on $T^{-1}(G)$ subset of the polyhedron $(K,g_0)$. On the other hand, they embody indices of generic intersections between embedded squares and bones, and this makes sense in $(M,g)$ by considering the embedded polyhedron $T(K)$.
\vskip2mm
{\em The reader may compare the {\em Regge curving} with the development in {\rm\cite{R}}, looking at the formulas written on the bottom of page {\rm 565}.}
\vskip2mm
Proposition \ref{Riesqnb2} gives bounds in terms of the mesh $\rho$ (for $\rho$ close to $0$) on the number of intersections between a generic square $G$ and all bones $\xi$. It is a prerequisite for the proofs of the main convergence theorems \ref{theorem F} and \ref{theorem G}, stated and proved in section \ref{main33}. 
\vskip2mm
In section \ref{pulling33}, definitions \ref{Reggemeas} and \ref{Gaussmeas} of {\em Regge} and {\em Gauss curvings} allow to formulate theorem \ref{theorem G} in a concise and conceptual way: indeed, theorem \ref{convergence00} roughly says that, relying on theorem \ref{fundametal}, {\em there exist polyhedral approximations of a Riemannian manifold for which the Regge curving converges to the Gauss curving}.
\vskip3mm
In section \ref{Reggetheorem}, the classical Regge theorem (see \cite{R}, \cite{C-M-S2}) is stated as theorem \ref{Regge1234} and proved using the work done before: given a polyhedral approximation $(K,g_0)$ of a domain $D\subset(M,g)$, the {\em Regge $g_0$-scalar curvature} is $\sum_\xi \alpha_\xi\,{\rm Vol}_{n-2}^{g_0}(\xi)$, where the sum is over all the bones $\xi\subset D$ of $K$ and $\alpha_\xi$ is the $g_0$-defect angle around $\xi$. The Regge theorem says that {\em the Regge $g_0$-scalar curvature converges, up to a universal constant, to the integral over $D$ of the $g$-scalar curvature as $\rho$ goes to $0$}.
\vskip1mm
The first mathematical proof of Regge's theorem was proposed by Cheeger,
M\"{u}ller and Schrader (see \cite{C-M-S2} and \cite{La}). Following Regge's original intuition, we give here a quite different proof.
\vskip1mm
To derive this theorem from theorem \ref{theorem G}, we introduce in section \ref{Riemannian1111} the notion of families of Riemannian squares $\underline G_r$, which are, around any $x\!\in\! M$, images under $\exp_x$ of Euclidean squares of side $r$ in $(T_xM,g_x)$, parametrised by $r>0$ and by the Stiefel manifold ${\bf St}_2(TM)$ of pairs of $g$-orthonormal vectors in $(T_xM,g_x)$, for any $x\,\in\,M$. Then, we switch to families of quasi-Riemannian squares $\bar G_r$ which are close to $\underline G_r$ and generic in the sense that theorem \ref{text-1133} holds for them.
\vskip2mm
Afterwards, useful integrals of curvatures, namely integrals of {\em curvings}, are defined in section \ref{someintegrals}. For $r>0$ fixed and small, define $I^g(\bar G_r)$ (define $I^{g_0}(\bar G_r)$) to be the $r^{-2}$-weighted mean of the {\em Gauss} (or the {\em Regge}) curving over ${\bf St}_2(TD)$. Theorem \ref{theorem G} implies that $I^{g_0}(\bar G_r)$ tends to $I^g(\bar G_r)$ as $\rho$ tends to $0$. Another technical intermediate integral $J^{g_0}(\bar G_r)$ is also introduced as a companion of $I^{g_0}(\bar G_r)$.
\vskip2mm
In section \ref{sketchRegge}, one mainly {\em states and admits} proposition \ref{CMS-7} and quotes (this is classical) that,  as $r\rightarrow0$ and up to universal constant, $I^g(\bar G_r)$ tends to the integral over $D$ of the scalar curvature (proposition \ref{CMS11}), proving the Regge theorem in section \ref{Regge333333} thanks to proposition \ref{CMS-7}. 
The proof of proposition \ref{CMS-7} goes through a description of $\bar G_r$ in Fermi coordinates around a given {\em bone} $\xi$ done in section \ref{bone1111}.
\vskip2mm
Proposition \ref{CMS-7} $(i)$ says that, uniformly in $r>0$ belonging to a given small interval and as $\rho$ tends to $0$, $I^{g_0}(\bar G_r)$ and $J^{g_0}(\bar G_r)$ near at least as fast as $\rho$; this is proved in section \ref{part(i)}.
\vskip2mm
Proposition \ref{CMS-7} $(ii)$ manages the computation of the limit of $J^{g_0}(\bar G_r)$ as $r$ tends to $0$, which involves a kinematic Euclidean integral (section
\ref{integraleuclidean}) and a tricky application of the Lebesgue theorem: it requires again a magnification argument producing a model situation (section \ref{magnificationagain1}), but also neglecting some boundary effects, see section \ref{negligeons}.
\vskip3mm
In part \ref{annexe1033}, we chose to outline some required statements and facts in singularity theory, since a development in the core of the text would have caused a detour. In section \ref{multiple1088}, we present briefly some material extracted from \cite{C-M}, in particular the useful {\em textures}, and outline properties we need. In section \ref{adaptThom}, we state and prove the central theorem \ref{text-1133bis}, which appeared in a simpler version in part \ref{M-T1} as theorem \ref{text-11}, then in its full strength in part \ref{M-T2} as theorem \ref{text-1133}.
Finally, section \ref{Whitneyandal} presents several proofs of technical lemmata involving extensions or restrictions of $C^\infty$-maps and their relations to Whitney $C^\infty$-topology.
\vskip2mm
A word about the semi-Riemannian case, which embodies relativity, the original scope of Regge's article \cite{R}. Part \ref{M-T1} and theorem \ref{theorem I} rely heavily on the choice of a Riemannian metric, through barycentric constructions, convexity results, triangle inequality, etc... In part \ref{M-T2}, given a semi-Riemannian metric of signature $(p,q)$ (with $p\!+\!q\!=\!n$, the ambient dimension), integrating around a loop which bounds a disk defines a holonomy transformation in the relevant group of isometries $SO_0(p,q)$: the changes seem there to be manageable, since there is a Levi-Civita
connexion associated to a semi-metric (see \cite{BO}, \cite{M-T-W} or
\cite{St}). In what concerns a Regge's type theorem, the situation needs to handle means under $SO_0(p,q)$, which no longer acts transitively on the tangent space. Maybe one could try to pair the semi-Riemannian metric with some related Riemannian metric in order to follow a route close to ours. We hope anyway that this text will encourage people to study the convergence problem in this larger setting. 

 \vskip2mm
Today, Alexandrov and Regge's above proposal may be thought as part of Discrete Differential Geometry (see \cite{B-S}) since it relies on a discretisation of a Riemannian manifold and gives an intrinsic interpretation of the curvature as an object in the frame of this discretised new object, addressing also to the question of a nice behavior under approximations. Here, we do not investigate metric extrinsic situations: the interested reader may for instance consult \cite{B-S}, \cite{Mo}. Giving a mathematical basis to \cite{R}, our work,  though restricted to Riemannian geometry, addresses the question: {\em may the today developped Quantum Loop Gravity Theory (see for instance  \cite{Chr0}, \cite{Ki}) be thought as an approximation of classical General Relativity?} This also points to the theory of elasticity and to its discrete version (see \cite{Chr0}, \cite{Chr1}).

\vskip2mm

This work was started by the first author during a long stay in Athens, at the School of Applied Mathematics and Physics of the NTUA (Polytechnio) in 1995. Its initial motivation has roots in exchanges with Lafontaine, at the time Jacques was preparing his Bourbaki lecture on the Cheeger-M\"{u}ller-Schrader theorem (see \cite{La}). Then, the first author benefited of fruitful advices and discussions with Marcel Berger from the middle of the nineties to the beginning of the two-thousands, at the first part of this project: Marcel indicated the importance of the Cayley-Menger determinants (see
\cite{Be1}, section 9.7) and, giving us the references \cite{Be3} and \cite{Va}, led to \cite{Sc} (see also sections 16, 17, 18 pp. 49-54 in \cite{Be5}). 
In the years 2004-5, an unexpected mathematical difficulty came out. 
While making the refinement of a simplicial complex {\it in a special way}, one has to ensure that intersections of a generic curve (of a generic surface) and any refined simplex of codimension $1$ (of codimension $2$) {\it stays of bounded multiplicity in the process} to gain some {\it convergence}. 
Asking Marcel Berger if he knew about such results was replied by an immediate address to the ``forgotten Th\' eor\`eme in G\' eom\' etrie finie'' of Ren\' e Thom  (\cite{Th}, \cite{Be2}). It was the beginning of a... new story !
Part \ref{M-T2} is dedicated to Marcel's memory.
\par
There, we owe to Jean-Jacques Risler his initial coaching and encouragements, when Thom's theorem showed up needing a proof. We dedicate part \ref{M-T1} to his memory.
We are heavily indebted to Alain Chenciner, who took much of his time, helping during two years to write a former complete proof of Thom's theorem with great care and gentleness, and to Marc Chaperon, who became, two years later, the coauthor of a seemingly not well-known but much improved paper, see \cite{C-M}. Part \ref{annexe1033} reflects those efforts.

Though the work had gone far, it would have been dropped if two events didn't prevent this failure.
\par On internet, the first author found a video taken from the nice conference given at the College de France by Snorre Christiansen \cite {Chr0} (on the dark day of the Bataclan (!)),  see his papers on the subject \cite{Chr1} and \cite{Chr2}. The first author understood that the story of approximating Riemannian manifolds by Polyhedra was still interesting people: the subject in those references is the twin problem, approximating Polyhedra by Riemannian manifolds. 
\par Then, after several talks with the first, the second author convinced the first that publishing a text would be valuable.
From a mess of sketches, unfinished drafts and paintings of ideas, the second author, by offering his will to understand, rigorous critics and his mathematical contributions, succeeded in that the efforts of both authors gave, four years later, this strange but seemingly achieved text.

\vskip2mm
We want to thank Sylvain Gallot for comments and contributions on the first part of this text,
 Marina Ville, who attended secret lectures in coffee shops and convinced us trying to express explanations in the language of Shakespeare and Harold Rosenberg for his interest, encouragements and his help in introducing this work in a decent way.
\vfill\eject

\part{Approximating a Riemannian metric and its parallel translation}\label{M-T1}
$$\hbox{\it In memory of Jean-Jacques Risler}$$

\section{Euclidean preliminaries}

\subsection{On how far a linear simplex is from a regular one}\label{2.1}
$ $

We follow Munkres' book \cite{Mu}, see definition 9.9 on page 90, and introduce 
the thickness of a 
linear $n$-simplex $\sigma$ having vertices $p_0,p_1,\dots,p_n$ sitting in a 
Euclidean space $(E,\langle\cdot,\cdot\rangle_E)$ of dimension $\geq n\,$.

\begin{defn} \label{1a}
Denoting by $r(\sigma)$ the radius of the biggest $n$-ball 
centered at the gravity 
center $\tilde p=(p_0+p_1+\dots+p_n)/(n+1)$ and inscribed in $\sigma\,$, and by 
$\delta(\sigma)$ the 
diameter of $\sigma$, the {\it thickness} of $\sigma$
is defined to be the ratio $t(\sigma)=r(\sigma)/\delta(\sigma)$.
Observe that $p_0,p_1,\dots,p_n$ are affinely independent if and only if $t(\sigma)>0$. In this case $\sigma$ is {\em nondegenerate}.
\end{defn}

Though this may be not explicit, the $n$-simplices are often assumed to be nondegenerate.
We shall make a wide use of the following
\begin{lem}\label{1b}
Setting $U_1=p_1-p_0,\dots,U_n=p_n-p_0\,$, one has for any scalars 
$u_1,\dots,u_n$
\begin{equation*}
r(\sigma)\,\sqrt{\sum_{i=1}^n u_i^2}\,\leq\, 
\left\Vert\,\sum_{i=1}^n u_i\,U_i\,\right\Vert_E\,  .
\end{equation*}
\end{lem}

\begin{proof} 
The full ellipsoid 
\begin{equation*}
{\mathcal E}=\left\{\sum_{i=1}^n u_i\,U_i\mid\,
\forall (u_1,\dots,u_n)\in\R^n : \sum_{i=1}^n u_i^2\leq 1 \right\} 
\end{equation*}
contains the 
translated simplex $\sigma-p_0\,$, thus also the ball $B(\tilde p-p_0,r(\sigma))\,$. 
Being an ellipsoid, ${\mathcal E}$ is preserved by $-\hbox{\rm Id}$ : by 
convexity, 
${\mathcal E}$ contains also the ball $B(0,r(\sigma))$. If 
$e_1,\dots,e_n$ is an orthonormal basis of the $n$-vector subspace generated in $E$ by $U_1,\cdots,U_n\,$, one 
derives from this inclusion for any 
$u=\sum_{i=1}^n u_i\,U_i$ 
with $\sum_{i=1}^n u_i^2= 1$ (i. e. for any boundary point $u\in\partial {\mathcal E}$) the inequality
\begin{equation*}
r^2(\sigma)=\sum_{i=1}^n u_i^2\,r^2(\sigma)=\Vert\,\sum_{i=1}^n 
u_i\,(r(\sigma)\,e_i)\,\Vert_E^2\,
\leq\, \Vert\,\sum_{i=1}^n u_i\,U_i\,\Vert_E^2\  ,
\end{equation*}
and we are done.
\end{proof}

To deal later with differentiable images of simplices, a more robust notion, 
closely related, is the openness. 

\begin{defn} \label{openness1} Recall that in a {\em regular} simplex all edges are equal.
$ $
\begin{enumerate}
 \item The {\it openness} of an $n$-simplex $\sigma \subset E$ is the 
ratio $\omega(\sigma)=\hbox{\rm vol}_n(\sigma)/\delta^n(\sigma)$, where 
$E$ is assumed to be of dimension $\geq n\,$.

The isoperimetric inequality for the $n$-simplex \cite[Inequality (1)]{Rez} says 
$\ \omega(\sigma)\leq \omega(\sigma_n)=\vert \sigma_n\vert\,$, where $\sigma_n$ is the $n$-dimensional regular simplex of diameter $1$ and $\vert \sigma_n\vert$ its $n$-volume.

\item In a Riemannian manifold $(M,g)$ of dimension $\geq n\,$,
a {\it differentiable $n$-simplex} is the image of an $n$-simplex $\sigma \subset E$ 
by an 
embedding $T$ from an open $U\subset E$, containing $\sigma$, into $(M,g)$. 
The openness of such a differentiable $n$-simplex $T(\sigma) \subset (M,g)$ is the ratio 

\centerline{$\omega(T(\sigma))=\hbox{\rm vol}_{n,g}(T(\sigma))/\delta_g^n(T(\sigma))\,,$}

\noindent where $\hbox{\rm vol}_{n,g}$ is the 
$n$-dimensional volume in $(M,g)$ and $\delta_g(T(\sigma))$ is the 
$g$-diameter. 
\end{enumerate}
\end{defn}

\begin{rem}\label{R.regular}
The $n$-dimensional regular simplex $\sigma_n$ of diameter $1$  verifies
\begin{equation*} t(\sigma_n)=r(\sigma_n)=\frac{1}{\sqrt{2\,n\,(n+1)}} \ \ \   
\hbox{and} \ \ \ \omega(\sigma_n)=\vert\sigma_n\vert=
\frac{\sqrt{n+1}}{2^{\frac{n}{2}}\,n!}\  .
\end{equation*} 
\end{rem}

\begin{lem}\label{1c}
The thickness and openness of a linear $n$-simplex $\sigma$ 
having $n+1$ 
vertices in a Euclidean space satisfy the following inequalities
\begin{equation}\label{1c1033}
\hskip4mm\frac{n}{n+1}\dfrac{\ \vert\sigma_n\vert^{\frac{n-1}{n}}} 
{\vert\sigma_{ n-1}\vert}\ \omega^{\frac1n}(\sigma)\,\geq\, t(\sigma)\,\geq\,
\frac{n}{n+1}\frac{1}{\vert\sigma_{n-1}\vert}\ \omega(\sigma)\,,
\end{equation}
where $\vert\sigma_l\vert$ stands for the $l$-volume of the regular 
$l$-dimensional simplex of 
diameter $1$ and where equality is achieved on the left if $\sigma$ is isometric 
to $\sigma_n$ up to a 
scaling factor. Thus a bound on $\omega(\sigma)$ gives a bound on $t(\sigma)$ and conversely.
More, $\,t(\sigma)\leq 1/\sqrt{2n(n+1)}$, where equality occurs if and only if $\sigma$ is isometric to $\sigma_n\,$, up to a scaling factor.
\end{lem}

\begin{proof} 
Call $\partial_0\sigma,\partial_1\sigma,\dots,\partial_n\sigma$ 
the $(n-1)$-faces of 
$\sigma$ opposite to the vertices $p_0,p_1,\dots,p_n$ and $H_0,H_1,\dots,H_n$ the hyperplanes generated by those faces. One has $B(\tilde p,r(\sigma))\subset\sigma$ (definition of $r(\sigma)$) and $B(\tilde p,r(\sigma))$ is tangent in its interior to some face $\sigma_n\subset\sigma\,$, thus $r(\sigma)=d_E(\tilde p,\partial_n\sigma)=d_E(\tilde p,H_n)$. With respect to each $\partial_i\sigma\,$, as $d(\tilde p,H_i)=\frac{d(p_i,H_i)}{n+1}\,$,
one has 
\begin{equation}\label{1ca} 
d_E(\tilde p,H_i)\,
\frac{\hbox{\rm vol}_{n-1}(\partial_i\sigma)}{n}=\frac{\hbox{\rm vol}_n(\sigma)}{n+1}\ 
\end{equation}
and (in particular)
\begin{equation}\label{1cb}
\frac{1}{n}\ r(\sigma)\,\hbox{\rm vol}_{n-1}(\partial_n\sigma)=\frac{1}{n+1}\ 
\hbox{\rm vol}_n(\sigma)\,.
\end{equation}

\par
As $d_E(\tilde p,H_i)\geq r(\sigma)\,$,
summing up all equations (\ref{1ca}) gives
\begin{equation*}
\hbox{\rm vol}_n(\sigma)=\frac{1}{n}\ \sum_{i=0}^n
d_E(\tilde p,\partial_i\sigma)\,\hbox{\rm vol}_{n-1}(\partial_i\sigma)\geq
\frac{1}{n}\,r(\sigma)\  
\sum_{i=0}^n\hbox{\rm vol}_{n-1}(\partial_i\sigma)\,,
\end{equation*}
which may be rewritten
\begin{equation}\label{1cc}
\frac{\ (\hbox{\rm vol}_n(\sigma))^{\frac{n-1}{n}}} 
{\hbox{\rm vol}_{n-1}(\partial\sigma)}\geq
\frac{r(\sigma)}{n\,(\hbox{\rm vol}_n(\sigma))^{\frac{1}{n}}}\  .
\end{equation}
The already cited isoperimetric inequality for the $n$-simplex \cite{Rez} tells
\begin{equation}\label{1iso}
\frac{\vert\sigma_n\vert^{\frac{n-1}{n}}}{(n+1)\,\vert\sigma_{n-1}\vert}=
\frac{\ (\hbox{\rm vol}_n(\sigma_n))^{\frac{n-1}{n}}}{\hbox{\rm vol}_{n-1}
(\partial\sigma_n)}\geq\frac{\ (\hbox{\rm vol}_n(\sigma))^{\frac{n-1}{n}}} 
{\hbox{\rm vol}_{n-1}(\partial\sigma)}\  ,
\end{equation}
which gives with (\ref{1cc}) 
\begin{equation*}
\frac{n}{n+1}\dfrac{\ 
\vert\sigma_n\vert^{\frac{n-1}{n}}}{\vert\sigma_{n-1}\vert}\ 
(\hbox{\rm vol}_n(\sigma))^{\frac{1}{n}}\,\geq\, r(\sigma)\,,
\end{equation*}
and the left-hand side is established. Equality holds if and only if equality 
holds in the 
isoperimetric inequality (\ref{1iso}): this occurs if and only if $\sigma$ is 
a 
regular simplex.

On the other hand, one always has $\hbox{\rm vol}_{n-1}
(\partial_n\sigma)\leq \vert\sigma_{n-1}\vert\ (\delta(\sigma))^{n-1}\,$ and as 
$\partial_n\sigma$ 
was assumed to be the closest face to $\tilde p\,$, starting from equation 
(\ref{1cb}), one gets with 
the inequality just quoted
\begin{equation*} 
\frac{\hbox{\rm vol}_n(\sigma)}{n+1}\ \leq\frac{1}{n}\ 
r(\sigma)\,\vert\sigma_{n-1}\vert\ (\delta(\sigma))^{n-1}\  ,
\end{equation*}
which gives the right-hand side inequality to prove. 
The critical left-hand side inequality (\ref{1c1033}) and remark
\ref{R.regular} give the last claim.
\end{proof}

\begin{defn}  Given $t_0>0\,$, denote by ${\mathfrak S}_{t_0}$ the set of all 
linear Euclidean $n$-simplices in 
$\R^{n}$ which are of thickness $\geq t_0\,$.
\end{defn}  

\begin{lem} \label{1fo} 
There exist $t_{0,k}>0$ depending only on $k,n$ and 
$t_0>0$ such that any $k$-face of 
an $n$-simplex in ${\mathfrak S}_{t_0}$ has thickness $\geq t_{0,k}>0\,$.
\end{lem}

\begin{proof} 
By homogeneity, it is enough to prove this claim on all simplices 
of diameter $1$ which 
are in ${\mathfrak S}_{t_0}$: they form (up to isometries) a compact set. If the claim was not 
true, one could find 
a sequence of simplices of diameter $1$ and of thickness $\geq t_0$ converging 
to some limit $n$-simplex 
sharing the same properties, showing at the same time a sequence of $k$-faces 
tending to a 
limit $k$-face of thickness $0\,$. This is a contradiction, for a $k$-simplex having 
thickness $0$ is degenerate, 
and so would be the limit $n$-simplex too.
\end{proof}

\subsection{A tool to study the proximity of Euclidean simplices}\label{2.2}

\begin{lem}\label{1d} Let $\sigma_{\bf v},\,\sigma_{\bf w}$ be simplices 
with vertices 
${\bf v}=\{v_0,v_1,\cdots,v_n\}$ and ${\bf w}=\{w_0,w_1,\cdots,w_n\}$ in $E$ and $F\,$,
$n$-dimensional Euclidean 
spaces. 
Assume there exists $\nu \in \, ]0,1]$ such that
\begin{equation}\label{1da}  
\forall \,i,j=0,1,\cdots,n\hskip8mm
\vert\,\Vert \,v_i-v_j\,\Vert_E-\Vert\,w_i-w_j\,\Vert_F\,\vert\leq\nu\ 
\delta(\sigma_{\bf v})\,. 
\end{equation}

{\rm 1}. Assume $v_0=0$ and $w_0=0$ and consider the linear map ${\mathcal A}$ 
sending $v_i$ to $w_i$ for all $i=0,1\dots,n\,$ ; given 
$\alpha_1,\alpha_2\in E$, call $\beta_1={\mathcal A}(\alpha_1),\,
\beta_2={\mathcal A}(\alpha_2)\in F\,$. One has
\begin{equation*}
\vert\,\langle\alpha_1,\alpha_2\rangle_E-\langle\beta_1,\beta_2\rangle_F\,
\vert\leq
\frac{9\,n\,\nu}{2\,t^2(\sigma_{\bf v})}\,\Vert\,\alpha_1\,\Vert_E\,\Vert\,\alpha_2\,\Vert_E\ ,
\end{equation*}
so that
\begin{equation}\label{1db}
\vert\,\Vert\,\alpha_1\,\Vert_E-\Vert\,\beta_1\,\Vert_F\,\vert\leq
\frac{9\,n\,\nu}{2\,t^2(\sigma_{\bf v})}\,\Vert\,\alpha_1\,\Vert_E\ .
\end{equation}

{\rm 2}. The vertices $v_0$ and $w_0$ may now be $\not=0\,$, one has
\begin{equation*}
\vert r(\sigma_{\bf v})-r(\sigma_{\bf w})\vert\leq \frac{9n\,\nu}{2\,t^2(\sigma_{\bf v})}\,r(\sigma_{\bf v})\,,
\end{equation*}
\begin{equation*}
(1-\nu)\,\delta(\sigma_{\bf v})\leq \delta(\sigma_{\bf w})\leq 
(1+\nu)\,\delta(\sigma_{\bf v})\,,
\end{equation*}
\begin{equation*}
\vert\,t(\sigma_{\bf v})-t(\sigma_{\bf w})\,\vert\leq 
\frac{\nu\,t(\sigma_{\bf v})}{1-\nu}\, 
(1+\frac{9n}{2\,t^2(\sigma_{\bf v})})\, .
\end{equation*}
\end{lem}

\begin{proof} 
First, proving part 1.,
from (\ref{1da}), we get for any $i,j=1,\dots,n$ (use the definition of 
$\delta(\sigma_{\bf v})$ 
and the inequality $\nu(2+\nu)\leq 3\,\nu$)
\begin{equation}\label{1dc}
\begin{split}
\vert\,\Vert v_i\Vert_E^2-\Vert w_i\Vert_F^2\,\vert  \leq 3\,\nu\ 
\delta^2(\sigma_{\bf v}) \,;\\
\vert\,\Vert v_i-v_j\Vert_E^2-\Vert w_i-w_j\Vert_F^2\,\vert \leq 
3\,\nu\ \delta^2(\sigma_{\bf v})\, .
\end{split}
\end{equation}
One gets, writing $\alpha_1=\sum_{i=1}^n\lambda_i\,v_i$ and $\alpha_2=\sum_{i=1}^n\mu_i\,v_i$
\begin{equation*}
\langle\alpha_1,\alpha_2\rangle_E
=\sum_{i,j=1}^n
\,\frac{\lambda_i\,\mu_j}{2}\,\,
\big(\Vert v_i\Vert_E^2+\Vert v_j\Vert_E^2-\Vert v_i-v_j\Vert_E^2\big)\,,
\end{equation*}
thus, using the inequalities (\ref{1dc}) above and Schwarz's 
inequality
\begin{multline*}
\big\vert\langle\alpha_1,\alpha_2\rangle_E-\langle\beta_1,\beta_2\rangle_F\big\vert
= \big\vert\sum_{i,j=1}^n\lambda_i\,\mu_j
\,(\langle v_i,v_j\rangle_E-
\langle w_i,w_j\rangle_F) \big\vert \leq \\
\leq(\sum_{i=1}^n\vert\lambda_i\vert)\,(\sum_{j=1}^n\vert\mu_j\vert)
\,\frac{9\,\nu\,\delta^2(\sigma_{\bf v})}{2}\leq \sqrt{(\sum_{i=1}^n\lambda_i^2)
\,(\sum_{j=1}^n\mu_j^2)}
\ \frac{9n\,\nu\,\delta^2(\sigma_{\bf v})}{2}\  .
\end{multline*}
Finally, inequality $r(\sigma_{\bf 
v})\,\sqrt{\sum_{i=1}^n\lambda_i^2}\leq \Vert \alpha_1\Vert_E$ 
(lemma \ref{1b}) together with definition \ref{1a}
give the 
estimates of part $1$.

Now, proving part 2., the vector $\alpha\,$, pointing from the gravity center in $\sigma_{\bf v}$ to a 
point on the hyperplane supporting the face in 
$\sigma_{\bf v}$ opposite to the vertex $v_i\,$, corresponds by ${\mathcal A}$ to 
the vector $\beta\,$, 
pointing from the gravity center of $\sigma_{\bf w}$ to the corresponding point on 
the hyperplane supporting the face in 
$\sigma_{\bf w}$ opposite to the vertex $w_i\,$. Thus, from the definition of 
$r(\sigma_{\bf v})\,$ and (\ref{1db}), one gets
\begin{equation*}
\vert r(\sigma_{\bf v})-r(\sigma_{\bf w})\vert\leq 
\frac{9n\,\nu}{2\,t^2(\sigma_{\bf v})}\,r(\sigma_{\bf v})\leq 
\frac{9n\,\nu}{2\,t(\sigma_{\bf v})}\,\delta(\sigma_{\bf v})\,,
\end{equation*}
and also
\begin{equation*}
t(\sigma_{\bf v})-\frac{9n\,\nu}{2\,t(\sigma_{\bf v})} \leq 
\frac{r(\sigma_{\bf w})}{\delta(\sigma_{\bf v})} 
\leq t(\sigma_{\bf v})+
\frac{9n\,\nu}{2\,t(\sigma_{\bf v})}\  .
\end{equation*}
From the hypotheses (\ref{1da}) follows directly
\begin{equation*}
(1-\nu)\,\delta(\sigma_{\bf v})\leq\delta(\sigma_{\bf w})
\leq(1+\nu)\,\delta(\sigma_{\bf v})\,,
\end{equation*}
so that 
\begin{equation*}
-\frac{\nu}{(1+\nu)}\,(t(\sigma_{\bf v})+\frac{9n}{2\,t(\sigma_{\bf v})})\leq 
t(\sigma_{\bf w})-t(\sigma_{\bf v})\leq
\frac{\nu}{(1-\nu)}\,(t(\sigma_{\bf v})+
\frac{9n}{2\,t(\sigma_{\bf v})})\ ,
\end{equation*}
which completes the result.
\end{proof}

\subsection{When are $n(n\!+\!1)/2$ positive reals the distances between the vertices 
of a Euclidean $n$-simplex ?}\label{2.3}
$ $

We state the following theorem, elementary but useful, central in our development.
\begin{thm}[Theorem of Schoenberg] \label{2Schoenberg} 
A symmetric matrix having $(n\!+\!1)\times(n\!+\!1)$ entries $\alpha_{i,i}=0$ and 
$\alpha_{i,j}=\alpha_{j,i}\geq 0$ for $i,j=0,1,\dots,n$ and $i\not=j$
can be viewed as giving the distances between two of $n+1$ points in an 
Euclidean space of 
dimension $n\ ($but not $n-1)$ if and only if the matrix with $n\times n$ entries 
$($now $i$ and $j$ run from 
$1$ to $n)$
\begin{equation*}
\frac{1}{2}\ [\alpha_{i,0}^2+\alpha_{j,0}^2-\alpha_{i,j}^2]
\end{equation*}
is positive definite.
\end{thm}
The reader may consult I. J. Schoenberg's article \cite{Sc}. 
\section{Riemannian background}
References on Riemannian geometry in connection with developments in 
this text are \cite{BGM}, \cite{B-K}, \cite{C-E}, \cite{Mi}, \cite{Pe}, \cite{St}, \cite{Spi}, \cite{dC2}.

\subsection{Local properties link Riemannian to Euclidean geometry}\label{EuclRiem}
$  $

Let $(M,g)$ be an open $n$-dimensional Riemannian manifold (with $n\geq2$) and 
$\partial M=\bar M\setminus M$ its boundary, which may be empty.

\begin{defn}\label{D.convex}$ $
 For any $x\in M$ and any positive real $r\,$, denote by $B(x,r)\subset M$ the 
open geodesic ball centered at $x$ of radius $r\,$. Denote by $B(0_x,r)$  
the open ball in the Euclidean tangent space $(T_xM,g_x)$.
A geodesic ball $B(x,r)$ is said to be {\em convex} (in fact strongly convex, see \cite{G-K-M} pages 159ff) if 
\begin{itemize}
 \item the exponential map $\exp_x : T_xM \rightarrow M$ is defined and one-to-one on 
 $B(0_x,r)$;
 
 \item for any two distinct points in $B(x,r)$, there exists in $(M,g)$ a unique minimising geodesic 
arc joining those points which is contained in $B(x,r)$;

\item the above property is true for any ball $B(x',r')\subset B(x,r)$.
\end{itemize}
For any $x\in M\,$, the {\em convexity radius} $r_x$ is defined to be the supremum of the 
nonempty set $\{ r>0 \mid \text{the geodesic ball}  \  B(x,r) \ \text{is convex}\}\,$.
\end{defn}
\begin{rem} \label{strictconv} A {\it special convexity} of the distance function $\varrho:=d(x,\bullet)$ holds in $B(x,r)$: at any $x'\in B(x,r)$, the $g$-Hessian $(Dd\varrho)_{x'}$ is a positive definite symmetric bilinear form on $\{\nabla\varrho\}_{x'}^\perp\subset T_{x'}M\,$ (appendix \ref{strictconv1}). In $B(x,r)$, the special convexity of $\varrho$ implies the {\it strict convexity} of $\varrho^2$.
\end{rem}

\begin{lem}\label{L.constante}
 For any relatively compact open set $W$ verifying $\bar W\subset(M,g)\,$, there exists a real number 
 ${\mathcal R}_W>0$ such that 
 \begin{equation}\label{ineqs}
2\,{\mathcal R}_W\leq \inf_{x\in W_{2{\mathcal R}_W}}r_x \hskip6mm 
\hbox{and}  \hskip6mm d_g(W,\partial M)\geq 4{\mathcal R}_W \  .
\end{equation}
\end{lem}
\begin{rem} {\it Denote} by $W_r$ the set of points in $M$ which are at distance $\leq r$ from $W\,$. The number $2{\mathcal R}_W$ works as a positive lower bound for the 
convexity radius at any point in $W_{2{\mathcal R}_W}\,$.
\end{rem}
\begin{proof}
For any subset $\Omega\subset M\,$, define
$r(\Omega) \!=\! \inf_{x\in \Omega}r_x$ and 
${\mathcal R}_W\!=\!r(W_{r(W)})/2 \!\leq\! r(W)/2\,$. 
Since $r(W_{r(W)})\leq r(W)$ one has $2{\mathcal R}_W \leq  r(W)$, thus
$W_{2{\mathcal R}_W}\subset W_{r(W)}\,$. Therefore
$2\,{\mathcal R}_W= \inf_{x\in W_{r(W)}}r_x \leq \inf_{x\in W_{2\,{\mathcal R}_W}}r_x\,$, 
which is the first stated inequality in (\ref{ineqs}).

Observe that no convex ball in $(M,g)$ contains points in $\partial M$, this proves that 
$d(W_{2{\mathcal R}_W},\partial M)= d(\partial W_{2{\mathcal R}_W}, \partial M)$ 
is $\geq\inf_{x\in \partial 
W_{2{\mathcal R}_W}}r_x\geq\inf_{x\in W_{2{\mathcal R}_W}}r_x\geq 
2{\mathcal R}_W\,$, this gives 
the second inequality in (\ref{ineqs}).
\end{proof}
The three propositions to come are linking local Riemannian and Euclidean geometries through a notion widely popularized by Gromov (see \cite{GLP}):
\begin{defn}\label{quasi-is} Given $C\in[0,1[\,$, a differentiable map $Q$ from a Riemannian manifold $(M_1,g_1)$ to another $(M_2,g_2)$ is $C$-quasi-isometric if, for any $p\in M_1$ and $u\in T_pM_1\,$, one has
$$(1-C)\,\Vert u\Vert_1\leq \Vert Q(u)\Vert_2\leq (1+C)\,\Vert u\Vert_1\  .
$$
More, if $Q$ is a one-to-one and $C$-quasi-isometric function from $(M_1,g_1)$ onto $(M_2,g_2)$, two connected Riemannian manifolds, then $Q$ is a $C$-quasi-isometry, i. e. for any $p$ and $q\in M_1\,$, one also has
$$(1-C)\,d_{g_1}(p,q)\leq d_{g_2}(Q(p),Q(q))\leq (1+C)\,d_{g_1}(p,q)\  .
$$
\end{defn}
\begin{convent}\label{convention1} Given a relatively compact open 
set  $W\,$, with $\bar W\subset(M,g)\,$, 
we now find bounds ${\mathcal R}_l>0$ for $l= 0,1,2$ verifying 
$$0< {\mathcal R}_2\leq {\mathcal R}_1\leq {\mathcal R}_0\leq {\mathcal R}_W\ 
$$ 
and such that, if $\rho$ belongs to $]0, {\mathcal R}_l]\,$, each of the three propositions to come holds. So, one single choice of $\rho>0$ small enough fits with all those three propositions being true.
\end{convent}

\begin{prop} \label{proposition A}
Given $C\in \,]0,1/2]\,$, one can find 
a positive real number 
${\mathcal R}_0={\mathcal R}_0(C) \leq  {\mathcal R}_W\,$, such that, for 
any $x\in W$, one has for any 
$q$ in $B(x,{\mathcal R}_0)$ and $u$ in $B(0_q,2{\mathcal R}_0)$, lifting $g_q$ on $T_u(T_qM)$
\begin{equation}\label{2f} 
\forall v \in T_uT_qM\hskip2mm (1-C)g_q(v,v)\leq (\exp_q^\ast g)_u(v,v)\leq 
(1+C)g_q(v,v)\,,
\end{equation}
thus $\exp_q$ is $C$-quasi-isometric from
$(B(0_q,2{\mathcal R}_0),g_q)$ to $(B(q,2{\mathcal R}_0),g)$
\begin{equation} \label{Eq.norme}
\forall v \in T_uT_qM\hskip2mm (1-C)\,\Vert \,v\,\Vert_q\leq 
\Vert \,d\exp_q(u)(v)\,\Vert_g\leq (1+C)\,\Vert \,v\,\Vert_q\   .
\end{equation}
As the open balls $(B(q,2{\mathcal R}_0),g),(B(0_q,2{\mathcal R}_0),g_q)$ are convex, $exp_q$ is also a metric $C$-quasi-isometry.
Furthermore, for any finite subset of $C\in]0,1/2]\,$, we can choose 
${\mathcal R}_0(C)$ to be non decreasing on this subset.
\end{prop}
\begin{proof} 
 If $u$ is in $TM\,$, call $\pi(u)\!\in M$ the point such that $u\!\in T_{\pi(u)}M\,$. 
The map $\pi :TM\! \rightarrow \!M$ is $C^2$ 
(even $C^\infty$ according to our hypotheses). Given $q$ and $C\in \,]0,1/2]\,$, we know 
there exists 
${\mathcal R}_q\leq {\mathcal R}_W$ such that, for any $u$ in 
$B(0_q,2{\mathcal R}_q)$, assertion (\ref{2f}) 
holds since $d\exp_q(0_q)=\hbox{\rm Id}\,$.

For a given $x\!\in\!M\,$, set 
${\mathcal R}(x)\!=\!\inf_{q\in \overline{B(x, {\mathcal R}_x })} {\mathcal R}_q\,$. 
By compactness of  $\overline{B(x,{\mathcal R}_x)} $ in $ M\,$, by continuity of 
$(q,v)\in \overline{B(x,{\mathcal R}_x )}\times TT_qM\mapsto 
d\exp_q(\pi(v))(v)\in T_{\exp_q(\pi (v))}M\,$, 
we know that ${\mathcal R}(x)$ is $>0\,$. Put ${\mathcal R}_0=\inf_{x\in 
\overline{W}}{\mathcal R}(x)/2\,$, this 
gives the result.
\end{proof}

\begin{prop} \label{proposition B}
 Given $t_1\in \,]0,\frac{1}{\sqrt{2n(n+1)}}]\,$, let  
 ${\mathcal R}_1={\mathcal R}_1(t_1)$ be 
 any positive number  
 satisfying ${\mathcal R}_1\leq {\mathcal R}_0(\frac{t_1^2}{9n})\,$. 
Then,  for any $x\in W$ and 
every 
${\bf p}=\{p_0,p_1,\dots,p_n\}\subset B(x,{\mathcal R}_1)\,$, if the linear 
$n$-simplex $\sigma_x\subset T_xM$ having vertices 
$\exp_x^{-1}p_0$, 
$\exp_x^{-1}p_1,\dots,\exp_x^{-1}p_n\,$, is of $g_x$-thickness 
 $\geq t_1\,$, one can build 
a 
Euclidean $n$-simplex whose edges have lengths the Riemannian distances 
$d(p_i,p_j)$.
\end{prop}

\begin{proof} 
 Introduce the vectors $U_{i,j}:=\exp_x^{-1}p_j-\exp_x^{-1}p_i$ and $U_i:=U_{0,i}$ in $T_xM\,$. 
The $U_i$ build a basis, due to the hypothesis on the $g_x$-thickness of $\sigma_x\,$. 
So, the symmetric matrix 
$D_x$ with entries
\begin{equation*}
D_{x;i,j}=\frac{1}{2}\,\left[\Vert\,U_i\,\Vert^2+\Vert\,U_j\,\Vert^2-\Vert\,U_j-U_i\,
\Vert^2\right] = g_x(U_i,U_j)\,,
\end{equation*}
where $i,j=1,\dots,n\,$, 
is positive definite.

Choose $C\leq \frac{t_1^2}{9n}$ and  
a positive number ${\mathcal R}_1$ smaller than the corresponding 
${\mathcal R}_0(C)>0$ given in proposition \ref{proposition A}.  Let $D$ be the symmetric 
matrix having $n\times n$ entries given by
\begin{equation*}
D_{i,j}=\frac{1}{2}[\,d^2(p_i,p_0)+d^2(p_j,p_0)-d^2(p_i,p_j)]\  \hskip2mm 
\hbox{where}\hskip2mm i,j=1,\dots,n\ .
\end{equation*}
If $D$ is known to be positive definite, one  
gets the conclusion of proposition \ref{proposition B} from the theorem of 
Schoenberg \ref{2Schoenberg} (see \cite{Sc}).
\par
Observe that  the $g_x$ and $g$-distances are achieved in 
$B(0_x,{\mathcal R}_1)$ and $B(x,{\mathcal R}_1)$ by minimizing length of 
arcs contained in those convex balls. Thus,  from the inequality given by (\ref{Eq.norme})
\begin{equation*}
\forall v \in T_uT_xM, \hskip6mm (1-C)\Vert \,v\,\Vert_x\,\leq 
\Vert \,d\exp_x(u)(v)\,\Vert_g\leq (1+C)\Vert \,v\,\Vert_x \  ,
\end{equation*}
one gets by integration (for any $i,j=0,1,\cdots,n$)
\begin{equation*}
(1-C)\,\Vert \,U_{i,j}\,\Vert_x\,\leq\, d(p_i,p_j)\, \leq\, (1+C)\,
\Vert \,U_{i,j}\,\Vert_x\  .  
\end{equation*}
Therefore, for any $i,j=0,1,\dots,n$, we get (since $C\,(2+C)\leq 3C$)
\begin{equation*}
\vert\,d^2(p_i,p_j)-\Vert\,U_{i,j}\,\Vert_x^2\,\vert\leq\,3C\,\delta^2(\sigma_x)\,. 
\end{equation*}
So that, by the above stated polarization
\begin{equation*}
\vert\,D_{i,j}-D_{x;i,j}\,\vert\leq \,\frac{9C\,\delta^2(\sigma_x)}{2}\  .
\end{equation*}
If ${\mathcal U}$ denotes a column matrix having entries $u_1,\dots,u_n\in{\mathbb R}\,$,
one writes
\begin{multline*}
^{t}{\mathcal U}\,D\,{\mathcal U}=\,^{t}{\mathcal U}\,D_x\,
{\mathcal U}+\,^{t}{\mathcal U}\,(D-D_x)\,{\mathcal U}\geq \\
\geq\Vert\,\sum_{i=1}^n u_i\,U_i\,\Vert^2-\sum_{i,j=1}^n
\vert u_i\vert\,\vert u_j\vert\,\vert\,D_{i,j}-D_{x;i,j}\,\vert\geq \\ 
\geq r^2(\sigma_x)\sum_{i=1}^n u_i^2-\frac{9n\,C
\delta^2(\sigma_x)}{2}\sum_{i=1}^n u_i^2\ ,\hskip4mm
\end{multline*}
using lemma \ref{1b} and Schwarz's inequality.

Since 
$t(\sigma_x)=r(\sigma_x)/\delta(\sigma_x)\geq t_1$ by assumption and $C\leq\frac{t_1^2}{9n}\,$, 
get $t^2(\sigma_x)-\frac{9nC}{2} >0\,$, so is 
$D$ positive definite, 
completing the proof.
\end{proof}

\begin{prop} \label{proposition C}\label{calRbnd4}
Given $t_1\in \,]0,\frac{1}{\sqrt{2n(n+1)}}]$ and any positive number 
${\mathcal R}_2={\mathcal R}_2(t_1)$ satisfying 
${\mathcal R}_2 \leq {\mathcal R}_0(\frac{t_1^2}{45n})\,$, 
given any $x\in W$ and  
${\bf p}=\{p_0,p_1,\dots,p_n\}\subset B(x,{\mathcal R}_2)\,$, 
if the linear $n$-simplex having vertices 
$\exp_x^{-1}p_0, 
\exp_x^{-1}p_1,\dots,\exp_x^{-1}p_n$ in $B(0_x,{\mathcal R}_2)$ 
has $g_x$-thickness $\!\geq\! t_1\,$, then, 
for any 
$q\!\in\!B(x,{\mathcal R}_2)$, the linear $n$-simplex $\!\subset \!T_qM$ built on the vertices
$\exp_q^{-1}p_0,\exp_q^{-1}p_1,\dots,\exp_q^{-1}p_n$ has 
$g_q$-thickness $\geq\! t_1/4$.
\end{prop}

\begin{proof} 
 By proposition \ref{proposition A}, given
$C \in \,]0,1/2]\,$ and ${\mathcal R}_2$ smaller than the 
corresponding 
${\mathcal R}_0(C)>0\,$, one already knows for any $x\in W\,$, 
$q\in B(x,{\mathcal R}_2)$, $u\in B(0_x,{\mathcal R}_2)$ 
and $v\in T_uT_xM\,$, $v\not= 0_x$
\begin{equation}\label{Eq.exp}
(1-\frac{2C}{1-C})\,\leq\, 
\dfrac{\Vert\,d(\exp_q^{-1}\circ
\exp_x)(u)(v)\Vert_q}{\Vert\,v\,\Vert_x} \,
\leq\,(1+\frac{2C}{1-C})\,.
\end{equation}
Thus, if $C= t_1^2/(45 n)<1/12$ (recall $t_1<1$), one knows $Q_{x,q}=\exp_q^{-1}\circ\exp_x$ is $1/4$-quasi-isometric from $(T_xM,g_x)$ to
$(T_qM,g_q)$, since $2C/(1-C)\,\leq \,3\,C \leq 1/4\,$. Let $ p_0,\dots,p_n \in B(0_x,{\mathcal R}_2)\,$, 
 set $v_i=\exp_x^{-1} p_i$, $w_i=\exp_q^{-1} p_i\,$, 
thus $Q_{x,q}(v_i)=w_i\,$, $i=0,\dots,n\,$. Integrating (\ref{Eq.exp}), one gets the inequality
\begin{equation*}
(1-\frac{2C}{1-C})\,\Vert\,v_i-v_j\,\Vert_x
\leq \Vert\,w_i-w_j\,\Vert_q\leq (1+\frac{2C}{1-C})\,\Vert\,v_i-v_j\,\Vert_x\ .
\end{equation*}
 Denote by $\sigma_v \subset T_xM$ (by $\sigma_w \subset T_qM$) the linear 
$n$-simplex constructed on the vertices $v_0,\dots,v_n$ (on $w_0,\dots,w_n$).

In view of
$2C/(1-C)\,\leq \,3\,C \leq 1/4\,$, assumption (\ref{1da}) is satisfied.
The second part of lemma \ref{1d} applies with 
 $\nu=3\,C\,$, so that 
\begin{equation*}
\vert\,t(\sigma_{\bf v})-t(\sigma_{\bf w})\,\vert\leq
\frac{\nu\,t(\sigma_{\bf v})}{1-\nu}\, (1+\frac{9n}{2\,t^2(\sigma_{\bf v})})\,,
\end{equation*}
and the conclusion follows if
\begin{equation}\label{garant1}
\frac{\nu\,t(\sigma_{\bf v})}{1-\nu}\, (1+\frac{9n}{2\,t^2(\sigma_{\bf v})})\leq 
\frac{3}{4}\,t(\sigma_{\bf v})\,.
\end{equation}
But $\nu=3C\leq 1/4\,,\, C=t_1^2/(45 n)$ and $t^2(\sigma_{\bf v})\geq t_1^2$ give (\ref{garant1}).
\end{proof}

\subsection{Some bounds arising from variations by geodesics}\label{bdvargeod}
$ $

A real $C\in\,]0,1/2]$ is given and  
${\mathcal R}_0={\mathcal R}_0(C)>0$ is chosen according to 
proposition \ref{proposition A}. By $\rho>0$ we denote some real number 
which is 
smaller than ${\mathcal R}_0\,$.
Thus we know from (\ref{Eq.norme}), that if $p\in W$ and $u\in T_pM\,$, for each 
$s\in\R$ with 
$\Vert\,s u\,\Vert\leq \rho\,$, one has 
\begin{equation}\label{bvga} 
\forall w\in T_{su}T_pM\hskip4mm\Vert\,d\exp_p(su)(w)\,\Vert_{\exp_p (su)}\leq 
(1+C)\,\Vert\,w\,\Vert_p\ ,
\end{equation}
where, for a fixed $p\in M$, all Euclidean tangent fibres $(T_{su}T_pM,\Vert\cdot\Vert_{su})$ are identified, using the parallelism, to the same $(T_pM,\Vert\cdot\Vert_{p})$.
\begin{defn} \label{D.variation}
Consider $c(s,t)=\exp_p(s(u+tw))$, a variation by geodesics 
with $s$ such that 
$\Vert\,s u\,\Vert\leq \rho\,$. Denote by $Y$ the Jacobi field $Y(s)=\frac{\partial c}{\partial t}(s,0)$ along 
the geodesic $c_0(s)=c(s,0)$ and  
by ${\mathcal P}_s$ the parallel translation along the geodesic $c_0$ from 
$c_0(0)=p$ to $c_0(s)\,$. 

Since $Y(s)=d\exp_p(su)(sw)$, inequality (\ref{bvga}) reads 
\begin{equation} \label{bvga1}
\Vert\,Y(s)\,\Vert_{\exp_p (su)}\leq (1+C)\,\Vert\,sw\,\Vert_p\  .
\end{equation}
If $s\not=0\,$, observe one has $\frac{Y(s)}{s}=d\exp_p(su)(w)$, so that 
$\frac{\partial Y}{\partial s}(0)=w\,$.
\par
Throughout this section, the fields $X_1,\dots,X_n$ {\it build}
 a parallel orthonormal frame along $c_0\,$.
\end{defn}
\par
\begin{defn} \label{Riembnd}
Define  ${\mathfrak R}_0 >0$ setting
\begin{equation*}
{\mathfrak R}_0=\sup_{p\in W}\  \  \sup_{u_1,u_2,u_3\in T_pM\setminus\{0\}}\  
\  \frac{\Vert\,R(u_1,u_2) u_3\,\Vert_p}
{\Vert\,u_1\,\Vert_p\,\Vert\,u_2\,\Vert_p\,\Vert\,u_3\,\Vert_p}\  .
\end{equation*}
\end{defn}

\begin{lem}\label{lemma 1st} 
Consider the above variation $c$ given in definition {\rm \ref{D.variation}}. The constant ${\mathcal C}_1={\mathfrak R}_0\,(1+C)/\sqrt{3}\,$ , 
is such that
\begin{equation*}
\Vert\,\frac{\partial Y}{\partial s}(s)-{\mathcal P}_s\,
\frac{\partial Y}{\partial s}(0)\,\Vert\leq
{\mathcal C}_1\,\Vert\,su\,\Vert^2\,\Vert\,w\,\Vert\  .
\end{equation*}
\end{lem}

\begin{proof} 
From the Jacobi equation, one gets
\begin{equation}\label{jacob1}
\langle\frac{\partial Y}{\partial s},X_i\rangle_s -
\langle\frac{\partial Y}{\partial s},X_i\rangle_0=-\int_0^s 
\langle R(Y,\frac{\partial c_0}{\partial s})
\frac{\partial c_0}{\partial s},X_i\rangle_\sigma \,d\sigma\  .
\end{equation}
Observe
\begin{multline*}
\int_0^s \sum_i\langle R(Y,\frac{\partial c_0}{\partial s})
\frac{\partial c_0}{\partial s},X_i\rangle^2_{\mid\sigma} \,d\sigma= \\ 
=\int_0^s \Vert\, R(Y,\frac{\partial c_0}{\partial s})
\frac{\partial c_0}{\partial s}\Vert^2 (\sigma) \,d\sigma
\leq ({\mathfrak R}_0)^2\,\Vert\,
\frac{\partial c_0}{\partial s}\,
\Vert^4\,\int_0^s \Vert\, Y(\sigma)\,\Vert^2\,d\sigma\  .
\end{multline*}
Use the estimate (\ref{bvga1}) and get
\begin{multline}\label{jacob2}
\int_0^s \sum_i\langle R(Y,\frac{\partial c_0}{\partial s})
\frac{\partial c_0}{\partial s},X_i\rangle^2_{\mid\sigma} \,d\sigma 
\leq \frac{({\mathfrak R}_0)^2\,(1+C)^2\,s^3}{3}\Vert\,u\,\Vert^4\,\Vert\, w\,\Vert^2\  .
\end{multline}

Since ${\mathcal P}_s\frac{\partial Y}{\partial s}(0)={\mathcal P}_s w 
=\sum_i \langle w,X_i(0)\rangle X_i (s)\,$, using Schwarz's inequality, putting together (\ref{jacob1}) and 
(\ref{jacob2}), one gets
\begin{multline*}
\Vert\!\sum_i\!\bigg{(}\!\langle\frac{\partial Y}{\partial s},X_i\rangle_s -
\langle w,X_i(0)\rangle\!\bigg{)}\!X_i(s)\Vert^2 
\!=\! \sum_i\!\bigg{(}\!\langle\frac{\partial Y}{\partial s},X_i\rangle_s -
\langle w,X_i(0)\rangle\!\bigg{)}^{\!\!2}\!=\!\\=\!\! \sum_i\!\bigg{(}\! 
\int_0^s \! \langle R(Y,\frac{\partial c_0}{\partial s})
\frac{\partial c_0}{\partial s},X_i\rangle_\sigma \,d\sigma\! \bigg{)}^{\!\!2} 
\! \leq \! s\int_0^s \! \sum_i\langle R(Y,\frac{\partial c_0}{\partial s})
\frac{\partial c_0}{\partial s},X_i\rangle^2_{\mid\sigma} \,d\sigma \leq \\
\leq\frac{({\mathfrak R}_0)^2\,(1+C)^2}{3}\Vert\,s u\,\Vert^4\,\Vert\, w\,\Vert^2\  ,
\end{multline*}
thus the result.
\end{proof}

The next lemma relies on the bound ${\mathfrak R}_0$ on the curvature, it was pointed to us by Sylvestre Gallot, with the nice proof below.
\begin{lem}\label{lemma 2nd}
If $c$ is the variation of definition {\rm \ref{D.variation}}, set ${\mathcal C}_2={\mathfrak R}_0\,2/3\,$. Given $\rho\in]0,\sqrt{\sqrt{3}/(2{\mathfrak R}_0)}]$ (a fortiori $\rho\in]0,\sqrt{1/2{\mathcal C}_2}]$, see remark {\rm\ref{mathcalC1,2}}), one has if $\Vert su\Vert\leq \rho$
\begin{equation*}
\Vert\,{\mathcal P}_s^{-1}\,Y(s)-sw\,\Vert\leq
{\mathcal C}_2\,\Vert\,su\,\Vert^2\,\Vert\,sw\,\Vert\  ,
\end{equation*}
or equivalently
\begin{equation*}
\Vert\,{\mathcal P}_s^{-1}\,d\exp_p(su)(w)-w\,\Vert\leq
{\mathcal C}_2\,\Vert\,su\,\Vert^2\,\Vert\,w\,\Vert\  .
\end{equation*}
\end{lem}
\begin{proof} 
Write the following Taylor expansion with integral remainder
\begin{multline*}
\langle Y(s)-s {\mathcal P}_s(w), X_i(s)\rangle=
\langle Y(s), X_i(s)\rangle-s\langle Y'(0), X_i(0)\rangle=\\=
\int_0^s(s-\sigma) \langle Y''(\sigma), X_i(\sigma)\rangle \,d\sigma=
\\=-\int_0^s(s-\sigma) \langle R(Y(\sigma),\frac{\partial c_0}{\partial s}(\sigma))\frac{\partial c_0}{\partial s}(\sigma), X_i(\sigma)\rangle \,d\sigma\  .
\end{multline*}
Thus, using Schwarz's inequality
\begin{multline*}
\Vert Y(s)-s {\mathcal P}_s(w)\Vert^2=
\sum_i\langle Y(s)-s \,{\mathcal P}_s(w), X_i(s)\rangle^2\leq\\\leq
\int_0^s(s-\sigma)^2 \,d\sigma\ \int_0^s\sum_i\langle R(Y(\sigma),\frac{\partial c_0}{\partial s}(\sigma))\frac{\partial c_0}{\partial s}(\sigma), X_i(\sigma)\rangle^2 \,d\sigma=\\=
\frac{s^3}{3}\,\int_0^s\Vert R(Y(\sigma),\frac{\partial c_0}{\partial s}(\sigma))\frac{\partial c_0}{\partial s}(\sigma)\Vert^2 \,d\sigma\  .
\end{multline*}
Using definition \ref{Riembnd}, we get the estimate
\begin{equation}\label{step1} \Vert Y(s)-s\,{\mathcal P}_s(w)\Vert^2\leq \frac{s^3}{3}\ {\mathfrak R}_0^2
\ \Vert u\Vert^4\int_0^s\Vert Y(\sigma)\Vert^2\,d\sigma\  .
\end{equation}
With the help of the triangle inequality, derive, for any $s$ with $\Vert su\Vert\leq \rho$
\begin{equation*} \Vert Y(s)\Vert \leq s\,\Vert{\mathcal P}_s(w)\Vert+ \frac{s^{2}}{\sqrt{3}}\ {\mathfrak R}_0
\ \Vert u\Vert^2\,\sup_{\sigma\in[0,s]}\Vert Y(\sigma)\Vert\  ,
\end{equation*}
which gives in turn, for any $\sigma\in[0,s]$
\begin{equation*} \sup_{\tau\in[0,\sigma]}\Vert Y(\tau)\Vert\leq \sigma\ \Vert w\Vert +
\frac{\sigma^{2}}{\sqrt{3}}\ {\mathfrak R}_0
\ \Vert u\Vert^2\,\sup_{\tau\in[0,\sigma]}\Vert Y(\tau)\Vert\  ,
\end{equation*}
and finally
\begin{equation} \label{step2} \sup_{\tau\in[0,\sigma]}\Vert Y(\tau)\Vert\leq (\frac{1}{1-\frac{s^{2}}{\sqrt{3}}\ {\mathfrak R}_0\ \Vert u\Vert^2})\  \sigma\ \Vert w\Vert 
\  .
\end{equation}
Bringing the last (\ref{step2}) in (\ref{step1}), one derives
\begin{multline*}\Vert Y(s)-s {\mathcal P}_s(w)\Vert^2\leq \frac{s^{6}}{9}\ {\mathfrak R}_0^2
\ \Vert u\Vert^4\ (\frac{1}{1-\frac{s^{2}}{\sqrt{3}}\ {\mathfrak R}_0\ \Vert u\Vert^2})^2\  \Vert w\Vert^2\leq\\\leq\frac{{\mathfrak R}_0^2}{9}
\ (\frac{1}{1-\frac{\rho^{2}}{\sqrt{3}}\ {\mathfrak R}_0})^2\ \Vert s\,u\Vert^4\  \Vert \,s\,w\Vert^2\  ,
\end{multline*}
which, defining ${\mathcal C}_2=2\,{\mathfrak R}_0/3\,$, gives the expected result
\begin{equation*}\Vert {\mathcal P}_s^{-1}\,Y(s)-s \,w\Vert=\Vert Y(s)-s {\mathcal P}_s(w)\Vert\leq {\mathcal C}_2\ \Vert s\,u\Vert^2\  \Vert s\,w\Vert\  .\qedhere
\end{equation*}
\end{proof}

\begin{rem} \label{mathcalC1,2} Given $C\!\in]0,1/2]$, setting ${\mathcal C}_2\!=\!{2\,{\mathfrak R}_0}/{3}$ implies $\sqrt{1/2\,{\mathcal C}_2}\!\leq\!\sqrt{\!\sqrt{3}/(2{\mathfrak R}_0)}$
(used in lemmas \ref{lemma 2nd}, \ref{lemma 3rd}) and ${\mathcal C}_1\!+\!{\mathcal C}_2\!<\! \sqrt{3}{\mathfrak R}_0$ (used later).
\end{rem}
\begin{lem}\label{lemma 3rd} 
With the hypotheses and notations of lemma {\rm \ref{lemma 2nd}}, given $C\in\,]0,1/2]$ and  
${\mathcal R}_0={\mathcal R}_0(C)>0$ according to 
proposition {\rm\ref{proposition A}},
if $\Vert su\Vert\leq \rho$ with 
$\rho\leq \min( {\mathcal R}_0,\sqrt{1/2\,{\mathcal C}_2})$, one has
\begin{equation*}
\Vert\,((d\exp_p (su))^{-1}\circ{\mathcal P}_s)\,(w)-w\,\Vert\leq
2\,{\mathcal C}_2\,\Vert w \Vert \,\Vert s u\Vert^2\ ,
\end{equation*}
and consequently (with $w_s={\mathcal P}_s w$)
\begin{equation*}
\Vert\,{\mathcal P}_s\circ (d\exp_p (su))^{-1}\   w_s-w_s\,\Vert\leq 
2{\mathcal C}_2\,\Vert w_s \Vert \,\Vert s u\Vert^2\ .
\end{equation*}
\end{lem}

\begin{proof} 
 The inverse of ${\mathcal L}=(d\exp_p (su))^{-1}\circ{\mathcal P}_s$ is 
${\mathcal L}^{-1}={\mathcal P}_s^{-1}\circ d\exp_p(su)$.
 Setting $A={\mathcal L}^{-1}-\hbox{\rm Id}\,$, if (using the operator norm) $\Vert\,A\,\Vert <1\,$, as 
${\mathcal L}-\hbox{\rm Id}=-A\circ ({\rm Id}+A)^{-1}=-A\circ ({\rm Id}-A+A^2-A^3+\dots)\,$, one has
\begin{equation*}
 \Vert\,{\mathcal L}-\hbox{\rm Id}\,\Vert \leq 
 \frac{\Vert\,{\mathcal L}^{-1}-\hbox{\rm Id}\,\Vert}
 {1-\Vert\,{\mathcal L}^{-1}-\hbox{\rm Id}\,\Vert}\ .
\end{equation*}
Remark \ref{mathcalC1,2} 
says $\rho^2\leq \sqrt{3}/(2{\mathcal R}_0)$
and we know $\Vert\,{\mathcal L}^{-1}-\hbox{\rm Id}\,\Vert\leq 
{\mathcal C}_2\Vert s u\Vert^2$ from lemma \ref{lemma 2nd},
 so that  
\begin{equation*}
\Vert\,{\mathcal L}-\hbox{\rm Id}\,\Vert\leq\frac{{\mathcal C}_2\,
\Vert\,s u\,\Vert^2}{1-{\mathcal C}_2\,\Vert\,s u\,\Vert^2}\leq 
2\,{\mathcal C}_2\,\Vert\,s u\,\Vert^2\  .
\end{equation*}
The second inequality follows immediately since 
$\Vert {\mathcal P}_s w\Vert = \Vert w\Vert$ and 
\begin{align*}
 \Vert\,{\mathcal P}_s\circ (d\exp_p (su))^{-1}\   w_s-w_s\,\Vert &=
 \Vert\,{\mathcal P}_s\circ (d\exp_p (su))^{-1}\ {\mathcal P}_s w- 
 {\mathcal P}_s w\,\Vert \\
 &= \Vert\,((d\exp_p (su))^{-1}\circ{\mathcal P}_s)\,(w)-w\,\Vert\ .\qedhere
 \end{align*}
\end{proof}

\subsection{On the difference along a curve between the parallel translations  
defined by two Riemannian metrics}\label{diffparll}
$ $

We consider two quasi-isometric metrics $g$ and $\hat g$ on $M\,$, i.e. there exists 
$C_1\in [0,1[\,$ such that
\begin{equation}\label{Eq.quasi}
(1-C_1)\,g\leq \hat g\leq (1+C_1)\,g\  .
\end{equation}
Observe that  equation (\ref{Eq.quasi}) is 
equivalent to 
$\vert\Vert\cdot\Vert_g^2-\Vert\cdot\Vert_{\hat g}^2\vert \leq 
C_1\,\Vert\cdot\Vert_g^2\,$, 
which in turn implies
$\vert\Vert\cdot\Vert_g-\Vert\cdot\Vert_{\hat g}\vert \leq C_1\,\Vert\cdot\Vert_g\,$.

Call $h=h^{g,\hat g}$ the {\it difference connection tensor} $h=\hat\nabla-\nabla\,$, built by the difference of the two 
Levi-Civita connections associated with the two Riemannian metrics $\hat g$ and $g$ 
on the $n$-manifold $M$ (see for instance \cite{MeD}). Thus 
$h(U,V)= \hat \nabla_U V- \nabla_U V$ for any vector fields $U$ and $V\,$.
 
If $\gamma : [0, 1]\rightarrow M$ is a differentiable curve starting at $p\,$, 
denote by ${\mathcal P}_t$ and $\hat{\mathcal P}_t$ the $g$ and $\hat 
g$-parallel translations along $\gamma$ 
from $p$ to $\gamma(t)$.

\begin{defn} \label{hbnd} 
If $X_i, i=1,\dots,n$ is a $g$-parallel orthonormal frame along $\gamma\,$, write 
$h_i=\langle h,X_i\rangle_g$ 
and define a $C^0$ $g$-norm of $h$ along $\gamma$
\begin{equation*}
H=\sqrt{\sum_i\Vert\,h_i\,\Vert^2}\  \  \hbox{where}\  \   
\Vert\,h_i\,\Vert=\sup_{t\in[0,1]}\ \sup_{u,v\in T_{\gamma(t)}M\setminus\{0\}}\ 
\frac{\vert\,h_i(u,v)\,\vert}{\Vert\,u\,\Vert_g\,\Vert\,v\,\Vert_g}\  .
\end{equation*}
\end{defn}
\begin{prop}\label{proposition 1st} One has the following control on
the difference between the 
two parallel translations along $\gamma$ 
\begin{equation*}
\forall t\in[0,1]\hskip1cm\Vert\,\hat{\mathcal P}_t-
{\mathcal P}_t\,\Vert_{p,\gamma(t)}\leq (\frac{1+C_1}{1-C_1})\ H\ 
\hbox{\rm length}_g\,(\gamma_{\mid [0,t]})\, .
\end{equation*}
\end{prop}

\begin{proof} 
Equation (\ref{Eq.quasi}) implies, for any $w\in T_{\gamma(t)}M\,$, both 
\begin{equation} \label{quasiso2}
\begin{split}
 (1-C_1)\,\Vert\,w\,\Vert_g &\leq \Vert\,w\,\Vert_{\hat g} \leq 
(1+C_1)\,\Vert\,w\,\Vert_g\  ,\\
\frac{1}{1+C_1}\,\Vert\,w\,\Vert_{\hat g} &\leq \Vert\,w\,\Vert_g\leq 
\frac{1}{1-C_1}\,\Vert\,w\,\Vert_{\hat g}\  .
\end{split}
\end{equation}

From the last line of (\ref{quasiso2}), we get, for any $v\in T_{\gamma(0)}M\,$, denoting by $\hat{\mathcal P}_{\cdot}(v)$ the field $t\mapsto \hat{\mathcal P}_t(v)$ along $\gamma$
\begin{equation*}
\Vert\,\hat{\mathcal P}_{\cdot}(v)\,\Vert_g
\leq \frac{1}{1-C_1}\Vert\,\hat{\mathcal P}_{\cdot}(v)\,\Vert_{\hat g}\ .
\end{equation*}
Since  
$\Vert\,\hat{\mathcal P}_{\cdot}(v)\,\Vert_{\hat g}=
\Vert\,v\,\Vert_{\hat g}\,$, 
using the first line of (\ref{quasiso2}) we get 
\begin{equation}
\label{quasiso3}\Vert\,\hat{\mathcal P}_{\cdot}(v)\,\Vert_g
\leq (\frac{1+C_1}{1-C_1})\Vert\,v\,\Vert_g\  .
\end{equation}
The frame $X_i(\cdot), i=1,\dots,n$ is $g$-parallel along $\gamma\,$, 
thus
\begin{equation*}
\frac{d}{dt}\!\bigg{(}\!\!\langle (\hat{\mathcal P}_\cdot-
{\mathcal P}_\cdot)(v),X_i\rangle_g\!\!\bigg{)}_{\!\!\mid t}\!\!=\!
\langle (\nabla_{\frac{d\gamma}{dt}}
\hat{\mathcal P}_\cdot(v))_t,X_i\rangle_g\!
= \!-\langle h(\frac{d\gamma}{dt}(t),\!\hat{\mathcal P}_t(\!v\!)),\!X_i\rangle_g\,.
\end{equation*}
Definition \ref{hbnd} reads 
\begin{equation*}
\Big\vert\langle 
h(\frac{\frac{d\gamma}{dt}(t)}{\Vert\frac{d\gamma}{dt}(t)\Vert_g},
\frac{\hat{\mathcal P}_t(v)}{\Vert \hat{\mathcal P}_t(v)\Vert_g}),
X_i\rangle_g\Big\vert\leq \Vert\,h_i\,\Vert\  ,
\end{equation*}
so that with (\ref{quasiso3}), one gets for any $v$
\begin{multline}\label{quasiso4}
\Big\vert \frac{d}{dt}\,\bigg{(}\langle (\hat{\mathcal P}_\cdot-
{\mathcal P}_\cdot)(v),X_i\rangle_g\bigg{)}_{\mid t}\Big\vert 
\leq
(\frac{1+C_1}{1-C_1})\ \Vert\,h_i\,\Vert \ 
\Vert\frac{d\gamma}{dt}(t)\Vert_g\,\Vert\,v\,\Vert_g\  .
\end{multline}
Then, one computes and gets the result claimed, bringing in (\ref{quasiso4})
\begin{equation*}
\begin{split}
 \Vert\, (\hat{\mathcal P}_t-{\mathcal P}_t)(v)\,\Vert^2 &
 =\sum_{i}\bigg{(} \int_0^t(\frac{d}{ds}\,\langle (\hat{\mathcal P}_\cdot-
{\mathcal P}_\cdot)(v),X_i\rangle_g)_s\  ds\bigg{)}^2 \leq\\
&\leq
(\frac{1+C_1}{1-C_1})^2 \big(\sum_i\  \Vert\,h_i\,\Vert^2 \big)\  
\hbox{\rm length}^2_g(\gamma_{\mid [0,t]})\,\Vert\,v\,\Vert_g^2 = \\
&= (\frac{1+C_1}{1-C_1})^2\ H^2\  
\hbox{\rm length}^2_g\,(\gamma_{\mid [0,t]})\ \Vert\,v\,\Vert_g^2\  .\qedhere
\end{split} 
\end{equation*}
\end{proof}

\subsection{Riemannian barycentric simplices}\label{rbs} 
$ $

Here we present the elementary piece of a polyhedral 
approximation of a Riemannian manifold.
See related developments in \cite{B-K}.

By $B\!\subset \!M$ denote an open convex ball of $(M,g)$, 
see definition \ref{D.convex} ($B$ is in fact a strongly convex ball in the sense of \cite{G-K-M}), verifying $\ov B \!\subset \!M\,$.

\begin{prop} \label{rbs1} 
Consider $n+1$ points, $p_0,\dots,p_n\,$, in the convex ball $B$ of $(M,g)$. 
If $\lambda_0,\cdots,\lambda_n$ are reals $\geq 0$ 
summing up to 
$1\,$, the function $\sum_{i=0}^n\lambda_i\,d^2(p_i,\cdot)$ defined on the compact 
set 
$\overline B$ achieves its minimum 
at a unique point inside $B\,$.
\end{prop}

\begin{proof} 
Assume the function achieves its minimum at $x\!\in\partial B\,$. Let $n_x$ be the unit normal pointing inside $\partial B$ at 
$x\,$. Due to the 
convexity of $B$, for each $i=0,1,\cdots,n\,$, the angle between the vectors 
$-\exp_x^{-1}p_i$ 
and $n_x$ is in $\,]\pi/2,\pi]\,$ (the gradient at $x$ of $\frac{1}{2}\,d^2(p_i,\cdot)$ is $-\exp_x^{-1}p_i\,$, see \cite{B-K}, proposition 6.4.6). The first variation formula of the length 
implies 
that all distances $d(p_i,x)$ are strictly decreased by pushing $x$ to $x'$ along a curve tangent to $n_x$ 
at $x\,$. So $x$ cannot be the minimum, contradiction. 

The uniqueness follows from the strict convexity of 
$\sum_{i=0}^n\lambda_i\,d^2(p_i,\cdot)$, 
which is a sum of strictly convex functions on $\overline{B}\,$, see remark \ref{strictconv}. Actually, along a 
geodesic $c\,$, the 
Hessian bilinear form of a function $f$ satisfies
\begin{equation*}
(\nabla df) (\frac{\partial c}{\partial t},\frac{\partial c}{\partial t})=
\frac{\partial^2 f\circ c}{\partial t^2}\   . \qedhere
\end{equation*}
\end{proof}
\begin{rem} \label{minind} The minimum $x$ of $f=\sum_{i=0}^n\lambda_i\,d^2(p_i,\cdot)$ for fixed $\lambda_i\geq0,\,i=0,\cdots,n$ does not depend 
on the convex ball $B$ containing ${\bf p}=p_0,p_1,\dots,p_n$ since $f$ is stricly convex and uniquely defined by ${\bf p}\,$. 
\end{rem}
\begin{defn} \label{D.spread.simplex}
$ $
\begin{enumerate}
\item\label{spread12}  
A set of $n+1$ points 
 ${\bf p} =\{p_0,\cdots,p_n\}\subset B$
 is said to be 
 {\em spread} if, for any $x\in B\,$, the linear simplex of $T_xM$ having vertices 
$\exp_x^{-1}p_0,\cdots,\exp_x^{-1}p_n$ is non 
degenerate.

\item \label{spread12bise}Denote by $\Delta$ the set 
\begin{equation*}
\Delta=\{\lambda=(\lambda_0,\cdots,\lambda_n)
\in\R^{n+1}\mid 
\lambda_0,\cdots,\lambda_n\geq 0 \hskip1mm\hbox{and}\hskip1mm\sum_i\lambda_i=1 
\}\ .
\end{equation*}
If ${\bf p} =\{p_0,\cdots,p_n\}\subset B$  is a set, define
\begin{equation*} 
\hat\tau=\hat\tau_{\bf p}=\{x\in B\mid \exists 
\lambda\in\Delta \hskip2mm\hbox{with}\hskip2mm\sum_i\lambda_i\,d^2(p_i,\cdot) 
\hskip2mm\hbox{is minimum at}\, x\}\  .
\end{equation*}
If ${\bf p}$ is {\it spread}, we say that $\hat\tau=\hat\tau_{\bf p}$ is the
{\em Riemannian barycentric simplex} built on ${\bf p}\,$.

  \item Define the geometric map ${\mathscr S}$, called {\it Status}, from 
$B$ to $\R^{n+1}$ by setting 
\begin{equation*} 
{\mathscr S} : x \in B \longmapsto \sum_i\,d^2(p_i,x)\,e_i\in \R^{n+1}
\end{equation*}
where $e_i, i=0,\dots,n$ is the canonical basis of $\R^{n+1}\,$. 

Consequently, if $\psi$ is an isometry between $(M,g)$ and $(M',g')$, writing ${\mathscr S}^\psi$ 
the Status map 
defined by the set of points $\psi({\bf p})$ (which is also spread in a convex 
ball), one has for any 
$x\in M$ making sense
\begin{equation*}
{\mathscr S}^\psi(\psi(x))={\mathscr S}(x)\,.
\end{equation*}

In the sequel, we often write 
$d_i(\cdot)$ or $d_i$ short for $d(p_i,\cdot)$.

\end{enumerate}
\end{defn}

{\it Instead of defining a spread set of points ${\bf p}=\{p_0,p_1,\dots,p_n\}$ as above}, we might have focused on sets ${\bf p}$ such that,
for any $q$ in $\hat\tau$,
the linear 
simplex with vertices $\exp_
q^{-1}p_0,\dots,\exp_q^{-1}p_n$ is $n$-dimensional (this is slightly more general and independent of $B$). Yet, this is close and the given definition of spread set is natural in our exposition.

We now prove that such a Riemannian barycentric simplex $\hat\tau_{\bf p}$ is embedded 
in $B$ and 
diffeomorphic to {\em the standard regular Euclidean $n$-simplex}
$\Sigma_n$, which is, up to 
isometry, the simplex generated by the vertices $e_0,e_1,\dots,e_n\,$, the canonical 
orthonormal basis in $\R^{n+1}\,$.

\begin{lem} \label{rbs4}
Given $\lambda$ in $\Delta\,$, {\em then} $\upsilon\!=\!\sum_i\lambda_i\,e_i$ is 
orthogonal to 
${\mathscr S}(B)$ at ${\mathscr S}(x)$ {\em if and only if}  $x$ 
is in $\hat\tau$ and
is the unique minimum of $\sum_i\lambda_i\,d_i^2\,$. 
\end{lem}
\begin{proof}
If $\upsilon=\sum_i\lambda_i\,e_i\,$, one has
\begin{equation*}
\sum_i\lambda_i\,d_i^2(\cdot)=\langle{\mathscr S}(\cdot),\upsilon \rangle
\hskip2mm\hbox{and}\hskip2mm\sum_i\lambda_i\,d\,d_i^2(x)(\cdot)=
d\langle {\mathscr S},\upsilon \rangle(x)(\cdot) \,.
\end{equation*}
Thus, if  $x \in \hat\tau$ is the minimum 
of 
$\sum_i\lambda_i\,d_i^2\,$, then $x$ is a critical point of 
$\langle{\mathscr S}(\cdot),\upsilon \rangle\,$, and this 
implies that $\upsilon$ is orthogonal to ${\mathscr S}(B)$ at ${\mathscr S}(x)$.

Conversely, if $\upsilon$ is orthogonal to ${\mathscr S}(B)$ at ${\mathscr S}(x)\,$, then 
$x$ is a critical point of $\langle{\mathscr S}(\cdot),\upsilon \rangle\,$. 
Since $\lambda \in  \Delta\,$, the function $\langle{\mathscr S}(\cdot),\upsilon \rangle$ 
is strictly convex in $B$ with $x$ as unique minimum (proposition 
\ref{rbs1}), 
and $x$ sits in $\hat\tau$ (definition of $\hat\tau$).
\end{proof}

\begin{cor} \label{remarquable} The map ${\mathscr S}$ is one-to-one from $\hat\tau$ onto the set of ``positive directions'' $\{v=\sum_{i=0}^n\lambda_i\,e_i \mid \lambda\in\Delta\}\,$.
\end{cor}
\begin{proof} For any $v=\sum_{i=0}^n\lambda_i\,e_i$ with $\lambda\in\Delta$ exists a {\it unique} point $x$ (which belongs to $\hat\tau$) such that $T_{{\mathscr S}(x)}{\mathscr S}(B)$ is orthogonal to $v\,$.
\end{proof}

\begin{prop}\label{rbs5} 
There exists an open set\, ${\mathcal U}\subset B\,$,  
containing $\hat\tau$ such that 
${\mathscr S} : {\mathcal U} \rightarrow \R^{n+1}$ is an embedding having a 
range  ${\mathscr S}({\mathcal U})$ which is a locally strictly convex embedded 
hypersurface in 
$\R^{n+1}$, 
meaning its second fundamental form is definite (see \cite{K-NII}, page {\rm40}, proposition {\rm5.5}). 
\end{prop}

\begin{proof}
First, Status is an immersion from $B$ into $\R^{n+1}\,$. In fact, its 
differential at any $x\in B$ is
\begin{equation*}
d {\mathscr S} (x)(\cdot)=\sum_i\,d d_i^2(x)(\cdot)\,e_i=-2\,\sum_i\,
\langle\exp_x^{-1}p_i,\cdot\rangle\ e_i\  ,
\end{equation*}
so it has kernel $\equiv\{0\}$ since ${\bf p}$ is spread. Let $x$ be a point in 
$\hat\tau$ and $U$ a small open set containing $x$ such that ${\mathscr S}(U)$ 
is embedded. Denote by $\upsilon$ 
a unit vector which is orthogonal to ${\mathscr S}(U)$ at ${\mathscr S}(x)$. 
The second fundamental form $l_\upsilon$ of this piece of embedded hypersurface 
in $\R^{n+1}$ is given by (see \cite{Mi}, pages 32ff.)
\begin{equation*}
\forall u\in T_xB \hskip1cm l_\upsilon(u,u)=\sum_i\,
\frac{\partial^2d_i^2}{\partial u^2}(x)\  \langle e_i,\upsilon \rangle\  .
\end{equation*}
Actually, for the functions $\langle{\mathscr S}(\cdot),\upsilon \rangle$ under concern, writing second order partial 
derivatives 
(at the critical point $x$ for $\langle{\mathscr S}(\cdot),\upsilon \rangle$) gives the 
Hessian bilinear form at $x\,$. 
The local strict convexity of ${\mathscr S}(U)$ follows now from the strict convexity of
$\sum_i\lambda_i\,d_i^2$ for $\lambda\in\Delta\,$. This is an open condition and 
it 
remains true in an open set containing $\hat\tau\,$.

Suppose that there exist two distinct points $x\in\hat\tau$ and $x'\in B$ 
such that 
${\mathscr S}(x)={\mathscr S}(x')$, this implies $d(p_i,x)=d(p_i,x')$ for 
$i=0,1,\cdots,n\,$. So, there exists a 
$\lambda\in\Delta$ such that $\sum_i\lambda_i\,d_i^2$ achieves its minimum at 
two distinct 
points in $B$, which is impossible by proposition \ref{rbs1}. Thus Status is an 
injective immersion from 
the compact set $\hat\tau$ and it is as well an embedding from an open set
${\mathcal U}\subset B$ containing $\hat\tau$ into $\R^{n+1}\,$.
\end{proof}
\begin{rem} One should refer to Apollonius in the case of Euclidean metrics 
(see \cite{Be1} 9.7.6) 
\begin{equation*}
x=\sum_{i=0}^n\lambda_i\,p_i \hskip2mm\hbox{is the minimum of}\hskip2mm 
\sum_{i=0}^n\lambda_i\,\Vert\,p_i-\cdot\,\Vert^2\ .
\end{equation*}
\end{rem}

Linked to Status, there is a kind of Legendre transform belonging to the 
range of 
convex duality (see \cite{Ar1} and \cite{Ar2}), we now present. 
\begin{defn} 
$ $

\begin{enumerate}
 \item Denote by $\Pi$ the one-to-one map which sends by central 
projection (from $0_{\R^{n+1}}$) a point in the open ``Northern Hemisphere'' ${\bf S}_+^n$ (``having North Pole'' ${\mathcal N}=\frac{e_0+\cdots +e_n}{\sqrt{n+1}}$ on the unit standard sphere 
${\bf S}^n\subset\R^{n+1}$) to a point on the affine 
hyperplane
\begin{equation*}
{\mathcal H}=\{\lambda=(\lambda_0,\cdots,\lambda_n)\in\R^{n+1}\mid 
\sum_{i=0}^n\lambda_i=1\} \  .
\end{equation*}

\item We know by proposition (\ref{rbs5}) that  
${\mathscr S}({\mathcal U})$ is embedded. Therefore, the Gauss map, denoted by $G\,$,  
 is well defined everywhere on 
the whole of ${\mathscr S}({\mathcal U})$. 

\item Recall that  $\Sigma_n$ is  the linear $n$-simplex having vertices the points $e_i$ 
in ${\mathcal H}\subset\R^{n+1}\,$. In fact, $\Sigma_n$ 
is rescaled by a factor $\sqrt{2}$ from the regular $n$-simplex $\sigma_n$ of diameter 
$1\,$.

\item The set $\hat\tau$ was defined in definition \ref{D.spread.simplex} (\ref{spread12bise}). We denote by ${\mathscr B}$ the composed mapping, called the {\it barycentric 
coordinates}
\begin{equation*}
{\mathscr B}=\Pi\circ G\circ {\mathscr S} \ .
\end{equation*}
By lemma \ref{rbs4}  and corollary \ref{remarquable}, for any $x\in\hat\tau\,$,
there is a unique $\lambda\in\Delta$ such that $x$ is the unique minimum of $\sum_i\lambda_i\,d_i^2\,$, with $v=\sum_i\lambda_i\,e_i$ a ``positive direction'' orthogonal to ${\mathscr S}(U)$ at ${\mathscr S}(x)$. One has $G\circ{\mathscr S}(x)=v/\Vert v\Vert\in {\bf S}_+^n$ and (since $\sum_i\lambda_i=1$) 
$$\forall x\in\hat\tau\ \ \ \ \ \ \ \ {\mathscr B}(x)=\Pi (G\circ {\mathscr S}(x))=\lambda=(\lambda_0,\cdots,\lambda_n)\in\Delta\subset{\mathcal H}\ .
$$

More, observe that the barycentric coordinates are not necessarily defined on the whole of 
${\mathcal U}$ because $G\circ {\mathscr S}$ has no reason to take all its 
values in ${\bf S}_+^n\,$. But we can 
restrict ${\mathscr S}$ to a connected open subset of 
$(G\circ{\mathscr S})^{-1} {\bf S}_+^n$ 
in ${\mathcal U}$ containing $\hat\tau$, its image by $G\circ{\mathscr S}$ is 
contained in ${\bf S}_+^n\,$, and ${\mathscr B}$ is then well 
defined.

\end{enumerate}
\end{defn}

\begin{thm}\label{rbs0}
The barycentric coordinates build a diffeomorphism from an 
open 
set $U\subset{\mathcal U}\,$, with
$\hat\tau\subset U$ (see definition {\rm{\ref{D.spread.simplex}} ({\ref{spread12bise}})} for $\hat\tau$), to an open image ${\mathscr B}(U)\subset {\mathcal H}$ 
containing $\Sigma_n\,$. More, one has ${\mathscr B}(\hat\tau)=\Sigma_n\,$.
\end{thm} 

\begin{proof} 
First we deduce from proposition \ref{rbs5} that there exists an open subset 
$U_1\,$, with $\hat\tau\subset U_1\subset{\mathcal U}\,$, such that 
$G$ is a diffeomorphism from ${\mathscr S}(U_1)$ 
to an open subset of ${\bf S}^n\,$.

Set $U=U_1\cap (G\circ{\mathscr S})^{-1}{\bf S}^n_+\,$, then the 
composed 
map $\Pi\circ G$ is a diffeomorphism well defined from 
${\mathscr S}(U)\subset{\mathscr S}(U_1)$ to an open set 
${\mathscr B}(U)\subset{\mathcal H}\,$.
Lemma \ref{rbs4}  and corollary \ref{remarquable} ensure the 
inclusion $\hat\tau\subset U$. The map ${\mathscr B}$ is thus a diffeomorphism 
from $U$ to ${\mathscr B}(U)\,$, because ${\mathscr S}$ is a 
diffeomorphism from $U$ to its image by proposition \ref{rbs5}. More, the 
set 
$\hat\tau\subset U$ is sent on $\Sigma_n$ by construction, this was the purpose of defining ${\mathscr S}$ and ${\mathscr B}$ by exploiting the geometric context emphasized in corollary \ref{remarquable}.
\end{proof}

\begin{cor}\label{rbs10} 
A set of $n+1$ points ${\bf p}\,$, which is spread in a convex 
ball $B$ of a Riemannian manifold $(M,g)$, generates an embedded differentiable 
$n$-simplex $\subset B\,$, the {\em Riemannian barycentric 
simplex} $\hat\tau=\hat\tau_{\bf p}$ constructed on the set of ``vertices'' 
${\bf p}$. A $k$-face of $\hat\tau$ is defined as a subset of $\hat\tau$, e. g. 
the $k$-face supported by $p_0,p_1,\cdots,p_
k$ is the set
\begin{gather*}
\{x\in \hat\tau\mid \exists \lambda\in\Delta_{0,1,\cdots,k}  \hskip1mm 
\hbox{with}\hskip1mm \sum_{i=0}^n\lambda_i\,d_i^2 \hskip1mm 
\hbox{is minimum at}\hskip1mm x\}\  ,\\
\hbox{where}\hskip4mm \Delta_{0,1,\cdots,k}=
\{\lambda\in\Delta\mid \lambda_{k+1}=\cdots=\lambda_n=0\}\  ,
\end{gather*}
thus the definition of this $k$-face depends only upon $p_0,p_1,\cdots,p_k$ and not upon
$p_{k+1},\cdots,p_n\,$.
\end{cor}

\begin{proof} 
This is just rephrasing the essence of the previous results.
\end{proof}

\begin{prop} \label{edges} 
The edges of Riemannian barycentric simplices are (up to a constant factor) natural arc-length parametrized geodesics .
\end{prop}

\begin{proof} 
As before $d_i=d(p_i,\cdot)$ for $ i=0,\dots,n\,$. 
Choose two vertices $p_0$ and $p_1$ and $\lambda\in\,]0,1[\,$. Let $x$ be the 
point 
where the function $(1-\lambda)d_0^2+\lambda d_1^2$ achieves its minimum and 
set $\xi=d(x,p_0)$. If $c(t)
$ is the arc-length parametrized geodesic from $p_0$ to $p_1\,$, introduce the 
point $c(\xi)$. One has
\begin{equation*} 
(1-\lambda)d^2(p_0,c(\xi))+\lambda d^2(p_1,c(\xi))\geq
(1-\lambda)d^2(p_0,x)+\lambda d^2(p_1,x)\,,
\end{equation*}
but as $d(p_0,x)=d(p_0,c(\xi))=\xi\,$, one infers (from the above inequality) 
$\lambda=0$ (excluded) or $d(p_1,c(\xi))\geq d(p_1,x)$, which implies 
$c(\xi)=x\,$, by uniqueness of the geodesic arc from $p_
0$ to $p_1$ in the convex ball $B\,$. So, the edge of $\hat\tau$ between $p_0$ 
and 
$p_1$ is geometrically identical with the image of the geodesic from $p_0$ to 
$p_1\,$.
At the critical point $x$ of $(1-\lambda)d_0^2+\lambda d_1^2\,$, compute
\begin{equation*}
(1-\lambda)\,\langle \exp_{c(\xi)}^{-1}p_0,\frac{dc}{dt}(\xi)\rangle+
\lambda\,\langle \exp_{c(\xi)}^{-1}p_1,\frac{dc}{dt}(\xi)\rangle=0\  ,
\end{equation*}
which implies in turn that 
$(1-\lambda)\,d(p_0,c(\xi))=\lambda\,d(p_1,c(\xi))$, and thus
\begin{equation*}
 (1-\lambda)\,\xi=\lambda\,(d(p_0,p_1)-d(p_0,c(\xi))) = \lambda\,(d(p_0,p_1)-\xi)\, ,
\end{equation*}
so $\lambda=\xi/d(p_0,p_1)$. The edge is parametrized proportionally to the 
arc of length, as expected. 
\end{proof}

\begin{rem} 
The $2$-face supported by $p_0,p_1,p_2$ is not totally geodesic in a 
general Riemannian manifold. The common situation is that totally geodesic 
submanifolds of dimension $\not=0,1$ or $n$ do not necessarily  exist in Riemannian manifolds of dimension 
$n\geq3\,$.  
\end{rem}

\section{Rough Riemannian estimates through the magnifying glass}\label{4}
$ $

According to Gromov, the ``magnifying glass'' involved in this section goes back to 
Riemann and Cartan, see \cite{GLP} proposition 3.15. We develop this tool for two reasons. In general, as a geometrical understanding of approximation, focusing from Riemannian towards Euclidean geometry, it is behind a construction to come up later, which is the cornerstone in the proof of our main theorem. More, {\it making a kind of an exercise in style}, it allows us to state rough estimates on the Hessian of the squared distance to a given point near this point (this is classical, see \cite{B-K}) and, {\it almost together}, on its covariant derivative.
\par
The manifold $(M,g)$ is $C^{k+1}$
(the metric $C^k$), with $k$ to be specified.

\begin{defn} $ $

\begin{enumerate}
\item 
The mapping $\pi\times\exp$ below is $C^{k-1}$
\begin{equation*}
\pi\times\exp:v\in TM\longmapsto (\pi(v),c_v(1))\in M\times M\ ,
\end{equation*}
here $c_v(t)$ is the geodesic verifying 
$\frac{dc}{dt}(0)=v$ and 
$\pi(v)\in M$ is such as $v\in T_{\pi(v)}M\,$.
We restrict $\pi\times\exp$ to an open neighborhood $U$ of the zero-section in 
$TM\,$, such that
\begin{enumerate}
 \item[(i)] the map $\pi\times\exp$ is invertible from $U$ to its image ;
 
 \item[(ii)] the open $U$ is fit to verify the equality
$\pi\times\exp(U)= B\times B\,$, where $B=B(x,\rho)$ is a convex ball centered at 
$x\in W\,$, so  $\rho\leq{\mathcal R}_W\,$, see lemma \ref{L.constante}.
\end{enumerate}

\item Given a function $f$ on $(M,g)$, write ${}^\sharp\hbox{\rm hess}_q f=(\nabla \grad f)_q$ for the $g$-paired $(1,1)$-tensor 
to the $(0,2)$-tensor $\hbox{\rm hess}_q f=\nabla df$
\begin{equation*}
\forall u,v \in T_qM\hskip4mm\ \  \hbox{\rm hess}_q f(u, v)=
\langle{}^\sharp\hbox{\rm hess}_q 
f(u), v\rangle_g=\langle\nabla_{\!\!u}\, \grad f, v\rangle_g
\ .
\end{equation*}
 
\item Given $p,q\in(M,g)\,$ and $L\in\hbox{End}((T_pM,g_p),(T_qM,g_q))$, write $\Vert L\Vert_{p,q}$ for
the operator norm of $L\,$ (sup-norm on $T_qM \otimes T^\ast_pM $) and $\Vert L\Vert_{\rm Eucl}$ for the natural Euclidean norm, one has
$$\Vert L\Vert_{p,q}\leq \Vert L\Vert_{\rm Eucl}\ .
$$

\item Denoting by $\Vert \cdot\Vert_{q,q,q}$ the sup-norm and by $\Vert\cdot\Vert_{\rm Eucl}$ the natural Euclidean norm on $T_qM \otimes T^\ast_qM\otimes T^\ast_qM \,$, one also has
$$\forall {\mathfrak S}\in T_qM\otimes T^\ast_qM \otimes T^\ast_qM\ \ \ \ \ 
\Vert {\mathfrak S}\Vert_{q,q,q}\leq\Vert{\mathfrak S}\Vert_{\rm Eucl}\  .
$$
\end{enumerate}
\end{defn}

The section is devoted to introduce the magnifying glass, then to prove the two 
theorems 

\begin{thm}\label{theorem 1st} Assume $k\geq7$. 
There exists a 
constant ${\mathcal M}_1$, uniform on $W\subset(M,g)$, such that, 
given a positive $\rho\leq{\mathcal R}_W\,$, 
for any $x$ in 
$W$ and any $p,q$ in $B(x,\rho)$, one has
\begin{equation*}
\Vert\,{}^\sharp\hbox{\rm hess}_q(\frac{d^2(p,\cdot)}{2})-
\hbox{\rm Id}_q\,\Vert_{q,q}
\leq {\mathcal M}_1\,\max (d^2(p,x),d^2(q,x))\leq {\mathcal M}_1\,\rho^2\  .
\end{equation*}
\end{thm}

\begin{thm}\label{theorem 2nd} Assume $k\geq 8\,$. 
There exists a constant ${\mathcal M}_2\,$, 
such that, given 
a positive $\rho\leq{\mathcal R}_W\,$,  
for any 
$x$ in $W$ and any $p,q$ in $B(x,\rho)$, one has
\begin{equation*}
\Vert\,(\nabla {}^\sharp\hbox{\rm hess}(\frac{d^2(p,\cdot)}{2}))_q\,\Vert_{q,q,q}
\leq {\mathcal M}_2\,\max (d(p,x),d(q,x))\leq {\mathcal M}_2\,\rho\  .
\end{equation*}
\end{thm}

\subsection{The magnifying glass and its metric properties.}\label{gloups33}

\begin{defn} \label{D.family}
$ $
\begin{enumerate}
 \item \label{mgsc1} 
Define the $C^{k-1}$-family ($\exp$ is $C^{k-1}$) of {\it scalings} on $B=B(x,\rho)$
\begin{equation*}
\psi_x: (t,q)\in[-1,1]\times B\longmapsto\exp_x(t\,\exp_x^{-1}q)\in B \  .
\end{equation*}
Set $\psi_t (q)$ short for $ \psi_{t,x}(q)=\psi_x (t,q)$.
Given $t\in[-1,1]\setminus\{0\}\,$, the maps $\psi_t : B(x,\rho) \longrightarrow B(x,t\rho)$ and 
$\psi_{\frac{1}{t}}: B(x,t\rho)\longrightarrow  B(x,\rho)$ are diffeomorphisms verifying
$(\psi_t)^{-1} = \psi_{\frac{1}{t}}\,$.

\item \label{mgmet1} For any $t\in [-1,1]\setminus \{0 \}\,$, define also the family of {\it rescaled} Riemannian metrics 
on $B$
\begin{equation*}
g_{t,x}=\frac{1}{t^2}\  \psi_{t,x}^\ast g \ ,
\end{equation*}
and the corresponding {\it rescaled} squared distances 
$d^2_{t,x}$ on $B\,$.

\end{enumerate}
\end{defn}
\begin{rem} We often write $g_t$  short for $g_{t,x}$ and  $d_t$ 
instead of $d_{t,x}\,$.
\end{rem}

\begin{lem}\label{mgsc2} 
Let $k$ be $\geq 2\,$. The family of mappings $\frac{1}{t}\,d\psi_t$ extends -  {\rm define} it to 
be $d\exp_x^{-1}$ at $t=0$ -  to a 
$C^{k-2}$-family on the whole $[-1,1]\times B\,$. Moreover, the family of inverses 
$t\,d\psi_{\frac{1}{t}}\circ \psi_t$  also extends 
to a 
$C^{k-2}$-family when defined to be $d\exp_x \circ \exp_x^{-1}$ at $t=0\,$.
\end{lem}

\begin{proof} 
If $t\not=0\,$, one has
\begin{equation} \label{mgsc3}
\frac{1}{t}\ d\psi_t(\cdot)=d\exp_x(t\exp_x^{-1}(\cdot))\circ 
d\exp_x^{-1}(\cdot)\,.
\end{equation}
The direct claim follows. If $t\not=0\,$, the claim on the inverses results from  
$(\frac{1}{t}\ d\psi_t)^{-1}=t\,d\psi_{\frac{1}{t}} \circ \psi_t\,$, then the 
implicit 
function theorem gives the $C^{k-2}$ regularity including the value $t=0\,$.

Observe that $t\,d\psi_{\frac{1}{t}}\circ \psi_t$  
tends to $d\exp_x \circ \exp_x^{-1}$ 
when $t$ tends to $0\,$.
\end{proof}

The main step towards the magnifying glass, its metric properties:
\begin{prop}\label{mg4} 
Let $k$ be $\geq3\,$. Given $t\in[-1,1]\setminus\{0\}\, , \ p,q\in B\ ,\ v,w\in T_pM\,$, the families
\begin{equation*} 
(t,x,v,w)\longmapsto g_{t,x}(v,w)\hskip4mm\hbox{and}\hskip4mm
(t,x,p,q)\longmapsto d^2_{t,x}(p,q)
\end{equation*}
extend at $t\!=\!0$ in $C^{k\!-\!2}$ and $C^{k\!-\!3}$ families respectively defined at $t\!=\!0$ by
$g_{0,x}(v,w)\!:=\!(\exp_x^{-1})^\ast g_x(v,w)$ and by $d^2_{0,x}(p,q)\!=\!\Vert\exp_x^{-1}p-\exp_x^{-1}q\Vert^2_x\,$.
\par {\rm Moreover,} if $k\geq4\,$,
one has, for any $p\in B$ 
\begin{equation} \label{mghess}
\frac{\partial }{\partial t}(d_{t,x}^2(p,\cdot)_{\mid t=0}\equiv 0\,.
\end{equation}
\end{prop}
The variation by geodesics $\gamma$ below plays a central role.
\begin{defn}\label{vargamma}
Given the geodesics $p(t)\!=\!\exp_x(t\exp_x^{-\!1}p)$ and $q(t)\!=\!\exp_x(t\exp_x^{-\!1}q)$, define, for each fixed $t$ in $[0,1]\,$, the 
unique geodesic $\gamma(s,t)$ joining the point 
$\gamma(0,t)=p(t)$ to the point $\gamma(1,t)=q(t)$. 
Setting $v=\exp_x^{-1} p$ and  $w=\exp_x^{-1} q\,$, the variation $\gamma$ 
reads
\begin{equation}\label{variation}
\begin{split}
\gamma: 
(s,t)\in[0,1]\times[0,1]\longmapsto\gamma(s,t)=\exp_{p(t)}(s\exp_{p(t)}^{-1}
q(t))\,,\\
 \gamma(0,t) =p(t)=\exp_x(tv)=\psi_t(p)\ ,   \\ 
 \gamma(1,t)=q(t)= 
\exp_x(tw)=\psi_t(q) \,  .
\end{split}
\end{equation}
\end{defn}

\begin{proof}
If $t\!\not=\!0\,$, taking care of (\ref{mgsc3}),
define the metric $\check g_{t,x}\! =\!\exp_x^\ast g_{t,x}$ on $B(0_x,\rho)$. 
The metric 
$\check g_{t,x}=g(d\exp_x(t\exp_x^{-1}(\cdot)),d\exp_x(t\exp_x^{-1}(\cdot)))$ 
tends to $g_x$ as 
$t$ 
goes to $0\,$. With $d\exp\,$, the family $\check g_{t,x}$ has $C^{k-2}$ 
regularity 
on $[-1,1]\times B(0_x,\rho)$, this gives through $d\exp_x$ the conclusion for 
$g_{t,x}\,$.
\par
Near the diagonal in $M\times M\,$, the $C^{k-1}$ regularity of  $(p,q)\mapsto d_g^2(p,q)$ follows from  
$d^2(p,q)=\Vert\exp_p^{-1}q\Vert_p^2=\Vert\exp_q^{-1}p\Vert_q^2\,$. This formula 
reads for the squared $g_{t,x}$-distance
\begin{equation} \label{mg-1}
d_{t,x}^2(p,q)=\frac{1}{t^2}\  
\Vert\exp_{\psi_t(p)}^{-1}{\psi_t(q)}\Vert_{\psi_t(p)}^2\,.
\end{equation}
We know (see above) that for any $t\not=0$ the family 
$d_{t,x}^2(p,q)$ is $C^{k-1}$ because $g$ is $C^k\,$ and {\it check the regularity 
including the value} $t=0\,$.
Write 
$\delta(t,x,p,q)=t^{2}\,d_{t,x}^{2}(p,q)$. It is well defined and $C^{k-1}$ in 
$(t,x,p,q)\,$, 
including the value $t=0\,$. One has 
$\delta(0,x,p,q)=\Vert\exp_x^{-1}x\Vert_x^2=0\,$. Writing $Y=\frac{\partial \gamma}{\partial t}$ the Jacobi 
field associated to $\gamma$ (definition \ref{vargamma}), one 
computes from (\ref{mg-1})
\begin{multline*}
\frac{\partial \delta}{\partial t}(0,x,p,q)=
2\langle\frac{\partial }{\partial 
t}(\exp_{\psi_t(p)}^{-1}\psi_t(q))_{\vert t=0}, 
\exp_x^{-1}x\rangle=\\=2\langle\nabla_{\frac{\partial \gamma}{\partial 
s}}Y(0,0), 
\frac{\partial \gamma}{\partial s}(0,0)\rangle=0\ .
\end{multline*}
Consequently, one gets  
$\,
\delta(t,x,p,q)=t^{2}\,\int_0^1(1-\tau)\ 
{\partial^2 \delta}/{\partial t^2}(\tau t,x,p,q)\,d\tau\,$,
which also reads (if $t\not=0$)
$\ 
d_{t,x}^{2}(p,q)=\int_0^1(1-\tau)\,{\partial^2 \delta}/{\partial t^2} 
(\tau t,x,p,q)\,d\tau\,$.
The function defined on the right hand-side is $C^{k-3}$ in $t,x,p,q$, 
{\it including} 
$t=0\,$, so $d_{t,x}^{2}(p,q)$ extends with the same regularity at $t=0\,$.

The proof of the remaining claims is given in appendix \ref{app2}. 
\end{proof}

\noindent{\bf A precise description of the {\it magnifying glass} at $x\,$.}\label{magnglass}
\begin{lem}\label{mg5} For any $t\in [-1,1]\,$, the two inverse $C^{k-2}$ 
families 
of scalings $\frac{1}{t}\,d\psi_t$ and $t\,d\psi_\frac{1}{t}\, $
act isometrically on distinguished slices of $[-1,1]\times TB$ through
\begin{gather*}
\frac{1}{t}\,d\psi_t : (1,q,u)\in (\{1\}\times TB,g_t) \longmapsto (t,q,u)\in 
(\{t\}\times TB,g) \,,\\
t\,d\psi_{\frac{1}{t}} :  (t,q,u)\in (\{t\}\times TB,g) \longmapsto (1,q,u)\in 
(\{1\}\times TB,g_t)\, .
\end{gather*}
The integral curves of $\frac{\partial }{\partial t}$ are geodesics 
in
$(\tilde B, \tilde g)=([-1,1]\times B,dt^2+g_t)$.
\end{lem}

\begin{proof} 
The {\it first} statement is putting together definition \ref{mgsc1} of 
the scalings, 
lemma \ref{mgsc2}, definition \ref{mgmet1} and 
proposition \ref{mg4}, while the {\it second}
comes with the definition of the metric 
$\tilde g$ on $\tilde B\,$.
\end{proof}

\begin{rem}
On the product $\tilde B=[-1,1]\times B\,$, consider the two metrics 
$\tilde g=dt^2+ g_t$ and $\bar g =
dt^2+ g\,$. The action of the ``magnifying glass''  can be caught in one picture. 
Define the map 
\begin{equation*} 
{\bf \Psi}:(t,q)\in (\tilde B, \tilde g)\longmapsto (t,\psi_t(q))\in 
(\tilde B, \bar g)\, .
\end{equation*}
For each $t\!\in\![-\!1,1]\!\setminus\!\{0\}$, set $B_t\!=\!\{t\} \!\times \!B$. The map 
${\bf \Psi}_t\!:\!(t,q)\!\in\! B_t\longmapsto (t,\psi_t(q))\!\in\! B_t$ is an isometry 
composed with a dilation (magnification) of factor $\frac{1}{t}\,$.
\end{rem}

Not so surprising is the following result
\begin{lem}\label{mg7} In $(\tilde B,\tilde g)$ the Euclidean slice $B_0=\{0\}\times B$ is 
totally 
geodesic: if 
$S_t(u)=\tilde\nabla_u\frac{\partial }{\partial t}$ is
the shape operator of the slice $B_t=\{t\}\times B$ in direction $u\in TB_t\,$, one has $S_0\equiv 0$ (as well as $\frac{\partial g}{\partial t}(0)\equiv 0\,$).
\end{lem}
The proof is given in appendix \ref{app20}

\subsection{The Hessian and its derivatives in the magnifying glass.}
\begin{defn}
Let $\nabla^t$ be the Levi-Civita connection of $g_t=g_{t,x}\,$. 
Denote by $\hbox{\rm grad}_t\,f$ and by $\hbox{\rm hess}^t\,f$ the $g_t$-gradient and the $g_t$-Hessian (i. e. $\hbox{\rm hess}^t f=\nabla^t\,df $)
of a function $f$ on $B\,$.
\end{defn}
\begin{rem} In the lines to come, for fixed $x,t\,$, we denote by $\hbox{\rm hess}^t\,d_{t,x}^2$ the $g_t$-Hessian of the function
$q\in B\mapsto d_{t,x}^2(p,q)\in\mathbb{R}\,$, also depending on a given $p\in B\,$. So, $\hbox{\rm hess}^t\,d_{t,x}^2$ depends on $(p,q)\in B\times B$ through
$$(\hbox{\rm hess}^t\,d_{t,x}^2)(p,q):=\hbox{\rm hess}_q^t\,d_{t,x}^2(p,\cdot_q)\  .
$$
\end{rem}

\begin{cor} 
If $k\geq 6\,$, at $t=0$ the family $\hbox{\rm hess}^t\,d_{t,x}^2$ 
$($defined for any $t\in\,]0,1])$ extends $C^{k-5}$ in a family defined for 
$t\in[-1,1]\,$, while $(\nabla^t\,\hbox{\rm hess}^t)\,d_{t,x}^2$ extends 
$C^{k-6}$ in the corresponding sense.
\end{cor}

\begin{proof} 
This is a direct consequence of previous proposition \ref{mg4}.
\end{proof}

\begin{lem} Given $t$ in $[-1,1]\,$, write $q'$ and $p'$ for 
$\psi_t(q)$ 
and $\psi_t(p)$. Define $d_p$ and $d_{p'}$ to be the functions
$d_p:q\mapsto d(p,q)$ and $d_{p'}:q'\mapsto d(p',q')$ (thus $d_{p'}(\psi_t(\cdot))=d_{\psi_t(p)}\circ \psi_t(\cdot)$).
\par
One has ($\sharp_t$ is paired with $g_t$ as $\sharp$ with $g$)
\begin{equation} \label{mgA} 
{}^{\sharp_t}\hbox{\rm hess}_q^t(d_{t,x}^2(p,\cdot))=d\psi_t^{-1}(q')\circ
{}^\sharp\hbox{\rm hess}_{q'}\,d_{p'}^2\circ d\psi_t(q)\,;
\end{equation}
\begin{equation} \label{mgB} 
(\nabla^t\,{}^{\sharp_t}\hbox{\rm hess}^t(d_{t,x}^2(p,\cdot)))_q=
d\psi_t^{-1}(q')\circ\nabla\,{}^\sharp
\hbox{\rm hess}_{q'}\,d_{p'}^2\circ d\psi_t(q)\, .
\end{equation}
\end{lem}

\begin{proof} 
If $\psi:(\breve M,\breve g)\rightarrow(M,g)$ is an isometry and $f$ a function on $(M,g)$, the 
$\breve g\,$-Hessian 
of $f\circ\psi$ verifies
$\   
\breve{\hbox{\rm hess}}_{\cdot}\, f\circ \psi=
\hbox{\rm hess}_{\psi(\cdot)}\, f ( d\psi (\cdot) ,d\psi (\cdot))\,.
$
One derives the claimed formulas, using  
$d_{t,x}^2(p,\cdot)=\frac{d^2(p',\psi_t(\cdot))}{t^2}\,$.
\end{proof}

\par

\begin{defn} \label{mgtan} 
A tensor $T$ of type $(m,l)$ on $\tilde B\,$ is  
an element of
$(\otimes^mT\tilde B)\otimes(\otimes^lT^\ast \tilde B)\,$.
It is tangential if, for any family of 
$1$-forms
$\varphi_1,\dots,\varphi_m$ 
among which one is colinear to 
$\langle\frac{\partial }{\partial t},\cdot\rangle\,$, respectively 
for any family of 
vectors $v_1,\dots,v_l$ among which one is colinear to $\frac{\partial 
}{\partial t}\,$, one has
\begin{equation*}
\varphi_1\otimes\dots\otimes\varphi_m(T)=0\  , 
\hskip3mm \hbox{respectively} 
\hskip3mm T(v_1\otimes\dots\otimes v_l)=0\  .
\end{equation*}
\end{defn}
\begin{rem} \label{tangt1} Any family of tensors $T_t\in\bigotimes B_t$ 
(differentiable 
in the parameter $t$) {\it defines} a tangential tensor $T$ in $\bigotimes \tilde B$ 
in a canonical way by {\it asking 
that} $T$ contracted with the field $\frac{\partial }{\partial t}$ (respectively 
with the 
$1$-form 
$\langle\frac{\partial }{\partial t},\cdot\rangle$) is $0$ {\it and also that} $T$ 
and $T_t$ have the same components  
along $B_t\,$: so, the {\it defined} $T$ is a tensor having non-zero components 
only along 
$B_t\,$. 
\end{rem}
\begin{lem}\label{tangt2} If $T$ is tangential, then $\nabla^t T$ and   
$\ \tilde{\nabla}_{\frac{\partial}{\partial t}} T$ are tangential.
\end{lem}
\begin{proof} The first claim is clear from the definitions. For a tangential field $X$ or $1$-form $\omega\,$, one has 
\begin{gather*}\langle\tilde{\nabla}_{\frac{\partial}{\partial t}}X,\frac{\partial}{\partial t}\rangle=
\frac{\partial}{\partial t}\langle X,\frac{\partial}{\partial t}\rangle-\langle X,\tilde{\nabla}_{\frac{\partial}{\partial t}}\frac{\partial}{\partial t}\rangle=0\  \ \hbox{and}\ \ \\
(\tilde{\nabla}_{\frac{\partial}{\partial t}}\omega)(\frac{\partial}{\partial t})=
\frac{\partial}{\partial t}(\omega(\frac{\partial}{\partial t}))-\omega(\tilde{\nabla}_{\frac{\partial}{\partial t}}\frac{\partial}{\partial t})=0\  .
\end{gather*}
As the covariant derivative acts as a {\it derivation} on a tensor $T$ of type $(m,l)$, it follows that $\tilde{\nabla}_{\frac{\partial}{\partial t}}T$ is also tangential in the general case.
\end{proof}
\begin{lem}\label{mgcom}\label{mgcom1} 
Let $T$ be a tangential tensor on $\tilde B\,$. 
If $S_{t_0}\equiv0$ on $B_{t_0}\,$, one has
\begin{equation*} \forall (t_0,q_0)\in B_{t_0}\ \ \ \ \ \ \ \ 
(\tilde{\nabla}_{\frac{\partial}{\partial t}}\nabla^t T)_{(t_0,q_0)}=
(\nabla^t\tilde{\nabla}_{\frac{\partial}{\partial t}} T)_{(t_0,q_0)}\  .
\end{equation*}
\end{lem}
The proof is given in appendix \ref{app22}.

Our next goal is to compute the $\tilde\nabla$ derivative of $
{}^\sharp\hbox{\rm hess}_q^t(d_{t,x}^2(p,\cdot))$ inserted in the 
``magnifying glass'' $(\tilde B,\tilde g)$ as a $(1,1)$-tensor (remark \ref{tangt1}). This allows to 
make Taylor expansions around $t=0$ and, 
with (\ref{mgA}) and (\ref{mgB}), will produce 
the 
estimates claimed in the theorems.

\begin{defn}
Let $f :\tilde B\longrightarrow \R$ be a $C^2$ function. 
Given any $t\in [-1,1]$ and $q\in B\,$, set $f_t (q):=f(t,q)$. 
The tensor  $^{\sharp_t} {\rm hess}^t f$ on $\tilde B$ is defined (see remark \ref{tangt1}) to be the 
unique tangential $(1,1)$-tensor satisfying, for any $(t,q)\in \tilde B$ and $u\in T_{(t,q))} B$
\begin{equation*}
 ^{\sharp_t} {\rm hess}_{(t,q)}^t f (u) = 
 \, ^{\sharp_t} {\rm hess}_q^t f_t (u) 
 =\nabla_u^t \,\grad_t f_t \ \ \ {\rm and}\ \ \ 
  ^{\sharp_t} {\rm hess}_{(t,q)}^t f (\frac{\partial}{\partial t}) = 0\,.
\end{equation*}
\end{defn}

\begin{lem} \label{mg8} 
Let $f:\tilde B\rightarrow \R $ be a differentiable 
function of $(t,q)$. 
One has on $B_0$
\begin{equation*}
\tilde{\nabla}_{\frac{\partial }{\partial t}}{}^{\sharp_t}\hbox{\rm hess}^t f=
\nabla^t\,\tilde{\nabla}_{\frac{\partial }{\partial t}}\,\grad_t f_t\  .
\end{equation*}
One also has
$\ 
\langle \tilde{\nabla}_{\frac{\partial }{\partial t}}\grad_t f_t, 
u\rangle=u(\frac{\partial f}{\partial t})\,$ and
$\ 
\langle \tilde{\nabla}_{\frac{\partial }{\partial t}}\grad_t f_t, 
\frac{\partial }{\partial t}\rangle=0\,$.
\end{lem}

\begin{proof} 
Taking $X=\grad_t f_t$ in lemma \ref{mgcom} ($\grad_t f_t$ is tangential, see remark \ref{tangt1}), one 
gets 
in $B_0$ (since $S_0=0$ by lemma  \ref{mg7})
\begin{equation*}
\tilde{\nabla}_{\frac{\partial }{\partial t}}{}^{\sharp_t}\hbox{\rm hess}^t f=
\tilde{\nabla}_{\frac{\partial }{\partial t}}\,\nabla^t\,\grad_t f_t
=\tilde{\nabla}_{\frac{\partial }{\partial t}}\,\nabla^t\,X 
=\nabla^t\,\tilde{\nabla}_{\frac{\partial }{\partial t}}\,X
=\nabla^t\,\tilde{\nabla}_{\frac{\partial }{\partial t}}\,\grad_t f_t\  .
\end{equation*}

Let $U$ be
a  $\tilde g$-parallel field along the integral lines of 
$\frac{\partial }{\partial t}\,$, everywhere tangent to the slices $B_t$
and such as $U_{\pi(u)}=u\,$.
As 
$[\frac{\partial }{\partial t},U]=
\tilde\nabla_{\frac{\partial }{\partial t}}U-
\tilde\nabla_U\frac{\partial }{\partial t}=0$ 
(use $\tilde\nabla_{\frac{\partial}{\partial t}} U=0$ and $S_0=0$ by lemma  \ref{mg7})
\begin{equation*}
\langle \tilde{\nabla}_{\frac{\partial }{\partial t}}\grad_t f_t, 
u\rangle_{\tilde g}=
\frac{\partial }{\partial t}\langle\grad_t f_t, U\rangle
=\frac{\partial }{\partial t}(d f(\cdot) (U))
=\frac{\partial }{\partial t}(U(f))=
u(\frac{\partial f}{\partial t})\, .
\end{equation*}
As $\grad_t f_t$ is tangential, lemma \ref{tangt2} implies
$\langle \tilde{\nabla}_{\frac{\partial }{\partial t}}\grad_t f_t, 
\frac{\partial }{\partial t}\rangle=0\,$.
\end{proof}

\begin{lem} \label{mghess1}
For any $p\in B\,$, one has at $t=0$
\begin{gather*}
\tilde\nabla_{\frac{\!\partial}{\partial t}_{\mid 
t=0}}\, \grad_t d_{t,x}^2(p,\cdot)=0 
 \ ,\ \ 
\tilde\nabla_{\!\frac{\partial }{\partial t}_{\mid t=0}}{}^{\sharp_t}\hbox{\rm hess}_\cdot^t 
(d_{t,x}^2(p,\cdot))= 0\ ,\\ \ \tilde\nabla_{\!\frac{\partial }{\partial t}_{t=0}}\nabla^t{}^{\sharp_t}
\hbox{\rm hess}^t_{\cdot}(\frac{d_{t,x}^2(p,\cdot)}{2})
= 0 \ \hbox{and} \ \ \tilde\nabla_{\!\frac{\partial }{\partial t}_{t=0}}\nabla^t\nabla^t\cdots\nabla^t(\frac{d_{t,x}^2(p,\cdot)}{2})
= 0\  .
\end{gather*}
\end{lem}
\begin{proof} Of course, one applies repeatedly remark \ref{tangt1} and lemma \ref{tangt2}, even if not mentioning this
fact explicitly.
Apply lemma \ref{mg8} to $f=d_{t,x}^2(p,\cdot)$, the first equality to prove 
is a consequence of (\ref{mghess}) (proposition \ref{mg4}). The second is a consequence of the 
first one and lemma \ref{mg8}.
\par
Then, $T={}^{\sharp_t}\hbox{\rm hess}_q^t (d_{t,x}^2(p,\cdot))/2$ 
defines a 
$(1,1)$-tangential 
tensor (remark \ref{tangt1}) on 
$\tilde B\,$. Since $S_0=0\,$ (lemma \ref{mg7}), one 
applies lemma \ref{mgcom1} 
to get $\ (\tilde{\nabla}_{\frac{\partial }{\partial t}}\nabla^t T)_{(0,\cdot)}=
(\nabla^t\tilde{\nabla}_{\frac{\partial }{\partial t}} T)_{(0,\cdot)}\,$ on $B_0\,$.
The second equality of lemma \ref{mghess1} tells 
${\tilde\nabla}_{\frac{\partial }{\partial t}} T\equiv0$ on 
$B_0\,$, which completes the proof.
\end{proof}

\subsection{Proving the theorems}
\noindent {\bf Proof of theorem \ref{theorem 1st}.}
\begin{proof} 
Set $w\!=\!\exp_x^{\!-1}q$ and, along the geodesic $q(t)\!=\! \exp_xtw$ in $(B,g)$, 
choose a $g$-parallel orthonormal frame 
${\mathcal V}_i$ and denote by 
${\mathcal U}_i^t$ the family of $\tilde g$-parallel orthonormal frames in 
$\tilde B$ defined by
the collection
$\  {\mathcal U}_i^t=t\,d\psi_t^{-1}(\psi_t(q))
({\mathcal V}_i)\,$  together with
$\ \frac{\partial }{\partial t}\,$. So, one has 
$g_t({\mathcal U}_i^t,{\mathcal U}_j^t)=g({\mathcal V}_i,
{\mathcal V}_j)=\delta_i^j$ 
and the 
${\mathcal U}_i^t$ are all tangent to $B_t\,$.
From work done before, using lemma \ref{mghess1} (based on (\ref{mghess})),
we can make the
following Taylor
expansion in $t$ around $t=0$ of
${}^{\sharp_t}\hbox{\rm hess}_\cdot^t(\frac{d_{t,x}^2(p,\cdot)}{2})$
\begin{multline} \label{mg11}
\langle{}^{\sharp_t}\hbox{\rm hess}_\cdot^t(\frac{d_{t,x}^2(p,\cdot)}{2}) 
({\mathcal U}_i^t),{\mathcal U}_j^t\rangle_{\tilde g}
=\langle{\mathcal U}_i^0,{\mathcal U}_j^0\rangle_{\tilde g} + \\ 
+t^2\int_0^1(1-\tau)\,\frac{d^2}{d t^2}
(\langle{}^{\sharp_t}\hbox{\rm hess}_\cdot^t(\frac{d_{t,x}^2(p,\cdot)}{2}) 
({\mathcal U}_i^t),
{\mathcal U}_j^t\rangle_{\tilde g})_{t\tau}\,d\tau\  .
\end{multline}

Setting $p'=p(t)= \psi_t (p)$ and $q'=q(t)=\psi_t (q)$ as before, and 
$H={}^\sharp\hbox{\rm hess}_{q'}(\frac{d_{p'}^2}{2})$ (recall $d_{p'}$ is the function $q'\mapsto d(p',q')$),
we get 
\begin{align*}
\langle{}^\sharp\hbox{\rm hess}_{q'}(\frac{d_{p'}^2}{2})& ({\mathcal V}_i),
{\mathcal V}_j\rangle_g = 
 \langle H (\frac{1}{t}\,d \psi_t(q) ({\mathcal U}_i^t)),
(\frac{1}{t}\,d \psi_t(q) ({\mathcal U}_j^t))\rangle_{g} \\
=& \langle t \,d\psi_t^{-1}(q')\circ H \circ 
\frac{1}{t} d \psi_t(q) ({\mathcal U}_i),t \,d\psi_t^{-1}(q')\circ 
\frac{1}{t} d \psi_t(q) ({\mathcal U}_j)\rangle_{g_t} \\
=&\langle t \,d\psi_t^{-1}(q')\circ H \circ 
\frac{1}{t} d \psi_t(q)({\mathcal U}_i)
,{\mathcal U}_j\rangle_{g_t} 
=\\=&\langle d\psi_t^{-1} (q')\circ H \circ 
d \psi_t(q) ({\mathcal U}_i),{\mathcal U}_j\rangle_{g_t}\ .
\end{align*}
Using equation (\ref{mgA}), we have 
\begin{equation} \label{mg12}
 \langle{}^\sharp\hbox{\rm hess}_{q'}(\frac{d_{p'}^2}{2}) ({\mathcal V}_i),
{\mathcal V}_j\rangle_g = \langle{}^{\sharp_t}\hbox{\rm 
hess}_q^t(\frac{d_{t,x}^2(p,\cdot)}{2}) 
({\mathcal U}_i^t), {\mathcal U}_j^t\rangle_{\tilde g}\ .
\end{equation}

\noindent Thus (\ref{mg11}) and (\ref{mg12}) produce at $q'$ in $B(x,t\rho)$ an estimate 
of 
$d_{p'}^2\,$: thanks to (\ref{mg12}), bounding 
\begin{equation*}
\sup_{t:\vert t\vert\leq 1} \big\vert \frac{d^2}{d t^2}
(\langle{}^{\sharp_t}\hbox{\rm hess}_q^t(\frac{d_{t,x}^2(p,\cdot)}{2}) 
({\mathcal U}_i^t),{\mathcal U}_j^t\rangle_{\tilde g})_t \big\vert
\end{equation*}
is equivalent to bound, along the geodesic $p(t)\times q(t)$ in $(B\times B,g\times g)$, the equal term
\begin{equation*}
{\mathfrak M}_{1,i,j}(p,q)=
\sup_{t:\vert t\vert\leq 1} \big\vert \frac{d^2}{d t^2}
(\langle{}^{\sharp}\hbox{\rm hess}_{\psi_t(q)}
(\frac{d_{\psi_t(p)}^2}{2}) 
({\mathcal V}_i),{\mathcal V}_j\rangle_g)_t \big\vert\   .
\end{equation*}

Denote by ${\mathcal V}_j^\ast$ the dual basis of 
${\mathcal V}_j\,$ and by ${\mathfrak T}$ a tensor verifying 
${\mathfrak T}_{p,q}\in \bigotimes^{1,1} (\{0_p\}\times T_qM)$ (with $(p,q)\in B\times B$). 
Set  ${\mathfrak T}_{ij} = 
\langle {\mathfrak T}({\mathcal V}_i), {\mathcal V}_j\rangle\,$. Thus 
\begin{equation}\label{mg14}
{\mathfrak T}=\sum_{i,j}{\mathfrak T}_{ij}\,{\mathcal V}_j\otimes 
{\mathcal V}_i^\ast\ .
\end{equation}
Here we have 
$\ 
{\mathfrak T}={}^\sharp
\hbox{\rm hess}_{\cdot_q}(\frac{d_{\cdot_p}^2}{2})\,$ and 
$\ {\mathfrak T}_{ij}=\langle{}^\sharp\hbox{\rm hess}_{\cdot_q}
(\frac{d_{\cdot_p}^2}{2}) ({\mathcal V}_i),{\mathcal V}_j\rangle_g\,$.

One deduces from  (\ref{mg14}) 
\begin{equation} \label{mg15}
\vert\,{\mathfrak T}_{ij}\,\vert\leq \Vert\,{\mathfrak T}\,\Vert_{\rm Eucl} \hskip2mm 
\hbox{in coherence with}
\hskip2mm \Vert\,{\mathfrak T}\,\Vert_{\rm Eucl}^2=\sum_{i,j}{\mathfrak T}^2_{ij} \ .
\end{equation}

\noindent {\bf Claim} {\em 
Consider a $(1,1)$-tensor ${\mathfrak T}$ on $(M\times M,g\times g)$ (as above) and its restriction 
along a product geodesic $t\mapsto (p(t),q(t))$, i. e. consider $t\mapsto {\mathfrak T}\circ (\psi_t,\psi_t)(p,q)={\mathfrak T}_{(p(t),q(t))}\,$. One has}
\begin{equation*}
\Vert\,\nabla_{\frac{\partial }{\partial t}}^{g\times g}{\mathfrak T}(p(\cdot),q(\cdot))\,\Vert_{\rm Eucl}\leq 
\sqrt{\Vert\,\frac{d p}{d t}\,\Vert^2+\Vert\,\frac{d q}{d t}\,\Vert^2}\  \ 
\Vert\,\nabla^{g\times g} {\mathfrak T}\,\Vert_{\rm Eucl} \  .
\end{equation*}

\begin{proof}
Indeed, one has $\nabla_{\frac{\partial }{\partial t}}^{g\times g}{\mathfrak T}
=\sum_{i,j}\nabla_{\frac{\partial }{\partial t}}^{g\times g}{\mathfrak T}_{ij}
\,{\mathcal V}_i\otimes {\mathcal V}_j^\ast$ and 
\begin{equation*}
\big{(}\frac{d}{dt}{\mathfrak T}_{ij}\circ(\psi_t,\psi_t)\big{)}_t=
\langle(\nabla^{g\times g} {\mathfrak T}_{ij})_{(p(t),q(t))},(\frac{d p(t)}{d t},\frac{d q(t)}{d t})\rangle_{g\times g}\,. \qedhere
\end{equation*}
\end{proof}

Writing $(\nabla^{g\times g})^2$ the second order covariant derivative in the metric $g\times g$ and 
interpreting the field $\frac{\partial }{\partial t}$ as a field along the geodesic $(p(t),q(t))$ and the $d/dt$-derivative
as $\nabla^{g\times g}_{\!\!(\frac{d p}{d t},\frac{d q}{d t})}\,$, we have, for any $p$ and $q$ in $B(x,t\rho)$
\begin{equation}\label{mg13}
{\mathfrak M}_{1;i,j}(p,q)\leq 2\ \mu(p,q)\!\!\!\sup_{p,q\in B(x,t\rho)} \Vert\,
\big((\nabla^{g\times g})^2({\mathfrak T}_{ij})\big)_{p,q}\,\Vert_{\rm Eucl}\  ,
\end{equation}
with $\mu(p,q)=\max(d^2(p,x),d^2(q,x))\,$.

In (\ref{mg13}), the right-hand side is bounded by compactness of $B$ and $C^4$ 
regularity of $d^2\,$. 
We need $C^5$ regularity of $g$ at this stage, but $C^7$ regularity of $g$ was 
necessary to state 
expansion (\ref{mg11}) :  there,
${}^\sharp\hbox{\rm hess}^t(d_{t,x}^2)$ (respectively $g_{t,x}$) are at 
least 
$C^4$ ($C^5$), so $g$ should at least be $C^7\,$, see proposition \ref{mg4}. Since 
$\rho\leq {\mathcal R}_W\,$, the proof is 
complete  if, in view of (\ref{mg11}),(\ref{mg12}), (\ref{mg15}), (\ref{mg13}), we set ${\mathcal M}_1=2\,{\mathcal M}$ with 
\begin{equation}\label{onlyyou}
{\mathcal M}:=\sup_{x\in W}\, \sup_{p,q\in B(x,{\mathcal R}_W)}\Vert\,
\big((\nabla^{g\times g})^2({}^\sharp\hbox{\rm hess}_{\cdot_q}
(\frac{d^2_{\cdot_p}}{2}))\big)_{p,q}\,\Vert_{\rm Eucl}  \  .\qedhere
\end{equation}
\end{proof}

\begin{rem} \label{mg16} See \cite{B-K} section 6 for a far more
precise version of this theorem through comparison of second order differential 
equations.
\end{rem}

To prove the last theorem, we use the classical result stated below (for a proof, see the appendix \ref{app222}).
\begin{prop}\label{mgapp1} 
Let $f$ be a $C^k$ function of $t$ varying in $I=[-T,T]$ such 
that 
$f(0)=0\,$,
with $T>0$ and $k\geq1\,$. Define $g(t)=f(t)/t\,$. The function $g$ is $C^{k-1}$ and the 
derivatives of $g$ in 
$I$ are bounded in terms of those of $f\,$, more precisely, for any $l$ with 
$0\leq l\leq k-1\,$,
 one has
\begin{equation*} 
\sup_{t\in I} \,\big\vert \frac{d^l g}{d t^l} (t) \big\vert\leq 
\frac{1}{l+1}\,\sup_{t\in I}\, \big\vert \frac{d^{l+1} f}{d t^{l+1}} 
(t)\big\vert\ .
\end{equation*}
\end{prop}
\noindent {\bf Proof of theorem \ref{theorem 2nd}}
\begin{proof} 
With the same notations as in the above proof of theorem \ref{theorem 1st} 
and the help of (\ref{mgB}), we get, setting $p'=p(t)$ and $q'=q(t)$ as before
\begin{equation}\label{mg18} 
\langle \!(\nabla {}^\sharp\hbox{\rm 
hess}\frac{d_{p'}^2}{2})_{q'}({\mathcal V}_l,{\mathcal V}_i),
{\mathcal V}_j\rangle_g 
\!=\!\frac{1}{t}\langle \!(\nabla^t\,{}^{\sharp_t}
\!\hbox{\rm hess}^t\frac{d_{t,x}^2(p,\cdot)}{2} )_{q}
({\mathcal U}_l^t,{\mathcal U}_i^t),
{\mathcal U}_j^t\rangle_{\tilde g}.
\end{equation}
The right-hand side is $0$ at $t=0$ (by lemma \ref{mghess1}) and at least 
$C^1$, 
as $d_{t,x}^2(p,\cdot)$ is at least $C^5$ if $k\geq8$ (apply proposition 
\ref{mg4} and proposition \ref{mgapp1}). Thus, we may write (recall $\tilde g=dt^2+g_t$)
\begin{multline*}
\frac{1}{t}\langle(\nabla^t\,{}^{\sharp_t}
\hbox{\rm hess}^t(\frac{d_{t,x}^2(p,\cdot)}{2}))_q ({\mathcal U}_l^t,
{\mathcal U}_i^t), {\mathcal U}_j^t\rangle_{\tilde g} 
=\\=t\!\int_0^1\!\!\frac{\partial }{\partial t}(\frac{1}{t}\,
\langle(\nabla^t\,{}^{\sharp_t}\hbox{\rm hess}^t(\frac{d_{t,x}^2(p,\cdot)}{2}))_q 
({\mathcal U}_l^t,{\mathcal U}_i^t),
{\mathcal U}_j^t\rangle_{\tilde g})_{t\tau}\,d\tau\,.
\end{multline*}

Thinking to the next steps, introduce 
\begin{multline*}
{\mathfrak S}_{\mid(\cdot_{p'},\cdot_{q'})}=\!(\nabla {}^\sharp\hbox{\rm hess}
(\frac{d_{\cdot_{p'}}^2}{2}))_{\cdot_{q'}}\!=\!
\sum_{l,i,j}{\mathfrak S}_{l,i,j}(\cdot_{p'},\cdot_{q'})\  
{\mathcal V}_j\otimes{\mathcal V}_l^\ast\otimes{\mathcal V}_i^\ast \!=\\=\!
\sum_{l,i,j}
\langle(\nabla {}^\sharp\hbox{\rm hess}(\frac{d_{\cdot_{p'}}^2}{2}) )_{\cdot_{q'}}
({\mathcal V}_l,{\mathcal V}_i),
{\mathcal V}_j\rangle_g\  
{\mathcal V}_j\otimes{\mathcal V}_l^\ast\otimes{\mathcal V}_i^\ast\  .
\end{multline*}

Define ${\mathfrak M}_{2;lij}(p,q)$ to be
\begin{equation*} 
{\mathfrak M}_{2;lij}(p,q)=\sup_{t:\vert t\vert\leq 1}\left\vert 
\frac{\partial }{\partial t}(\frac{1}{t}\,\langle (\nabla^t\,{}^{\sharp_t}
\hbox{\rm hess}^t(\frac{d_{t,x}^2(p,\cdot)}{2}) )_q
({\mathcal U}_l^t,{\mathcal U}_i^t),
{\mathcal U}_j^t\rangle_{\tilde g})\right\vert\ , 
\end{equation*}
which may be written with the help of equation (\ref{mg18})
\begin{equation*} 
{\mathfrak M}_{2;lij}(p,q)=\sup_{t:\vert t\vert\leq 1}\left\vert 
\frac{\partial }{\partial t}\,\langle (\nabla {}^\sharp\hbox{\rm 
hess}(\frac{d_{p(t)}^2}{2}))_{q(t)}
({\mathcal V}_l,{\mathcal V}_i),
{\mathcal V}_j\rangle_g\right\vert\  \ ,
\end{equation*}
which is thus bounded. Again, interpret the
$\frac{d}{dt}$-derivative  
as $\nabla^{g\times g}_{\!\!(\frac{d p}{d t},\frac{d q}{d t})}\,$. 
For any $p$ and $q$ in 
$B(x,\rho)$ (with $x\in W$ and $\rho\leq{\mathcal R}_W$) this leads to
\begin{equation*}
{\mathfrak M}_{2;lij}(p,q)\leq 
 \sqrt{2\,\mu(p,q)}\,\sup_{p,q\in B(x,\rho)} \Vert\,
(\nabla^{g\times g}{\mathfrak S}_{l,i,j})_{p,q}\,\Vert\  ,
\end{equation*}
with $\mu(p,q)=\max(d^2(p,x),d^2(q,x))$ as in the proof of theorem \ref{theorem 1st}.
Defining the norm of a tensor in the way we did proving 
theorem \ref{theorem 1st} and choosing ${\mathcal M}_2=\sqrt{2}\,{\mathcal M}$ (see (\ref{onlyyou})) completes the proof of theorem \ref{theorem 2nd},
because of the inequality (due to the nature of a product metric)
\begin{equation*}
\sup_{x\in W}\,\sup_{p,q\in B(x,{\mathcal R}_W)} \Vert\,
\big( \nabla^{g\times g} ((\nabla \ {}^\sharp\hbox{\rm hess})_{\cdot_q}
(\frac{d^2_{\cdot_p}}{2}))\big)_{p,q}\,\Vert_{\rm Eucl}\leq {\mathcal M}\ .\qedhere  
\end{equation*}
\end{proof}

\section{Riemannian barycentric simplex close to Euclidean}
\label{S.close}
In this section ${\mathcal R}_0\leq{\mathcal R}_0(1/2) \leq {\mathcal R}_W $ refers to section \ref{EuclRiem}, where, in
proposition \ref{proposition A}, it guarantees that $\exp_x$ is for any 
$x \in W$ a 
$1/2$-quasi-isometry in $B(x,{\mathcal R}_0)$. 
Moreover, ${\bf p}=\{p_0,p_1,\dots,p_n\}$ is a spread set of 
points in a convex 
ball $B(x,{\mathcal R}_0)$, which, by definition \ref{D.spread.simplex}-\ref{spread12}, implies that, for any $q\in\hat\tau$, the linear 
simplex $\hat\tau_q$ with vertices $w_0(q)=\exp_q^{-1}p_0,\dots,w_n(q)=\exp_q^{-1}p_n$ has $g_q$-thickness $t_q(\hat\tau_q)>0\,$. 
Thus ${\bf p}$ 
builds a true Riemannian barycentric simplex $\hat\tau_{\bf p}=\hat\tau\,$, 
see theorem \ref{rbs0} and its corollary. 
 
\begin{defn} \label{genthick}
The (generalized) thickness of a Riemannian barycentric simplex is 
\begin{equation*}
t_g(\hat\tau)=\inf_{q\in\hat\tau}t_q(\hat\tau_q)\,,
\end{equation*}
where $\hat\tau_q\subset T_qM$ is the linear 
simplex built on the vertices $w_i(q)=\exp_
q^{-1}p_i\,$. 
The fields $w_i(q)$ are central in the computations to go.
\end{defn}

Given some $t_0>0\,$, we {\it focus} in this section on barycentric simplices 
$\hat\tau=\hat\tau_{\bf p}$ of thickness $t_g(\hat\tau)\geq t_0\,$, 
while we {\it a priori restrict the study} to those $\hat\tau=\hat\tau_{\bf p}$ whose vertices are at mutual Riemannian distances $d(p_i,p_j)$ and allow to exist an actual linear
$n$-dimensional Euclidean
simplex with vertices at mutual distances $d(p_i,p_j)$
(we assume this Euclidean
simplex to exist in section \ref{S.close}).

On such a $\hat\tau\,$, a flat metric $\hat g$ for which the vertices are at mutual distances $d(p_i,p_j)$ can be defined 
by means of the ``barycentric Riemannian coordinates'': we now describe a 
way to do this.

\begin{defn} \label{D.euclidianisable}
Equip $T_xM$ with the Euclidean structure $\bar g$ (in general $\not=g_x$) such 
that, writing $v_0=\exp_x^{-1}p_0, \dots,v_n=\exp_x^{-1}p_n$  
\begin{equation*}
\forall i,j=0,1,\dots,n\hskip1cm\Vert\,\,v_i-v_j\,\Vert_{\bar g}=d(p_i,p_j)\,.
\end{equation*}
Call $\phi$ the diffeomorphism that sends a point $v=\sum_{i=0}^n\lambda_i\,v_i$ 
varying 
in a neighborhood of $\hat\tau_x$ in $T_xM$ to the unique minimum of  
$\sum_{i=0}^n\lambda_i\,d^2(p_i,\cdot)$ in $\hat\tau\!=\!\phi(\hat\tau_x)\!\subset \!B(x,{\mathcal R}_0)$ (see section \ref{rbs}, ${\bf \lambda}=\{\lambda_0,\dots,\lambda_n\}$ varies in some open 
subset of $\R^{n+1}$ containing 
$\Delta=\{\lambda\mid\sum_{i=0}^n\lambda_i=1\ ,\  \hbox{all}
\ \lambda_i\geq0\}$).
\end{defn}

\begin{defn}\label{bareucl} $ $
\begin{enumerate}
 \item The metric $\hat g$ on $\hat\tau$ such that $\phi^\ast\hat g=\bar g$ defines a 
{\it Euclidean realization $(\hat\tau,\hat g)$}, so $(\hat\tau,g)$ {\it admits a Euclidean realization}. 
Call $h$ the tensor $h=\hat\nabla-\nabla$ difference of the Levi-Civita 
connections associated with $\hat g$ and $g\,$ (see section \ref{diffparll}).

\item Given  a $C^1$-curve $\gamma \!:\! [0, 1]\rightarrow \!\hat\tau$ with 
$\gamma(0)=q\,$, call ${\mathcal P}_\gamma^g$ and ${\mathcal P}_\gamma^{\hat g}$ the 
$g$ and $\hat g$-parallel 
translations along $\gamma\,$, from $q$ to $\gamma (1)$.
\end{enumerate}
\end{defn}

Recall that ${\mathfrak R}_0,{\mathcal M}_1,{\mathcal M}_2$ are positive reals introduced in definition \ref{Riembnd}, 
theorem \ref{theorem 1st} and theorem \ref{theorem 2nd}.
\begin{defn}\label{constanteC_7R_3} A real $t_0>0$ is given. As it plays everywhere the role of a lower bound on the thickness of Euclidean $n$-simplices as in definition \ref{genthick}, it also satisfies $t_0<1\,$, and even $t_0<1/\sqrt{n(n+1)}\,$).
\begin{enumerate}
\item Define the constant ${\mathcal C}_7$ to be
$${\mathcal C}_7= 2^63^2n^2t_0^{-4}({\mathcal M}_2+3{\mathcal M}_1+2^33\sqrt{3n}\,{\mathfrak R}_0/t_0)\ .
$$
\item
Then, denote by ${\mathcal R}_3$ the positive real number  
$${\mathcal R}_3 = \min (\frac{{\mathcal R}_0}{2},\frac{t_0}{5\sqrt{n\,{\mathcal M}_1}}, \frac{1}{\sqrt{{\mathcal C}_7}})\,.
$$ 
\end{enumerate}
\end{defn}
\begin{rem}\label{need33}
For further need, 
observe (as $ t_0 <1$) one has 
$\frac{t_0}{5\sqrt{n\,{\mathcal M}_1}}<\frac{1}{\sqrt{2\,{\mathcal M}_1}}<1$ and $\frac{1}{\sqrt{{\mathcal C}_7}}<\frac{1}{\sqrt{2{\mathcal C}_2}}$ (indeed $\frac{1}{\sqrt{{\mathcal C}_7}}\leq\frac{t_0^{5/2}}{2^{9/2}3^{7/4}n^{5/4}\sqrt{{\mathfrak R}_0}}<\frac{\sqrt{3}}{2\sqrt{{\mathfrak R}_0}}=\frac{1}{\sqrt{2{\mathcal C}_2}}\,$, see lemma \ref{lemma 2nd}) and this allows to apply lemma \ref{lemma 3rd}.
 \end{rem}

We now state the main theorem in this section.

\begin{thm} \label{theorem 3} Given $t_0>0$ and ${\mathcal C}_7$ as in definition {\rm\ref{constanteC_7R_3}}, for any 
$\rho\in]0,{\mathcal R}_3]\,$ one has ${\mathcal C}_7\,\rho^2<1$ 
{\em and}, for any $x$ in 
$W$, each Riemannian barycentric simplex $\hat\tau\subset B(x,\rho)$ 
of thickness $t_g(\hat\tau)\geq t_0$ admitting a Euclidean realization $(\hat \tau,\hat g)$ verifies, for 
all $q\in \hat\tau$ and $u\in T_q M$
\begin{align*}
(i) &   
 \hskip6mm & \vert\,\Vert\,u\,\Vert_{\hat g}^2-\Vert\,u\,\Vert_g^2\,\vert & 
 \leq \Vert\,u\,\Vert_g^2 \  {\mathcal C}_7\,\rho^2\ ; \\ 
 (i') &   
 \hskip6mm & \vert\,\Vert\,u\,\Vert_{\hat g}-\Vert\,u\,\Vert_g\,\vert & 
 \leq \Vert\,u\,\Vert_g \  {\mathcal C}_7\,\rho^2\ ; \\
(ii) & \hskip6mm & H=\sup_{q\in\hat\tau}\Vert\,h_q\,\Vert_g & \leq 
{\mathcal C}_7\,\rho\  ; \\
(iii) &\hskip6mm & \Vert\,({\mathcal P}_\gamma^{\hat g}-{\mathcal 
P}_\gamma^g)(u)\,\Vert_g 
& \leq \Vert\,u\,\Vert_g\ \hbox{\rm length}_g (\gamma)\ {\mathcal C}_7\,\rho\ ,
\end{align*}
where
$\gamma : [0, 1]\rightarrow \hat\tau$ is any differentiable curve starting at 
$q\,$.
\end{thm}
\begin{rem}
The rest of the section is devoted to the proof of this theorem. Most of the work
done consists in controlling the tensor $h\,$ and apply proposition 
\ref{proposition 1st} 
to the $g$ and $\hat g$-parallel translations. Around the next lemma, see also \cite[section 8]{B-K}.
\end{rem}

\begin{defn} \label{formule1022}
Let $q\in \hat\tau$, $v\!=\!\sum_{i=0}^n\lambda_i\,v_i \!\in\! \hat\tau_x\!\subset\! T_xM$ with
$v_i\!=\!\exp_x^{-1}p_i$ (definition \ref{D.euclidianisable}), verify $\phi(v)\!=\! q$.
Define the $g$-symmetric operator $A_q$
\begin{equation*} 
A_q=\sum_{i=0}^n\lambda_i\ {}^\sharp\hbox{\rm hess}_q\, 
(\frac{d^2(p_i,\cdot)}{2})\,.
\end{equation*}
\end{defn}

\begin{lem} \label{inv} 
The operator $A_q$ is invertible for all $q\in\hat\tau\,$. Set $v_{ij}=v_j-v_i\in T_xM$ and 
$w_{ij}(q)=w_j(q)-w_i(q)
=\exp_q^{-1}p_j-\exp_q^{-1}p_i\in T_qM\,$. One has for any $v\in \hat\tau_x$ such that $\phi(v)=q$
\begin{equation} \label{formule}
d\phi(v)\,v_{ij}=A_q^{-1}\,w_{ij}(q)\ .
\end{equation}
Furthermore, one has
\begin{equation}\label{inv333} 
\Vert\,A_q^{-1}-\hbox{\rm Id}_q\,\Vert\leq 2\,{\mathcal M}_1\,\rho^2\  .
\end{equation}
\end{lem}

\begin{proof} 
So, we have $x\in W$, 
$\hat\tau \subset B(x,\rho)$, where $\rho \leq {\mathcal R}_3 \leq {\mathcal R}_0/2 \leq {\mathcal R}_W/2$ (here $\rho\leq{\mathcal R}_W$ is enough), and $v_i=\exp_x^{-1} p_i\,$. 
We know (by theorem \ref{theorem 1st}) 
a uniform 
bound ${\mathcal M}_1$ to exist on $W\subset(M,g)$, such that, for any $x\in W$ 
and $p,q\in B(x,\rho)\subset B(x,{\mathcal R}_W)$, the following holds
\begin{equation*}
\Vert\,{}^{\sharp}\hbox{\rm hess}_q\,(\frac{d^2(p,\cdot)}{2})-
\hbox{\rm Id}_q\,\Vert\leq {\mathcal M}_1\,\rho^2\  .
\end{equation*}
Thus, for any ${\bf \lambda}$ with $\sum_{i=0}^n\lambda_i\!=\!1$ and $\lambda_i\!\geq \!0\,$, writing $f_i\!=\!d^2(p_i,\cdot)/2$
\begin{equation}\label{invert A}
\Vert\,A_q-\hbox{\rm Id}_q\,\Vert\leq \sum_{i=0}^n
\lambda_i\Vert\,{}^{\sharp}\hbox{\rm hess}_q\,f_i-\hbox{\rm Id}_q\,\Vert\leq  {\mathcal M}_1\,\rho^2\  .
\end{equation}
Since $\rho \leq {\mathcal R}_3 \leq \frac{1}{\sqrt{2{\mathcal M}_1 }}\,$ (remark \ref{need33}), one has 
$\Vert\,A_q-\hbox{\rm Id}_q\,\Vert \leq 1/2$ and it follows that $A_q$ is invertible.

For sake of notational simplicity, we prove formula (\ref{formule}) in the 
case $i=0$ and $j=1\,$. At 
any point $v\in \hat\tau_x\,$, the constant field $v_{01}$ is the $\alpha$-derivative at $\alpha=0$ of the 
rectilinear curve 
\begin{equation*} 
v(\alpha)=v+\alpha\,(v_1-v_0)=(\lambda_0-\alpha)\,v_0+(\lambda_1+\alpha)\,
v_1+\sum_{i=2}^n\lambda_i\,v_i\  .
\end{equation*}
The image by $\phi$ of this line $v(\alpha)$ is a curve $q(\alpha)$ (we study 
$q=\phi(v)$) 
given by the 
implicit equation $F_\alpha(q(\alpha))=0\,$, where
\begin{equation*} 
F_\alpha(\cdot)=(\lambda_0-\alpha)\,df_0(\cdot)+(\lambda_1+\alpha)\,
df_1(\cdot)+\sum_{i=2}^n\lambda_i\,df_i(\cdot)\,.
\end{equation*}
Dualizing in the metric $g$ the form $F_\alpha$ into a field $F_\alpha^{\sharp}$ and 
taking the derivative
of this field in $\alpha$ at $\alpha=0\,$, one gets
\begin{equation*} 
0=\nabla_{\frac{\partial q}{\partial 
\alpha}(0)}(F_\alpha^{\sharp}(q(\alpha)))=\exp_{q}^{-1}p_0-
\exp_{q}^{-1}p_1+A_q\,(\frac{\partial q}{\partial \alpha}(0))\,,
\end{equation*}
and this proves formula (\ref{formule}) since $A_q^{-1}$ exists.

Finally, one has 
\begin{equation*}
\Vert\,A_q^{-1}\,\Vert=\bigg{(}\inf_{u\in T_qM\setminus\{0\}}
\frac{\Vert\,A_q (u)\,\Vert_q}{\Vert\, u\,\Vert_q}\bigg{)}^{-1}\,.
\end{equation*}
As one knows by (\ref{invert A})
\begin{equation*}
\inf_{u\not=0}\frac{\Vert\,A_q (u)\,\Vert_q}{\Vert\, u\,\Vert_q}\geq
1-{\mathcal M}_1\,\rho^2 \hskip2mm 
\end{equation*}
one concludes as wished (recall $2{\mathcal M}_1\,\rho^2$ is $<1$)
\begin{equation*}
\Vert\,A_q^{-1}-\hbox{\rm Id}_q\,\Vert\leq 
\Vert\,A_q-\hbox{\rm Id}_q\,\Vert\ \Vert\,A_q^{-1}\,\Vert
\leq\frac{{\mathcal M}_1\,\rho^2}{1-{\mathcal M}_1\,\rho^2}\leq
2\,{\mathcal M}_1\,\rho^2\   . \qedhere
\end{equation*}
\end{proof}

The next proposition is central in proving theorem \ref{theorem 3}.
\begin{prop} \label{lemma 2} 
Given $t_0\!>\!0$ as in definition {\rm\ref{constanteC_7R_3}}, there exists ${\mathcal C}_5$ such that, 
for any positive $\rho\!\leq\!{\mathcal R}_3$ (here $\rho\!\leq\!{\mathcal R}_2/2$ matters), for any $x\!\in\! W$, each Riemannian barycentric simplex $\hat\tau\!=\! \hat\tau_{\bf p}\! \subset\! B(x,\rho)$ 
of thickness $t_g(\hat\tau)\geq t_0$ verifies, for 
any $q\in\hat\tau$ and $i,j\!=\!0,1,\dots,n$ with $i<j\,$ 
\begin{equation*}
\Vert\,\nabla_{A_q^{-1}w_{ij}}A_\cdot^{-1}w_{ij}(\cdot)\,\Vert_q\leq 
{\mathcal C}_5\,\delta_q^2(\hat\tau_q)\,\rho\  ,
\end{equation*}
where $\delta_q(\hat\tau_q)$ is the $g_q$-diameter of $\hat\tau_q\,$ (thus $\delta_q(\hat\tau_q)\leq 2\rho$).
\par
Actually, ${\mathcal C}_5= 2^3\,({\mathcal M}_2+3{\mathcal M}_1+2^3\,3\sqrt{3n}\,{\mathfrak R}_0/t_0)$ works.
\end{prop}
 
The proof is divided in several lemmas and relies on some computations involving 
variations by geodesics we now introduce. As before, we specialize $i=0$ and $j=1$ and 
look at the two variations 
\begin{equation*}
\begin{split}
 c_0(s,t)=\exp_{p_0}\,(s \exp_{p_0}^{-1}\circ\exp_q (t\exp_q^{-1}p_1))\,; \\
 c_1(s,t)=\exp_{p_1}\,(s \exp_{p_1}^{-1}\circ\exp_q (t\exp_q^{-1}p_0))\,.
\end{split}
 \end{equation*}

Recall $w_i(q)=\exp_q^{-1}p_i$ (definition \ref{genthick}).
 
 \begin{lem}\label{lemma 3} 
 The following equalities are true
 \begin{align*} 
 (i) &  \ \ w_0(c_0(s,t)) = -s\,\frac{\partial c_0}{\partial s}(s,t) \,; 
\ w_1(c_1(s,t))=-s\,\frac{\partial c_1}{\partial s}(s,t)  \,; \\
 (ii)  & \ \ w_1(c_0(1,t)) = (1-t)\,\frac{\partial c_0}{\partial t}(1,t) \,; 
   \  w_0(c_1(1,t))=(1-t)\,\frac{\partial c_1}{\partial t}(1,t)  \, ; \\
(iii) & \ (\nabla_{w_0}  w_0)_{\mid c_0(s,t)} = -w_0(c_0(s,t)) \,  ;  
\ (\nabla_{w_1}  w_1)_{\mid c_1(s,t)}=-w_1(c_1(s,t))  \,; \\
 (iv)  &  \ (\nabla_{w_0}  w_1)_{\mid c_1(1,t)} =
 -(1-t)\,\frac{\partial^2 c_1}{\partial t\partial s}(1,t) = 
 -(1-t)\,\frac{\partial Y_1}{\partial s}(1,t)  \,; \\
 (v)  & \ (\nabla_{w_1}  w_0)_{\mid c_0(1,t)}  =-(1-t)\,
  \frac{\partial^2 c_0}{\partial t\partial s}(1,t) 
= -(1-t)\,\frac{\partial Y_0}{\partial s}(1,t) \,,
 \end{align*}
  where $Y_0$ and $Y_1$ are the Jacobi fields defined by $c_0$ and $c_1\,$.
 \end{lem}
 
\begin{proof} 
 If $c(s)$ is the geodesic from $p$ to $q\,$, as $\  
d\exp_p(\exp_p^{-1}q)(\exp_p^{-1}q)=\frac{\partial c}{\partial s}(1)$, this 
equals $-\frac{\partial \tilde c}{\partial s}(0)$ for the 
reversed geodesic $\tilde c(s)=c(1-s)$ and one has 
$\frac{\partial c}{\partial s}(1)=-\exp_q^{-1}p\,$. 
The curve $\gamma(\sigma)=c(s\sigma)=\exp_p(\sigma s \exp_p^{-1}q)$ is the geodesic from 
$\gamma(0)=p$ to $\gamma(1)=c(s)$. Thus 
\begin{equation*}
s\,\frac{\partial c}{\partial s}(s)\!=\!
d\exp_p(s\exp_p^{-1}q)(s\exp_p^{-1}q)\!=\!
\frac{\partial \gamma}{\partial \sigma}(1)
\!=\!-\exp_{\gamma(1)}p\!=\!-\exp_{c(s)}^{-1}p\,.
\end{equation*}
Setting $c=c_0$ and $c=c_1$ gives $(i)$.

As for $(ii)$, the initial velocity vector of the geodesic running from
$c_0(1,t)$ to $p_1$ is
$w_1(c_0(1,t))=\exp^{-1}_{c_0(1,t)}p_1\,$. But 
$\sigma\in[0,1]\rightarrow c_0(1, \sigma)$ is the 
geodesic 
from $q$ to $p_1\,$. By uniqueness, the vectors $\frac{\partial c_0}{\partial t}(1,t)$ 
and $w_1(c_0(1,t))$ are colinear and point in the same direction : checking the lengths 
gives $w_1(c_0(1,t))=(1-t)\, \frac{\partial c_0}{\partial t}(1,t)$.

Formulas $(iii)$ follow directly from formulas $(i)$  and from the 
definitions of $c_0$ and $c_1\,$. Formulas $(iv)$ and $(v)$  
follow by taking derivative in $t$ in 
both formulas $(i)$, while $s=1$ is fixed, and by using  
formulas $(ii)$ to relate this derivative to the derivative along $w_1$ in the first case 
(to $w_0$ in the second).
\end{proof}

 \begin{lem} \label{lemma 5} One has
\begin{equation*} 
\frac{\partial Y_1}{\partial s}(0,0)=
d\exp_{p_1}^{-1}(q)\,(w_0(q))\,, \ \ \  
\frac{\partial Y_0}{\partial s}(0,0)=
d\exp_{p_0}^{-1}(q)\,(w_1(q))\,.
\end{equation*}
\end{lem}

\begin{proof} Compute (recall $w_0(q)=\exp_q^{-1}p_0$)
\begin{multline*}\  
\frac{\partial c_1}{\partial s}(0,t)
\!=\!\exp_{p_1}^{-1}\circ\exp_q(t\,w_0(q))\, ,\ \ 
\frac{\partial Y_1}{\partial s}(0,t)\!=\!
\frac{\partial^2 c_1}{\partial t\partial s}(0,t)=\\
=(d\exp_{p_1}^{-1}(\exp_q(t \,w_0(q)))\circ d\exp_q(t \,w_0(q)))\,(w_0(q))\, , 
\end{multline*}
and put $t=0\,$.
The second equality 
follows along analogous lines.  
\end{proof}

Write ${\mathcal W}_i$ for $w_{0i}=w_i-w_0\,$, in details 
 ($p_0,\dots,p_n \in B(x,\rho)$ being the vertices of $\hat\tau$)
\begin{equation} \label{Wi}
{\mathcal W}_i(q)=w_{0i}(q)=\exp^{-1}_q p_i-\exp^{-1}_qp_0\ .
\end{equation}

\begin{lem} \label{lemma 6} 
There exists a constant ${\mathcal C}_3$
such that, 
for any $x\in W$, any $\rho\in]0,{\mathcal R}_3]\,$,
each Riemannian barycentric simplex  
$\hat\tau  =\hat\tau _{\bf p} \subset B(x,\rho)$, one has for any $q\in\hat\tau\,$ and $i,j=0,1,\dots,n$
\begin{equation*}
\Vert\,\nabla_{{\mathcal W}_i}  {\mathcal W}_j\,\Vert_q\leq 
{\mathcal C}_3\ \delta_q^3(\hat\tau_q)\   .
\end{equation*}
The constant ${\mathcal C}_3= 2^23^{3/2}\,{\mathfrak R}_0$ fits. 
\end{lem}

Observe that if $i$ or $j=0\,$, there is nothing to prove since  
${\mathcal W}_0 \equiv 0\,$.

\begin{proof} 
From previous lemma \ref{lemma 3}, one has at $q=c_1(1,0)$
\begin{multline*}
(\nabla_{w_0}  (w_1-w_0))_{\mid q}=-\,\frac{\partial Y_1}{\partial s}(1,0)+w_0(q)= \\
= -\,\frac{\partial Y_1}{\partial s}(1,0)+{\mathcal P}\,
\frac{\partial Y_1}{\partial s}(0,0)-{\mathcal P}\,
\frac{\partial Y_1}{\partial s}(0,0)+w_0(q)\,,
\end{multline*}
where ${\mathcal P}$ stands for the $g$-parallel translation along $c_1(\cdot,0)$ from 
$p_1$ to $q\,$. 
Consider the variation 
\begin{equation*} 
\tilde c_1(s,t)=\exp_{p_1}
(s(\frac{\partial c_1}{\partial s}(0,0)+t\,d\exp_{p_1}^{-1}(q)(\exp_q^{-1} p_0)))\,,
\end{equation*}
This variation has the same Jacobi field as $c_1$ along the curve 
$c_1(s,0)=\tilde c_1(s,0)$: indeed, a computation shows 
${\partial \tilde c_1}/{\partial t}(s,0)=Z_1(s,0)=Y_1(s,0) $. Using 
definition \ref{bareucl} and 
proposition \ref{proposition A} (making only $C=1/2$), get
\begin{equation*}
\Vert\,d\exp_{p_1}^{-1}(q)(\exp_q^{-1} p_0)\,\Vert\leq 
\Vert\,d\exp_{p_1}^{-1}(q)\,\Vert\ \Vert\,w_0\,\Vert_q \leq 2\ 
\Vert\,w_0\,\Vert_q\ , 
\end{equation*}
and (by lemma \ref{lemma 1st})
\begin{equation*}
\Vert\,\frac{\partial Y_1}{\partial s}(1,0)-{\mathcal P}\,
\frac{\partial Y_1}{\partial s}(0,0)\,\Vert_q\leq 2\,{\mathcal C}_1\  
\Vert\,w_1\,\Vert_q^2\ 
 \Vert\,w_0\,\Vert_q\  .
\end{equation*}
With the help of lemma \ref{lemma 5} and previous lemma \ref{lemma 3rd} (see remark \ref{bdvargeod}), one has a uniform 
${\mathcal C}_2$ on $W\subset(M,g)$ such that (recall $\rho$ is smaller than the 
convexity 
radius in $W$ and satisfies 
${\mathcal M}_1\,\rho^2\leq 1/2$)
\begin{equation*}
\Vert {\mathcal P}  \frac{\partial Y_1}{\partial s}(0,0)- w_0\Vert_q =
\Vert\,{\mathcal P}\,d\exp^{-1}_{p_1}(q) (w_0)-w_0\,\Vert_q\leq 
2\,{\mathcal C}_2\  \Vert\,w_1\,\Vert_q^2\  \Vert\,w_0\,\Vert_q\  .
\end{equation*}
We thus obtain, putting together the previous lines
\begin{equation} \label{aster1}
\Vert\, \nabla_{w_0}  (w_1-w_0) \,\Vert_q\leq 2\,({\mathcal C}_1+ {\mathcal C}_2)\  
\Vert\,w_1\,\Vert_q^2\  \Vert\,w_0\,\Vert_q\  ,
\end{equation}
so that, going from $0$ and $1$ to $i$ and $j=0,1,\dots,n$
\begin{equation*}
\Vert \nabla_{(w_j-w_i)}  (w_j-w_i) \Vert_q\leq 2\,({\mathcal C}_1+ {\mathcal C}_2)\Vert w_i\Vert_q
\Vert w_j\Vert_q(\Vert w_i \Vert_q+\Vert w_j\Vert_q)\,.
\end{equation*}
Now, as ${\mathcal W}_i=w_{0i}\,$, writing
\begin{equation*}
\nabla_{{\mathcal W}_i}  {\mathcal W}_j+\nabla_{{\mathcal W}_j}  {\mathcal W}_i=
\nabla_{{\mathcal W}_i}  {\mathcal W}_i+\nabla_{{\mathcal W}_j}  {\mathcal W}_j-
\nabla_{({\mathcal W}_i-{\mathcal W}_j)}  ({\mathcal W}_i-{\mathcal W}_j)
\end{equation*}
we get from the last inequalities the following ones
\begin{equation*}
\Vert \nabla_{{\mathcal W}_i}  {\mathcal W}_j+\nabla_{{\mathcal W}_j} 
{\mathcal W}_i \Vert_q
\leq 12\,({\mathcal C}_1+ {\mathcal C}_2)\  \delta_q^3(\hat\tau_q)\,.
\end{equation*}
But, one may also write
\begin{gather*}
\nabla_{{\mathcal W}_i}  {\mathcal W}_j-\nabla_{{\mathcal W}_j}  {\mathcal W}_i 
=\nabla_{w_i}  (w_j-w_0)-\nabla_{w_j}  (w_i-w_0)-\nabla_{w_0}  (w_j-w_i)\,,
\\
\nabla_{w_i}  (w_j-w_0)=\nabla_{w_i}  (w_j-w_i)+
\nabla_{w_i}  (w_i-w_0)\,,
\end{gather*}
and the two corresponding equations obtained by cycling in $0,i,j\,$. This, with the help 
of (\ref{aster1}), allows to write 
\begin{equation*}
\Vert \nabla_{{\mathcal W}_i}  {\mathcal W}_j-\nabla_{{\mathcal W}_j}  
{\mathcal W}_i \Vert_q
\leq 12\,({\mathcal C}_1+ {\mathcal C}_2)\  \delta_q^3(\hat\tau_q)\,.
\end{equation*}
Remark \ref{mathcalC1,2} says ${\mathcal C}_3=12\,\sqrt{3}\,{\mathfrak R}_0>12\,({\mathcal C}_1+ {\mathcal C}_2)\,$, hence the result .
\end{proof}

\begin{lem} \label{lemma 7}
Given $t_0>0$, there exists a real ${\mathcal C}_4$ such that, 
for any positive $\rho\leq{\mathcal R}_3\,$, any $x\in W\,$, each Riemannian barycentric simplex $\hat\tau= \hat\tau_{\bf p} \subset B(x,\rho)$ 
of thickness $t_g(\hat\tau)\geq t_0$ verifies, for 
any $q\in\hat\tau$ and $i,j=0,1,\dots,n$ with $i<j\,$ 
\begin{equation*}
\Vert\,A_q^{-1}(\nabla_{A_q^{-1}w_{ij}}  w_{ij})\,\Vert_q\leq 
{\mathcal C}_4\ \delta_q^2(\hat\tau_q)\,\rho\  .
\end{equation*}
One may choose ${\mathcal C}_4=  2^6\,3\,\sqrt{3n}\,{\mathfrak R}_0/t_0\,$.
\end{lem}

\begin{proof} 
In the basis ${\mathcal W}_l\,$, $l=1,\dots,n\,$, write 
$A_q^{-1}w_{ij}=\sum_{l=1}^n\alpha_l^{ij}\  {\mathcal W}_l\,$ and
\begin{multline*}
\Vert\,\nabla_{A_q^{-1}w_{ij}}  w_{ij}\,\Vert \leq \sum_l\vert\alpha_l^{ij}\vert\  
\Vert\,\nabla_{{\mathcal W}_l}  w_{ij}\,\Vert = \\
=\sum_l\vert\alpha_l^{ij}\vert\  
\Vert\,\nabla_{{\mathcal W}_l}  {\mathcal W}_j-\nabla_{{\mathcal W}_l}  
{\mathcal W}_i\,\Vert\leq\sum_l\vert\alpha_l^{ij}\vert\  
(\Vert\,\nabla_{{\mathcal W}_l}  {\mathcal W}_j\,\Vert+\Vert\,\nabla_{{\mathcal W}_l}  
{\mathcal W}_i\,\Vert)\  .
\end{multline*}
At any point $q\in\hat\tau\,$, get with lemma \ref{lemma 6} and 
Schwarz's inequality
\begin{equation*}
\Vert\,\nabla_{A_q^{-1}w_{ij}}  w_{ij}\,\Vert_q \leq 2\,\sqrt{n}\,
{\mathcal C}_3\,\delta_q^3(\hat\tau_q)\,\sqrt{\sum_l\vert\alpha_l^{ij}\vert^2}\ \,.
\end{equation*}
The thickness $t_g(\hat\tau)$ is $\geq t_0>0$ (see 
definition \ref{genthick}), thus $t_q(\hat\tau_q)\geq t_0\,$. Recalling 
$r(\hat\tau_q)\,\sqrt{\sum_l\vert\alpha_l^{ij}\vert^2}\leq
\Vert\sum_{l=1}^n\alpha_l^{ij}\  {\mathcal W}_l\Vert$ (see lemma \ref{1b}), we get
\begin{equation}\label{ast3}
\Vert\,\nabla_{A_q^{-1}w_{ij}}  w_{ij}\,\Vert_q \leq 
\frac{2\,\sqrt{n}\,{\mathcal C}_3\,\delta_q^2(\hat\tau_q)}{t_0}\,
\Vert\,A_q^{-1}w_{ij}\,\Vert\ .
\end{equation}
Lemma \ref{inv} gives
\begin{equation} \label{estime1}
\Vert\,A_q^{-1}w_{ij}\,\Vert_q\leq (1+2\,{\mathcal M}_1\,\rho^2)\,\Vert\,w_{ij}\,\Vert_q 
\leq 2\,\delta_q(\hat\tau_q)\, .
\end{equation}
Putting (\ref{estime1}) into inequality (\ref{ast3}) gives
$
\Vert\,\nabla_{A_q^{-1}w_{ij}}  w_{ij}\,\Vert_q \leq 
\frac{4\,\sqrt{n}\,{\mathcal C}_3}{t_0}\,\,\delta_q^3(\hat\tau_q)\,.
$
As
${\mathcal C}_4=16\,\sqrt{n}\,{\mathcal C}_3/t_0\,$, the inequalities $\Vert A_q^{-1}\Vert\leq 2$ and $\delta_q(\hat\tau_q)\leq 2\rho$ complete the proof, recalling ${\mathcal C}_3=12\,\sqrt{3}\,{\mathfrak R}_0$ (lemma \ref{lemma 6}).
\end{proof}
\noindent {\bf Completing the proof of proposition \ref{lemma 2}.}
\begin{proof} 
Theorem \ref{theorem 2nd} gives a uniform 
${\mathcal M}_2$ on $W\subset(M,g)$ such that, for any $\rho\leq{\mathcal R}_3$ 
and $x\in W\,$,
for any $p,q\in B(x,\rho)$, one has the inequality 
$\Vert(\!\nabla{}^\sharp\hbox{\rm hess})_q ({d^2(p,\cdot)}/{2})\Vert
\leq {\mathcal M}_2\rho\,.$
As
$A_q\!=\!\sum_k\!\lambda_k(\phi^{-1}(\cdot)){}^{\sharp}\hbox{\rm hess}_q ({d^2(p_k,\cdot)}/{2})$ (definition \ref{formule1022}), setting for $v\in \hat\tau_x$ with $\phi(v)\!=\!q$ (proof of lemma \ref{inv}) $ 
v_{ij}(\alpha)\!=\!v\!+\!\alpha(v_j\!-\!v_i)\!=\!\sum_{k}^n\!\lambda_k (v_{ij}(\alpha))\,v_k
$, the $\alpha$-derivative at $\alpha\!=\!0$ of $\lambda_k(v_{ij}\!(\alpha))\!=\!\lambda_i\!-\!\alpha$ if $k\!=\!i$, $=\!\lambda_j\!+\!\alpha
$ if $k\!=\!j$, $=\!\lambda_k$ if $k\!\not=\!i,j$, is the derivative of $\lambda_k(\phi^{-1}(\cdot))$ in the direction $
d\phi(v)v_{ij}\!=\!A_q^{-1}w_{ij}(q)
$ (see (\ref{formule})). Thus
\begin{multline*}
(\nabla_{A_q^{-1}  w_{ij}}A)(A_q^{-1}  w_{ij}) 
\!= \!\sum_{k=0}^n\!\lambda_k\circ\phi^{-1}\!(q)(\nabla_{A_q^{-1}  w_{ij}}\!{}^{\sharp}\hbox{\rm hess}_q 
\frac{d^2(p_k,\cdot)}{2})(A_q^{-1}  w_{ij})+\hfill\\
\hfill+ ({}^{\sharp}\hbox{\rm hess}_q 
\frac{d^2(p_j,\cdot)}{2} 
-{}^{\sharp}\hbox{\rm hess}_q 
\frac{d^2(p_i,\cdot)}{2})\,(A_q^{-1}  w_{ij})\, .
\end{multline*}
The $\lambda_k$ sum up to $1$ and are $\geq0\,$, so by 
(\ref{estime1}), applying theorem \ref{theorem 1st} while {\it choosing the $x$ {\em there} (we name it $y$ {\em here}) in $\hat\tau$} (the choice ${\mathcal R}_3\leq {\mathcal R}_0/2$ done in definition \ref{constanteC_7R_3} is used here to guarantee
$\hat\tau\subset B(y,2\rho)$ with $2\rho\leq{\mathcal R}_0$), applying proposition \ref{proposition A} (getting $\delta(\hat\tau)\leq3\,\delta_q(\hat\tau_q)/2$) and using $\delta(\hat\tau)\leq2\,\rho\,,\,\delta_q(\hat\tau_q)\leq2\,\rho\,$, one gets
\begin{multline}\label{estimetralala}
\Vert\,(\nabla_{A_q^{-1}  w_{ij}}A)\,(A_q^{-1}  w_{ij})\Vert_q\leq \\ 
\leq \sum_{k=0}^n\lambda_k\,\Vert\,\nabla({}^{\sharp}\hbox{\rm hess}_q 
\frac{d^2(p_k,\cdot)}{2})\Vert\, \Vert\,A_q^{-1}  w_{ij}\Vert^2_q+\hfill\\
\hfill+\Vert ({}^{\sharp}\hbox{\rm hess}_q 
\frac{d^2(p_j,\cdot)}{2} 
-{}^{\sharp}\hbox{\rm hess}_q 
\frac{d^2(p_i,\cdot)}{2})\,(A_q^{-1}  w_{ij})\Vert\leq \\
\leq {\mathcal M}_2\,\rho\,\Vert\,A_q^{-1}  w_{ij}\Vert^2_q+2
\,{\mathcal M}_1\,\delta^2(\hat\tau)\,\Vert\,A_q^{-1}  w_{ij}\Vert_q
\leq\\\leq (4\,{\mathcal M}_2+12\,{\mathcal M}_1)\ \rho\ \delta_q^2(\hat\tau_q)\, .
\end{multline}
But, since 
$\nabla A^{-1}=-A^{-1}\,(\nabla A)\,A^{-1}\,$, one has
\begin{equation*} 
\Vert\nabla_{A_q^{-1}w_{ij}}  A_\cdot^{-1}  w_{ij}\Vert_q\leq\Vert A_q^{-1}(\nabla_{A_q^{-1}  w_{ij}}A)(A_q^{-1}  w_{ij})\Vert_q +
\Vert A_q^{-1}\,\nabla_{A_q^{-1}  w_{ij}} w_{ij}\Vert_q\  .
\end{equation*}
Estimate (\ref{estimetralala}) and the inequality stated in lemma \ref{lemma 7}, 
thus establish for ${\mathcal C}_5= 8\,{\mathcal M}_2+24\,{\mathcal M}_1+{\mathcal C}_4\,$ and $i,j=0,1,\dots,n\,$ 
 (use $\Vert\,A_q^{-1}\,\Vert\leq 2$)
\begin{equation*}
 \Vert\,\nabla_{A_q^{-1}w_{ij}}  A_\cdot^{-1}  w_{ij}\,\Vert_q\leq 
 {\mathcal C}_5\,\delta_q^2(\hat\tau_q)\,\rho\  .\qedhere
 \end{equation*}
\end{proof}
\begin{defn} \label{bartauq} A Euclidean simplex $(\hat\tau_q,g_q)$ whose vertices are the points
$\exp_q^{-1}p_0,\dots,\exp_q^{-1}p_n$ satisfies definition \ref{genthick}. Call $(\overline{\tau}_q,g_q)$ the simplex with vertices 
${\mathcal A}_i\!=\!A_q^{-1}  {\mathcal W}_i(q)\!\in\!(T_qM,g_q)$, where (according to (\ref{Wi})) ${\mathcal W}_i(q)\!=\!\exp_q^{-1}p_i\!-\!\exp_q^{-1}p_0=w_i-w_0\,$.
\end{defn}
The following geometrical link between $\hat\tau_q$ and $\overline{\tau}_q$ will be needed
\begin{lem} \label{tricky} Given $t_0>0$ as in definition {\rm\ref{constanteC_7R_3}}, 
for any positive $\rho\leq{\mathcal R}_3\,$ (thus 
$\rho < \frac{t_0}{5\sqrt{n\,{\mathcal M}_1}}$), any $x\in W$, each Riemannian barycentric simplex $\hat\tau= \hat\tau_{\bf p} \subset B(x,\rho)$ 
of thickness $t_g(\hat\tau)\geq t_0$ is such that, for 
any $q\in\hat\tau\,$, the 
Euclidean 
simplex $\overline{\tau}_q$ (see definition {\rm\ref{bartauq}}) verifies
\begin{equation*} 
r_q(\overline{\tau}_q) \geq 
\frac{r_q(\hat\tau_q)}{2}\hskip4mm\hbox{and}\hskip4mm
t_q(\overline{\tau}_q) \geq \frac{t_q(\hat\tau_q)}{2}\geq \frac{t_0}{2}\  .
\end{equation*}
\end{lem}

\begin{proof} 
Going back to lemma \ref{inv} (we know  
$2\,{\mathcal M}_1\,\rho^2<1$), one has
\begin{equation*}
\forall q\in \hat\tau\hskip2cm
\Vert\,A_q^{-1}  {\mathcal W}_i-{\mathcal W}_i\,\Vert_q\leq 2\,{\mathcal M}_1\rho^2\,
\Vert\,{\mathcal W}_i\,\Vert_q\  .
\end{equation*}
Applying lemma \ref{1d} (take $\nu=2\,{\mathcal M}_1\rho^2$), one gets from its second 
part 
\begin{equation*}
\vert\,r_q(\overline{\tau}_q)-r_q(\hat\tau_q)\,\vert\leq
\frac{9\,n\,{\mathcal M}_1}{t_0^2}\,\rho^2\,r_q(\hat\tau_q)\, ,
\end{equation*}
and the first claim follows as $\rho<{t_0}/{(5\sqrt{n\,{\mathcal M}_1})}\,$. 
Another application of the same lemma implies that the second claim follows if
\begin{equation*} 
\frac{2\,{\mathcal M}_1\,\rho^2}{1-2\,{\mathcal 
M}_1\,\rho^2}(1+\frac{9\,n}{2\,t_0^2})\leq\frac{1}{2}\  ,
\end{equation*}
true if 
$\ t_0^2\,{\mathcal M}_1\,\rho^2+3\,n\,{\mathcal M}_1\,\rho^2\leq {t_0^2}/{6}\,$, thus if
$4\,n\,{\mathcal M}_1\,\rho^2\leq {t_0^2}/{6}\,$ (the thickness of a simplex is always $\leq1\,$)
and {\it a fortiori} if $\rho\leq {t_0}/{(5\sqrt{n\,{\mathcal M}_1})}\,$.
\end{proof}

{\bf Proof of theorem \ref{theorem 3}.}

\begin{proof}
We start by proving $(ii)$. Given $v\in\hat\tau_x$ such that $\phi(v)=q\,$, the construction of the metric $\hat g$ and the 
definition of 
the fields ${\mathcal A}_i=A_q^{-1}  {\mathcal W}_i=d\phi(v)(v_i-v_0)$ imply that 
$\hat\nabla_{{\mathcal A}_i}{\mathcal A}_j\equiv 0$ for all 
$i,j=1,\dots,n\,$ (this is true for 
$v\in\hat\tau_x$, see lemma \ref{inv}). 
So, the difference connection tensor $h$ reads in the frame ${\mathcal A}_i$ for 
$i=1,\dots,n$ 
\begin{equation*}
h({\mathcal A}_i,{\mathcal A}_j)=-\nabla_{{\mathcal A}_i}{\mathcal A}_j\ .
\end{equation*}
By polarization, we get for all $i,j=1,\cdots,n$
\begin{equation*}
2\,h({\mathcal A}_i,{\mathcal A}_j)=h({\mathcal A}_i,{\mathcal A}_i)+h({\mathcal A}_j,{\mathcal A}_j)-h(A_q^{-1} w_{ij},A_q^{-1} 
w_{ij})\, .
\end{equation*}
From proposition \ref{lemma 2} we have
\begin{equation*} 
\Vert\,h({\mathcal A}_i,{\mathcal A}_j)\,\Vert\leq \frac{3}{2}\,
{\mathcal C}_5\,\rho\,\delta^2_q(\hat\tau_q)\  .
\end{equation*}
Taking ${\mathcal U}=\sum_{i=1}^n\lambda_i\,{\mathcal A}_i\,$, once again Schwarz's inequality plus lemma \ref{1b} give, with $r_q(\overline{\tau}_q) \geq 
{r_q(\hat\tau_q)}/{2}$ taken from lemma \ref{tricky}
\begin{multline*}
\Vert\,h({\mathcal U},{\mathcal U})\,\Vert\leq
\sum_{i,j=1}^n(\vert\lambda_i\vert\,\vert\lambda_j\vert)\,
\frac{3\,{\mathcal C}_5\,\rho\ \delta^2_q(\hat\tau_q)}{2}
\leq(\sum_{i=1}^n\lambda_i^2)\,\frac{3\,n\,{\mathcal C}_5\,\rho\ 
\delta^2_q(\hat\tau_q)}{2} \leq \\
\leq \Vert\,{\mathcal U}\,\Vert^2\,\frac{3\,n\,{\mathcal C}_5\,\rho\ 
\delta^2_q(\hat\tau_q)}{2\,r^2_q(\overline{\tau_q})}\leq\Vert\,{\mathcal U}\,\Vert^2\,\frac{6\,n\,{\mathcal C}_5\,\rho\ 
\delta^2_q(\hat\tau_q)}{r^2_q(\hat\tau_q)} 
\leq \Vert\,{\mathcal U}\,\Vert^2\,\frac{6\,n\,{\mathcal C}_5\,\rho\ }{t_0^2}\ .
\end{multline*}
Setting ${\mathcal C}_6=6\,n\,{\mathcal C}_5\,$, we have established, at any point 
$q\in\hat\tau\subset B(x,\rho)$ with $x\in W\,$, the estimate $(ii)$, true for any 
$\hat\tau$ 
of thickness $t_g(\hat\tau)\geq t_0\,$
\begin{equation} \label{ast5}
\sup_{q\in\hat\tau}\Vert\,h_q\,\Vert\leq \frac{{\mathcal C}_6}{t_0^2}\,\rho\ .
\end{equation}

We now prove $(i)$ ($(i')$ follows directly). Actually, it can be seen as a 
consequence of $(ii)$. For all $i,j=0,1,\dots,n\,$, introduce the function
\begin{equation*}
q\in\hat\tau\longmapsto \hat g(A_q^{-1} w_{ij},A_q^{-1} w_{ij})-g(A_q^{-1} 
w_{ij},A_q^{-1} 
w_{ij})\, .
\end{equation*}
This function vanishes along the $g$-geodesic $c$ joining 
$c(0)=p_i$ to $c(1)=p_j$ : indeed, we know 
from the barycentric construction (proposition \ref{edges}) that $c$ is also the
$\hat g$-geodesic running from $c(0)=p_i$ to $c(1)=p_j\,$, so that 
$\hat g(A_q^{-1} w_{ij},A_q^{-1} w_{ij})=
g(A_q^{-1} w_{ij},A_q^{-1} w_{ij})=d^2(p_i,p_j)$.

So, since $A_q^{-1} w_{ij}$ is $\hat g_{\hat\tau}$-parallel (in 
$(\hat\tau,\hat g)$) 
along a $g$-geodesic $q(\cdot)$ that goes from 
$q(0)=p_i$ to $q(L)=q$, one gets
\begin{multline*} 
(\hat g-g) (A_{q}^{-1} w_{ij},A_{q}^{-1} w_{ij}) 
= \int_0^L\frac{\partial}{\partial t}[(\hat g-g) (A_{q(\cdot)}^{-1} 
w_{ij},A_{q(\cdot)}^{-1} w_{ij})]_t\,dt = \\ 
=-2\,\int_0^L[g (\nabla_{\frac{\partial q}{\partial t}(\cdot)}(A_{q(\cdot)}^{-1} 
w_{ij}),A_{q(\cdot)}^{-1} w_{ij})]_t\,dt = \\ 
=2\,\int_0^L[g (h_{q(\cdot)}(\frac{\partial q}{\partial t},A_{q(\cdot)}^{-1} 
w_{ij}),A_{q(\cdot)}^{-1} w_{ij})]_t\,dt\  .
\end{multline*}
Using $\Vert h_q(U_1,U_2)\Vert_q\!\leq\!
H\,\Vert U_1 \Vert_q \,\Vert U_2\Vert_q\,$, 
with $H\!=\!\sup_{q\in\hat\tau}\Vert h_q \Vert\,$, thanks to $(ii)$, (\ref{estime1}),
(\ref{ast5}),
one may write at any $q\in\hat\tau$ (see definition \ref{bartauq})
\begin{multline} \label{ast2}\!
\vert(\hat g-g) (A_{q}^{-1} w_{ij},A_{q}^{-1} w_{ij})\vert\! \leq \!\!
2H\,d(p_i,q)\Vert\,A_q^{-1} 
w_{ij}\,\Vert_q^2\!\leq\!\frac{2\rho^2{\mathcal C}_6}{t_0^2}\,
\delta^2_q(\overline{\tau}_q)\,.
\end{multline}
Taking $\,{\mathcal U}\in T_q\hat\tau$ and writing 
${\mathcal U}=\sum_{i=1}^nu_i\,{\mathcal A}_i\,$, we have
\begin{equation*}
\vert\,(\hat g-g) ({\mathcal U},{\mathcal U})\,\vert\leq 
(\sum_{i,j=1}^n\vert u_i\vert\,\vert u_j\vert)\ \vert\,(g-\hat g) 
({\mathcal A}_i,{\mathcal A}_j)\,\vert\  ,
\end{equation*}
polarizing gives with (\ref{ast2})
\begin{equation*}
\vert\,(\hat g-g) ({\mathcal U},{\mathcal U})\,\vert\leq (\sum_{i=1}^n\vert 
u_i\vert)^2 
\ \frac{3\,\rho^2\,{\mathcal C}_6\,\delta^2_q(\overline{\tau}_q)}{t_0^2}\  .
\end{equation*}
Again using Schwarz's inequality, lemma \ref{1b} and  lemma 
\ref{tricky}, we 
get
\begin{multline*}
\vert\,(g-\hat g) ({\mathcal U},{\mathcal U})\,\vert \leq 
\rho^2\,\frac{3\,n\,{\mathcal C}_6}{t_0^2}\,
\frac{\delta^2_q(\overline{\tau}_q)}{r_q^2(\overline{\tau}_q)}\,
\Vert\,{ \mathcal U}\,\Vert_q^2\leq \rho^2\,\frac{12\,n\,
{\mathcal C}_6}{t_0^4}\,\Vert\,{\mathcal U}\,\Vert_q^2\  .
\end{multline*}
Thus $(i)$ and $(ii)$ hold with (remember $t_0<1$)
$${\mathcal C}_7=\max(\frac{{\mathcal C}_6}{t_0^2},\frac{12\,n\,
{\mathcal C}_6}{t_0^4})=\frac{2^23\,n\,{\mathcal C}_6}{t_0^4}=\frac{2^33^2n^2\,{\mathcal C}_5}{t_0^4}
$$ 
and, in view of 
${\mathcal C}_5= 2^3\,({\mathcal M}_2+3{\mathcal M}_1+2^3\,3\sqrt{3n}\,{\mathfrak R}_0/t_0)$ (see proposition \ref{lemma 2}), one gets 
$${\mathcal C}_7= 2^63^2n^2t_0^{-4}({\mathcal M}_2+3{\mathcal M}_1+2^33\sqrt{3n}\,{\mathfrak R}_0/t_0)\ .
$$

Now $(i)$, $(ii)$ and proposition \ref{proposition 1st} (set $C_1={\mathcal C}_7\rho^2$) imply $(iii)$.
\end{proof}

\section{A pattern fit to build a 
Euclidean 
approximation}

For all the material on complexes, triangulations and differential topology, we refer to 
\cite{Mu}, \cite{G-G}, \cite{Mi} and \cite{Gre}.

\subsection{Controlled refinements of embedded simplicial complexes}\label{6.1}$ $

\begin{convent}\label{convention2} Let $W$ be a relatively compact open 
set in $(M,g)\,$. 
Part of the general logical scheme which is followed in this section is to rely on successive bounds ${\mathcal R}_l>0$ for $l= 0,\cdots,5$, those for $l=0,1,2$ referring to convention {\rm\ref{convention1}}, for $l=3\,$ to section {\rm\ref{S.close}}, all verifying 
$$0<{\mathcal R}_5\leq{\mathcal R}_4\leq{\mathcal R}_3\leq {\mathcal R}_2\leq {\mathcal R}_1\leq {\mathcal R}_0\leq {\mathcal R}_W\ .
$$ 
Those ${\mathcal R}_l\,$, if $\rho$ belongs to $]0, {\mathcal R}_l]\,$, guarantee a series of $l$--statements (indexed by $l$) - depending on an integer $E=E(\rho)$ (increasing as $\rho$ decreases) - to hold true. {\em A choice of $\rho>0$ small enough} (of $E=E(\rho)$ large enough) ensures all $l$-statements are true for $l=0,\cdots,5\,$.
\end{convent}

Recall that $W$ is a chosen relatively compact open set in $(M,g)$ and ${\mathcal R}_W$ 
is 
defined according to lemma \ref{L.constante}. 

In this section \ref{6.1}, theorem \ref{2n} shows that {\it subdividing} enough {\it in a ``regular way'' a finite simplicial complex $K$ embedded through a simplicial mapping $T$ onto $T(K)\subset W\subset(M,g)\,$} and {\it using $g$-barycentric coordinates} enables to replace each simplex $\tau_E$ of the subdivided $T(K_E)$ by a Riemannian barycentric simplex $\hat\tau_E$ having the same vertices $p_0,\cdots,p_n$ which carries a well defined Euclidean metric $\hat g$ that verifies $d_g(p_i,p_j)=d_{\hat g}(p_i,p_j)$ for any $i,j\in\{0,1,\cdots,n\}\,$. Moreover, {\it the central theorem {\rm\ref{fundametal}} of this paper is true for any $(\hat\tau_E,\hat g)$}. The ``order'' of the subdivision process depends on a pair $(\rho,E)$ expressing dual faces of the same thing: the thinness of the subdivided simplicial complex is governed by a real $\rho>0\,$, the {\em mesh}, while the corresponding increase (order) of the simplices is ruled by the coupled integer $E=E(\rho)$, see definition \ref{linkint}. {\it A bound ${\mathcal R}_5\leq {\mathcal R}_3\leq \cdots\leq {\mathcal R}_W$ (depending on $T,K, W\subset(M,g)$) ensures that, in section {\rm \ref{6.1}}, all results on $K_E$ hold if $\rho\in]0,{\mathcal R}_5/4]\,$, as is announced in convention {\rm\ref{convention2}}.}

\par
Denote by $K$ a {\it finite} $n$-dimensional simplicial complex equipped with the 
standard 
Euclidean metric ${\rm std}$, induced from the standard Euclidean structure in some 
$\R^{N+1}$ (with $N$ large enough) by viewing $K$ as a sub-complex of the standard 
regular simplex $\Sigma_N\subset \R^{N+1}$ whose vertices are $e_0,e_1,\dots,e_N\,$, the vectors of the canonical 
basis in $\R^{N+1}\,$: {\em each simplex in $K$ is standard regular} (replace the ``standard simplices'' in the proof of lemma 9.4, page 92, \cite{Mu}).
More, we restrict our 
concern to a simplicial complex having underlying topological space $\vert K\vert$ 
homeomorphic to 
a connected $n$-dimensional manifold, whose boundary may be not empty.

\begin{defn}\label{D.Simplicial.Embedding}
 A {\em differentiable simplicial mapping} $T:K\rightarrow 
(M,g)$ is a continuous map from 
$\vert K\vert$ to its range $T(\vert K\vert)$ which also is 
a differentiable 
diffeomorphism from each $n$-simplex 
$\sigma$ of $K$ 
to its image $T(\sigma)$, in the sense that $T$ can 
separately be extended as a 
diffeomorphism from a neighborhood of $\sigma$ (in the affine $n$-space supported by each 
$\sigma$) to a neighborhood of $T(\sigma)$ in $M\,$. The map $T$ is {\it a differentiable simplicial embedding}
if it also is a homeomorphism from 
$\vert K\vert$ to $T(\vert K\vert)$.
\end{defn}
\begin{lem} \label{2a} \label{2a-bis}$ $
\begin{enumerate}
\item Given a Riemannian manifold $(M,g)$ and a differentiable simplicial embedding
$T:(K,{\rm std})\rightarrow (M,g)\,$, there exist two positive best constants 
$\alpha_T$ and $\beta_T$ such that, for any $n$-simplex $\sigma$ in $K$
\begin{equation*}
\forall y\in\sigma\,,\ \forall u \in T_y\sigma\setminus\{0\}\,,\ 
\alpha_T\leq \frac{\Vert\,dT(y)(u)\,\Vert_g}{\Vert\,u\,\Vert_{{\rm std}}}\leq\beta_T\   .
\end{equation*}
\item \label{2a-bis33} Denoting by $K_E$ a {\rm regular subdivision of integer} 
$E$ of $K$ (the meaning of a {\rm regular subdivision of integer $E$ is explained in course of the 
proof below),} there exist positive reals ${\mathscr D}_n,{\mathcal D}_n,\Theta_n$ 
depending only on $n$,  such 
that, for any $E\in \N^*$ and 
any $n$-simplex $\sigma_E$ in $K_E\,$, one has
\begin{gather*}
(i) \ \, \frac{{\mathscr D}_n}{E}\leq \hbox{\rm diam}_{{\rm std}} (\sigma_E)\leq 
\frac{{\mathcal D}_n}{E}\ , \ \  
 (ii)\ \ \frac{\sqrt{n+1}}{2^{\frac{n}{2}}\, n!}\geq\omega_{{\rm std}}(\sigma_E)\geq \Theta_n \  ,\\   \hbox{and} \ \ \ \  (ii') \ \ 
 t_{{\rm std}}(\sigma_E)\geq {\bf t}_n:=\frac{n}{n+1}\, 
 \frac{1}{\lvert \sigma_{n-1}\vert}\,\Theta_n\ .
\end{gather*}
\end{enumerate}
\end{lem}
\begin{proof}
 Denote by $K_n\,$ the finite set of all $n$-simplices $\sigma\in K$. The first statement holds making the following definitions of $\alpha_T$ and $\beta_T$
 \begin{gather*}
 \alpha_T = \inf
\left\{ \frac{\Vert\,dT(y)(u)\,\Vert_g} {\Vert\,u\,\Vert_{{\rm std}} }\ \mid
\sigma \in K_n\ ,\,  y\in \sigma\ ,\ u\in T_y\sigma\setminus \{0\}  \right\}\ , \\
 \beta_T = \sup
\left\{ \frac{\Vert\,dT(y)(u)\,\Vert_g} {\Vert\,u\,\Vert_{{\rm std}} }\ \mid
 \sigma \in K_n\,,\  y\in \sigma\,,\  y\in \sigma\,,\ \ u\in T_y\sigma\setminus \{0\} \right\} \  .
\end{gather*}

It is enough to prove $(i),(ii)$ in a standard regular $n$-simplex of $K\,$ since $K$ is
a collection of such simplices. The reader may consult \cite{Mu} lemma 9.4, we adapt the proof given there for the {\it 
standard simplex $\subset{\mathbb R}^n$} (having vertices $0,e_1,\dots,e_n$) to the standard regular 
$n$-simplex $\Sigma_n\subset\{x_0+x_1+\dots+x_n=1\}$ ($\Sigma_n\subset{\mathbb R}^{n+1}$ has vertices $e_0,e_1,\dots,e_n$). {\em We now describe a regular simplicial subdivision of $K$ of integer} $E\,$.

First, the hyperplanes $\{x_0+x_1+\dots+x_n=l\}\,$, with $l\in \Z\, $, which intersect non 
trivially the unit cube $[0,1]^{n+1}$ correspond to the values $l=0,1,\dots,n+1$, so 
this 
kind of intersections arises 
in a finite number of polyhedral shapes. 

Another observation is that cutting 
$\Sigma_n$ 
by all  
cubes $(k_0/E,\dots,k_n/E)+[0,(1/E)]^{n+1}$ with $k_i\in \Z $, produces, up to a 
homothety of dilation ratio $E$, the same polyhedral shapes as cutting some unit cube 
$(k_0,\dots,k_n)+[0,1]^{n+1}$ with 
$k_i\in \Z $, by a hyperplane $\{x_0+x_1+\dots+x_n=m\}$ with $m$ integer, which is 
itself, 
up 
to translations, the thing observed in the first place. 

So, whatever $E$ is, dividing in this way a regular simplex produces only a 
finite number of 
polyhedral shapes. A simple 
barycentric subdivision (as in the proof of \cite[Lemma 7.8]{Mu}) then produces the {\em regular simplicial decomposition of integer} $E\,$, whose simplices also come out in a finite number of geometric shapes. Therefore $(i)$ and $(ii)$ follow.
\par As for $(ii')$, let $\sigma_E$ be an $n$-simplex in $(K_E,{{\rm std}})$. Perform a homothety $H_\lambda$ of dilation ratio $\lambda$ on 
$\sigma_E$ in order to get $\diam (H_\lambda(\sigma_E))=1\,$. 
Thus 
$\omega_{{\rm std}} (\sigma_E)=\omega_{{\rm std}} (H_\lambda(\sigma_E))=\vol_{{\rm std}} (H_\lambda(\sigma_E))$. The 
isoperimetric inequality given in \cite[inequality (4)]{Rez} reads (see remark \ref{R.regular})
$\vol_{{\rm std}} (H_\lambda(\sigma_E)) \leq \vert \sigma_n\vert=\omega(\sigma_n)$. Thus, one has 
$\ \Theta_n \leq \omega_{{\rm std}}(\sigma_E)\leq\omega(\sigma_n) =
\vert \sigma_n\vert={\sqrt{n+1}}/({2^{\frac{n}{2}}\, n!})\,$,
which implies $(ii')$ by lemma \ref{1c}.
\end{proof}
\begin{defn} \label{2j}  Keep the notations of lemma \ref{2a}.
\begin{enumerate}
 
 \item \label{2j-1}Define two uniform positive coupled constants 
(depending only upon $\alpha_T,\beta_T$ and $n$) 
\begin{equation} \label{2c}
\overline{A}:=\bigg{(}\frac{\alpha_T}{\beta_T}\,
\frac{{{\mathscr D}}_n}{{\mathcal D}_n}\bigg{)}^{\!\!n}\,\Theta_n\ \ \ \hbox{and}
\ \ \  \overline{t}:=\frac{n}{n+1}\ \frac{\overline{A}}{2^n\vert\sigma_{n-1}\vert}\ ,
\end{equation}
where $\vert\,\sigma_{n-1}\,\vert$ is the $(n-1)$-volume of the regular 
$(n-1)$-simplex having edges of  length $1\,$.
 In connection with the above $\bar t>0\,$, introduce two auxiliary positive reals 
\begin{equation*}
t_1:=\frac{\overline{t}}{12^n}\  ,\  \   t_0:=\frac{\overline{t}}{4\,(6^n)}= \frac{n}{4(n+1)}\ 
\frac{\overline A}{12^n\,\vert\sigma_{n-1}\vert }\ .
\end{equation*}

Since 
${\overline A} \leq \Theta_n \leq \vert 
\sigma_n\vert =\frac{\sqrt{n+1}}{2^{\frac{n}{2}}\, n!}\,$, 
one has ($t_1<t_0$ if $n\geq3$)
\begin{equation}\label{curiosity} 0<t_1\leq t_0 < \frac{1}{\sqrt{2n(n+1)}}<1/2\ .
\end{equation}

\item  \label{flop}As $t_1 \in \, ]0,\frac{1}{\sqrt{2n(n+1)}}[\,$ (a direct consequence of (\ref{curiosity})),   
Euclidean $n$-simplices with thickness $\geq t_1$ exist (see lemma \ref{1c}). 
Let $C>0\,$, denoted here by $\check C\,$, be such that (recall $n\geq2$)
$$\check C\leq  
\frac{t_1^2}{45n}
\ , \ \ \ \hbox{thus also} \ \ \ \check C<\frac{2\,t_1^2}{3^3n}<\frac{t_1^2}{9n}<\frac{1}{12} \ .
$$
Denote by ${\mathcal R}_0 (\check C)$ the bound related to the above $C\!=\!\check C$ through proposition \ref{proposition A} and define a new positive real bound ${\mathcal R}_4$
\begin{equation*}
{\mathcal R}_4 =  \min({\mathcal R}_0 (\check C), {\mathcal R}_3, \frac{1}{\sqrt{3\,{\mathcal C}_7}})\,,
\end{equation*}
here ${\mathcal R}_3,{\mathcal C}_7$ 
are linked to theorem \ref{theorem 3}.
From ${\mathcal R}_4\leq{\mathcal R}_3$, one still has
$
{\mathcal R}_4 \leq \min(\frac{1}{\sqrt{2\,{\mathcal C}_2}}, \frac{t_0}{5\sqrt{n\,{\mathcal M}_1}})\,$,
where ${\mathcal C}_2$ is related to lemma \ref{lemma 2nd} and 
${\mathcal M}_1$ to theorem \ref{theorem 1st}. 
\par\noindent {\it If $\rho \leq {\mathcal R}_4\leq  {\mathcal R}_3\,$, one has ${\mathcal C}_7\,\rho^2\leq1/3$
and theorem {\rm\ref{theorem 3}} applies.}

While $C=\check C\leq \frac{t_1^2}{45n} < \frac{t_1^2}{9n}\,$ is chosen as above, and as (by convention \ref{convention2}) one has ${\mathcal R}_4  \leq {\mathcal R}_W\,$,
{\it we may apply propositions {\rm \ref{proposition A}}, {\rm \ref{proposition B}} and 
{\rm\ref{proposition C}} 
on each ball $B(x,{\mathcal R}_4)$ for any $x\in W\,$.}
\end{enumerate}
\end{defn}
\begin{rem}\label{convex333} As ${\mathcal R}_4 \leq {\mathcal R}_0(\check C) \leq {\mathcal R}_W\, $,  
so is actually ${\mathcal R}_4$ 
smaller than {\it half} of the convexity radius at all 
points in $\overline{W}_{2{\mathcal R}_4}$, see lemma \ref{L.constante}.
In fact, as announced
$${\mathcal R}_4\leq\min ({\mathcal R}_3,{\mathcal R}_0(\check C))\leq{\mathcal R}_2\leq{\mathcal R}_1\leq {\mathcal R}_0(1/2)\ .
$$
\end{rem}

\begin{defn}\label{m} 
An embedded $n$-simplex $\tau$, in a convex ball of a 
Riemannian $n$-dimensional manifold $(M,g)$, is {\em Euclideable} if 
\begin{enumerate}
 \item  its vertices build up a $g$-Riemannian barycentric simplex 
$\hat{\tau}$, i.e. 
the vertices of $\tau$ determine unique barycentric coordinates defining 
$\hat{\tau}$;

\item by means of the barycentric coordinates, one can determine a unique 
Euclidean 
structure $\hat g$ on $\hat{\tau}$ for which the vertices are at the same 
$\hat g$ and $g$-distances.
\end{enumerate}
\end{defn}
\begin{defn} \label{linkint}
Given $\rho\in]0,\beta_T\,{\mathcal D}_n]\,$, denote by  $E(\rho)$ the (positive) integer part of the real $1+(\beta_T\,{\mathcal D}_n/\rho)$. Thus one has the following coupling 
between $\rho\,$ and $E(\rho)$ 
\begin{equation}\label{2l}
\frac{\beta_T\,{\mathcal D}_n}{\rho}\, < \,E\,\leq\, 1+\frac{\beta_T\,
{\mathcal D}_n}{\rho}\,\leq\,\frac{2\,\beta_T\,{\mathcal D}_n}{\rho} \  .
\end{equation}
\end{defn}

The main result in this section.
\begin{thm}\label{2n} Let $T\!:\!(K,{\rm std})\!\rightarrow W\!\subset\!(M,g)$ be a differentiable simplicial embedding. There exists ${\mathcal R}_5\!>\!0$ with ${\mathcal R}_5\!\leq\!\min(4\,\beta_T\,{\mathcal D}_n\,,\,{\mathcal R}_4)\!\leq\! {\mathcal R}_W$ such that, for any $\rho \!\in ]0, {\mathcal R}_5/4]$ and $E\!=\! E(\rho)$, there is a regular subdivision $K_E$ of $K$ such that any 
$n$-simplex $\sigma_E\!\in\!K_E$ gives a {\rm Euclideable}
$g$-Riemannian barycentric simplex $\hat\tau_E$; more, theorem 
{\rm\ref{theorem 3}} applies to $(\hat\tau_E,\hat g)$. If $\varpi$ is the $g$-gravity center 
of $\!\hat\tau_E\,$, if $(\hat\tau_E)_q\subset T_qM$ is the linear simplex of vertices those of $\exp_q^{-1} \hat\tau_E\,$, then (see 
definition \rm{\ref{genthick}} for $(ii)$)
\begin{align*}
(i) &   \hskip4mm & \hat\tau_E\subset B(\varpi,\rho)  & \hskip3mm \hbox{and} \hskip4mm \hbox{\rm diam}_g(\hat\tau_E)\leq 2\,\rho\  ; & \\
(ii) &  \hskip2mm & t_g(\hat\tau_E)\geq t_0 \  \  &\hbox{meaning}\ \ \ \forall q\in \hat\tau_E \ \ \ t_{g_q}((\hat\tau_E)_q)\geq t_0 \ ;&   \\
(iii) &  \hskip6mm & \hbox{\rm vol}_g(\hat\tau_E)\geq A_0\,\rho^n & \ \ \  \hbox{where} \ \ \ \    A_0=\frac{\bar A}{12^n}\ .&\\
\end{align*}
{\it One can cover $W$ by finitely many convex balls $B(x_i,{\mathcal R}_4)$ and define,
for any $i$, a sub-complex $K_{E,i}$ of $K_E$ such that $T(K_{E,i})\!\subset\! B(x_i,{\mathcal R}_4)$ and $K_E\!=\!\cup_i K_{E,i}$.}
\end{thm}
The results to follow build a proof of the above theorem.
\begin{prop}\label{2o}
Let $T:(K,{\rm std})\rightarrow W\subset(M,g)$ be 
a differentiable simplicial embedding. 
For any  
$\rho\in\,]0,\beta_T\,{\mathcal D}_n]$ verifying  $\rho\leq {\mathcal R}_4\,$, let
$E = E(\rho) \in \N$. 
Given any $n$-simplex $\sigma_E$ in $K_E$, the regular simplicial subdivision of integer $E$, 
denoting by $\tau_E$ the image $T(\sigma_E)$, one has the following control on 
the embedding given by $T$ 
\begin{align*}
(i) &  \hskip4mm&\hbox{\rm diam}_g(\tau_E) < \rho\  ;& \hskip3cm & \hskip9mm \\
(ii) &  \hskip4mm& \omega_g(\tau_E)\geq \frac{\overline{A}}{2^n}\ ; & & \\
(iii) &  \hskip4mm & \hbox{\rm vol}_g(\tau_E)\geq \frac{\bar A}{2^n}\,\rho^n\ .&\\
\end{align*}
\end{prop}
\begin{rem}\label{markit} Thus $\rho\leq{\mathcal R}_4$ implies $\tau_E \!\subset \!B(x,{\mathcal R}_4)$ for some $x\in W\,$.
\end{rem}

\begin{proof}
\hskip1mm The extrinsic diameter $\hbox{\rm diam}_g(\tau_E)$ of $\tau_E$ in 
$(M,g)$ is smaller than its intrinsic diameter 
as metric subspace. Using lemma \ref{2a} part 1. get
\begin{equation*}
\hbox{\rm diam}_g(\tau_E)\leq \beta_T\ \hbox{\rm diam}_{{\rm std}}(\sigma_E)\,,
\end{equation*}
and by the change of variables in a volume integral 
(definition of $\alpha_T$)
\begin{equation*}
\hbox{\rm vol}_g(\tau_E)\geq\alpha_T^n\  \hbox{\rm vol}_{{\rm std}}(\sigma_E)\,.
\end{equation*}
Using lemma \ref{2a} part 2. and  (\ref{2l}) of definition \ref{linkint}, we get
\begin{equation*}
\hbox{\rm diam}_{{\rm std}}(\sigma_E)\leq \frac{{\mathcal D}_n}{E}
< \frac{\rho}{\beta_T} \  ,
\end{equation*}
thus $\hbox{\rm diam}_g(\tau_E) < \beta_T\,\frac{\rho}{\beta_T}=\rho\ $, 
proving $(i)$.
\par
Finally, one has (using (\ref{2l}) of definition \ref{linkint}, as $\rho\leq \beta_T\,{\mathcal D}_n$)
\begin{multline*}
\hskip4mm\hbox{\rm vol}_g(\tau_E)\geq\alpha_T^n\  
\hbox{\rm vol}_{{\rm std}}(\sigma_E)\geq
\alpha_T^n\  \Theta_n\,(\hbox{\rm diam}_{{\rm std}}(\sigma_E))^n\geq \hfill \cr
\hfill\geq\alpha_T^n\  \Theta_n\,(\frac{{\mathscr D}_n}{E})^n\geq
(\frac{\alpha_T}{\beta_T}\,\frac{{\mathscr D}_n}{{\mathcal D}_n})^n\,
\frac{\Theta_n}{2^n}\,\rho^n\geq \frac{\overline{A}}{2^n}\,\rho^n\geq 
\frac{\overline{A}}{2^n}\, (\hbox{\rm diam}_g(\tau_E))^n\  ,
\end{multline*}
and the last two inequalities produce $(iii)$ and $(ii)$.
\end{proof}

\begin{defn} \label{secantmap}
Given $x$ in $W$, 
for any $q\in B(x,{\mathcal R}_4)$ {\it define $T_q$} to be (remark \ref{convex333}) 
$$T_q:=\exp_q^{-1}\circ\  T:T^{-1}B(q,{2\mathcal R}_4)\longrightarrow B(0_q,{2\mathcal R}_4)\subset T_qM\ .
$$
{\it Denote} by $L_x$ the affine map, called the {\it secant map} (of $T$ on $K_E$), see \cite{Mu} page 90, sending 
any $\sigma_E\in K_E$ with $\tau_E=T(\sigma_E)\subset B(x,{\mathcal R}_4)$ to the linear simplex $L_x(\sigma_E)\subset B(0_x,{\mathcal R}_4)\subset T_xM$ of same vertices as $T_x(\sigma_E)$. 
\end{defn}

\begin{prop}\label{2p} 
Let $T:(K,{\rm std})\rightarrow W\subset(M,g)$ be 
a differentiable simplicial embedding 
and let  $\rho \in \,]0, \beta_T\,{\mathcal D}_n]$,
also chosen such that $\rho\leq {\mathcal R}_4\,$. 
With 
$E = E(\rho)$ (see definition {\rm\ref{linkint}}), 
given an $n$-simplex $\sigma_E$ in $K_E$,  
and $x\in W$ such that $\tau_E= T(\sigma_E) \subset B(x,{\mathcal R}_4)$ (see 
remark {\rm \ref{markit}}), 
one 
has, for any $q$ in $B(x,{\mathcal R}_4)$
\begin{gather*}
(i)\hskip4mm \hbox{\rm diam}_{g_q}(T_q(\sigma_E))\leq 2\,\rho\, ;
\hskip1mm 
(ii) \hskip4mm \omega_{g_q}(T_q(\sigma_E))\geq 
\frac{\overline{A}}{6^n}\,;\\\hskip1mm 
(iii) \hskip4mm \hbox{\rm vol}_{g_q}(T_q(\sigma_E))\geq 
\frac{\overline{A}}{3^n}\,\rho^n\,.
\end{gather*}
\end{prop}

\begin {proof}
By remark \ref{convex333}, proposition \ref{proposition A} applies at all 
points $q\in \overline{W}_{2{\mathcal R}_4}\,$.
From proposition \ref{proposition A} and definition \ref{2j} \ref{flop}, as we chose $C<1/2\,$ and ${\mathcal R}_4\leq {\mathcal R}_0(C)$, one 
has for any $q$ in $B(x,{\mathcal R}_4)$ and $u$ in $B(0_q,2{\mathcal R}_4)$ 
\begin{equation}\label{2pa}
 \forall v \in T_uT_qM\,, \quad \frac{1}{2}\,\Vert \,v\,\Vert_q \leq 
\Vert \,d\exp_q(u)(v)\,\Vert_g\leq \frac{3}{2}\,\Vert \,v\,\Vert_q\  .
\end{equation}
Thus,
defining the distance with respect to $g$ or to $g_q$ by 
minimizing 
the length of arc is here straightforward because the balls 
$B(q,2{\mathcal R}_4)$ or $B(0_q,2{\mathcal R}_4)$ are 
both convex in the respective metrics.  The conclusions are then directly derived from 
proposition \ref{2o}, using the infinitesimal control given in (\ref{2pa}).
\end{proof}

\begin{rem} \label{secantmap1033}
In our development, in order to handle the soon occurring secant maps defined by simplicial mappings, a notion of $C^1$-approximation for differentiable simplicial maps is needed.
\end{rem}
Close to Munkres' text \cite{Mu} is the
\begin{defn}\label{Argh55}
Let $f : (K,\hbox{std}) \rightarrow ({\mathbb R}^n,\Vert\,\cdot\,\Vert)$ 
be a differentiable simplicial mapping. 
Given $\eta>0$ and a subdivision $K'$ of $K\,$, a map $g:\vert K\vert \rightarrow \R^n$ is a {\em strong} {\em $\eta$-approximation} to $f$ on $K'$  if $g$ is a differentiable simplicial mapping $g:
 K' \rightarrow \R^n$ and  
for all $\sigma'\in K',y\in\sigma',u\in 
T_y\sigma'$
\begin{gather*}
(a) \hskip4mm \Vert\,g(y)-f(y)\,\Vert\leq \eta\  ,\\
 (b) \hskip4mm\Vert\,dg(y)(u)-df(y)(u)\,\Vert\leq \eta\,\Vert\,u\,\Vert_{{\rm std}}\ .
\end{gather*}
Skipping $K'\,$, one also says losely {\em $g$ is a strong $\eta$-approximation to $f\,$}.
\end{defn}

\begin{rem} 
In our context, this definition is equivalent to the one given by 
Munkres, 
(see \cite[definition 8.5]{Mu} and the lines that follow it).
\end{rem}

\begin{defn}\label{constant99} Denote by $K_{E,x}$ the subcomplex of those $\sigma_E\in K_E$ contained in $B(x,{\mathcal R}_4)$ with all their faces. Define
\begin{equation*}{\mathscr T}_x(E) :=\sup_{\sigma_E\in K_{E,x}} \,\sup_{y,z\in\bar\sigma_E}\ \Vert dT_x(z)-dT_x(y)\Vert_{{\rm std},x}\ .
\end{equation*}
Set also
$${\bf T}_x(E):=\frac{2\,\sqrt{n}\,{\mathscr T}_x(E)}{{\bf t}_n} \ \ \ \hbox{and}\ \ \  {\bf T}(E):=\sup_{x\in W} \ {\bf T}_x(E)\ .
$$
\end{defn}

\begin{rem} If $\rho$ is small enough ($\rho\leq {\mathcal R}_4$), thus $E=E(\rho)$ sufficiently large, the bound ${\mathscr T}_x(E)$ exists for any $x\in \bar W$. It tends to $0$ as $E$ tends to $\infty\,$, because $\hbox{diam}_{\rm std}(\sigma_E)\leq \frac{{\mathcal D}_n}{E}<\frac{\rho}{\beta_T}$ (by definition \ref{linkint}). More, the bound ${\bf T}(E)$ is finite and tends to $0$ with $E$ tending to $\infty$ since the $T_x$ constitute a continuous piecewise differentiable family of mappings continuously indexed by $x$ varying in the compact $\bar W\,$. This allows to state the following lemma, very close to Lemma 9.3 in Munkres' book, see \cite{Mu}, which here is {\it uniform} as $x$ varies in $W\,$.
\end{rem}

\begin{lem}\label{2q} Let $T:(K,{\rm std})\rightarrow W\subset(M,g)$ be 
a differentiable simplicial embedding. Define $\rho(\eta)\in]0,{\mathcal R}_4]$ to be such that ${\bf T}(E)\leq \eta\,$ for $E\geq E(\rho(\eta))$.
For any $x\in W\,$, $\rho\in]0,\rho(\eta)]$, with $E=E(\rho),$ and any
$n$-simplex $\sigma_E\in K_E$, one knows the {\rm secant map} $L_x$ to be a {\em strong} 
$\eta$-{\em approximation to} $T_x$, meaning,
for all $y$ in 
$\sigma_E$ and 
$u$ in $T_y\sigma_E$ 
\begin{gather*}
(a)\hskip4mm \Vert\,T_x(y)-L_x(y)\,\Vert_x\leq \eta\  ,\\
 (b) \hskip4mm\Vert\,dT_x(y)(u)-dL_x(y)(u)\,\Vert_x\leq 
 \eta\,\Vert\,u\,\Vert_{{\rm std}}\ . 
\end{gather*} 
\end{lem} 
\begin{rem} Observe that the smaller $\eta\,$, thus the smaller $\rho$, the larger $E=E(\rho)$ must be (definition \ref{linkint}).
\end{rem}
\begin{proof} See appendix \ref{preuveMunkres}.
\end{proof}

\begin{prop}\label{2r}
Let $T:(K,{\rm std})\rightarrow W\subset(M,g)$ be 
a differentiable simplicial embedding. Choose some $x\in W\,$.
Let  $\rho \in \,]0, \beta_T\,{\mathcal D}_n]$ be small enough so that (in accordance with previous lemma {\rm\ref{2q}})
$\rho \leq \min({\mathcal R}_4,\rho(\alpha_T/4))$ and
$E=E(\rho)$.
For any $n$-simplex in $\sigma_E$ in $K_E\,$ such that $T(\sigma_E)\subset B(x,{\mathcal R}_4)$,
one has
\begin{align*} 
 & (i) & \hbox{\rm diam}_{g_x}(L_x(\sigma_E)) & \leq & 
\hbox{\rm diam}_{g_x}(T_x(\sigma_E)) & \leq & \,2\,\rho\ ; & \\
 & (ii)  &\omega_{g_x}(L_x(\sigma_E)) & \geq & \frac{1}{2^n}\ 
\omega_{g_x}(T_x(\sigma_E)) 
& \geq &\,\frac{\overline{A}}{12^n}\ ;&  \\
& (iii) &  \hbox{\rm vol}_{g_x}(L_x(\sigma_E)) & \geq 
& \frac{1}{2^n}\,\hbox{\rm vol}_{g_x}(T_x(\sigma_E)) & \geq 
& \frac{\overline{A}}{6^n}\,\rho^n\ .& 
\end{align*}
\end{prop}

\begin{proof} 

The left inequality in $(i)$ expresses the Euclidean character of 
the metric $g_x$, while the right inequality is proposition \ref{2p} $(i)$.

Given $x\in W$, let $E$ be chosen as in the statement to prove and $\sigma_E$ in $K_E$ satisfy 
$T(\sigma_E) \subset B(x,{\mathcal R}_4)$. Let $y\in \sigma_E$ and 
$u\in T_y \sigma_E\,$, $u\not= 0\,$. 
Set  $v=dT_x(y)(u)$. 
Applying lemma \ref{2q} (b) with $\eta=\alpha_T/4$ (we assumed $\rho\leq \rho(\alpha_T/4)$), one gets (with $\Vert\cdot\Vert_y= \Vert\cdot\Vert_{\rm std}$)
\begin{equation} \label{2ra} 
\hskip4mm\Vert\,v-d L_x(y)\circ (dT_x(y))^{-1}(v)\,\Vert_x\leq  
\frac{\alpha_T}{4}\,\Vert\,(dT_x(y))^{-1}(v)\,\Vert_{\rm std} \  .
\end{equation}
But, one may write 
\begin{equation*}
\frac{\Vert\,dT_x(y)(u)\,\Vert_x}{\Vert\,u\,\Vert_{{\rm std}}} = 
\frac{\Vert\,d\exp_x^{-1}(T(y))(dT(y)(u))\,\Vert_x}{\Vert\,dT(y)(u)\,\Vert_{T(y)}} \, 
\frac{\Vert\,dT(y)(u)\,\Vert_{T(y)}}{\Vert\,u\,
\Vert_{{\rm std}}} \ ,
\end{equation*}
thus, bringing in the inequalities (see proposition 
\ref{proposition A} and lemma \ref{2a})
\begin{gather*}
(1-C)\,\leq \Vert \,d\exp_x(T_x(y))\,\Vert_{x,T(y)}\leq \,(1+C)\,,\\
\alpha_T\leq \Vert\,dT(y)\,\Vert_{y,T(y)}\leq \beta_T\, ,
\end{gather*}
one derives
\begin{gather}
\frac{\alpha_T}{1+C}\leq \Vert\,dT_x(y)\,\Vert_{y,x}\leq 
\frac{\beta_T}{1-C}\ , \label{Ineg1} \\ 
\frac{1-C}{\beta_T}\leq \Vert\,(dT_x(y))^{-1}\,\Vert_{x,y}\leq 
\frac{1+C}{\alpha_T}  . \label{Ineg2}
\end{gather}
Lines (\ref{Ineg2}) and  (\ref{2ra}) give
\begin{equation*}
1-\frac{(1+C)}{4}\leq\,\Vert\,d L_x(y)\circ 
(dT_x(y))^{-1}\,\Vert_{x,x}\leq  1+\frac{(1+C)}{4}\  .
\end{equation*}
As 
$1-\frac{\eta\,(1+C)}{\alpha_T}\geq 1/2\,$, we finally 
deduce

{\it $(\ast)$ the map $L_x \circ T_x^{-1}$, which sends $T_x(\sigma_E)$ to 
$L_x(\sigma_E)$, 
has infinitesimal contraction factor not less than $1/2\,$.}

The first inequality in $(iii)$ follows from $(\ast)$, the second 
from proposition \ref{2p} $(iii)$. The first 
inequalities in (i)  
and $(iii)$ give the first inequality in $(ii)$.  
Proposition \ref{2p} $(ii)$ implies the second in $(ii)$.
\end{proof}

\begin{prop}\label{2s}
Let $T:(K,{\rm std})\rightarrow W\subset(M,g)$ be 
a differentiable simplicial embedding. Choose some $x\in W$ and 
let  $\rho \in \,]0, \beta_T\,{\mathcal D}_n]$ be such that 
$\rho \leq \min({\mathcal R}_4,\rho(\alpha_T/4))$ and  $E=E(\rho)$ (so proposition {\rm\ref{2r}} applies). 
Then, 
any $n$-simplex $\sigma_E$ of $K_E$ with $\tau_E:= T(\sigma_E)\subset B(x,{\mathcal R}_4)$
determines a 
$g$-Riemannian barycentric simplex $\hat\tau_E$ which is Euclideable, one has (see 
definition \rm{\ref{genthick}} for $t_g(\hat\tau_E)$)
\begin{equation*}
t_g(\hat\tau_E)\geq  t_0 \hskip4mm \hbox{(see definition {\rm\ref{2j} (\ref{2j-1})} for the positive real number $t_0$)}.
\end{equation*}
\par The proof \em{even} shows (definitions \rm{\ref{genthick}} and \rm{\ref{secantmap}} give $(\hat\tau_E)_q=L_q(\sigma_E)$)
$$\forall q\in B(x,{\mathcal R}_4) \ \ \ \ t_{g_q}(L_q(\sigma_E))=t_{g_q}((\hat\tau_E)_q)\geq t_0\  .
$$
\end{prop}

\begin{proof}
From lemma \ref{1c} we know  
\begin{equation} \label{Eq.lineaire}
 t_{g_x}(L_x(\sigma_E))\geq \frac{n}{n+1}\, \frac{1}{\vert\sigma_{n-1}\vert}\,
  \omega_{g_x}(L_x(\sigma_E))\,.
\end{equation}
We saw in definition \ref{2j} (\ref{2j-1}) that $t_1=\overline{t}/12^n \leq t_0< \frac{1}{\sqrt{2n(n+1)}}\,$.
Apply propositions \ref{proposition B} and 
\ref{proposition C} 
with $t_1$ to be $t_1=\overline{t}/12^n$ (the value introduced in definition \ref{2j} (\ref{2j-1})).
Proposition \ref{2r} $(ii)$, 
proposition \ref{proposition C} and (\ref{Eq.lineaire}) give $t_{g_x}(L_x(\sigma_E))\geq  \frac{\overline{t}}{6^n}\,$, so that
\begin{equation}\label{deadoralive}
 \forall q\in B(x,{\mathcal R}_4) \hskip18mm t_{g_q}(L_q(\sigma_E))\geq  
t_0\  .
\end{equation}
So, the vertices of $\tau_E$ are spread, see definition \ref{D.spread.simplex}, and thus generate a Riemannian barycentric simplex $\hat\tau_E$ by corollary \ref{rbs10}.
The inequalities (\ref{deadoralive}) imply that $\hat\tau_E$ 
has a $g$-thickness (in the generalized sense) $\geq t_0>0$
(see definition  \ref{genthick}). 
Proposition \ref{proposition B} and the fact that $B(x,{\mathcal R}_4)$ is a 
convex ball guarantee $\hat\tau_E$ to be Euclideable.
\end{proof}

\begin{prop}\label{2t} 
Let $T:(K,{\rm std})\rightarrow W\subset(M,g)$ be an embedding of a finite 
$n$-simplicial complex $K\,$.  Choose $x\in W$ and let  $\rho \in \,]0, \beta_T\,{\mathcal D}_n]$ be such that 
$\rho \leq \min({\mathcal R}_4,\rho(\alpha_T/4))$ and  $E\!=\!E(\rho)$ (so proposition {\rm\ref{2r}} applies). 
Then, 
any $n$-simplex $\sigma_E$ of $K_E$ with $\tau_E\!=\! T(\sigma_E)$ satisfying $\tau_E\subset B(x,{\mathcal R}_4)\,$
determines a 
$\,g$-Riemannian barycentric simplex $\hat\tau_E$ which is Euclideable and, calling  $\varpi$ the $g$-gravity center 
 of $\hat\tau_E\,$, one has
\begin{align*}
(i) &   \hskip4mm & \hat\tau_E\subset B(\varpi,\rho)\,,  & \hskip3mm 
\hbox{\rm so that} 
\hskip4mm \hbox{\rm diam}_g(\hat\tau_E)\leq 2\,\rho\  ; & \\
(ii) &  \hskip4mm & t_g(\hat\tau_E)\geq t_0 \,, &\hskip3mm\hbox{\rm and even}\ \ \forall q\in B(x,{\mathcal R}_4) \ \ t_{g_q}((\hat\tau_E)_q)\geq t_0 \ ;&   \\
(iii) &  \hskip4mm & \hbox{\rm vol}_g(\hat\tau_E)\geq A_0\,\rho^n\ ,&\hskip3mm 
\hbox{\rm where} 
\hskip4mm A_0= \overline{A}/12^n\  .&  
\end{align*}
So, apart from its last statement, theorem {\rm\ref{2n}} is established, setting ${\mathcal R}_5:=\inf (4{\mathcal D}_n\,\beta_T,{\mathcal R}_4,4\rho(\alpha_T/4))$.
\end{prop}

\begin{proof}
The existence of a Riemannian barycentric simplex $\hat\tau_E$ which is Euclideable and satisfies $(ii)$ was already proven in proposition 
\ref{2s}.

 We know $\,\hbox{\rm diam}_g(\tau_E)< \rho\,$ 
(from proposition \ref{2o}), so that we have $\tau_E  \subset B(p_i,\rho)\,$ for 
any vertex $p_i$ of $\tau_E\,$. By convexity of the balls $B(p_i,
\rho)\,$ (due to $\rho\leq {\mathcal R}_4$) and by definition of $\hat\tau_E\,$, we also have 
$\hat\tau_E  \subset B(p_i,\rho)\,$ for any $p_i\,$, see proposition \ref{rbs1}. 
So, the distance from $p_i$ 
to any point inside $\hat\tau_E $ is smaller than $\rho$, thus 
the ball $B(\varpi,\rho)$ contains all vertices $p_i$ of $\hat\tau_E $ : as 
$B(\varpi,\rho)$ is convex, it contains also $\hat\tau_E $, by definition of 
such a simplex. This gives $(i)$.

It remains to prove $(iii)$.
If $\,p_0,p_1,\dots,p_n$ are the vertices of $\tau_E$, set 
$v_0=\exp_x^{-1}p_0,v_1=\exp_x^{-1}p_1,\dots,v_n=\exp_x^{-1}p_n\in T_xM$ (see section \ref{S.close}).
Bring back the map $\phi$ from definition \ref{D.euclidianisable}, related to the ``$g$-Riemannian barycentric 
coordinates'', 
that sends a point $u=\sum_0^n\lambda_i\,v_i\in L_x(\sigma_E)\subset T_xM\,$ with 
$\lambda_0,
\lambda_1,\dots,\lambda_n\geq 0\,$ and $\sum_0^n\lambda_i=1$,  to the unique 
point 
$\phi(u)$ in $\hat\tau_E $ which minimizes the function 
$\sum_0^n\lambda_i\,d_g^2(p_i,{\cdot})$. 
As we already 
quoted in previous proposition \ref{2s}, the map $\phi$ provides $\hat\tau_E $ with a 
Euclidean metric $\hat g$, which is in fact deduced from 
a Euclidean metric $\bar g$ (see definition \ref{bareucl}) on 
$T_xM$ verifying $d_{\bar  g}(v_i,v_j)=d_g(p_i,p_j)$ (in general 
$\bar g\not=g_x$). 
Use $(i),(ii)$ to bring in theorem \ref{theorem 3} $(i')$ which asserts the 
existence of a positive constant
${\mathcal C}_7$ such that for 
$\rho \leq {\mathcal R}_4$
\begin{equation*}
\forall q\in \hat\tau_E \hskip4mm\forall u\in T_qM \hskip6mm
\vert\,\Vert\,u\,\Vert_g - \Vert\,u\,\Vert_{\hat g}\,\vert\leq 
{\mathcal C}_7 \ \rho^2\, \Vert\,u\,\Vert_g\  , 
\end{equation*}
which implies 
\begin{equation}
\label{Eq.1}\hbox{\rm vol}_g(\hat\tau_E)\geq (\frac{1}{1+{\mathcal C}_7 \ \rho^2})^n\,
\hbox{\rm vol}_{\hat g}(\hat\tau_E)= 
(\frac{1}{1+{\mathcal C}_7 \ \rho^2})^n\,\hbox{\rm vol}_{\bar g}(L_x(\sigma_E))\,,
\end{equation}
the last equality resulting from the definition of $\phi$. Using the 
definition 
of the Euclidean metric $\bar g$ through the equalities 
$d_g (\exp_xv_i,\exp_xv_j)=\Vert v_i-v_j\Vert_{\bar g}$ (for $i,
j=0,\dots,n$), recalling that ${\mathcal R}_4\leq {\mathcal R}_0(\check C)$ and that ${\mathcal R}_0(\check C)$ was chosen so that 
proposition \ref{proposition A} applies with $\check C$ given in 
definition \ref{2j} \ref{flop}, one has
\begin{equation}\label{2tb}
\vert\,\Vert\,v_i-v_j\,\Vert_{\bar g} - \Vert\,v_i-v_j\,\Vert_{g_x}\,\vert\leq
\check C\, \Vert\,v_i-v_j\,\Vert_{g_x}\ .
\end{equation}
That is, the mapping $\exp_x$ is $\check C$-quasi-isometric from a  
convex ball 
$B(0_x,\rho)\subset(T_xM,g_x)$ 
to another $B(x,\rho)\subset(M,g)$. Thus, integrating the infinitesimal 
control, it follows (again)
 that the $g_x$ and $g$-distances show a corresponding equivalence.
Using (\ref{1db}) in 
lemma \ref{1d}, the above line (\ref{2tb}) tells that the identity 
mapping 
sending $(L_x(\sigma_E), g_x)$ to $(L_x(\sigma_E), \bar g)$ is a $C_1$-quasi-isometric 
mapping
\begin{equation*} 
(1-C_1)\,\leq\,\Vert\,\,\hbox{\rm Id}\,\Vert_{g_x,\bar g}\,\leq \,(1+C_1)\,,
\end{equation*}
with $C_1=9n\check C/2\,t_1^2\,$, and we get
\begin{equation}\label{Eq.2}
\hbox{\rm vol}_{\bar g}(L_x(\sigma_E))\geq (1-C_1)^n\, 
\hbox{\rm vol}_{g_x}(L_x(\sigma_E))\, .
\end{equation}
From equations (\ref{Eq.1}), (\ref{Eq.2}) and proposition \ref{2r} $(iii)$, we derive $(iii)$
\begin{equation}\label{Eq.3}\hbox{\rm vol}_g(\hat\tau_E)\geq
(\frac{1-C_1}{1+{\mathcal C}_7 \ \rho^2})^n\,\hbox{\rm vol}_{g_x}(L_x(\sigma_E))\geq\frac{\overline{A}}{12^n}\,\rho^n=A_0\,\rho^n\  .
\end{equation}
Indeed, one knows that
$\,(1-C_1)/(1+{\mathcal C}_7 \ \rho^2)\geq 1/2\,$, for $\ 1/(1+{\mathcal C}_7 \ \rho^2)\geq 3/4$
(see definition \ref{2j} (\ref{flop})) and
$\,(1-C_1)\geq 2/3\,$ (since $C_1=9n\check C/2\,t_1^2\leq 1/3\,$, since $\check C$ was chosen $\leq 2t_1^2/(3^3n)$, again by definition \ref{2j} (\ref{flop})).
\end{proof}

{\bf Completing the proof of theorem \ref{2n}}
\begin{proof} Consider an embedded simplicial complex $T:(K,{\rm std})\rightarrow W\subset(M,g)$ and ${\mathcal R}_4>0$ satisfying the assumptions of definition \ref{2j}. Take a maximal family of distinct points $\{x_i\}_{i=1,\dots,N}$ in $\overline{W}$ such that any pair of points $\{x_i,x_j\}$ for $i\not=j$ satisfies $d(x_i,x_j)\geq {\mathcal R}_4/2\,$: thus, any point of $W$ sits in some ball $B(x_i,{\mathcal R}_4/2)$. Take $\rho\leq{\mathcal R}_5/4\,$ with ${\mathcal R}_5:=\inf (4{\mathcal D}_n\,\beta_T,{\mathcal R}_4,4\rho(\alpha_T/4))$, so $\rho\in ]0, {\mathcal R}_4/4]\,$, set $E=E(\rho)$ and consider the regular simplicial subdivision $K_E$ of integer $E$. Proposition \ref{2t} applies to $K_E$ and all its $n$-simplices. Any simplex $T(\sigma_E)$ has diameter less than $\rho$ by proposition \ref{2o} $(i)$, thus is included in some $B(x_i,{\mathcal R}_4)$. Call $K_{E,i}$ the sub-complex of $K_E$ obtained by collecting those $n$-simplices $\sigma_E$ for which $\tau_E=T(\sigma_E)$ meets $B(x_i,{\mathcal R}_4/2)$, together with all their faces. By proposition \ref{2o} $(i)$, $T(K_{E,i})$ is contained in $B(x_i,{\mathcal R}_4)$. Since any $n$-simplex in $K_E$ meets some of the $B(x_i,{\mathcal R}_4/2)$, the proof of the theorem is complete.
\end{proof}

\subsection{From a refined simplicial complex embedded in a Riemannian manifold to its Euclidean companion}\label{6.2}
$ $

From what has already been done, starting from a finite simplicial complex $(T,K)$ that 
triangulates a piece of a Riemannian manifold and refining it enough in a nice way, 
one can replace each 
embedded simplex $\tau_E=T(\sigma_E)$ by the Riemannian barycentric simplex $\hat\tau_E$ having same vertices which admits a Euclidean realization 
$(\hat\tau_E,\hat g_E)$ and verifies all claims of theorem \ref{theorem 3}.

\begin{defn} \label{hatT1}A real $\rho\!\in ]0,{\mathcal R_5}/4]$ is given with $E\!=\!E(\rho)$, according to definition \ref{linkint}, so that the results of section \ref{6.1} hold. Call $\hat T_E : \lvert K \lvert \rightarrow M$ the differentiable simplicial mapping (definition \ref{D.Simplicial.Embedding}) defined to be, on each $\sigma_E 
\in K_E\,$, the barycentric map 
that sends the vertices of $\sigma_E$ to the paired vertices of 
$\tau_E=T(\sigma_E)$ (see theorem \ref{2n}): a point which has affine barycentric coordinates in 
$\sigma_E$ is sent to the unique corresponding one 
having same {\it Riemannian} barycentric coordinates in $(M,g)$.
\end{defn}

\begin{thm} \label{hatT} If the integer $E$ is large enough, the regular simplicial subdivision $K_E$ of $K$ is such that the mapping $\hat T_E:K_E\rightarrow M$ is a differentiable simplicial embedding close to $T\,$. 
\par Precisely, for any 
$\delta\!>\!0$ small enough exists a positive integer $E_1(\delta) $ such that, for any 
integer $E\!\geq\! E_1(\delta)\,$, the
map 
$\hat T_E : (K_E,{\rm std}) \longrightarrow (M,g)$ is a $2\delta$-approximation to $T\,$, one has
\begin{equation}\label{soclosetoyou}\forall p\in \lvert K \lvert\ \ \ \ \ \ \ \ \ \ \ \ \ \ 
 d_g(\hat T_E (p), T(p)) \leq 2\delta \  .
\end{equation}
\end{thm}

\begin{rem} By its very definition, the new embedded simplicial complex 
$\hat T_E:K_E\rightarrow M$ has the same combinatorics as $T_E:K_E\rightarrow (M,g)$. If 
$E$ 
is large enough ($\rho$ is small enough) (\ref{soclosetoyou}) says {\it the embedded images 
$\hat T_E(K_E)$ and $T_E(K_E)$ cover almost the same piece of $(M,g)$.}
\end{rem}

\begin{proof} 
The proof is given in appendix \ref{Whitehead}. It relies on the following 
\begin{thm} \label{James1} Let $K$ be a finite simplicial complex. Let 
$f:K\rightarrow {\mathbb R}^n$ be a differentiable simplicial immersion $($or embedding$)$. 
There 
exists $\epsilon>0$ such that any strong $\epsilon$-approximation to $f$ is a differentiable simplicial immersion 
(or embedding). 
\end{thm}

See \cite[theorem 8.8]{Mu} for a proof of the last theorem, a central tool to establish J. H. C. Whitehead's existence theorem of {\em intrinsic} triangulations on manifolds (for an {\em extrinsic} approach, using Whitney's embedding theorem, see \cite{Cai2}). The fact that theorem \ref{James1} is also true if one replaces the target ${\mathbb R}^n\,$ by a manifold $M$ is left to the reader as (part) of an exercise in Munkres' book, see \cite[(d) page 89]{Mu}: theorem \ref{hatT} follows then from this ``extended'' version. Yet, in appendix \ref{Whitehead},
using a refined version of the last statement given in theorem \ref{2n}, which allows to deal with finitely many maps $\exp_{x_i}^{-1}\circ T$ taking values in the vector spaces $T_{x_i}M$, we mix theorem \ref{James1} with the quasi-isometric nature of the maps $\exp_x$ in $(M,g)$ to prove theorem \ref{hatT} .
\end{proof}

{\it From now on,} our concern turns exclusively towards the new simplicial complex embedding (definition \ref{hatT1})
$$\hat T_E:K_E\longrightarrow W\subset(M,g)\,.
$$

\begin{rem}
There are two natural metrics on  $\hat T_E(K_E)\subset W$. The one 
is 
the smooth Riemannian metric $g$, given from the beginning. The other is the 
piecewise flat metric $\hat g_E$ 
introduced by means of the Riemannian barycentric coordinates, see definitions 
\ref{bareucl} and \ref{m}. 
\end{rem}
And, to be complete, singular over the $(n-2)$-skeleton $K_{E,n-2}$ of $K_E\,$, there is a 
third natural {\it analytic and flat metric on} $K_E\setminus K_{E,n-2}\,$
\begin{defn}\label{D.metric1}
Denote by $g_0$ the metric 
obtained {\it by requiring the map 
$\hat T : (\sigma,  g_0) \rightarrow (\hat T(\sigma), \hat g_E)$ to be an isometry} when restricted to each $n$-simplex $\sigma\in K_E$ {\it and by viewing the metric union of two $n$-simplices $(\sigma_1,g_0)$ and $(\sigma_2,g_0)\,$ which are adjacent in $K_E$} along a common 
$(n-1)$-face {\it to sit isometrically in a same nice Euclidean (flat) $({\mathbb R^n},g_0)\,$}. This defines the analytic flat $g_0$ on $\sigma_1\cup\sigma_2\,$. Propagate $g_0$ from an $n$-simplex to all its adjacent neighbors in this way, this defines the expected global analytic flat Riemannian metric $g_0$ over $K_E\setminus K_{E,n-2}\,$.

Moreover, for any simplex $\eta\!\in\! K$, the metric $g_0$ induces a Euclidean metric on the tangent space $T(\eta\setminus\partial\eta)$.

Observe again that $g_0$ depends on $E\,$, though we mostly skip $E\,$.

Observe also that $(K,g_0)$ is a {\em polyhedron} with respect to 
the first instance of convention \ref{D.metric133}.
\end{defn}

To enable a comparison of the parallel translation in the metric $g$ with the 
one in 
the metric $\hat g_E\,$, we have to study and control the jump caused by the 
$C^1$-rip existing in the metric $\hat 
g_E$ at an interface between two different $n$-simplices $\hat\tau_1$ and 
$\hat\tau_2$ having a common $(n-1)$-face. 

\begin{defn} \label{D.metric}
Denote by $\phi_i$ (for  $i=1,2$) the map $\hat T_{\mid \sigma_i}$ and, 
for 
$q\in\hat\tau_1\cap\hat\tau_2\,$, define $\Phi_q=(d\phi_2\circ d\phi_1^{-1})(q)$. 
The union $\sigma_1\cup\sigma_2$ is equipped with the everywhere analytic Euclidean 
metric $g_0$ 
making of $\phi_1,\phi_2$ isometries from $(\sigma_1, g_0),(\sigma_2, g_0)$  
onto $(\hat\tau_1,\hat g_1),(\hat\tau_2,\hat g_2)$, where 
$\hat g_i=\hat g_{\mid \hat\tau_i}\,$.
\end{defn}

\begin{rem} 
Take $y\in\sigma_1\cap\sigma_2\,$, one has canonically 
$T_y\sigma_1\equiv T_y\sigma_2\,$. 
In general, if $u\in T_y\sigma_1=T_y\sigma_2\,$, the images $d\phi_1(y)(u)$, 
$d\phi_2(y)(u)$ differ. 
But one has $d\phi_1(y)=d\phi_2(y)$ on $T_y(\sigma_1\cap\sigma_2)$, i. e. 
one 
has 
$\Phi_q-\hbox{\rm Id}=0$ on $T_q(\hat \tau_1\cap \hat\tau_2)$, where 
$q= \phi_1(y)=\phi_2(y)$.
\end{rem}

\begin{prop}\label{rip1}
Given $\rho\in]0,{\mathcal R}_3]$ such that ${\mathcal C}_7\,\rho^2<1/9$ (definition {\rm \ref{constanteC_7R_3}}), 
if $\sigma_1$ and $\sigma_2$ have a common $(n-1)$-face, one has, for $i=1,2$
\begin{gather}\label{rip331}\forall q\in\hat\tau_1\cap\hat\tau_2  \hskip1cm
\Vert\,\Phi_q-\hbox{\rm Id}\,\Vert_{\hat g_i,\hat g_i}\leq 4\,\sqrt{2}\,\,
{\mathcal C}_7\,\rho^2\,(\frac{1-2\,{\mathcal C}_7\, \rho^2}{1-3\,
{\mathcal C}_7\,\rho^2})\,,
\\
\label{rip332}
\forall q\in\hat\tau_1\cap\hat\tau_2  \hskip1cm \Vert\,\Phi_q-
\hbox{\rm Id}\,\Vert_{g,g}\leq 12\,\sqrt{2}\,{\mathcal C}_7\,\rho^2\  .
\end{gather}
\end{prop}

\begin{proof}
Given $q\in \hat\tau_1\cap\hat\tau_2\,$, let $a,u_1,u_2$ be the $g,\hat g_1$ and 
$\hat g_2\,$-unit normal 
vectors to $T_q(\hat\tau_1\cap\hat\tau_2)$, pointing from 
$\hat\tau_1$ to $\hat\tau_2 \,$. Thus $\Phi_q(u_1)=u_2\,$ and 
one can state and prove the following
\begin{lem} \label{rip1-1} 
One has $(\Vert\,u\,\Vert_1$ stands for 
$\Vert\,u\,\Vert_{\hat g_1}\,$, etc...$)$
\begin{equation*}
\Vert\,\Phi_q-\hbox{\rm Id}\,\Vert_{1,1}=\Vert\,u_2-u_1\,\Vert_1\  .
\end{equation*}
\end{lem}

\begin{proof}
Any $v\in T_qM$ such that $\Vert\,v\,\Vert_1=1$ may be written
\begin{equation*}
v=\cos\theta\,u_1+\sin\theta\,v_1\hskip2mm\hbox{with}\hskip2mm
\langle u_1,v_1\rangle_1=0 \  , \Vert\,u_1\,\Vert_1=\Vert\,v_1\,\Vert_1=1\  .
\end{equation*}
As $\Phi_q-\hbox{\rm Id}=0$ on $T_q(\hat\sigma_1\cap\hat\sigma_2)$, one gets
\begin{equation*}
\Vert\,\Phi_q(v)-v\,\Vert_1=\Vert\,\cos\theta\,(u_2-u_1)\,\Vert_1\leq
\Vert\,u_2-u_1\,\Vert_1\  ,
\end{equation*}
the conclusion follows.
\end{proof}

One also needs the following estimate, a consequence of theorem \ref{theorem 3}.
\begin{lem} \label{rip1-2} 
Let $\rho\in]0,{\mathcal R}_3]$ be such that ${\mathcal C}_7\,\rho^2<1/9$ (see proposition {\rm\ref{rip1}}). One has 
\begin{equation*}
\max (\Vert\,a-u_1\,\Vert_1, \Vert\,a-u_2\,\Vert_2)\leq 
2\sqrt{2}\,{\mathcal C}_7\, 
\rho^2\  .
\end{equation*}
\end{lem}

\begin{proof}
Take $v\!\in \!T_q(\hat\tau_1\cap\hat\tau_2)$ with $\Vert v\Vert_g\!=\!1$, so (use $g(a,v)\!=\!\hat g_1(u_1,v)\!=\!0$)
\begin{equation}\label{rip2}
\vert \hat g_1(a,v)- \hat g_1(u_1,v) \vert=\vert \hat g_1(a,v)\vert=\vert \hat 
g_1(a,v)- g(a,v) \vert\  .
\end{equation}
One writes
\begin{multline}\label{rip3}
\vert \hat g_1(a,v)- g(a,v) \vert \leq
\frac{1}{2}(\vert \hat g_1(a+v,a+v)- g(a+v,a+v) \vert + \\ 
+\vert \hat g_1(a,a)- g(a,a) \vert  + \vert \hat g_1(v,v)- g(v,v) \vert) \  .
\end{multline}
Thanks to theorem \ref{theorem 3}, one knows
\begin{equation*}
\vert \hat g_1(a,a)- g(a,a) \vert
\leq {\mathcal C}_7\,\rho^2\,g(a,a)={\mathcal C}_7\,\rho^2\  ,
\end{equation*}
so, one has
\begin{equation} \label{rip5}
\vert \hat g_1(a,a)- \hat g_1(u_1,u_1) \vert
\leq {\mathcal C}_7\,\rho^2\  .
\end{equation}
From theorem \ref{theorem 3} and (\ref{rip2}), (\ref{rip3}), one gets (use $g(a,v)=0$)
\begin{equation}\label{rip6}
\vert \hat g_1(a,v)- \hat g_1(u_1,v) \vert\leq \frac{1}{2}\,
{\mathcal C}_7\,\rho^2(2\,g(a,a)+2\,g(v,v))=2\,{\mathcal C}_7\,\rho^2\  .
\end{equation}
Writing 
$a=\alpha_1\,u_1+\beta_1\,v_1$ with $\Vert u_1\Vert_1=\Vert v_1\Vert_1=1$ 
and $\langle u_1,v_1 \rangle_1=0$ (hence $v_1$ is in 
$T_q\hat\tau_1\cap\hat\tau_2$), one derives from (\ref{rip6}) if 
$\ 9\,{\mathcal C}_7\,\rho^2\leq 1$
\begin{equation*}
\vert \hat g_1(a-u_1,v_1) \vert=\beta_1\leq 2\,{\mathcal C}_7\,\rho^2
\Vert v_1\Vert_g
\leq\frac{2\,{\mathcal C}_7\,\rho^2}{1-{\mathcal C}_7\,\rho^2}\leq 
(\frac{3}{2})^2\,{\mathcal C}_7\,\rho^2\  ,
\end{equation*}
which implies $\beta_1^2\leq (\frac{3}{2})^4\,{\mathcal C}_7^2\,\rho^4$,
and from (\ref{rip5}), one deduces
\begin{multline}\label{rip8}
\vert \hat g_1(\alpha_1u_1+\beta_1v_1,\alpha_1u_1+\beta_1v_1)\!-\! 
\hat g_1(u_1,u_1) \vert
\!=\!\vert\alpha_1^2+\beta_1^2-1\vert\!\leq\! {\mathcal C}_7\,\rho^2\, .
\end{multline}
As $\ 9\,{\mathcal C}_7\,\rho^2\leq 1\,$, from 
$\beta_1^2\leq (\frac{3}{2})^4\,{\mathcal C}_7^2\,\rho^4\leq 
\frac{3^2}{2^4}\,{\mathcal C}_7\,\rho^2$ and (\ref{rip8}), one has 
$\vert1-\alpha_1^2\vert\,\leq \beta_1^2+{\mathcal C}_7\,\rho^2\leq 
(1+\frac{3^2}{2^4})\,{\mathcal C}_7\,\rho^2\,$;
using 
$\vert1-\alpha_1^2\vert=\vert1-\alpha_1\vert\,
(1+\alpha_1)\geq \vert1-\alpha_1\vert$ 
(the choice of $a$ and $u_1$ implies $\alpha_1\geq0$), one gets
\begin{multline*}
\Vert \, a- u_1\, \Vert_1=\sqrt{(1-\alpha_1)^2+\beta_1^2}\leq \\
\leq
\big((1+\frac{3^2}{2^4})^2+(\frac{3}{2})^4\big)^{\frac{1}{2}}\ {\mathcal C}_7\,\rho^2=
(\frac{25^2+3^4\,2^4}{2^8})^{\frac{1}{2}}\ 
{\mathcal C}_7\,\rho^2=\\=(\frac{1921}{256})^{\frac{1}{2}}\ 
{\mathcal C}_7\,\rho^2\leq 2\,\sqrt{2}
 \ {\mathcal C}_7\,\rho^2\  ,
\end{multline*}
and in the same way
\begin{equation*}
\Vert \, a- u_2\, \Vert_2 \leq 2\,\sqrt{2}\,{\mathcal C}_7\, \rho^2\ . \qedhere
\end{equation*}
\end{proof}

Applying lemma \ref{rip1-1}, write
\begin{equation} \label{rip1-5}
\Vert\,\Phi_q-\hbox{\rm Id}\,\Vert_{1,1}=\Vert\,u_2-u_1\,\Vert_1\leq 
\Vert\,u_2-a\,\Vert_1+\Vert\,u_1-a\,\Vert_1\  .
\end{equation}
Applying again theorem \ref{theorem 3}, one gets
\begin{equation} \label{rip1-4}
\Vert\,u_2-a\,\Vert_1\leq (\frac{1-{\mathcal C}_7\,\rho^2}
{1-3\,{\mathcal C}_7\,\rho^2})\,\Vert\,u_2-a\,\Vert_2\  .
\end{equation}
{\it Indeed,} theorem \ref{theorem 3} $(i')$ allows to write
\begin{gather}
\vert\,\Vert\,u_2-a\,\Vert_2-\Vert\,u_2-a\,\Vert_g\,\vert\leq 
{\mathcal C}_7\, \rho^2\,\Vert\,u_2-a\,\Vert_g\  ,\notag \\
\label{rip1-3}\vert\,\Vert\,u_2-a\,\Vert_1-\Vert\,u_2-a\,\Vert_g\,\vert \leq  
{\mathcal C}_7\, \rho^2\,\Vert\,u_2-a\,\Vert_g\  ,
\end{gather}
from which we get
\begin{equation}\label{rip10}
\vert\,\Vert\,u_2-a\,\Vert_2-\Vert\,u_2-a\,\Vert_1\,\vert\leq 
2\,{\mathcal C}_7\, \rho^2\,\Vert\,u_2-a\,\Vert_g\  ;
\end{equation}
but one also has from (\ref{rip1-3})
\begin{equation}\label{rip11}
(1-{\mathcal C}_7\,\rho^2)\,\Vert\,u_2-a\,\Vert_g\leq \Vert\,u_2-a\,\Vert_1\  ,
\end{equation}
and (\ref{rip10}), (\ref{rip11}) give
\begin{equation*}
\vert\,\Vert\,u_2-a\,\Vert_2-\Vert\,u_2-a\,\Vert_1\,\vert\leq 
(\frac{2\,{\mathcal C}_7\, \rho^2} 
{1-{\mathcal C}_7\,\rho^2})\,\Vert\,u_2-a\,\Vert_1\  ;
\end{equation*}
this, in turn, implies (\ref{rip1-4}).
One deduces (\ref{rip331}) from (\ref{rip1-4}), (\ref{rip1-5}) and lemma \ref{rip1-2}.
Apply once more theorem \ref{theorem 3} and state
\begin{equation*}
\Vert\,\Phi_q-\hbox{\rm Id}\,\Vert_{g,g}\leq (\frac{1+{\mathcal C}_7\, 
\rho^2}{1-{\mathcal C}_7\,\rho^2})\,\Vert\,\Phi_q-\hbox{\rm Id}
\,\Vert_{1,1}\,\ .
\end{equation*}
Get 
$(1+{\mathcal C}_7\rho^2)\!/\!(1-{\mathcal C}_7\rho^2)\!\leq\! 5/4$ and 
$(1\!-\!2\,{\mathcal C}_7\rho^2)\!/\!(1-3\,{\mathcal C}_7\rho^2)\!\leq\! 7/3$ from ${\mathcal C}_7\rho^2\!\leq\! 1/9\,$. To prove (\ref{rip332}) of proposition \ref{rip1}, write
\begin{multline*}
\Vert\,\Phi_q-\hbox{\rm Id}\,\Vert_{g,g}\leq 4\,\sqrt{2}\,
{\mathcal C}_7\,\rho^2\,
(\frac{1+{\mathcal C}_7\, \rho^2} {1-{\mathcal C}_7\,\rho^2})
\,(\frac{1-2\,{\mathcal C}_7\, \rho^2}{1-3\,{\mathcal C}_7\,\rho^2})\leq \\ 
\leq 4\,\sqrt{2}\,\frac{5}{4}\,\frac{7}{3}\,{\mathcal C}_7\,\rho^2\leq
12\,\sqrt{2}\,{\mathcal C}_7\,\rho^2\  .\qedhere
\end{multline*}
\end{proof}

\subsection{Riemannian barycentric textures}\label{6.3}
$ $

To ensure the increasing closeness (with $E$) of the $g$ and $\hat g_E$-parallel 
translations along a 
given curve in $M$, as one regularly refines the source $K$ (of the 
simplicial embedding $\hat T
$ in $W\subset(M,g)$) to $K_E$,  {\it one 
must control the increase of intersections between the curve and the 
$(n-1)$-faces of Riemannian barycentric $n$-simplices,} 
each crossing of which may cause a rip. No jump in the {\it Euclidean} plane, but, in what concerns the alluded bound, using an upper bound on the {\it local degree of intersection} of a generic curve with a straight line (well known to be $\leq 3$), one can find, in terms of $E\,$, a bound on the number of intersections of a generic curve of finite length with all linear sides of Euclidean triangles building an increasing $E$-regularly refined linear triangulation. In the general case, a work on singularities founded on a theorem of 
Ren\' e Thom plays 
the central role
(see the book of Marcel Berger \cite{Be2}, V.16, pages 390-1, who indicated 
\cite{Th} to us, see also \cite{C-M}). To be able to apply the results we need, generalizing the role played by the set of straight lines in the Euclidean plane, we have to 
construct a ``texture'' whose 
pieces, roughly speaking, are all Riemannian barycentric $k$-simplices built on any 
``non-degenerating'' 
set of $k+1$ points of a convex ball of $(M,g)$. This brings in a new kind 
of magnifying glass 
around a point, which formalizes the passage to the limit done by making the mesh of a 
barycentric 
triangulation thinner and thinner while making the integer $E$ larger and 
larger.

\subsubsection{The manifold of linear simplices centered at $p$ in $(T_pM,g_p)$}
\label{linsimp}
\begin{defn} 
Write ${\bf v}$ short for a set of $k+1$ tangent vectors $(v_0,\dots,v_k)$ in 
$T_pM$ for  $p\in M$ and $k\in\{0 , \dots , n\}$.

 If $k\geq 1$ and $(v_1-v_0, \dots , 
v_k-v_0)\not=(0,\dots,0)$, the standard ratio of ${\bf v}$ is defined to be 
\begin{equation*}
{\mathcal{S}}({\bf v})=\frac{\Vert(v_1-v_0)\wedge\dots\wedge(v_k-v_0)\Vert^2}
{(\frac{1}{k}(\Vert v_1-v_0\Vert^2 +\dots +\Vert v_k-v_0\Vert^2))^k}\  .
\end{equation*}
Denote by ${\mathcal N}({\bf v})$ the numerator and by ${\mathcal D}({\bf v})$ 
the denominator
of ${\mathcal{S}}({\bf v})$. 
\end{defn}

\begin{defn}
Define the manifold of centered $k$-affine frames of $M$
\begin{gather*}
{\mathcal F}(k,TM):=\bigcup_{p\in M}{\mathcal F}(k,T_pM)\ 
\hskip2mm \hbox{where} \hskip2mm\\ 
{\mathcal F}(k,T_pM)\!:=\!\{{\bf v}\in (T_pM)^{k+1}\!\mid\!
{\bf v} \hskip1mm \hbox{is affinely free and} \hskip1mm v_0 + \dots + v_k=0_p\}\,.
\end{gather*}
\end{defn}

\begin{rem} 
One can naturally relate a given ${\bf v}$ to the linear simplex 
$\sigma_{\bf v}$ in $T_pM$ having vertices $(v_0,\dots,v_k)$. If ${\bf v}$ is in 
${\mathcal F}(k,T_pM)$, this $\sigma_{\bf v}\subset T_pM$ is a nondegenerate 
$k$-simplex {\it having barycentre at the origin $0_p$} (we also say that is is {\it centered at $p$}). 
\end{rem}

\begin{defn}
Recall that the {\it $k$-standard simplex} is 
the simplex having 
vertices 
$0,e_1,\dots,e_k$ in $\R^k$ (where $e_1,\dots,e_k$ is the canonical basis).  
{\it It is not regular} !  
\end{defn}

The standard ratio ${\mathcal{S}}({\bf v})={\mathcal{S}}(\sigma_{\bf v})$ of a 
simplex $\sigma_{\bf v}$ (this name is justified by lemma \ref{critical0}, just below) is another measure close to ``openness'' and 
``thickness'', a link is given by the
\begin{prop} \label{srop} 
If ${\bf v}$ is in ${\mathcal F}(k,T_pM)$, the following inequalities hold for 
$\sigma_{\bf v}$
$($see definition {\rm \ref{openness1}} of $\omega)$
\begin{equation*} 
(k!)^2\,\omega^2(\sigma_{\bf v})\leq {\mathcal{S}}({\bf v})\leq 
(k!)^2\,k^k\,\omega^2(\sigma_{\bf v})\  .
\end{equation*}
\end{prop}

\begin{proof} 
This is a consequence of the inequalities 
\begin{equation*}
\forall \alpha_1,\dots,\alpha_k\in{\mathbb R}\hskip2cm\sup_{i=1}^k\,\alpha_i^2\leq\sum_{i=1}^k \alpha_i^2\leq 
k\,\sup_{i=1}^k\,\alpha_i^2
\end{equation*}
and of the equality 
$\Vert(v_1-v_0)\wedge\dots\wedge(v_k-v_0)\Vert=k!\,\hbox{\rm vol}_k(\sigma_{\bf v})$.
\end{proof}

\begin{lem} \label{critical0}
The function ${\mathcal{S}}:{\mathcal 
F}(k,T_pM)\rightarrow \R$ is submersive outside ${\mathcal{S}}^{-1}\{1\}$ 
and $1$ is its maximum value. In fact, one has the 
equivalence between the following three properties 

 $(i)$ \hskip4mm the point ${\bf v}$ is critical for ${\mathcal{S}}$ ;

 $(ii)$ \hskip4mm ${\mathcal{S}}({\bf v})=1$ ;

 $(iii)$ \hskip4mm $(v_1-v_0,\dots,v_k-v_0)$ are orthogonal and share the same norm.
\end{lem}

\begin{proof}
The $\log$-concavity (the reals $\alpha_i$ are $>0$)
\begin{equation*}
\log(\frac{\alpha_1+\dots+\alpha_k}{k})\geq\sum_{i=1}^k\frac{\log\alpha_i}{k}
\end{equation*}
and Hadamard's inequality give
\begin{equation}\label{Hadamard}
{\mathcal{S}}({\bf v})\leq\frac{\Vert(v_1-v_0)\wedge\dots\wedge(v_k-v_0)\Vert^2}
{\Vert v_1-v_0\Vert^2 \cdots \Vert v_k-v_0\Vert^2}\leq 1\  .
\end{equation}
Actually, the maximum of ${\mathcal{S}}$ is $1$ and is achieved precisely on the set of 
those $\sigma_{\bf v}$ which differ from the standard 
simplex 
by an isometry that sends $v_0$ on $0\in\R^k$ ($\sigma_{\bf v}$ is said to be 
``standard at'' $v_0$) followed by a homothety of center $0\,$. 
From this follows $(iii)$ {\it implies} $(ii)$.

Looking at the definition of the standard ratio, one computes directly 
${\mathcal{S}}(\sigma_{\bf v})=1$ on a ${\bf v}$ such that $\sigma_{\bf v}$ is 
``standard at'' $v_0$, so ${\mathcal{S}}$ achieves there a maximum (because of 
(\ref{Hadamard})). Moreover, if $\sigma_{\bf v}$ is a point where 
${\mathcal{S}}(\sigma_{\bf v})=1\,$,  so is $\sigma_{\bf v}$ a point where ${\mathcal{S}}$ 
achieves its maximum and which is interior in ${\mathcal F}(k,T_pM)$ (by definition of 
this set), so $\sigma_{\bf v}$ is a critical point of ${\mathcal{S}}\,$, thus 
$(ii)$ {\it implies} $(i)$.

To show 
 $(i)$ {\it implies} $(iii)$, one computes the derivative of ${\mathcal{S}}$ at a point 
${\bf v}\in{\mathcal F}(k,T_pM)$ 
acting on ${\bf h}\in(T_pM)^{k+1}$ such that $h_0+\dots+h_k=0\,$  
\begin{multline}\label{critical}
\frac{1}{2}\,d\,{\mathcal{S}}({\bf v})({\bf h})\,{\mathcal D}^2({\bf v}) = \\ 
= \sum_{i=1}^k\langle 
(v_1-v_0)\wedge\dots\wedge(v_k-v_0),
(v_1-v_0)\wedge\dots\wedge(h_i-h_0)\wedge\dots\wedge(v_k-v_0)\rangle{\mathcal D}
({\bf v}) \\ 
-(\sum_{i=1}^k\langle v_i-v_0 , h_i-h_0\rangle)\,{\mathcal D}^{\frac{k-1}{k}}
({\bf v})\,{\mathcal N}({\bf v})\,.
\end{multline}
From this equality, one checks
\begin{lem}\label{critical1} 
Critical ${\bf v}$ implies $v_i-v_0$ orthogonal to $v_j-v_0$ for $i\not=j\,$.
\end{lem}
\begin{proof} The subspace generated by 
$v_1-v_0,\dots,v_k-v_0\,$ is denoted by $V_{\bf v}$. 
Choose
${\bf h}=(h_0,\dots,h_0,-k\,h_0)$ with $h_0$ being any vector in $V_{\bf v}$ orthogonal 
to $v_k-v_0$ and read what (\ref{critical}) tells at a critical ${\bf v}\,$. One 
gets (as 
$\,{\mathcal D}({\bf v})\not=0$ and the last term in (\ref{critical}) is $0$)
\begin{equation} \label{critical2}
0\!=\! 
\langle \!(v_1-v_0)\wedge\dots\wedge\!(v_k-v_0),(v_1-v_0)\wedge\dots\wedge\!(v_{k-1}
-v_0)\wedge(\!-k-1)h_0\rangle \,.
\end{equation}
So (\ref{critical2}) says that the orthogonal subspace to $v_k-v_0$ in 
$V_{\bf v}$ 
is contained in the subspace generated by
$v_1-v_0,\dots,v_{k-1}-v_0\,$. Caring about dimensions, this proves lemma \ref{critical1} by doing 
the same deduction exchanging the role just played by index $k$ with index $1\,$, etc. 
up to $k-1\,$.
\end{proof}
\begin{lem}\label{critical3} 
At a critical ${\bf v}\,$, all $v_i-v_0$ have the 
same norm.
\end{lem}
\begin{proof} 
Making the choice $2\,{\bf h}=(v_0-v_1,v_1-v_0,0,\dots,0)\,$, one gets this time, 
injecting this value of ${\bf h}$ in (\ref{critical}) (use $v_1-v_0$ 
orthogonal to $v_i-v_0$ with $i\not=1\,$, from lemma \ref{critical1})
\begin{equation*}
0={\mathcal N}({\bf v})\,{\mathcal D}({\bf v})-\Vert v_1-v_0\Vert^2\,
{\mathcal D}^{\frac{k-1}{k}}({\bf v})\,
{\mathcal N}({\bf v})\,,
\end{equation*}
and this gives $\Vert v_1-v_0\Vert^2={\mathcal D}^{\frac{1}{k}}({\bf v})$, hence 
lemma \ref{critical3} is proven (do the same thing for any index in place of $1$).
\end{proof}

Lemmas \ref{critical1}, \ref{critical3} show $(i)\Rightarrow(iii)$ in lemma \ref{critical0}, 
now proven.
\end{proof}

\begin{defn}\label{gloups1}$ $
\begin{enumerate}
\item Define ${\mathbb F}_s^k(p)\!:=\!{\mathcal S}^{-1}\{s\}\cap {\mathcal F}(k,T_pM)$, then set  
${\mathbb F}_s^k\!:=\!\cup_{p\in M}{\mathbb F}_s^k(p)$ and ${\mathbb F}_{\geq s}^k\!:=\!\cup_{s'\geq s}{\mathbb F}_{s'}^k\!=\!{\mathcal S}^{-1}[s,+1]\cap {\mathcal F}(k,T_pM)\,$.
Thus, one has $\ {\mathbb F}_{>0}^k=\cup_{s>0}{\mathbb F}_s^k={\mathcal F}(k,TM)$.
\par {\it If not confusing, we skip $k\,$, writing
${\mathbb F}_s$ instead of ${\mathbb F}_s^k\,$.}

\item Denote by ${\mathcal H}_0(p)\subset (T_pM)^{k+1}$ the  $n$-codimensional
subspace of ${\bf v}$ whose components $v_i$
sum
up to $0\,$.
\end{enumerate}
\end{defn}

A way to rephrase the previous results is the matter of the next

\begin{prop} \label{gloups2}
Given $s\in\,]0,1[\,$, the set 
${\mathbb F}_{\geq s}(p)\subset {\mathcal H}_0(p)\subset (T_pM)^{k+1}$ is a conical 
submanifold with boundary of maximal dimension $nk$ and vertex $0\,$. The 
manifold ${\mathbb F}_{\geq s}$ is the total space of a fibre-bundle with fibre 
${\mathbb F}_{\geq s}(p)$ which retracts by deformation $($the retraction is performed fibrewise$)$ onto a 
subset which is itself a fibre-bundle fibered over $M$, with conical fibres 
$\left((O(T_pM)/O(n-k))\times\R_+\right)/\sim$, where the equivalence relation $\sim$ identifies all points $({\bf v},0)\,$, i. e. for all ${\bf v},{\bf v}'$ one has
$({\bf v},0)\sim({\bf v}',0)$. For any $p\in M$, the group 
$\R_+^\ast$ 
acts on ${\mathbb F}_{\geq s}(p)$ (thus on ${\mathbb F}_{\geq s}$) through the 
homotheties of center $0_p$ in $T_pM\,$.
\end{prop}

\begin{proof} 
The previous lemma \ref{critical0} says precisely that ${\mathcal S}$ is a 
submersion outside the critical set
${\mathcal S}^{-1}\{1\}\,$, so the first part follows by noticing that ${\mathcal F}(k,T_pM)$ 
is an open set in ${\mathcal H}_0(p)\subset(T_pM)^{k+1}$.

The second part is deduced by 
following the gradient 
lines of ${\mathcal S}$ over ${\mathbb F}_{\geq s}\,$, which necessarily lead to a critical 
point  in ${\mathcal S}^{-1}\{1\}\,$, note that, by homogeneity, one can restrict the study 
to a compact level 
of ${\mathcal D}\,$. Therefore this provides (according to lemma \ref{critical0}), 
 up to homothetical factor, an orthonormal system of $k$ vectors, hence it gives a 
deformation that retracts ${\mathbb F}_{\geq s}$ on the expected total space.
\end{proof}

Yet we need another result related to the previous.
\begin{lem} \label{submersion} For any $s\in\,]0,1[\,$, the smooth function 
$\nu^2$ defined by $\nu^2:={\mathcal D}^{1/k}$ is submersive on each 
${\mathcal S}^{-1}\{s\}$ and thus has level sets which are submanifolds transversal to 
${\mathcal S}^{-1}\{s\}\,$.
\end{lem}

\begin{proof} 
It is enough to look upon the generatrix of the cone ${\mathcal S}^{-1}\{s\}$ 
passing through a point 
${\bf v}$, more precisely to compute the effect of ${\mathcal D}({\bf v})$ on 
the vector ${\bf h}={\bf v}\,$. One gets
\begin{equation*}
d {\mathcal D}({\bf v})({\bf v})=k\,{\mathcal D}({\bf v})\not=0\  ,
\end{equation*}
thus the result.
\end{proof}

\begin{defn}
\label{cylindrification}
Let $s$ be in $\,]0,1[\,$. Choose some fixed real 
$\alpha>0\,$, {\it which shall be given a precise value later, see convention {\rm\ref{convention3}}} below. Define 
${\mathbb F}_{\geq s,\alpha}^k:={\mathbb F}_{\geq s}^k\cap\nu^{-1}\{\alpha\}\,$.  
Define the {\it cylinders} $\Gamma_{\geq s}^k$ and $\Gamma_{>0}^k$ to be\begin{equation} \label{Kylinder} 
\Gamma_{\geq s}^k:=({\mathbb F}_{\geq s}^k\cap\nu^{-1}\{\alpha\})\times
\R={\mathbb F}_{\geq s,\alpha}^k\times\R\  \ \ \hbox{and}\ \ \ 
\Gamma_{>0}^k:=\cup_{s > 0}\Gamma_{\geq s}^k\ .
\end{equation}
{\it Recall that we skip $\,k\,$ if its mention seems redundant.} \end{defn}

\begin{rem} \label{srth} $ $
\begin{enumerate}
\item If ${\bf v}\in {\mathbb F}_{\geq s}\cap\nu^{-1}\{\alpha\}$ and $t>0$ are given, then
$t\,{\bf v}\in {\mathbb F}_{\geq s}$ gives rise to  $({\bf v},t)\in\Gamma_{\geq s}\,$, so 
$({\bf v},t)$ {\it 
can be now thought} for any $t\in\R\,$. This allows to extend 
mappings defined on ${\mathbb F}_{\geq s}$ as mappings on 
the cylinder $\Gamma_{\geq s}\,$, which is a new manifold with boundary, fibered over $M\,$, see proposition
\ref{text0}.
 
\item If, for a simplex $\sigma_{\bf v}\,$, one has  ${\mathcal S}(\sigma_{\bf v})\geq s>0$, 
then one also has (see definition \ref{1a})
$t(\sigma_{\bf v}) \geq t >0$  for 
some $t= t(s) >0\,$.

Indeed, by lemma \ref{1c} and 
proposition \ref{srop}, one knows 
that
\begin{equation*}
t(\sigma_{\bf v})\geq\frac{k}{k+1}
\frac{\omega(\sigma_{\bf v})}{\vert \sigma_{k-1}\vert}\geq 
\frac{k}{(k+1)!\,k^{k/2}}\frac{\sqrt{{\mathcal S}(\sigma_{\bf v})}}{\vert 
\sigma_{k-1}\vert}\  ,
\end{equation*}
hence one can choose $t=\frac{k}{(k+1)!\,k^{k/2}}\frac{s^{1/2}}{\vert \sigma_{k-1}\vert}\,$.

\item In the same way, using lemma \ref{1c} and 
proposition \ref{srop}, it follows that,  if the thickness $t(\sigma_{\bf v})$ of a 
simplex $\sigma_{\bf v}$ is $\geq t>0\,$, then  the 
standard ratio ${\mathcal S}(\sigma_{\bf v})$ will be $\geq s>0\,$, where 
$s=s(t)$ can be given in terms of $t\,$.
\end{enumerate}

\end{rem}

\begin{convent} \label{convention3} We now fix $\alpha:=1/2\sqrt{k}$ (see definition {\rm\ref{cylindrification}}). 
\end{convent}

\begin{lem}\label{text2}
Given $s_0>0\,$, put
$ {\mathbb F}_{>0,\alpha}={\mathbb F}_{>0}\cap\nu^{-1}\{\alpha\}\,$.
One has 
\begin{equation*}
\forall {\bf v}\in{\mathbb F}_{>0,\alpha}(p) \hskip6mm  \hskip6mm
\sigma_{\bf v}
\subset  B(0_p,1)\,.
\end{equation*}
Consequently, for any $\rho >0\,$, one has  
\begin{equation*}
\forall {\bf v}\in {\mathbb F}_{>0,\alpha\rho}
\hskip6mm  \hskip6mm
\sigma_{\bf v}
\subset B(0_p,\rho)\, .
\end{equation*}
\end{lem}

\begin{proof} 
If ${\bf v}$ is in ${\mathbb F}_{>0,\alpha}(p)$, then
\begin{equation*}
\Vert v_1-v_0\Vert^2+\dots+\Vert v_k-v_0\Vert^2= \alpha^2k\  .
\end{equation*}
Since ${\mathcal S}(\sigma_{\bf v}) >0\,$, we have $v_i\not= v_j$ if $i\not= j\,$. 
Thus (convention \ref{convention3})
\begin{equation*}
\max_i \Vert v_i-v_0\Vert < \alpha\sqrt{k}  \qquad \hbox{and} \qquad
\max_{i,j} \Vert v_i-v_j\Vert < 2\alpha\sqrt{k}=1\ .
\end{equation*}
Hence, the gravity center $0_p$ of $\sigma_{\bf v}$ is at distance $< 1$ 
from 
each $v_i\,$. Thus, for any ${\bf v}\in{\mathbb F}_{>0,\alpha}(p)$ one has 
$\sigma_{\bf v}\subset B(0_p,1)$. 
\end{proof}

\subsubsection{Organization of a new magnifying glass.}\label{newmagnglass}
$ $
\begin{defn} \label{barcoord}$ $
\begin{enumerate}
\item\label{first}
Consider the general case where ${\bf v}$ is the set of $k+1$ affinely free vertices 
of a linear simplex $\sigma_{\bf v}$ in the Euclidean space $(T_pM, g_p)$ which has barycentre $\tilde v$ and define $v_i(t)\!:=\!t\,(v_i-\tilde v)+\tilde v$. 
Performing the {\it first} homothety $h_t$ which sends (for $i\!=\!0,\cdots,k$) each $v_i\in{\bf v}$ to $h_t(v_i)\!:=\! v_i(t)$
produces the simplex $\sigma_{{\bf v},t}$ whose 
vertices are the points $v_i(t)$. 
\item \label{homothetie}
The barycentric coordinates 
${\bf \lambda}(t)\!=\!(\lambda_0(t),\dots,\lambda_k(t))$ linked to $\sigma_{{\bf v},t}$ (coordinates in the $k$-subspace generated by 
$\sigma_{\bf v}$) may be expressed in terms of those linked to $\sigma_{\bf v}\,$,
namely ${\bf \lambda}=(\lambda_0,\dots,\lambda_k)$, one has for $i=0,1,\dots,k$
\begin{gather*}
\lambda_i(t)\!=\!(\lambda_i-\frac{1}{k+1})\,\frac{1}{t}+\frac{1}{k+1}
\,, \ \hbox{defining the homothety}\\ 
\ \ h_t(\lambda):=\mu=\lambda(t)=(\lambda_0(t),\dots,\lambda_k(t))\ .
\end{gather*}
Notice that this {\it second} homothety $h_t:\lambda\mapsto \lambda(t)$ is the dilation of center $(\frac{1}{k+1},\dots,\frac{1}{k+1})$ and ratio ${1}/{t}\,$ and is ``dual'' to the first, which is a dilation of center $\tilde v$ and ratio $t\,$.
\item
The {\it space of barycentric coordinates in a $k$-affine space} is 
\begin{equation*}
{\mathcal H}^k=\{{\bf \lambda\in\R^{k+1}\mid 
\lambda_0+\dots+\lambda_k=1}\}\  .
\end{equation*}
\end{enumerate}
\end{defn}
Apollonius's device to recover the point 
having new barycentric coordinates ${\bf \lambda}(t)$ take the form of a kind of
``magnifying glass'' around $\tilde v$ of scaling ratio $1/t$
\begin{multline*}
\hbox{\it Find the point ${\mathcal Q}_{{\bf \lambda},{\bf v}}(t):=q(t)=q$ where}
\hskip3mm \\
\psi_{{\bf \lambda},{\bf v},t}(\cdot):=\sum_{i=0}^k\lambda_i(t)\,\Vert v_i(t)-\cdot\Vert^2 \\
\hbox{\it achieves its minimum ($\lambda_i(t)$ is given in definition {\rm\ref{barcoord} (\ref{homothetie})}).}
\end{multline*}

In the Riemannian case, with ${\bf v}\in{\mathcal F}(k,T_pM)$ and  
$p_i(t)=\exp_p(tv_i)$, 
\begin{multline} \label{fontionAR}
\hbox{\it Find the point ${\mathcal Q}_{{\bf \lambda},{\bf v}}(t):=q(t)=q$ 
where}\hskip3mm \\  
\psi_{{\bf \lambda}\,,
{\bf v},t}(\cdot):=\sum_{i=0}^k\lambda_i(t)\,d^2(p_i(t),\cdot) \\
\hbox{\it achieves its minimum ($\lambda_i(t)$ given in definition {\rm\ref{barcoord} (\ref{homothetie})}).}
\end{multline}
\begin{rem} \label{newglasses}
Thus, in a given convex ball $B$, the Riemannian device is organized  around each of its point $p\,$: for each given Riemannian barycentric simplex $\hat\tau\subset B$ centered at $p$, one considers the whole family of Riemannian barycentric simplices whose vertices are, for $t\in]0,1]\,$, ``$t$-homothetic'' to the vertices of the given one and relies on the ``associated'' barycentric coordinates ${\bf \lambda}(t)$ to follow the evolving simplex. This process mimics sort of a ``haze magnifying glass'' around $p$ of ``scaling ratio'' $1/t\,$. Supported on the parallel description of the centered linear simplices in $T_pM$ done in section \ref{linsimp}, one cares about the $t$-redundancy arising in this way to describe the $p$-centered Riemannian barycentric simplices $\hat\tau\,$ (see lemma \ref{simpletext} and theorem \ref{text-11}): {\it for a good choice of a fixed $\rho>0$,} every $p$-centered Riemannian barycentric $k$-simplex $\hat\tau\,$, of thickness $\geq$ some given $ t_{0,k}>0$ and diameter small enough, is described through its vertices $q_i$ in linking $\hat\tau$ to ${\bf v}\in {\mathbb F}_{\geq s_{0,k},\rho}^k(p)$ through ${\bf w}\!=\!(\exp_p^{-1}(q_0),\dots,\exp_p^{-1}(q_k))$ and ${\bf w}\!=\!t {\bf v}$ with $t\in]0,1]$ (see remark \ref{psi}),
where $s_{0,k}$ is small enough compared to $t_{0,k}$ (see proposition \ref{srop}, remark \ref{srth}, lemma \ref{1c}). So, $\hat\tau$ {\it is seen in a tilted way} (through the evolving $\lambda(t)$) {\it as locus of minimizers of 
$\psi_{\lambda,{\bf v},t}$}. As $t$ goes to $0\,$, one can cylindrify and pass to the limit as done in section \ref{magnglass}.
\end{rem}

\begin{rem} \label{barybary} As $t$ varies, $p$ remains the barycentre of the Riemannian barycentric
simplices $\hat\tau_{t{\bf v}}$ ({\it if they exist}) having vertices 
$p_i(t)\!=\!\exp_p(v_i(t))$: 
{\it indeed,} the barycentre of $\hat\tau_{t{\bf 
v}}$ has coordinates $\lambda_i(t)\!=\!\frac{1}{k+1}$, so, taking the $\cdot$ derivative of 
$\psi_{{\bf \lambda},{\bf v},t}(\cdot)$ in this case, the gradient of 
$\psi_{{\bf \lambda},{\bf v},t}$ defined in (\ref{fontionAR}) at the minimizing $p$ is (proposition \ref{rbs1} and its proof)
\begin{equation*}
\sum_{i=0}^k\frac{1}{k+1}\exp_p^{-1}p_i(t)=\sum_{i=0}^k\frac{t}{k+1}\,v_i=0_p\ .
\end{equation*}
\end{rem}

\begin{rem} 
The goal of this new ``magnifying glass'' is to take advantage of 
the identification of
$\frac{1}{t}\exp_p^{-1}\hat\tau_{t{\bf v}}$ with $\sigma_{\bf v}\subset T_pM$ in the limit, as $t$ goes to $0\,$, see section \ref{limit33}.
\end{rem}

For fixed $p\in M,\,{\bf w}\in{\mathcal F}(k,T_pM)$, barycentric coordinates ${\bf \mu}\in{\mathcal H}^k\,$,
{\it we want to put a hand on
the point minimizing} 
$$\psi_{{\bf \mu},{\bf w}}(\cdot)=\sum_{i=0}^k\mu_i\,d^2(\exp_p(w_i),\cdot),
$$ 
existence and uniqueness are to be studied later.

\begin{rem} \label{crx}
Given ${\bf w}\in{\mathcal F}(k,T_pM),\,\mu\in {\mathcal H}^k$ ($\sigma_{\bf w}$ is centered at $0_p$) 
 \begin{equation*}
\psi_{{\bf \mu},{\bf w}}(\cdot):=\sum_{i=0}^k\mu_i\,d^2(\exp_p(w_i),\cdot)
\end{equation*} 
can be written 
$\psi_{{\bf \lambda},{\bf v},t}(\cdot)$ for some 
$t\in\R_+\setminus\{0\}$ and 
${\bf v}\in{\mathcal F}(k,T_pM)\cap \nu^{-1}\{\alpha\}\,$ 
($\alpha =1/2\sqrt{k}\,$, convention \ref{convention3}). 
\par {\it Indeed,} in view of definition \ref{barcoord} (\ref{first}), (\ref{homothetie}), one has ($\sigma_{\bf v}$ centered at $0_p$)
\begin{equation*}
h_t:(\lambda,{\bf v})\longmapsto (h_t(\lambda),h_t({\bf v}))=(\lambda(t),t\,{\bf v})\,.
\end{equation*}
Taking $t=\nu({\bf w})/\alpha,{\bf v}={\bf w}/t$ and $\lambda$ solution of $\mu=h_t(\lambda)$, one has
\begin{equation*}\label{crx0}\psi_{{\bf \mu},{\bf w}}(\cdot)=\psi_{h_t({\bf \lambda},{\bf v})}(\cdot)=\psi_{\lambda(t),t{\bf v}}(\cdot)=\psi_{{\bf \lambda},{\bf v},t}(\cdot)\,.
\end{equation*}
\end{rem}
See a comment in appendix \ref{pseudohomothety}.
\vskip1mm

Of course, the question of the existence and 
uniqueness of the point where $\psi_{{\bf \lambda},{\bf v},t}$ achieves a minimum 
will occupy us for quite a few 
pages. 

Recall $W$ is a relatively compact open set 
$\subset (M,g)$ 
and, by lemma \ref{L.constante}, ${\mathcal R}_W$ is less than half of the convexity 
radius at all points in 
$W_{2{\mathcal R}_W}$.

\begin{rem} \label{psi} 
Given $p\in M$ and ${\bf v}\in{\mathcal F}(k,T_pM)$, 
if ${\bf v}$ is in $\nu^{-1}\{\rho\}$, lemma \ref{text2} implies 
$\sigma_{\bf v}\subset B(0_p,2\sqrt{k}\rho)$. Hence the barycentre $0_p$ 
is at distance $< 2\sqrt{k}\rho$ from any $v_i$ and $p_i(t)=\exp_p(t\,v_i)$ are points in  
$B(p,{\mathcal R}_W)$ if $t\in\,]0,\alpha\,{\mathcal R}_W/\rho]\,$. If 
$ {\bf \lambda}\,,{\bf v}\in{\mathcal H}^k\times{\mathbb F}_{>0,\rho}$ and  
$t\in\,]0,\alpha{\mathcal R}_W/\rho]$,  the function $\psi_{{\bf \lambda},{\bf v},t}$ is thus defined and {\rm differentiable} on the 
convex ball $B(p,{\mathcal R}_W)$ (recall $\alpha =1/2\sqrt{k}\,$, convention \ref{convention3}). Since the vertices are spread in a convex ball, theorem \ref{rbs0} and its corollary imply the existence of the Riemannian barycentric $k$-simplex $\hat\tau_{t{\bf 
v}}\,$ (and $p$ is the barycentre, see remark \ref{barybary}). So, if the barycentric coordinates are $\geq0\,$, the question of the existence and unicity of a minimizer of $\psi_{{\bf \lambda},{\bf v},t}$ has already be settled with the existence of $\hat\tau_{t{\bf v}}\,$ in section \ref{rbs}.

On the other hand, if $\sigma_{\bf v}$ is of diameter 
$\leq \rho\,$, then  $\nu(\sigma_{\bf v})\leq \rho\,$.
Indeed, if ${\rm diam}\,\sigma_{\bf v}\leq \rho$, one has 
$\max_i\Vert v_i-v_0\Vert\leq \rho\,$ and $\nu({\bf v})\leq \rho\,$. 
Thus, any Riemannian barycentric $k$-simplex of barycentre $p$ and diameter small enough is described along the lines just above.
\end{rem}

We now formalize the cylindrification of the family $\psi_{{\bf \lambda},{\bf v},t}\,$, taking advantage of the cylindrification done in section \ref{linsimp}.
\begin{prop} \label{text0} 
Given $\rho>0$, let $V$ be an open set in ${\mathcal H}^k\times {\mathbb F}_{>0,\rho}\,$.  
Denote by $\pi({\bf v})$ the point $p$ such that ${\bf v}$ belongs to
${\mathcal F}(k,T_pM)$. Consider the set
\begin{equation*}
\Delta=\{({\bf \lambda},{\bf v},q)\mid (\lambda,{\bf v})\in V
\   \hbox{\rm and} \  q\in  B(\pi({\bf v}),{\mathcal R}_W)\}\ .
\end{equation*}
The functions $\{\psi_t\mid t \in  \,]0,\alpha{\mathcal R}_W/\rho]\}$ 
defined through
\begin{equation*}
\psi_t : ( \lambda,{\bf v},q)\in \Delta\longmapsto 
\psi_{\lambda, {\bf v},t} (q) \in \R
\end{equation*}
constitute a differentiable family.
This family 
extends in a differentiable family over an open interval $I$ containing
$[0,\alpha\,{\mathcal R}_W/\rho]\,$,  
defining at $t=0$
\begin{equation} \label{text1}
\psi_0 : ( \lambda,{\bf v},q) \longmapsto 
\psi_{\lambda, {\bf v},0} (q) :=
\langle-2\,\exp_p^{-1}q,\sum_{i=0}^k\lambda_i\,v_i\rangle+d^2(p,q)\,  .
\end{equation}
\end{prop}

\begin{rem} ``Cylindrifying'' in the way described, one may observe that the set 
$V\times I$ becomes an open set in ${\mathcal H}^k\times\Gamma_{>0}\,$, see  
(\ref{Kylinder}), definition \ref{cylindrification}
and related remark. Notice also $\sum_{i=0}^k\lambda_i\,v_i=\sum_{i=0}^k\lambda_i(t)\,tv_i\,$.
\end{rem}

\begin{proof} 
As $t$ goes to $0\,$, the points $p_i(t)=\exp_p(tv_i)$ converge to the fixed 
center of gravity $p\,$.
In view of $\sum_{i=0}^k\lambda_i=\sum_{i=0}^k\frac{1}{k+1}=1\,$, expressing $\lambda_i(t)$ thanks to definition \ref{barcoord} (\ref{homothetie}), we restate 
$\psi_{{\bf \lambda},{\bf v},t}$ in the following way
\begin{equation*}
\psi_{{\bf \lambda},{\bf v},t}(\cdot)
=\sum_{i=0}^k(\lambda_i-\frac{1}{k+1})\,(\frac{d^2(p_i(t),\cdot)-
d^2(p,\cdot)}{t})+\sum_{i=0}^k \frac{d^2(p_i(t),\cdot)}{k+1}\  .
\end{equation*}
As $t\not=0$ tends to $0\,$, we know that 
\begin{equation*} 
\frac{d^2(p_i(t),q)-d^2(p,q)}{t} \hskip2mm  \hbox{tends to}\hskip2mm 
\langle \nabla_p\,d^2(\cdot_p,q),v_i  \rangle=-2\langle \exp_p^{-1}q,v_i  \rangle\  . 
\end{equation*}
To manage the question of the differentiability (see proposition \ref{mgapp1} or 
apply Taylor's formula with integral remainder), observe that, if 
$f(t)=d^2(p_i(t),\cdot)-d^2(p,\cdot)$ is $C^{r+1}$ and 
$f(0)=0\,$, then $g(t)=f(t)/t$ extended at $t=0$ by $g(0)=f'(0)$ is $C^r\,$.
\end{proof}
\begin{defn}\label{uniqmini0}
Given $p\in W\subset M$ and ${\bf v}\in{\mathcal F}(k,T_pM)$, we know $\psi_{{\bf \lambda},{\bf v},0}\,$ is defined in a convex ball $B(p,{\mathcal R}_W)\,$ for $\lambda\in {\mathcal H}^k\,$ (proposition \ref{text0}, remark \ref{psi}). {\it If the minimum of $\psi_{{\bf \lambda},{\bf v},0}\,$ exists and is unique,
define} the map ${\mathcal Q}_{{\bf v},0}\,$ sending $\lambda\in {\mathcal H}^k\,$ 
to this unique minimum ${\mathcal Q}_{{\bf v},0}(\lambda)\in B(p,{\mathcal R}_W)\,$.
\end{defn}
\begin{lem}\label{uniqmini0-bis}  {\rm On its set of definition}, i. e. for any $p\in W\subset M\,$, ${\bf v}\in{\mathcal F}(k,T_pM)$ and $\lambda\in{\mathcal H}^k$ such that $\psi_{{\bf \lambda},{\bf v},0}\,$ achieves a unique minimum (in the convex ball $B(p,{\mathcal R}_W)$), the map ${\mathcal Q}_{{\bf v},0}\,$ satisfies
\begin{equation} \label{identity}
{\mathcal Q}_{{\bf v},0}({\bf \lambda})=\exp_p(\sum_{i=0}^k\lambda_i\,v_i) \, .
\end{equation}
\end{lem}
\begin{proof} 
The proof holds for $\lambda\in{\mathcal H}^k$ such that
a unique minimum of $\psi_{{\bf \lambda},{\bf v},0}$ exists.
The minimum point $q$ is critical for $\psi_{{\bf \lambda},{\bf v},0}\,$,
thus one has 
\begin{equation}\label{text5}
\forall h \in T_qM \hskip8mm 0=\sum_{i=0}^k\lambda_i\langle 
d(\exp_p^{-1}\cdot_q)(q)(h),v_i\rangle_p+\langle\exp_q^{-1}p, h\rangle_p\  . 
\end{equation}
We deduce from (\ref{text5})
\begin{equation*}
\exp_q^{-1}p=-\sum_{i=0}^k\lambda_i \,{}^t(d(\exp_p^{-1}\cdot)_q(q))(v_i)\, , 
\end{equation*}
or 
\begin{equation} \label{transpose} 
{}^t(d(\exp_p\cdot))(\exp_p^{-1}q)(\exp_q^{-1}p)=
-\sum_{i=0}^k\lambda_i \,v_i\  .
\end{equation}
On one hand we know that
\begin{multline}\label{forgotten-25}
\langle {}^t(d(\exp_p\cdot))(\exp_p^{-1}q)(\exp_q^{-1}p), \exp_p^{-1}q\rangle = \\ 
=\langle \exp_q^{-1}p, d(\exp_p\cdot)(\exp_p^{-1}q)\exp_p^{-1}q\rangle=
-\langle \exp_q^{-1}p, \exp_q^{-1}p\rangle = \\
=-\langle \exp_p^{-1}q, \exp_p^{-1}q\rangle\  .
\end{multline}
On the other hand, if $w$ verifies $\langle w,\exp_p^{-1}q\rangle=0\,$, one has
\begin{multline}\label{orthogonal}
\langle {}^t(d(\exp_p\cdot))(\exp_p^{-1}q)(\exp_q^{-1}p), w\rangle = \\
= \langle \exp_q^{-1}p, d(\exp_p\cdot)(\exp_p^{-1}q)(w)\rangle=0\  ,
\end{multline}
since $d\exp$ preserves orthogonality to the geodesic rays from $p\,$, by Gauss' lemma. 

Hence (\ref{orthogonal}) shows that the vectors 
$-{}^t(d(\exp_p\cdot))(\exp_p^{-1}q)(\exp_q^{-1}p)$ and 
$\exp_p^{-1}q$ are proportional, and we conclude with 
(\ref{forgotten-25}) they are  equal.
Using (\ref{transpose}) we get $\exp_p^{-1}q=\sum_{i=0}^k\lambda_i \,v_i\,$, 
implying the claim.
\end{proof}

Now a  
tool used to derive the existence and uniqueness of minimizers for an open set of functions 
$\psi_{{\bf \lambda},{\bf v},t}\,$ neighboring the limit function $\psi_{{\bf \lambda},{\bf v},0}\,$.

\begin{defn} 
Let $K$ be a compact set in a manifold ${\mathcal M}\,$, let 
${\mathcal O}\subset\R^l$ 
be an open set containing $0$ and $(M,g)$ a Riemannian manifold. 
Denote by ${\mathcal B}\subset {\mathcal M}\times M$ a relatively compact manifold 
which is fibered over the compact $K$ with fibres $B_\mu\,$, $\mu \in K$, each $B_\mu$ being an 
open convex and relatively 
compact geodesic ball contained in $(M,g)$.
\end{defn}

\begin{prop} \label{strategy} 
With the above notations, define ${\mathcal D}$ to be the set
\begin{equation*}
{\mathcal D}=\{(\mu,\nu,q)\mid (\mu,q)\in \bar{{\mathcal B}} \hskip2mm 
\hbox{\rm and} 
\hskip2mm \nu\in {\mathcal O}   \}\  .
\end{equation*}
Let $\varphi$ be a $C^2$ mapping
\begin{equation*} 
\varphi  : (\mu,\nu,q)\in {\mathcal D}  \longmapsto    \varphi(\mu,\nu,q)\in\R \  . 
\end{equation*}
Denote by $B'_\mu$ $($by $B''_\mu)$ the open ball in $(M,g)$ having same center as 
$B_\mu$ and 
radius divided by two $($radius equal to three-fourth of the original radius$)$.
Assume that for any $\mu\in K$, the following two properties hold

$(a) \hskip1cm \forall q\in \bar B_\mu,\  
{}^{\sharp}\hbox{\rm hess}_q\,\varphi_{\mu,0} >>0\  ,$

$(b)$ \hskip2mm the function $\varphi_{\mu,0}$ achieves a minimum in $B_\mu$ ({\em unique}, in view of {\em (a)}) 
which actually belongs to $\bar B'_\mu\,$.

Then, there exists a relatively compact open neighborhood 
${\mathcal V}$ 
of $0$ in ${\mathcal O}$ such that,  for any $\mu\in K$ and any $\nu\in{\mathcal V}\,$, 
one also has

$(i) \hskip1cm \forall q\in \bar B_\mu, \ 
 {}^{\sharp} \hbox{\rm hess}_q\,\varphi_{\mu,\nu} >> 0\  ,$

$(ii)$ \hskip2mm the function $\varphi_{\mu,\nu}$ achieves a minimum in $B_\mu$ ({\em unique}, in view of {\em (i)}) 
which actually belongs to $ B''_\mu\,$.
\end{prop}

\begin{proof} 
Consider on the space of $g$-symmetric $(1,1)$-tensors $h$ on $TM$ the operator 
norm $\Vert\cdot\Vert\,$. From $(a)$ and the $C^2$-regularity of the family 
$\varphi_{\mu,\nu}$ (as $\nu$ 
tends to $0\,$, the function $\varphi_{\mu,\nu}$ tends to $\varphi_{\mu,0}$ in the $C^2$ 
sense on $\bar B_\mu$), there exists a $C>0$ and a relatively compact open neighborhood 
${\mathcal V}_1$ of $0$ 
in ${\mathcal O}$ such that, for any $\mu\in K, q\in \bar B_\mu,\nu\in{\mathcal V}_1\,$, one has
\begin{equation*}
{}^{\sharp}\hbox{\rm hess}_q\,\varphi_{\mu,0}\geq C\,\hbox{\rm Id}_q \hskip6mm 
\hbox{\rm and} \hskip6mm 
\Vert {}^{\sharp}\hbox{\rm hess}_q\,\varphi_{\mu,\nu}-{}^{\sharp}
\hbox{\rm hess}_q\,\varphi_{\mu,0}\Vert\leq\frac{C}{2}\  .
\end{equation*}
So, one has for any $\mu\in K,q\in\bar B_\mu$
\begin{equation*}
{}^{\sharp}\hbox{\rm hess}_q\,\varphi_{\mu,\nu}\geq\frac{C}{2}\,\hbox{\rm Id}_q\  ,
\end{equation*}
which proves $(i)$.

 Define 
\begin{equation*} 
a= \inf_{\mu\in K}(\inf_{\bar B_\mu \setminus 
B''_\mu}(\varphi_{\mu,0})-\inf_{\bar B'_\mu}(\varphi_{\mu,0}))\  .
\end{equation*}
From $(b)$ we have $a>0\,$. The $C^0$-proximity of $\varphi_{\mu,\nu}$ and 
$\varphi_{\mu,0}$ implies the existence of an open neighborhood ${\mathcal V}$ 
of $0$ included in ${\mathcal V}_1$ such that for any $\nu \in {\mathcal V}$ we 
have
\begin{equation*}
\sup_{\mu \in K}\,\sup_{q\in \bar B_\mu} \vert \varphi_{\mu,\nu}(q)- 
\varphi_{\mu,0}(q)\vert \leq \frac{a}{3}\ .
\end{equation*}

Thus, if $\varphi_{\mu,\nu}$ achieves a minimum at  
$q_1\in \bar B_\mu \setminus B''_\mu\,$ (this minimum is unique since $\varphi_{\mu,\nu}$ verifiying (i) is strictly convex), one can write
\begin{equation*}
\inf_{\bar B_\mu \setminus B''_\mu}\varphi_{\mu,\nu}=\varphi_{\mu,\nu}(q_1)\geq 
\varphi_{\mu,0}(q_1)-\frac{a}{3}\geq
\inf_{\bar B_\mu \setminus B''_\mu}\varphi_{\mu,0}-\frac{a}{3}\  ;
\end{equation*}
if $\varphi_{\mu,0}$ (verifying (i)) achieves its minimum at $q_2\in \bar B'_\mu\,$, 
one also has
\begin{equation*}
\inf_{\bar B'_\mu}\varphi_{\mu,0}=\varphi_{\mu,0}(q_2)\geq 
\varphi_{\mu,\nu}(q_2)-\frac{a}{3}\geq
\inf_{\bar B_\mu}\varphi_{\mu,\nu}-\frac{a}{3}\  .
\end{equation*}
Putting together the above inequalities gives
\begin{equation*}
0 = \inf_{\bar B_\mu \setminus B''_\mu}\varphi_{\mu,\nu}-
\inf_{\bar B_\mu}\varphi_{\mu,\nu}\geq
\inf_{\bar B_\mu \setminus B''_\mu}\varphi_{\mu,0}-
\inf_{\bar B'_\mu}\varphi_{\mu,0}-\frac{2a}{3}\geq\frac{a}{3}\  ,
\end{equation*}
a contradiction. This proves $(ii)$.
\end{proof}
\subsubsection{The limit $\psi_0$ fits to apply proposition {\rm\ref{strategy}}}
\label{limit00}
\begin{defn} 
For any $x\in W \subset (M,g)$, denoting by $B_x$  the 
open, relatively compact and convex ball $B(x,{\mathcal R}_W)$,
for any $p\in B_x\,$, the vector valued map $\exp_p^{-1}:B_x\rightarrow T_pM$
can be seen as $\exp_p^{-1}(\cdot)=\sum_i \langle \exp_p^{-1}(\cdot),a_i\rangle_p \, a_i\,$, 
where $a_i$ is an orthonormal basis in 
$(T_pM,g_p)$.
The following definition makes sense for any $q\in B_x$
\begin{equation*} {}^{\sharp}\hbox{\rm hess}_q
\exp_p^{-1} (\cdot_q)=\sum_i \langle {}^{\sharp}\hbox{\rm hess}_q\exp_p^{-1}(\cdot_q),a_i\rangle_p \, a_i\  .
\end{equation*}
Set 
\begin{equation*} \Vert{}^{\sharp}\hbox{\rm hess}_q\exp_p^{-1} (\cdot_q)\Vert^2=
\sum_i\Vert\langle {}^{\sharp}\hbox{\rm hess}_q\exp_p^{-1} (\cdot_q),a_i\rangle_p\Vert_{q,q}^2
\end{equation*}
and
\begin{equation*}
{\mathscr H}^g=
\sup_{x\in W}\sup_{p,q\in\bar B_x}
\Vert{}^{\sharp}\hbox{\rm hess}_q\exp_p^{-1} (\cdot_q)\Vert\  .
\end{equation*}
\end{defn}

\begin{lem}\label{text}
There exists $\rho_1\in\,]0,{\mathcal R}_W/3]$ 
such that, for any $x\in W\,$,
any $p\in B(x,\rho_1)$ and $q\in B(x,3\rho_1)$, and for all 
$({\bf \lambda},{\bf v})\in{\mathcal H}^k\times{\mathcal F}(k,T_pM)$ verifying 
$\Vert \sum_{i=0}^k\lambda_i\,v_i\Vert_p\leq 1/(4\,{\mathscr H}^g)$, one has 
\begin{equation*}
{}^{\sharp}\hbox{\rm hess}_q(\psi_{{\bf \lambda},{\bf v},0})\geq 
\frac{1}{2}\,\hbox{\rm Id}_q>>0\ .
\end{equation*}
\end{lem}
\begin{rem} \label{text-33}
The hypothesis $\Vert \sum_{i=0}^k\lambda_i\,v_i\Vert_p\leq 1/(4\,{\mathscr H}^g)$ is 
satisfied for all 
$({\bf \lambda},{\bf v})\in([0,1]^{k+1}\cap {\mathcal H}^k)\times 
{\mathbb F}_{\geq s_0,\frac{\alpha}{4{\mathscr H}^g}}$
(see lemma \ref{text2}).
\end{rem}

\begin{proof} 
One knows there exists $\rho_1\in\,]0,{\mathcal R}_W/3]$ such that, for any 
$p\in B(x,\rho_1)$, one has for any 
$q\in B(p,2\rho_1)\subset B(x,3\rho_1)\subset B(x,{\mathcal R}_W)$
\begin{equation*}
{}^{\sharp}\hbox{\rm hess}_q(d^2(p,\cdot))\geq \hbox{\rm Id}_q>>0\ .
\end{equation*}
As the inequality below also holds
\begin{equation*}
\Vert{}^{\sharp}\hbox{hess}_q\langle-2\exp_p^{-1}(\cdot_q),
\sum_{i=0}^k\lambda_i\,v_i\rangle_p\Vert_p
\leq 2\,{\mathscr H}^g \Vert \sum_{i=0}^k\lambda_i\,v_i\Vert_p\  ,
\end{equation*}
one gets the result using (\ref{text1}).
\end{proof}

\begin{defn} For any $\rho>0\,$,
 $p\in M$ and ${\bf v}\in{\mathbb F}_{>0}(p)$, define the relatively compact open set
\begin{equation} \label{U_v0}
{\mathcal U}_{\bf v,\rho}:= \{{\bf \lambda} \in{\mathcal H}^k\mid 
\Vert\sum_{i=0}^k\lambda_i\,v_i\Vert_p< \rho\}\subset{\mathcal H}^k\ .
\end{equation}
\end{defn}

\begin{lem} \label{text6} 
Given $x\in W$ and $s_0>0$, choose $\rho>0$ such that $\rho\leq \rho_2:=\min(\rho_1,1/(4\,{\mathscr H}^g))\,$ (where $\rho_1$ is defined in lemma {\rm\ref{text}}).
For any $p\in B(x,\rho),\,{\bf v}\in{\mathbb F}_{\geq s_0,\alpha\rho}(p)$, the properties below are true. 
\par $(a)$ One has the inclusion $([0,1]^{k+1}\cap{\mathcal H}^k)\subset{\mathcal U}_{\bf v,\rho}\ .$

\par $(b)$ The two properties below are satisfied for any $\lambda\in {\mathcal U}_{\bf v,\rho}\,$:

$\hskip1cm(i) \ \ \forall q \in \overline{B(p,2\rho)} \hskip1.5cm {}^{\sharp}
\hbox{\rm hess}_q(\psi_{{\bf \lambda},{\bf v},0}) \geq 
\frac{1}{2}\,\hbox{\rm Id}_q>>0 \ ;$
\par
$\hskip1cm (ii) \ \  \psi_{{\bf \lambda},{\bf v},0}$ achieves a minimum at a 
point 
${\mathcal Q}_{{\bf v},0}({\bf \lambda})$ belonging to $B(p,\rho)$ which is unique.

\par $(c)$
The map  
${\mathcal Q}_{{\bf v},0}$ is diffeomorphic from ${\mathcal U}_{\bf v,\rho}$ onto $B(p,\rho)\cap \exp_p(V_{\bf v})$, where $V_{\bf v}$ is the vector subspace generated by ${\bf v}$.

\par $(d)$
For any ${\bf v},\lambda$ such that $\psi_{{\bf \lambda},{\bf v},0}$ makes sense and for any $t\in]0,1]$, one has
the following crux invariance  
\begin{equation} \label{crx3}\psi_{{\bf \lambda},{\bf v},0}(\cdot)=\psi_{\lambda(t),t{\bf v},0}(\cdot)\  .
\end{equation}
\par One has for any $t\in]0,1]$, with $\lambda(t)=h_t(\lambda)$  (definition {\rm\ref{barcoord} (\ref{homothetie})})
\begin{gather}\label{crx1} {\mathcal U}_{t{\bf v},t{\rho}}={\mathcal U}_{\bf v,\rho}\ , \ \ \ {\mathcal Q}_{t{\bf v},0}(\lambda(t))\!=\!{\mathcal Q}_{{\bf v},0}({\bf \lambda})\ ,\\
\label{crx11}\,B(p,\rho)\cap \exp_p(V_{\bf v})={\rm rge}({\mathcal Q}_{{\bf v},0})={\rm rge}({\mathcal Q}_{t{\bf v},0})\,.
\end{gather}

\par $(e)$ The map ${\mathcal Q}_{t{\bf v},0}\,$ is defined on ${\mathcal U}_{{\bf v},\rho/t}\,$, ``$1/t$-dilated'' from the initial ${\mathcal U}_{\bf v,\rho}={\mathcal U}_{t{\bf v},t{\rho}}$ and properties (b) and (c) still hold for ${\mathcal Q}_{t{\bf v},0}\,$ while replacing ${\mathcal U}_{\bf v,\rho}\,, B(p,\rho)$ and $B(p,2\rho)$ occurring in their statements by the ``$1/t$-dilated'' ${\mathcal U}_{{\bf v},\rho/t}\,,B(p,\rho/t)$ and $B(p,2\rho/t)$.

\par $(f)$  One has, for any $\rho\in]0,\rho_2]\,,\, p\in B(x,\rho),\,{\bf v}\in{\mathbb F}_{\geq s_0,\alpha\rho}(p)$
\begin{equation}\label{range11}B(p,\rho_2)\cap \exp_p(V_{\bf v})\subset {\rm rge}({\mathcal Q}_{{\bf v},0})\ .
\end{equation}
\end{lem}

\begin{proof}
Given ${\bf v}\in {\mathbb F}_{\geq s_0,\alpha\rho}\,,\,{\bf \lambda}\in[0,1]^{k+1}\cap {\mathcal H}^k\,$, one has (see lemma \ref{text2})
\begin{equation*}
\Vert\sum_{i=0}^k\lambda_i\,v_i\Vert_p\leq \sup_{i=0}^k\Vert 
v_i\Vert_p<\rho\  ,
\end{equation*}
and (a) follows.

To prove (b), as $\rho\leq \rho_1\,$, one gets from (\ref{text1})
 $\psi_{0,{\bf v},0}=d^2(p,\cdot) 
=4\,\rho^2$
on $\partial B(p,2\rho)$. Since $\Vert\sum_{i=0}^k\lambda_i\,v_i\Vert_p\leq \rho\leq\frac{1}{4{\mathscr H}^g}$ one also gets, according 
to lemma \ref{text},  
 ${}^{\sharp}\hbox{\rm hess}_q(\psi_{{\bf \lambda},{\bf v},0}) \geq 
\frac{1}{2}\,\hbox{\rm Id}_q>>0 $ on $B(p,2\rho)$.
But one also has on $B(p,2\rho)$ (see (\ref{text1})) 
\begin{equation*}
\vert\psi_{{\bf \lambda},{\bf v},0}-d^2(p,\cdot)\vert\leq 
4\,\rho\,\Vert\sum_{i=0}^k\lambda_i\,v_i\Vert_p\  .
\end{equation*}
As $\Vert\sum_{i=0}^k\lambda_i\,v_i\Vert_p<\rho$ and since 
$\psi_{{\bf \lambda},{\bf v},0} >0$ on $\partial B(p,2\rho)$ and 
$\psi_{{\bf \lambda},{\bf v},0}(p)=0$, the function $\psi_{{\bf \lambda},{\bf v},0}$ achieves a minimum in $B(p,2\rho)$.
With assumption $\rho\leq \min(\rho_1,1/(4{\mathscr H}^g))$ (remark \ref{text-33}) and lemma \ref{text}, one gets 

$ (i) \hskip6mm {}^{\sharp}\hbox{\rm hess}_q(\psi_{{\bf \lambda},{\bf v},0})
\geq\frac{1}{2}\,\hbox{\rm Id}_q>>0 \ $ on $\overline{B(p,2\rho)}$;

$(ii)$ the minimum of $\psi_{{\bf \lambda},{\bf v},0}$ is achieved at a 
unique point in $B(p,2\rho)$. 

So, {\it any $\psi_{{\bf \lambda},{\bf v},0}$ with 
${\bf v}\in {\mathbb F}_{\geq s_0,\alpha\rho}\,$ and 
${\bf \lambda}\in {\mathcal U}_{\bf v,\rho}$ determines a unique minimum 
${\mathcal Q}_{{\bf v},0}({\bf \lambda})$ in $ B(p,2\rho)$, thus ${\mathcal Q}_{{\bf v},0}$ is defined on each ${\bf \lambda}\in {\mathcal U}_{\bf v,\rho}$.} That ${\mathcal Q}_{{\bf v},0}({\bf \lambda})$ even belongs to $B(p,\rho)$ results from (c).

The proof of (c) follows from  
(\ref{identity}) (see lemma \ref{uniqmini0-bis}). Indeed, in view of 
(\ref{U_v0}), one has ${\mathcal Q}_{{\bf v},0}({\mathcal U}_{\bf v,\rho})=\exp_p({\mathcal U}_{\bf v,\rho})=B(p,\rho)\cap \exp_p(V_{\bf v})\,$.

Then, (d) follows from (\ref{text1}), which implies, for any ${\bf v}$ and $\lambda$ such that $\psi_{{\bf \lambda},{\bf v},0}$ makes sense  
\begin{equation*} \psi_{{\bf \lambda},{\bf v},0}(\cdot)=\psi_{\lambda(t),t{\bf v},0}(\cdot)\ \ \hbox{in view of} \ \ \sum_i\lambda_i\,v_i=\sum_i\lambda_i(t)\,t\,v_i\  .
\end{equation*}
Indeed,  (\ref{crx11}) and the last equality in (\ref{crx1}) follow from this invariance (\ref{crx3}).
Moreover, one computes for any $\rho'>0$
\begin{multline} \label{crx2}{\mathcal U}_{t{\bf v},\rho'}\!=\!
\{\mu \mid \Vert\sum_i\mu_i \,t\,v_i\Vert< \rho'\}\!=\!\{h_t(\lambda) \mid \Vert\sum_i\lambda_i(t) \,t\,v_i\Vert< \rho'\}\\=
h_t(\{\lambda \mid \Vert\sum_i\lambda_i\,v_i\Vert< \rho'\})=h_t({\mathcal U}_{{\bf v}, \rho'})={\mathcal U}_{{\bf v}, \frac{\rho'}{t}}\  ,
\end{multline}
using the definition of $h_t$ given in definition \ref{barcoord} (\ref{homothetie}) and the first equality in (\ref{crx1}) follows making $\rho'=t\rho\,$.

Given 
${\bf \lambda}\in{\mathcal U}_{{\bf v}, \rho}\,$, (b) says ${\mathcal Q}_{{\bf v},0}({\bf \lambda})$ is defined and, in view of (d) (\ref{crx1}), ${\mathcal Q}_{t{\bf v},0}(h_t({\bf \lambda}))$ is also defined.
The first statement of (e) follows from $h_t({\mathcal U}_{{\bf v}, \rho})={\mathcal U}_{{\bf v}, \rho/t}\,$, while (\ref{crx3}) implies the other. 

All previous statements hold for any $\rho\in]0,\rho_2]\,$, implying (f).
\end{proof}
\begin{lem}\label{text-66} Assumptions of lemma {\rm\ref{text6}}; ${\mathcal Q}_{{\bf v},0}$ (${\mathcal Q}_{t{\bf v},0}$ if $t\in]0,1]$) is diffeomorphic from $\overline{\mathcal U}_{\bf v,\rho_2}$ (from $\overline{\mathcal U}_{{\bf v}, \rho_2/t}$) onto
$\overline{B(p,\rho_2)}\cap\exp_p(V_{\bf v})$.
\end{lem}
\begin{proof} Taking $\rho=\rho_2\,$, if $\lambda$ belongs to $\partial\,{\mathcal U}_{\bf v,\rho_2}\,$, by relaxing strict inequalities to large in the above proof, one gets $\psi_{{\bf \lambda},{\bf v},0}$ achieves a minimum inside $\overline{B(p,\rho_2)}$, which, in view of (b) (i), implies that it is unique. Thus ${\mathcal Q}_{{\bf v},0}$ is also defined on $\partial\,{\mathcal U}_{\bf v,\rho_2}$ and the rest follows as above.
\end{proof}
\subsubsection{Relying on the limit $\psi_0$}\label{limit33}
\begin{lem} \label{text7}
Define $\rho_3=\rho_2/2\,$ (see lemma {\rm \ref{text6}} for $\rho_2$). Given $x\in W\,$, $s_0>0\,$, denote by
${\mathcal K}_{\rho_3}$ the set of $({\bf \lambda},{\bf v})$, where 
${\bf \lambda}$ varies in ${\mathcal U}_{{\bf v},2\rho_3}$ as ${\bf v}$ varies in 
${\mathbb F}_{\geq s_0,\alpha\rho_3}(p)$ and $p$ in $\overline{B(x,\rho_3)}\,$. There  
exists $\eta>0$ such that, for any 
$\,t\,$ in $[-\eta,\eta]\,$ and any $({\bf \lambda},{\bf v})\in {\mathcal 
K}_{\rho_3}\,$, 
 all functions $\psi_{{\bf \lambda},{\bf v},t}\,$ satisfy

$ (i) \hskip1mm \hskip3mm   \forall q \in \overline{B(p,4\rho_3)} \hskip1cm 
{}^{\sharp}\hbox{\rm hess}_q\,\psi_{{\bf \lambda},{\bf v},t}>>0\  ;$

$(ii)$ \hskip1mm $\psi_{{\bf \lambda},{\bf v},t}$ achieves a minimum at a 
unique point ${\mathcal Q}_{{\bf v},t}(\lambda)\!\in B(p,3\rho_3)$.
\end{lem}
\begin{proof} 
From lemma \ref{text6} (b), the choice of $\rho_3\,$ implies for
$({\bf \lambda},{\bf v})\in \overline{\mathcal K}_{\rho_3}$:

 $ (i) \hskip6mm {}^{\sharp}\hbox{\rm hess}_q(\psi_{{\bf \lambda},{\bf v},0}) 
\geq \frac{1}{2}\,\hbox{\rm Id}_q>>0 \ $ on $\overline{B(p,4\rho_3)}\,$;

 $(ii)$ \hskip2mm $\psi_{{\bf \lambda},{\bf v},0}$ achieves a minimum 
at a unique point in $\overline{B(p,2\rho_3)}\,$.

In the frame of proposition \ref{strategy}, make
$\mu=({\bf \lambda},{\bf v})$ and $B_\mu=B(p,4\rho_3)$
\begin{equation*}
K=\overline{\mathcal K}_{\rho_3}=\{({\bf \lambda},{\bf v})\mid {\bf \lambda}\in  
\overline{{\mathcal U}_{{\bf v},2\rho_3}} \hskip2mm \hbox{and} \hskip2mm {\bf v}\in 
{\mathbb F}_{\geq s_0,\alpha\rho_3}(p)\ 
 \hbox{and} \hskip2mm p\in \overline{B(x,\rho_3)}\}\  ,
\end{equation*}
set $\nu=t$ and $\varphi_{\mu,\nu}=\psi_{{\bf \lambda},{\bf v},t}\,$, thus 
$\varphi_{\mu,0}=\psi_{{\bf \lambda},{\bf v},0}\,$. 
The first lines in this proof precisely say: {\it the hypotheses of 
proposition {\rm \ref{strategy}} are fulfilled}. Hence, this proposition \ref{strategy} tells there exists an $\eta>0$ such that the 
claim of lemma \ref{text7} holds. 
\end{proof}

\begin{lem} \label{text8} Given $x\in W,s_0>0\,$, choose $\rho_3$ and $\eta>0$ 
according to lemma {\rm \ref{text7}}. There exists $\gamma\in\,]0,\eta]$ 
such that the mappings 
\begin{equation*}
 \left\{  {\mathcal Q}_{{\bf v},t} \ \mid\  t\in 
[-\gamma,\gamma]\,,\ 
 {\bf v} \in {\mathbb F}_{\geq s_0,\alpha\rho_3}(p),\ 
 p\in \overline{B(x,\rho_3)}  \right\}
\end{equation*}
constitute a 
{\rm differentiable 
family of embeddings} of $\overline{\mathcal U}_{\bf v,\frac{3}{2}\rho_3}$ into $M$, where 
we have defined ${\mathcal Q}_{{\bf v},t}(\lambda):={\mathcal Q}_{\lambda, {\bf v}} (t)$, 
see {\rm (\ref{fontionAR})} for ${\mathcal Q}_{\lambda, {\bf v}} (t)$.
\end{lem}

\begin{proof} 
Set ${\mathcal Q}({\bf v},t,\lambda):= {\mathcal Q}_{\lambda, {\bf v}} (t)={\mathcal Q}_{{\bf v},t} (\lambda)$.
We know from lemma \ref{text7} that, for any ${\bf v}\in{\mathbb F}_{\geq s_0,\alpha\rho_3}$ and $t\in[-\eta,\eta]\,$, the mappings ${\mathcal Q}_{{\bf v},t}$ are well 
defined from ${\mathcal U}_{{\bf v},2\rho_3}$ into $B(p,3\rho_3)$. 

 More, 
${\mathcal Q}$ is defined on $\overline{\mathcal K}_{\rho_3}\,$ (see lemma \ref{text-66}), on which holds the implicit equation
\begin{equation*}
d\psi_{{\bf \lambda},{\bf v},t}({\mathcal Q}({\bf v},t,\lambda))=
d\psi_{{\bf \lambda},{\bf v},t}({\mathcal Q}_{{\bf v},t}(\lambda))=0\  .
\end{equation*}
We know further, from  lemma \ref{text7}, that one has on 
$\overline{B(p,4\rho_3)}$ 
\begin{equation*} 
 {}^{\sharp} \hbox{\rm hess}_q\,\psi_{{\bf \lambda},{\bf v},t}>>0 \ .
\end{equation*}
Thus the implicit function theorem can be applied and the mapping 
${\mathcal Q}: ({\bf v},t,\lambda)\rightarrow {\mathcal Q}({\bf v},t,\lambda)$ as well as
the family ${\mathcal Q}_{{\bf v},t}$ are indeed differentiable. 

For any 
${\bf v}\in{\mathbb F}_{\geq s_0,\alpha\rho_3}\,$, we know (from lemma \ref{text6}(f)(\ref{range11})) 
that the 
mapping ${\mathcal Q}_{{\bf v},0}$ 
embeds ${\mathcal U}_{{\bf v},2\rho_3}$ onto $B(p,2\rho_3)\cap\exp_p(V_{\bf v})$. 
Therefore, for any ${\bf v}\in {\mathbb F}_{\geq s_0,\alpha\rho_3}$ exists $\gamma_{\bf v}>0$ (see \cite{G-G} chapter III or \cite{Mu}, theorem 3.10) 
such that, in $C^\infty({\mathcal U}_{{\bf v},2\rho_3}, M)$, any $\gamma_{\bf v}$-close $C^1$-approximation
of ${\mathcal Q}_{{\bf v},0}\,$ is {\it a proper 
embedding} from $\overline{\mathcal U}_{{\bf v},\frac{3}{2}\rho_3}$ into $M\,$, neighbouring the proper embedding (see lemma \ref{text-66})
${\mathcal Q}_{{\bf v},0}: \overline{\mathcal U}_{{\bf v},\frac{3}{2}\rho_3}\rightarrow M$. 

Set $\underline\gamma=\inf_{{\bf v}}\gamma_{\bf v}\,$. We claim  
$\underline\gamma$ is {\it positive}.

If this was not true, one could find a sequence $t_i$ tending to $0$ and 
another  sequence ${\bf v}_i\in {\mathbb F}_{\geq s_0,\alpha\rho_3}\,$, such that 
the ${\mathcal Q}_{{\bf v}_i,t_i}: \overline{\mathcal U}_{{\bf v},\frac{3}{2}\rho_3}\rightarrow M$ are not embeddings. Since 
${\mathbb F}_{\geq s_0,\alpha\rho_3}$ is compact, one could 
extract a  
sequence $({\bf v}_i,t_i)$ 
tending to $({\bf v},0)$, where 
${\bf v} \in {\mathbb F}_{\geq s_0,\alpha\rho_3}\,$.
This is {\it a contradiction}, because   
${\mathcal Q}_{{\bf v}_,0}$ shares the property to 
have, in $C^1(\overline{\mathcal U}_{{\bf v},\frac{3}{2}\rho_3}, M)$ a $C^1$-neighborhood of embeddings.
As $\underline\gamma$ is $>0$ and 
$({\bf v},t)\mapsto\Vert {\mathcal Q}_{{\bf v},t}-{\mathcal Q}_{{\bf v},0} 
\Vert_{C^1}$ is $C^0$, there exists $\gamma>0$ 
such that
\begin{equation*}
\forall t\in [-\gamma,\gamma]\hskip1.5cm \sup_{{\bf v}}
\Vert {\mathcal Q}_{{\bf v},t}-{\mathcal Q}_{{\bf v},0} \Vert_{C^1}\leq 
\underline\gamma\  .\qedhere
\end{equation*}
\end{proof}

\begin{lem} \label{text9} 
 Given $x$ in $W$ and $s_0>0\,$, choose $\rho_3$ and $\gamma$
according to lemma {\rm \ref{text8}}. There exists a 
positive $\beta \leq \min(\gamma,1)$ such 
that, for any $t\in[-\beta,\beta]\,,\,p\in  \overline{B(x,\beta\rho_3)}$ and 
${\bf v} \in {\mathbb F}_{\geq s_0,\alpha\rho_3}(p)\,$, the 
differentiable mappings 
$[{\mathcal Q}]_{{\bf v}_,t}: =({\mathcal Q}_{{\bf v}_,t})_{\mid_{{\mathcal U}_{{\bf v},\frac{3}{2}\rho_3}}}$
have properly embedded transversal images to $\partial B(x,\rho_3)$, all diffeomorphic to $k$-embedded disks.
\end{lem}
\begin{proof} 
By lemma \ref{text8}, the images $\hbox{\rm rge}\,([{\mathcal Q}_{{\bf v},t}])$ are 
embedded if $\vert t\vert\leq \gamma\,$. From lemmas \ref{uniqmini0-bis},
\ref{text6} (c) and (e), for any ${\bf v} \in {\mathbb F}_{\geq s_0,\alpha\rho_3}(x)\,$, we know that the embedded 
$\hbox{\rm rge}\,([{\mathcal Q}_{{\bf v},0}])$ have 
transversal intersections with 
$\partial B(x,\rho_3)$. Thus there exists $\beta_0\in]0,1]$ such that, for any $p\in\overline{B(x,\beta_0\rho_3)}$ and 
${\bf v} \in {\mathbb F}_{\geq s_0,\alpha\rho_3}(p)\,$, the embedded 
$\hbox{\rm rge}\,([{\mathcal Q}_{{\bf v},0}])$ still cut 
$\partial B(x,\rho_3)$ {\it transversally}.
As $\,\cup_{p\in\overline{B(x,\beta_0\rho_3)}}{\mathbb F}_{\geq s_0,\alpha\rho_3}(p)\,$ is compact, the classical stability of transversal properties imply the existence of $\beta\in ]0,\min(\beta_0,\gamma)]\,$ such that for any $t\in[-\beta,\beta]\,,\,p\in  \overline{B(x,\beta_0\,\rho_3)}\,,\,{\bf v}\in{\mathbb F}_{\geq s_0,\alpha\rho_3}(p)\,$, 
the  embedded 
$\hbox{\rm rge}\,([{\mathcal Q}_{{\bf v},t}])$ cut $\partial B(x,\rho_3)$ transversally. 
So, as $t$ describes $[-\beta,\beta]$, in view of the implicit function theorem, all $[{\mathcal Q}_{{\bf v},t}]^{-1}(\partial B(x,\rho_3))\subset{\mathcal U}_{{\bf v},\frac{3}{2}\rho_3}$ are seen to be sets diffeomorphic to a $(k-1)$-dimensional sphere (case $t=0$), each of which bounds an open $[{\mathcal Q}_{{\bf v},t}]^{-1} B(x,\rho_3)$ diffeomorphic to a $k$-dimensional ball in ${\mathcal U}_{{\bf v},\frac{3}{2}\rho_3}\,$.
\end{proof}

\begin{rem}\label{bon-bon}
Given ${\bf v} \!\in \!{\mathbb F}_{\geq s_0,\alpha\rho_3}(x)$, define ${\mathcal U}^{\bf v}_t$ (see lemmas \ref{text8}, \ref{text9})
\begin{equation*} 
{\mathcal U}^{\bf v}_t=
[{\mathcal Q}_{{\bf v},t}]^{-1}{B(x,\rho_3)}\, .
\end{equation*}
\begin{enumerate}
 \item\label{bon-bon102} {\it The inverse images $\overline{\mathcal U}^{\bf v}_t\!=\![{\mathcal Q}_{{\bf v},t}]^{-1}\overline{B(x,\rho_3)}$ 
are},
by virtue of the last two lemmas, $k$-{\it dimensional 
differentiable sub-manifolds, all 
diffeomorphic to a standard $k$-dimensional ball}, and the ``leaves'' $[{\mathcal Q}_{{\bf v},t}](\overline{\mathcal U}^{\bf v}_t) $ are 
{\it properly embedded in} $\overline{B(x,\rho_3)}$.

\item \label{bon-bon333}Recall $\rho_3$ is uniform over $ W\,$. Take into account that all developments were done for $p$ 
in 
a ball centered at a point $x$ and write $\beta_x$ for $\beta\,$. Denote by $\rho_
x$ the product $\rho_x=\alpha\,\beta_x\,\rho_3\leq \rho_3\,$. 
\end{enumerate}
\end{rem}

\subsubsection{Riemannian barycentric textures}\label{RBT}

\begin{defn}\label{texture3333}
Set $L_{{\bf v},t}:=[{\mathcal Q}_{{\bf v},t}] ({\mathcal U}^{\bf v}_t)$ and define $X:=X_{x,\rho_x,s_0}^k$
(of dimension $\kappa\!+\!k\!:=\!(k\!+\!1)n\!+\!k\,$, in accordance with \cite{C-M}, 
section 7) to be
\begin{equation*} 
X_{x,\rho_x,s_0}^k:=
\{({\bf v},t)\times L_{{\bf v},t}\mid{\bf v}\in 
{\mathbb F}_{\geq s_0,\rho_x}^k(p),\,t\in[-1,1], \,p\in B(x,\rho_x)\}\ .
\end{equation*}
\end{defn}
We can state the theorem which was our scope in this section. 
\begin{thm}\label{texture333}
The sets $L_{{\bf v},t}$ are the {\rm leaves} of the foliation ${\mathcal L}$ 
which composes a {\rm texture} in the sense of section {\rm 7} of {\rm \cite{C-M}}, the 
total space is
$X=X_{x,\rho_x,s_0}^k$, the {\rm morphism} is $\ {\mathpzc p}:({\bf v},t,q)\longmapsto q$ 
and the {\rm basis} $B(x,\rho_3)$. This texture is {\rm clean}.
\end{thm}
\begin{note} A somehow sufficient for our purpose and simpler definition of a {\em field fresh texture} is given later, see definition \ref{Texturedefn} in section \ref{multiple1088}.
\end{note}

\begin{proof} 
Recall the cylinder of definition \ref{cylindrification} and 
the way we introduced it:
given ${\bf v}\!\in\! {\mathbb F}_{\geq s}^k\cap\nu^{-1}\{\alpha\}$ and $t\!>\!0\,$, 
any $t\,{\bf v}\!\in\! {\mathbb F}_{\geq s}^k$ gives rise to  
$({\bf v},t)\!\in\!({\mathbb F}_{\geq s}^k\cap\nu^{-1}\{\alpha\})\times
\R_+^\ast\subset \Gamma_{\geq s}^k\,$. For any $p\!\in\! B(x,\rho_x)$, one thus can view any element  
$({\bf v},t)\in{\mathbb F}_{\geq s_0,\alpha\,\rho_3}^k(p)\times [-\beta_x,\beta_x]$ as $({\bf v}',t')=(\beta_x\,{\bf v},t/\beta_x)\in{\mathbb F}_{\geq s_0,\rho_x}^k(p)\times [-1,1]\,$. Apply
lemmas \ref{text8}, \ref{text9} to the maps $[{\mathcal Q}_{{\bf v},t}]$ defined on 
${\mathcal U}^{\bf v}_t\,$. It follows that the {\it leaves} $\,L_{{\bf v},t}$ are properly embedded connected 
submanifolds of $B(x,\rho_3)$. According to section $7$ in \cite{C-M}, this proves that 
$X_{x,\rho_x,s_0}^k$ 
fulfills the conditions defining a {\it clean texture}.
\end{proof}

Given $k\,$, one has around any point $x$ in the compact set $\bar W$ a 
texture $X_{x,\rho_x,s_0}^k$ with basis $B(x,\rho_3)$. 
Cover $\bar W$ with balls $B(x,\rho_x/2)$ and extract a finite 
subcovering $B(x_i,\rho_{x_i}/2)$. 

\begin{defn} \label{texture33333} $ $
\begin{enumerate} 
 \item  Denote by $\underline r$ the positive number 
 $\underline r=\min_i \rho_{x_i}\,$.

\item Denote by ${\mathfrak S}_{t_0}(W)$ the set of all Riemannian 
barycentric $n$-simplices of thickness $\geq t_0$ embedded in the relatively compact open 
subset $W\,$.

\item Recall lemma \ref{1fo} and the definition of $t_{0,k}>0$ as positive lower 
bound 
(depending on $k$ and $t_0$) for the thickness of any $k$-face in a linear Euclidean 
$n$-simplex of thickness $\geq t_0>0\,$.
Fix $s_0=s_{0,k}>0$ ensuring that any $\sigma_{\bf v}$ of thickness 
$\geq t_{0,k}>0$ has standard ratio $\geq s_{0,k}>0\,$, see remarks \ref{srth}, \ref{newglasses}.
\end{enumerate}
\end{defn}

Keeping in mind remark \ref{psi}, we also need the following
\begin{lem} \label{simpletext} 
For any Riemannian barycentric $n$-simplex $\hat\tau$ of thickness $\geq t_0$ and diameter 
$\leq \underline{r}/4\,$, each $k$-face is a $k$-Riemannian barycentric simplex 
$\hat\tau_k\,$, whose barycentre we call $p\,$. This $\hat\tau_k$ gives rise to a unique 
${\bf v}\in{\mathbb F}_{\geq s_{0,k},\underline{r}}^k(p)$ such that $\hat\tau_k\subset L_{{\bf v},t}\,$ for some $t\in[-1,1]\,$.
\end{lem}
\begin{proof} 
Since $\hbox{diam}_g(\hat\tau)\leq \underline{r}/4\,$, the simplices $\hat\tau_k$ and $\hat\tau$ are included in $B(p,\underline{r}/4)$. As $p$ belongs to some $B(x_i,\rho_{x_i}/2)$, one also has (use $\underline{r}\leq \rho_3$)
$\hat\tau_k\subset B(x_i,3\rho_3/4)$. The thickness of $\hat\tau_k$ is $\geq t_{0,k}$ 
(definition \ref{genthick}, lemma \ref{1fo}),
its vertices $p_0,\dots,p_k$ determine 
${\bf w}\!=\!(\exp_p^{-1}p_0,\dots,\exp_p^{-1}p_k)$. One has 
$d(0_p,w_i)=\Vert w_i\Vert\leq \underline{r}/2\,$, thus $\max_i  \Vert w_i-w_0\Vert\leq \underline{r}$ and
one gets $\nu({\bf w})\leq \underline{r}\,$. A $1/t$-homothety (with $t\in]0,1]$) on ${\bf w}$ provides the 
wanted ${\bf v}={\bf w}/t\in{\mathbb F}_{\geq s_{0,k},\underline{r}}^k(p)\,$. That $\hat\tau_k$ is included in $L_{{\bf v},t}$ follows from $\hat\tau_k\subset\hbox{rge}({\mathcal Q}_{{\bf v},t})$ and $\hat\tau_k\subset B(x_i,\rho_3)$ (definition \ref{texture3333}, remark \ref{bon-bon}(\ref{bon-bon333})).
\end{proof}

Now, a well fit version of Thom's theorem related to our development. 
One may also read in parallel section \ref{multiple1088} of part \ref{annexe1033}.

In the sequel, we consider submanifolds $N\subset M$ which are compact $(n-k)$-dimensional 
and written in the form
$N=f^{-1}\{0\}\,$, where $f$ is a differentiable 
map from $M$ to $\R^k\,$, regular over $N\,$, fitting with the theorem
applied below, as it is stated in \cite{C-M}.

\begin{defn} Set ${\bf n}_A$ for the number (finite or infinite) of points in a set $A\,$.
Recall the definition (see \cite{C-M}) of the local $({\mathcal L},{\mathpzc p})$-degree at a point 
$a\in f^{-1}\{0\} \cap \bar W$ of $f^{-1}\{0\}$ with respect to the 
texture ${\mathcal L}$ 
with 
leaves $L$ (here ${\mathcal L}=X$ is foliated by $L$)
\begin{equation*} 
\hbox{deg}_{{\mathcal L},{\mathpzc p},a} f^{-1}\{0\}:=
\inf_{U,U_1} \sup_L\  {\bf n}_{\,f^{-1}\{0\}\cap U \cap {\mathpzc p}(L)}\ ,
\end{equation*}
where the infimum is over all open neighborhoods $U$ in $f^{-1}\{0\}$ and $U_1$ in $M$ of the same point $a$ and the supremum is over any leaf $L$ of the foliation 
${\mathcal L}$ induced in ${\mathpzc p}^{-1}(U_1)$.
\end{defn}
\begin{rem}\label{multipl1033} Actually, the degree of intersection $\hbox{deg}_{{\mathcal L},{\mathpzc p},a} f^{-1}\{0\}$ is bounded from above by another integer (see lemma \ref{degleqmult}),
the 
{\it multiplicity} $\mu_{{\mathcal L},{\mathpzc p},a}\,f^{-1}\{0\}\,$, to be defined 
later, see definition \ref{multiple1044} in 
section \ref{multiple1088} (see also section 2 of \cite{C-M}, the version used here is actually developed in lemma 7.2 of \cite{C-M}). 
The fact that this multiplicity is also upper semi-continuous under perturbations 
(see \cite{G-G} pp. 167-8-9, see lemma \ref{semicontmult}) is used in the proof of theorem \ref{text-11} below.
\end{rem}
For future use, we not only consider maps from $M$ to $\R^k$ defining 
$(n-k)$-dimensional submanifolds, but also differentiable families of such maps, i. e.  
maps $\bar f\in C^\infty({\mathscr T}\times M,\R^k)$, 
where ${\mathscr T}$ is a $c$-dimensional manifold. 

\begin{exmp}
Let $\Delta$ be a line in $\R^4\,$. Any one-parameter $C^\infty$-family of curves in $\R^4$ which
meets $\Delta$ may be $C^\infty$-perturbed in another one-parameter family of curves in $\R^4$  which doesn't meet $\Delta$ and is as $C^\infty$-close as wished from the initial one. This cannot be done for a two-parameter $C^\infty$-family of curves in $\R^4$ as a
family of curves constituted by lines all parallel to a given
line in a hyperplane hitting $\Delta\,$.
\end{exmp}

We shall need the following lemma
\begin{lem}\label{unlemme} Given an $n$-manifold $M$ and a $c$-manifold ${\mathscr T}$, call $B^{n}$ the unit ball in ${\mathbb R}^{n}$ and $\varphi$ a diffeomorphism from an open ${\mathcal O}\subset M$ onto an open of ${\mathbb R}^{n}$ that contains $\bar B^n$. Set $V=\varphi^{-1}B^{n}$. The restriction mapping
$$ {\mathcal R}_{\bar V}: \bar f\!\in\! C^\infty({\mathscr T}\times M,{\mathbb R}^l)\longmapsto
{\mathcal R}_{\bar V}(\bar f)\!=\!\bar f_{\mid{\mathscr T}\times\bar V}\!\in\! C^\infty({\mathscr T}\times\bar V,{\mathbb R}^l)
$$
is onto, open and continuous in the Whitney topology.
Thus, if ${\mathcal U}$ is open dense in $C^\infty({\mathscr T}\times M,{\mathbb R}^l)$, so is ${\mathcal R}_{\bar V}({\mathcal U})$ in $ C^\infty({\mathscr T}\times\bar V,{\mathbb R}^l)$. If ${\mathcal V}$ is open dense in $C^\infty({\mathscr T}\times\bar V,{\mathbb R}^l)$, so is $ {\mathcal R}_{\bar V}^{-1}({\mathcal V})$ in $C^\infty({\mathscr T}\times M,{\mathbb R}^l)$.
\end{lem}
\begin{rem}\label{unlemme444} A function on $\bar V$ is $C^k$ ($k$ may be $\infty$) if it has a $C^k$-extension in some open containing $\bar V$ (see \cite{D}, p. 18 ff).
\end{rem}
\begin{proof} This lemma is proved as lemma \ref{unlemme333} in part \ref{M-T2}.
\end{proof}
\begin{thm}[Thom revisited, local form] \label{text-11}
Given $t_0\!>\!0,\,c\!\in \!\N\,$, let $(M,g)$ be a Riemannian $n$-manifold and ${\mathscr T}$ a $c$-manifold.
One can find an open dense subset ${\mathcal V}_c$ of the open set in $C^\infty({\mathscr T}\times M,\R^k)$ of families $\bar f$
of regular proper maps $f_{\vartheta}$, where
$f_\vartheta(\cdot)\!:=\!\bar f(\vartheta,\cdot)$, for which, given a compact
$\Theta\subset{\mathscr T}\,$, there exists $\rho_{\bar f}\!>\!0$ $($depending also
on $\Theta)$
such that, for any $\rho\!\in\,]0,\rho_{\bar f}]\,$, if
$\vartheta$ varies in $\Theta\,$, the sets $f_\vartheta^{-1}\{0\}$ intersect all $k$-faces
of Riemannian barycentric simplices $\hat\tau\!\in\!{\mathfrak S}_{t_0}(W)$ of diameter
$\leq \!2\rho$ in at most $m_k(\kappa\!+\!c)$ points ($m_k(\kappa\!+\!c)$ is an integer depending only on $k$ and $\kappa+c$, see remark {\rm{\ref{codimension1066}}}, and here $\kappa\!=\!(k\!+\!1)n$, see definition {\rm\ref{texture3333}}).
\end{thm}
\begin{note} We later present and prove a more precise version of the above theorem, see theorem \ref{text-1133}. This local Thom theorem, which relies on theorem \ref{thm7.1C-Mbis}, is a byproduct of theorem 7.1 in \cite{C-M} and its proof.
\end{note}

\begin{proof}
Take any Riemannian barycentric $n$-simplex $\hat\tau$ of diameter $\leq 2\rho\,$, where $\rho$ is chosen such that
$4\rho\leq \underline r\leq \rho_{x_i}\leq \rho_3\,$. The barycentre $p$ of a $k$-face $\hat\tau_k\,$  is in some $B(x_i,\rho_{x_i}/2)$ and $\hat\tau_k\subset B(x_i,3\rho_3/4)$ (lemma
\ref{simpletext}), so (lemma \ref{1fo}, remarks \ref{srth}, \ref{newglasses})
there exist $i$, $t\in\,]0,1]\,,\,{\bf v}\in{\mathbb F}_{\geq s_{0,k},\rho_{x_i}}^k(p)$ such that
$\hat\tau_k$
is a subset of $L_{{\bf v},t}$, {\it a leaf of the texture}
$X_i^k=X_{x_i,\rho_{x_i},s_{0,k}}^k\,$.

Let $\bar f\!\in\!C^\infty({\mathscr T}\times M,\R^k)$ be a family of regular proper maps $f_\vartheta$.
From the definition of a texture and the nature of the map ${\mathpzc p}\,$, since the
leaves $L$ of $X_i^k$ are properly embedded
in $B(x_i,\rho_3)$ with the sets $f_\vartheta^{-1}\{0\}$ (for any $\vartheta\!\in\!{\mathscr T}$), calling ${
\mathcal N}_a$ for a neighborhood basis of $a$ in $M\,$, one may take as definition
of the local $(X_i^k,{\mathpzc p})$-degree at a point $a\in f_\vartheta^{-1}\{0\} \cap \bar W\cap B(x_i,\rho_3)$
\begin{equation*}
\hbox{deg}_{X_i^k,{\mathpzc p},a} f_\vartheta^{-1}\{0\}:=\inf_{U\in{\mathcal N}_a}
\sup_L  \,{\bf n}_{\,f_\vartheta^{-1}\{0\}\cap U \cap {\mathpzc p}(L)}\  .
\end{equation*}
Define the restriction maps
\begin{gather*}{\mathcal R}'_i:C^\infty({\mathscr T}\times M,\R^k)\rightarrow C^\infty({\mathscr T}\times\overline{B(x_i,3\rho_3/4)},\R^k)\\ {\mathcal R}''_i:C^\infty({\mathscr T}\times B(x_i,\rho_3),\R^k)\rightarrow C^\infty({\mathscr T}\times\overline{B(x_i,3\rho_3/4)},\R^k)\,.
\end{gather*}
Theorem \ref{thm7.1C-Mbis} says: for any $i$ exists an open and dense set ${\mathcal V}_{c,i}\subset C^\infty({\mathscr T}\times B(x_i,\rho_3),\R^k)$ of $\bar f$ such that, for any
$\vartheta\in{\mathcal T}$, the following bound $\hbox{deg}_{X_i^k,{\mathpzc p},a} f_\vartheta^{-1}\{0\}\leq m_k(\kappa+c)$ holds. Then ${\mathcal V}_c:=\cap_i{\mathcal R}_i'^{-1}({\mathcal R}_i''({\mathcal V}_{c,i}))$ is {\em open and dense} since {\em all ${\mathcal R}_i'^{-1}({\mathcal R}_i''({\mathcal V}_{c,i}))$ are so} (lemma \ref{unlemme}). As $\bar W\subset\cup_i \overline{B(x_i,3\rho_3/4)}$, the same theorem implies that for any
$\bar f\in {\mathcal V}_c\subset C^\infty({\mathcal T}\times M,\R^k)$, for any
$\vartheta\in{\mathcal T}\,$, all $(X_i^k,{\mathpzc p})$-degree of $f_\vartheta$ at a point
$a\in f_\vartheta^{-1}\{0\} \cap \bar W$ are
$\leq m_k(\kappa+c)$.

Actually, for any family $\bar f\in {\mathcal V}_c\,$, theorem \ref{thm7.1C-Mbis} implies that, for any $i$ and any parameter $\vartheta$, any $a\in f_\vartheta^{-1}\{0\} \cap \bar W$, the {\it multiplicity} verifies $\mu_{X_i^k,{\mathpzc p},a}\,f_\vartheta^{-1}\{0\}\leq m_k(\kappa+c)\,$. In view of its upper semi-continuity (see remark \ref{multipl1033}), the multiplicity is even seen to satisfy this bound
on an open set
$\theta\times U$ around each $\vartheta\times a$ in ${\mathscr T}\times M$,
and thus, for any $\vartheta'\in \theta$, one has
$\sup_L \,{\bf n}_{\, L\,\cap\, U \,\cap\, f_{\vartheta'}^{-1}\{0\}}\leq m_k(\kappa+c)$
(lemma 7.2 of \cite{C-M}).
If $\vartheta $ varies in a compact $\Theta\subset {\mathscr T}\,$, one can find a finite
subcovering of the compact $\Theta\times \bar W$ by such open sets
$\theta_l\times U_l, l=1,\cdots,m\,$, and
this implies that the claim must be true. Indeed, {\it if this was not the case},
one could find a sequence $\rho_j>0$ tending to $0$ and
$\vartheta_j\times a_j\in\Theta\times \bar W$ converging to a
point $\vartheta\times a$ verifying
\begin{equation} \label{absurd}
\sup_L \,{\bf n}_{\, L\,\cap\, B(a_j,\rho_j) \,\cap\, f_{\vartheta'}^{-1}\{0\}}\geq
m_k(\kappa+c)+1\  .
\end{equation}
But $\vartheta\times a$ belongs to some $\theta_l\times U_l$ and, by its very
definition, the local degree of $f_\vartheta$ for any $\vartheta \in \theta_l$ is
$\leq m_k(\kappa+c)$ around any point
of $U_l\,$; this gives a contradiction: for sufficiently large $j\,$, one would get
$B(a_j,2\rho_j)\subset U_l$ and $\vartheta_j\in\theta_l$ together with
(\ref{absurd}).
\end{proof}

\begin{thm}[Thom revisited, global form] \label{text11}
Given $t_0>0\,$, let $(M,g)$ be a Riemannian $n$-dimensional manifold. There exists $\rho>0$
and, in the open set of maps $f\in C^\infty(M,\R^k)$ such that $f^{-1}\{0\}$ is
compact, there exists a {\rm dense open subset} consisting of maps $f$ for which
$f^{-1}\{0\}$ intersects
all $k$-faces of Riemannian barycentric simplices $\hat\tau\in{\mathfrak S}_{t_0}(W)$
of diameter $\leq 2\rho$ in a {\rm finite} number of points.
\end{thm}
\begin{rem}
This finite number depends on the map
$f$ $($but {\it not} on
the simplices under consideration$)$.
\end{rem}
\begin{proof}
After the same initial considerations as in the proof of theorem \ref{text-11}, one applies here directly $(i)$ of
theorem
7.6 of \cite{C-M}.
\end{proof}

\subsection{On the number of disjoint simplices meeting a curve}\label{6.4}
$ $

Recall that $W$ is an open relatively compact subset of $(M,g)$.
Consider a real number
$r>0$ such that  $2r\leq{\mathcal R}_W$ (defined in lemma \ref{L.constante}) and even 
$2r\leq {\mathcal R}_0(1/2)$, where
the last constant refers to proposition \ref{proposition A}, making $C=1/2\,$.
Let $\gamma$ be a differentiable curve of finite length sitting in $W\,$.
\begin{defn} \label{bds1} 
Call $r$-tube around $\gamma$ the set 
\begin{equation*}
\hbox{\rm tub}_r(\gamma)=\{q\in M\mid \exists p\in \hbox{\rm rge}(\gamma) 
\hskip2mm \hbox{such that} \hskip4mm d(p,q)\leq r\} \  .
\end{equation*} 
Thus, any subset of $M$ of diameter $\leq r$ which meets $\gamma$ is included 
in
the $r$-tube around $\gamma\,$.
\end{defn}
\begin{lem} \label{bds2} 
For any $r >0\,$, defining $N-1$ to be the integer part of 
$\hbox{\rm length}(\gamma)/2r\,$,
one can choose a {\em net} of $N$ points $p_i$ on $ \hbox{\rm rge}(\gamma)$ such that 
any point $x$ in $ \hbox{\rm rge}(\gamma)$ verifies 
\begin{equation*} 
\exists i\in \{1 , \dots , N\} \quad \hbox{such that} \quad
d(p_i,x)\leq r\  .
\end{equation*}
Moreover, calling $b_n$ the volume of the unit Euclidean ball of 
radius $1\,$, one has
\begin{equation*}
\hbox{\rm vol}\, \hbox{\rm tub}_r(\gamma)\leq 3^nb_n\,
N\,r^n\   .
\end{equation*}
\end{lem}

\begin{proof} 
The first claim is even true for the induced metric structure on the 
curve $\gamma\,$, so it 
holds {\it a fortiori} while considering the ambient metric.

With the help of proposition \ref{proposition A}, one gets a uniform estimate 
for 
the volume of a ball having center $x\in W$ and radius $2r$
\begin{equation*}
\hbox {\rm vol} \,B(x,2r) \leq 
(1+\frac{1}{2})^n \,{\rm vol} \,B(0_x,2r)=3^n b_n r^n\  .
\end{equation*}
On the other hand, any point of $\hbox{\rm tub}_r(\gamma)$ is at a distance 
$\leq 2r$ from a point 
$p_i$ in the {\it net} chosen in the statement, hence
$\hbox{\rm tub}_r(\gamma)\subset\cup_i \,B(x_i,2r)$. The second statement follows.
\end{proof}

\begin{defn}\label{retourtheo2n} Now,
$T:K\rightarrow  W\subset M$  is a finite embedded simplicial complex 
(of the type seen in definition \ref{A.approx}) and 
we assume that for any $n$-simplex $\sigma \in K\,$, the 
diameter of $\tau = T(\sigma)$ is $\leq 2\rho$ and its volume is 
$\geq A_0 \rho^n\,$, where $2\rho=r\leq {\mathcal R}_0(1/2)/2$ and $A_0>0$ are 
given.

Actually, $\tau,K$ and $A_0$ point at $\hat\tau_E,K_E$ and $A_0$ in
theorem \ref{2n}. 
\end{defn}
{\it Without counting the degrees} of intersections of a given simplex 
with a curve
$\gamma$ {\em now running in}  $T(K)$, i. e. {\it counting $1$} if there is {\it at least} such an intersection and $0$ 
{\it if not}, one can give a 
useful bound.

\begin{prop}\label{bds4} With the assumptions and notations of lemma {\rm\ref{bds2}} and definition {\rm\ref{retourtheo2n}},
setting ${\mathcal B}\!=\!6^n\,b_n/(2A_0)$, the number $\nu$ of $n$-simplices 
of $T(K)$ meeting 
$\gamma\,\subset \,T(K)$ at least once is bounded, actually
\begin{gather*}
\nu\leq {\mathcal B} \,\hbox{\rm length}(\gamma)\,\frac{1}{\rho}\ \ \ 
\hbox{if} 
\ \  \ \hbox{\rm length}(\gamma)\geq 4\rho \ \ \ \hbox{and} \\ 
\  \ \nu< 2{\mathcal B}\ \ 
\ \hbox{if} \  \ \ \hbox{\rm length}(\gamma)< 4\rho\,.
\end{gather*}
\end{prop}

\begin{proof}
If $\hbox{\rm length}(\gamma)\!\geq\! 4\rho\,$, as the interiors of the $n$-simplices $\tau$ are
disjoint and those meeting $\gamma$ are $\subset\hbox{\rm tub}_{2\rho}(\gamma)$ 
(definition \ref{bds1}), one gets (using lemma \ref{bds2})
\begin{multline*}
\nu\,A_0\,\rho^n\leq \sum_{\tau \mid \tau\cap \gamma \not=\emptyset}
\hbox{\rm vol}\,\tau\leq 
\hbox{\rm vol}\,\hbox{\rm tub}_{2\rho}(\gamma) \leq \\ 
\leq  {3^nb_n}\,(\hbox{\rm int}
(\frac{\hbox{\rm length}(\gamma)}{4\rho})+1)\,(2\rho)^n
 \leq \frac{6^nb_n}{2}\,\hbox{\rm length}(\gamma)\,\rho^{n-1}\  .
\end{multline*}
Thus, defining ${\mathcal B}=6^n\,b_n/(2A_0)$, one gets the result.

If $\hbox{\rm length}(\gamma)< 4\rho\,$, take one point $x_1$ on $\gamma$ which 
is at distance $< 2\rho$ from any other point on the curve ; all simplices meeting 
$\gamma$ and of diameter $\leq2\rho$ 
are contained in the same open ball $B(x_1,4\rho)$: one gets 
(use proposition \ref{proposition A})
$\nu\,A_0\,\rho^n < 3^n\,b_n\,(2\rho)^n\,$,
thus $\nu<2\,{\mathcal B} \,$.
\end{proof}

\begin{cor} \label{bds45} Same assumptions and notations as in proposition {\rm\ref{bds4}}. The amount of $(n-1)$-simplices in $T(K)$
cut at 
least once by $\gamma$ is at most 
$(n+1)\,{\mathcal B} \,\hbox{\rm length}(\gamma)/(2\rho)$, 
if $\hbox{\rm length}(\gamma)\geq 4\rho\,$, and at most $(n+1)\,{\mathcal B}\,$, if 
$\hbox{\rm length}(\gamma)< 4\rho\,$.
\end{cor}

\begin{proof} Each $(n-1)$-simplex in $T(K)$ is the $(n-1)$-face of an $n$-simplex.
The upper bound then results from definition \ref{retourtheo2n} which assumes the lower bound $\hbox{vol}_g(\tau)\geq A_0\,\rho^n$ for each $n$-simplex $\tau\in T(K)$ (or, applying theorem \ref{2n} in the setting of this text,  $\hbox{vol}_g(\hat\tau_E)\geq A_0\,\rho^n$ for each $n$-simplex $\hat\tau_E\in\hat T(K_E)$), proposition \ref{bds4} and both crude observations:  
{\it each $n$-simplex has $n+1$ 
faces of dimension $n-1$} and {\it each interface crossed by $\gamma$  is the 
face of two adjacent simplices.} 
\end{proof}

\section{Picking up the pieces}\label{7}

Below $D$ denotes a compact domain in $(M,g)$,
$T$ a piecewise differentiable (proper) embedding of a finite simplicial complex $K$ into $M$ and $W$ a relatively
compact open set of $M$. More, they are assumed to verify $D\subset T(K)\setminus\partial T(K)$ and $T(K)\subset W$.

\begin{rem} \label{bibopalula}
In this section, the interplay between the following three metrics is essential (go back to definition \ref{D.metric1}):

1. the Riemannian metric $g$ on the open manifold $M\,$;

2. the piecewise linear flat metric $\hat g$ on $D\subset M$ pieced together by requiring, on each Riemannian barycentric $n$-simplex $\hat\tau\subset W$ with vertices those of $\tau=T(\sigma)$, that $\hat g_{\mid\hat\tau}$ is, through the Riemannian barycentric coordinates, the Euclidean metric such that,
between the vertices, the Euclidean distances equal the Riemannian distances;

3. the analytic flat metric $g_0$ on $K\setminus K_{n-2}$, 
derived on each $n$-simplex with $\hat g$ from $g\,$, then pasted over $K\setminus K_{n-2}\,$.

By construction, each
$(\hat\tau,\hat g)$ is isometric through $\hat T$ to $(\sigma,g_0)$.

\end{rem}

\begin{lem} \label{calLD} 
Let $W_1$ and  $W_2$ be relatively compact open sets contained in the open manifold $M$ with $W_1\subset W_2$  and $D\subset W_1\,$.
Any embedded simplicial complex $T:K\rightarrow W_2$ such that 
$D\subset T(K)\subset W_2$ gives rise to an embedded simplicial complex 
$T':K'\rightarrow W_1$ such that $D$ is contained in the 
interior of $T'(K')$ and one has, in the notations of lemma {\rm \ref{2a}}
\begin{equation*}
\alpha_T\leq \alpha_{T'}\leq \beta_{T'}\leq \beta_T\  .
\end{equation*}
\end{lem}

\begin{proof} 
As $W_1\subset W_2\,$, one can refine $T$ by performing barycentric subdivisions 
until the simplices of a refined $T'$ meeting $D$ are all contained in $W_1\,$. 
This gives the result.
\end{proof}
\begin{rem} For any regular curve 
$\gamma:[0,1]\rightarrow D\setminus T(K_{n-2})$, denote by 
${\mathcal P}_\gamma^g$ (respectively by ${\mathcal P}_{T^{-1}\gamma}^{g_0}$) the 
parallel 
translation along $\gamma$ (from $\gamma(0)$ to $\gamma(1)$) in the metric $g$ 
(respectively the 
corresponding one along $T^{-1}\gamma$ in the metric $g_0$).
\end{rem}
{\bf Caution}. {\it From now on, the symbol $T$ often readily designates the refined triangulation we used to name $\hat T\,$, or something of the same kind. This is quite natural since, whatever is the initial trangulation we consider, we shall soon be willing to work with triangulations sharing all properties of definition  {\rm\ref{convpolyhedra}} and theorem {\rm\ref{theorem I}} below, so why should those guys always wear a hat, like $\hat T\,$ pointing to an ancestor $T\,$?}
\begin{defn} \label{convpolyhedra} 
Given an integer $c\geq0\,$, given a compact domain $D$ in an $n$-dimensional open Riemannian 
manifold $(M,g)$,
a {\it strong metric approximation by polyhedra of $(M,g)$ over $D$} requires 
constants ${\mathfrak C}_1$, ${\mathfrak C}_2(c)$, ${\mathfrak C}_3>0$ and 
${\mathfrak C}_4>0\,$, determined only by $D\subset(M,g)$, such that, for any
$\rho>0$ exists a triplet $(T,K,g_0)$, depending on $\rho\,$, consisting of a polyhedron (see definition \rm{\ref{D.metric133}}) $(K,g_0)$ and a  
differentiable simplicial embedding (definition \ref{D.Simplicial.Embedding}) $T$ of $K$ into $M$ verifying

$(i)$ $T(K)$ contains $D$ in its interior.

$(ii)$ One has 
for any $n$-simplex $\sigma$ of $K$ having vertices 
$s_0,s_1,\dots,s_n\,$, 
\begin{gather*}\hskip6mm
(a)\hskip6mm\forall i,j=0,1,\dots,n
\hskip6mm d_g(T(s_i),T(s_j))=d_{g_0}(s_i,s_j)\,;\\
\hskip6mm(b)\hskip6mm \forall v\in T\sigma\hskip6mm\vert T^\ast g(v,v)-g_0(v,v)\vert 
\leq {\mathfrak C}_1\ \rho^2\ T^\ast g(v,v)\,. 
\end{gather*}

$(iii)$ Given a $c$-manifold ${\mathscr T}$ of parameters
$\vartheta$, one can find  an open dense subset ${\mathcal V}_c$ of the open set in $C^\infty({\mathscr T}\times M,\R^k)$ of families $\bar f$
of regular proper maps $f_{\vartheta}$, where
$f_\vartheta(\cdot):=\bar f(\vartheta,\cdot)$, for which, if $\Theta$ is a
compact subset of ${\mathscr T}$, one can find
$\rho_{\bar f}>0$ such that,
for any $\rho\in\,]0,\rho_{\bar f}]\,$, any $\vartheta\in\Theta\,$, and any regular
curve $\gamma:[0,1]\rightarrow M$ with $\gamma([0,1])\subset f_\vartheta^{-1}\{0\}\,$, one has, for any $u\in T_{\gamma(0)}M$
\begin{equation*}
\hskip1mm (c)\hskip1mm \Vert{\mathcal P}_{\gamma}^g(u)-
dT ({\mathcal P}_{T^{-1}\gamma}^{g_0}(dT^{-1}(u)))\Vert_g\leq
{\mathfrak C}_2(c)\rho\,\hbox{\rm length}(\gamma)\,\Vert u\Vert_g\  .
\end{equation*}

$(iv)$ For any $n$-simplex $\sigma$ of $K$ and any $p\in
T(\sigma)$, 
the thickness $t_{g_p}(L_p(\sigma))$ (see definition \ref{secantmap}) and the 
$g$-volume of $T(\sigma)$ verify
\begin{equation*} 
t_{g_p}(L_p(\sigma))\geq {\mathfrak C}_3 \hskip6mm
\hbox{and} \hskip6mm 
\hbox{\rm vol}_g(T(\sigma))\geq {\mathfrak C}_4\,\rho^n\ . 
\end{equation*}
\end{defn}

\begin{rem} \label{remarquable1} In view of definition \ref{genthick}, if $T$ is a former $\hat T$, the above property $(iv)$ also reads: for any $n$-simplex $\sigma$ of $K$, 
the generalised thickness $t_g(\hat\tau)$ and the 
$g$-volume $\hbox{\rm vol}_g(\hat\tau)$ of the Riemannian barycentric simplex $\hat\tau=\hat T(\sigma)$ (see definition \ref{hatT1}; use $A_0\geq {\mathfrak C}_4$, see (\ref{glouhips})) verify
\begin{equation*} 
t_g(\hat\tau)\geq {\mathfrak C}_3 \hskip6mm
\hbox{and} \hskip6mm \hbox{\rm vol}_g(\hat\tau)\geq {\mathfrak C}_4\,\rho^n\ . 
\end{equation*}
\end{rem}

\begin{thm} \label{theorem I} 
Let $D$ be a compact domain in an open $n$-dimensional Riemannian manifold $(M,g)$. 
For any given integer $c\!\in\!\N$, one can find constants 
${\mathfrak C}_1, {\mathfrak C}_2(c),{\mathfrak C}_3,{\mathfrak C}_4$, 
only determined by $c$ and $D\subset(M,g)$, such that there exists a strong metric 
approximation by 
polyhedra of $(M,g)$ over $D$ related to those constants as stated in the definition
{\rm \ref{convpolyhedra}} above.
\end{thm}
{\bf Caution}. {\it At the beginning of this proof $T,K$ deserve once more the roles of an initial triangulation, but at its end $\hat T_E, K_E$ (the integer $E$ is paired with $\rho$, the mesh $\rho$ goes to $0$ as $E$ goes to infinity, see definition \rm{\ref{linkint}} and around) will be renamed $T,K$, in accordance with the new roles they play in the statement.}
\begin{proof}
Choose a relatively compact open set $W$ that contains $D\,$.
Cover $D$ by a finite number of relatively compact open sets $C_i\subset M$ that are the 
interior of 
cubical images of charts $\varphi_i:[0,1]^n\rightarrow M\,$, all included in $W\,$. Applying 
theorem $10.4$ of \cite{Mu}, one can build up a finite differentiable simplicial embedding 
from the 
``union'' of those simplicial embeddings $\varphi_i:K_i\rightarrow C_i$ (based on a 
simplicial decomposition of the cube $[0,1]^n$) and 
even manage this ``union'' to contain $D$ in its interior and to be contained 
in $W\,$, as the union of the $C_i$ does. Here, call $T: K\rightarrow M$ the resulting 
{\em initial} embedded 
simplicial complex with $D$ 
contained in the interior of $ T(K)$ and $T(K)$ contained in $W\,$. 

Consider the non empty set ${\mathfrak T}_{D,W}$ of all embedded finite simplicial complexes
$(T,K)$ of the kind:

$T:(K,{\rm std})\rightarrow W\subset (M,g)$ is such that $T(K)$ contains $D$ in its 
interior.

Refer to lemma \ref{2a} and define the positive number 
${\mathcal L}_{{\mathfrak T}_{D,W}}$
\begin{equation*}
0<{\mathcal L}_{{\mathfrak T}_{D,W}}=\sup_{T\in {\mathfrak T}_{D,W}}
\frac{\alpha_T}{\beta_T}\leq 1\  .
\end{equation*}
Though $\frac{\alpha_T}{\beta_T}$ is always positive, it can be made as small as 
wanted while choosing $(T,K)\in {\mathfrak T}_{D,W}\,$. Now, using lemma \ref{calLD}, one 
gets 
for $W_1\subset W_2$ the 
inequality 
${\mathcal L}_{{\mathfrak T}_{D,W_1}}\geq {\mathcal L}_{{\mathfrak T}_{D,W_2}}\,$. 
But, one also has - by definition of ${\mathfrak T}_{D,W}$ - the inclusion 
${\mathfrak T}_{D,W_1}\subset{\mathfrak T}_{D,W_2}\,$, thus 
${\mathcal L}_{{\mathfrak T}_{D,W_1}}\leq {\mathcal L}_{{\mathfrak T}_{D,W_2}}\,$. 
Consequently, ${\mathcal L}_{{\mathfrak T}_{D,W}}$ depends only on $D\,$: we denote it 
by ${\mathcal L}_{{\mathfrak T}_D}\,$.

Choose a $(T ,K )\in {\mathfrak T}_{D,W} $ such that 
$\frac{\alpha_{T }}{\beta_{T }} \geq {\mathcal L}_{{\mathfrak T}_D}/2\,$. 
This will produce constants $t_0$ and $A_0$ associated with that 
$(T ,K )$ (see 
formula (\ref{2c}),  
proposition \ref{2s}) which are bounded 
from below only in terms  of 
${\mathcal L}_{{\mathfrak T}_D}$  
\begin{equation}\label{glouhips}
A_0\geq{\mathfrak C}_4=\bigg{(}\frac{{\mathcal L}_{{\mathfrak T}_D}}{24}\,
\frac{{{\mathscr D}}_n}{{\mathcal D}_n}\bigg{)}^{\!\!n}\,\Theta_n\hskip4mm\hbox{and}
\hskip4mm t_0\geq {\mathfrak C}_3=\frac{n}{n+1}\ 
\frac{{\mathfrak C}_4}{4\,\vert\sigma_{n-1}\vert}\  .
\end{equation}
Now, theorem \ref{2n} applied to $(T ,K )$ implies $(iv)$. {\it Indeed}, it implies that, for any sufficiently small $\rho>0\,$, the regular subdivision of integer 
$E=E(\rho)$ is such
that each $n$-simplex $\sigma_E\in K_E$ gives rise 
to a Riemannian barycentric simplex $\hat\tau_E=\hat T(\sigma_E)$ sharing
the following properties
\begin{align*}
\dag(i) &   \hskip4mm & \hat\tau_E\subset B(\varpi,\rho)  & \hskip3mm \hbox{and} \hskip4mm
\hbox{\rm diam}_g(\hat\tau_E)\leq 2\,\rho\, ; & \\
\dag(ii) &  \hskip4mm & t_g(\hat\tau_E)\geq t_0 \geq {\mathfrak C}_3 \,; & &    \\
\dag(iii) &  \hskip4mm & \hbox{\rm vol}_g(\hat\tau_E)\geq A_0\,\rho^n\geq 
{\mathfrak C}_4\,\rho^n \,,&  &   
\end{align*}
where $\varpi$ is the gravity center of 
$\hat\tau_E\,$, thus giving $(iv)$ of definition \ref{convpolyhedra}.

Notice that, in the above situation, theorem \ref{2n} also implies that each $\hat \tau_E$ shares the properties of theorem \ref{theorem 3}, where $t_0$ is given above. Thus, there exists a 
constant ${\mathcal C}_7\,$, depending on $t_0\,$,  such that for 
any $x$ in $W$ and $\rho\in]0,{\mathcal R}_3]\,$, each Riemannian barycentric simplex 
$\hat\tau\subset B(x,\rho)$ 
having thickness $t_g(\hat\tau)\geq t_0$ and Euclidean realization $(\hat\tau,\hat g)$ verifies, for 
any $q\in \hat\tau\,$, 
any $u\in T_q M\,$, any differentiable curve $\gamma\subset\hat\tau$ starting at $q$
\begin{align*}
(i)&   
 \hskip6mm & \vert\,\Vert\,u\,\Vert_{\hat g}^2-\Vert\,u\,\Vert_g^2\,\vert & 
 \leq \Vert\,u\,\Vert_g^2 \  {\mathcal C}_7\,\rho^2\,; \\ 
 (i') &    \hskip6mm & \vert\,\Vert\,u\,\Vert_{\hat g}-\Vert\,u\,\Vert_g\,\vert & \leq 
\Vert\,u\,\Vert_g \  {\mathcal C}_7\,\rho^2\,; \\
(ii) & \hskip6mm & H=\sup_{q\in\hat\tau}\Vert\,h_q\,\Vert_g & \leq 
{\mathcal C}_7\,\rho\,;\\
(iii) &\hskip6mm & 
\Vert\,({\mathcal P}_\gamma^{\hat g}-{\mathcal P}_\gamma^g)(u)\,\Vert_g & 
\leq  \hbox{\rm length}_g (\gamma)\ {\mathcal C}_7\,\rho\ \Vert u\Vert_g \, .
\end{align*}
To summarize, if $\rho$ is small enough (to be able to apply theorems 
\ref{theorem 3}, \ref{2n} and \ref{hatT}), 
all simplices $\tau_E=T_E (\sigma_E)$ are 
Euclideable by theorem \ref{2n} and share also the results of theorem 
\ref{theorem 3}.
By theorem  \ref{hatT}, the map $\hat T_E$ which replaces each 
$T_E(\sigma_E)=\tau_E$ 
by its barycentric Riemannian companion $\hat \tau_E$  - if $\rho$ is small enough - 
($\hat T_E$ still goes 
from $K_E$ to $M$) is a differentiable simplicial embedding such that $\hat 
T_E(K_E)$ contains $D$ 
in 
its interior and $\overline{\hat T_E(K_E)}\subset W\,$: for $E$ large, theorem \ref{hatT} implies
that $\hat T_E$ and $T_E$ are 
close in a way 
uniformly controlled by $\rho$ on each $\sigma_E$ since $(i)$ holds. 
Consider $K_E$ equipped with the metric $g_0$ given in definition \ref{D.metric1}. Recall that 
$g_0$ is an analytic flat metric singular over the $(n-2)$-skeleton of $K_E$ and depends 
on $E\,$. Moreover  for any $n$-simplex $\sigma\in K_E\,$, the map 
$ \hat T_E : (\sigma, g_0) \rightarrow (\hat T_E(\sigma), \hat g)$ is an isometry.

Deduce from definition \ref{constanteC_7R_3} that the constant ${\mathcal C}_7$  of theorem \ref{theorem 3} verifies
$${\mathcal C}_7\leq \tilde{\mathcal C}_7:=2^6\,3^2\,n^2\,t_0^{-6}({\mathcal M}_2+2^33\sqrt{3n}\,{\mathfrak R}_0)\ ,
$$
and as $\tilde{\mathcal C}_7\,t_0^6$ only depends on $D\subset(M,g)$ (restricting the Riemannian bounds ${\mathfrak R}_0$ and ${\mathcal M}_2$ to be estimated only over $D$), {\bf defining}
$\tilde{\mathfrak C}_1:=2^6\,3^2\,n^2\,{\mathfrak C}_3^{-6}({\mathcal M}_2+2^33\sqrt{3n}\,{\mathfrak R}_0)\,$ {\bf and setting} $(T,K):=(\hat T_E,K_E)$, all statements given $(i),(ii)$ in 
definition \ref{convpolyhedra} are 
established with $\tilde{\mathfrak C}_1$ in place of ${\mathfrak C}_1\,$, except $(iii)$ which concerns parallel translations.

As for this last point, thanks to $(iii)$ of theorem \ref{theorem 3} just recalled, 
we know that {\it the result holds separately on each $\hat \tau_E\,$} (with $\tilde{\mathfrak C}_1$ in place of ${\mathfrak C}_1$), and if there 
wouldn't be the jumps as $\gamma$ crosses interfaces between  adjacent simplices, setting 
${\mathfrak C}_2=\tilde{\mathfrak C}_1$ would end the job.
So, we only have to care about those jumps, i. e. about the discrepancies caused by the 
crossings.

Thanks to theorem \ref{text-11} (Thom revisited, local form), applied in the case $k=n\!-\!1$,
bring in the generic set ${\mathcal V}_c$ of $C^\infty$-families of curves, choosing families of curves
$\bar\gamma$ as in definition \ref{convpolyhedra} $(iii)$ and restricting the parameters to vary in a compact 
$\Theta \subset {\mathscr T}$. There exists $\rho_{\bar f}>0$ 
(or $\rho_f>0$ if $c\!=\!0$) such that, for $\rho\in\,]0,\rho_{\bar f}]$ (for 
$\rho\in\,]0,\rho_f]$) any $\gamma_\vartheta\in\bar\gamma\in{\mathcal V}_c$ 
(any $\gamma\in{\mathcal V}$) cuts, {\it counting degrees}, every $(n-1)$-face (even the open leaf containing this face, {\it see remark {\rm \ref{n-2?...}}}) of any 
barycentric simplex in $ W$ with diameter $\leq 2\rho$ and thickness 
$\geq t_0\geq{\mathfrak C}_3$ only a finite number 
of times $\leq m_{n-1}(\kappa+c)$ (of course, $c=0$ if we consider a single curve).  

Hence, bringing in corollary \ref{bds45} (assuming 
$\rho$ such that $\hbox{\rm length}(\gamma)\geq 4\rho$), the number (counted 
with degrees) of $(n-1)$-faces cut by $\gamma$ (by each $\gamma_\vartheta$ in 
$\bar\gamma$) is bounded by
$(n+1)\, m_{n-1}(\kappa+c)\,\hbox{\rm length}(\gamma)\,{\mathcal B} /(2\rho)$, where 
we even know that ${\mathcal B}$ is smaller than $6^n\,b_n/(2{\mathfrak C}_4)$ (still 
write ${\mathcal B}$ for this bound $6^n\,b_n/(2{\mathfrak C}_4)$, depending only on 
$D\subset(M,g)$).

Apply proposition \ref{proposition 1st} to the $g$ and $\hat g$-parallel translations along a curve, considering first the parallel fields to be continuous at any point $q$ in some $\hat\tau_1\cap\hat\tau_2$.
Control each succssive rip $\Phi_q$ thanks to proposition \ref{rip1} and (\ref{rip331}). Conclude in the way (\ref{rip332}) is established from (\ref{rip331}).

Setting ${\mathfrak C}_2(c)=\tilde{\mathfrak C}_1+6\,\sqrt{2}\, (n+1)\,
m_{n-1}(\kappa+c)\,\tilde{\mathfrak C}_1\,{\mathcal B}\,$ gives the generic asymptotic uniform 
statement $(iii) (c)$ in definition 
\ref{convpolyhedra}.
\end{proof}
\begin{rem}\label{n-2?...} A generic single curve $\gamma$ does not meet any $(n-2)$-face of one of the Riemannian barycentric simplices $\hat\tau_E\,$. But, in general,  some curves of a generic family of curves $\bar\gamma$ (depending on at least one real parameter) will cut $(n-2)$-faces. Along such curves the $\hat g$-parallel translation is generally ambiguous, we have to deal with such situations in part \ref{M-T2}. 

We make another concluding observation: given any curve $\gamma$ in a family $\bar\gamma$ indexed by a compact set of parameters, the difference between the $g$ and $\hat g$-parallel translations, if this one is defined, along any curve $\gamma\in\bar\gamma\,$ is uniformly controlled in terms of $\rho$ and on $\sup_{\gamma\in\bar\gamma}\rm{length}(\gamma)$.
\end{rem}

\section{Appendices}
\subsection{Proof of remark \ref{strictconv}}\label{strictconv1}
$ $
\begin{proof} Given $x'\in B(x,r)$ and $u\in \{\nabla\varrho\}_{x'}^\perp\setminus\{0\}\,$, observe $(Dd\varrho)_{x'}(u,u)\,$ is the second derivative of the length integral computed for the variation 
$$V_\epsilon:t\in [0,\varrho(x')]\mapsto V_\epsilon(t)=V(t,\epsilon):=\exp_x(t(v+\epsilon\,w))\in M\ ,
$$
where the parameter of variation $\epsilon$ describes an interval around $0\,$ and $c\!=\!\exp_x(\cdot\,v)$ is the unique geodesic $c$ minimising the distance between $x\!=\!c(0)$ and $x'\!=\!c(\varrho(x'))$, while $w\!\in\! T_xM$ is sent to $u$ by parallel translation along $c$ from $x$ to $x'\,$. So, $(Dd\varrho)_{x'}(u,u)\leq 0$ would imply that $Y(t)\!=\!(\partial V/\partial \epsilon)(t,0)$ is a non trivial Jacobi field vanishing at $t\!=\!0$ and $t\!\leq\!\varrho(x')$; thus, given any small $\eta\!>\!0$, the geodesic $c$ wouldn't minimise the distance between $c(0)\!=\!x$ (between $c(\varrho(x'))\!=\!x'$) and any point $c(\varrho(x')+\eta)$ (any point $c(-\eta)$), see for instance \cite{G-K-M}, page 143, Lemma 1 or \cite{K-NII}, page 87, theorem 5.7. Actually, the cut-radius is never greater than the conjugate radius (while the last is always $\geq \pi/\sqrt{\sup(0,K_{\rm max})}$, where $K_{\rm max}$ is the supremum of the sectional curvatures in $(M,g)$). But this does not fit with the fact that $c$ minimises the distance between any two points belonging to a diameter of $B(x,r)$.
\end{proof}

\subsection{Some geometrical considerations about Status}\label{app1}
$ $
\begin{thm}\label{rbs7}
Given $(M,g),(M',g')$ two Riemannian 
surfaces and ${\bf p}$ and ${\bf p}'$ spread sets of points, 
sitting in convex balls $B\subset M,\,B'\subset M'\,$,
let ${\mathscr S}$, ${\mathscr S}'\,$ be the status maps associated to 
${\bf p}$ and ${\bf p}'$, let 
$U\supset {\bf p}\,,\,U'\supset {\bf p}'$ be open sets included in the domains 
of injectivity of 
${\mathscr S}$ and ${\mathscr S}'\,$ (see proposition {\rm \ref{rbs5}}). 
{\em Assume that} there is a bijective pairing of $x\in U$ with $x'\in U'$ such that $x$ and $x'$ share the $g$ and $g'$-same Riemannian barycentric coordinates and verify $d_g(p_i,x)=d_{g'}(p'_i,x')$ for $i=0,1,2\,$, {\em or assume, equivalently, that} the images ${\mathscr S}(U)$ and ${\mathscr 
S}'(U')$ are identical in 
$\R^{3}\,$. {\em Then} the uniquely defined map $\psi\,$ that verifies ${\mathscr S}(x)={\mathscr S}'(\psi(x))$ for $x\in U$ and $x'=\psi(x)\in U'$ (paired in the assumed pairing) is a Riemannian isometry from $(U,g)$ onto $(U',g')$.
\end{thm}

\begin{rem}\label{rbs8} 
From this result, one can infer for instance that a flat 
region in a 
Riemannian {\it surface} is characterized by the fact that Status, for admissible choices of 
points, has an image 
containing a piece of an elliptic paraboloid. Indeed, it can be shown that images of Satus-maps defined in flat spaces are contained in elliptic paraboloids (quadrics in ${\mathbb R}^{n+1}$ defined, in well-chosen coordinates, as a locus $z_0=z_1^2+\dots+z_n^2$).
\end{rem}

\begin{proof} 
First make some observations implied by the same data in any dimension $n\,$, i. e. not only for surfaces. By lemma \ref{rbs4}, the point $p_i$, being the unique point 
such 
that the tangent plane to ${\mathscr S}(U)$ at ${\mathscr S}(p_i)$ is orthogonal 
to $e_i\,$, corresponds 
to $p'_i$, which shares the same property, thus $\psi(p_i)=p'_i\,$. Moreover 
$\psi(x)$ and $x$ have the same ``barycentric coordinates''
$(\lambda_0(x),\dots,\lambda_n(x))\,$: indeed,
$x$ and $\psi(x)$ are sent to the same point in the identical images by ${\mathscr 
S}$ and ${\mathscr S}'\,$, and the normal at ${\mathscr 
S}(x)={\mathscr S}'(\psi(x))$ determines $(\lambda_0(x),\dots,\lambda_n(x))\,$; the definition 
of Status gives $d_g(p_i,x)=d_{g^\prime}(p'_i,\psi(x)))$ for $i=0,1,\cdots,n$,  
so
$\psi(x)=x'\,$.
\begin{lem}\label{auxiliaireinfinit}
\par Set $a_i=-\lambda_i(x) (\exp^g_x)^{-1}p_i$ and $a'_i=-\lambda_i(x) (\exp^{g'}_{x'})^{-1}p'_i$ 
for $i=0,1,...,n\,$. For $i=0,1,\dots,n$ and any $v\in T_xM\,$, one has
\begin{equation}\label{barycorr}\sum_i \!a_i \!=\! \sum_i \!a'_i=\!0\, ,\  \Vert a_i\Vert_g=\Vert a'_i\Vert_{g^\prime}\    {\rm and}
\ \langle a_i,v\rangle_g= 
\langle a'_i,d\psi(x)(v)\rangle_{g'}\, .
\end{equation}
\end{lem}
\begin{proof} Take the $x$-derivative of $d_g(p_i,x)=d_{g^\prime}(p'_i,\psi(x))$ to get the last equalities, the first are direct features of the stated Riemannian barycentric $g$ to $g'$-correspondence ($x$ and $x'$ are minimum loci).
\end{proof}
One can make a step further, identifying through $\psi$ the Riemannian barycentric simplex defined by ${\bf p}$ in $(M,g)$ with the one defined by ${\bf p}'$ in $(M',g')\,$, thus getting on the same simplex, say the one
generated by ${\bf p}$ in $(M,g)$, another metric $\psi^* g'$, which we still denote by $g'$: now $\psi$ reads $\psi={\rm Id}\,$ and the result to prove reduces to $g=g'\,$. From the equalities (\ref{barycorr}), one gets for {\it surfaces} and $i,j=0,1,2$
\begin{equation}\label{barycorr2} \langle a_i,a_j\rangle_g=\langle a'_i,a'_j\rangle_{g'}\  .
\end{equation}
Indeed, in view for instance of $a_0+a_1+a_2=a'_0+a'_1+a'_2=0\,$, write
\begin{multline*}
2\langle a_0,a_1\rangle_g=\Vert a_0+a_1\Vert_g^2-\Vert a_0\Vert_g^2-\Vert a_1\Vert_g^2=\Vert a_2\Vert_g^2-\Vert a_0\Vert_g^2-\Vert a_1\Vert_g^2\\
=\Vert a'_2\Vert_{g'}^2-\Vert a'_0\Vert_{g'}^2-\Vert a'_1\Vert_{g'}^2=\Vert a'_0+a'_1\Vert_{g'}^2-\Vert a'_0\Vert_{g'}^2-\Vert a'_1\Vert_{g'}^2=2\langle a'_0,a'_1\rangle_{g'}\, .
\end{multline*}
Then, mixing (\ref{barycorr2}) with the last equalities in (\ref{barycorr}), which here read (for any $i=0,1,2$) $ \langle a_i,v\rangle_g= 
\langle a'_i,v\rangle_{g'}\,$, one gets for any $i,j=0,1,2$
\begin{equation*}
\langle a'_i,a'_j\rangle_{g'}=\langle a_i,a_j\rangle_g=\langle a'_i,a_j\rangle_{g'}\ ,\ \ {\rm thus}
\ \ \langle a'_i,a'_j-a_j\rangle_{g'}=0\, ,
\end{equation*}
which implies $a_j=a'_j$ for $j=0,1,2$ and (as expected) $g=g'$.
\end{proof}

{\bf Caution}. {\it This result does not generalise in dimension $\geq 3\,$. 
}
\begin{prop} \label{ctrexple} Given any $3$-simplex with vertices $p_0,\dots,p_3\,$ in ${\mathbb R}^3$, one can find a Riemannian metric $g'$ {\em which is distinct} from the standard Euclidean structure $g$, {\em but} $d_g(p_i,\cdot)$ and $d_{g'}(p_i,\cdot)$ have the same level sets, i. e. for any $x\in{\mathbb R}^3$ one has $d_g(p_i,x)=d_{g'}(p_i,x)$ for $i=0,1,2,3\,$. So the bijective map  $\psi$ of theorem {\rm\ref{rbs7}} here reduces to 
$\psi\!=\!{\rm Id}$ and {\em is not} an isometry.
\end{prop}

\begin{proof}
Consider ${\bf p}=p_0,p_1,\dots, p_n$ a $g$ and $g'$-spread set in a manifold equipped with two Riemannian metrics $g$ and $g'\,$, with ${\bf p}$ included  in  ``convex'' $g$ and $g'$-balls. Call $\sigma$ the $g$-Riemannian barycentric simplex having vertices ${\bf p}\,$. Assume $d_i(x)=d_g(p_i,x)=d_{g'}(p_i,x)=d'_i(x)$ for any $x$ in $\sigma$ (thus for $x$ near $\sigma$), thus one has $\psi=\rm{Id}$. Setting
\begin{equation}\label{ail}a_i=a_i(x)\!:=\!-\lambda_i(x)(\exp^g_x)^{-1}p_i\ ,\ a'_i:=a_i'(x)\!=\!-\lambda_i(x)(\exp^{g'}_{x})^{-1}p_i\ ,
\end{equation}
then (\ref{barycorr}) in lemma \ref{auxiliaireinfinit}
here reads (for $i\!=\!0,1,\dots,n,x\in\sigma,v\in T_xM$)
\begin{equation}\label{barycorr3}\sum_i \!a_i = \sum_i \!a'_i=0\, ,\  \Vert a_i\Vert_{g_x}=\Vert a'_i\Vert_{g_x^\prime}\    {\rm and}
\ \langle a_i,v\rangle_{g_x}= 
\langle a'_i,v\rangle_{g'_x}\, .
\end{equation}
\begin{rem} If $x=p_i$ for some $i\,$, observe that (\ref{barycorr3}) is still verified. 
\end{rem}
\begin{rem} If $x\not=p_i$ (for any $i$), the affine mapping $L_x$ which sends $a'_i(x)$ to $a_i(x)\,$ is linear and $g_x$-symmetric positive definite.
{\em Indeed}, the linearity is straightforward from the first equality in (\ref{barycorr3}), while the fact that $L_x$ is $g_x$-symmetric positive definite comes from the last identities, which are {\it equivalent to} saying $
g'_x(\bullet,\ast)=g_x(L_x(\bullet),\ast)\,$.
\end{rem}

\begin{lem} \label{contre1} Choose $n+1$ affinely independent points $a_i$ in ${\mathbb R}^n$ (where $i=0,1,\dots,n$) verifying as vectors $\sum_i \!a_i =0\,$. Any $g$-symmetric positive definite linear map $L$ for which the corresponding $a'_i:=L^{-1}(a_i)$ satisfy (for any $i=0,1,\dots,n$)
\begin{equation}\label{contre2} \langle a'_i-a_i,a_i\rangle_{g}=0
\end{equation} 
defines a scalar product $g'(\bullet,\ast)=g(L(\bullet),\ast)\,$ such that {\rm (\ref{barycorr3})} holds.
\end{lem}
\begin{proof} The identity $\sum_i \!a'_i=0$ results from the linearity of $L\,$. As $L$ is
$g$-symmetric positive definite, $g'(\bullet,\ast)$ is a scalar product. The way $L$ defines $g'$ from $g$ and $a'_i$ from $a_i$ (for any $i$) gives the last equalities in
{\rm (\ref{barycorr3})}, while (\ref{contre2}) allows to write for any $i$
\begin{equation*}\Vert a_i\Vert_{g}^2=\langle a_i,a_i\rangle_{g}=\langle a'_i,a_i\rangle_{g}=\langle a'_i,L_x(a'_i)\rangle_{g}=
\Vert a'_i\Vert_{g'}^2\ ,
\end{equation*}
producing the remaining central equalities stated in {\rm (\ref{barycorr3})}.
\end{proof}
\begin{lem}\label{contre3} Given affinely independent $a_i\in {\mathbb R}^3,\,i\!=\!0,1,2,3$ verifying $\sum_i \!a_i \!=\!0$ (as vectors), there exists around the standard scalar product $g$ a two-dimensional embedded small disk ${\bf D}$ of scalar products $h$ paired, as in lemma {\rm\ref{contre1}}, with a symmetric positive definite linear map $L$ such that ({\rm\ref{contre2}}) and ({\rm\ref{barycorr3}}) hold (with $a'_i:=L^{-1}(a_i),\,i=0,1,\dots,3$). 
\par If all $a_i=a_i(x)$ depend in a differentiable way on a variable $x\,$, so is the change of the disk
${\bf D}_x$ that one can define relative to the $a_i(x)$.
\end{lem}
\begin{rem} \label{poof} The disk ${\bf D}$ shows up as the image of a {\it planar affine disk} ${\mathcal D}$ under the mapping $B\mapsto (B+{\rm Id})^{-1}\,$. 
\end{rem}
\begin{proof} {\it Temporarily}, $h$ deserves the role devoted to $g'\,$.
As in lemma \ref{contre1}, consider any scalar product $h$ in ${\mathbb R}^3$ to be given by a linear $g$-symmetric positive definite mapping $L$ through $h(\bullet,\ast)=g(L(\bullet),\ast)\,$. We now describe the set of scalar products $h$ that satisfy {\rm (\ref{barycorr3})} and are near enough to $g\,$. 
\par First, for $i=1,2,3$, the three identities (\ref{contre2}) imply that a mapping $L$ of interest is linked to 
$L^{-1}={\rm Id}+B\,$ where $B$ is $g$-symmetric and uniquely coupled (through $g$) with a bilinear symmetric form $\beta$ having, in the basis $a_1,a_2,a_3\,$, a symmetric matrix $\beta_{i,j}=\langle B(a_i),a_j\rangle_g=\langle B(a_j),a_i\rangle_g=\beta_{j,i}$ with diagonal coefficients $\beta_{i,i}=\beta(a_i,a_i)=$ all equal to zero. Indeed, taking $a'_i=L^{-1}(a_i)$, one gets $\langle a'_i- a_i,a_i\rangle_g=0$ since
\begin{multline}\label{condn1!!!}0=\beta(a_i,a_i)=\langle B(a_i),a_i\rangle_g=\langle (L^{-1}-{\rm Id})(a_i),a_i\rangle_g=\\=\langle L^{-1}(a_i),a_i\rangle_g-\langle a_i,a_i\rangle_g=\langle a'_i,a_i\rangle_g-\langle a_i,a_i\rangle_g=\langle a'_i- a_i,a_i\rangle_g\  .
\end{multline} 
Then (\ref{contre2}) for $i=0$ reads $\langle B(a_1+a_2+a_3),a_1+a_2+a_3\rangle_g=0\,$ or
\begin{equation}\label{condn!} \beta_{1,2} +  \beta_{1,3} + \beta_{2,3}=0\ .
\end{equation}
\par Summarizing the description, we get, in a $2$-dimensional vector subspace ${\mathcal B}$ of the space of $g$-symmetric bilinear forms, {\it a relevant affine planar disk} ${\mathcal D}$ {\it consisting} of the symmetric bilinear forms $\beta$ whose (symmetric) matrices $\beta_{i,j}$ in the basis $a_1,a_2,a_3\,$ share the constraints
\par 1. $\beta_{1,2}=\beta_{2,1}$ and $\beta_{1,3}=\beta_{3,1}\,$ are given by selecting a point $(\beta_{1,2},\beta_{1,3})$ in the {\it disk} ${\mathcal D}=\{\beta^2_{1,2} +  \beta^2_{1,3}\leq \epsilon\}$ (with $\epsilon>0$ small); 
\par 2. $\beta_{i,i}=0$ and $\beta_{2,3}=\beta_{3,2}=-(\beta_{1,2} +  \beta_{1,3})$ are imposed. 
\par The $g$-symmetric mappings $B$ uniquely determined by $\beta\in{\mathcal B}$ (meaning $\beta(\bullet,\ast)=g(B\bullet,\ast)$) provides a unique $L=({\rm Id}+B)^{-1}$ for which $h=g(L(\bullet),\ast)$ satisfies {\rm(\ref{contre2})}, thus {\rm (\ref{barycorr3})}, and is near enough to $g\,$.
Hence, the affine planar small disk ${\mathcal D}$ around $0$ of $g$-symmetric mappings $B$ is sent, through the transformation $B\mapsto L=(B+Id)^{-1}\,$, to the disk ${\bf D}$ satisfying the required properties of the statement. 
\par Furthermore, the correspondence sending a bilinear form $\beta$ to the unique $g$-symmetric linear map $B$ defined by $\beta(\bullet,\ast)=g(B\bullet,\ast)$ is differentiable, while the correspondence sending a real symmetric matrix $(\beta_{i,j})$ to the symmetric bilinear form $\beta_x\,$, whose matrix is $(\beta_{i,j})$ in the basis $a_1(x),a_2(x),a_3(x)\,$, is also differentiable because the $a_i(x)$ are assumed to be differentiable. This proves the last statement.
\end{proof}
\begin{rem} The changing of the plane ${\mathcal B}_x$ as $x$ varies may be described by writing, for a relevant symmetric bilinear form $\beta_x\in {\mathcal B}_x\,$,  its matrix $B_{i,j}(x)$ in the {\it fixed} canonical basis $e_1,e_2,e_3$ of ${\bf R}^3\,$ (denoting by $P(x)$ the matrix with coefficients $P_{l,j}(x):=\alpha_j^l(x):=\langle a_j(x),e_l\rangle_g$)
$$\beta_{i,j}=\sum_{k,l=1}^3({}^tP(x))_{i,k}\,B_{k,l}(x) \,P_{l,j}(x)\  .
$$
In the basis $e_1,e_2,e_3\,$, the coefficients $B_{k,l}(x)$ corresponding to a relevant $\beta_x$ must build a symmetric matrix (i. e. $B_{k,l}(x)=B_{l,k}(x)$) and satisfy a well defined system of four independent linear equations, the transcription in terms of the $B_{k,l}(x)$ of the four independent equations $\beta_{1,1}=\beta_{2,2}=\beta_{3,3}=\beta_{1,2}+\beta_{1,3}+\beta_{2,3}=0\,$, produced by skipping from the basis $a_1(x),a_2(x),a_3(x)$ to $e_1,e_2,e_3\,$. The first three ($i=1,2,3$) are
$$\sum_{k,l=1}^3  \alpha_i^k(x) \,\alpha_i^l(x)\,B_{k,l}(x)=0\ .
$$
The last one is
$$\sum_{k,l=1}^3 (\alpha_1^k(x)\,\alpha_2^l(x)+\alpha_1^k(x)\,\alpha_3^l(x)+\alpha_2^k(x)\,\alpha_3^l(x))\,B_{k,l}(x)=0\  .
$$
\end{rem}
\par
In ${\mathbb R}^3\,$, consider some nondegenerate simplex $\sigma\,$. Near the Euclidean canonical scalar product $g$ and defined in the interior $\sigma^{\mathrm{o}}$ of $\sigma$, {\it choose} a {\it distinct} Riemannian structure $h\,$: actually, along the lines of the above lemma \ref{contre3}, choose at any $x\in\sigma^{\mathrm{o}}$ a Euclidean structure $h_x\!\in\!{\bf D}_x$ in a $C^\infty$-coherent way (this can be done since $\sigma^{\mathrm{o}}$ is a contractible set on which the hypotheses of lemma \ref{contre3} are fulfilled). 
Actually, we choose those distinct natural ``Riemannian'' metrics $g$ and $h\,$ so as to share all properties of lemma {\rm\ref{contre1}} ``with respect to'' $\sigma\,$, meaning, for $x\in\sigma^{\mathrm{o}}\,$ and $a_i(x)=-\lambda_i(x)\,(p_i-x)$ (see lemma \ref{auxiliaireinfinit}), we can rely on lemmas \ref{contre1}, \ref{contre3}.
From lemma \ref{contre3} and {\it its proof}, $h$ is given at any $x\in \sigma^{\mathrm{o}}$ through some linear $g_x$-symmetric $B_x$ in a small affine disk ${\mathcal D}_x$ (remark \ref{poof}): one takes $L^h_x=(B_x+{\rm Id})^{-1}\,$ and gets $h_x(\bullet,\ast)=g_x(L^h_x(\bullet),\ast)$. Take a small $g$-ball $b$ contained in $\sigma^{\mathrm{o}}$ and a very tiny ``bump function'' on $b\,$, i. e. a $C^\infty$-function $\varphi:{\mathbb R}^3\rightarrow {\mathbb R}\,$ with support in $b\,$, everywhere $\geq0\,$, somewhere $>0$ and having a very small $C^2$-norm (with respect to $g$). 
 \begin{defn} On ${\mathbb R}^3$, define the Riemannian metric $\, g'_x(\bullet,\ast)\!=\!g_x(L_x(\bullet),\ast)\,$ with $\, L_x:=(\varphi(x)\,B_x+{\rm Id})^{-1}$ if $x\in b$ and $g'_x=g_x$ if $x\notin b\,$.
\end{defn}
The metric $g'$ is by definition equal to $g$ outside of $b\,$, distinct from $g$ somewhere on $b\,$, {\it although very close} to $g$ in the $C^2$-sense. Thus the $g'$-geodesics are on a huge distance very close to the $g$-geodesics and the curvature of $g'$ remains very close to $0\,$ (which is the curvature of $g$) in a $g$-big set containing $\sigma\,$. By a theorem of Whitehead (see \cite{C-E}, p. 103), the $g'$-convexity radius is very large, thus $\sigma$ is contained in a big $g'$-convex ball (the $C^2$-smallness of $\varphi$ is taken very small ``compared'' to the given size of $\sigma\,$, producing a large enough $g'$-convexity radius).
 \par As $g'$ shares with $g$ all equalities (\ref{barycorr3}) {\it with respect to the simplex $\sigma$} (and $a_i\,$, at any $x\in b\,$, is
$a_i(x)=-\lambda_i(x)\,(p_i-x)$, compare with lemma \ref{auxiliaireinfinit}), we know from the last equalities in (\ref{barycorr3}), for any $i$ and $x\in b\,$, that the tangent space $T_x(\partial B^g(p_i,d_i(x)))$ {\it is both $g$-orthogonal to $a_i(x)$ and $g'$-orthogonal to} $a'_i(x)=L_x^{-1}(a_i(x))\,$ (if $x\notin b\,$, then $a'_i(x)=a_i(x)$). As $x$ varies, the distribution of hyperplanes $\{v\in T_xM\mid
\langle v, a'_i(x)\rangle_{g'}=0\}=T_x(\partial B^g(p_i,d_i(x)))$ already integrates into the family of $g$-spheres
$\partial B^g(p_i,d_i(\bullet))\,$.
Moreover, the $g$-distance function $d_i$ to $p_i$ verifies $d_i(x)=\Vert a_i(x)\Vert_g/\lambda_i(x)\,$. As the $g$-spheres $\partial B^g(p_i,d_i(\bullet))\,$ are the level sets of $d_i\,$, the $g'$-gradient of $d_i$ is $g'$-perpendicular to $\partial B^g(p_i,d_i(\bullet))$, thus proportional to $a'_i(\bullet)\,$. Consider the curve $c(t)$ such that 
$c(0)\!=\!x $ and ${d c}/{dt}\!=\!{a'_i\circ c}/{\Vert a'_i\circ c\Vert_{g'}}\,$. 
The proof of lemma \ref{contre3} shows that $\varphi(x)\,B_x\!+\!{\rm Id}$ belongs to ${\mathcal D}_x$. Uses that the gradient of $\frac{1}{2}\,d_i^2$ is $a_i\,$ (\cite{B-K}, proposition 6.4.6) and (\ref{contre2}) $\langle a_i,a'_i\rangle_g\!=\!\langle a_i,a_i\rangle_g\,$ to get
\begin{multline}\label{contre4}  \frac{d}{dt}(d_i\circ c)\!=\! \frac{1}{(d_i\circ c)}\,\frac{d}{dt}
\frac{(d_i\circ c)^2}{2}\!=\!\frac{\lambda_i\circ c}{\Vert a_i\circ c\,\Vert_g}\,\langle \frac{a_i\circ c}{\lambda_i\circ c}, \frac{a'_i\circ c}{\Vert a'_i\circ c\Vert_{g'}}\rangle_g\!=\\=\langle \frac{a_i\circ c}{\Vert a_i\circ c\,\Vert_g}, \frac{a'_i\circ c}{\Vert a_i\circ c\Vert_g}\rangle_g=\langle \frac{a_i\circ c}{\Vert a_i\circ c\,\Vert_g}, \frac{a_i\circ c}{\Vert a_i\circ c\Vert_g}\rangle_g=1\  ,
\end{multline}
so, since $\Vert{d c}/{dt}\Vert_{g'}\!=\!1$, {\it the $g'$-gradient of $d_i$ is of constant $g'$-norm} $1\,$. Thus {\it the $d_i$-level $g$-spheres are $g'$-parallel hypersurfaces.} This implies that {\it the integral curves of the $g'$-gradient of $d_i$ are arc-$g'$-length parametrised $g'$-geodesics starting at $p_i\,$}, as we recall just below. 
\par 
Indeed, if $f$ is a function defined on a Riemannian manifold $(M,h)$ verifying $\Vert\nabla^hf\Vert_h=1\,$, the local one-parameter subgroup of diffeomorphism $\phi_s$ induced by $\nabla^hf$ preserves the ``level sets'' $f^{-1}(c)$ of $f\,$: this follows from
$$\frac{d}{ds}(f\circ\phi_s(x))=df(\phi_s(x)) (\frac{d\phi_s(x)}{ds})=\Vert\nabla_{\phi_s(x)}^hf\Vert^2_h=1\ .
$$
This implies (for any $x\in M$) 
\begin{equation}\label{chapeau!}f\circ\phi_s(x)-f\circ\phi_{s_0}(x)=s-s_0\ .
\end{equation}
So, given a field $X$ tangent to a ``level set'' $f^{-1}(c)$ ($c$ is a regular value), define the extended (in $M$) vector field of the type $X_{\phi_s(x)}=d \phi_s(x)(X)$, where $x$ is in $f^{-1}(c)\,$; one has $[X,\nabla^hf]=0$ and thus one computes
$$\langle D^h_{\nabla^hf}\nabla^hf,X\rangle_h=-\langle \nabla^hf, D^h_{\nabla^hf}X\rangle_h=-\langle \nabla^hf, D^h_X\nabla^hf\rangle_h=0\  ,
$$
Hence, an integral curve of $\nabla^hf$ is an $h$-geodesic and, in a convex region, it is a minimising arc-$h$-length parametrised geodesic.
\par\smallskip
\noindent\emph{Conclusion}
\par
{\it In both the distinct metrics} $g$ and $g'\,$, the function $d_i$ is the distance function to $p_i\,$, one has $d_i=d'_i\,$ (this is deduced from the equality $d_i\circ\phi_s(x)-d_i\circ\phi_{s_0}(x)=d'_i\circ\phi'_s(x)-d'_i\circ\phi'_{s_0}(x)=s-s_0\,$, see the above equality (\ref{chapeau!}) - of course, $\phi_s$ and $\phi_s'$ are the local one parameter subgroups of diffeomorphisms associated to $\nabla^g d_i$ and $\nabla^{g'} d'_i$). As this is true for any $i$, it also shows that any point $x$ has the same barycentric coordinates in both the distinct Riemannian metrics $g$ and $g'$ with respect to ${\bf p}=p_0,p_1,p_2,p_3\,$.
This completes the construction.
\end{proof}
\subsection{Completing the proof of proposition \ref{mg4}}\label{app2}
$ $
\begin{proof}
At $t=0\,$, we show the extensions are
$\ 
d_{0,x}^2(p,q)=\Vert\,v-w\,\Vert_x^2 \,$, where $p=\exp_xv,\ 
  q=\exp_xw\,$.
Taking twice the $t$-derivative of 
$\ 
t^2\,d_{t,x}^2(p,q)= \Vert\exp_{p(t)}^{-1}{q(t)}\Vert_{p(t)}^2=
\Vert\frac{\partial \gamma}{\partial s}(0,t)\Vert_{p(t)}^2\,
$,
one gets 
\begin{equation*} 
2\,t\,d_{t,x}^2(p,q)+t^2\,\frac{\partial}{\partial t}(d_{\cdot,x}^2(p,q))= 
2\,\langle\nabla_{\frac{\partial \gamma}{\partial s}}Y(0,t),
\frac{\partial \gamma}{\partial s}(0,t)\rangle_{p(t)}\  \  \ \ \ \ \hbox{and}
\end{equation*}
\begin{multline} \label{mg3}
\ \ \ \ 2\,d_{t,x}^2(p,q)+4\,t\,
\frac{\partial}{\partial t}(d_{\cdot,x}^2(p,q))+t^2\,
\frac{\partial^2}{\partial t^2}(d_{\cdot,x}^2(p,q))= \\
=2\,\langle\nabla_{\frac{\partial \gamma}{\partial t}}
\nabla_{\frac{\partial \gamma}{\partial s}}Y(0,t),
\frac{\partial \gamma}{\partial s}(0,t)\rangle_{p(t)}+
2\,\Vert\,\nabla_{\frac{\partial \gamma}{\partial s}}Y(0,t)\,\Vert_{p(t)}^2\  .
\end{multline}

The Jacobi field equation 
$\nabla_{\frac{\partial \gamma}{\partial s}}
\nabla_{\frac{\partial \gamma}{\partial s}}Y=
-R(Y,\frac{\partial \gamma}{\partial s})\frac{\partial \gamma}{\partial s}$ gives, for all $t\in]0,1]\,$, in a $g$-orthonormal parallel frame $X_i(s)$ 
along $s\mapsto \gamma(s,t)$
\begin{equation}\label{variation1033}
\langle\nabla_{\!\frac{\partial \gamma}{\partial s}}Y,X_i\rangle (s,t)\!=\!
-\!\!\!\int_0^s\!\!\langle R(Y,\frac{\partial \gamma}{\partial s})
\frac{\partial \gamma}{\partial s},X_i\rangle_{\!\gamma(\sigma,t)} d\sigma\!+\!
\langle\nabla_{\frac{\!\partial \gamma}{\partial s}}Y,X_i\rangle (0,t)\,.
\end{equation}
Everything for $t=0$ sits in $(T_xM,g_x)$, thus one infers from (\ref{variation1033}) and from $\frac{\partial \gamma}{\partial s}(s,0)=0$ (by definition of $\gamma$ (\ref{variation}))
$\lim_{t\rightarrow0}\nabla_{\!\frac{\partial \gamma}{\partial s}}Y(s,t)\!={\frac{\partial Y}{\partial s}}(s,0)=
\frac{\partial Y}{\partial s}(0,0)$, so ${\frac{\partial Y}{\partial s}}(s,0)$ is 
constant and 
$Y(s,0)\!-\!Y(0,0)\!=\!s{\frac{\partial Y}{\partial s}}(0,0)$. Set $s=1$ and get 
$Y(1,0)-Y(0,0)=\frac{\partial Y}{\partial s}(0,0)$, this gives 
\begin{equation*}
\frac{\partial Y}{\partial s}(0,0)=
\frac{\partial q}{\partial t}(0)-\frac{\partial p}{\partial t}(0)=w-v\  ,
\end{equation*}
and the conclusion we wanted, making $t=0$ in (\ref{mg3}).
\par To prove (\ref{mghess}), again
work with 
$\gamma(s,t) = 
\exp_{p(t)}(s\exp_{p(t)}^{-1}q(t)) $ (with $s,t\in [0,1]\,$, see definition \ref{vargamma})  and 
the associated Jacobi field  $Y_t(s)=\frac{\partial\gamma}{\partial t}(s,t)$.
Bringing in (\ref{mg-1}), we want to compute
\begin{multline}\label{mg10}
\frac{\partial }{\partial 
t}(\Vert\,\frac{\exp_{p(t)}^{-1}q(t)}{t}\,\Vert_{g_t}^2) 
=2\langle \frac{\exp_{p(t)}^{-1}q(t)}{t},
\nabla_{\frac{\partial\gamma}{\partial t}}(
\frac{\exp_{p(t)}^{-1}q(t)}{t})\rangle_{g_t}=\\
 =-2\,\langle \frac{1}{t}\frac{\partial\gamma}{\partial s}(0,t),
\frac{1}{t^2}(\frac{\partial\gamma}{\partial s}(0,t)
-t\,\frac{\partial Y_t}{\partial s}(0))\rangle_{g_t} \ .
\end{multline}
Now, perform an order $3$ Taylor expansion (in $t\,$, around $t=0$) 
of 
$\frac{\partial\gamma}{\partial s}(0,t)$ and 
$t\,\frac{\partial Y_t}{\partial s}(0)$. Choose a $g$-parallel 
${\mathcal Z}_i(t)$ orthonormal frame  along the geodesic $p(t)=\gamma(0,t)$. 
One gets
\begin{multline*}
\langle \frac{\partial\gamma}{\partial s}(0,t),{\mathcal Z}_i(t)\rangle_g\!=\!
t\,\langle \frac{\partial Y_0(0)}{\partial s},{\mathcal Z}_i(0)\rangle_g\! 
+\!\frac{t^2}{2}
\langle (\nabla_{\frac{\partial}{\partial t}}
\nabla_{\frac{\partial}{\partial s}}Y)_{\!(0,0)},{\mathcal Z}_i(0)\rangle_g \,+\\ 
+ t^3\int_0^1\frac{(1-\tau)^2}{2}\,
\langle (\nabla_{\frac{\partial}{\partial t}}
\nabla_{\frac{\partial}{\partial t}}\nabla_{\frac{\partial}{\partial s}}Y)_{(0,t\tau)}, 
{\mathcal Z}_i(t\tau)\rangle_g \  d\tau\  ,
\end{multline*}
and 
\begin{multline*}
t\,\langle \frac{\partial Y_t}{\partial s}(0,t),{\mathcal Z}_i(t)\rangle_g=
t\,\langle \frac{\partial Y_0(0)}{\partial s},{\mathcal Z}_i(0)\rangle_g 
+ t^2\, 
\langle (\nabla_{\frac{\partial}{\partial t}}
\nabla_{\frac{\partial}{\partial s}}Y)_{(0,0)}, {\mathcal Z}_i(0)\rangle_g + \\
+ t^3\int_0^1(1-\tau)\,
\langle (\nabla_{\frac{\partial}{\partial t}}
\nabla_{\frac{\partial}{\partial t}}
\nabla_{\frac{\partial}{\partial s}}Y)_{(0,t\tau)},
{\mathcal Z}_i(t\tau)\rangle_g \  d\tau\  ,
\end{multline*}
so that
\begin{multline}\label{Eq.Taylor}
\langle \frac{\partial\gamma}{\partial s}(0,t)
-t\,\frac{\partial Y_t}{\partial s}(0,t),{\mathcal Z}_i(t)\rangle_g 
=\\= t^3\int_0^1\frac{\tau^2-1}{2}\,
\langle 
(\nabla_{\frac{\partial}{\partial t}}\nabla_{\frac{\partial}{\partial t}}
\nabla_{\frac{\partial}{\partial s}}Y)_{(0,t\tau)},
{\mathcal Z}_i(t\tau)\rangle_g\  d\tau\  .
\end{multline}
{\it Indeed,} a non-obvious cancellation occurs for this variation 
$\gamma\,$: the $t^2$-term in the Taylor expansion (\ref{Eq.Taylor})
vanishes with $ 
(\nabla_{\frac{\partial}{\partial t}}
\nabla_{\frac{\partial}{\partial s}}Y)_{(0,0)}\,$.
\begin{lem} \label{Topono}The variation $\gamma$ verifies
\begin{equation} \label{Eq.t2} 
(\nabla_{\frac{\partial}{\partial t}}
\nabla_{\frac{\partial}{\partial s}}Y)_{(s,0)}\equiv 0\  .
\end{equation}
\end{lem}
The proof below is taken from Shlomo Sternberg's book \cite[section 10.1, pages 217ff]{St}, who refers to Wolfgang Meyer 
 \cite{MeW}, see section 1.3, pages 5-8, in his text dedicated to a theorem of Toponogov. 
 \begin{proof}Indeed, we know that 
$(\nabla_{\frac{\partial}{\partial t}}Y)_{(0,0)}=
(\nabla_{\frac{\partial}{\partial t}}Y)_{(1,0)}=0\,$, 
the curves 
$p(t)$ and $q(t)$ being geodesics by definition of $\gamma\,$. It will be 
proved just below that 
$s\mapsto(\nabla_{\frac{\partial}{\partial t}}Y)_{(s,0)}$ is an affine function 
of $s\,$. Both facts give together
$(\nabla_{\frac{\partial}{\partial t}}Y)_{(s,0)}\equiv0\,$. So, we get 
$(\nabla_{\frac{\partial}{\partial s}}
\nabla_{\frac{\partial}{\partial t}}Y)_{(s,0)}\equiv 0$ and the conclusion 
$(\nabla_{\frac{\partial}{\partial t}}
\nabla_{\frac{\partial}{\partial s}}Y)_{(s,0)} \equiv 0$ 
follows then from $(R(\frac{\partial p(\cdot)}{\partial 
t},\frac{\partial\gamma}{\partial s})Y)_{t=0}\equiv0\,$, which holds 
for $t=0$ because $\frac{\partial\gamma}{\partial s}_{t=0}\equiv0$ by 
definition 
of the variation.

It remains to show that 
$s\mapsto(\nabla_{\frac{\partial}{\partial t}}Y)_{(s,0)}$ is an affine function.
Take the second derivative of this function 
 with respect to $s$ 
and check that the result is in fact identically $0\,$.
\begin{multline*}
\nabla_{\frac{\partial}{\partial s}}\nabla_{\frac{\partial}{\partial s}}
\nabla_{\frac{\partial}{\partial t}}Y
=\nabla_{\frac{\partial}{\partial s}}R(\frac{\partial \gamma}{\partial s},Y)Y
+\nabla_{\frac{\partial}{\partial s}}
\nabla_{\frac{\partial}{\partial t}}\nabla_{\frac{\partial}{\partial s}}Y=\\
=\nabla_{\frac{\partial}{\partial s}}(R(\frac{\partial \gamma}{\partial s},Y)Y) 
+ \nabla_{\frac{\partial}{\partial s}}
\nabla_{\frac{\partial}{\partial t}}
\nabla_{\frac{\partial}{\partial t}}\frac{\partial \gamma}{\partial s} = \\
=\nabla_{\frac{\partial}{\partial s}}(R(\frac{\partial \gamma}{\partial s},Y)Y)+
R(\frac{\partial \gamma}{\partial s},\frac{\partial \gamma}{\partial t})
(\nabla_{\frac{\partial }{\partial t}}\frac{\partial \gamma}{\partial s})
+\nabla_{\frac{\partial}{\partial t}}\nabla_{\frac{\partial}{\partial s}}
\nabla_{\frac{\partial}{\partial t}}\frac{\partial \gamma}{\partial s} = \\
=\nabla_{\frac{\partial}{\partial s}}(R(\frac{\partial \gamma}{\partial s},Y)Y) 
+
R(\frac{\partial \gamma}{\partial s},\frac{\partial \gamma}{\partial t})
(\nabla_{\frac{\partial }{\partial t}}\frac{\partial \gamma}{\partial s})
+\nabla_{\frac{\partial}{\partial t}}(R(\frac{\partial \gamma}{\partial s},
\frac{\partial \gamma}{\partial t})\frac{\partial \gamma}{\partial s})+\\+
\nabla_{\frac{\partial}{\partial t}}\nabla_{\frac{\partial}{\partial t}}
\nabla_{\frac{\partial}{\partial s}}\frac{\partial \gamma}{\partial s}\  \,.
\end{multline*}
The last term of the last expression obtained 
vanishes for 
$s\mapsto\gamma(s,t)$ are geodesics and the three others vanish too at 
$t=0\,$, since all 
these terms, when written in 
fully tensorial expressions, are all linear in 
$\frac{\partial \gamma}{\partial s}(s,0)\equiv0$ 
(since $\gamma(s,0)\equiv x$) or 
$\nabla_{\frac{\partial \gamma}{\partial s}(s,0)}
\frac{\partial \gamma}{\partial s}(s,0)=0\,$.
Thus (\ref{Eq.t2}) holds.
\end{proof}
In view of (\ref{Eq.Taylor}), the last term in
(\ref{mg10}) 
is seen to tend 
to $0$ with $t$, this shows 
$\ \frac{\partial }{\partial t}(d_{t,x}^2(p,\cdot)_{\mid t=0}\equiv 0\,$.
\end{proof}

\subsection{Proof of lemma \ref{mg7}}\label{app20}
$ $
\begin{proof} 
Choose $u\in T_{(t,q)}B_t$ and extend it around in a field $U\,$, first in $B_t\,$, 
and then in $\tilde B$ by defining $U$ as the constant field in the 
$\frac{\partial }{\partial t}$ direction generated 
by $U$ in $B_t$ (thus $U=\psi_{\tau-t}(U_{B_t})$ remains tangent to all slices $B_\tau$). By construction, one 
has 
$[\frac{\partial }{\partial t},U]=0\,$. One also has
\begin{equation*}
\langle S_t(u),u \rangle_{g_t} = \langle\tilde\nabla_u
\frac{\partial }{\partial t},u\rangle_{\tilde g} = \frac{1}{2}\frac{\partial 
}{\partial t}
\langle U,U\rangle_{\tilde g} = 
\frac{1}{2}\frac{\partial }{\partial t} g_t(U,U)\,.
\end{equation*}
Thus, we want to compute  
$\nabla_{\frac{\partial }{\partial t}}(\frac{d\psi_t(q)(u)}{t})$ 
at $t=0\,$: 
if ${\mathcal P}^{-1}$ is the parallel translation along the 
geodesic ray from 
$q$ to $x\,$, this leads us to find
\begin{equation} \label{mglim}
\lim_{t\rightarrow0}\frac{({\mathcal P}^{-1}\circ d\exp_x(t\exp_x^{-1}(q))-{\rm Id})\circ 
d\exp_x^{-1}(q)(u)}{t}\ .
\end{equation}
As $S_0$ only matters here, the radius of the ball
$B(x,\rho)$
may be chosen small {\it ad libitum} while applying lemma \ref{lemma 2nd} to 
$\ 
c(\sigma,\tau)\!=\!\exp_x(\sigma(\exp_x^{\!-1}q\!+\!\tau\,(d\exp_x^{\!-1}(q)(u))))\,
$.
The Jacobi field along $c(\sigma,0)$ is
$$Y(\sigma)\!=\!\frac{\partial c}{\partial 
\tau}(\sigma)\!=\!d\exp_x(\sigma(\exp_x^{\!-1}q))(d\exp_x^{\!-1}(q)(u))
$$
and one gets the following estimate (in fact, $\sigma$ plays here the role 
attributed to 
$s$ in the statement of lemma \ref{lemma 2nd}, while $\tau=0$)
\begin{multline}\label{Eq.sigma}
\Vert({\mathcal P}^{-1}\!\circ d\exp_x(\sigma \exp_x^{-1}(q))\!-
\hbox{\rm Id}_x) 
(d\exp_x^{-1}(q)(u))\Vert_x  
\leq\\\leq
{\mathcal C}_2\,\Vert d\exp_x^{-1}(q)(u)\Vert_x\  
\Vert\sigma\exp_x^{-1}(q)\Vert_x^2\,.
\end{multline}
Making $t=\sigma$ in (\ref{Eq.sigma}) shows that the limit (\ref{mglim}) is $0\,$. One gets the conclusion $S_0\equiv0$ and
$\frac{\partial g_t}{\partial t}(0)\equiv0\,$.
\end{proof}

\subsection{Proof of lemma \ref{mgcom1}}\label{app22}
$ $
\begin{proof} 
First, consider the case $T=X$ is a tangential vector field on $\tilde B\,$, thus $\langle X, \frac{\partial }{\partial t}\rangle=0\,$. At a value $t_0$ where the shape operator 
$S_{t_0}$ vanishes identically, we want to prove
\begin{equation*} 
(\tilde{\nabla}_{\frac{\partial }{\partial t}}\nabla^tX)_{(t_0,q_0)}=
(\nabla^t\tilde{\nabla}_{\frac{\partial }{\partial t}}X)_{(t_0,q_0)}\,.
\end{equation*}
Take any field $U$ tangent to $B_{t_0}$ and extend it in a $\tilde g$-parallel 
field $U$ (on $\tilde B$) along the geodesic integral lines of 
$\frac{\partial }{\partial t}\,$, 
so $U$ remains tangent to 
all slices $B_t\,$.
At $t=t_0\,$, one has $(\nabla_U\frac{\partial }{\partial t})_{t_0}=S_{t_0}(U)=0$ 
and, 
by construction, one also has 
$(\nabla_{\frac{\partial }{\partial t}}U)_{t_0}=0$, so one has 
$[U,\frac{\partial }{\partial t}]_{B_{t_0}}=0\,$.  

Thus  
$(\tilde{\nabla}\tilde{\nabla}_{\frac{\partial }{\partial t}}X)(U)\!=\!\tilde{\nabla}_U\tilde{\nabla}_{\frac{\partial }{\partial t}}X$ 
and $(\tilde{\nabla}_{\frac{\partial }{\partial t}}\tilde{\nabla}X)(U)\!=\!\tilde{\nabla}_{\frac{\partial }{\partial t}}\tilde{\nabla}_UX$ have in fact equal tangential parts to $B_{t_0}$, 
up to the 
tangential part of 
$\tilde R(\frac{\partial }{\partial t},U)X\,$. But, due to the fact that 
$S_{t_0}=0\,$, one 
deduces 
from Codazzi's equation that precisely this tangential part of 
$\tilde R(\frac{\partial }{\partial t},U)X$ vanishes along $B_{t_0}\,$.
Now, the tangential part to $B_{t_0}$ of 
$\tilde{\nabla}_U\tilde{\nabla}_{\frac{\partial }{\partial t}}X$ 
is 
equal to 
$\nabla^t_U\tilde{\nabla}_{\frac{\partial }{\partial t}}X$ (actually, 
$\tilde{\nabla}_U\tilde{\nabla}_{\frac{\partial }{\partial t}}X$ 
is tangential because $\langle\tilde{\nabla}_{\frac{\partial }{\partial t}}X, 
\frac{\partial }{\partial t}\rangle=0$ and $S_{t_0}=0$). On the other hand, the 
tangential part to $B_{t_0}$ of 
$\tilde{\nabla}_{\frac{\partial }{\partial t}}\tilde{\nabla}_UX$ coincides with 
$\tilde{\nabla}_{\frac{\partial }{\partial t}}\nabla^t_UX\,$, since 
$\nabla^t_UX$ is the 
tangential part of $\tilde{\nabla}_UX\,$. The proof is 
complete in case of a tangential field $X\,$.

Any tangential $1$-differential form $\omega$ is canonically associated to a (dual) tangential field $X\,$, namely
$X$ is the field such that $g(X,\cdot)=\omega\,$. One has
for tangential $1$-differential form and fields $\omega,U,W$
\begin{multline*} (\tilde{\nabla}_{\frac{\partial}{\partial t}}\nabla^t \omega)(U,W)=\\=
\frac{\partial}{\partial t}((\nabla^t \omega)(U,W))-(\nabla^t \omega)(\tilde{\nabla}_{\frac{\partial}{\partial t}}U,W)-(\nabla^t \omega)(U,\tilde{\nabla}_{\frac{\partial}{\partial t}}W)=\\=\frac{\partial}{\partial t}(U(\omega(W)))-\frac{\partial}{\partial t}(\omega(\nabla^t_UW))-\nabla^t_{\tilde{\nabla}_{\frac{\partial}{\partial t}}U}(\omega(W))+\\+\omega(\nabla^t_{\tilde{\nabla}_{\frac{\partial}{\partial t}}U}W)-
\nabla^t_U(\omega(\tilde{\nabla}_{\frac{\partial}{\partial t}}W))+\omega(\nabla^t_U\tilde{\nabla}_{\frac{\partial}{\partial t}}W)\,,
\end{multline*}
while
\begin{multline*} (\nabla^t\tilde{\nabla}_{\frac{\partial}{\partial t}} \omega)(U,W)=U((\tilde{\nabla}_{\frac{\partial}{\partial t}}\omega)(W)))-
(\tilde{\nabla}_{\frac{\partial}{\partial t}}\omega)(\nabla^t_UW)=
\\=U(\frac{\partial}{\partial t}(\omega(W)))-U(\omega(\tilde{\nabla}_{\frac{\partial}{\partial t}}W))-\\
-\frac{\partial}{\partial t}(\omega(\nabla^t_UW))+\omega(\tilde{\nabla}_{\frac{\partial}{\partial t}}\nabla^t_UW)\,.
\end{multline*}
Choosing $U$ and $W$ to be $\frac{\partial}{\partial t}$-parallel (thus tangent to all leaves $B_t$), the two equations above become simpler. The first reads  
\begin{equation}\label{mgcom22} (\tilde{\nabla}_{\frac{\partial}{\partial t}}\nabla^t \omega)(U,W)\!=\!\frac{\partial}{\partial t}(U(\omega(W)))-\frac{\partial}{\partial t}(\omega(\nabla^t_UW))\,,
\end{equation}
and the second
\begin{equation} \label{mgcom33}(\nabla^t\tilde{\nabla}_{\!\frac{\partial}{\partial t}} \omega)(U,W)\!=\!U(\frac{\partial}{\partial t}(\omega(W)))\!-\!\frac{\partial}{\partial t}(\omega(\nabla^t_{\!U}W))\!+\omega(\tilde{\nabla}_{\!\frac{\partial}{\partial t}}\tilde{\nabla}_{\!U}W)\,.
\end{equation}
At any point $(t_0,q_0)\in B_{t_0}\,$, since $S_{t_0}\equiv0\,$, one gets $[\frac{\partial}{\partial t},U]_{(t_0,q_0)}=0\,$: the two first terms of the right-hand sides in (\ref{mgcom22}) and (\ref{mgcom33}) are equal. The second are obviously equal. Again, $S_{t_0}\equiv0$ and the Codazzi equation imply $\tilde{\nabla}_{\frac{\partial}{\partial t}}\tilde{\nabla}_UW=0\,$. The statement is established in case $T=\omega\,$.

If $T$ is a tangential tensor of type $(l,m)\,$, the operators $\nabla^t,\,\tilde{\nabla}_{\frac{\partial}{\partial t}}$ being derivations, the result follows from a straightforward computation.
\end{proof}

\subsection{Proof of proposition \ref{mgapp1}}\label{app222}
$ $
\begin{proof} Taylor's expansion with integral remainder implies that $g=f/t$ is $C^{k-1}\,$. 
If $l\!=\!0,1,\dots, k-1,\,t\!\not=\!0\,$, {\it set}
$f_l (t)\!:=\! t^{l+1} \,\frac{d^l g}{d t^l}(t)\, 
$. One has
 
\begin{equation}\label{mgapp2}
f_0= f\  , 
\ f_1=t\,\frac{d f}{d t}-f\  ,\  \dots\  ,\  f_l= t \,\frac{d f_{l-1}}{d t}-l\,f_{l-1}\  .
\end{equation}
This implies immediately $f_l(0)=0$ for $l=0,1,\dots, k-1\,$, which follows by induction 
from $f(0)=f_0(0)=0\,$.
If $l=0\,$, the inequality to be proved results from 
\begin{equation*} 
g(t)=\frac{1}{t}\,\int_0^t \frac{d f}{d t}(s)\,ds \ \ \  
\hbox{which gives}\ \ \ 
\sup_{t\in I} \,\vert g (t) \vert\leq \sup_{t\in I}\, \vert \frac{d f}{d t} (t)\vert\  .
\end{equation*}
Assume the induction formula below (which holds if $l=1$) to be true for $l\leq k-2$
\begin{equation}\label{mgapp3} 
\frac{d f_l}{d t}(t)=t^l\,\frac{d^{l+1} f}{d t^{l+1}} (t)\ .
\end{equation}
A computation gives, starting from (\ref{mgapp2}) 
and inserting 
the 
induction formula (\ref{mgapp3}) for $l$ to get the same for $l+1$
\begin{multline*} 
\frac{d f_{l+1}}{d t}(t)=\frac{d}{dt}(t\,\frac{d f_l}{d t}(t)-(l+1)\,f_l(t)) 
= \frac{d}{dt}(t^{l+1}\,\frac{d^{l+1} f}{d t^{l+1}} (t))-(l+1)\,\frac{d f_l}{d t}(t) = \\
=(l+1)\,t^l\,\frac{d^{l+1} f}{d t^{l+1}} (t)-(l+1)\,\frac{d f_l}{d t}(t)+t^{l+1}\,
\frac{d^{l+2} f}{d t^{l+2}} 
=t^{l+1}\,\frac{d^{l+2} f}{d t^{l+2}} \  .
\end{multline*}
Now, with the help of (\ref{mgapp2}), (\ref{mgapp3}) and $f_{l-1}(0)=0\,$, write 
\begin{multline*}
f_l(t)=t^l\,\frac{d^l f}{d t^l} (t)-l\,f_{l-1}(t)=\\=(l\int_0^ts^{l-1}ds)
(\frac{d^l f}{d t^l}) (t)-l\,\int_0^ts^{l-1}(\frac{d^l f}{d t^l})(s)ds = \\
=l\int_0^ts^{l-1}(\frac{d^l f}{d t^l}(t)-\frac{d^l f}{d t^l}(s)) ds 
=l\,\int_0^ts^{l-1}(t-s)(\frac{\frac{d^l f}{d t^l}(t)-\frac{d^l f}{d t^l}(s)}{t-s}) ds\  .
\end{multline*}
Thus, one has
\begin{equation*}
\vert \frac{d^l g}{d t^l}(t)\vert=\vert \frac{f_l (t)}{t^{l+1}}\vert\leq l\,
\bigg{\vert}\int_0^t\vert\frac{s^{l-1}(t-s)}{t^{l+1}}\vert ds\bigg{\vert}\,
\sup_I\,\vert \frac{d^{l+1} f}{d t^{l+1}}(t)\vert\  ;
\end{equation*}
computing the last integral, the result follows.
\end{proof}

\begin{cor} Let $f$ be a $C^k$ function of $t$ varying in $I=[-T,T]$ such that 
$f(0)=\frac{d f}{d t}(0)=0\,$,
with $T>0$ and $k\geq2\,$. The function $h(t)=f(t)/t^2$ is $C^{k-2}$ and the derivatives of 
$h$ 
in $I$ are bounded in terms of those of $f\,$, more precisely, for any $l$ with 
$0\leq l\leq k-2\,$, one has
\begin{equation*} 
\sup_{t\in I} \,\vert \frac{d^l h}{d t^l} (t) \vert\leq 
\frac{1}{(l+1)(l+2)}\,\sup_{t\in I}\, \vert \frac{d^{l+2} f}{d t^{l+2}} 
(t)\vert\  .
\end{equation*}
\end{cor}

\begin{proof} 
Apply proposition \ref{mgapp1} to $g(t)\!=\!f(t)/t$ and $h(t)\!=\!g(t)/t\,$.
\end{proof}

\subsection{Proof of lemma \ref{2q}}\label{preuveMunkres}
$ $
\begin{proof} As $E\geq E({\mathcal R}_4)$, given any $\sigma_E\in K_E\,$, one has for any $x\in T(\sigma_E)$ the inclusion $T(\sigma_E)\subset B(x,{\mathcal R}_4)$. As $T_x$ is well defined from any such $\sigma_E\subset T^{-1}(B(x,{\mathcal R}_4))$ into $B(0_x,{\mathcal R}_4)\subset T_xM$, so is the secant map $L_x$ too, from the same $\sigma_E$ into $T_xM\equiv {\mathbb R}^n\,$. One has for any $y,z\in\sigma_E$ (using the mean value equality on the linear arc from $y$ to $z$)
\begin{equation}\label{Mun1}
\Vert T_x(z)-T_x(y)-dT_x(y)(z-y)\Vert_x\leq {\mathscr T}_x(E) \Vert z-y\Vert_{\rm std}\  .
\end{equation}
One also has, for any $y\in\sigma_E$ and
$u\in T_y\sigma_E$
\begin{equation}\label{Mun2}
\Vert dT_x(z)(u)-dT_x(y)(u)\Vert_x\leq {\mathscr T}_x(E)\ \Vert u\Vert_{\rm std}\  .
\end{equation}
Observe that $\hbox{diam}_{\rm std}(\sigma_E)\leq 1$ if the standard metric ${\rm std}$ on the initial $K$(before refinement) verifies that any regular $n$-simplex of $K$ has diameter $1$. If ${\rm std}$ is homothetic to the previous by the factor $\sqrt{2}$, then $\hbox{diam}_{\rm std}(\sigma_E)\leq 1$ is true if $E\geq 2$. Anyway, we assume $\hbox{diam}_{\rm std}(\sigma_E)\leq 1$ to be true. Thus {\em the map $T_x$ is a strong 
${\mathscr T}_x(E)$-approximation (on $\sigma_E$) to the {\em affine} tangent map $D_y(T_x)(\cdot):=T_x(y)+dT_x(y)(\cdot-y)$.}
\par
Moreover, the secant map $L_x$ verifies $L_x(u_i)=T_x(u_i)$ on each vertex $u_i$ and one has for $i,j\in\{0,1,\cdots,n\}\,$ (as in (\ref{Mun1}))
\begin{equation}\label{Mun3}
\Vert L_x(u_j)-L_x(u_i)-dT_x(y)(u_j-u_i)\Vert_x\leq  {\mathscr T}_x(E)\ \Vert u_j-u_i\Vert_{\rm std}\  .
\end{equation}
If $u=\sum_{i=1}^n\lambda_i\,(u_i-u_0)$, get $\sum_{i=1}^n\vert\lambda_i\vert\,\hbox{diam}_{\rm std}(\sigma_E)\leq \sqrt{n}\,\Vert u\Vert_{\rm std}/t_{\rm std}(\sigma_E)\,$ from lemma \ref{1b}. So (\ref{Mun3}) gives for any $y\in\sigma_E\,,\, u\in T_y\sigma_E$ ($L_x$ is affine)
\begin{equation}\label{Mun4}
\Vert (dL_x(y)-dT_x(y))(u)\Vert_x\leq \frac{\sqrt{n}}{t_{\rm std}(\sigma_E)}\ {\mathscr T}_x(E)\ \Vert u\Vert_{\rm std}\  .
\end{equation}
Finally, in view of the inequality (true, by (\ref{Mun1}), for any $i$) 
$$\Vert L_x(u_i)-D_yT_x(u_i)\Vert_x=\Vert T_x(u_i)-D_yT_x(u_i)\Vert_x\leq 
{\mathscr T}_x(E)\ \hbox{diam}_{\rm std}(\sigma_E)\  ,
$$
writing $z=\sum_{i=0}^n\mu_i\,u_i\,$ with $\sum_{i=0}^n\mu_i=1$ and all $\mu_i\geq 0\,$, one has
\begin{multline}\label{Mun5}\Vert L_x(z)-D_yT_x(z)\Vert_x=
\Vert \sum_i\mu_i\  (L_x(u_i)-D_yT_x(u_i))\Vert_x\leq\\
\leq 
\sup_i \Vert L_x(u_i)-D_yT_x(u_i)\Vert_x\leq {\mathscr T}_x(E)\ \hbox{diam}_{\rm std}(\sigma_E)\ .
\end{multline}
Then, as $\hbox{diam}_{\rm std}(\sigma_E)\leq 1\,$, (\ref{Mun4}) and (\ref{Mun5}) tell that {\em $L_x$ is a strong $\frac{\sqrt{n}}{{\bf t}_n}\,{\mathscr T}_x(E)$-approximation to $D_yT_x$ on $\sigma_E\,$.} Since $1\leq\frac{\sqrt{n}}{{\bf t}_n}\,$, the conclusion follows in view of definition (\ref{constant99}) and the remark after.
\end{proof}

\subsection{Proof of theorem \ref{hatT}}\label{Whitehead}
$ $ 

Here $E$ is chosen such as $\rho\leq{\mathcal R}_4/4$ holds.
We shall prove:
\par {\it For any 
$\delta>0$ small enough exists $E_1(\delta)\!\in\!\N^\ast$ such that, for any 
integer $E\geq E_1(\delta)\,$, the
mapping 
$\hat T_E : (K_E,{\rm std}) \rightarrow (M,g)$ is a simplicial embedding and a $2\delta$-approximation to $T$, one has (see ({\rm\ref{soclosetoyou}}))}
\begin{equation}\forall p\in \lvert K \lvert\ \ \ \ \ \ \ \ \ \ \ \ \ \ 
 d_g(\hat T_E (p), T(p)) \leq 2\,\delta \  .
\end{equation}

Proposition \ref{James2} below makes the bridge between theorem \ref{James1} and our situation. First, some objects. 

\begin{defn}\label{D,covering}$ $
With the notations and definitions of section \ref{6.1},  
we enforce here the requirement on the covering done in the proof of theorem \ref{2n} and cover this time $\bar W$ by a finite number of balls $B(x_i,{\mathcal R}_4)$ such 
that the balls $B(x_i,{\mathcal R}_4/4)$ still cover $\bar W$, where all $x_i$ are in $ W$.
Recall or define various sub-complexes of $K_E$

\begin{itemize} 
\item $K_E^i$ is constituted by
all $n$-simplices $\sigma_E$ of $K_E$ such that 
$T(\sigma_E)\subset B(x_i,{\mathcal R}_4)$ plus all their faces;

\item $K_{E,i}$ is constituted by all 
$n$-simplices 
$\sigma_E$ such that 
$T(\sigma_E)\cap B(x_i,{\mathcal R}_4/2)\not=\emptyset$ plus all their faces;

\item $k_{E,i}$ is constituted by all  
$n$-simplices 
$\sigma_E$ such that 
$T(\sigma_E)\cap B(x_i,{\mathcal R}_4/4)\not=\emptyset$ plus all their faces.
\end{itemize}

\end{defn}

\begin{rem}\label{batifol} 
If $\sigma_E\,\in\! K_{E,i}\,$, one has $T(\sigma_E)\!\subset\! B(x_i,{\mathcal R}_4)$. On the other hand, for any $\sigma_E\in K_E$ exists an $ i\in{\mathbb N}$ such that $T(\sigma_E)\subset B(x_i,{\mathcal R}_4/2)\  $.

{\it Indeed},
given $i$, 
observe that, by definition, any $\sigma_E\in K_{E,i}$ is such that 
$T(\sigma_E)$ has a point $x$ belonging to $B(x_i,{\mathcal R}_4/2)$. By the choice of $\rho$  (implying $\rho \leq {\mathcal R}_4/4\,$) and of $E=E(\rho)$, one has  
 $\diam (T(\sigma_E)) < \rho$ (proposition \ref{2o}$(i)$). Thus $T(\sigma_E)\subset B(x,{\mathcal R}_4/4)\subset B(x_i,3{\mathcal R}_4/4)$. 
 
{\it As for the second claim},
$T(\sigma_E)$ meets a ball of the covering $B(x_i,{\mathcal R}_4/4)$, thus the conclusion since $\diam (T(\sigma_E)) < \rho$ and
$\rho \leq {\mathcal R}_4/4\,$.

Since  $\rho \leq {\mathcal R}_4/4\,$, it follows from the lines above that one has
\begin{equation*}
 k_{E,i}\subset K_{E,i}\subset K_E^i\subset K_E\  \  \  \hbox{and}\  \  \  
 \cup_i k_{E,i}=\cup_i K_{E,i}=\cup_i K_E^i=K_E\  .
\end{equation*}
\end{rem}

\begin{defn} \label{entier22} For a given $\rho_0 :=\min (\beta_T\,{\mathcal D}_n, {\mathcal R}_4/32)$, as in definition \ref{linkint} define the integer $E(\rho_0)$ which is the integer part of $1+(\beta_T/{\mathcal D}_n\rho_0)$.
\end{defn}
\begin{prop} \label{James2} 
Consider the differentiable simplicial embedding $f=T:K\rightarrow W\subset (M,g)$. 
There exists $\delta\in \,]0,{\mathcal R}_4/32]$ small enough (depending on $f$)
for which, given any integer  $E\geq E(\rho_0)$, 
the property 
\par {\em a map $h : K_E \rightarrow M$ such that, for any $i\,$,
the map $h_i=\exp_{x_i}^{-1}\circ h$ is a strong $\delta$-approximation to 
$f_i=T_{x_i}=\exp_{x_i}^{-1}\circ T:K_{E,i}\rightarrow(T_{x_i}M,g_{x_i})$}
\par implies the property
\par {\em the map $h\,$ itself is a differentiable simplicial embedding into $M$}
\par (for the notion of a strong approximation, see definition \rm{\ref{Argh55}})
\end{prop}

\begin{proof} 
First, a lemma.
\begin{lem}\label{L.Approx}
Any map $h :K_E \longrightarrow M$ satisfying all hypotheses given in proposition 
{\rm\ref{James2}} is a $2\delta$-approximation of $T$, one has
\begin{equation*}
 \forall p\in \lvert K \lvert\ \ \ \ \ \ \ \ d_g(h(p), T(p)) < 2 \delta\ .
\end{equation*}
\end{lem}

\begin{proof}
By proposition \ref{proposition A}, for any $x\in W$, $\exp_x$ is $\frac{1}{2}$-quasi-isometric on 
$B(0_x, {\mathcal R}_4) \!\subset \!T_x M$, so one has 
\begin{equation*}
\forall i\,, \ \forall u,v \in B(0_{x_i}, {\mathcal R}_4) \subset 
T_{x_i}M\,,\ \ \  d\left( \exp_{x_i} (u), \exp_{x_i} (v)\right) \leq 
2 \Vert v-u\Vert_{x_i}\ .
\end{equation*} 
The claim follows, exploiting $K_E= \cup_i K_{E,i}$ and going back to the definitions of 
$h_i$ and $T_i$ in terms of $h$ and $T$ through $\exp_{x_i}\,$.  
\end{proof}

We proceed to the proof of proposition \ref{James2}. Choose an integer $E\geq E(\rho_0)\,$.
Theorem \ref{James1} tells us that, for any $i$, 
there exists $\delta_i>0$ 
such that, for  any integer  $E\geq E(\rho_0)$, 
any strong $\delta_i$-approximation $h_i$ to $f_i=T_{x_i}$ on $K_{E,i}$ 
is also an embedding into 
$T_{x_i}M\,$. So $\exp_{x_i}\circ\, h_i : K_{E,i} \rightarrow M$, which is $h$ 
restricted 
to $K_{E,i}\,$, is also an 
embedding, as $\hbox{rge} (h_i)$ is included in the ball $ B(x_i,2{\mathcal R}_4)$ 
(recall ${\mathcal R}_4\leq{\mathcal R}_W$). 

Set $\delta_0 :=\min_i ({\mathcal R}_4/32 ;\, \delta_i)$ and choose 
$\delta \in \, ]0,\delta_0]\,$. Consider 
the differentiable 
simplicial mappings $h$ such that, for any $i$, the mapping $h_i$ is a 
strong $\delta$-approximation 
to $T_{x_i}$ on $K_{E,i}$ and thus an embedding from $K_{E,i}$ into $(T_{x_i}M,g_{x_i})$ 
by Whitehead's theorem \ref{James1}.

 As $E$ is $\geq E(\rho_0)\,$, any $n$-simplex $\sigma_E$ in $k_{E,i}$ verifies 
${\rm diam}_g \,T(\sigma_E) < \rho_0$ (proposition \ref{2o}) and the choice (see definition \ref{entier22})
$\rho_0 \leq {\mathcal R}_4/32$ gives the inclusion
$T(\sigma_E)\subset B(x_i,9{\mathcal R}_4/32)$.

 We now analyse how $h$ may fail to be injective. Take any $\sigma_E\in K_E$, thus 
$\sigma_E$ is in some $k_{E,i}\,$. Suppose that $\sigma'_E$ is a simplex in $K_E\,$, 
verifying 
$h(\sigma_E)\cap h(\sigma'_E)\not=\emptyset\,$. Consider a point $y \in h(\sigma_E)\cap h(\sigma'_E)$ and set 
$y= h(p)= h(p^\prime)$ with $p\in \sigma_E$ and $p^\prime\in \sigma_E^\prime\,$.
In view of the choice 
$\delta\leq {\mathcal R}_4/32\,$, we get (apply  
lemma \ref{L.Approx}) $d(y, T(p))< {\mathcal R}_4/16$ and 
$d(y, T(p^\prime))< {\mathcal R}_4/16\,$, implying 
$T(\sigma'_E) \cap B(x_i, 7{\mathcal R}_4/16)\not= \emptyset\,$. Thus $\sigma_E$ and 
$\sigma'_E$ are both simplices in $K_{E,i}\,$. As $\delta\leq {\mathcal R}_4/32\,$, one even has $h(\sigma_E)$ and $h(\sigma'_E)\subset B(x_i,{\mathcal R}_4)$.

So, any differentiable simplicial map $h$ as in proposition \ref{James2} verifies that $\exp_{x_i}\circ h_i$, for 
any $i\,$, is a differentiable simplicial 
embedding of $K_{E,i}$ into 
$M$; moreover, by the above lines, it is also injective. Being a local embedding
which is injective, it is embedding.
\end{proof}

To prove theorem 
\ref{hatT}, we make some further steps.
\begin{defn}\label{new map} {\it Introduce the new map} $\vartheta$ (given $(x,q)\in W\times B(x,2\rho)$, it sends an open set $\subset T_{v'}T_xM\equiv T_xM\,$, with $v'\in T_xM\,$, into $T_qM$) 
\begin{gather*}\vartheta : W\times B(x,2\rho)\times B(0_x,2\rho)\longrightarrow T_qM\hskip5mm \\
\hskip2cm (x,q,w) \longmapsto \vartheta(x,q,w)=(Q_{x,q}-d\exp_x(\exp_x^{-1}(q)))(w)\ ,\\
\hbox{where} \ \ \ Q_{x,q}=\exp_q^{-1}\circ\exp_x\ 
\end{gather*}
and the corresponding upper-bound $\Theta$ on the $C^2$-norm of $\vartheta$ achieved over the set $\overline{W\times B(x,2\rho)\times B(0_x,2\rho)}\,$ (where $\vartheta$ is well-defined and $C^2$-regular).
\end{defn}
\begin{defn}
Recall that $\phi$ (see definition \ref{D.euclidianisable}, section \ref{S.close}) is the barycentric map which sends the 
linear simplex 
$\hat\tau_x$ in $(T_xM,g_x)$ built on the vertices  
$v_0\!=\!\exp_x^{-1}p_0\,,\dots ,
v_n=\exp_x^{-1}p_n$ 
to the Riemannian barycentric simplex $\hat\tau$ built on the vertices $p_0,\dots,p_n$ 
(in a convex ball $B(x, {\mathcal R}_W)\subset(M,g)$), and {\em define} $\phi_x$ to be the 
mapping 
$\exp_x^{-1}\circ\, \phi$ (sending a neighborhood of $\hat\tau_x\subset(T_xM,g_x)$ into 
$(T_xM,g_x)$).
\end{defn}

\begin{prop}\label{Argh33}  
Given $t_0$ as in definition {\rm{\ref{2j}}}, one can find $\rho_1 (\mu)$ belonging to $]0,\min({\mathcal R}_5/32,1/2)]$ 
such that, for any $\rho \in \, ]0, \rho_1 (\mu)]$, 
any $x\in W$ and any
Riemannian barycentric simplex $\hat\tau\subset B(x,2\rho )$ whose 
thickness 
verifies $\forall q\in B(x,{\mathcal R}_4) \ \ t_{g_q}(\hat\tau_q)\geq t_0$ (see theorem {\rm \ref{2n}}), one has,  for 
all $v\in\hat\tau_x$ and $w\in T_x M$ 
\begin{equation*}
 (a) \  \Vert\,\phi_x(v)-v\,\Vert_x\leq \mu\  \  \  \hbox{\rm and}\  \  \  (b)\  
 \Vert\,d\phi_x(v)(w)-w\,\Vert_x\leq \mu\,\Vert\,w\,\Vert_x\  .
\end{equation*}
\end{prop}
\begin{rem} In other words, setting $E= E(\rho)$, for any $\sigma_E\in K_E$ verifying 
$T(\sigma_E) \subset B(x, 2\rho)$, the map
$\phi_x$ is a strong
$\mu$-approximation to $\hbox{Id}_x$ on 
$ \hat \tau_x = L_x(\sigma_E)\subset B(0_x, 2\rho) \subset (T_xM,g_x)$, see 
definition \ref{Argh55}. Or, if this makes sense,  {\it when lifted by} $\exp_x^{-1}$, the Riemannian barycentric coordinates in $M$ relative to $p_0,\cdots,p_n$ constitute a strong $\mu$-approximation to the affine coordinates in $T_xM$ relative to
$v_0,\cdots,v_n\,$.
\end{rem}

\begin{proof} 
First $(a)$ follows from $(b)$. Indeed, one has $\phi_x(v_i)= v_i$ for $i=0,\dots,n\,$, so  $(b)$ and
the mean value inequality along the linear arc from $v_i$ to $v\,$ give, using $\rho\!\leq\! \rho_1(\mu) \!\leq\! 1/2$ and ${\rm diam}_{g_x}(\hat\tau_x)\!\leq\!2\,{\rm diam}_{g}(\hat\tau)\!\leq\!2\,\rho$
\begin{equation}
\Vert\,\phi_x(v)-v\,\Vert_x\leq \Vert\,v-v_i\,\Vert_x\  (\sup_{v'\in \hat\tau_x} 
\Vert\, d\phi_x(v')-\hbox{Id}_x\ \Vert_{x,x})\leq \,2\rho\, \mu\leq \mu\  .
\end{equation}

In order to prove $(b)$, having in mind to apply lemma \ref{1b}, 
the idea is to make a choice of $\rho\,$ enabling us to get (for $i<j$ and $i,j\in\{0,1,\dots,n\}$) the following {\it ``likely''} inequalities (linked with $\mu$)
\begin{equation}\label{C}\Vert d\phi_x(v) (v_{i,j})-v_{i,j}\Vert_x\leq \frac{\mu\,t_0}{\sqrt{n}}\,\Vert v_{i,j}\Vert_x \ ,\ \ \ \hbox{where} \ \ \ v_{j,i}=v_j-v_i\ ,
\end{equation}
which urges to play with 
$Q_{x,q}=\exp_q^{-1}\circ\exp_x$ and estimate ($q=\exp_xv'$)
$$\Vert d\phi_x(v) (v_{i,j})-v_{i,j}\Vert_x=
\Vert (d\exp_x(v'))^{-1} \circ A_q^{-1} (Q_{x,q}(v_j)- Q_{x,q}(v_i))-v_{i,j}\Vert_x\ ,
$$
(for $A_q\,$, see definition \ref{formule1022} and lemma \ref{inv}).
Writing $w_i=Q_{x,q}(v_i)$ and $w_{i,j}=w_j-w_i\,$, this expression is less than $A+B$ with
\begin{gather*}A:=\Vert (d\exp_x(v'))^{-1} \circ A_q^{-1}\,w_{i,j}-
(d\exp_x(v'))^{-1} w_{i,j}\Vert_x \ ,\\ B:=\Vert (d\exp_x(v'))^{-1} w_{i,j}-v_{i,j}\Vert_x\ .
\end{gather*}

Note that $Q_{x,q}$ is for sure a $1/3$-quasi-isometry since $\exp_q$ and $\exp_x$ are $1/12$-quasi-isometries and $\exp_x^{-1}$ is a $1/11$-quasi-isometry.
The bound (\ref{inv333}) on $\Vert A_q^{-1}-\hbox{Id}_q\Vert_{q,q}$ in lemma  \ref{inv}, the fact that $(d\exp_x(v'))^{-1}$ and $Q_{x,q}\,$ are quasi-isometries give (here, one has $q\in\tau\subset B(x,2\rho)$)
\begin{multline}\label{A}A=\Vert (d\exp_x(v'))^{-1} \circ (A_q^{-1}-\hbox{Id}_q) (w_{i,j})
\Vert_x\leq\\\leq \frac{12}{11}\,\Vert (A_q^{-1}-\hbox{Id}_q) (w_{i,j})\Vert_q\leq  
\frac{12}{11}\,\Vert A_q^{-1}-\hbox{Id}_q\Vert_{q,q}\,\Vert w_{i,j}\Vert_q\leq \\\leq \frac{96}{33}\,{\mathcal M}_1\rho^2\,\Vert v_{i,j}\Vert_x=\frac{32}{11}\,{\mathcal M}_1\,\rho^2\,\Vert v_{i,j}\Vert_x\  .
\end{multline}

As $B$ vanishes if $q=x\,$, as $x$ and $q$ are close, applying twice the mean-value theorem gives the estimate below (since $d_g(x,q)\leq 2\rho$)
\begin{multline}\label{B}B=\!\Vert (d\exp_x(v'))^{-1} (\vartheta(x,q,v_j)- \vartheta(x,q,v_i))\Vert_x\!\leq\\\leq
\!
\frac{12}{11}\ \Vert \vartheta(x,q,v_j)- \vartheta(x,q,v_i)\Vert_q\leq\\
\leq \frac{12}{11}\ \sup_{(x,q,\omega)\in W\times B\times B_x}\Vert d_\omega\vartheta (x,q,\omega)\Vert_{x,q}\,\Vert v_j-v_i\Vert_x\leq \\\leq \frac{12}{11}\ \sup_{(x,q,\omega)\in W\times B\times B_x}\Vert d_q\,d_\omega\vartheta (x,q,\omega)\Vert_{x,q}\,\Vert \exp_x^{-1}q\Vert_x\,\Vert v_{i,j}\Vert_x\leq\\\leq
 \frac{24}{11}\ \Theta\,\rho\,\Vert v_{i,j}\Vert_x\  .
\end{multline}
Setting $\rho<\rho_1(\mu):=(t_0\,\mu)/\frac{8}{11}(\frac{3}{2}\,{\mathcal M}_1+2\ \Theta)\sqrt{n})\,$, one gets (\ref{C}) for $i<j\,,\, i,j\in\{0,1,\dots,n\}$ from (\ref{A}) and (\ref{B})  ($\rho$ is chosen $\leq 1/2$)
\begin{equation*}\Vert d\phi_x(v) (v_{i,j})-v_{i,j}\Vert_x\leq \frac{8}{11}(3\,{\mathcal M}_1\,\rho+2\,\Theta)\,\rho\,\Vert v_{i,j}\Vert_x\leq \frac{\mu\,t_0}{\sqrt{n}}\,\Vert v_{i,j}\Vert_x\ .
\end{equation*}
Use $t_{g_q}(\hat\tau_q)\geq t_0\,$, lemma \ref{1b} and (\ref{C}) to get $(b)$; for any $\,u=\sum_{j\not=0}\lambda_j \,v_{0,j}\,$
\begin{multline*}\Vert (d\phi_x(v)-\hbox{Id}_x)(u)\Vert_{x}\leq \sum_{j\not=0}\vert\lambda_j \vert\,\Vert d\phi_x(v)(v_{0,j})-v_{0,j}\Vert_{x}\leq\\
\leq \frac{\mu\,t_0}{\sqrt{n}}\,\sum_{j\not=0}\vert\lambda_j \vert\,\Vert v_{0,j}\Vert_x
\leq \mu\,(t_0\,\sup_{j\not=0}\Vert v_{0,j}\Vert_x)\,\sqrt{\sum_{j\not=0}\lambda_j^2}\ 
\leq\\\leq \mu\ r(\hat\tau_x)\,\sqrt{\sum_{j\not=0}\lambda_j^2}\leq \mu\ \Vert \sum_{j\not=0}\lambda_j \,v_{0,j}\Vert_x=\mu\ \Vert u\Vert_x\ .\qedhere
\end{multline*}
\end{proof}

\begin{defn} \label{ro2mu} Given $\mu>0\,$, define 
$\rho_2 (\mu):=\min(\rho(\mu),\rho_1(\mu))$, where $\rho(\mu)$ (where $\rho_1(\mu)$) refers to lemma \ref{2q} 
(proposition {\rm\ref{Argh33}}). 
\end{defn}

\begin{prop}\label{James-33}  
Given $\mu\in]0,\min(1/2,\beta_T/(1-C)]$, $\rho \in ]0, \rho_2 (\mu)]$ and $E=E(\rho)$, 
for any $x\in W$ and any $n$-simplex  $\sigma_E$ in $K_E$ verifying 
$T(\sigma_E) \subset B(x,2\rho)$,  the map
$\hat T_{E,x}=\exp_x^{-1}\circ \hat T_E$ is a strong ${\bf c}_1\mu$-approximation to the secant 
map $L_x$ 
on $\sigma_E$, where ${\bf c}_1:=\sup(1,\frac{2\beta_T}{1-C})$.
\end{prop}

\begin{proof} 
Let  $\rho \in \, ]0, \rho_1 (\mu)]\,$. 
By proposition \ref{Argh33}, for any simplex $\sigma_E$ such that  
$L_x(\sigma_E)\subset B(0_x, 2\rho)\subset (T_xM, g_x) $, the mapping $\phi_x$ is a 
strong $\mu$-approximation to $\hbox{Id}_x$ on $L_x(\sigma_E)$.

As 
$\Vert\,dT_x(y)\,\Vert_{y,x}\leq \beta_T/(1-C)$ (see (\ref{Ineg1})), 
one infers, using lemma \ref{2q} 
(in view of the choice $\mu \leq \beta_T/(1-C)$), for any $y\in \sigma_E$ and $u\in T_y \sigma_E$
\begin{multline}\label{James-55} 
\Vert\,dL_x(y)(u)\,\Vert_x\leq \Vert\,dT_x(y)(u)\,\Vert_x+\Vert\,dL_x(y)(u)-dT_x(y)(u)\,\Vert_x\leq\\ 
\leq  (\frac{\beta_T}{1-C} +\mu)\,
\Vert\,u\,\Vert_{{\rm std}} \leq  \frac{2\beta_T}{1-C} \, \Vert\,u\,\Vert_{{\rm std}}\  .
\end{multline}

The mapping $\hat T_{E,x}=\exp_x^{-1}\circ \hat T_E$ is equal to $\phi_x\circ L_x$ on 
each 
simplex $\sigma_E\in K_E$ (if $T(\sigma_E)$ is contained in $B(x,2\rho)$). 
Since
$\phi_x$ is a strong $\mu$-approximation to $\hbox{Id}_x$ on $L_x(\sigma_E)$, one 
can write (for any $v=L_x(y)$)
\begin{equation}  \label{James66}
\Vert\,\hat T_{E,x}(y)-L_x(y)\,\Vert_x=\Vert\,(\phi_x-\hbox{Id}_x)(v)\,\Vert_x
\leq \mu  \  ,
\end{equation}
and also (here (\ref{James-55}) also comes in)
\begin{multline}  \label{James55}
\Vert(d\hat T_{E,x}(y)-dL_x(y))(u)\Vert_x\!=\!
\Vert(d(\phi_x\!-\!\hbox{Id}_x)\circ (dL_x(y)))(u)\Vert_x\!\leq \\
\leq \Vert\,d(\phi_x-\hbox{Id}_x)\,\Vert_{x,x} \Vert\,dL_x(y)(u)\,\Vert_x
 \leq \frac{2\beta_T}{1-C}\,\mu\,
\Vert\,u\,\Vert_{{\rm std}}\  .
\end{multline}
Lines (\ref{James66}) and (\ref{James55}) express the claim given in proposition 
\ref{James-33}.
\end{proof}

\begin{lem}\label{James3} 
If $Q$ is a $\tilde C$-quasi-isometric 
$C^2$-map from $({\mathbb R}^n,\Vert\,\cdot\,\Vert_1)$ to 
$({\mathbb R}^n,\Vert\,\cdot\,\Vert_2)$, $\bar B\subset {\mathbb R}^n$ a compact convex set, 
define ${\mathfrak Q}=\sup_{x\in \bar B} \Vert d^2 Q(x)\Vert_{}$. 
\par
Consider a differentiable simplicial mapping 
$f:(K,\Vert\cdot\Vert_{{\rm std}})\rightarrow 
 \bar B\subset ({\mathbb R}^n,\Vert\cdot\Vert_1)$. 
 Set ${\mathfrak F}:=\sup_{y\in K}\Vert\,df(y)\,\Vert_{{\rm std},1}$ and 
 $\tilde{\bf c}:=1+\tilde C+ \wt{\mathfrak Q}\,{\mathfrak F}\,$. 
 
 Given $\mu>0$, if 
 $g :(K,\Vert\cdot\Vert_{{\rm std}})\rightarrow  \bar B$
 is a strong $\mu$-approximation to $f$ on $K$, then $Q\circ g$ 
is a strong $\tilde{\bf c}\,\mu$-approximation to $Q\circ f$ on $K$. 
\end{lem}

\begin{rem} \label{R.quasi}
We compute $\tilde{\bf c}\,$ in case $f=T_{x'}$ with $x'\in W\,$ and  
$Q:=Q_{x'',x'}=\exp_{x'}^{-1}\circ\exp_{x''}$ with $x''\in B(x',{\mathcal R}_4)$. Remembering $C\leq 1/12$ and using proposition 
\ref{proposition A}, 
we see that the mapping $Q\,$, well defined on
$B(x'',{\mathcal R}_4)$, is $\tilde C$-quasi-isometric, with $\tilde C=4\,C\,$. 
Line (\ref{Ineg1}) 
tells
$\Vert\,dT_x(y)\,\Vert_{y,x}\leq {\beta_T}/{1-C}\,$,
thus one has ${\mathfrak F}\leq\beta_T/(1-C)$ and gets 
$\tilde{\bf c}={\bf c}_2:=1+ 4 C+{\mathfrak Q}\,\beta_T/(1-C)$. 
\end{rem}

\begin{proof} 
By integration along the arc $h(t)=(1-t)\,f(y)+t\,g(y)$, $t\in [0,1]\,$, a mean value inequality gives $(a)$ in definition \ref{Argh55} for $Q\circ g$ and $Q\circ f\,$.
\begin{multline} \Vert\,Q\circ g(y)-Q\circ f(y)\,\Vert_2 =
 \\ \Vert\int_0^1
 \frac{dQ}{dt}(h(t))\,dt\,\Vert_2 =
\Vert\,\int_0^1
 dQ\,(h(t))\,\frac{dh}{dt}(t)\,dt\,\Vert_2\leq  \\ 
 \leq\int_0^1
 \Vert\,dQ(h(t))\,\Vert_{1,2}\  \Vert\, g(y)-f(y)\,\Vert_1\,dt\leq  (1+\tilde C)\,\mu 
\  .
\end{multline}

To prove property $(b)$ of definition \ref{Argh55} for $Q\circ g$ and $Q\circ f$, observe, using the 
mean 
value inequality applied to $d Q$
\begin{multline}  \Vert\,(d(Q\circ g)(y)-d(Q\circ f)(y))(u)\,\Vert_2
\leq \\
\leq
 \Vert\,dQ(g(y))\,\Vert_{1,2}\  \Vert\,dg(y)(u)-df(y)(u)\,\Vert_1+\\  
 + \Vert\,dQ(g(y))-dQ(f(y))\,\Vert_{1,2} \ \Vert\,df(y)(u)\,\Vert_1\leq \\ 
\leq(1+\tilde C+ \wt{\mathfrak Q}\,{\mathfrak F})\,\mu\,\Vert\,u\,\Vert_{{\rm std}}=\tilde{\bf c}\,\mu\,\Vert\,u\,\Vert_{{\rm std}}
\ . \qedhere
\end{multline}
\end{proof}

One also has the following
\begin{lem} \label{James-66} 
Given $\mu\in]0, \min (\frac{1}{2}, (\frac{\beta_T}{1-C}))]$, bring in the corresponding $\rho_2 (\mu)$ 
of definition {\rm\ref{ro2mu}}.
Given $\rho\in]0,\rho_2 (\mu)]\,$ and 
any 
$i$, setting $E= E(\rho)$, the 
mappings $h_i=\hat T_{E,x_i}=\exp_{x_i}^{-1}\circ \hat T_E\,$ are 
strong ${\bf c}_3\,\mu$-approximations to $f_i=T_{x_i}=\exp_{x_i}^{-1}\circ T$ 
on $K_{E,i}\,$,  where 
${\bf c}_3:= {\bf c}_2 (1+{\bf c}_1)$.
\end{lem}

\begin{proof} 
The integer $E$ being chosen according to the hypothesis, given any $i\,$, let $\sigma_E$ be any $n$-simplex 
$\sigma_E\in K_{E,i}\,$. By proposition 
\ref{James-33}, 
one knows 
$h_x=\hat T_{E,x}=\exp_x^{-1}\circ T$ is  a strong ${\bf c}_1 \mu$-approximation to $L_x$ on any simplex $\sigma_E$ 
such that  
$T(\sigma_E)$ is contained in $B(x, 2\rho)$. 

By lemma \ref{2q}, for any $i\,$ and  
$x\!\in\! B(x_i, {\mathcal R}_4/2)$, 
the map $T_x$ is a strong $\mu$-approximation to $L_x$
on any simplex $\sigma_E$ of $K_E$ 
with $T(\sigma_E)\! \subset\! B(x,2\rho)$.

Thus, for any $x\in B(x_i,{\mathcal R}_4/2)$ and any 
$i \,$, the map $\hat T_{E,x}$ is a 
strong $(1+{\bf c}_1)\,\mu$-approximation to $T_x$ on any 
$\sigma_E$ with 
$T(\sigma_E)\subset B(x,2\rho)$.

For any $x$ in $B(x_i,{\mathcal R}_4/2)$, the mapping 
$Q_i=\exp_{x_i}^{-1}\circ\exp_x$ is a $\tilde C$-quasi-isometric map from the 
closed ball 
$\bar B=\bar B(0_x,{\mathcal R}_4)\subset (T_x M,g_x)$ to $(T_{x_i}M,g_{x_i})$ 
(see the remark \ref{R.quasi}).

Thus, by  lemma \ref{James3} and remark \ref{R.quasi}, one knows that 
$h_i=\hat T_{E,x_i}=Q_i\circ \hat T_{E,x}$ 
is a strong
${\bf c}_3\,\mu$-approximation to 
$T_{x_i}= Q_i\circ T_x\,$ on any simplex $\sigma_E$ with 
$T(\sigma_E) \subset B(x,2\rho)$, for any $i$ and
$x\in B(x_i,{\mathcal R}_4/2)$, thus 
{\it $h_i$ 
is a 
$ {\bf c}_3\,\mu$-approximation to 
$T_{x_i}$ on $K_{E,i}$ and for any $i\,$}. 
\end{proof}

{\bf Completing the proof of theorem {\rm \ref{hatT}}}
\begin{proof}
Recall (from definition \ref{entier22}) $\rho_0 :=\min (\beta_T\,{\mathcal D}_n, {\mathcal R}_4/32)$. \par
According to proposition \ref{James2}, for 
$\delta>0$ small enough and 
$E\geq E(\rho_0)$,  
if, for any differentiable 
simplicial mapping $h : K_E \longrightarrow M\,$, the property below is true
(recall $h_i=\exp_{x_i}^{-1}\circ h$ and $f_i=T_{x_i}=\exp_{x_i}^{-1}\circ T$)
\par
{\it for any $i\,$ the map $h_i$ is a strong $\delta$-approximation to 
$f_i$ on $K_{E,i}\,$,}
\par\noindent then \hskip4mm
{\it $h$ is, with $T\,$, an embedding into $M\,$.}

Let $\delta>0$ be chosen small enough so that proposition \ref{James2} applies. {\it Set $h=\hat T_E\,$ and
 $\mu (\delta) := \min (1/2, \beta_T/(1-C) 
,\delta/{\bf c}_3)$.}
Lemma \ref{James-66} implies that, for any integer 
$E\geq E(\rho_2 (\mu(\delta)))$ and any $i\,$, the map $h_i=\exp_{x_i}^{-1}\circ 
h$ is a strong $\delta$-approximation to 
$f_i=T_{x_i}$ on $K_{E,i}\,$. 

Setting  
$E_1 (\delta):= \max (E(\rho_2 (\mu (\delta))), E(\rho_0))$, proposition \ref{James2} now implies that, for any $E\geq  E_1 (\delta)$, the mapping  
$\hat T_E :  K_E \rightarrow M$ is a differentiable simplicial 
 embedding. 

More, lemma \ref{L.Approx} says $d( \hat T_E (p), T(p)) \leq 2\delta\,$ for any $p\in \lvert K \lvert\,$.
\end{proof}

\subsection{Pseudo-homotheties}\label{pseudohomothety}
$ $
The geometric meaning of the last equalities in remark \ref{crx} is the existence of {\it pseudo-homotheties of ratio $t$ centered at the gravity center $p$} deforming a Riemannian barycentric simplex with vertices $\exp_pv_i=p_i$ (and $\sum_iv_i=0_p$) in a twisted geometric way into
the Riemannian barycentric simplex with vertices $\exp_p(tv_i)=p_i(t)$. Indeed,  one has 
$${\mathcal Q}_{{\bf \lambda},{\bf v}}(t)\!={\mathcal Q}_{h_t({\bf \lambda},{\bf v})}(1)=\!{\mathcal Q}_{\lambda(t),t{\bf v}}(1)\!=\!{\mathcal Q}_{\mu,{\bf w}}(1)\,.
$$
\vfill\eject

\part{Approximating Riemannian curvature}\label{M-T2}

\hskip4.3cm {\em To Marcel Berger}

\section{Holonomy}\label{holonomy}
\subsection{Path's decompositions and restricted holonomy}\label{holonomy1}

\par For all the basic material on holonomy in Riemannian manifolds and path's decompositions, see \cite{Be4}, \cite{Bes}, \cite{G-K-M}, \cite{K-N}, \cite{Li} and \cite{Pe}, for the basic material on Lie groups which is needed, see \cite{C-E}, \cite{I-L} and \cite{W}, for a discussion closer to our development see
\cite{Chr2}, \cite{N} and \cite{S}.
\begin{nota} Denote by $(M,g)$ a Riemannian differentiable manifold and by $(E,\nabla)$ a differentiable vector bundle over $M$ equipped with a linear connection whose associated covariant derivative is denoted $\nabla\,$, while $p_E:E\rightarrow M$ is the bundle projection. 

Let $\tilde g$ be a field of scalar products on $E$, independent of $\nabla\,$. 
\begin{rem} In the {\em Riemannian case} a Riemannian manifold $(M,g)$ is given, $E=TM$ and $\nabla$ is the Levi-Civita connection of $g=\tilde g\,$.
\end{rem}
\end{nota}
\begin{defn}\label{psquare1}$ $
\begin{itemize}\item
A {\em parametrised square} is the range 
$G([0,1]\times[0,1])$ of a differentiable mapping into a manifold $M$ 
$$G:[0,1]\times[0,1]\longrightarrow M\ .
$$
Such a square has a natural orientation (the one induced by $[0,1]\times[0,1]$ through $G$). The boundary curve $\Gamma\,$ is a naturally oriented loop. We often call $x:=G(0,0)$ its {\em base point} $G(0,0)$.

{\em Parametrised squares having source a more general square, for instance restrictions of the initial $G$
$$G:(s,t)\in[a,b]\times[c,d]\subset [0,1]\times[0,1]\longrightarrow G(s,t)\in M \ ,
$$
but also extensions, will tacitly arise in this text.}

\item Denote by $\gamma_\tau^{\varrho,\sigma}$ and $\delta_\varrho^{\tau,\upsilon}$ the paths
\begin{equation*} \label{paths1}\gamma_\tau^{\varrho,\sigma}\!:\!s\in[\varrho,\sigma]\!\mapsto\! G(s,\tau)\in M \  \ \hbox{and}\  \ \delta_\varrho^{\tau,\upsilon}\!:\!t\in[\tau,\upsilon]\!\mapsto\! G(\varrho,t)\in M
\  .
\end{equation*}
\item One also defines
\begin{equation*} \gamma_\tau^{\sigma,\varrho}=(\gamma_\tau^{\varrho,\sigma})^{-1}\  ,\  \  
\delta_\varrho^{\upsilon,\tau}=(\delta_\varrho^{\tau,\upsilon})^{-1}\  \   \hbox{and}\  \  \gamma_\tau^\sigma=\gamma_\tau^{0,\sigma}\  ,\  \  \delta_\varrho^\upsilon=\delta_\varrho^{0,\upsilon}\ .
\end{equation*}
\end{itemize}
\end{defn}
\begin{lem} \label{psquare2} A loop $\Gamma\subset M$ is homotopycally trivial if and only if it is homotopic to the boundary $\partial \,{\rm rge}(G)$ of a parametrised square $G\,$. 
\end{lem}
\begin{proof} If $\Gamma$ is homotopic to some $\partial \,{\rm rge}(G)\,$, one uses $G$ to define a homotopy retracting $\Gamma$ onto the base point $G(0,0)$. There are many ways to do this. For instance, if $\delta\vee\gamma$ denotes the composed path first following $\gamma$ then $\delta$ (so $\delta$ has for origin the end point of $\gamma$), define
\begin{equation} \label{paths0} \Gamma_t=\delta_0^{t,0}\vee\gamma_t^{1,0}\vee\delta_1^{0,t}\vee\gamma_0^{0,1}\  
 \hbox{where}\ \Gamma_1\ \hbox{is homotopic to}\  0\  \hbox{and to}\ \Gamma \,.
\end{equation}
The reciprocal property is almost the definition of a homotopy: there exists a continuous mapping
\begin{gather*} H:[0,1]\times[0,1]\longrightarrow M\  \  \hbox{such that} \cr 
H(0,t)=H(1,t)=H(s,0)=x
\  \  \  \  \hbox{and}\  \  \  \  \hbox{\rm rge}(H(\cdot,1))=\Gamma\  .
\end{gather*}
As any continuous $H:[0,1]\times[0,1]\rightarrow M$ can be approximated by a differentiable homotopy $G\,$, one gets the result.
\end{proof}

\begin{defn}\label{paraltrspt33} Consider a parametrised square $G$ defined in an open square $I\times I$ containing $[0,1]\times[0,1]$. Denote by $X$ a vector field along $G$, i. e. a mapping $X:I\times I\rightarrow E$ such that $p_E\circ X=G$. Parallel translation ${\mathcal P}$ along the lines $\gamma$ and $\delta$ that constitute $G$ (see definition {\rm\ref{psquare1}}) enables to define new vector fields, namely ${\mathcal P}_{\gamma;s} X$ and ${\mathcal P}_{\delta;t} X$ which satisfy (for any $\varrho,\tau,s,t$ making sense)
\begin{equation*} {\mathcal P}_{\gamma;s} X (\varrho,\tau)\!:=\! {\mathcal P}_{\gamma_\tau^{\varrho-s,\varrho}}(X(\varrho-s,\tau))\ 
\hbox{and}\ \, {\mathcal P}_{\delta;t} X (\varrho,\tau)\!:=\! {\mathcal P}_{\delta_\varrho^{\tau-t,\tau}}(X(\varrho,\tau-t))\,.
\end{equation*}
\end{defn}
\begin{lem} \label{fait12} The fields ${\mathcal P}_{\gamma;s} X$ and ${\mathcal P}_{\delta;t} X$ are (with $X$) differentiable in $(\varrho,\tau,s)$ and $(\varrho,\tau,t)$. 
\end{lem}
\begin{proof} This results from the differentiable dependence of the solutions upon the initial conditions of differentiable equations.
\end{proof}
\begin{nota}\label{nota121}
The loop $\Gamma_{\sigma,\upsilon}:=(\delta_0^\upsilon)^{-1}\vee(\gamma_\upsilon^\sigma)^{-1}\vee\delta_\sigma^\upsilon\vee\gamma_0^\sigma$ has inverse
\begin{equation} \label{nota12}
\Gamma_{\sigma,\upsilon}^{-1}:=(\gamma_0^\sigma)^{-1}\vee(\delta_\sigma^\upsilon)^{-1}\vee\gamma_\upsilon^\sigma\vee\delta_0^\upsilon\  
\end{equation}
and call ${\mathcal P}_{\Gamma_{\sigma,\upsilon}}^{-1}$ the parallel translation along this loop: ${\mathcal P}_{\Gamma_{\sigma,\upsilon}}^{-1}$ belongs to the connected component ${\bf GL}^+(E_x)\subset{\bf GL}(E_x)$ that contains $\hbox{\rm Id}$, and even to ${\bf O}_g^+(T_xM)\subset{\bf O}_g(T_xM)$ in the Riemannian case (see \cite{C-E}, \cite{W}).
\par In a group ${\bf G}$, the connected component of $\hbox{\rm Id}$ is denoted by ${\bf G}^+$. 
\end{nota}
\begin{lem} (see {\rm\cite{G-K-M}, A $(iii)$} pp. {\rm 53-4})\label{squcurvparrtrslt} If $G$ is a parametrised square of $M$, the mappings $(\sigma,\upsilon)\mapsto \Gamma_{\sigma,\upsilon}$ or $(\sigma,\upsilon)\mapsto \Gamma_{\sigma,\upsilon}^{-1}$ are differentiable and (with $x=G(0,0)$) one has in ${\bf gl}(E_x)$ (in ${\bf o}_g(T_xM)$ if $E\!=\!TM, \tilde g\!=\!g$)
\begin{equation} \lim_{\sigma\rightarrow0}\lim_{\upsilon\rightarrow0}\frac{{\mathcal P}_{\Gamma_{\sigma,\upsilon}}^{-1}\!-\!\hbox{\rm Id}}{\sigma\upsilon}\!=\!R_x(dG(\frac{\partial}{\partial \sigma}),dG(\frac{\partial}{\partial \upsilon}))\!=\!\lim_{\upsilon\rightarrow0}\lim_{\sigma\rightarrow0}\frac{{\mathcal P}_{\Gamma_{\sigma,\upsilon}}^{-1}\!-\!\hbox{\rm Id}}{\sigma\upsilon}\,.
\end{equation}
\end{lem}
The proof is given in appendix \ref{squcurvparrtrslt1}.
\begin{rem}\label{singsquare} If the images of 
$\frac{\partial}{\partial \sigma}$ and $\frac{\partial}{\partial \upsilon}$ are colinear at some point of $I\times I$, the map $G$ becomes singular and $R_x(dG(\frac{\partial}{\partial \sigma}),dG(\frac{\partial}{\partial \upsilon}))$ vanishes.
\end{rem}
 
The way $(\sigma,\upsilon)\in I\times I\mapsto {\mathcal P}_{\Gamma_{\sigma,\upsilon}}^{-1}\in{\bf GL}^+(E_x)\subset{\bf End}(E_x)={\bf gl}(E_x)$ is defined through the definition of the loop (\ref{nota12}) implies
\begin{equation}\label{blue3}\frac{\partial {\mathcal P}_{\Gamma_{\sigma,\upsilon}}^{-1}}{\partial \sigma}(0,0)\!=\!
\frac{\partial {\mathcal P}_{\Gamma_{\sigma,\upsilon}}^{-1}}{\partial \upsilon}(0,0)\!=\!\frac{\partial^2 {\mathcal P}_{\Gamma_{\sigma,\upsilon}}^{-1}}{\partial \sigma^2}(0,0)\!=\!\frac{\partial^2 {\mathcal P}_{\Gamma_{\sigma,\upsilon}}^{-1}}{\partial \upsilon^2}(0,0)\!=\!0
\end{equation}
and we know from lemma \ref{squcurvparrtrslt}
\begin{equation}\label{blue4}\frac{\partial^2 {\mathcal P}_{\Gamma_{\sigma,\upsilon}}^{-1}}{\partial \sigma\partial \upsilon}(0,0)=\frac{\partial^2 {\mathcal P}_{\Gamma_{\sigma,\upsilon}}^{-1}}{\partial \upsilon\partial \sigma}(0,0)=R_x(dG(\frac{\partial}{\partial \sigma}),dG(\frac{\partial}{\partial \upsilon}))\  .
\end{equation}
We thus get in ${\bf gl}(E_x)$ from (\ref{blue3}) and (\ref{blue4}) 
\begin{equation} \label{formula12} {\mathcal P}_{\Gamma_{\sigma,\upsilon}}^{-1}=\hbox{\rm Id}+
\sigma \upsilon\  R_x(dG(\frac{\partial}{\partial \sigma}),dG(\frac{\partial}{\partial \upsilon}))+o_x(\sigma^2+ \upsilon^2)\  .
\end{equation}
 \par From now, compare with \cite{S}, \cite{N}, \cite{Be4}, \cite{Bes}, \cite{K-N}, \cite{Li}.
\begin{nota} Define the loop $\Gamma_{\varrho,\tau;\sigma,\upsilon}^{-1}$ to be
\begin{equation} \label{nota1233}
\Gamma_{\varrho,\tau;\sigma,\upsilon}^{-1}=(\gamma_\tau^{\varrho,\sigma})^{-1}\vee(\delta_\sigma^{\tau,\upsilon})^{-1}\vee\gamma_\upsilon^{\varrho,\sigma}\vee\delta_\varrho^{\tau,\upsilon}\  ,
\end{equation}
so that $\Gamma_{0,0;\sigma,\upsilon}^{-1}=\Gamma_{\sigma,\upsilon}^{-1}\,$.
\end{nota}
\begin{lem}\label{movebase} Features leading to analogous {\rm(\ref{blue3}), (\ref{blue4}), (\ref{formula12})} hold for the ``square'' (definition {\rm\ref{psquare1}}) $G_{\mid [\varrho,\sigma]\times[\tau,\upsilon]}$, where $(\varrho,\tau)$ and $(\sigma\!-\!\varrho,\upsilon\!-\!\tau)$ now play the roles of $(0,0)$ and $(\sigma,\upsilon)$. There exists 
$\epsilon_{G(\varrho,\tau)}(\sigma\!-\!\varrho,\upsilon\!-\!\tau)$
{\em continuous} in $\varrho,\tau,\sigma,\upsilon$, tending to $0$ with $\Vert(\sigma\!-\!\varrho,\upsilon\!-\!\tau)\Vert$, such that 
\begin{equation*}\Vert o_{G(\varrho,\tau)}((\sigma\!-\!\varrho)^2\!+\!(\upsilon\!-\!\tau)^2)\Vert\leq \epsilon_{G(\varrho,\tau)}(\sigma\!-\!\varrho,\upsilon\!-\!\tau)\,((\sigma\!-\!\varrho)^2+(\upsilon\!-\!\tau)^2)\ ,
\end{equation*}
where $\Vert\cdot\Vert$ is the operator norm $\Vert\cdot\Vert_{\tilde g,\tilde g}$ defined in ${\bf gl}(E_{G(\varrho,\tau)})$ by $\tilde g\,$.
\end{lem}
A proof is given in appendix \ref{Taylor2}.
\begin{nota} \label{notation007}Write $R_G(\sigma,\upsilon)$ in place of $R_{G(\sigma,\upsilon)}(dG(\frac{\partial}{\partial \sigma}),dG(\frac{\partial}{\partial \upsilon}))\,$.
\end{nota}
\begin{prop} \label{formula14} The following equality holds in ${\bf gl}(E_x)$ (respectively in ${\bf o}(T_xM)$)
\begin{equation*}  
\frac{\partial  {\mathcal P}_{\Gamma_{\sigma,\upsilon}}^{-1}}{\partial \upsilon}(\sigma,0)=\int_0^\sigma {\mathcal P}_{\gamma_0^s}^{-1}\circ R_G(s,0) \circ {\mathcal P}_{\gamma_0^s}\ \ d\,s\  .
\end{equation*}
\end{prop}
\begin{proof}
The next expression ($N$ is an integer) will be a useful key
\begin{equation}\label{blue5}
\frac{\partial {\mathcal P}_{\Gamma_{\sigma,\upsilon}}^{-1}}{\partial \upsilon}(\sigma,0)=
\lim_{\upsilon\rightarrow0}\frac{{\mathcal P}_{\Gamma_{\sigma,\upsilon}}^{-1}-\hbox{\rm Id}}{\upsilon}=
\lim_{N\rightarrow\infty}\frac{{\mathcal P}_{\Gamma_{\sigma,\frac{\sigma}{N}}}^{-1}-\hbox{\rm Id}}{\frac{\sigma}{N}}\  .
\end{equation}
So, introduce the decomposition of the rectangle loop $\Gamma_{\sigma,\frac{\sigma}{N}}^{-1}$ into ``square'' loops $\Delta_{i,N}$, all having base point $G(0,0)$
\begin{gather}\label{blue6}
\Gamma_{\sigma,\frac{\sigma}{N}}^{-1}=\Delta_{N-1,N}\vee\dots\vee \Delta_{1,N}\vee \Delta_{0,N}\  \  \  \hbox{where, for}\\
\label{blue7}
i=0,1,\dots,N-1 \hskip10mm \Delta_{i,N}=(\gamma_0^{\frac{i\sigma}{N}})^{-1} \vee\Gamma_{\frac{i\sigma}{N};0,\frac{(i+1)\sigma}{N},\frac{\sigma}{N}}^{-1}\vee\gamma_0^{\frac{i\sigma}{N}}\  .
\end{gather}
Apply (\ref{formula12}) (lemma \ref{movebase}) to $\Gamma_{i,N}^{-1}=\Gamma_{\frac{i\sigma}{N};0,\frac{(i+1)\sigma}{N},\frac{\sigma}{N}}^{-1}$ to get in 
${\bf gl}(E_{G(\frac{i\sigma}{N},0)})$
\begin{equation} \label{blue8}
{\mathcal P}_{\Gamma_{i,N}}^{-1}=\hbox{\rm Id}+
(\frac{\sigma}{N})^2\  R_G(\frac{i\sigma}{N},0)+o_i((\frac{\sigma}{N})^2)\  .
\end{equation}
One has (apply lemma \ref{movebase}, here $\Vert\cdot\Vert$ is the operator norm $\Vert\cdot\Vert_{\tilde g,\tilde g}$ defined in ${\bf gl}(E_{(\frac{i\sigma}{N},0)})$ by $\tilde g$)
\begin{equation}\label{blue9}\Vert\,o_i((\frac{\sigma}{N})^2)\,\Vert\leq (\frac{\sigma}{N})^2\,\epsilon_i(\frac{\sigma}{N})\  ,
\end{equation}
with $\epsilon_i(\frac{\sigma}{N})\rightarrow0$ as $N\rightarrow\infty$. 
Formulas (\ref{blue7}) and (\ref{blue8}) give
\begin{equation}\label{formule1}
{\mathcal P}_{\Delta_{i,N}}=\hbox{\rm Id}+
(\frac{\sigma}{N})^2\  {\mathcal P}_{\gamma_0^{\frac{i\sigma}{N}}}^{-1}\circ R_G(\frac{i\sigma}{N},0)\circ{\mathcal P}_{\gamma_0^{\frac{i\sigma}{N}}}+o'_i((\frac{\sigma}{N})^2)\  ,
\end{equation}
where $o'_i((\frac{\sigma}{N})^2)$ takes the role of $o_i((\frac{\sigma}{N})^2)$, this time in ${\bf gl}(E_x)$.

\par By the uniform continuity of the $\epsilon_{G(\varrho,\tau)}(\sigma-\varrho,\upsilon-\tau)$ (lemma \ref{movebase}), for any $\epsilon>0$  exists $\eta$ such that $(\sigma-\varrho)^2+(\upsilon-\tau)^2<\eta$ implies $\vert\epsilon_{G(\sigma-\varrho,\upsilon-\tau)}\vert<\epsilon\,$. For any $\epsilon>0$, the  compacity of the path $\gamma_0^\sigma$ thus allows to choose a large enough $N$ such that, for $i=0,1,\dots,N-1$
\begin{equation}\label{formule2}\Vert\,o'_i((\frac{\sigma}{N})^2)\,\Vert\leq (\frac{\sigma}{N})^2\,\epsilon\  .
\end{equation}
Setting $\ K=\sup_{s\in[0,\sigma]}\Vert\,{\mathcal P}_{\gamma_0^s}^{-1}\circ R_G(s,0)\circ{\mathcal P}_{\gamma_0^s}\,\Vert$ and
$P_{i,N}={\mathcal P}_{\Delta_{i,N}}-\hbox{\rm Id}$, one gets
from (\ref{formule1}) and (\ref{formule2}) 
\begin{equation} \label{formule3} \Vert\,P_{i,N}\,\Vert=\Vert\,{\mathcal P}_{\Delta_{i,N}}-\hbox{\rm Id}\,\Vert\leq (\frac{\sigma}{N})^2(K+\epsilon)\  .
\end{equation}
One also has from (\ref{blue6})
\begin{equation} \label{formule4} {\mathcal P}_{\Gamma_{\sigma,\frac{\sigma}{N}}}^{-1}={\mathcal P}_{\Delta_{N-1,N}}\circ\dots\circ {\mathcal P}_{\Delta_{0,N}}=
(\hbox{\rm Id}+P_{N-1,N})\circ\dots\circ(\hbox{\rm Id}+P_{0,N})\  .
\end{equation}
Expand (\ref{formule4})
\begin{multline} \label{formule5} {\mathcal P}_{\Gamma_{\sigma,\frac{\sigma}{N}}}^{-1}-\hbox{\rm Id}=
\sum_{i=0}^{N-1}P_{i,N}+\sum_{0\leq i<j\leq N-1}P_{j,N}\circ P_{i,N}+\dots\\
\dots+P_{N-1,N}\circ \dots\circ P_{1,N}\circ P_{0,N}\  ,
\end{multline}
and set $a_N:=({\sigma}/{N})^2(K+\epsilon)$, so that (\ref{formule3}) and (\ref{formule5}) give
\begin{equation} \label{formule6} \Vert\,{\mathcal P}_{\Gamma_{\sigma,\frac{\sigma}{N}}}^{-1}-\hbox{\rm Id}-\sum_{i=0}^{N-1}P_{i,N}\,\Vert\leq
\binom{N}{2}\,a_N^2+\binom{N}{3}\,a_N^3+\dots+a_N^N\  .
\end{equation}
As, for any $l\leq N$
\begin{equation} \label{blue11}\binom{N}{l}\,a_N^l= \frac{N(N-1)\dots(N-l+1)}{l!}\,a_N^l\leq \frac{(Na_N)^l}{l!}\  ,
\end{equation}
one gets from (\ref{formule6}) and (\ref{blue11}) 
\begin{equation} \label{formule7} 
\Vert\,{\mathcal P}_{\Gamma_{\sigma,\frac{\sigma}{N}}}^{-1}-\hbox{\rm Id}-\sum_{i=0}^{N-1}P_{i,N}\,\Vert\leq
e^{N\,a_N}-1-N\,a_N=N^2\,a_N^2 \ O(1)\  ,
\end{equation}
where $O(1)$ stays bounded as $N$ tends to $\infty$ ($N\ a_N$ tends to $0$). Thanks to the {\em key} (\ref{blue5}), to (\ref{formule7}) and to $a_N:=({\sigma}/{N})^2(K+\epsilon)$, we get
\begin{multline}\label{blue12}\Vert\,\frac{\partial {\mathcal P}_{\Gamma_{\sigma,\upsilon}}^{-1}}{\partial \upsilon}(\sigma,0)-\lim_{N\rightarrow\infty}\frac{\sum_{i=0}^{N-1}P_{i,N}}{\frac{\sigma}{N}}\,\Vert=\\=
\lim_{N\rightarrow\infty}\Vert\,\frac{{\mathcal P}_{\Gamma_{\sigma,\frac{\sigma}{N}}}^{-1}-\hbox{\rm Id}-\sum_{i=0}^{N-1}P_{i,N}}{\frac{\sigma}{N}}\,\Vert=\lim_{N\rightarrow\infty} \ O(\frac{N^3\,a_N^2}{\sigma})=0\  .
\end{multline}
But, from (\ref{formule1}) we also have
\begin{multline}\label{blue13}
\lim_{N\rightarrow\infty}\frac{\sum_{i=0}^{N-1}P_{i,N}}{\frac{\sigma}{N}}=\\=
\lim_{N\rightarrow\infty}  \sum_{i=0}^{N-1} \frac{\sigma}{N}\  
({\mathcal P}_{\gamma_0^{\frac{i\sigma}{N}}}^{-1}\circ R_G(\frac{i\sigma}{N},0)\circ{\mathcal P}_{\gamma_0^{\frac{i\sigma}{N}}})+
\lim_{N\rightarrow\infty} \sum_{i=0}^{N-1} o_i'(\frac{\sigma}{N})=\\=
\int_0^\sigma {\mathcal P}_{\gamma_0^s}^{-1}\circ R_G(s,0) \circ {\mathcal P}_{\gamma_0^s}\ \ d\,s+
\lim_{N\rightarrow\infty}  \sum_{i=0}^{N-1} o'_i(\frac{\sigma}{N})\  .
\end{multline}
Bringing in (\ref{formule2}), the last term in (\ref{blue13}) is seen to tend to $0$ as $N\rightarrow\infty\,$. Putting together (\ref{blue12}) and (\ref{blue13}), the proof is complete.
\end{proof}
\par The next theorems state equalities in ${\bf gl}(E_x)$ (actually in ${\bf o}(T_xM)$ in the Riemannian case).
\begin{nota} \label{AlalaBlala} Define $A_{\varrho,\tau}$ and $B_{\varrho,\tau}$ to be the paths
\begin{equation} \label{blue14}A_{\varrho,\tau}=\gamma_\tau^\varrho\vee\delta_0^\tau \hskip5mm  \hbox{and} \hskip5mm
B_{\varrho,\tau}=(\gamma_\tau^{\varrho,1})^{-1}\vee\delta_1^\tau\vee\gamma_0^1=A_{\varrho,\tau}\vee\Gamma_{1,\tau}\  .
\end{equation}
\end{nota}
\begin{thm}  \label{theorem A} Given any smooth parametrised square $G\,$, one has
\begin{equation} \label{formule8}
{\mathcal P}_{\Gamma_{1,1}}^{-1}-\hbox{\rm Id}=\int_0^1\!\!\int_0^1 {\mathcal P}_{B_{\varrho,\tau}}^{-1}\circ R_G(\varrho,\tau) \circ {\mathcal P}_{A_{\varrho,\tau}}\ \ d\,\varrho\,d\,\tau\  .
\end{equation}
Thus, for a parametrised square defining a homotopy through paths $\gamma_t^1$ running from a fixed origin $x$ to a fixed extremity $y$
\begin{gather*} G:(s,t)\in[0,1]\times[0,1]\longmapsto  G(s,t)\in M\   ,\\
t\mapsto G(0,t)=\delta_0^1(t)=x\  ,\  \  \  \  t\mapsto G(1,t)=\delta_1^1(t)=y\  \\ 
s\mapsto G(s,t)=\gamma_t^1(s)\  ,
\end{gather*}
the following formula holds true
\begin{equation} \label{formule9} {\mathcal P}_{\gamma_1^1}-{\mathcal P}_{\gamma_0^1}=
\int_0^1 \!\!\int_0^1 {\mathcal P}_{\gamma_\tau^{\varrho,1}}\circ R_G(\varrho,\tau) \circ {\mathcal P}_{\gamma_\tau^\varrho}\ \ d\,\varrho\,d\,\tau\  .
\end{equation}
\end{thm}
\begin{proof} By proposition \ref{formula14} (change the variable $\upsilon$ to $\upsilon-\tau$), one has
\begin{equation*}
\frac{\partial  {\mathcal P}_{\Gamma_{0,\tau ; 1,\upsilon}}^{-1}}{\partial \upsilon}(1,\tau)=\int_0^1 {\mathcal P}_{\gamma_\tau^s}^{-1}\circ R_G(s,\tau) \circ {\mathcal P}_{\gamma_\tau^s}\ \ d\,s\  ,
\end{equation*}
so that, writing the decomposition
\begin{equation*}\Gamma_{1,\upsilon}^{-1}=(\gamma_0^1)^{-1}\vee (\delta_1^\tau)^{-1}\vee\gamma_\tau^1 \vee  \Gamma_{0,\tau ; 1,\upsilon}^{-1} \vee \delta_0^\tau
\end{equation*}
we get
\begin{equation} \label{formulatralala}
\frac{\partial  {\mathcal P}_{\Gamma_{1,\upsilon}}^{-1}}{\partial \upsilon}(1,\tau)=\int_0^1 {\mathcal P}_{B_{\varrho,\tau}}^{-1}\circ R_G(\varrho,\tau) \circ {\mathcal P}_{A_{\varrho,\tau}}\ \ d\,\varrho\  .
\end{equation}
Integrating (\ref{formulatralala}) in $\tau$ from $0$ to $1$, one expresses 
$
\int_0^1({\partial  {\mathcal P}_{\Gamma_{1,\upsilon}}^{-1}}/{\partial \upsilon})(1,\tau)\,d\,\tau$ as to give (\ref{formule8}).
Since $\delta_0^1(t)=x$ and $\delta_1^1(t)=y$, one gets with (\ref{blue14})
\begin{equation}  {\mathcal P}_{A_{\varrho,\tau}}= {\mathcal P}_{\gamma_\tau^\varrho} \hskip5mm  \hbox{and} \hskip5mm
 {\mathcal P}_{\gamma_0^1}\circ{\mathcal P}_{B_{\varrho,\tau}}^{-1}= {\mathcal P}_{\gamma_\tau^{\varrho,1}} \  , 
\end{equation}
so that (\ref{formule8}) can be written
as formula (\ref{formule9}) of the claim.
\end{proof}
\begin{defn}  \label{normcurv} 
If $G$ is a parametrised {\em immersed} square, put on $\hbox{\rm rge}(G)$ the orientation induced by $G$ and denote by $d\,{\rm area}$ the element of surface induced by the Riemannian metric $g\,$.
\par Define, for a parametrised {\em immersed} square $G$ in $(M,g)$, the {\em normed curvature} $R^G$ (as an operator acting on $E_{G(\varrho,\tau)}$) 
\begin{equation} \label{Courbnorm} R^G(\varrho,\tau)=\frac {R_G(\varrho,\tau)}{\Vert dG(\frac{\partial}{\partial \varrho})\wedge dG(\frac{\partial}{\partial \tau})\Vert_g}=\frac {R( dG(\frac{\partial}{\partial \varrho}), dG(\frac{\partial}{\partial \tau})}{\Vert dG(\frac{\partial}{\partial \varrho})\wedge dG(\frac{\partial}{\partial \tau})\Vert_g}\  .
\end{equation} 

\end{defn}
\begin{cor}\label{corTheorem A} Thanks to definition {\rm\ref{normcurv}}, with the same assumptions, one can rephrase {\rm(\ref{formule9})} for a parametrised {\em immersed} square $G$ in $(M,g)$
\begin{equation} \label{formule10} {\mathcal P}_{\gamma_1^1}-{\mathcal P}_{\gamma_0^1}=
\int\!\!\int_{[0,1]\times[0,1]} {\mathcal P}_{\gamma_\tau^{\varrho,1}} \circ R^G(\varrho,\tau) \circ {\mathcal P}_{\gamma_\tau^\varrho}\ \ G^\ast d\,{\rm area}\,.
\end{equation}
\end{cor}
\begin{proof} Rewrite (\ref{formule9}) using (\ref{Courbnorm}).
\end{proof}
\begin{nota} If $\delta_0\!=\!\delta_1\!=\!\gamma_0\!\equiv\! y\!=\!x\,$, call $\Gamma^{\varrho,\tau}$ the loop (based at $G(\varrho,\tau)$)
\begin{equation} \label{extraloop1}  \Gamma^{\varrho,\tau}:= B_{\varrho,\tau}\vee A_{\varrho,\tau}^{-1}=
(\gamma_\tau^{\varrho,1})^{-1}\vee(\gamma_\tau^\varrho)^{-1}\  .
\end{equation}
\end{nota}
\begin{rem} In the case above, thanks to ({\rm\ref{extraloop1}}) formula 
({\rm\ref{formule10}}) reads 
\begin{equation} \label{formule10biss} {\mathcal P}_{\gamma_1^1}-{\mathcal P}_{\gamma_0^1}=
\int\!\!\int_{[0,1]\times[0,1]} {\mathcal P}_{\gamma_\tau^\varrho}^{-1} \circ {\mathcal P}_{\Gamma^{\varrho,\tau}}^{-1}\circ R^G(\varrho,\tau) \circ {\mathcal P}_{\gamma_\tau^\varrho}\ \ G^\ast d\,{\rm area}\,.
\end{equation}
\end{rem}
We now give a second development of this way of thinking and another theorem, which embodies a classical elementary piece of the two-dimensional Gauss-Bonnet theorem.
\begin{defn} \label{defn14} Two parametrised squares $G_1$ and $G_2$ sharing the same parametrised boundary loop $\partial G_1=\partial G_2=\Gamma$ and the same base point $G_1(0,0)=G_2(0,0)=x$ are {\em homotopic fixing the boundary} if there exists a map $GG$ defined on the cube with values into $M$ which reduces to $G_1$ on a face, to $G_2$ on the opposite, which on the other faces describes a homotopy leaving $\Gamma$ globally invariant and having fixed base point $x$. We denote by $[G]$ the homotopy class defined in this sense for a parametrised square $G$.
\end{defn}
\begin{nota} \label{MC1033} For any parametrised square $G$, call ${\mathfrak G}$ the curve in ${\bf G}={\bf GL}^+(E_x)$ (in ${\bf G}={\bf O}_{g_x}^+(T_xM)$ in the Riemannian case) ${\mathfrak G}: \tau\in[0,1]\mapsto{\mathcal P}_{\Gamma_{1,\tau}}^{-1}\in {\bf G}$, running from $\hbox{\rm Id}$ to ${\mathcal P}_{\Gamma_{1,1}}^{-1}={\mathcal P}_\Gamma^{-1}\,$ (see ({\rm\ref{nota12}})).
\end{nota}
\begin{defn} \label{MC} The Maurer-Cartan vector valued differential $1$-form $\omega$ is defined on ${\bf G}$ subgroup of ${\bf GL}(E_x)$ (see \cite{I-L} section {\rm1.6}, page {\rm16}ff or \cite{K-N} page {\rm41}) by requiring, on each left invariant field ${\mathpzc b}\,$, i. e. ${\mathpzc b}_A:=A\,b\in T_A{\bf G}$ with $A\in {\bf G}$ and $b$ in the Lie algebra $\underline {\bf G}\equiv T_{\hbox{\rm Id}}{\bf G}$
\begin{equation}  \  \label{MC0} \omega_A({\mathpzc b})=\omega_A(A\,b)=b\in \underline {\bf G}\subset {\bf gl}(E_x)\equiv T_{\hbox{\rm Id}}{\bf GL}(E_x)\, .
\end{equation}
\end{defn}
\begin{lem} \label{MC1} Consider the basis $E_i^j=e_j\otimes e_i^\ast$ of ${\bf End}(E_x)={\bf gl}(E_x)$ defined by a chosen basis $e_1,\dots,e_l$ in $E_x$ (the dual is $e_1^\ast,\dots,e_l^\ast$). Setting $A=\sum_{i,j} A_i^j E_i^j$ and thus
$dA=\sum_{i,j} dA_i^j E_i^j$, the form $\omega$ reads in the embedding $\ \mu: A\in {\bf G} \mapsto (A_i^j)\in {\bf Mat}_{l,l}({\mathbb R})$ (as usual, the $d A_i^j$ denote the coordinate forms $(E_i^j)^\ast$ in ${\bf End}(E_x)$ evaluated on $dA$)
\begin{equation} \label{MC2} \omega_A\equiv(\mu(A))^{-1}\,d\mu_A= (\sum_{k} (A^{-1})_k^j\,(d A)_i^k)\  ,
\end{equation}
and, in this context, it is usually denoted by $A^{-1}\,dA$.
\end{lem}
\begin{proof} Indeed, (\ref{MC0}) and (\ref{MC2}) both define left invariant differential $1$-forms on ${\bf G}$ coinciding on $\underline {\bf G}$, since $\omega_{\hbox{\rm Id}}(b)\,=\,b$ for any $ b\,\in\, \underline {\bf G}\,$.
\end{proof}
\begin{lem} \label{MC3-2}The exterior product and exterior derivative of $\underline {\bf G}$-valued differential $1$-forms like $\omega=A^{-1}\,dA$ make sense and the Maurer-Cartan equation holds
\begin{equation} \label{MC3} 
d\omega= d\,(A^{-1}\,dA)=-(A^{-1}\,dA)\wedge (A^{-1}\,dA)=-\,\omega\wedge\omega\  .
\end{equation}
\end{lem}
For a proof, see appendix \ref{MC30}.
\begin{thm} \label{theorem B} Recall from {\rm(\ref{blue14})} that ${\mathcal P}_{A_{\varrho,\tau}}={\mathcal P}_{\gamma_\tau^\varrho}\circ{\mathcal P}_{\delta_0^\tau}$. The following formula (in ${\bf gl}(E_x)$ or in ${\bf o}_g(T_xM)$, in the Riemannian case) holds true for any parametrised square $G$
\begin{equation} \label{formula15} \int_0^1\!\!{\mathcal P}_{\Gamma_{1,\tau}}\circ 
\frac{\partial  {\mathcal P}_{\Gamma_{1,\upsilon}}^{-1}}{\partial \upsilon}(1,\tau)\ d\,\tau=
\int_0^1\!\!\int_0^1\!\!{\mathcal P}_{A_{\varrho,\tau}}^{-1}\circ R_G(\varrho,\tau)\circ{\mathcal P}_{A_{\varrho,\tau}} \ \ d\,\varrho\,d\,\tau\, .
\end{equation}
The left-hand side (called the {\rm generalised angle} associated to $G$) equals the integral of the Maurer-Cartan form $A^{-1}\,dA$ over the path ${\mathfrak G}$, i. e.
\begin{equation} \label{genangle} \int_0^1{\mathcal P}_{\Gamma_{1,\tau}}\circ 
\frac{\partial  {\mathcal P}_{\Gamma_{1,\upsilon}}^{-1}}{\partial \upsilon}(1,\tau)\ \ d\,\tau=\int_{\mathfrak G}A^{-1}\,dA\  .
\end{equation}
If $GG$ is a homotopy fixing the boundary between two parametrised squares $G_1$ and $G_2$, the difference of generalised angles is given by
\begin{equation} \label{formula16}
\int_{{\mathfrak G}_2}A^{-1}\,dA-\int_{{\mathfrak G}_1}A^{-1}\,dA=-\int\!\!\int_{{\mathfrak G}{\mathfrak G}}(A^{-1}\,dA)\wedge(A^{-1}\,dA)\  ,
\end{equation}
where ${\mathfrak G}{\mathfrak G}$ is the parametrised surface equal to the range of $GG$ in ${\bf gl}(E_x)$ (in ${\bf o}_g(T_xM)$, in the Riemannian case)  and having boundary ${\mathfrak G}_1\cup{\mathfrak G}_2$ (notice that those two curves have the same end-points).
\par If $M$ is $2$-dimensional, $E=TM$ and $\nabla$ is the covariant derivative defined by a Riemannian metric on $M$, if $G$ is an embedding, the above theorem reduces to part of the classical Gauss-Bonnet theorem for 
\hbox{\rm rge}(G). In this case, the generalised angle reduces to the usual angle associated to the parallel translation along $\Gamma^{-1}$.
\end{thm}
\begin{rem} \label{GaussBonnet:2} Actually, the last claim is only half of the Gauss-Bonnet theorem for an elementary piece of surface: the other part consists in 
\vskip1mm
1. the turning tangent theorem: the angle between the unit tangent vector to $\partial G$ and a non zero fixed field over ${\rm rge}(G)$ (say $\partial/\partial s$) after accumulation during a full loop, including the jumps at corners, is $2\pi$;
\vskip1mm
 2. the fact that the geodesic curvature of an arc-length's parametrised curve $c(\alpha)$ is the $\alpha$-derivative of the angle at $c(\alpha)$ between a parallel field along $c$ and ${d c}/{d \alpha}$, see \cite{C} or \cite{dC1}, and this can be applied to the curve $\partial G$ (see \cite{C} pages 175-9 or \cite{dC1} page 264ff).
\end{rem}
\begin{proof}
Start with the following decomposition (see (\ref{nota12}) and (\ref{nota1233}))
\begin{equation}\Gamma_{1,\tau}\vee\Gamma_{1,\upsilon}^{-1}=(\delta_0^\tau)^{-1}\vee\Gamma_{0,\tau,1,\upsilon}^{-1}\vee\delta_0^\tau\ , 
\end{equation}
and derive in $\upsilon$ the corresponding equality between parallel translations. Applying proposition \ref{formula14}, one has
\begin{equation}
{\mathcal P}_{\Gamma_{1,\tau}}\circ 
\frac{\partial  {\mathcal P}_{\Gamma_{1,\upsilon}}^{-1}}{\partial \upsilon}(1,\tau)=
\int_0^1{\mathcal P}_{\delta_0^\tau}^{-1}\circ{\mathcal P}_{\gamma_\tau^\varrho}^{-1}\circ R_G(\varrho,\tau)\circ{\mathcal P}_{\gamma_\tau^\varrho}\circ{\mathcal P}_{\delta_0^\tau} \  \ d\,\varrho\, .
\end{equation}
which gives formula (\ref{formula15}) by integrating in $\tau$ from $0$ to $1$.
\par Formula (\ref{genangle}) is given by a change of variables using notation \ref{MC1033}.
\par Concerning the next statement on the {\em generalised angle}, if a homotopy like $GG$ exists, it can be realised in a differentiable way from ${\mathbb R}^3$ to $M$, and it describes in ${\bf G}$ a homotopy of paths ${\mathfrak G}_u$
for $u \in [0,1]$ having fixed extremities $\hbox{\rm Id}$ and ${\mathcal P}_\Gamma^{-1}\,$ and linking ${\mathfrak G}_1$ to  ${\mathfrak G}_2$. 
Applying (\ref{MC3}) (lemma \ref{MC3-2}) and Stokes' formula
\begin{equation*}
\int_{{\mathfrak G}_2}A^{-1}\,dA-\int_{{\mathfrak G}_1}A^{-1}\,dA=\int\!\!\int_{{\mathfrak G}{\mathfrak G}}d(A^{-1}\,dA)\  ,
\end{equation*}
one gets (\ref{formula16}).
\par To derive the last statement, first observe that ${\bf O}_g^+(T_xM)$ is commutative because $(M,g)$ is now $2$-dimensional (and $E=TM$). Consider on $\hbox{\rm rge}(G)$ the complex structure induced by $G$ from the one which tells
$i \, \frac{\partial}{\partial \varrho}=\frac{\partial}{\partial \tau}$ on $[0,1]\times[0,1]\,$. In this setting, a rotation in ${\bf O}_g^+(T_xM)$ is naturally associated to a complex number $e^{i\theta}$. The angle $\theta(0)=0$ being associated with the rotation ${\mathcal P}_{\Gamma_{1,0}}^{-1}=\hbox{\rm Id}\,$, a continuous angle $\theta(\tau)$ (as $\tau$ describes $[0,1]$) naturally corresponds 
to the rotation ${\mathcal P}_{\Gamma_{1,\tau}}^{-1}$. Take a vector $w\in T_xM$ of $g_x$-norm $1$ and compute the following expression derived from the left-hand side of (\ref{formula15}) 
\begin{multline} \label{formula17} \int_0^1\langle{\mathcal P}_{\Gamma_{1,\tau}}\circ 
\frac{\partial  {\mathcal P}_{\Gamma_{1,\upsilon}}^{-1}}{\partial \upsilon}(1,\tau)\,(w),i\,w\rangle\,d\,\tau=\cr=
\int_0^1\langle e^{-i\theta(\tau)}\, 
i\,\frac{\partial  \theta(\upsilon)}{\partial \upsilon}(\tau)\,e^{i\theta(\tau)} \,w,i\,w\rangle\,d\,\tau
=\cr=\langle i\,w,i\,w\rangle\int_0^1\, 
\frac{\partial  \theta(\upsilon)}{\partial \upsilon}(\tau)\,d\,\tau
=\theta(1)-\theta(0)=\theta(1)\, ,
\end{multline}
and $\theta(1)$ is the angle corresponding to the rotation ${\mathcal P}_{\Gamma_{1,1}}^{-1}$ (the parallel translation in the reverse direction along the boundary $\Gamma=\partial \,\hbox{\rm rge}(G)$). Starting from the right-hand side of (\ref{formula15}), one gets thanks to (\ref{Courbnorm})
\begin{multline} \label{formula18} 
\int_0^1\!\!\int_0^1 \langle R_G(\varrho,\tau) {\mathcal P}_{A_{\varrho,\tau}}(w),i\,{\mathcal P}_{A_{\varrho,\tau}}(w)\rangle  d\,\varrho\,d\,\tau=\cr=\int\!\!\int_{\hbox{\rm rge}(G)} \frac{\langle R_G(\varrho,\tau) {\mathcal P}_{A_{\varrho,\tau}}(w),i\,{\mathcal P}_{A_{\varrho,\tau}}(w)\rangle}{\Vert dG(\frac{\partial}{\partial \varrho})\wedge dG(\frac{\partial}{\partial \tau})\Vert} \ d\,\hbox{\rm area}_{\,\hbox{\rm rge}(G)}=\cr=-\,\int\!\!\int_{\hbox{\rm rge}(G)} K(G(\varrho,\tau)) \ d\,\hbox{\rm area}_{\,\hbox{\rm rge}(G)}\, ,
\end{multline}
where $K(G(\varrho,\tau))$ is the Gauss curvature of the surface $M$ at the point $G(\varrho,\tau)\in\hbox{\rm rge}(G)$. Quoting that the boundary is parametrised in the reverse direction (with respect to the trigonometric orientation given by the complex structure), putting (\ref{formula17}) and (\ref{formula18}) together, one gets the claimed part of the classical Gauss-Bonnet theorem for a rectangular piece of surface (as explained in remark \ref{GaussBonnet:2}).
\end{proof}

\subsection{On proximity of paths in the linear group.}
\label{proximity33}$ $

We present here various metric equivalences (they arise while comparing different holonomies along paths) and a technical tool (proposition \ref{prox1}), all used in the proof of the main theorem (see section \ref{main33}).
As the applications of the results of this subsection concern the special case $(E_x,\tilde g_x)=(T_xM,g_x)$, we write $(E_x,g_x)$ meaning $(E_x,\tilde g_x)$, for sake of notational simplicity. 

\begin{nota} Denote by $\Vert\,\cdot\,\Vert_x$ the norm on $(E_x, g_x)$. Consider ${\bf G}$ a Lie subgroup (verifying $\hbox{dim}\,{\bf G}\!>\!0$) of ${\bf GL}^+(E_x)\!\subset\!\hbox{\bf End}(E_x)$ (see notation \ref{nota121}). Special cases are 
${\bf G}={\bf O}_{g_x}^+(E_x)$ or ${\bf G}={\bf GL}^+(E_x)$. 
\par If $e_i$ is a $g_x$-orthonormal basis on $E_x\,$, define the scalar product $\overline g_x$ on $\hbox{\bf End}(E_x)$ by asking $E_i^j=e_j\otimes e_i^\ast$ to be $\overline g_x$-orthonormal. It induces on ${\bf G}\subset \hbox{\bf End}(E_x)$ a natural Riemannian structure
$({\bf G},\overline g_x)$.
\par Besides the above structure of Riemannian sub-manifold, consider on ${\bf G}$ the homogeneous Riemannian metric $\underline g_x$ defined by left translating (over ${\bf G}$) the restriction of $\overline g_x$ on $T_{\hbox{\rm Id}}{\bf G}=\underline {\bf G}\subset \hbox{\bf End}(E_x)\  $.  
\end{nota}
\begin{rem}\label{isomm} One has $\underline g_x\!=\!\overline g_x$ on ${\bf G}\!=\!{\bf O}_g^+(E_x)$. {\em A direct proof:} the left-action on ${\bf O}_g(E_x)$ extends to the natural left-action $L\in\hbox{\bf End}(E_x)\mapsto h\circ L\in \hbox{\bf End}(E_x)$, sending a $\overline g_x$-orthonormal basis of $\hbox{\bf End}(E_x)$ to another $\overline g_x$-orthonormal basis, and this implies that the notion of $\overline g_x$ and $\underline g_x$-orthonormal basis coincide all over $T{\bf O}_g(E_x)$. 
\par Furthermore, for any $x,y\!\in\! M\,$, the group $({\bf GL}(E_x),\underline g_x)$ {\em is isometric to} $({\bf GL}(E_y),\underline g_y)$ while $({\bf GL}(E_x),\overline g_x)$ {\em is isometric to} $({\bf GL}(E_y),\overline g_y)$. Thus ${\bf c}_1,{\bf c}_2,{\bf c}_3,{\bf r}_1$ defined below {\em are absolute constants.}
\end{rem}
\begin{rem} \label{inv1} For any $A\in {\bf O}_g(E_x)$ and $B\in {\bf End}(E_x)$ one has
\begin{equation*} \Vert\,B-A\,\Vert_{x,x}=\Vert \,A(A^{-1}B-\hbox{\rm Id})\,\Vert_{x,x}=\Vert \,A^{-1}B-\hbox{\rm Id}\,\Vert_{x,x}\  ,
\end{equation*}
and the $\underline g_x$-invariance reads for $A,B\in {\bf GL}(E_x)$
\begin{equation*} d_{\underline g_x}(A,B)=d_{\underline g_x}(A^{-1}B,\hbox{\rm Id}) \  ,\  \  \  d_{\underline g_x}(A^{-1},B^{-1})=d_{\underline g_x}(BA^{-1},\hbox{\rm Id})\,.
\end{equation*}
\end{rem}
\begin{lem} \label{opnorm1} The operator norm of $h\in{\bf GL}(E_x)$, which is acting on $(\hbox{\bf End}(E_x),\overline g_x)$ through
$L\mapsto h\circ L$, is equal to the operator norm $\Vert h\Vert_{x,x}$ of $h$ acting on $(E_x, g_x)$. In case 
$h$ is in ${\bf O}_g(T_xM)$, one has $\Vert h\Vert_{x,x}=1$.
\end{lem}
\begin{proof} Indeed, in the orthonormal basis $E_i^j$, one computes (for $h\in{\bf G}$)
\begin{multline*}\Vert\,h\circ L\,\Vert_{\overline g_x}^2=\sum_{i,j}\sum_k(h_k^j\, L_i^k)^2=\sum_i\sum_{j,k}(h_k^j\, L_i^k)^2=\cr=
\sum_i\Vert\,h\circ L(e_i)\,\Vert_x^2 \leq \sum_i \Vert h\Vert_{x,x}^2\,\Vert\,L(e_i)\,\Vert_x^2=\Vert h\Vert_{x,x}^2
\Vert\,L\,\Vert_{\overline g_x}^2\  .
\end{multline*}
This implies $\Vert h\Vert_{x,x}\geq \Vert h\Vert_{\overline g_x,\overline g_x}\,$. Thus both operator norms are equal for one may consider any vector $u$ to be the image of the linear map $\lambda_u$ sending $e_1$ to $u$ and each other $e_i$ on $0$ and one has $\Vert u\Vert_{x}=\Vert \lambda_u\Vert_{\overline g_x}\,$.
\end{proof}
\begin{rem}\label{convex33} By {\em convex ball}, we mean {\em a strongly convex ball} (see definition {\rm \ref{D.convex}}).
\end{rem}
\begin{lem} \label{nequiv1}
Choose ${\bf r}_1>0$ (remark {\rm\ref{isomm}}) so that $\bar B_1=\overline{B_{\underline g_x}(\hbox{\rm Id},{\bf r}_1)}$ is a convex ball in $({\bf G},\underline g_x)$. 
Set ${\bf c}_1:=\max(a^{-1},b)$ where 
\begin{equation} \label{blue40} a^{-1}:=\sup\{\Vert\,h\,\Vert_{x,x}\mid h\in \bar B_1\}\  \hbox{and}\  b:=\sup\{\Vert h^{-1}\Vert_{x,x}\mid h\in \bar B_1\}\  ,
\end{equation}
so ${\bf c}_1\!\geq\!1$ (and ${\bf G}\!\subset\!{\bf O}_g(E_x)\!\Leftrightarrow\!{\bf c}_1\!=\!1$).
For any $h\!\in\! \bar B_1$ and $u\!\in\! T_h{\bf GL}(E_x)$ 
\begin{equation} \label{ineq1} {\bf c}_1^{-1}\,\Vert\,u\,\Vert_{\overline g_x}\leq \Vert\,u\,\Vert_{\underline g_x}\leq {\bf c}_1\,\Vert\,u\,\Vert_{\overline g_x}\  ,
\end{equation}
thus, the Riemannian distances $d_{\underline g_x}\!$  of $(\!\bar B_1\!,\!\underline g_x\!)$ and $d_{\overline g_x}\!$ of $(\!\bar B_1\!,\!\overline g_x\!)$ are equivalent (see definition {\rm \ref{quasi-is}} and its concluding remark).
\end{lem}
\begin{rem} Again, ({\rm\ref{ineq1}}) implies that $\underline g_x$ and $\overline g_x$ coincide
on ${\bf O}_g(E_x)$.
\end{rem}
\begin{proof} If $u$ belongs to $T_h{\bf G}$, one has by lemma \ref{opnorm1}
\begin{gather}\label{blue41}
\Vert\,u\,\Vert_{\underline g_x}=\Vert\,h^{-1}\circ u\,\Vert_{\overline g_x} \leq \Vert h^{-1}\Vert_{x,x}\, \Vert\,u\,\Vert_{\overline g_x}\  ,
\\ \label{blue42}
\Vert u \Vert_{\overline g_x}\!=\!\Vert h\circ h^{-1}\circ u \Vert_{\overline g_x} \leq \Vert\,h\,\Vert_{x,x}\,\Vert h^{-1}\circ u\,\Vert_{\overline g_x}\!=\!\Vert\,h\,\Vert_{x,x}\,\Vert u\,\Vert_{\underline g_x}\,,
\end{gather}
and (\ref{blue41}),(\ref{blue42}) together with (\ref{blue40}) produce (\ref{ineq1}). 
\par Integrating (\ref{ineq1}) and computing the Riemannian distance in each metric space from its arc-length's definition, the last claim follows.
\end{proof}
\begin{lem}\label{nequiv1022} There exists ${\bf c}_1'\!>\!0$ verifying
for any $h_1,h_2\!\in\! \bar B_1$
\begin{equation} \label{ineq2} ({\bf c}'_1)^{-1}\,d_{\underline g_x}(h_1,h_2) \leq  \Vert\, h_2-h_1\,\Vert_{\overline g_x} \leq {\bf c}_1\,d_{\underline g_x}(h_1,h_2)\  .
\end{equation}
\end{lem}
\begin{proof}
The right-hand side of (\ref{ineq2}) comes from the isometric embedding of $({\bf G},\overline g_x)$ in $(\hbox{\bf End}(E_x),\overline g_x)$ plus lemma \ref{nequiv1}, last claim. 
\par As for the left-hand side in (\ref{ineq2}), the function 
\begin{equation*}q(h_1,h_2):=d_{\overline g_x}(h_1,h_2)/\Vert\, h_2-h_1\,\Vert_{\overline g_x}\ ,
\end{equation*} 
defined to be $1$ if $h_1=h_2\,$, is continuous on $\bar B_1\times\bar B_1$ and achieves a maximum $a_1\geq1\,$: this relies on the classical interplay between intrinsic and extrinsic metrics on an embedded compact submanifold in a Euclidean space.
Setting ${\bf c}'_1= a_1\,{\bf c}_1$ and using the last claim of lemma \ref{nequiv1} completes the proof of lemma \ref{nequiv1022}.
\end{proof}
\begin{lem} \label{nequiv2}$ $ Here ${\bf G}$ is ${\bf O}_{g_x}^+(E_x)$ or ${\bf GL}^+(E_x)$ and ${\bf r}_1\!>\!0$ is as in lemma {\rm \ref{nequiv1}} ($\bar B_1\!=\!\overline{B_{\underline g_x}(\hbox{\rm Id},{\bf r}_1})\!\subset\! ({\bf G},\underline g_x)$ is convex).
\begin{itemize}
\item
There exists ${\bf r}'_1\!\in]0,{\bf r}_1]$ such that, for any $h\!\in\! {\bf G}$ 
\begin{equation} \label{ineq+} \Vert h-\hbox{\rm Id}\Vert_{\overline g_x} \leq {\bf r}'_1\ \ \hbox{implies}\ \ 
d_{\underline g_x}(h,\hbox{\rm Id}) \leq {\bf r}_1\ ,\ \hbox{thus}\ \ h\in \bar B_1\  .
\end{equation}
\item One has for any $h\in {\bf G}$
\begin{equation}\label{ineq2bisbis}\Vert h-\hbox{\rm Id}\Vert_{x,x}\leq \frac{{\bf r}'_1}{\sqrt{n}}\ \ \ \hbox{implies}\ \ \ 
d_{\underline g_x}(h,\hbox{\rm Id})\leq {\bf r}_1\,, \ \hbox{thus}\ h\in\bar B_1\ .
\end{equation}
\item
There exists ${\bf c}_2\!\geq\!1$ such that,
for any $h_1,h_2\!\in\!\bar B_1$ one has
\begin{equation} \label{ineq2bis} {\bf c}_2^{-1}\,d_{\underline g_x}(h_1,h_2) \leq  \Vert\, h_2-h_1\,\Vert_{x,x} \leq {\bf c}_2\,d_{\underline g_x}(h_1,h_2)\  .
\end{equation}
\item
And {\rm(\ref{ineq2bis})} holds for $h_1\!\in\! {\bf O}_{g_x}^+(E_x),h_2\!\in\!{\bf GL}^+(E_x)$ assuming {\em this time} $h_1^{-1}h_2\!\in\!\bar B_1\,$, where $\bar B_1$ is the ${\bf r}_1$-ball in $({\bf GL}^+(E_x),\underline g_x)\,$.
\end{itemize}
\end{lem}
\begin{rem} \label{clim} If $h_1,h_2$ are in $\bar B_1$ convex, then $d_{\underline g_x}(h_1,h_2)$ is achieved by a unique $\underline g_x$-geodesic $\subset\bar B_1$, so $d_{\underline g_x}(h_1,h_2)=d_{\underline g_x}(h_1^{-1}h_2,{\rm Id})$ equals the distance between $h_1$ and $h_2$ in the metric subspace $(\bar B_1,{\underline g_x})$.
\end{rem}
\begin{proof} The first claim is clear since $\underline g_x$ and $\Vert\cdot-\ast\Vert_{\bar g_x}$ both induce metrics on ${\bf G}$ (the first is Riemannian, the second is by restriction), thus their respective metric balls define equivalent neighborhood systems. 
\par 
The second claim follows from the first and
from the classical equivalence of the norms $\Vert \,\cdot\,\Vert_{x,x}$ and $\Vert \,\cdot\,\Vert_{\overline g_x}$ on $\hbox{\bf End}(E_x)$
\begin{equation}\label{equiv222}\frac{1}{n}\Vert\cdot\Vert_{\overline{g}_x}^2\leq\Vert\cdot\Vert_{x,x}^2\leq\Vert\cdot\Vert_{\overline{g}_x}^2\ .
\end{equation}
\par  
Lemma \ref{nequiv1022} and (\ref{equiv222}) give the third claim (set ${\bf c}_2:=\sqrt{n}\,{\bf c}'_1\,$).
\par
If  $h_1,h_2\in\bar G$ verify  $h_1^{-1}\,h_2$ belongs to $\bar B_1$,
the third claim implies
\begin{equation*} \label{ineq2bis1033} {\bf c}_2^{-1}\,d_{\underline g_x}(h_1^{-1}\,h_2,\hbox{\rm Id}) \leq  \Vert\, h_1^{-1}\,h_2-\hbox{\rm Id}\,\Vert_{x,x} \leq {\bf c}_2\,d_{\underline g_x}(h_1^{-1}\,h_2,\hbox{\rm Id})\  ,
\end{equation*}
and the last claim follows from remark \ref{inv1}.
\end{proof}
Related to the continuity of $h\in {\bf G}\mapsto h^{-1}\in{\bf G}$, recall the
\begin{lem}\label{continv} If $h_1, h_2\in {\bf G}$ verify
$\Vert h_1^{-1}\Vert_{x,x}\,\Vert h_2-h_1 \Vert_{x,x}<1$, one has 
\begin{equation*} \Vert\, h_2^{-1}-h_1^{-1}\,\Vert_{x,x}\leq \,\frac{\Vert h_1^{-1}\Vert_{x,x}^2\,\Vert h_2-h_1 \Vert_{x,x}}{1-\Vert h_1^{-1}\Vert_{x,x}\,\Vert h_2-h_1 \Vert_{x,x}}\  .
\end{equation*}
\end{lem}
\begin{proof} First, write
\begin{multline} \label{continv33}h_2^{-1}-h_1^{-1}=(h_1+(h_2-h_1))^{-1}-h_1^{-1}=\\=
((\hbox{\rm Id}+h_1^{-1}(h_2-h_1))^{-1}-\hbox{\rm Id})\,h_1^{-1}\,.
\end{multline}
As $\Vert h_1^{-1}( h_2-h_1) \Vert_{x,x} <1$, one gets lemma \ref{continv} from (\ref{continv33}), writing
\begin{multline*} \Vert h_2^{-1}-h_1^{-1}\Vert_{x,x}\leq
\Vert h_1^{-1}\Vert_{x,x}^2\Vert h_2-h_1\Vert_{x,x}(\sum_{l=0}^\infty 
\Vert h_1^{-1}\Vert_{x,x}^l\Vert h_2-h_1\Vert_{x,x}^l)
=\\= \frac{\Vert h_1^{-1}\Vert_{x,x}^2\Vert h_2-h_1\Vert_{x,x}}{1-\Vert h_1^{-1}\Vert_{x,x}\Vert h_2-h_1\Vert_{x,x}}\,.\qedhere
\end{multline*}
\end{proof}
\begin{prop} \label{prox1} Given $\epsilon\in]0,{\bf r}_1/2]$, with ${\bf r}_1>0$ like in lemma {\rm \ref{nequiv1}} (thus $\bar B_1=\overline{B_{\underline g_x}(\hbox{\rm Id},{\bf r}_1)}$ is convex), there exists ${\bf c}_3>0$ such that, if ${\mathfrak G}_1,{\mathfrak G}_2 : [0,1] \rightarrow {\bf G}$ are two paths such that ${\mathfrak G}_1(0)={\mathfrak G}_2(0)$ and 
\begin{equation} \label{prox3} 
 \forall t\in[0,1]\   \  \  \  \  d_{\underline g_x}({\mathfrak G}_1(t),{\mathfrak G}_2(t))\leq \epsilon\  ,
\end{equation}
then one has
\begin{equation*} \Vert\!\int_{{\mathfrak G}_2}\!\!\!A^{-1}\,dA-\!\int_{{\mathfrak G}_1}\!\!\!A^{-1}\,dA\,\Vert_{x,x}
\leq ({\bf c}_3\,(\hbox{\rm length}_{\underline g_x}\!({\mathfrak G}_1)+\hbox{\rm length}_{\underline g_x}\!({\mathfrak G}_2))+1)\,\epsilon\  .
\end{equation*}
\end{prop}
\begin{proof} Thanks to (\ref{prox3}), each point ${\mathfrak G}_1(t)$ is linked to ${\mathfrak G}_2(t)$ in $({\bf G},\underline g_x)$
by a unique $\underline g_x$-geodesic parametrised by $s\in[0,1]$ (with constant speed), constituting in this way a homotopy ${\mathfrak G}{\mathfrak G}$ that links ${\mathfrak G}_1$ to ${\mathfrak G}_2$ with fixed origin ${\mathfrak G}{\mathfrak G}(s,0)={\mathfrak G}_1(0)={\mathfrak G}_2(0)$. Define ${\mathfrak G}_3(s)$ to be ${\mathfrak G}{\mathfrak G}(s,1)$ (so ${\mathfrak G}_3(0)\!=\!{\mathfrak G}_1(1),{\mathfrak G}_3(1)\!=\!{\mathfrak G}_2(1))$. This curve ${\mathfrak G}_3$, as any other $\underline g_x$-geodesic ${\mathfrak G}{\mathfrak G}(s,t)={\mathfrak G}_t(s)$, has $\underline g_x$-length $\leq\epsilon$.
Bring in (\ref{formula16}) from theorem \ref{theorem B}
\begin{equation*} 
\int_{{\mathfrak G}_2}\!A^{-1}dA-\!\int_{{\mathfrak G}_1}A^{-1}dA\!=\!\int_{{\mathfrak G}_3}\!A^{-1}dA -\!\int\!\!\int_{{\mathfrak G}{\mathfrak G}}(A^{-1}dA)\wedge(A^{-1}dA)\, .
\end{equation*}
\begin{lem} \label{prox5}
One has
\begin{equation*}\Vert\,\int_{{\mathfrak G}_3}A^{-1}\,dA\,\Vert_{x,x}\leq \hbox{\rm length}_{\underline{g}_x}({\mathfrak G}_3)\leq\epsilon    \  .
\end{equation*} 
\end{lem}
\begin{proof} Start from
\begin{equation*} \int_{{\mathfrak G}_3} A^{-1}\,dA=
\int_0^1 {\mathfrak G}_3^{-1}(t)\,\frac{\partial {\mathfrak G}_3}{\partial t}(t)\,d\,t\  .
\end{equation*}
As $(T_{\hbox{Id}}{\bf G},\bar g_x)=(T_{\hbox{Id}}{\bf G},\underline{g}_x)$, lemma \ref{prox5} follows from (\ref{equiv222}), using the fact that ${\mathfrak G}_3^{-1}(t)$ is an isometry
on $({\bf G},\underline g_x)$
\begin{multline*} \Vert\int_{{\mathfrak G}_3} A^{-1}\,dA\,\Vert_{x,x}\leq
\int_0^1 \Vert{\mathfrak G}_3^{-1}(t)\,\frac{\partial {\mathfrak G}_3}{\partial t}(t)\Vert_{\overline g_x}\,d\,t
=\\ =\int_0^1 \Vert{\mathfrak G}_3^{-1}(t)\,\frac{\partial {\mathfrak G}_3}{\partial t}(t)\Vert_{\underline g_x}\,d\,t
=
\int_0^1 \Vert\frac{\partial {\mathfrak G}_3}{\partial t}(t)\Vert_{\underline g_x}\,d\,t\leq \epsilon\  .\qedhere
\end{multline*}
\end{proof}
\begin{lem} \label{prox6}  There exists ${\bf c}_3>0$ such that
\begin{equation*}\Vert\!\int\!\!\!\int_{{\mathfrak G}{\mathfrak G}}\!(A^{-1}\,dA)\wedge(A^{-1}\,dA)\Vert_{x,x}\!\leq\! {\bf c}_3\,(\hbox{\rm length}_{\underline g_x}\!({\mathfrak G}_1)+\hbox{\rm length}_{\underline g_x}\!({\mathfrak G}_2))\,\epsilon    \  .
\end{equation*}
\end{lem}
\begin{proof} 
Use (see appendix \ref{diantre44bis})
\begin{equation} \label{diantre44} \Vert \,{\mathfrak G}{\mathfrak G}^\ast (A^{-1}\,dA)\wedge(A^{-1}\,dA)\,\Vert_{x,x}\leq2\sqrt{n}\ 
\Vert\frac{\partial {\mathfrak G}{\mathfrak G}}{\partial s}\wedge\frac{\partial {\mathfrak G}{\mathfrak G}}{\partial t}\Vert_{\underline g_x}\ ,
\end{equation}
and write
\begin{multline*} 
\Vert\!\int\!\!\!\int_{{\mathfrak G}{\mathfrak G}}\!(A^{-1}\,dA)\wedge(A^{-1}\,dA)\,\Vert_{x,x}
\leq\\\leq2\sqrt{n}\int_0^1\!\!\!\!\int_0^1\!\!\Vert \frac{\partial {\mathfrak G}{\mathfrak G}(s,t)}{\partial s}\wedge\frac{\partial {\mathfrak G}{\mathfrak G}(s,t)}{\partial t}\Vert_{\underline g_x} ds\,dt= 2\sqrt{n}\ \hbox{\rm area}_{\underline g_x}\!\!(\hbox{\rm rge}({\mathfrak G}{\mathfrak G}))\  .
\end{multline*}
The proof of lemma \ref{prox6} then results from the following lemma \ref{prox8} applied to the homotopy ${\mathfrak G}{\mathfrak G}$ in $(N,h)=({\bf G},\underline g_x)$.
\begin{lem} \label{prox8} Let $D\subset(N,h)$ be a compact domain in a Riemannian manifold and ${\mathcal R}>0$ be smaller than the convexity radius of the ${\mathcal R}$-neighborhood of $D\,$. A real $\epsilon\in]0,{\mathcal R}/2]$ is given. Let ${\mathfrak G}_1,{\mathfrak G}_2: [0,1]\rightarrow D$  be paths such that, for any $t\in[0,1]$, one has $d_h({\mathfrak G}_1(t),{\mathfrak G}_2(t))\leq \epsilon$.
There exist ${\bf c}>0$ depending upon $D\subset(N,h)$, independent of ${\mathfrak G}_1,{\mathfrak G}_2$ and a homotopy 
$H: (s,t)\in[0,1]\times[0,1]\mapsto H(s,t)\in M$
from ${\mathfrak G}_1$ to ${\mathfrak G}_2$ through $h$-geodesics ${\mathfrak G}_t=H(\cdot,t)$ such that
\begin{equation*} \hbox{\rm area}_{h}(\hbox{\rm rge}(H))\leq {\bf c}\,(\hbox{\rm length}_h({\mathfrak G}_1)+\hbox{\rm length}_h({\mathfrak G}_2))\,\epsilon\  .
\end{equation*}
\end{lem}
\begin{proof} Because of the $\epsilon$-proximity of ${\mathfrak G}_1$ and ${\mathfrak G}_2$ (in $(N,h)$), for any given $t$ there exists a unique $g$-geodesic ${\mathfrak G}_t$ joining ${\mathfrak G}_1(t)$ to ${\mathfrak G}_2(t)$: this determines $H$ such that $H(0,t)\!=\!{\mathfrak G}_1(t),H(1,t)\!=\!{\mathfrak G}_2(t), H(\cdot,t)\!=\!{\mathfrak G}_t(\cdot)$. Cut this homotopy into a finite number $\nu$ of pieces $H_i=H_{[0,1]\times[t_i,t_{i+1}]}:[0,1]\times [t_i,t_{i+1}]\rightarrow M$ ($i=0,1,\dots,\nu-1$) whose union is $H$ and which are such that
the sets $\hbox{\rm rge}(H_i)$ are contained in balls $B(p_i,2\epsilon)\subset M$. This can be realised in a way we now explain. 
Take the mid points $p_t$ of the geodesics ${\mathfrak G}_t$ and consider the family of balls $B(p_t,2\epsilon)$. The union of those balls covers $\hbox{\rm rge}(H)$. Moreover, by uniform continuity, there exists $\eta>0$ such that, if $t$ and $t'\in[0,1]$ verify $\vert t-t'\vert\leq \eta$, the mid points $p_t$ and $p_{t'}$ satisfy  $d(p_t,p_{t'})\leq \epsilon$, thus $H_{[0,1]\times[t,t']}: [0,1]\times[t,t']\rightarrow M$ has a range contained in both the {\em convex} balls $B(p_t,2\epsilon)$ and $B(p_{t'},2\epsilon)$. Choosing, for $\nu=\hbox{\rm int}(1/\eta)+1$, a family of  $\nu+1$ points $t_0=0<t_1<\dots<t_{\nu-1}<t_\nu=1$ such that $t_{i+1}-t_i\leq \eta$ completes the construction.

The end of the proof of lemma \ref{prox8}, and thus of the proposition, relies on two more lemmata.
\begin{lem}  \label{prox9} Lemma {\rm \ref{prox8}} is true if $(N,h)$ is a Euclidean vector space.
\end{lem}
\begin{proof} If $(N,h)$ is a Euclidean vector space, the homotopy $H$ reads
\begin{equation*} H(s,t)=(1-s)\,{\mathfrak G}_1(t)+s\,{\mathfrak G}_2(t)\  .
\end{equation*}
So, one has (recall $s\in[0,1]$)
\begin{equation} \label{blue66}\Vert\,\frac{\partial H}{\partial t}(s,t)\,\Vert_h\leq (1-s)\,\Vert\,\frac{\partial {\mathfrak G}_1}{\partial t}(t)\,\Vert_h+
s\,\Vert\,\frac{\partial {\mathfrak G}_2}{\partial t}(t)\,\Vert_h\  ,
\end{equation}
and also 
\begin{equation} \label{blue77}\Vert\,\frac{\partial H}{\partial s}(s,t)\,\Vert_h= \Vert\,{\mathfrak G}_2(t)-{\mathfrak G}_1(t)\,\Vert_h\leq \epsilon\  .
\end{equation}
Thus, by Hadamard's inequality, the area of the range of $H$ satisfies (use (\ref{blue66}) and (\ref{blue77}))
\begin{multline*} \hbox{\rm area}_h(\hbox{\rm rge}(H)) \leq \int_0^1\!\!\int_0^1\Vert\,\frac{\partial H}{\partial s}(s,t)\,\Vert_h\  \Vert\,\frac{\partial H}{\partial t}(s,t)\,\Vert_h\,d\,s\,d\,t
\leq \\\leq \epsilon\,(\int_0^1\!\!\int_0^1(1-s)\Vert\,\frac{\partial {\mathfrak G}_1}{\partial t}(t)\,\Vert_h\,d\,s\,d\,t+
\int_0^1\!\!\int_0^1s\,\Vert\,\frac{\partial {\mathfrak G}_2}{\partial t}(t)\,\Vert_h\,d\,s\,d\,t)=\\=
\frac{1}{2}\,(\hbox{\rm length}_h({\mathfrak G}_1)+\hbox{\rm length}_h({\mathfrak G}_2))\,\epsilon\  .\qedhere
\end{multline*}
\end{proof}
\begin{lem}  \label{prox10} Lemma  {\rm \ref{prox8}} is true on each $B(p_i,2\epsilon)$.
\end{lem}
\begin{proof} Recall from section \ref{EuclRiem}, definition \ref{quasi-is}, proposition \ref{proposition A}: for any $p\in D$, the balls $(B(p,2\epsilon),h)$ and $(B(0_p,2\epsilon),h_p)$ show a metric uniform equivalence through $\exp_p$, i. e. there exists ${\mathcal C}\in[0,1[$ depending only on $D\subset(M,g)$ such that
$(1-{\mathcal C})\,h_p \!\leq\! \exp_p^\ast h \!\leq\!   (1+{\mathcal C})\,h_p\,$, and more  
\begin{equation}  \label{prox12}
(1-{\mathcal C})\,\Vert\,\cdot\,\Vert_p \leq \Vert\,\cdot\,\Vert_{\exp_p^\ast h} \leq   (1+{\mathcal C})\,\Vert\,\cdot\,\Vert_p\  .
\end{equation}
\par
For each $B(p_i,2\epsilon)$, consider $\tilde H_i=\exp_{p_i}^{-1}\circ H_i$, and apply lemma \ref{prox9} to
$\tilde H_i$ having values in the Euclidean space $(T_{p_i}M,h_{p_i})$. One gets with the help of (\ref{prox12}), setting
$\epsilon'=\epsilon/(1-{\mathcal C})$
\begin{multline*} \hbox{\rm area}_{h_{p_i}}(\hbox{\rm rge}(\tilde H_i)) \leq\\\leq
\frac{1}{2}\,(\hbox{\rm length}_{h_{p_i}}(\exp_{p_i}^{-1}\circ {{\mathfrak G}_1}_{\mid [t_i,t_{i+1}]})+\hbox{\rm length}_{h_{p_i}}(\exp_{p_i}^{-1}\circ {{\mathfrak G}_2}_{\mid [t_i,t_{i+1}]}))\,\epsilon'\  ,
\end{multline*}
which gives in turn (again with the help of (\ref{prox12}))
\begin{equation*} \label{prox11} \hbox{\rm area}_h\,\!(\hbox{\rm rge}(\,H_i\,)) \!\leq\!
\frac{(1\!+\!{\mathcal C})^2}{2(1\!-\!{\mathcal C})^2}\,(\hbox{\rm length}_h({{\mathfrak G}_1}_{\mid [t_i,t_{i+1}]})+\hbox{\rm length}_h({{\mathfrak G}_2}_{\mid [t_i,t_{i+1}]}))\,\epsilon\  ,
\end{equation*}
hence the result, putting $\,{\bf c}=(1+{\mathcal C})^2/2(1-{\mathcal C})^2\,$.
\end{proof}
Lemma \ref{prox8} follows by collecting applications of lemma \ref{prox10} to each $H_i$ and by the additive character of area and length.
\end{proof}
Lemma \ref{prox8} can be applied to $(N,h)\!=\!({\bf G},\underline g_x)$ and to ${\mathfrak G}_1,{\mathfrak G}_2,\,H\!=\!{\mathfrak G}{\mathfrak G}$ as in the beginning of the proof of proposition \ref{prox1}. Setting ${\bf c}_3\!:=\!2\sqrt{n}\,{\bf c}$ (see lemma \ref{prox8} for ${\bf c}$) completes the proof of lemma \ref{prox6}. 
\end{proof}
Proposition \ref{prox1} is established.
\end{proof}

\section{Holonomy of a piecewise flat polyhedron and its Regge curvature}
\label{holonomie}\par Basic material in topology can be found in \cite{Gre}, \cite{Mu}, our development points at least to \cite{C-M-S2}, \cite{C-M-S3}, \cite{C-M-S1}, \cite{Chr1}, \cite{Chr2}, \cite{F}, \cite{Ha}, mainly to the fundatorial \cite{R} (see also the chapter 42 on Regge Calculus in \cite{M-T-W}).

\subsection{Polyhedra} \label{polyedres}
\begin{assump} \label{full} Recall that we restrict our concern to finite simplicial complexes $K$ whose topological underlying spaces are homeomorphic to compact connected $n$-dimensional manifolds with boundary (which may be empty), see definition \ref{D.metric133}. We refer to \cite{Mu}, \cite{Gre}, \cite{Ha}. {\em We denote by $\sigma,K...$ a simplex, a simplicial complex... as well as the corresponding underlying topological spaces $\vert\sigma\vert, \ \vert K\vert=\cup_{\sigma\in K}\vert\sigma\vert...$} 
\end{assump}
Recall from definition \ref{D.metric133}
that a {\em polyhedron} is an $n$-dimensional simplicial complex $K$ satisfying assumption {\rm\ref{full}} equipped with an analytic locally Euclidean metric $g_0$ on $K\setminus K_{n\!-\!2}$ for which any $n$-simplex may be naturally seen as isometric to a linear simplex in $\R^n$. 

Denote by $\dagger\!K_k$ the set of all $k$-simplices in $K$ (thus $\dagger\! K_k\subset K_k$).
\begin{defn} Define the {\em pulp} of $K$ to be ${\mathcal K}=K\!\setminus\! (K_{n-2}\cup\partial K)$. Denote by $\nabla$ the Levi-Civita connection of $g_0$ on ${\mathcal K}\,$.
\end{defn}
The local ``affine character'' of $K$ stems from its definition, it provides us with a canonical flat connection $\nabla$ on ${\mathcal K}\,$. But one may vary the metric $g_0$ on ${\mathcal K}$ while the same flat Levi-Civita connection $\nabla$ defines the same $\nabla$-parallel translation along a given curve.
\begin{rem} Observe that, as topological set, ${\mathcal K}\!=\!K\!\setminus\! (K_{n-2}\cup\partial K)$ is connected and is the union of all interiors of $n$-simplices in $K$ and of all relative interiors of the $(n-1)$-simplices sitted in the interior of $K\,$.
\end{rem}
\begin{lem} \label{trpar1} Along two paths homotopic in the flat $({\mathcal K},g_0)$ sharing the same fixed end-points, the $\nabla$-parallel translations coincide. 
\end{lem}
\begin{proof} This is a direct consequence of theorem \ref{theorem A}.
\end{proof}
\begin{defn} In a Riemannian manifold $(M,g)$, given a loop $\gamma$ based at a point $x\in M\,$, the {\rm holonomy} (created by $\gamma$ at $x$) is the parallel translation along $\gamma$ acting as an element of ${\bf O}_{g_x}(T_xM)$ on $(T_xM,g_x)$.
\end{defn}
\begin{rem} The holonomy transformations at $x$ build a subgroup of ${\bf O}_{g_x}(T_xM)$ while those stemming from loops homotopic to $x$ build a closed Lie-subgroup of ${\bf O}_{g_x}^+(T_xM)$ (\cite{Bes} chapter X or \cite{Be4} chapter 13).
\end{rem}
\begin{lem}\label{trpar2} The homotopy group of ${\mathcal K}$, based at a point $p\in {\mathcal K}$, admits a canonical surjective representation onto the holonomy group of ${\mathcal K}$ at $p\,$: one sends a homotopy class of loops based at $p$ to the parallel translation along a loop belonging to this class.
\end{lem}
\begin{proof} This is a direct consequence of lemma \ref{trpar1}.
\end{proof}
\begin{defn} Denote by $\varpi_\sigma$ (by $\varpi_\eta$) the center of gravity of an $n$-simplex $\sigma\,$ (of an $(n-1)$-face $\eta$). The {\rm tree} (of $K$ or of the pulp ${\mathcal K}$) is the graph which is the collection of {\em all} linear arcs joining in $\sigma\,$, for each $\sigma$ and $(n-1)$-face $\eta$ in $\sigma\,$, $\varpi_\sigma$ to $\varpi_\eta\,$. The {\rm graph root} is the (sub)-graph which is the collection of {\em all} linear arcs joining in $\sigma$ the point $\varpi_\sigma$ to $\varpi_\eta$ for any $(n-1)$-face $\eta\subset\sigma$ not included in $\partial K\,$.
\end{defn}
\begin{lem} \label{racine} A full homotopic description of ${\mathcal K}$ is given by the graph root, actually this graph root is a deformation retract of ${\mathcal K}\,$.   
\end{lem}
\begin{proof} Indeed, one retracts in a natural way ${\mathcal K}$ onto the {\rm tree}, and then the {\rm tree} onto the graph root, justifying the assumption above on the homotopic identification of ${\mathcal K}$ with its {\em root}. The first retraction may be performed by a piecewise linear procedure and can be defined on each $n$-simplex $\sigma\in K$ in a simple way, {\em as follows}.
\par Consider each $n$-sub-simplex in $\sigma\in K$ built as the {\em join of an $(n\!-\!2)$-face} $\xi$ and {\em of the segment} $[\varpi_\sigma,\varpi_\eta]$, where $\xi\!\subset\!\partial\eta$ and $\eta$ is an $(n\!-\!1)$-face $\subset\!\sigma\,$. One can define the retraction pushing linearily each $n$-sub-simplex of this type onto $[\varpi_\sigma,\varpi_\eta]$ (gaze at $(1\!-\!t)\xi\!+\!t[\varpi_\sigma,\varpi_\eta]$, for $t\not=0$). As all $\xi\!\in\! K_{n-2}, \xi\!\subset\!\sigma$ are removed from all $\sigma\!\in\! K$ while defining ${\mathcal K}\,$, the linear retractions defined by each $\xi\!\subset\!\sigma$ for any $\sigma$ pile up over ${\mathcal K}$ in the desired retraction onto the tree. The second retraction is clear.
\end{proof}
\begin{lem}  
The complex
 $K$ has {\em a simply connected pulp} ${\mathcal K}$
 if and only if the graph root of $K$ is simply connected.
As a generic loop in $K$ avoids $K_{n-2}\cup \partial K\,$, if $K$ has a simply connected pulp, then $K$ is simply connected (the converse is not true in general). 
\end{lem}
\begin{proof} This is a direct consequence of lemma \ref{racine}.
\end{proof}
\begin{defn} Given an $n$-polyhedron $(K,g_0)$ (definition {\rm\ref{D.metric133}}), call {\em cutting} of $K$ a {\em connected} $n$-polyhedron $(\tilde K,\tilde g_0)$ for which exists a list of triplets $\eta\in \dagger\!K_{n\!-\!1},\eta',\eta''\!\in\!\dagger\!\tilde K_{n\!-\!1}\,$, with $\eta',\eta''\!\subset\!\partial \tilde K$ $\tilde g_0$-isometric through $\phi_\eta\!:\!\eta'\!\rightarrow\!\eta''$, such that, for each triplet $\!\eta,\eta',\eta''\!$, identifying $\eta'$ with $\eta''$ through $\phi_\eta$ gives $\eta\!\not\subset\!\partial K$, building $(K,g_0)\!$ out of $\!(\tilde K,\tilde g_0)$ ($\tilde g_0$ induces $g_0$). 
\vskip2mm
\par {\bf Caution}, we now {\em abusively} write $g_0$ instead of $\tilde g_0\,$.
\end{defn}
\begin{lem} \label{1connpulp1}Every $n$-polyhedron $(K,g_0)$ has a cutting $(\tilde K,\tilde g_0)$ with simply connected pulp. Thus $\tilde {\mathcal K}=\tilde K\setminus (\tilde K_{n-2}\cup\partial\tilde K)$ inherits a $\tilde g_0$-orthonormal parallelisation and $\tilde K_{n-2}$ is included in the boundary $\partial \tilde K\,$. 
\end{lem}
\begin{proof} As the complex $K$ is finite, the homotopy group in ${\mathcal K}$ is finitely generated and its generators can be viewed as loops in the {\em graph root} (${\mathcal K}$ is homotopically equivalent to the connected graph root, lemma \ref{racine}). By cutting such a loop at some interface $\eta$, which means ``cutting'' $K$ along $\eta\,$, i. e. by creating two $(n\!-\!1)$-boundary simplices $\eta'$ and $\eta''$ out of $\eta\,$, producing in this way $\tilde K\,$, the first {\em cutting} of $K$, one kills the corresponding generator in the homotopy group. Precisely, the complex $\tilde K$ is the closure of the metric space $(K\setminus\eta,d)$ for the distance function $d(p,q)$ defined between two points $p,q\in K\setminus\eta$ to be the shortest $g_0$-length of rectifiable paths contained in $K\setminus\eta$ and joining $p$ to $q$, see Gromov \cite{GLP}. The (connected) graph root of $\tilde K$ has strictly less loops than the one of $K\,$: it has been modified (by a cut) only along $\eta$.
This procedure is done until the initial homotopy group of ${\mathcal K}$ {\em has disappeared} (after a finite number of steps), proving the first statement.
\par As for the second statement, because the manifold $\tilde {\mathcal K}$ is simply connected and flat (it is equipped with the flat connection $\nabla$), the result follows from the triviality of the holonomy group, a consequence of lemma \ref{trpar2}. {\em Indeed}, one can parallelise the tangent bundle to $\tilde {\mathcal K}$ by parallel translating a basis of $T_p\tilde {\mathcal K}$ from any point $p$ to any point $q$ along any path from $p$ to $q$, this will produce a well-defined smooth field of tangent bases over the whole of $\tilde {\mathcal K}$. 
\par Finally, observe that any $s\in \dagger\!\tilde K_{n-2}$ with $s\not\subset\partial \tilde K$ would be surrounded by $n$-simplices having $s$ as a face, producing a non trivial loop in the graph root, which is impossible.
\end{proof}
\begin{rem} \label{remarquons} Suppose a $g_0$-orthonormal parallelisation is fixed on a simply connected pulp $\tilde {\mathcal K}\,$, so that, given $p\in \tilde {\mathcal K}\,$, the tangent bundle $T\tilde {\mathcal K}$ decomposes as a true product $T\tilde {\mathcal K}= \tilde {\mathcal K}\times T_p\tilde {\mathcal K}\,$. As in ${\mathbb R}^n\,$, the parallel translation along any curve from $p$ to any $q\in\tilde {\mathcal K}$ identifies in a canonical way $T_p \tilde {\mathcal K}$ with $T_q \tilde {\mathcal K}\,$, {\em this results from the above proof}.
\end{rem}
\begin{rem}\label{orientmess0} A simply connected pulp $\tilde {\mathcal K}$ can be given an orientation, which induces a canonical orientation on each boundary $(n-1)$-simplex (relative to the outward normal, for instance).
\end{rem}
Below, we describe in a cut $\tilde K$ of $K$ with simply connected pulp the $g_0$-parallel translation produced along a curve $\gamma\!\subset\! ({\mathcal K},g_0)\,$. All parts of $\gamma$ that do not cross an $(n\!-\!1)$-simplex along which a cut is done (while building $\tilde K$ out of $K$) contribute through the canonical parallelism depicted in remark \ref{remarquons}, and, for simplicity, we do not mention it explicitly (or as ${\rm Id}$). The non-trivial description is at the cuts along $(n\!-\!1)$-simplices $\eta\,$.
\begin{defn} \label{phieta}  \label{majphieta} Starting from a cutting $\tilde K$ that has a simply connected pulp and rebuilding $K\,$, {\em we first define the mapping $\phi_\eta$ producing an identification along $\eta\,$.} {\em A cut} (proof of lemma {\rm\ref{1connpulp1}}) along $\eta$ (with vertices $p_1,\dots,p_n$) creates two distinct $g_0$-isometric $(n\!-\!1)$-faces $\eta'\!\subset\!\partial \tilde K$ (with vertices $p'_1,\dots,p'_n$) and $\eta''\!\subset\!\partial\tilde K$ (with vertices $p''_1,\dots,p''_n$) that border two different $n$-simplices $\sigma',\sigma''\!\in\!\tilde K\,$. {\em The $g_0$-isometry which identifies back $\eta'$ with $\eta''$ is $\phi_\eta$.} So $\phi_\eta$ sends each $p'_i$ on $p''_i$ and is extended near $\eta'$ so as to be a $g_0$-isometry and to send the $g_0$-unit outward normal $\in T_{\varpi_{\eta'}}\sigma'$ to the $g_0$-unit inward normal $\in T_{\varpi_{\eta''}}\sigma''$. If no cutting was done along $\eta\!=\!\sigma'\cap\sigma''\in K$ to produce $\tilde K\,$, {\em set} $\phi_\eta\!:=\!\hbox{\rm Id}\,$.
\par Identifying  $T_{\varpi_{\eta'}}\sigma'$ with $T_{\varpi_{\eta''}}\sigma''$ through the canonical parallelism in $(\tilde {\mathcal K},\tilde g_0)$ (remark {\rm\ref{remarquons}}), $\Phi_\eta\!:=\!d\phi_\eta(\varpi_{\eta'})$ may be viewed as an element of ${\bf O}_{{g_0},{\varpi_{\eta'}}}(T_{\varpi_{\eta'}}\sigma')$. Due to the canonical parallelism, this element $\Phi_\eta$ can also be viewed, for any $p,q\in\tilde {\mathcal K}\,$, as a $g_0$-linear isometry from $(T_p\tilde {\mathcal K},g_0)$ to $(T_q\tilde {\mathcal K},g_0)$, {\em and we shall now do this tacitly}.
\end{defn}
\begin{rem}\label{orient mess}  The (local) $g_0$-isometry  $\phi_\eta$ in $\tilde {\mathcal K}\,$ is direct if and only if it reverses the {\em induced} orientation from $\eta'$ to $\eta''$ (see remark {\rm\ref{orientmess0}}).
\end{rem}
\begin{lem}\label{trpar3} Consider a cutting $\tilde K$ of $K$ having a simply connected pulp and $\tilde{\mathcal K}=\tilde K\setminus  \partial\tilde K$. A differentiable path in ${\mathcal K}$ crossing once a single $(n-1)$-face $\eta$ produces a parallel translation in $({\mathcal K},g_0)$ which reads $\Phi_\eta$ in
a $g_0$-orthonormal parallelisation of $\tilde{\mathcal K}\,$.
\end{lem}
\begin{proof} Such a path starts in an $n$-simplex $\sigma'$ and ends in another $\sigma''\,$, with $\eta=\sigma'\!\cap\!\sigma''\,$. In each $\sigma'$ and $\sigma''\,$, one picks up the restriction of a chosen $g_0$-orthonormal parallelisation $\chi$ over $\tilde {\mathcal K}\,$. One infers from definition \ref{majphieta} that the $({\mathcal K},g_0)$-parallel translation along such a path reads in  $\chi$ as the $g_0$-isometry $\Phi_\eta$ which goes from $T_{\omega_{\eta '}}\sigma'$ to $T_{\omega_{\eta ''}}\sigma''$ while identifying $\eta'\!\subset \!\sigma'$ to $\eta''\!\subset\!\sigma''$ through $\phi_\eta\,$. Indeed $\chi$ is translated to itself as long as the path is interior to $\sigma'$ or $\sigma''\,$. The change in parallel translation occurs only at the crossing with $\eta\,$. Choose the first $n-1$ vectors of the parallelisation $\chi$ tangent to $T_{\omega_{\eta '}}\eta'\,$ at $\omega_{\eta '}$. In $T_{\omega_{\eta '}}\sigma'\equiv T_{\omega_{\eta ''}}\sigma''$, look at  $T_{\omega_{\eta ''}}\eta''=
\Phi_\eta(T_{\omega_{\eta '}}\eta')$ to infer the conclusion.
\end{proof}
\begin{rem} \label{Deltacomplex1}Thanks to appendix {\rm\ref{Deltacomplex}},
one can also handle the case where two $g_0$-isometric $(n\!-\!1)$-faces $\eta',\eta''$ of the same $n$-simplex are identified to produce the same $\eta$ in $\sigma$.
\end{rem}
\begin{lem} \label{trpar4} Let $(K,g_0)$ be an $n$-polyhedron. Each path $\gamma:[0,1]\rightarrow {\mathcal K}$ is homotopic in ${\mathcal K}$ to a path $\gamma'\!\subset\! {\mathcal K}$ having the same extremities and transverse to $K_{n\!-\!1}$. Let $\eta_1,\dots,\eta_r$ be the $(n\!-\!1)$-faces crossed by $\gamma'$.
In a cutting $\tilde K$ with simply connected pulp, the parallel translation ${\mathcal P}^0_\gamma$ along $\gamma$ from $\gamma(0)$ to $\gamma(1)$ (in $({\mathcal K},g_0)$) has a counterpart in 
$(\tilde{\mathcal K},g_0)$ such that 
$${\mathcal P}^0_\gamma \ \ \ \ {\rm reads}\ \ \ \ \Phi_{\eta_r}\circ\cdots\circ\Phi_{\eta_2}\circ\Phi_{\eta_1}$$ 
in a $g_0$-orthonormal parallelisation of $\tilde{\mathcal K}=\tilde K\setminus (\tilde K_{n-2}\cup\partial \tilde K)\,$.\end{lem}
\begin{proof}
The path $\gamma$ can be retracted (with the whole of ${\mathcal K}$, by the piecewise linear retraction described in the proof of lemma \ref{racine}) onto a path in the graph root, and this - choosing, in the homotopy class, an economic parametrisation of the retracted path, then completing this path in a path joining $\gamma(0)$ to $\gamma(1)$ by connecting its origin by a segment to $\gamma(0)$ and its end in the same way to $\gamma(1)$ - produces a piecewise differentiable path which is itself close to a differentiable $\gamma'$ transverse to $K_{n-1}$, homotopic to $\gamma$ and with same extremities,
so $\gamma$ and $\gamma'$ define the same parallel translation (lemma \ref{trpar1}).
Applying lemma \ref{trpar3} to the successive pieces of
$\gamma'$ meeting a single cut completes the proof.
\end{proof}
\begin{rem}\label{Moebius} Parallel translation, when considered along a loop $\gamma$ in a non orientable $K$, may eventually create an isometry reversing orientation: the Moebius strip is built by the twisted gluing of two opposite sides of a rectangle (which is the  union of two symmetric triangles).
\end{rem}

\subsection{Restricted holonomy in polyhedra}\label{restrholo}
$ $

We want here to describe holonomy transformations which, in the representation sending a loop $\subset{\mathcal K}$ to the $\nabla$-parallel translation it induces (see lemma \ref{trpar2}), correspond to the {\em part} of the homotopy group in ${\mathcal K}$ which is trivial in $K\,$.
In the Moebius strip of remark \ref{Moebius} - or in the Klein bottle which requires a further identification of the two edges left, building the Moebius strip - one can observe that a loop running in the middle cannot be obtained as a composition of loops winding around vertices.
Actually, the {\em part} we shall describe is a subgroup of the holonomy group which can be thought as the ``restricted holonomy group'' of $K$ (see \cite{Bes} chapter X or \cite{Be4} pages {\rm637-ff}). 
\begin{defn} \label{essentialbone} A {\em bone} is any $(n\!-\!2)$-simplex whose relative interior sits in the interior of $K\,$. {\em We borrow this terminology to Regge \cite{R}.} Since we can not circle around $(n\!-\!2)$-simplices sitting in the boundary of $K\,$, we leave them out of this study. 
\end{defn}
\begin{lem} All non-trivial holonomy transformations created by loops homotopically trivial in $K$ are by loops winding around {\em ``bones''}. 
\end{lem}
\begin{proof} A transformation under concern is created by $\nabla$-parallel translating along a loop in ${\mathcal K}$ based at $p$ which is trivial in the homotopy group of $K$ at $p\,$. And ${\mathcal K}$ has the homotopy type of its graph root, i. e. of a {\em ``bouquet'' of circles} (winding around {\em ``bones''}).  
\end{proof}
\begin{defn} \label{ast1} If $s$ is a simplex of $K$, call {\em open star} the set $\hbox{\rm st} (s)$ which is the union of relative interiors of all simplices containing $s$.
\end{defn}
\begin{lem} \label{aster3} The open star of a simplex $s$ whose relative interior sits in the interior of $K$ is homeomorphic to an open ball in ${\mathbb R}^n\,$, thus admits a parallelisation (and an orientation). In particular, this holds for the open star $\hbox{\rm st} (\eta)$ of $\eta\!=\!\sigma'\cap\sigma''\!\in\! \dagger\!K_{n\!-\!1}\,$, where $\sigma',\sigma''\!\in\!K$ are two distinct {\em contiguous} $n$-simplices (see also appendix {\rm\ref{Deltacomplex}}).
\end{lem}
\begin{proof} Recall that we consider simplicial complexes satisfying assumption \ref{full}. Thus, any relative interior of $s$ sitting in the interior of $K$ has an open neighborhood in $K$. On the other hand, $\hbox{\rm st} (s)$ is open in $K$ since it is the complement in $K$ of the union of all (closed) simplices whose relative interiors are not in $\hbox{\rm st} (s)$. 
Hence $\hbox{st}(s)$ is itself an open connected neighborhood in $K$ of the relative interior of $s$. As it retracts on any interior point in $s\,$, it is homeomorphic to an open ball.
\end{proof}
\par We depict the situation around a bone $\xi\,$, the study is done in ${\rm st}(\xi)$ and we first define a basic geometric structure coming out.
\begin{defn} \label{conicdisk}  A {\em conic disk} of {\em angle} $\beta\in{\mathbb R}_+$ is the set $[0,r[\times {\mathbb R}$ ($r$ is the {\em radius}) in which points $(s,\theta)$ and $(s',\theta')$ are identified if $s=s'=0$ (defining the {\em apex} $0$) or if $\theta=\theta'\ {\rm mod}\,(\beta)$, equipped with the locally Euclidean metric $ds^2+s^2d\theta^2\,$, singular at $0$ if $\beta\not=2\pi\,$. A {\em conic disk} is thus a piece of a
$2$-dimensional flat cone.
\end{defn}
\begin{assump} \label{essential333}\label{xi-orientation}  
An arbitrary orientation is fixed on each bone.
\end{assump}
The remark below pieces together facts coming from the inner $g_0$-geometry of successive pairs of contiguous $n$-simplices containing $\xi\,$.
\begin{rem}\label{parallstar}$ $ If $\xi$ is a bone, there exists along paths $\subset \!\hbox{\rm st}(\xi)\!\setminus \!\xi$ a $\nabla$-parallel translation, which is the one defined in $({\mathcal K},g_0)$. The $g_0$-geometry of $\hbox{\rm st}(\xi)\!\setminus \!\xi$ is nice (see \cite{C-M-S2} section 1 for more details): 
\begin{itemize}
\item Around any point $q\in\xi\!\setminus\!\partial\xi\,$, exists a $g_0$-orthogonally splitted neighborhood $U_q\!=\!V_q\!\times\! W_q$ where $V_q$ is an open in $(\xi\setminus\partial\xi,{g_0}_{\mid\xi})$ and $(W_q,{g_0}_{\mid W_q})$ is a small conic disk (definition {\rm\ref{conicdisk}}) of apex $q$ and angle $\beta_\xi=\sum_i \beta_i$, where $\beta_i$ is the interior dihedral angle along $\xi$ of any $n$-simplex $\sigma_i$ having $\xi$ as a face. 
\item Any point $p$ in $\hbox{\rm st}(\xi)\setminus \xi$ belongs to an $n$-simplex $\sigma$ containing $\xi$ and can be connected by a path {\em sitting in} $\sigma$ to a point $q\in\xi\,$. This allows to define in $T_p(\hbox{\rm st}(\xi)\setminus \xi)$ a vector $v$ (a subspace $\xi$) ``parallel'' to $v\in T_q\xi$ (to $\xi\equiv T_q\xi$) and its orthogonal complement $\xi^\perp$: if two contiguous $n$-simplices $\sigma_1,\sigma_2$ contain $\xi$, so do their interface $\sigma_1\cap\sigma_2$. Given $q\in\xi\setminus\partial\xi$, there exists $\xi_q^{\perp}\,$, with $\xi_q^{\perp}\subset\hbox{\rm st}(\xi)\,$, a {\em slice  $g_0$-orthogonal to} $\xi$ {\em at} $q$ {\em extending} $W_q$, such that the collection of $\xi_q^{\perp}$ constitute a foliation of $\hbox{\rm st}(\xi)\,$.
\end{itemize}
One defines the {\em defect angle} $\alpha_\xi\!:=\!2\pi\!-\!\beta_\xi\!\in\!{\mathbb R}\,$. 
{\em For instance} $\alpha_\xi=2\pi/3$ around any vertex in the regular octahedron.
\end{rem}
\begin{lem} \label{Regge1} A loop $\gamma\!:\![0,1]\!\rightarrow\! \hbox{\rm st}(\xi)\!\setminus \!\xi$ {\em circling once} around a bone $\xi$ defines a generator in the homotopy group of ${\mathcal K}$ and is homotopically trivial in $K$. {\em One may orient $\xi^\perp$ according to $\gamma$}. The parallel translation ${\mathcal P}^0_\gamma$ along $\gamma$ (with $p\!=\!\gamma(0)\!=\!\gamma(1)$) is a rotation in $(T_p{\mathcal K},g_0)$ which fixes $\xi$ and rotates in $\xi^\perp\,$ by the {\em defect angle} $\alpha_\xi\!\in\!{\mathbb R}$ independent of $\gamma\,$.
\end{lem}
\begin{rem}\label{renvparcours} Reversed orientations of  $\gamma^{\!-\!1}$ and $\xi^\perp$ give the same $\alpha_\xi\,$.
\end{rem}
\begin{proof} We focus the study on $k:=\hbox{\rm st}(\xi)\!\subset\! K$. First $\hbox{\rm st}(\xi)$ is homeomorphic to the open unit ball $B$ in ${\mathbb R}^n$ (lemma \ref{aster3}). Moreover $\hbox{\rm st}(\xi)\!\setminus \!\xi$ is itself homeomorphic to $B\!\setminus \!(B\!\cap\! {\bf X})$, where ${\bf X}$ is an $(n\!-\!2)$-vector subspace of ${\mathbb R}^n$. Thus $\hbox{\rm st}(\xi)\!\setminus \!\xi$ retracts onto a ``disk deprived from its center'' (in the retraction of $B\!\setminus\!( B\!\cap\! {\bf X})$ given by projecting orthogonally to $B\!\cap\! {\bf X}^\perp$, the point $0$ is left out of its target), this proves the first statement. 
\par  
Replace $\gamma$ by a homotopic loop $\gamma'$ contained in the orthogonal slice $\xi^{\perp}_p\!\subset\! \hbox{\rm st}(\xi)$  through $p$ (remark \ref{parallstar}), with $\gamma'$ such that its central projection on a small circle $\subset\xi^{\perp}_p$ of center $\{q\}\!=\!\xi^{\perp}_p\cap\xi$ {\em is bijective.} This allows to orient any slice $\xi^\perp$ in accordance with $\gamma\,$. Enumerate all $(n\!-\!1)$-simplices $\eta_i$ in $k$ containing $\xi$ ($i=0,\dots,r-1 \ \hbox{\rm mod}\, (r) $) according to their successive unique crossings by $\gamma'\,$. Then $\eta_i$ is the interface of two contiguous $n$-simplices $\sigma_i,\sigma_{i+1}$, which have the bone $\xi$ as common face. Those simplices are circling the bone $\xi$ (see also \cite{C-M-S2} section 1 for further developments). In the simply connected pulp of a cut $\tilde k$ stemming from $k$ (lemma \ref{1connpulp1}), one has $\xi \!\subset\!\partial \tilde k\,$, the cut occuring along a {\em single} $\eta_i\,$, say $\eta_0$ ({\em if not,} the graph root of $k$ {\em would not be} connected and simply connected). 
The parallel translation along $\gamma$ reads in a $g_0$-orthonormal parallelisation of the pulp of $\tilde k$ (lemma \ref{trpar4})
\begin{equation}\label{Blueblues} {\mathcal P}^0_\gamma=\Phi_{\eta_{r-1}}\circ\cdots\circ\Phi_{\eta_1}\circ\Phi_{\eta_0}=\Phi_{\eta_0}\ .
\end{equation} 
Actually,  for $i=1,\cdots,r-1$, one has $\Phi_{\eta_i}={\rm Id}$ if the only cut is done along $\eta_0\,$. More, ${\mathcal P}^0_\gamma$ preserves the global orientation of $\hbox{st}(\xi)$, since parallel translation preserves orientation from one simplex $\sigma_i$ to the next $\sigma_{i+1}$ (lemma \ref{trpar3}) and a global orientation exists in $\hbox{\rm st}(\xi)$.
\par
Observe that ${\mathcal P}^0_\gamma\,$, while computed {\em through a local procedure} in $g_0$-orthonormal parallelisations either of $\tilde{\mathcal K}$ (simply connected pulp of a cutting $\tilde K$ of $K$) or of the pulp of $\tilde k$ (through (\ref{Blueblues})), corresponds to the same element in the holonomy group of $({\mathcal K},g_0)$ at $p\,$. As ${\mathcal P}^0_\gamma={\rm Id}$ on $\xi\subset T_p(\hbox{\rm st}(\xi))$ (remark \ref{parallstar}), it remains to understand the part of ${\mathcal P}^0_\gamma$ which is a $g_0$-rotation in $\xi^\perp_p$. 
One has a notion of parallel vectors all over the simply connected pulp of $\tilde k\,$. Inside a ``developped'' planar cone bounded by two half-lines, identified in one generatrix in the cone (see definition \ref{conicdisk}), follow a curve that gives a loop around $0$ in the cone: using (\ref{Blueblues}), the claimed role of the defect angle $\alpha_\xi=2\pi-\sum_i\beta_i\,$ shows up ($\beta_i$ is the dihedral angle of $\sigma_i$ along $\xi$) in the $g_0$-rotation ${\mathcal P}^0_\gamma$ acting on $\xi^\perp_p$ in $({\mathcal K},g_0)$. {\em Indeed}, $\alpha_\xi\!=\!0$ if and only if the cone is flat.
\end{proof}
Exploiting the notion of direct basis, a choice of orientation on the $(n-2)$-subspace $\xi\subset{\mathbb R}^n$ and on $\xi^\perp$ determines an orientation on ${\mathbb R}^n\,$.
\begin{defn} \label{Reggeoperbone} Given a bone $\xi\!\in\! \dagger\!K_{n\!-\!2}\,, p\!\in\!\hbox{\rm st}(\xi)\!\setminus \xi$ and an orientation of the pair $(\hbox{\rm st}(\xi), \xi)$, {\em define} ${\bf r}_0^\xi$ on $T_p{\mathcal K}\,$ to be the $g_0$-rotation fixing pointwise $\xi\,$, which, in a slice $\xi^\perp\,$, is the $g_0$-parallel translation along a loop $\gamma\subset\xi^\perp\setminus\xi\subset\hbox{\rm st}(\xi)$ based at $p\,$, winding once around $\xi\,$, oriented to give back the orientation of $\hbox{\rm st}(\xi)$ through the orientation it induces on $\xi^\perp$ and the chosen orientation of $\xi$ (lemma {\rm\ref{Regge1}}, assumption {\rm\ref{xi-orientation}}). 
\end{defn} 
\begin{lem} \label{Reggeoperbonebis} 
Since the orthogonal splitting $T_p{\mathcal K}\!=\!\xi\!\oplus\!\xi^\perp$ has a canonical meaning in $\hbox{\rm st}(\xi)\!\!\setminus \!\xi$ (see remark {\rm\ref{parallstar}}), the $g_0$-rotation ${\bf r}_0^\xi$ depends on $\xi$ but not on  $p\!\in\!\hbox{\rm st}(\xi)\!\setminus \xi\,$, nor on $\gamma\,$. {\em Moreover}, given any homotopy 
\begin{gather*} H:(s,t)\in[0,1]\times[0,1]\rightarrow H(s,t)\in{\mathcal K}\ \ \hbox{such that} \\
 p'\!=\!H(0,0),H(s,0)\!=\!H_0(s)\!:=\!\gamma'(s), p\!=\!H(0,1), H(s,1)\!=\!H_1(s)\!:=\!\gamma(s)\,,
\end{gather*} 
the parallel translation along $\gamma'\,$ reads ${\mathcal P}_A^{\!-\!1}\!\circ\!{\bf r}_0^\xi\!\circ\!{\mathcal P}_A\,$, where  $A=H\circ \alpha$ and $\alpha\subset[0,1]\times[0,1]$ is any embedded path from $(0,0)$ to $(0,1)$.
\par Moreover, ${\bf r}_0^\xi$ rotates in the oriented plane $\xi^\perp\subset T_p{\mathcal K}$ by an angle $\alpha_\xi\in{\mathbb R}$ verifying that, changing $\xi$ to $-\xi$ and keeping the ambient orientation of $\hbox{\rm st}(\xi)$ fixed, one gets $\alpha_{-\xi}=\alpha_\xi\,$.
\end{lem}
\begin{proof} Indeed, the splitting $T_p{\mathcal K}=\xi\oplus\xi^\perp$ is unambiguously defined in the whole of $\hbox{\rm st}(\xi)\!\setminus \xi\,$. And if, as above, $H\!\subset\!{\mathcal K}$ is a free homotopy from $\gamma'$ to $\gamma\,$, reorganising $H$ in a homotopy from the loop $A\vee\gamma'\vee A^{-1}$ to the loop $\gamma$ through loops having fixed origin $p\,$, theorem \ref{theorem A} applies.

If one changes $\xi$ to $-\xi$ and keeps the ambient orientation of $\hbox{\rm st}(\xi)$ fixed, one passes from $\gamma$ to $\gamma^{-1}$ (the definition pairs $\gamma$ with ${\bf r}_0^\xi$ and $\gamma^{-1}$ with ${\bf r}_0^{-\xi}$) and from the oriented plane $\xi^\perp$ having orientation induced by $\gamma$ to the same plane with opposite orientation induced by $\gamma^{-1}$.
\end{proof}
\subsection{Playing with holonomy and parametrised squares}$ $\label{play}

This section requires material of part \ref{annexe1033}.
\vspace{-1mm}
\subsubsection{On contacts between parametrised squares and triangulation}$ $
\label{contacted}
\begin{assump} \label{assumptionX} $ $
\begin{itemize}
\item Consider a {\em differentiable simplicial embedding} $T$ (see section {\rm\ref{6.1}}) from $K$ into a manifold $M\,$, having {\em small} simplices. 
The word {\em small} refers implicitly to a Riemannian metric $g$ on $M\,$: for {\em small} simplex $s\in K$ and {\em small} parametrised square $G$ verifying $T(s)\cap  \hbox{\rm rge} (G)\not=\emptyset\,$, one can find a strongly $g$-convex ball containing $T(s)\cup \hbox{\rm rge} (G)$. 
\item For fixed $k\,$, all $k$-faces $\hat\eta$ of all $n$-simplices $\hat\sigma\!\in\! T(K_n)$ appear as subsets {\em of  leaves of a texture} $({\mathcal L},p)$
(see sections {\rm\ref{RBT}, \ref{multiple1088}}, the opening of section \ref{bone1111} or fact {\rm\ref{slip1033}} and the lines before). {\em This holds for all $\hat \xi\!:=\!T(\xi)$ with $\xi\!\in\!\dagger\! K_{n-2}\,$} (set $ \dagger\!\hat K_{n-2}\!:=\!T(\dagger\!K_{n-2})$).
\item From now on, any parametrised square $G$ (definition {\rm \ref{psquare1}}) is assumed to be {\em small} and we suppose that $G$ is, for some $\epsilon_G>0$, an {\bf embedding} of $]-\epsilon_G,1+\epsilon_G[\times
]-\epsilon_G,1+\epsilon_G[$ into $M$.
\end{itemize}
\end{assump}
\begin{defn} \label{psquare3} A parametrised square $G$ with $\Gamma\!:=\!\partial (\hbox{\rm rge} (G))$ (definition {\rm \ref{psquare1}}) {\em hits neatly} $\hat K_{n-2}\!:=\!T(K_{n-2})$ if $\hbox{\rm rge} (G)$ cuts transversally each intersected bone $\hat \xi\!\in\!\dagger\!\hat K_{n-2}$ {\em in finitely many distinct points $\notin\!\Gamma\!\cup \partial \hat\xi$}.
\end{defn}
\begin{defn} \label{psquare4} A square $G$ is {\em neatest} in the data $(T,K,M)$ if {\em there exists at most one $(s,t)$ in} $]0,1[\times]0,1[$ such that $G$ hits {\em neatly} $\hat K_{n-2}$ in the {\em neat point} $x=G(s,t)$ (thus $x\in\hat\xi\setminus(\Gamma\cup\partial\hat\xi)$ for {\em at most one} bone $\xi\in\dagger\!\hat K_{n-2}$). It is {\em trivial} if this intersection is empty.
\end{defn}
\begin{defn} \label{psquare41} Given $G\,$, a point $x\in (\hbox{\rm rge} (G)\cap\hat \xi)\setminus(\Gamma\cup\partial\hat \xi)$, with $\xi\in \dagger\!K_{n-2}\,$, is {\em temperate} if $m_x(\hbox{\rm rge} (G),\hat \xi)$ is finite (assumption {\rm \ref{assumptionX}} and definition {\rm \ref{multiple1044}}). 
A square $G$ is {\em temperate} if $\hbox{\rm rge} (G)\cap \hat K_{n-2}$
has {\em  at most one point} which is moreover {\em temperate}: a {\em neatest} square is {\em temperate}.
\end{defn}
\begin{assump} \label{Regge0} All definitions and statements given in connection with {\em Regge curvatures} (definition below) of a parametrised square $G$, sitting in the image $T(K)$ of a finite simplicial complex $(K,g_0)$ embedded in $M$, are in fact meant in the analytic flat space $({\mathcal K},g_0)$ for $T^{-1}\!\circ \!G([0,1]^2)\!=\!T^{-1} \!({\rm rge} (G))$. Generic {\em parametric families} of squares $G$ have contact points with $\hat K_{n-2}\!:=\!T(K_{n-2})$ {\em showing multiplicity}. 
\par {\em We often write} ${\mathcal P}^0_\Gamma$ {\em instead of} ${\mathcal P}^{g_0}_{T^{-1}\circ\Gamma}$ {\em and} $R_0^G$ {\em often means} $R_0^{T^{-1}\circ G}$.
\end{assump}
\begin{defn} \label{Regge2} The {\em Regge curvature} $R_0^G$ ({\em also denoted} $R_0(\Gamma)$) of a parametrised square $G$ with $x_G\!=\!G(0,0)$ and $\Gamma\!=\!\partial {\rm rge} (G)\!\subset\! T({\mathcal K})$ is, for $p_G\!=\!T^{-1}(x_G)$, the linear mapping below, from $T_p{\mathcal K}$ into itself 
\begin{gather*} R_0^G=R_0(\Gamma):={\mathcal P}^0_{\Gamma^{-1}}-\hbox{\rm Id}=({\mathcal P}^0_{\Gamma})^{-1}-\hbox{\rm Id}\ ,\hskip3mm\hbox{where}\hskip3mm\\
\Gamma^{-1}:=\Gamma_{1,1}^{-1}=\gamma_0^{1,0}\vee\delta_1^{1,0}\vee\gamma_1^{0,1}\vee\delta_0^{0,1}\  
\end{gather*} 
and ${\mathcal P}^0_{\Gamma}$ denotes the $g_0$-parallel translation along the boundary loop $\Gamma_{1,1}=\Gamma$ based at $x_G$ (in fact along $T^{-1}\circ\Gamma$ based at $p_G=T^{-1}(x_G)$), see {\rm(\ref{paths0})} and notation \ref{nota121}. 
\end{defn}
The hypothesis that $G$ {\em hits neatly} $\hat K_{n-2}$ can usefully be precised.
\begin{lem} \label{generic+1} If $G$ {\em hits neatly} $\hat K_{n-2}\,$, there exists a reparametrisation $\grave{G} $ of $G\,$, homotopic and as close as wished to $G\,$, such that, for any $(s,t)$ in $[0,1]\times[0,1]\,$, all curves $\grave{\gamma}_t$ or $\grave{\delta}_s$ meet at most once $\hat K_{n-2}\,$. 
\par If, more generally, $G$ has finitely many points of intersection with $\hat K_{n-2}$ which are all {\em temperate}, one can manage: for any $(s,t)$ in $[0,1]\times[0,1]\,$, all curves $\grave{\gamma}_t$ or $\grave{\delta}_s$ contain at most one of those {\em temperate points}.
\end{lem}
\begin{proof} Suppose that, for a fixed $\tau\,$, the curve $\gamma_{\tau}$ meets $\hat K_{n-2}$ in $m$ values of $s\,$, say at $s=\sigma_0,\dots,\sigma_{m-1}\,$. Around each $(\tau,\sigma_1),\dots,(\tau,\sigma_{m-1})$ exists a rectangle of values $(s,t)$ such that the only intersection of $\gamma_t(s)$ with $\hat K_{n-2}$ corresponds to the selected value. In the rectangle that corresponds to $(\tau,\sigma_1)$, one can $C^\infty$-perturb inside $G$ the family of curves $\gamma_t$ so that the perturbed family $\tilde\gamma_t$ verifies $\gamma_{\tau}(\sigma_1)\not=\tilde\gamma_{\tau}(\sigma_1)$ and remains parametrised according to the untouched family $\delta_s\,$, is unchanged outside of the rectangle, still embedded, and also homotopic to the previous one. Once enough analogous perturbations have carefuly been performed, one gets the desired property for all the curves $\gamma\,$. The same procedure applied to the curves $\delta\,$, keeping the curves $\gamma$ fixed, establishes the claim.
\par The second claim is proved along the same scheme, replacing neat intersections with $\hat K_{n-2}$ by temperate ones (definition \ref{psquare41}).
\end{proof}
\begin{nota} \label{mathcalL} Here ${\mathcal L}$ is associated to the texture whose leaves contain all $(n-2)$-simplices $\hat\xi\in\dagger\!\hat K_{n-2}$ (second point of assumption {\rm\ref{assumptionX}}).
\end{nota}
\begin{lem} \label{generic+2} If $\hbox{\rm rge} (G)$ cuts $\hat K_{\!n\!-\!2}$ out of $ \Gamma\cup\hat K_{\!n\!-\!3}$ and ${\mathcal M}\!\in\!{\mathbb N}$ verifies 
$$\sum_{\{x\in \hbox{\scriptsize\rm rge} (G)\cap\hat K_{n-2}\}}\,m_x( 
\hbox{\rm rge} (G),{\mathcal L})\leq {\mathcal M}\ ,
$$
the square $G$ can be $C^\infty$-approximated by $\acute{G}$ hitting {\em neatly} $\hat K_{n-2}$ in less than ${\mathcal M}$ points. Moreover $\acute{G}$ may be chosen {\em homotopic in $M\,$ to $G$ fixing $\Gamma=\partial G=\partial \acute{G}$} and sharing the enforced neatness of lemma {\rm \ref{generic+1}}.
\end{lem} 
\begin{proof} 
Let $B$ be an open disk in the plane centered at the origin $0$. Choose an open $U\!\subset ]0,1[^2$ containing all points $(s_x,t_x)$ sent to intersection points $G(s_x,t_x)\!=\!x$ in $(\hbox{\rm rge}(G)\!\setminus \!\Gamma)\cap(\hat \xi\!\setminus\!\partial\hat\xi)$, for some $\xi\!\in\! \dagger\!K_{\!n\!-\!2}\,$. Choose a mapping $\check G\in C^\infty(B\times I^2,M)$ verifying $\check G_0(s,t):=\check G(0,s,t)=G(s,t)$ and, for any $b\in B$, $\check G_b=G$ outside of $U\,$, {\em to be transverse} to all $\hat \xi\,$ hit by $\hbox{\rm rge}(G)$ on the compact $\bar B_0\times I^2\,$, where $B_0\subset \bar B_0\subset B$ is an open ball centered at the origin.
Consider the ``small perturbations'' $\check G_b$ of $G$ with $b$ near to $0\,$. By lemma \ref{degleqmult}, there exists an open $U_x\!\subset\! M$
near any $x\in(\hbox{\rm rge}(G)\!\setminus \!\Gamma)\cap(\hat \xi\!\setminus\!\partial\hat\xi)$ such that, {\em if} $\check G_b$ is transverse to $\hat \xi$ (this is true for a dense set of $b\in B$ by lemma 4.6, page 53 in \cite{G-G}), the {\em neat} points in ${\rm rge}(\check G_b)\cap\hat \xi\cap U_x$ are less than $m_x(\hbox{\rm rge}(G),{\mathcal L})$. More, if $b$ is small enough and tends to $0$, all points in ${\rm rge}(\check G_b)\cap\hat K_{n-2}$ congregate around the points $x\in (\hbox{\rm rge}(G)\setminus \Gamma)\cap(\hat K_{n-2}\setminus \hat K_{n-3})$. Setting $\acute G:=\check G_b$ for $b$ small enough, ${\rm rge}(\acute G)\cap\hat K_{n-2}$ has less than ${\mathcal M}$ neat points.
The rest follows from lemma {\rm \ref{generic+1}}.
\end{proof}
\begin{lem} \label{psquare5} Let $G:[0,1]\times[0,1]\longrightarrow M$ be a parametrised square hitting $\hat K_{n-2}$ neatly {\em in the enforced sense of lemma \ref{generic+1}} in ${\mathcal N}$ points, so ${\mathcal N}$ {\em is the number of points in} $({\rm rge}(G)\setminus\Gamma)\cap(\hat K_{n-2}\setminus\hat K_{n-3}$). Denote by $N$ the smallest integer such that ${\mathcal N}\leq N^2$. There exists a decomposition of $G$ into $N^2$ neatest (maybe trivial) parametrised squares ($G_{i,j}$ is the restriction of $G$ to the subsquare $[s_{i,j},s_{i+1,j}]\times[t_j,t_{j+1}]$)  
$$G_{i,j}:[s_{i,j},s_{i+1,j}]\times[t_j,t_{j+1}]\rightarrow M\hskip2mm\hbox{where}\hskip2mm i,j=0,\dots,N-1 \  .
$$
\par If $G$ has ${\mathcal N}'$ {\em temperate} points of intersection with $\hat K_{n-2}$, a corresponding decomposition holds replacing in their roles the neatest by temperate squares and $N$ by the smallest integer $N'$ such that ${\mathcal N}'\leq (N')^2\,$.
\end{lem}
\begin{proof} In view of the hypothesis and of the related assumption preceeding the statement of the lemma to be proved, one
may choose $t_0=0$ and $t_1$ such that $G([0,1]\times[0,t_1])$ contains $N$ points of intersections with $\hat K_{n-2}$, and none on its boundary curve. Then decompose $[0,1]\times[0,t_1]$ into squares
$[s_{i,0},s_{i+1,0}]\times[0,t_1]$ such that 
$$G_{i,0}:[s_{i,0},s_{i+1,0}]\times[0,t_1]\rightarrow M\hskip2mm\hbox{where}\hskip2mm i=0,\dots,N-1 \  
$$
is neatest (using the transversality of $G$ and $\hat K_{n\!-\!2}$ and the definitions \ref{psquare3} and \ref{psquare4}). 
\par Doing the same thing, one chooses $t_2\!>\!t_1$ such that $G([0,1]\!\times\![t_1,t_2])$ contains $N$ points of intersections with $\hat K_{n-2}$... After a finite number of such steps, the proof is complete.
\par The second case can be proved by following the same scheme while replacing neatest squares by temperate ones.
\end{proof}
\subsubsection{The case of $G$ hitting neatly $\hat K_{n-2}:=T(K_{n-2})$.}
\label{sssec:num1}
$ $

\begin{nota} Denote by $\Gamma_{i,j}$ the naturally oriented loop
\begin{equation} \label{paths2}\Gamma_{i,j}:=\partial G_{i,j}= \delta_{s_{i,j}}^{t_{j+1},t_j}\vee\gamma_{t_{j+1}}^{s_{i+1,j},s_{i,j}}\vee\delta_{s_{i+1,j}}^{t_j,t_{j+1}}\vee\gamma_{t_j}^{s_{i,j},s_{i+1,j}}\, ,
\end{equation} boundary of the square $G_{i,j}$ (see lemma {\rm\ref{psquare5}}) and by $\Gamma_j$ the oriented boundary loop of $G([0,1]\times[0,t_{j+1}])$. Finally, set 
\begin{equation}  \label{decloop1}\Delta_{i,j}=
A_{i,j}^{-1}\vee\Gamma_{i,j}^{-1}\vee A_{i,j} \hskip2mm\hbox{where}\hskip2mm A_{i,j}=\gamma_{t_j}^{0,s_{i,j}}\vee\delta_0^{0,t_j}\  .
\end{equation}
\end{nota}
\begin{lem} \label{decloop2} One has
\begin{equation*} \Gamma_{j-1}\vee\Gamma_j^{-1}=  \Delta_{N-1,j}\vee \dots\vee \Delta_{1,j}\vee\Delta_{0,j}\  .
\end{equation*}
\end{lem}
\begin{proof} The reader is invited to draw a scheme.
\end{proof}
\begin{nota} According to ({\rm\ref{paths2}), (\ref{decloop1}}), lemma {\rm\ref{Regge1}} and definition {\rm\ref{Regge2}}, the neatest square $G_{i,j}$ gives rise to the endomorphisms
\begin{equation}\label{decloop3} R_0^{G_{i,j}}\!=\!R_0(\Gamma_{i,j})\!=\!({\mathcal P}^0_{\Gamma_{i,j}})^{-1}-\hbox{\rm Id}\,,\, R_0(\Delta_{i,j})\!=\!({\mathcal P}^0_{A_{i,j}})^{-1}\!\circ\! R_0^{G_{i,j}} \!\circ\! {\mathcal P}^0_{A_{i,j}}\  ,
\end{equation}
where, in coherence with assumption {\rm\ref{Regge0}}, ${\mathcal P}^0_{A_{i,j}}$ stands for ${\mathcal P}^0_{T^{-1}\circ A_{i,j}}\,$. 
\end{nota}
\begin{rem} As the parallel translation with respect to $g_0$ is an isometry, one has ($\Vert \,\cdot\,\Vert_{0,0}$ is the operator norm in the metric $g_0$)
\begin{equation} \label{iso1} \Vert \,R_0^{G_{i,j}}\,\Vert_{0,0}=\Vert\,R_0(\Delta_{i,j})\,\Vert_{0,0} \  .
\end{equation}
\end{rem}
\begin{assump} \label{Blue1000} We assume that there exists a bound ${\mathscr C}$ such that, for all $i,j=0,\dots,N-1 $ one has
\begin{equation} \Vert\,R_0^{G_{i,j}}\,\Vert_{0,0} \leq {\mathscr C} \hskip3mm \hbox{and} \hskip3mm 
N{\mathscr C} \hskip3mm \hbox{belongs to} \hskip3mm [0,1]\  .
\end{equation}
\end{assump}
\begin{prop} \label{Regge4} If assumption {\rm\ref{Blue1000}} holds, one has
\begin{equation*} \Vert\,{\mathcal P}^0_{\Gamma_{j-1}}\circ ({\mathcal P}^0_{\Gamma_j})^{-1}-\hbox{\rm Id}-\sum_{i=0}^{N-1}R_0(\Delta_{i,j})\,\Vert_{0,0} \leq N^2{\mathscr C}^2\  .
\end{equation*}
\end{prop}
\begin{proof} Thanks to lemma \ref{decloop2} and to (\ref{decloop3}), one can write
\begin{multline}\label{binome1088} {\mathcal P}^0_{\Gamma_{j-1}}\circ ({\mathcal P}^0_{\Gamma_j})^{-1}-\hbox{\rm Id}=\cr=(R_0(\Delta_{N-1,j})+\hbox{\rm Id})\circ\cdots\circ(R_0(\Delta_{1,j})+\hbox{\rm Id})\circ(R_0(\Delta_{0,j})+\hbox{\rm Id})-\hbox{\rm Id}\  .
\end{multline}
Hence, one also has
\begin{multline*} \label{Regge-33} {\mathcal P}^0_{\Gamma_{j-1}}\circ ({\mathcal P}^0_{\Gamma_j})^{-1}-\hbox{\rm Id}
=\sum_{i=0}^{N-1}R_0(\Delta_{i,j})+\!\!\!\!\!\!\sum_{0\leq i_1<i_2\leq N-1}\!\!\!\!\!\!R_0(\Delta_{i_2,j})\circ R_0(\Delta_{i_1,j})+\cr+\dots +
R_0(\Delta_{N-1,j})\circ\dots\circ R_0(\Delta_{1,j})\circ R_0(\Delta_{0,j})\  .
\end{multline*}
Recalling (\ref{iso1}) and the assumption \ref{Blue1000}, one gets
\begin{equation*} \Vert{\mathcal P}^0_{\Gamma_{j-1}}\!\circ\! ({\mathcal P}^0_{\Gamma_j})^{-1}-\hbox{\rm Id}-\sum_{i=0}^{N-1}R_0(\Delta_{i,j})\Vert_{0,0}
\!\leq\!\! \sum_{k=2}^N\binom{N}{k}\,{\mathscr C}^k\leq \sum_{k=2}^\infty\,\frac{(N{\mathscr C})^k}{k!}\  .
\end{equation*}
As $z=N{\mathscr C}$ belongs to $[0,1]$, the result follows in view of
\begin{multline*} e^z-1-z=\frac{z^2}{2}+\frac{z^3}{3!}+\dots+\frac{z^{l+1}}{(l+1)!}+\dots \leq \cr
\leq z^2(\frac{1}{2}+\frac{1}{2^2}+\dots+\frac{1}{2^l}+\dots)=z^2\  .\qedhere
\end{multline*}
\end{proof}
\begin{nota}  Set 
\begin{equation}\label{Regge5} B_{i,j}=\gamma_{t_j}^{1,s_{i,j}}\vee\delta_1^{0,t_j}\vee\gamma^{0,1}_{0} \hskip4mm \hbox{and observe that} \hskip4mm \Gamma_{j-1}^{-1}\vee A_{i,j}^{-1}=B_{i,j}^{-1}\  .
\end{equation}
\end{nota}
\begin{prop} \label{Regge6} If assumption {\rm\ref{Blue1000}} holds, one has (notation {\rm (\ref{decloop3})})
\begin{equation*}\Vert\,({\mathcal P}^0_{\Gamma})^{-1}-\hbox{\rm Id}-\sum_{i,j=0}^{N-1}({\mathcal P}^0_{B_{i,j}})^{-1}\circ R_0^{G_{i,j}}\circ{\mathcal P}^0_{A_{i,j}}\,\Vert_{0,0}\leq N^3{\mathscr C}^2\  .
\end{equation*}
\end{prop}
\begin{proof} One can write (defining ${\mathcal P}^0_{\Gamma_{-1}}:= \hbox{\rm Id}$) 
\begin{multline*}({\mathcal P}^0_{\Gamma})^{-1}-\hbox{\rm Id}={\mathcal P}^0_{\Gamma^{-1}}-\hbox{\rm Id}={\mathcal P}^0_{\Gamma_{N-1}^{-1}}-\hbox{\rm Id}=\cr=
{\mathcal P}^0_{\Gamma_{N-1}^{-1}}-{\mathcal P}^0_{\Gamma_{N-2}^{-1}}+{\mathcal P}^0_{\Gamma_{N-2}^{-1}}-\dots
-{\mathcal P}^0_{\Gamma_{1}^{-1}}-{\mathcal P}^0_{\Gamma_{0}^{-1}}+{\mathcal P}^0_{\Gamma_{0}^{-1}}-\hbox{\rm Id}=\cr
=\sum_{j=0}^{N-1}({\mathcal P}^0_{\Gamma_{j-1}})^{-1}\circ({\mathcal P}^0_{\Gamma_{j-1}}\circ ({\mathcal P}^0_{\Gamma_j})^{-1}-\hbox{\rm Id})\  .
\end{multline*}
Thus, bringing in (\ref{Regge5}), because parallel translation with respect to $g_0$ is an isometry, one gets by applying proposition \ref{Regge4} 
\begin{multline*}\Vert ({\mathcal P}^0_{\Gamma})^{-1}-\hbox{\rm Id}-\sum_{i,j=0}^{N-1}({\mathcal P}^0_{B_{i,j}})^{-1}\circ R_0^{G_{i,j}} \circ {\mathcal P}^0_{A_{i,j}}\,\Vert_{0,0}\leq \cr
\leq \sum_{j=0}^{N-1}\Vert ({\mathcal P}^0_{\Gamma_{j-1}})^{-1}\!\circ\!({\mathcal P}^0_{\Gamma_{j-1}}\!\circ\! ({\mathcal P}^0_{\Gamma_j})^{-1}\!-\!\hbox{\rm Id}\!-\!\sum_{i=0}^{N-1}R_0(\Delta_{i,j}))\Vert_{0,0}\leq \!N^3{\mathscr C}^2\  .\qedhere
\end{multline*}
\end{proof}
\begin{prop} \label{ReggeGB} If assumption {\rm\ref{Blue1000}} holds, one has (notation {\rm (\ref{decloop3})})
\begin{equation*} \Vert\,\sum_{j=0}^{N-1}{\mathcal P}^0_{\Gamma_{j-1}}\circ (({\mathcal P}^0_{\Gamma_{j}})^{-1}-({\mathcal P}^0_{\Gamma_{j-1}})^{-1})-\sum_{i,j=0}^{N-1} R_0(\Delta_{i,j}) \,\Vert_{0,0}\leq N^3{\mathscr C}^2\  .
\end{equation*} 
\end{prop}
\begin{proof} Use again proposition \ref{Regge4} and sum up the inequalities derived from its conclusion for $i=0,\dots,N-1$ to get
\begin{multline*} \Vert\,\sum_{j=0}^{N-1}{\mathcal P}^0_{\Gamma_{j-1}}\circ (({\mathcal P}^0_{\Gamma_{j}})^{-1}-({\mathcal P}^0_{\Gamma_{j-1}})^{-1})-\sum_{i,j=0}^{N-1} R_0(\Delta_{i,j})\,\Vert_{0,0}\cr
\leq \sum_{j=0}^{N-1}\Vert\,{\mathcal P}^0_{\Gamma_{j-1}}\circ ({\mathcal P}^0_{\Gamma_j})^{-1}-\hbox{\rm Id}-\sum_{i=0}^{N-1}R_0(\Delta_{i,j})\,\Vert_{0,0} \leq N^3{\mathscr C}^2\  .\qedhere
\end{multline*}
\end{proof} 
\begin{rem} \label{commutativity}   The above propositions rely on shifting from non commutative decompositions of parallel translations along boundary loops of parametrised squares (hitting neatly $\hat K_{n-2}:=T(K_{n-2})$) to approximations by {\em sums} of $g_0$-rotations, making sense in $\hbox{\rm End}(T_{p}{\mathcal K})$ for some $p\in \mathcal K\,$. Beware that the explicit list of the involved rotations keeps the non-commutative character hidden in the paths $A_{i,j}, B_{i,j}$  involved, which allow expressing (via $T^{-1}$) a rotation generated by a curve $\Gamma_{i,j}$ {\em circling around and near} a given bone $\hat\xi:=T(\xi)$ as a rotation acting in $T_{p_G}{\mathcal K}\,$. Indeed, the way to connect a neighboring point of $\hat\xi$ in $\hbox{\rm rge}(G)$ to the base point $x_G=G(0,0)=T(p_G)$ is here by no means canonical if they are several intersection points of $\hbox{\rm rge}(G)$ with different bones, it is even an intimate and unavoidable part of the analysis we make.
\end{rem}

We now introduce a tool which allows expressing neatest (later temperate) Regge curvatures $R_0^{G_{i,j}}$ (later $R_0^{\backslash\!\!\!G_{i,j}}$) in terms of rotations around bones hit by the corresponding neatest $G_{i,j}$ (temperate $\backslash\!\!\!G_{i,j}$).
\begin{rem} \label{smallandneat} Recall that - from assumption {\rm\ref{assumptionX}} 
 - a parametrised square is assumed to be {\em embedded} in $M$. Thus the boundary curve of a neatest parametrised square winds at most once around a given bone $\hat\xi$. Our analysis leads to the definition given below of the {\em index of a neatest square $G$ with respect to $\xi$}, denoted by $\,\hbox{index}_{G,\xi}\,$, related to the {\em intersection number} $\#({\rm rge}\,(G),\hat\xi)$ of ${\rm rge}\,(G)$ with $\hat\xi$, see \cite{Hi}, page 132. Due to the {\em smallness} required in assumption {\rm\ref{assumptionX}},
one can choose an orientation in an open subset $U$ of $M$ containing $G\cup T(\hbox{\rm st}(\xi)\setminus \xi)$.
\end{rem}
\begin{defn} \label{indiceGxi}  Here $G$ is a neatest parametrised square with its canonical orientation and $\hat\xi$ a bone, both meeting at a single point, interior in $\hat\xi$ and $ \hbox{\rm rge}(G)$. If $U$ is an open ball containing $\hat\xi\cup \hbox{\rm rge}(G)$ in $M$ (remark {\rm\ref{smallandneat}}), orientations on $\hat\xi$ and $U$ are fixed, inducing an oriented $\hat\xi^\perp$. One compares the orientations of two local surfaces transverse to $\hat\xi$, like a slice tangent to $\hat\xi^\perp$ and $\hbox{\rm rge}(G)$. Define the {\em index of the neatest square $G$ versus $\xi$} (denoted by ${\rm index}_{G,\xi}$) to be
\par  $+1$ if $\hbox{\rm rge}(G)$ and $\hat\xi^\perp$ inherits the same induced orientation; 
\par  $-1$ if $\hbox{\rm rge}(G)$ and $\hat\xi^\perp$ inherits the opposite induced orientations.
\par\noindent If $\hbox{\rm rge}(G)$ and $\hat\xi$ do not meet or meet on $\partial G$ or $\partial \hat\xi\,$, set $\,{\rm index}_{G,\xi}:=0$. {\em Equivalently,} define ${\rm index}_{G,\xi}:=\#({\rm rge}(G),\hat\xi)$.
\end{defn}
\begin{lem} \label{indice34}
Given a parametrised square $G$ hitting neatly $\hat K_{n-2}$
(definition {\rm\ref{psquare3}}) in the enforced sense of lemma {\rm \ref{generic+1}}, using the decomposition $G_{i,j}$ of lemma {\rm\ref{psquare5}} one {\em extends} the above definition, {\em defining} the index of $G$ versus $\xi\!\in\! \dagger\!K_{n\!-\!2}$ (in an oriented neighborhood of $\hbox{\rm rge}(G)\!\subset\!M$)
\begin{equation} \label{indice332}\hbox{\rm index}_{G,\xi}:=\sum_{i,j}\hbox{\rm index}_{G_{i,j},\xi}\  .
\end{equation}
\par $(i)$This index of $G$ versus $\xi\in \dagger\!K_{n-2}$ is the {\em intersection number} of ${\rm rge}\,(G)$ versus $\hat\xi\,$, one has $\hbox{\rm index}_{G,\xi}=\#(\hbox{\rm rge}(G),\hat\xi)$. 
\par $(ii)$ If $G$ and $G'$ are parametrised squares hitting neatly $\hat K_{n-2}$ verifying $\hbox{\rm rge}(G)\!=\!\hbox{\rm rge}(G')$  sharing the same orientation, one has for $\xi\! \in \!\dagger\!K_{n-2}$
\begin{equation*} \hbox{\rm index}_{G,\xi}= \hbox{\rm index}_{G',\xi}\ .
\end{equation*}
\par $(iii)$ 
For any non-trivial $G_{i,j}$, one can find a loop $\partial G_l$ (with $l\!=\!l(i,j)$) bounding a non trivial neatest square $G_l$ homotopic in $\hbox{\rm rge}(G)\setminus \hat\xi$ to $G_{i,j}$ (sharing the same orientation), where $G_1,G_2,\dots,G_{\mathcal N}$ are enumerated in the right order so that $\partial G$ is homotopic in $\hbox{\rm rge}(G)\setminus \hat\xi$ to
\begin{equation}\label{indice333}\partial G_{\mathcal N}\vee\dots\vee\partial G_2\vee\partial G_1\ ,
\end{equation}
 (see lemma {\rm\ref{psquare5}} for ${\mathcal N}$) and  {\em then}
\begin{equation}\label{indice334} \hbox{\rm index}_{G,\xi}=\sum_{l=1}^{{\mathcal N}}\hbox{\rm index}_{G_l,\xi}\  .
\end{equation}
\end{lem}
\begin{proof} In the decomposition of lemma {\rm \ref{psquare5}}, the usual {\em infinitesimal} orientation on $G$ induces a canonical orientation on the neatest squares $G_{i,j}$. Bringing in the properties of the {\em intersection number} of $\hbox{\rm rge}(G)$ and $\hat\xi\,$ (see chapter 5, subsection 2 of \cite{Hi}) establishes claim $(i)$.
\par For $(ii)$, observe that $G$ and $G'$ are {\em embeddings} (assumption \ref{assumptionX}) sharing the same oriented images $\hbox{\rm rge}(G)\!=\!\hbox{\rm rge}(G')$, thus showing the same intersection points and intersection numbers.
\par As for $(iii)$, $\hbox{\rm rge}(G)\setminus \hat \xi$ has the homotopy type of a bouquet of circles, so, collecting the non trivial neatest $G_{i,j}$ achieves a composition (\ref{indice333}) of $\partial G\,$. One has $\hbox{\rm index}_{G_{i,j},\xi}\!=\!\hbox{\rm index}_{G_{l},\xi}$ (see definition \ref{indiceGxi}) since $G_l$ and $G_{i,j}$ are homotopic in $\hbox{\rm rge}(G)\setminus \hat\xi\,$, implying (\ref{indice334}) (use (\ref{indice332})).
 \end{proof}
The next lemma rephrases definition \ref{Reggeoperbone} (see also lemma \ref{Reggeoperbonebis}).
\begin{lem} \label{Reggeoperbone1} 
Given a bone $\xi\,$, given $x\in T(\hbox{\rm st}(\xi)\setminus \xi)$ with $p=T^{-1}(x)$, {\em there exists} a neatest square $G$ verifying $G(0,0)=x\,,\  \hbox{\rm index}_{G,\xi}=1$ and 
\begin{equation}\label{anytime}{\mathcal P}_\Gamma^0= {{\bf r}_0^\xi}\hskip3mm \hbox{and} \hskip3mm R_0^G=({{\bf r}_0^\xi})^{-1}-\hbox{\rm Id}\ ,
\end{equation}
where ${\bf r}_0^\xi$ is the $g_0$-rotation in $T_p{\mathcal K}$ of definition {\rm\ref{Reggeoperbone}}.
A square $G$ such that ${\rm rge}(G)\!\subset \!T(\hbox{\rm st}(\xi))),\,G(0,0)\!=\!x,\,\hbox{\rm index}_{G,\xi}\!=\!1$ verifies {\rm (\ref{anytime})}.
\end{lem}
\begin{proof} To prove the existence, select, $g_0$-orthogonal to $\xi\,$, a slice $\Sigma\subset\hbox{\rm st}(\xi)$ and a {\em well-oriented} circle $\Sigma$ centered at $\Sigma\cap\xi$ winding once around $\xi$ and included in $\hbox{\rm st}(\xi)\setminus\xi$. The image by $T$ of this circle is the boundary of a {\em well-oriented} disk contained in $T(\Sigma)$, which may be homotoped to a square $G\subset T(\hbox{\rm st}(\xi))$ such that $G(0,0)=x$ and $ \hbox{\rm index}_{G,\xi}=1$ (thus $G$ is {\em neatest}). In the set of neatest squares $G$ such that $G(0,0)=x$ and $ \hbox{\rm index}_{G,\xi}=1$ exists a homotopy connecting any pair of them: {\em indeed}, every such a square is homotopic in $T(\hbox{\rm st}(\xi))$ to a well oriented small disk centered at a point of $\hat\xi$ and contained in a slice orthogonal to $\hat\xi\,$. Coherently, $G$ such that ${\rm rge}(G)\!\subset \!T(\hbox{\rm st}(\xi))),\,G(0,0)\!=\!x,\,\hbox{\rm index}_{G,\xi}\!=\!-\!1$ verifies ${\mathcal P}_\Gamma^0\!=\! ({{\bf r}_0^\xi)^{\!-\!1}}$ and $ R_0^G\!=\!{\bf r}_0^\xi\!-\!\hbox{\rm Id}$ and (\ref{anytime}) shows that $R_0^G$ is independent of $p\,$ for any neatest $G$ such that ${\rm rge}(G)\subset T(\hbox{\rm st}(\xi)))$.
\end{proof}

\subsubsection{Now $G$ may hit $\hat K_{n-2}$ in finitely many temperate points.}
\label{sssec:num2}$ $

The importance of dealing with temperate points - i. e. {\em points of finite multiplicity} (in the sense of definition \ref{multiple1044}) - comes up from the need of a {\em good} control, for generic $G$, on the number ${\mathcal N}$ of points in $\hat K_{n-2}\cap {\rm rge}(G)$ while making smaller and smaller the mesh of $(T,K,M)$ (through the local theorem of Thom, see \cite{C-M}). This was central in theorem \ref{theorem I}, on which relies the convergence we shall see in section \ref{convhol1}, and is also related to the genericity of temperate points of intersection between {\em families} of parametrised squares and bones
(i. e. points of intersection which may show multiplicity). Using transversality, the results of section \ref{sssec:num1} are extended to this case. 

We first {\em define} the {\em index} of a {\em temperate} square.
\begin{lem}\label{singular00} Let $G$ be temperate, so there exists at most a bone $\xi_0$ such that ${\rm rge}(G)$ hits $\hat\xi_0$ in a single point outside of $\Gamma\cup\partial\hat\xi_0$ with {\em finite multiplicity} (definition {\rm\ref{psquare41}}). For generic squares $\acute{G}$ tending $C^\infty$ to $G$ with $\Gamma\!=\!\partial{\rm rge}(\acute{G})\!=\!\partial {\rm rge}(G)\,$, decomposed in neatest squares $\acute{G}\!=\!\acute{G}_{{\mathcal N}_{\acute{G}}}\vee\cdots\vee\acute{G}_1$ (lemma {\rm \ref{indice34}}), the indices $\hbox{\rm index}_{\acute{G},\xi}=\sum_l\hbox{\rm index}_{\acute{G}_l,\xi}$ tend to a limit {\em denoted by} $\hbox{\rm index}_{G,\xi}$ depending only on $G\,$.
\end{lem}
\begin{rem} This extension of the index may readily be adapted to parametrised squares $G$ showing {\em a finite number of temperate points}.
\end{rem}
\begin{proof} If $G$ is trivial, there is nothing to prove (in this case $R_0^G=0$). So, let $\xi_0$ be as in the statement. A square $\acute{G}$ transverse to  $\hat K_{n-2}$ (not meeting $\hat K_{n-3}$), with $\partial G=\partial \acute{G}\,$, is such that $\hbox{\rm rge}(\acute{G})$ hits neatly $\hat K_{n-2}$. Since, with respect to some given Riemannian metric on $M\,$, the distance between $\hbox{\rm rge}(G)$ and $\partial\hat\xi_0$ is always $\geq a>0\,$, one also has $\hbox{\rm rge}(\acute{G})\cap\partial\hat\xi_0=\emptyset\,$ for all $\acute{G}$ in some neighborhood ${\mathcal U}\subset C^\infty([0,1]^2,M)$ of $G\,$. 
More  (see \cite{Hi}, theorem 2.1, page 131-2, then lemmas 1.1, 1.2 and corollary 1.3 on pages 122-3), {\em any two squares $\acute{G}$ and $\acute{G}'$ near enough to $G$ are homotopic through a homotopy $H$ verifying} $\hbox{\rm rge}(H)\subset{\mathcal U}\,$ (no $H_t$ hits $\partial\hat\xi_0$) and one may choose $H:[0,1]\times[0,1]^2\rightarrow M$ to be {\em transverse} to $\hat\xi_0\,$, so $H^{-1}\hat\xi_0$ is a union of paths in $[0,1]^3$ with end-points in $\{0\}\times]0,1[^2\,\cup\,\{1\}\times]0,1[^2\,$;
checking the orientations, one concludes that $\hbox{\rm index}_{\acute{G},\xi_0}$ is constant for all $\acute{G}$ near enough to $G\,$. As, for $\xi\in \dagger\!K_{n-2}$ with $\xi\not\!=\!\xi_0$ and $\acute{G}$ near enough to $G\,$, one has $\hbox{\rm index}_{\acute{G},\xi}\!=\!\hbox{\rm index}_{G,\xi}\!=\!0$ and, for any $\xi\!\in\! \dagger\!K_{n-2}$, one thus can {\em set} $\hbox{\rm index}_{G,\xi}\,:=\,\lim_{\acute{G}\rightarrow G}\hbox{\rm index}_{\acute{G},\xi}\,$.  
Lemmas {\rm \ref{indice34}}, \ref{Reggeoperbone1} imply the claim
\end{proof}
\par
To express the Regge curvature of a temperate square $G$ with ${\rm rge}(G)\!\subset\! T({\rm st}(\xi))$ (definition \ref{psquare41}) in terms of ${\bf r}_0^\xi$, we introduce the 
\begin{defn}\label{eurk} If $G$ (${\rm rge}(G)\!\subset\! T({\rm st}(\xi))$) is a temperate square, set 
\begin{gather*}s_{G,\xi}\!:=\!0\ \ \hbox{if}\ \ {\rm index}_{G,\xi}\!=\!0\ ,\ s_{G,\xi}:={\rm index}_{G,\xi}/\vert {\rm index}_{G,\xi}\vert \ \hbox{if}\ \ {\rm index}_{G,\xi}\not=0\ ,\\ \ \hbox{and}\ \ {\bf r}_0^{\xi,G}:=s_{G,\xi}\ (({{\bf r}_0^\xi})^{-s_{G,\xi}}-\hbox{\rm Id})\ .
\end{gather*} 
\end{defn} 
\begin{lem} \label{banalbof} Given a neatest square $G$ such that ${\rm rge}(G)\!\subset\! T({\rm st}(\xi_0))$ for some bone $\xi_0\,$, one has $R_0^G\!\!=\!\!s_{G,\xi_0}\,{\bf r}_0^{\xi_0,G}\!\!=\!\!\sum_\xi {\rm index}_{G,\xi}\,{\bf r}_0^{\xi,G}$. 
\par If  $G$ is {\em temperate} and ${\rm rge}(G)\!\subset\! T({\rm st}(\xi_0))$, then $s_{G,\xi_0}\,{\bf r}_0^{\xi_0,G}$ is the Regge curvature $R_0^{G'}$ induced by a {\em neatest} $G'$ satisfying $s_{G',\xi_0}\!=\!s_{G,\xi_0}\,$.
\end{lem}
\begin{proof} Go back to definitions {\rm\ref{Reggeoperbone}}, {\rm\ref{eurk}}, then use lemma \ref{Reggeoperbone1}.
\end{proof}
In coherence with lemma \ref{Reggeoperbonebis}, one introduces the
\begin{defn} \label{cruelmanque} If $G_{i,j}$ is {\em temperate} and hits $\xi_0\,$, define 
\begin{equation} \label{temperate1033} {\bf r}_0^{\xi_0,G_{i,j}}:=s_{\xi_0,G_{i,j}}
\ ({\mathcal P}^0_A)^{-1}\!\circ(({\bf r}_0^{\xi_0})^{-s_{\xi_0,G_{i,j}}}-{\rm Id})\!\circ \!{\mathcal P}^0_A\ ,
\end{equation}
where $A$ is any embedded path in ${\rm rge}(G_{i,j})\!\setminus\!\hat\xi_0$ from $G(s_{i,j},t_j)$ to some point $x'\!=\!G(s,t)\in T(\hbox{\rm st}(\xi_0)\!\setminus \!\xi_0)$ (two such paths $A$ and $A'$ produce the same ${\bf r}_0^{\xi_0,G_{i,j}}$ through (\ref{temperate1033}) since they produce homotopic loops at $G(s_{i,j},t_j)$ in ${\rm rge}(G_{i,j})\!\setminus\!\hat\xi_0\,$, see also lemma \ref{Reggeoperbonebis} and its proof).
\end{defn}
\begin{lem} \label{geomversion} Given a square $G$ with $\hbox{\rm rge}(G)$ transverse to $\hat K_{n-2}$ and $\hat K_{n-3}$ and $\partial\,\hbox{\rm rge}(G)\cap \hat K_{n-2}\!=\!\emptyset$ (thus $\hbox{\rm rge}(G)\cap\hat K_{n-3}=\emptyset$ and $\hbox{\rm rge}(G)$ hits neatly $\hat K_{n-2}$, definition {\rm\ref{psquare3}}), then $R_0({\Delta_{i,j}})\!=\!({\mathcal P}^0_{A_{i,j}})^{-1}\!\circ\! R_0^{G_{i,j}} \!\circ\! {\mathcal P}^0_{A_{i,j}}$  (defined in {\rm(\ref{decloop3})}) acts on $T_{p_G}{\mathcal K}$ with $p_G=T^{-1}(x_G)$ and reads
\begin{equation} \label{babe1}R_0({\Delta_{i,j}})\!=\!\!\sum_{\xi\in \dagger\!K_{n-2}}{\rm index}_{G_{i,j},\xi}\ {\bf r}_0^\xi(\Delta_{i,j})
\ ,
\end{equation}
where, if $\ {\rm index}_{G_{i,j},\xi}\not=0\,$, one defines ${\bf r}_0^\xi(\Delta_{i,j})\!:=\!({\mathcal P}^0_{A_{i,j}})^{-1}\circ{\bf r}_0^{\xi,G_{i,j}}\!\circ \!{\mathcal P}^0_{A_{i,j}}$ (definition {\rm\ref{cruelmanque}}) and, if $\ {\rm index}_{G_{i,j},\xi}=0\,$, one sets ${\bf r}_0^\xi(\Delta_{i,j}):=0\,$.
\end{lem}

\begin{proof} As $G$ is transverse to $\hat K_{n-2}$ and to $\hat K_{n-3}$, its range has a finite number ${\mathcal N}$ of intersection points with bones $\hat\xi$ for $\xi\in \dagger\!K_{n-2}$, points which are all of multiplicity $1$. As the square $G_{i,j}$ is neatest (lemma {\rm \ref{psquare5}}), the result follows from definition \ref{indiceGxi}, lemma \ref{Reggeoperbonebis}, lemma \ref{banalbof} (a bone $\xi$ gives rise to an exponent $=-1,0$ or $1$ in (\ref{babe1})).
\end{proof}

\begin{lem}\label{singular1} Let $G$ be a temperate square (see definition {\rm\ref{psquare41}}), $\xi_0$ as in lemma {\rm\ref{singular00}} and $\{y\}\!=\!{\rm rge}(G)\cap\xi_0$. For $\acute{G}$ near to $G$
and $A\!\subset\!{\rm rge}(G)\!\setminus\!\xi_0$ embedded path from $x_G\!=\!G(0,0)$ to $x'\!=\!G(s,t)\!\in \!T(\hbox{\rm st}(\xi_0)\!\setminus \!\xi_0)$ (lemma {\rm\ref{Reggeoperbonebis}}), one has 
\begin{gather} \label{blue-333}
{\mathcal P}^0_{G}=\lim_{\acute{G}\rightarrow G}{\mathcal P}^0_{\acute{G}}
 = ({\mathcal P}^0_{A})^{-1}\circ({\bf r}_0^{\xi_0})^{\hbox{\tiny\rm index}_{G,\xi_0}}
 \circ{\mathcal P}^0_{A}\  ,
\\ \label{blue-334}R_0^G=({\mathcal P}^0_{A})^{-1}\circ(({\bf r}_0^{\xi_0})^{-\hbox{\tiny\rm index}_{G,\xi_0}}
 -\hbox{\rm Id})\circ{\mathcal P}^0_{A}\ .
\end{gather}
Recalling definition {\rm\ref{multiple1044}}, notation {\rm\ref{mathcalL}} and assumption {\rm\ref{assumptionX}} gives
\begin{equation}\label{bluette1033}\vert \hbox{\rm index}_{G,\xi_0}\vert\leq m_{y} (\hbox{\rm rge}(G),\hat\xi_0)\leq m_{y} (\hbox{\rm rge}(G),{\mathcal L})\,.
\end{equation} 
{\em Set} ${\mathscr C}_G \!:=\!\Vert {\bf r}_0^{\xi_0,G}\Vert_{0,0}$
(definition {\rm\ref{cruelmanque}}), $\backslash\!\!\!{\mathscr C}_G \!:=\!\vert \hbox{\rm index}_{G,\xi_0} \vert\,{\mathscr C}_G$, one has
\begin{equation}\label{bluette} 
\Vert R_0^G-\!\!\sum_\xi\hbox{\rm index}_{G,\xi}\,{\bf r}_0^{\xi,G}\Vert_{0,0}\!\leq\! \backslash\!\!\!{\mathscr C}_G^2\hskip3mm\hbox{if}\hskip3mm
\backslash\!\!\!{\mathscr C}_G\!\leq\!1\ .
\end{equation}
\end{lem}
\begin{proof} If $G$ is trivial, there is nothing to prove (in this case $R_0^G=0$). So, let $\xi_0$ be {\em the} (unique) bone hit by $G$ at {\em the} (unique) point $y$ (with $y\notin\Gamma\cup\partial\hat\xi_0\,$) having {\em finite multiplicity}.  
\par By lemma \ref{generic+2}, $G$ may be approached by $\acute G$ hitting neatly $K_{n-2}$ in less that 
$m_{y} (\hbox{\rm rge}(G),\hat\xi_0)$ neatest points and one gets (\ref{bluette1033}) in view of definitions \ref{multiple1044}, \ref{intdegree1044}, (\ref{indice332}) and lemma \ref{degleqmult}, since 
$$\vert \hbox{\rm index}_{G,\xi_0}\vert\leq\hbox{deg}_{y}(\hbox{\rm rge}(G),{\mathcal L})\leq m_{y} (\hbox{\rm rge}(G),\hat\xi_0)\leq m_{y} (\hbox{\rm rge}(G),{\mathcal L})\,.
$$
A binomial expansion (\ref{blue-333}) similar to the one done in the lines below (\ref{binome1088}) gives (\ref{bluette}). {\em Indeed}, write from (\ref{blue-334}) and definitions \ref{eurk}, 
\ref{cruelmanque}
\begin{multline*}{\mathcal P}^0_{A}\circ R_0^G\circ ({\mathcal P}^0_{A})^{-1}=({\bf r}_0^{\xi_0})^{-s_{G,\xi_0}\,\vert {\rm index}_{G,\xi_0}\vert}-{\rm Id}=\\=((({\bf r}_0^{\xi_0})^{-s_{G,\xi_0}}-{\rm Id})+{\rm Id})^{\vert {\rm index}_{G,\xi_0}\vert}-{\rm Id}={\rm index}_{G,\xi_0}\,{\bf r}_0^{\xi_0,G}+ \ {\it other}\ {\it terms}\,,
\end{multline*}
and the remaining terms verify $\Vert {\it other}\ {\it terms}\Vert_{0,0}\leq \backslash\!\!\!{\mathscr C}_G^2$ if $\backslash\!\!\!{\mathscr C}_G\!\leq\!1\,$.
\end{proof}
One now extends the preceeding propositions {\rm \ref{Regge4},  \ref{Regge6}, \ref{ReggeGB}}.
\begin{assump} \label{bounded22} From now on, we consider parametrised squares $G$ whose total amount of intersections with {\em any} single bone $\xi$, {\em  counting multiplicities}, is bounded by a fixed integer $\mu\,$, meaning we assume
$$\sum_{\{x\in \hbox{\scriptsize\rm rge}(G)\cap \hat \xi\}}\,m_x(\hbox{\rm rge}(G), \hat \xi)\leq \mu
$$
(thanks to theorem {\rm\ref{text-1133}}, $\mu\!=\!m_{n\!-\!2}(\!(n\!-\!1)n\!+\!c)$ fits definition \ref{testsquare1} $(ii)$).
\end{assump}
\begin{lem}\label{generalisation1} Assume $G$ is such that $\Gamma\cap \hat K_{n\!-\!2}\!=\!\emptyset$ and $\hbox{\rm rge}(G)\!\cap\!\hat K_{n\!-\!2}$ consists in ${\mathcal N}'$ temperate points, while 
{\em counting multiplicities} those points are ${\mathcal N}$ (lemma {\rm \ref{psquare5}}): if $N'$ is the smallest integer such that ${\mathcal N}'  \!\leq\! (N')^2 $, setting $\backslash\!\!\!\Delta_{i,j}\!:=\!A_{i,j}^{-1}\vee\backslash\!\!\!\Gamma_{i,j}^{-1}\vee A_{i,j}$ where $\backslash\!\!\!\Gamma_{i,j}$ is the boundary of the temperate square $\backslash\!\!\!G_{i,j}$ in a decomposition of 
$G$, if $A_{i,j}$ and $B_{i,j}$ keep their meanings (see {\rm(\ref{paths2})}, {\rm (\ref{decloop1})} and {\rm(\ref{Regge5})}), {\em one generalises} {\rm(\ref{decloop3})} 
\begin{equation}\label{decloop3333} R_0(\backslash\!\!\!\Delta_{i,j})=({\mathcal P}^0_{A_{i,j}})^{-1}\circ R_0^{\backslash\!\!\!G_{i,j}}\circ {\mathcal P}^0_{A_{i,j}} \ ,
\end{equation}
{\em and gets} for $i,j=0,1,\dots,N'-1$, according to {\rm (\ref{babe1})} and {\rm (\ref{bluette1033})} in lemmas {\rm\ref{geomversion}}, {\rm\ref{singular1}} and definitions {\rm\ref{eurk}}, {\rm\ref{cruelmanque}} (see lemma {\rm\ref{singular1}} for $\backslash\!\!\!{\mathscr C}_{\backslash\!\!\!G_{i,j}}$)
\begin{gather}\label{mygush}
\,{\bf r}_0^\xi({\backslash\!\!\!\Delta_{i,j}}):=({\mathcal P}^0_{A_{i,j}})^{\!-1}\!\circ\! {\bf r}_0^{\xi,{\backslash\!\!\!G_{i,j}}}\circ {\mathcal P}^0_{A_{i,j}}\ ,
\\ \label{mygush33}\Vert R_0({\backslash\!\!\!\Delta_{i,j}})-\!\!\sum_\xi\hbox{\rm index}_{\,\backslash\!\!\!G_{i,j},\xi}\,{\bf r}_0^\xi({\backslash\!\!\!\Delta_{i,j}})\Vert_{0,0}\leq \backslash\!\!\!{\mathscr C}_{\backslash\!\!\!G_{i,j}}^2\ \ \hbox{if}\ \ \backslash\!\!\!{\mathscr C}_{\backslash\!\!\!G_{i,j}}\leq 1\ ,\\
 \label{fantastic}\sum_{\xi\in K_{n-2}}\sum_{i,j=0}^{N'-1}
\!\!\vert \,\hbox{\rm index}_{\,\backslash\!\!\!G_{i,j},\xi}\,\vert \leq {\mathcal N} \leq {\mathcal N}'\,\mu
\ .
\end{gather}
\end{lem}
\begin{proof} Formula (\ref{fantastic}) follows from $\vert \,\hbox{\rm index}_{\,\backslash\!\!\!G_{i,j},\xi}\,\vert\leq {\bf n}_{\,\backslash\!\!\!G_{i,j}\,\cap\,\xi}\,$, where ${\bf n}_{\,\backslash\!\!\!G_{i,j}\,\cap\,\xi}\,$ is the global amount {\em with multiplicities} of points in $\backslash\!\!\!G_{i,j}\,\cap\,\xi\,$.
\end{proof}
\begin{assump}[see assumption \ref{Blue1000}] \label{hyp555} Set ${\mathscr C}'\!:=\!2\mu{\mathscr C}$. Assume $G$ is as in lemma {\rm\ref{generalisation1}} and one has (for $i,j\!=\!0,\dots,N'-1$ and any $\xi$)
\begin{equation*} \Vert {\bf r}_0^{\xi,\backslash\!\!\!G_{i,j}}\Vert_{0,0}\!=\!\Vert{\bf r}_0^{\xi}({\backslash\!\!\!\Delta_{i,j}})\Vert_{0,0}={\mathscr C}_{\backslash\!\!\!G_{i,j}} \!\leq\!{\mathscr C} \hskip2mm \hbox{and} \hskip2mm 
N'{\mathscr C}' \!\in\! [0,1] \hskip1mm (\hbox{thus} \hskip1mm {\mathscr C}\!\!\leq\!\! {\mathscr C}'\!\!\leq \!\!1)\,.
\end{equation*}
\end{assump}
\begin{rem} \label{plouc} Using $\Vert {\bf r}_0^{\xi,\backslash\!\!\!G_{i,j}}\Vert_{0,0} \!\!\leq\!\! {\mathscr C}$, ({\rm\ref{decloop3333}}), ({\rm\ref{mygush33}}) (assumptions {\rm\ref{bounded22}}, {\rm\ref{hyp555}}), one gets
$\Vert R_0(\backslash\!\!\!\Delta_{i,j})\Vert_{0,0}\! \leq \!\!\Vert \!\sum_\xi\!\!\hbox{\rm index}_{\,\backslash\!\!\!G_{i,j},\xi}\,{\bf r}_0^\xi({\backslash\!\!\!\Delta_{i,j}})\Vert_{0,0} \!+\! \backslash\!\!\!{\mathscr C}_{\backslash\!\!\!G_{i,j}}^2\!\!\leq\!\! {\mathscr C}'$.
\end{rem}
\par \noindent{\em Under assumptions {\rm\ref{bounded22}}, {\rm \ref{hyp555}}, the three propositions below are true.}
\begin{prop}({\bf extended}  {\rm\ref{Regge4}})\label{Regge4bis} One has
\begin{equation*} \Vert\,{\mathcal P}^0_{\backslash\!\!\!\Gamma_{j-1}}\circ ({\mathcal P}^0_{\backslash\!\!\!\Gamma_j})^{-1}-\hbox{\rm Id}-\sum_{i=0}^{N'-1}R_0(\backslash\!\!\!\Delta_{i,j})\,\Vert_{0,0} \leq N'^2{\mathscr C}'^2\  .
\end{equation*}
\end{prop}
\begin{prop}({\bf extended} {\rm\ref{Regge6}})\label{Regge6bis} One has
\begin{equation*}\Vert\,({\mathcal P}^0_{\Gamma})^{-1}-\hbox{\rm Id}-\sum_{i,j=0}^{N'-1}({\mathcal P}^0_{B_{i,j}})^{-1}\circ R_0^{\backslash\!\!\!G_{i,j}} \circ {\mathcal P}^0_{A_{i,j}}\,\Vert_{0,0}\leq N'^3{\mathscr C}'^2\  .
\end{equation*}
\end{prop}
\begin{prop}
({\bf extended}  {\rm\ref{ReggeGB}}) \label{ReggeGBbis} One has 
\begin{equation*} \Vert\!\sum_{j=0}^{N'-1}\!{\mathcal P}^0_{\backslash\!\!\!\Gamma_{j-1}}\circ (({\mathcal P}^0_{\backslash\!\!\!\Gamma_{j}})^{\!-\!1}-({\mathcal P}^0_{\backslash\!\!\!\Gamma_{j-1}})^{-1})-\!\sum_{i,j=0}^{N'-1} R_0(\backslash\!\!\!\Delta_{i,j})\Vert_{0,0}\!\leq\! N'^3{\mathscr C}'^2\,.
\end{equation*} 
\end{prop}
\begin{proof} All proofs of those propositions are entirely similar to the original ones. But this time, making a decomposition into $N'^2$ {\em temperate squares}, one uses assumption \ref{hyp555} and remark \ref{plouc}.
\end{proof} 
\subsubsection{A word on the case $\Gamma=\partial \,\hbox{\rm rge}(G)$ contains points of $\hat K_{n-2}$}\label{sssec:num3}
$ $

\begin{lem} \label{psquare44} If the square $G$ verifies $\Gamma$ contains points of $\hat K_{n-2}\,$, but $\hbox{\rm rge} (G)\setminus\Gamma$ has only finitely many temperate points of intersection with $\hat K_{n-2}\,$, {\em one can still define} the parallel translation along the curve $\Gamma\!=\!\partial\, \hbox{\rm rge}(G)$ {\em extended relatively to $G$}: for $\epsilon\in]0,1[$, $G_\epsilon:=G_{\mid [\epsilon,1-\epsilon]\times[\epsilon,1-\epsilon]}$ and $\Gamma_\epsilon:=$ the oriented $\partial G_\epsilon\!=\!
\delta_\epsilon^{1-\epsilon,\epsilon}
\gamma_{1-\epsilon}^{1-\epsilon,\epsilon}
\delta_{1-\epsilon}^{\epsilon,1-\epsilon}
\gamma_\epsilon^{\epsilon,1-\epsilon}$, {\em up to a conjugation} $\lim_{\epsilon\rightarrow0}
{\mathcal P}_{\Gamma_\epsilon}^0$ makes sense (see the proof). {\em Set} ${\mathcal P}_\Gamma^0\!:=\!
\lim_{\epsilon\rightarrow0}
{\mathcal P}_{\Gamma_\epsilon}^0\,$.
\end{lem}
\begin{rem}
One can still decompose $G_\epsilon$ along the lines of lemma {\rm\ref{generalisation1}},
choosing the smallest integer $N'$ such that ${\mathcal N}'\leq (N')^2$ and producing temperate squares $\backslash\!\!\!G_{\epsilon;i,j}\,$: the parallel translation along the curve $\backslash\!\!\!\Gamma_{\epsilon;i,j}$, boundary of a corresponding {\em temperate} parametrised square $\backslash\!\!\!G_{\epsilon;i,j}\,$, defines a corresponding Regge curvature $R_0^{\backslash\!\!\!G_{\epsilon;i,j}}\,$.
\end{rem}
\begin{proof} Lemma \ref{psquare44} consists in seeing $G$ as limit of a family of increasing parametrised squares $G_\epsilon$, having boundary curves $\Gamma_\epsilon\,$, which are loops starting and ending at $x(\epsilon)=G(\epsilon,\epsilon)$. The number of points in $\hat K_{n-2}\cap(\hbox{\rm rge}(G)\setminus\Gamma)$ is finite, so the parallel translation $({\mathcal P}_{x(\cdot)}^0)^{-1}{\mathcal P}_{\Gamma(\cdot)}^0{\mathcal P}_{x(\cdot)}^0$ is seen, up to a conjugation by the $g_0$-parallel translation along the curve $p(\cdot):=T^{-1}(x(\cdot))$ (say between any values $\epsilon_1$ and $\epsilon_2$ {\em small enough} with $0<\epsilon_1<\epsilon_2$),
to be {\em ``constant'' in} $\hbox{\rm End}(T_{p(\cdot)}{\mathcal K})\,$,  which proves that the definition of ${\mathcal P}_\Gamma^0$ given in the statement is coherent.
\end{proof}
\begin{rem} In the sequel, we often skip this case, dealing implicitly with the above extension of the parallel translation.
\end{rem}

\subsection{Dihedral angles and Regge curving}
\subsubsection{Bounds on dihedral angles in the Euclidean space}\label{dihedral33}
\begin{lem}\label{diedre1} Each Euclidean $n$-simplex of thickness $\geq t_0$ has dihedral angles in $[{\bf a},\pi-{\bf a}]$, where ${\bf a}\in ]0,\pi[$ depends only on $n$ and $t_0$.
\end{lem}
\begin{proof} If the lemma does not hold, a sequence of $n$-simplices of diameter $1\,$, thickness $\geq t_0\,$, has a dihedral tending to $0\,$ (or to $\pi\,$, but then a subsequence has a dihedral tending to $0$). A subsequence converges to an $n$-simplex of diameter $1\,$, thickness $\geq t_0\,$, with a dihedral angle equal to $0\,$, a contradiction.
\end{proof}
\begin{lem}\label{diedre2} Given ${\bf a}\in ]0,\pi[$ and $\alpha_1,\alpha_2$ in $[0,\pi]$, with
$\alpha_1$ or $\alpha_2$ in $[{\bf a},\pi-{\bf a}]$, set $\alpha=\alpha_2-\alpha_1$. If $\alpha$ belongs to $[-\pi/2,\pi/2]$, one has the inequalities
\begin{equation}\label{blue1984}\frac{\vert\alpha\vert}{\pi}\leq \frac{\vert\sin\alpha\vert}{2}\leq \frac{\vert\cos\alpha_2-\cos\alpha_1\vert}{\sin\alpha_1+\sin\alpha_2}\leq \frac{\vert\cos\alpha_2-\cos\alpha_1\vert}{\sin{\bf a}}\  .
\end{equation}
\end{lem}
\begin{proof} The inequality on the left in (\ref{blue1984}) comes from the concavity of $\sin$ on $[0,\pi/2]$, while the one on the right comes from $\sin\alpha_1+\sin\alpha_2\geq \sin{\bf a}$. As for the one in the middle of (\ref{blue1984}), one writes the inequality 
\begin{equation*}\frac{\vert\sin\alpha\vert}{2}\leq\frac{\vert\sin\alpha\vert}{1+\cos\alpha}\ ,
\end{equation*}
then derives the expected inequality (achieved if $\alpha=0$) using 
\begin{equation}\label{diedre9}\frac{-\sin^2\alpha}{1+\cos\alpha}=\cos\alpha-1=
(\cos\alpha_2-\cos\alpha_1)(\frac{\sin\alpha}{\sin\alpha_1+\sin\alpha_2})\  ,
\end{equation}
which itself follows from
\begin{multline*} \cos\alpha-1=\cos\alpha_1(\cos\alpha_2-\cos\alpha_1)+
\sin\alpha_1(\frac{\sin^2\alpha_2-\sin^2\alpha_1}{\sin\alpha_2+\sin\alpha_1})=\cr=(\cos\alpha_2-\cos\alpha_1)(\cos\alpha_1-
\sin\alpha_1(\frac{\cos\alpha_2+\cos\alpha_1}{\sin\alpha_2+\sin\alpha_1}))=
\cr=(\cos\alpha_2-\cos\alpha_1)
(\frac{\sin\alpha}{\sin\alpha_2+\sin\alpha_1})
\  .\qedhere
\end{multline*}
\end{proof} 
\begin{lem} \label{diedre3} Given ${\bf a}\!\in ]0,\pi[$, suppose $g_1,g_2$ are Euclidean metrics on a vector space such that there exists ${\mathscr C}_1\!\in\! [0,1/2]$ verifying
$$\vert\,\Vert\,\cdot\,\Vert^2_{g_1}\!-\!\Vert\,\cdot\,\Vert^2_{g_2}\,\vert\leq {\mathscr C}_1\Vert\,\cdot\,\Vert^2_{g_1}
\  .$$
Denoting by $\alpha_1$ (by $\alpha_2$) the $g_1$-value (the $g_2$-value) of an angle (in $[0,\pi]$) between two intersecting half-lines, one has
\begin{equation}\label{diedre11}\vert\cos\alpha_2-\cos\alpha_1\vert\leq 24\,{\mathscr C}_1\  .
\end{equation}
If one assumes $\alpha_1$ or $\alpha_2$ is in $[{\bf a},\pi-{\bf a}]$ and ${\mathscr C}_1\leq {\sin{\bf a}}/{24}$, then $\alpha=\alpha_2-\alpha_1$ belongs to $[-\pi/2,\pi/2]$ and  lemma {\rm \ref{diedre2}} implies
\begin{equation}\label{diedre0}\vert\alpha\vert\leq \pi\,\frac{\vert\cos\alpha_2-\cos\alpha_1\vert}{\sin{\bf a}}\leq
\frac{24\,\pi}{\sin{\bf a}}\,{\mathscr C}_1\  .
\end{equation}
\end{lem}
\begin{proof} 
Choose $v$ and $w$ having $g_1$-norm equal to $1$ and supporting the two given half-lines in the same half-plane, one has
\begin{multline}\label{diedre5}\vert\cos\alpha_2-\cos\alpha_1\vert\leq\cr\leq
 \frac{\vert g_2(v,w)-g_1(v,w) \vert}{\Vert\,v\,\Vert_2\Vert\,w\,\Vert_2}\,+
\,\vert g_1(v,w) \vert\,\frac{\vert 1-\Vert\,v\,\Vert_2\Vert\,w\,\Vert_2\vert}{\Vert\,v\,\Vert_2\Vert\,w\,\Vert_2}\  .
\end{multline}
One has from the hypothesis
\begin{multline}\label{diedre01}\vert g_2(v,w)-g_1(v,w) \vert\leq\cr\leq \frac{1}{2}(\vert \Vert\,v+w\,\Vert_2^2-\Vert\,v+w\,\Vert_1^2\, \vert+
\vert \Vert\,v\,\Vert_2^2-\Vert\,v\,\Vert_1^2\, \vert+\vert \Vert\,w\,\Vert_2^2-\Vert\,w\,\Vert_1^2\, \vert)\leq\cr
\leq \frac{{\mathscr C}_1}{2}(\Vert\,v+w\,\Vert_1^2+\Vert\,v\,\Vert_1^2+\Vert\,w\,\Vert_1^2)\leq 3{\mathscr C}_1\  .
\end{multline}
One also has 
\begin{equation} \label{diedre7}\vert \Vert\,v\,\Vert_2-\Vert\,v\,\Vert_1\, \vert=\vert \Vert\,v\,\Vert_2-1\, \vert=\frac{\vert \Vert\,v\,\Vert_2^2-1\, \vert}{\vert \Vert\,v\,\Vert_2+1\, \vert}\leq {\mathscr C}_1\  .
\end{equation}
By (\ref{diedre01}), (\ref{diedre7}), the first term on the right-hand side in (\ref{diedre5}) verifies 
\begin{equation}\label{diedre8bis}\frac{\vert g_2(v,w)-g_1(v,w) \vert}{\Vert\,v\,\Vert_2\Vert\,w\,\Vert_2}\leq
\frac{3{\mathscr C}_1}{(1-{\mathscr C}_1)^2}\  ,
\end{equation}
while the second satisfies
\begin{equation}\label{diedre8}\vert g_1(v,w) \vert\,\frac{\vert 1-\Vert\,v\,\Vert_2\Vert\,w\,\Vert_2\vert}{\Vert\,v\,\Vert_2\Vert\,w\,\Vert_2}
\leq \frac{{\mathscr C}_1(2+{\mathscr C}_1)}{(1-{\mathscr C}_1)^2}\  .
\end{equation}
{\em Indeed,} one has 
\begin{multline*}\vert g_1(v,w) \vert \  \vert 1-\Vert\,v\,\Vert_2\Vert\,w\,\Vert_2\vert\leq\\\leq 
\max((\vert 1-\Vert\,v\,\Vert_2\  \vert)(1+\Vert\,w\,\Vert_2),(1+\Vert\,v\,\Vert_2)(\vert \,1-\Vert\,w\,\Vert_2\,\vert))\ .
\end{multline*}
Putting together (\ref{diedre5}), (\ref{diedre8bis}), (\ref{diedre8}), one gets (\ref{diedre11}) (use $(1\!-\!{\mathscr C}_1)\!\geq\!1/2$).
\par
As $\alpha_1$ and $\alpha_2$ are in $[0,\pi]$, the difference $\alpha=\alpha_2-\alpha_1$ is in $[-\pi/2,\pi/2]$ {\em if and only} if $\cos\alpha\geq0$. Bring back (\ref{diedre9})
\begin{equation*}\hskip2cm\cos\alpha=1+\frac{(\cos\alpha_2-\cos\alpha_1)\,\sin(\alpha_2-\alpha_1)}{\sin\alpha_2+\sin\alpha_1}\  \,,
\end{equation*}
thus $\cos\alpha$ is $\geq0$ if
$\,\vert{(\cos\alpha_1-\cos\alpha_2)\,\sin(\alpha_2-\alpha_1)}/{\sin{\bf a}}\vert\leq1\  ,
$
thus {\em a fortiori} if the following inequality holds
\begin{equation} \label{diedre24} \vert\,\cos\alpha_1-\cos\alpha_2\,\vert\leq \sin{\bf a}\  ,
\end{equation}
and (\ref{diedre11}) implies that (\ref{diedre24}) holds if $24\,{\mathscr C}_1\leq\sin{\bf a}$, which was assumed in the statement. This also proves (\ref{diedre0}): actually, lemma \ref{diedre2} holds because $\vert\alpha\vert\leq \pi/2$ is true.
\end{proof}
\begin{lem} \label{diedre12} Given ${\bf a}\!\in ]0,\pi[$, suppose $g_1,g_2$ are Euclidean metrics on a vector space for which exists ${\mathscr C}_1\in [0,\sin{\bf a}/24]$ such that 
\begin{equation} \label{diedre-6} \vert\,\Vert\,\cdot\,\Vert^2_{g_1}-\Vert\,\cdot\,\Vert^2_{g_2}\,\vert\leq {\mathscr C}_1\Vert\,\cdot\,\Vert^2_{g_1}
\  .
\end{equation}
Let $\eta'$ and $\eta''$ be half-hyperplanes having common boundary $\xi\!=\!\eta'\cap\eta''$, an $(n\!-\!2)$-vector subspace. They define a prism of $g_1$-dihedral  (of $g_2$-dihedral) angle $\beta_1$  (angle $\beta_2$) in $]0,\pi[$. If $\beta_1$ or $\beta_2$ is in $[{\bf a},\pi\!-\!{\bf a}]$, then 
\begin{equation*}\vert\,\beta_1-\beta_2\,\vert\leq \frac{76\,\pi}{\sin{\bf a}}\,{\mathscr C}_1\ .
\end{equation*}
\end{lem}
\begin{proof} Assume $\beta_1$ is in $[{\bf a},\pi-{\bf a}]$. Choose $v_1$ and $w_1$ of $g_1$-norm equal to $1$, pointing $g_1$-orthogonally to $\xi$ into $\eta'$ and $\eta''$. Call $\alpha_1$ the dihedral angle $\beta_1=\angle_{g_1}(v_1,w_1)$ and set $\alpha_2=\angle_{g_2}(v_1,w_1)$ (the angle $\alpha_2$ is in general distinct from $\beta_2$). By previous lemma \ref{diedre3}, one gets from (\ref{diedre0})
\begin{equation}\label{diedre-1}\vert\,\alpha_2-\alpha_1\,\vert\leq \frac{24\,\pi}{\sin{\bf a}}\,{\mathscr C}_1\  .
\end{equation}
Introduce in the same way $v_2$ and $w_2$ of $g_2$-norm equal to $1$, pointing $g_2$-orthogonally to $\xi$ into $\eta'$ and $\eta''$. The dihedral angle $\beta_2=\angle_{g_2}(v_2,w_2)$, as the angle $\alpha_2=\angle_{g_2}(v_1,w_1)$, can be viewed as spherical distances between points on the $g_2$-unit sphere, thus one has
\begin{equation} \label{diedre-12} \vert\,\beta_2-\alpha_2\,\vert\leq \angle_{g_2}(v_2,v_1)+\angle_{g_2}(w_2,w_1)\  ,
\end{equation}
where $\vartheta'=\angle_{g_2}(v_2,v_1)$ and $\vartheta''=\angle_{g_2}(w_2,w_1)$ are absolute angles counted in $[0,\pi/2[$. {\it Indeed} $v_1,v_2$ (or $w_1,w_2$) point into the same half-hyperplane $\eta'$ (or $\eta''$).
We now manage a control on $\vartheta'$ (and $\vartheta''$).
\par If $\vartheta'=0$, this is the best situation possible. If not, $v_1$ and $v_2$ generate a plane $P'$ tangent to $\eta'$, thus there exists a vector $u\in(\xi\cap P')\setminus\{0\}$. One writes ($u,v_1$ and $v_2$ are in the same half-plane $P'\cap \eta'$)
\begin{multline}\label{diedre-2}\vartheta'=\angle_{g_2}(v_2,v_1)=\vert\,\angle_{g_2}(u,v_2)-\angle_{g_2}(u,v_1)\,\vert=\cr=
\vert\,\frac{\pi}{2}-\angle_{g_2}(u,v_1)\,\vert=\vert\,\angle_{g_1}(u,v_1)-\angle_{g_2}(u,v_1)\,\vert\  .
\end{multline}
Set now $\mu_1=\angle_{g_1}(u,v_1)$ and $\mu_2=\angle_{g_2}(u,v_1)$ and apply lemma \ref{diedre3} to the situation of the two half-lines supported by $u$ and $v_1\,$: this works because $\mu_1=\pi/2$ and $\mu_2$ is in $[0,\pi]$ (thus the lemma applies, even making ${\bf a}=\pi/2$ and $\sin{\bf a}=1$) and the hypotheses on $g_1$, $g_2$ and ${\mathscr C}_1$ are the same in both lemmas \ref{diedre3}, \ref{diedre12}. One gets from (\ref{diedre-2}) and (\ref{diedre0})
\begin{equation}\label{diedre-100} \vartheta'=\vert\,\mu_2-\mu_1\,\vert\leq 24\,\pi\,{\mathscr C}_1\  ,
\end{equation}
and the same for $\vartheta''$. Thanks to (\ref{diedre-1}), (\ref{diedre-12}) and (\ref{diedre-100}), the proof is complete if $\beta_1\in[{\bf a},\pi-{\bf a}]$, which was assumed in the beginning.
\par If $\beta_2\in[{\bf a},\pi-{\bf a}]$, inverting the roles of $g_1$ and $g_2$, the only thing which changes is the constant ${\mathscr C}_1$ in the inequality (\ref{diedre-6}), namely we only can guarantee (from (\ref{diedre-6}))
\begin{equation*} \vert\,\Vert\,\cdot\,\Vert_{g_1}^2-\Vert\,\cdot\,\Vert_{g_2}^2\,\vert\leq {\mathscr C}'_1\Vert\,\cdot\,\Vert_{g_2}^2\hskip3mm \hbox{where} \hskip3mm {\mathscr C}'_1=\frac{{\mathscr C}_1}{1-{\mathscr C}_1}\leq
\frac{24}{23}{\mathscr C}_1\  .
\end{equation*}
As all the rest works the same, a computation shows that the statement also holds in this case.
\end{proof}
\subsubsection{Bounding the defect angle created by a bone}
\label{desirable}$ $

\begin{assump}\label{assumpred} In an oriented open $U$, near and around an oriented bone $\xi\,$, the {\em defect angle} $\alpha_\xi\in{\mathbb R}$ defines a rotation on the oriented $\xi^\perp$ (this orientation fits with those of $U$ and $\xi$, lemma {\rm\ref{Regge1}}, remark {\rm\ref{renvparcours}}, definition {\rm\ref{Reggeoperbone}}, lemma {\rm\ref{Reggeoperbonebis}}). Given a small neatest square $G$ hitting $\xi$, with $\hbox{index}_{G,\xi}=1$ and $\hat\xi\cup {\rm rge}(G)\subset U$ (definition {\rm\ref{psquare4}}, assumption {\rm\ref{assumptionX}}), the orientations on ${\rm rge}(G)$ and $\xi^\perp$ are paired. 
\end{assump}
\begin{prop} \label{diedre4} Given a compact domain $D$ in a Riemannian $n$-manifold $(M,g)$ and $\rho\in]0,\sqrt{{\sin{\bf a}}/{24\,{\mathfrak C}_1}}\,]$ (where ${\bf a}\!\in]0,\pi[$ is as in lemma {\rm\ref{diedre1}}) small enough so that theorem {\rm \ref{theorem I}} applies, one has an embedded simplicial complex $T: K\rightarrow (M,g)$, with $D\subset T(K)$, and, paired with $\rho\,$ {\em as in the statement} of theorem {\rm \ref{theorem I}}, a piecewise flat metric $g_0$ on $K\,$. There exists ${\mathfrak C}_5>0$ depending on $D\subset(M,g)$ such that, for any 
bone $\xi\in K\,$, the Regge curvature of a neatest square $G$ satisfying assumption {\rm\ref{assumpred}} verifies (see lemma {\rm\ref{Reggeoperbone1}}) 
\begin{equation}\label{diedre4444}\Vert\,R_0^G\,\Vert_{0,0}=\Vert\,({\bf r}_0^{\xi})^{-1}-{\rm Id}\,\Vert_{0,0}\leq {\mathscr C}:={\mathfrak C}_5 \, \rho^2\  {\rm and}\ \ \vert \alpha_\xi\vert \leq {\mathfrak C}_5 \  \rho^2\ ,
\end{equation}
where $-\alpha_\xi\,$ is the angle of the rotation 
$R_0^G+\hbox{\rm Id}$ around $\xi$ and ${\mathfrak C}_5$ is
${\mathfrak C}_5\!:=\!6^n\,76\,\pi\,b_n\,{\mathfrak C}_1/({\mathfrak C}_4\sin{\bf a})
$ (for ${\mathfrak C}_1,{\mathfrak C}_4\!>\!0$, see theorem {\rm\ref{theorem I}}).
\par If $G$ is only temperate, one has (see assumptions {\rm\ref{bounded22}}, {\rm\ref{hyp555}})
\begin{equation*}\Vert\,R_0^G\,\Vert_{0,0}\leq {\mathscr C}':=2\,\mu\,{\mathfrak C}_5 \, \rho^2\  .
\end{equation*}
\end{prop}
\begin{proof} Recall that the number $\nu$ of $n$-simplices $\hat\sigma$ containing a point $p$ is bounded:
actually, in proposition \ref{bds4} (to apply it, use remark \ref{remarquable1} linked to theorem \ref{theorem I}) one establishes the bound $\nu\!\leq\!  6^n\,b_n/{\mathfrak C}_4$ (with ${\mathfrak C}_4>0$), where $b_n$ is the volume of the unit $n$-dimensional ball. 
\par
Consider an $n$-simplex $\sigma\in K$ and a bone $\xi\in\partial\sigma\,$. Let $\beta^{g_0}$ be the $g_0$-dihedral angle defined by $\xi$ in $\sigma\,$ and $\eta',\eta''$ be the $(n-1)$-faces bounding this angle (thus $\eta'\cap\eta''=\xi$). Let $p$ be a point interior to $\xi\,$. 
Pushing forward $g_0$ on $\hat\sigma=T(\sigma)$, at $y=T(p)\in\hat\xi=T(\xi)$, one reads the same dihedral angle $\beta^{g_0}$ in $(T_yM,(T^{-1})^\ast {g_0}_p)$ between $T_y\hat\eta'$ and $T_y\hat\eta''$ along $T_y\hat\xi$. From the machinery helping to establish theorem \ref{theorem I}, actually by proposition \ref{2o} we know that the openness of $\hat\sigma$ verifies
$\omega_g(\hat\sigma)\geq \bar A/2^n>0\,$ (for $\bar A$, see (\ref{2c}) in definition \ref{2j}). 
\par As $T$ is (by \ref{theorem I} $(ii),\,(b)$) a ${\mathfrak C}_1\rho^2$-isometry from $(\sigma,g_0)$ to $(\hat\sigma,g)$, we also know that for $\rho>0$ small enough
$\omega_{g_0}(\sigma)\geq \bar A/2^{n+1}>0\,$, and this implies by lemma \ref{1c} that $t_{g_0}(\sigma)\geq$ to some positive uniform constant.
Thus $\beta^{g_0}$ is in $[{\bf a},\pi-{\bf a}]$ for some uniform ${\bf a}>0$ by lemma \ref{diedre1}. 
Let  $\beta^{g_y}$ be the $g_y$-dihedral angle between $T_y\hat\eta'$ and $T_y\hat\eta''$ along $T_y\hat\xi$ in $(T_yM,g_y)$. 
Theorem \ref{theorem I} $(ii),(b)$ reads for any $v\in T_yM$ 
\begin{equation*}\vert\,g_y(v,v)-(T^{-1})^\ast g_0(v,v)\,\vert\leq {\mathfrak C}_1\,\rho^2\,g_y(v,v)\  .
\end{equation*}
Due to the hypothesis made on the choice of $\rho\,$, lemma \ref{diedre12} applies with ${\mathscr C}_1={\mathfrak C}_1\,\rho^2\leq \frac{\sin{\bf a}}{24}\,$, giving
\begin{equation} \label{diedre-3} \vert\, \beta^{g_0}-\beta^{g_y}\,\vert\leq \frac{76\pi}{\sin{\bf a}}\,{\mathscr C}_1=\frac{76\,\pi\,{\mathfrak C}_1}{\sin{\bf a}}\,\rho^2\  .
\end{equation}
\par Take $\xi$ as in the claim ($\xi$ sits in the interior of $K\,$, definition \ref{essentialbone}).
Bringing in lemmas \ref{Regge1}, \ref{Reggeoperbonebis}, we see that $R_0^G+\hbox{\rm Id}$ is $\hbox{\rm Id}$ along all directions parallel to $\xi$ and is a $g_0$-rotation in the plane $\xi^\perp$ by an angle $-\alpha_\xi$ verifying (the simplicial dihedral angles below $\beta_i^{g_0},\beta_i^{g_y}\,$ are in $[0,\pi]$) 
\begin{equation}\label{diedre-1003}\vert\alpha_\xi\vert = \vert 2\,\pi-\sum_{i=0}^{\nu-1} \beta_i^{g_0}\vert\leq \sum_{i=0}^{\nu-1} \vert\beta_i^{g_y}-\beta_i^{g_0}\vert\  .
\end{equation}
Putting (\ref{diedre-3}) and  (\ref{diedre-1003}) together with the bound on $\nu\,$, one gets
\begin{equation*}\vert\,\alpha_\xi\,\vert\leq\sum_{i=0}^{\nu-1} \vert\,\beta_i^{g_y}-\beta_i^{g_0}\,\vert\leq
\frac{6^n\,76\,\pi\,b_n\,{\mathfrak C}_1}{{\mathfrak C}_4\sin{\bf a}}\,\rho^2\  .
\end{equation*}
As the $g_0$-norm of $R_0^G$ verifies 
$\Vert\,R_0^G\,\Vert_{0,0}=\sqrt{2\,(1-\cos\alpha_\xi)}\leq \vert\,\alpha_\xi\,\vert\  ,$ (see appendix \ref{diedre-4bis}),
(\ref{diedre4444}) follows.
\par As for the last claim, use (\ref{diedre4444}), remark \ref{plouc} (thus lemma \ref{generalisation1}).
\end{proof}
\subsubsection{The Regge curving}\label{convenable} 
$ $

\begin{assump} \label{assumpredbis} Assumption {\rm\ref{assumpred}} holds as in section {\rm\ref{desirable}}. From now on, our concern is, applying theorem {\rm\ref{theorem I}}, to study a polyhedral approximation to $(M,g)$ as $\rho\rightarrow0\,$. Doing this, proposition {\rm\ref{diedre4}} shows that the defect angle in ${\mathbb R}$ verifies
$\vert\alpha_\xi\vert=O(\rho^2)$.
Thus, {\em we shall now consider $\alpha_\xi$ to be an element of} $]-\pi,\pi[\,$. 
\end{assump}
Take $p=T^{-1}x$ in $T^{-1}(U)\cap{\mathcal K}\,$. The Regge curvature $({\bf r}_0^\xi)^{-1}-{\rm Id}$ of a bone $\xi\in \dagger\!K_{n\!-\!2}$ (see lemma \ref{Reggeoperbone1}) is an exact measure of the non-flatness of $g_0$ around and near the singularity (located along $\xi\,$) of $g_0\,$. 
\par Given $x\in(M,g)$ and any
$X,Y\in T_xM$ with $\Vert X\wedge Y\Vert=1\,$ defining a given oriented plane $P\subset T_xM\,$, the curvature endomorphism 
$$R(P):Z\in T_xM\mapsto R(P)Z:=R(X,Y)Z\in T_xM
$$ is {\em $g_x$-skew-symmetric}. 
The {\em Regge curving} below is a {\em $g_0$-skew-symmetric} endomorphism of $T_pM$ which reflects enough $({\bf r}_0^\xi)^{-1}-{\rm Id}\,$, so that, in this polyhedral approximation, their difference is an $O(\rho^2)\,$.
\begin{defn} \label{Regge-200} Let $\Omega_0$ be the $g_0$-volume element in $T^{-1}(U)\cap{\mathcal K}\,$ and $\ast_0$ the $g_0$-Hodge homomorphism sending any $r$-vector $a=a_1\wedge\cdots\wedge a_r\,$, with $a_1,\dots,a_r\in T_p{\mathcal K}$, to the $(n-r)$-vector $\ast_0(a)$ uniquely defined by
\begin{equation*}\forall b\in \wedge^{n-r}T_p{\mathcal K}\  \    \   \   \  \  \  \langle \ast_0 (a), b\rangle=\Omega_0(a\wedge b)\, .
\end{equation*}
\end{defn}
\begin{defn} \label{Regge22} The {\em Regge curving} $\underline{R}_0^G$ of a {\em neatest parametrised square} $G$ with ${\rm rge}(G)\subset T(\hbox{st}(\xi))$ ($G$ plays the role of the above plane $P$) is at $p=T^{-1}(x)=T^{-1}(G(0,0))$ the linear $g_0$-skew-symmetric mapping from $T_p{\mathcal K}$ into itself defined by setting 
\begin{equation} \label{Regge10000}\underline{R}_0^G:=-\alpha_\xi \ \hbox{\rm index}_{G,\xi}\ \ast_0(\xi\wedge \cdot)\  ,
\end{equation} 
where $\xi\in\wedge^{n-2}(T_p{\mathcal K})$ denotes the $(n\!-\!2)$-vector of norm $1$ supporting $T_p\xi\!\subset \!T_p{\mathcal K}$ and orienting the bone $\xi\,$ (remark {\rm\ref{parallstar}}, lemmas {\rm\ref{Reggeoperbonebis}}, {\rm\ref{Reggeoperbone1}}).
\end{defn}
\begin{rem} \label{banalbane} A look at lemma {\rm\ref{Regge1}}, remark {\rm\ref{renvparcours}}, definition {\rm\ref{Reggeoperbone}} and lemma {\rm\ref{Reggeoperbonebis}}, is helpful. Observe ({\rm\ref{Regge10000}}) is unchanged by a change of orientation of $U$ or $\xi\,$: recall $\alpha_{-\xi}=\alpha_\xi$ from lemma {\rm\ref{Reggeoperbonebis}}, watch $\hbox{\rm index}_{G,\xi}$ is changed into $-\hbox{\rm index}_{G,\xi}\,$. But ({\rm\ref{Regge10000}}) is changed into its opposite by a change of orientation on $G$, as it must be, since $G$ plays the role of $P\,$. 
\end{rem}
\begin{lem} \label{banalbanane}
If $G$ is a {\em temperate} square meeting $\xi$, there exists a neatest square $G'$ with ${\rm rge}(G')\!\subset\! T(\hbox{st}(\xi)), R_0^{G'}\!=\!s_{G,\xi}\,{\mathcal P}^0_A\!\circ {\bf r}_0^{\xi,G}\!\circ\!(\!{\mathcal P}^0_A)^{\!-\!1}$ ($A$ is ``as'' in definition {\rm\ref{cruelmanque}}, see also below definition {\rm\ref{ReggeRegge!!}}), verifying $s_{G',\xi}=s_{G,\xi}$, thus definition {\rm\ref{Regge22}} allows to write
\begin{equation*}\underline{\bf r}_0^{\xi,G}\!:=\!s_{G,\xi}\,(\!{\mathcal P}^0_A)^{\!-\!1}\!\circ\underline R_0^{G'}\circ{\mathcal P}^0_A\!=\!-\alpha_{\xi}\,(\!{\mathcal P}^0_A)^{\!-\!1}\!\circ \ast_0 (\xi\wedge\cdot)\!\circ\!{\mathcal P}^0_A\, .
\end{equation*}
\end{lem}
\begin{proof} Go back to definitions \ref{eurk}, \ref{cruelmanque} and lemma \ref{banalbof}.
\end{proof}
\begin{lem} \label{Regge-1} Suppose one has $\vert\,\alpha_\xi\,\vert\leq 1$. If $G$  with ${\rm rge}(G)\!\subset\! T(\hbox{st}(\xi))$ is a neatest square
hitting the bone $\xi\,$, then $R_0^G$ and $\underline{R}_0^G$ coincide on $\xi$ and one has
\begin{equation} \label{toto1} R_0^G = \underline{R}_0^G +O(\alpha_\xi^2) \  , \ \ \hbox{in fact one has}\  
\ \ \Vert\, R_0^G -\underline{R}_0^G\,\Vert_{0,0} \leq \alpha_\xi^2 \  .
\end{equation}
\end{lem} 
\begin{proof} The case $\hbox{\rm index}_{G,\xi}=0$ is clear. Choose on $\xi$ the orientation fitting with $\hbox{\rm index}_{G,\xi}=1$. In a well chosen basis of $T_p{\mathcal K}$, i. e. a $g_0$-orthonormal basis built up by completing a $g_0$-orthonormal basis of $T_p\xi$, restricting the computations in the plane $T_p\xi^\perp$ that matters and taking advantage of $\alpha_\xi\ \hbox{\rm index}_{G,\xi}=\alpha_\xi$, one gets ($\theta_1$ and $\theta_2$ belong to the interval determined by $0$ and $\alpha_\xi$)
\begin{multline*}R_0^G +\alpha_\xi\ \hbox{\rm index}_{G,\xi}\ \ast_0(\xi\wedge \cdot) +\frac{\alpha_\xi^2}{2}\,\hbox{\rm Id}=\\
\begin{pmatrix} \cos \alpha_\xi -1+\frac{\alpha_\xi^2}{2}&+\sin\alpha_\xi-\alpha_\xi\\
-\sin\alpha_\xi+\alpha_\xi&\cos \alpha_\xi -1+\frac{\alpha_\xi^2}{2}\end{pmatrix}=\frac{\alpha_\xi^3}{3!}
\begin{pmatrix} \sin \theta_1 &-\cos\theta_2\\
\cos\theta_2&\sin \theta_1\end{pmatrix}\  ,
\end{multline*}
and thus (write $\sin^2\theta_1+\cos^2\theta_2\leq \theta_1^2+1\leq \alpha_\xi^2+1$)
\begin{equation}\label{toto}
\Vert\, R_0^G +\alpha_\xi\ \hbox{\rm index}_{G,\xi}\ \ast_0(\xi\wedge \cdot) +\frac{\alpha_\xi^2}{2}\,\hbox{\rm Id}\,\Vert_{0,0} \leq \frac{\vert\,\alpha_\xi^3\,\vert}{3!}\sqrt{1+\alpha_\xi^2}\  ,
\end{equation}
giving
\begin{equation*}
R_0^G = \underline{R}_0^G -\frac{\alpha_\xi^2}{2}\,\hbox{\rm Id}+O(\alpha_\xi^3)\  .
\end{equation*}
The rest of (\ref{toto1}) also follows from (\ref{toto}), using $\vert\,\alpha_\xi\,\vert\leq 1$.
\end{proof}
\begin{defn} \label{ReggeRegge!!} In view of the previous lines, define the {\em Regge curving} of a temperate $G$ (see lemmas {\rm\ref{banalbof}},  {\rm \ref{singular1}}, {\rm\ref{banalbanane}}) to be
\begin{equation}\label{ReggeRegge!}
\underline{R}_0^G:=({\mathcal P}_A^0)^{-1}\circ(\sum_{\xi\in\dagger\! K_{n-2}} \!\!\!\!-\alpha_\xi\ \hbox{\rm index}_{G,\xi}\,\ast_0\!(\xi\wedge \cdot))\circ{\mathcal P}_A^0\  ,
\end{equation}
where ${\rm rge}(G)\!\cap\!\hat K_{n\!-\!2}\!\subset\!\hat\xi_0\,$,$\  A\!\subset\!{\rm rge}(G)\!\setminus \!\hat  \xi_0$ is an embedded path from $G(0,0)$ to $x'\!=\!G(\varrho,\tau)\!\in\! T(\hbox{\rm st}(\xi_0)\!\setminus\! \xi_0)$ (see definition {\rm\ref{cruelmanque}}). Thanks to lemma {\rm\ref{banalbanane}}, one has
$\ \underline{R}_0^G=\sum_{\xi\in \dagger\!K_{n-2}} \hbox{\rm index}_{G,\xi}\ \underline{\bf r}_0^{\xi,G}\ .$
\par
{\em If} $G$ is as in lemma {\rm\ref{generalisation1}}, {\em extend} ({\rm \ref{ReggeRegge!}}) {\em setting} $\,
\underline{R}_0^{G}\!\!:=\!\!\sum_{i,j\!=\!0}^{N'\!-\!1}\underline{R}_0({\backslash \!\!\!\Delta}_{i,j})$, 
where ${\mathcal P}_{A_{i,j}}^0\circ\!\underline{R}_0({\backslash \!\!\!\Delta}_{i,j})\!\circ\,({\mathcal P}_{A_{i,j}}^0)^{-1}\!=\!
\underline{R}_0^{\backslash \!\!\!G_{i,j}}\!=\!\sum_{\xi\in  \dagger\!K_{n-2}} \hbox{\rm index}_{\backslash \!\!\!G_{i,j},\xi}\ \underline{\bf r}_0^{\xi,\backslash \!\!\!G_{i,j}}\,$.
\end{defn}
\begin{lem}  \label{Regge-2033}Assume the hypotheses of proposition {\rm \ref{diedre4}}, that $\rho$ is such that $\vert\,\alpha_\xi\,\vert\leq 1$ and the hypotheses of lemma {\rm\ref{generalisation1}} hold for $G\,$. Then
\begin{equation*}  
\Vert \!\sum_{i,j=0}^{N'-1}\!\!R_0({\backslash \!\!\!\Delta_{i,j}})-\underline R_0^G\Vert_{0,0}\!=\!\Vert\sum_{i,j=0}^{N'-1}\!(R_0({\backslash \!\!\!\Delta_{i,j}}) -\underline{R}_0({\backslash \!\!\!\Delta_{i,j}}))\Vert_{0,0}\! \leq \!{\mathcal N}\,(\mu+1)\,{\mathfrak C}_5^2\ \rho^4\  .
\end{equation*}
\end{lem}
\begin{proof} Write $\Vert\cdot\Vert$ short for $\Vert\cdot\Vert_{0,0}=\Vert\cdot\Vert_{g_0,g_0}\,$. Proposition {\rm \ref{diedre4}} says $\Vert R_0^\xi\Vert\!\leq\! {\mathscr C}\!=\!{\mathfrak C}_5\rho^2$ for any $\xi\!\in\!\dagger\! K_{n\!-\!2}$. {\em So} (\ref{mygush}), (\ref{mygush33}), (\ref{fantastic}), (\ref{toto1}), lemmas {\rm\ref{generalisation1}}, \ref{banalbanane}, \ref{Regge-1} and definition \ref{ReggeRegge!!} give (the penultimate inequality is derived using (\ref{mygush33}), the definition of $\backslash\!\!\!{\mathscr C}_{\backslash\!\!\!G_{i,j}}$ in lemma {\rm\ref{singular1}} and also (\ref{toto1}))
\begin{multline*} \label{toto2} \Vert\sum_{i,j=0}^{N'-1}(R_0(\backslash \!\!\!\Delta_{i,j}) -\underline{R}_0(\backslash \!\!\!\Delta_{i,j}))\Vert \leq\sum_{i,j=0}^{N'-1}\Vert R_0^{\backslash \!\!\!G_{i,j}} -\underline{R}_0^{\backslash \!\!\!G_{i,j}}\Vert\leq\\ \leq \!\!
\sum_{i,j=0}^{N'-1}\!(\Vert R_0^{\backslash \!\!\!G_{i,j}} \!\!-\!\!\!\!\sum_{\xi\in \dagger\!K_{n-2}}\!\!\!\!\hbox{\rm index}_{\,\backslash \!\!\!G_{i,j},\xi}\,{\bf r}_0^{\xi,\backslash \!\!\!G_{i,j}}\Vert  +\!\!\!\sum_{\xi\in \dagger\!K_{n-2}}\!\!\!\!\vert \hbox{\rm index}_{\,\backslash \!\!\!G_{i,j},\xi}\vert\,\Vert {\bf r}_0^{\xi,\backslash \!\!\!G_{i,j}} \!-\underline{\bf r}_0^{\xi,\backslash \!\!\!G_{i,j}}\!\Vert)\!\leq
\\\leq (\!\sum_{\xi\in \dagger\!K_{n-2}}\sum_{i,j=0}^{N'-1}[\vert \hbox{\rm index}_{\,\backslash \!\!\!G_{i,j},\xi}\vert^2\!\!
+\!\!\vert \hbox{\rm index}_{\,\backslash \!\!\!G_{i,j},\xi}\vert])\ {\mathscr C}^2\leq  {\mathcal N}\,(\mu+1)\,{\mathfrak C}_5^2\ \rho^4\,.\qedhere
\end{multline*}
\end{proof}

\section{Convergence in holonomy}\label{convhol1}

\subsection{Preparatory bounds}$ $\label{preparatorybds}$ $

First a definition (repeated as definition \ref{carrementbis}). Below, $\Pi=\bar I^2$ is the set of definition of all parametrised squares $G\,$.
\begin{defn} \label{carrement} Denote by $I$ a bounded open interval containing $[0,1]\,$. View $\bar I^2$ as $\Pi\!=\!\{\!(se_1\!+\!te_2\!)\!\in\! {\mathbb R}^2\mid (s,t)\!\in\!\bar I^2\}\,$, where $e_1,e_2$ is the canonical basis of ${\mathbb R}^2\,$. 
Given $w\in{\bf S}^1\subset{\mathbb R}^2\,$, the $2$-dimensional square
$\Pi$ may be thought as a $1$-{\it parameter} family of segments $l$ of straight lines describing a set of segments $\bar l_w=\{l\in\bar l_w \mid l \ \hbox{parallel to} \ w\}$  so that $\coprod_{l\in \bar l_w}l\!=\! \Pi\,$. Denote by $\bar l_1$ and $\bar l_2$ the sets $\bar l_{e_1}$ and $\bar l_{e_2}$ in $\Pi\,$.
\end{defn} 
The next theorem is the central tool to handle generic intersections. It relies on theorem 7.1 of \cite{C-M} (see also theorem \ref{thm7.1C-Mbis}). It allows to apply, for $l\in \bar l_w$ and $w\in {\bf S}^1$, theorem \ref{theorem I} $(iii)$ to {\it all curved lines} ${\bf G}_l\!:=\!G(l)$ of {\it embedded squares} $G$, {\it jointly} managing intersections of squares ${\bf G}\!:=\!G(\bar I^2)$ and their lines ${\bf G}_l$ with approximating polyhedra.

Recall from definition \ref{texture33333} that ${\mathfrak S}_{t_0}(W)$ is the space of all Riemannian barycentric $n$-simplices $\hat\sigma$ buildable in $W$ having thickness $\geq\! t_0$.
Any $k$-face of an $n$-simplex that belongs to ${\mathfrak S}_{t_0}(W)$ and has small diameter (say $\leq2\rho$ for some given $\rho>0$) slips into (see theorem \ref{texture333}) some leaf of a local field fresh texture ${\mathcal L}\equiv X_i^k$ (see section \ref{multiple1088}, definition \ref{Texturedefn}) of total space $X_i^k\!=\!X^k_{x_i,\rho_{x_i},s_{0,k}}$, each $X_i^k$ having as basis a convex (definition \ref{D.convex}) ball $B(x_i,\rho_3)$, where $\rho_3\geq \rho_{x_i}\!>\!0$, the $B(x_i,\rho_{x_i}/2)$ building a finite covering of $\bar W$. The leaves $L$ are embedded $k$-submanifolds in the balls $B(x_i,\rho_3)$ and $p$ verifies $p\circ \iota_L\!=\!i_L\,$, where $i_L$ embeds $L$ in $W\,$.
\par Observe that the construction leading to theorem \ref{texture333} was done for fixed
$k\,$, so $\rho_3\!=\!\rho_{3,k},\,\rho_{x}\!=\!\rho_{x,k}$ depend on $k$ (as $s_0\!=\!s_{0,k}$ does). To care about this, choose $\rho_3$ and $\rho_{x}$ to be the minima of $\rho_{3,k},\rho_{x,k}$ for $k\!=\!n-2,n-1$ (this fits with lemma \ref{text8} and its proof) and {\it jointly} choose the finite covering $B(x_i,\rho_{x_i}/2)$, dealing with those values of $\rho_3$ and $\rho_{x}$ ($k\!=\!n-2$ and $k\!=\!n-1$ only matter to us). Recall from lemma {\rm \ref{simpletext}} and its proof (see also the beginning of section {\rm\ref{bone1111}}) the 
\begin{fact}\label{slip1033}All $\hat\sigma\!\in\!{\mathfrak S}_{t_0}(W)$ with diameter $\!\leq \!\inf_i \rho_{x_i}/4\,$ have $k$-faces which are part of a leaf of some $X_i^{k}$ and $\hat\sigma\!\subset\!B(x_i,3\rho_3/4)$.
\end{fact}

Given a compact domain $\bar D$ in a relatively compact open $W\!\subset\! M\,$, select a polyhedral approximation according to theorem \ref{theorem I} and definition  \ref{convpolyhedra}, thus {\em a sequence} of $T:K\!\rightarrow \!W\!\subset \!(M,g)$  with $\bar D\subset T(K)$, paired with any $\rho\!>\!0$ small enough, realised by a Riemannian barycentric triangulation $\hat K\!:=\!T(K)$ of thickness $\geq t_0\!>\!0\,$ and mesh $\rho$. 

\begin{defn}\label{testsquare1} A polyhedral approximation $T,K,\hat K\!:=\!T(K)$ of $D\subset(M,g)$ having thickness $\!\geq\!t_0\!>\!0$ is given, as in the description above.
 \par A parametrised square $G$ is {\em testable} for $\rho$ if, for all $\hat\sigma\in{\mathfrak S}_{t_0}(W)$ of diameter $\!\leq \!2\rho$, all $(n\!-\!1)$ and $(n\!-\!2)$ faces $\hat\eta$ and $\hat\xi$ of $\hat\sigma\,$, {\em counting multiplicities} ($m_{n-1}$ and $m_{n-2}$ are defined in remark {\rm\ref{codimension1066}}) one has
\vskip1mm
\par $(i)$ \textbullet\  for any $l\in\bar l_w,,w\in{\bf S}^1$, the curves ${\bf G}_{l}\!\in\!{\bf G}_{w}\!:=\!\{{\bf G}_{l}\mid l\!\in\!\bar l_w\}$ intersect $\hat\eta$ in at most $m_{n-1}(n^2+\!2\!+\!c)$ points, actually one has
$$\sum_{x\in \hat\eta\,\cap\,{\bf G}_{l}} m_x({\bf G}_{\l},\hat\eta)\leq m_{n-1}(n^2+\!2\!+\!c)\, ;
$$ 
\hskip1cm \textbullet\  {\em moreover} $(iii)$ of definition \ref{convpolyhedra} holds along any ${\bf G}_{l}\!\in\!{\bf G}_{w}\,$;  
\vskip1mm
\par $(ii)$ ${\bf G}\!:=\!G(\bar I^2)$ cuts $\hat\xi$ in at most $m_{n-2}(\!(n\!-\!1)n\!+\!c)$ points, in fact
$$\sum_{x\in \hat\xi\,\cap\,{\bf G}} m_x({\bf G},\hat\xi)\leq m_{n-2}((n\!-\!1)n+c)\, .
$$
\vskip1mm
\par\noindent $(iii)$ For any $\rho'\!\in]0,\rho]$ and for $K$ paired with $\rho'$, one has 
$${\bf G}\!\cap\! \hat K_{n-3}=\emptyset\ \ \hbox{and}\ \ (\partial{\bf G})\!\cap\! \hat K_{n-2}=\emptyset\ .
$$
\end{defn}
\begin{thm}\label{text-1133}{\rm (companion to Thom's theorem {\rm\ref{text-11}}).} 
Let ${\mathscr T}$ be a manifold, $\Theta\!\subset\!{\mathscr T}$ a compact domain, both of dimension $c$.
A polyhedral approximation $K,T$ of $D\subset(M,g)$ is given. In the open set of ${\mathscr T}$-families of embedded squares exists a residual ${\mathcal V}_c\!:=\!\{\bar G\!=\! \{\bar G_\vartheta , \vartheta\!\in\!{\mathscr T}\}\}$ such that, for any $\bar G\!\in\!{\mathcal V}_c$, one can find $\rho_{\bar G}>0$ (also depending on $\Theta$) verifying that for any $\rho\in]0,\rho_{\bar G}]$ exists a dense open set $\Theta_\rho\subset\Theta$ of $\vartheta$ such that $G_\vartheta$ is testable for $\rho$.
\par {\em We say that} (among $\vartheta\!\in\!\!\Theta$) {\em almost any} $G_\vartheta$ in $\bar G\!\in\!\!{\mathcal V}_c$ is testable for $\rho$.
\end{thm}
\begin{proof} The theorem is proved as theorem \ref{text-1133bis} (in part \ref{annexe1033}).
\end{proof}
\begin{rem} \label{bounded222} Assumption {\rm\ref{bounded22}} is true for all
 squares verifying the hypotheses of theorem {\rm\ref{text-1133}} with $\mu=m_{n-2}(\!(n\!-\!1)n\!+\!c)$.
\end{rem}

\begin{assump}\label{embedded squares} All $G$ are defined on an open $U\!\subset \R^2$ containing $\bar I^2$, with $I$ open interval and $I\supset \![0,1]$, and they have embedded ranges in convex (definition \ref{D.convex}) balls of $W$. Set
$(S,g_S)\!:=\!(\hbox{\rm rge}(G),g_{\mid\hbox{\tiny\rm rge}(G)})\subset(M,g)$. If $\Vert \cdot\Vert_{g_S}$ is the operator norm $\Vert \cdot\Vert_{\hbox{\tiny\rm can},g_S}$ (linking $(\bar I^2,\hbox{\rm can})$ to $(S,g_S)$), {\it define} $\beta_G$ to be the positive Lipschitz constant such that
\begin{equation}\label{Lipschitzconsts} \beta_G:=\sup_{(s,t)\in \bar I^2}\Vert dG(s,t)\Vert_{g_S}\ .
\end{equation}
\end{assump}

\begin{defn} \label{undeterm} Let $G$ be a parametrised square embedded in $(M,g)$. Denote by $\Pi_\Delta$ any sub-square (actually, a sub-rectangle) included in $\Pi\equiv[0,1]^2$ and parallel to it. Such a $\Pi_\Delta$ is paired in a canonical way by $G$ with $\Delta:=G(\partial \Pi_\Delta)$, the image of its boundary in $S=\hbox{\rm rge}(G)$.
\end{defn}
\begin{prop} \label{sqRiem} There exists $c_1$ depending on $D\,\subset\!(M,g)$ such that, given $G\,$, $\Delta,\Pi_\Delta$ as in definition {\rm\ref{undeterm}}, for $r\!\in]0,1/2]$ small enough so that $G^{-1}(G(\Pi))_r\!\subset\!\bar I^2$ (here, $X_r$ is {\em the $r$-neighborhood} of $X\!\subset\!(M,g)$)
\begin{gather}
\label{caraf1}\hbox{\rm vol}_g\,((G(\Pi))_r)\!\leq\! c_1\,\beta_G^2\,r^{n-2}\, ;
\\\label{caraf2}
\hbox{\rm area}_g(G(\Pi_\Delta))\leq\hbox{\rm area}_g(G(\Pi))\leq \beta_G^2\  ;
\\
\label{caraf3}\hbox{\rm length}_g(\Delta)\leq 4\,\beta_G\ .
\end{gather}
\end{prop}
\begin{proof} First, choose a maximal net of $n_1$ points $x_i$ on $G(\Pi)=\Sigma$ such that, for any $i,j$ with $i\!\not=\!j\,$, one has $d_{(S,g_S)}(x_i,x_j)\!\geq\! r$  (with $S\!=\!\hbox{\rm rge}(G)$). Thus, for any point $x\!\in\!\Sigma$ exists $i$ such that $d_{(S,g_S)}(x,x_i)\!<\! r\,$. This implies $\Sigma_r\!\subset\!\cup_{i=1}^{n_1}B^g(x_i,2r)$. On the other hand, if $p_i\!\in\!\Pi$ is such that $G(p_i)\!=\!x_i\,$, the balls $B^{\rm can}(p_i,r/2\beta_G)\!\subset\!\R^2$ are disjointed since (by (\ref{Lipschitzconsts}))
$$\forall i,j \hskip1cm i\not=j\hskip1cm \beta_G\,d_{\rm can}(p_i,p_j)\geq d_{(S,g_S)}(x_i,x_j)\geq r\  .
$$
As $G^{-1}(\Sigma_r)\!\subset\![-1/2,3/2]^2$ for $r$ small enough, the $n_1$ points $p_i$ are less than $(16/\pi)\,(\beta_G/r)^2$. More, there exists a uniform constant $ \bar c_1$ on $D\subset(M,g)$ such that, for any ball of radius $2r$ having center in $D\,$, one has $\hbox{\rm Vol}_g B^g(x,2r)\leq \bar c_1\,r^n\,$. Setting $c_1\!:=\!(16/\pi)\,\bar c_1\,$, the above inequalities pile up
and give (\ref{caraf1}). The other results are easy to derive.
\end{proof}
\begin{rem}\label{sqRiem1033} If $\forall (s,t)\!\in \!\bar I^2$ one has $\Vert dG(s,t)\Vert_{g_S}\!\geq\!\alpha_G\!>\!0$, ({\rm\ref{caraf1}}) gives
$$\hbox{\rm Vol}_g (\Sigma_r)
\leq c_1\,\hbox{\rm area}_g(G(\Pi))\,({\beta_G}/{\alpha_G})^2\,r^{n-2}\ .
$$
From now on $\hbox{\rm rge}(G)$ denotes $G(\Pi)\!\subset\! W$, unless otherwise stated.
\end{rem}
\begin{lem} \label{Riesqnb} Given $\rho\!>\!0$ small enough, there exists ${\mathfrak C}_6$ depending only on $D\!\subset\!(M,g)$ such that, for any (small) square $G$ with $\hbox{\rm rge}(G)\!:=\!G(\Pi)\!\subset\! W$, the number $\nu$ of $n$-simplices $\hat\sigma\!\in\!{\mathfrak S}_{t_0}(W)$ meeting $\hbox{\rm rge}(G)$ {\em at least once} (the number of elements in $\{\sigma\mid \hat\sigma\cap\hbox{\rm rge}(G)\not=\emptyset\}$)
verifies
\begin{equation}\nu\leq {\mathfrak C}_6\,({\beta_G}/{\rho})^2\ .
\end{equation}
\end{lem}
\begin{proof} As any $\hat\sigma$ is of diameter $\leq 2\rho$ (theorem {\rm \ref{2n}}), if $\hat\sigma\cap\hbox{\rm rge}(G)\not=\emptyset$ one has $\hat\sigma\subset(\hbox{\rm rge}(G))_{2\rho}$. From (\ref{caraf1}), proposition \ref{sqRiem}, write
\begin{equation*} \hbox{\rm vol}_g\,((\hbox{\rm rge}(G))_{2\rho})\!\leq\! c_1\,\beta_G^2\,(2\rho)^{n-2}\  .
\end{equation*}
Theorem \ref{theorem I} and remark \ref{remarquable1} produce ${\mathfrak C}_4\!>\!0$ depending on $D\!\subset\!(M,g)$ such that
$\hbox{\rm vol}_g\,(\hat\sigma)\!\geq\! {\mathfrak C}_4\,\rho^n$. Setting ${\mathfrak C}_6\!:=\!2^{n-2}\, c_1{\mathfrak C}_4^{-1} $ gives the result. 
\end{proof}
\begin{nota} Define ${\mathscr C}$ to be ${\mathscr C}={\mathfrak C}_5\,\rho^2$ given in proposition {\rm \ref{diedre4}} and set, as in lemma {\rm\ref{generalisation1}}, ${\mathscr C}'= 2\,\mu\, {\mathscr C}\,$, where $\mu=m_{n\!-\!2}((n\!-\!1)n\!+\!c)$ (assumption {\rm\ref{bounded22}}, remark {\rm\ref{bounded222}}). Define also
\begin{equation}\label{red1}\backslash \!\!\!{\mathfrak C}_5:= \mu\, {\mathfrak C}_5\ \ \hbox{and}\ \ \backslash \!\!\!{\mathfrak C}_6:=\frac{1}{2}\,\binom{n+1}{2}\,\mu\,{\mathfrak C}_6\ .
\end{equation}
\end{nota}
\begin{nota} \label{betagbar}
Given $\Theta\!\subset\!{\mathscr T}$ compact and $\bar G\!\in\!{\mathcal V}_c$ as in theorem \ref{text-1133}, we make the abuse to denote by the same symbol $\bar G$ the former $\bar G=\{G_\vartheta\!\mid\! \vartheta\!\in\!{\mathscr T}\}$ and also its restriction $\{G_\vartheta\!\mid\! \vartheta\!\in\!\Theta\}$ to $\Theta$. We set
$$\beta_{\bar G}:=\sup_{\vartheta\in\Theta}\beta_{G_\vartheta}\ .
$$
\end{nota}
\begin{prop}  \label{Riesqnb2} Let $D$ be a compact domain in $(M,g)$. 
One can find ${\mathfrak C}_7,{\mathfrak C}_8,{\mathfrak C}_9\!>\!0$, depending only on $D\!\subset\!(M,g)$, such that, given any family $\bar G\in{\mathcal V}_c$ restricted to $\Theta$ compact domain in a $c$-dimensional manifold ${\mathscr T}$, for any $\rho\!\in]0,\rho_{\bar G}]\cap[0, \beta_{\bar G}]$, the polyhedral approximation $T\!:\! K\!\rightarrow \!(M,g)$ paired with $\rho$ and $g_0$ through theorem {\rm \ref{theorem I}}, induces, for {\em almost any} $\vartheta$ in $\Theta$ (theorem {\rm\ref{text-1133}} $(iii)$), a division into $N'^2$ temperate ``sub-squares'' $\backslash\!\!\! G_{i,j}$ indexed by $i,j \!=\!0,\dots, N'-1$, where $N'$ and ${\mathscr C}'$ verify
\begin{equation*} N'\rho\leq {\mathfrak C}_7\,\beta_{\bar G} \  ,\  \  N'\,{\mathscr C}'\leq {\mathfrak C}_8\,\beta_{\bar G}\,\rho\  \  \hbox{and}\  \  
{N'}^3\,{{\mathscr C}'}^2\leq {\mathfrak C}_9 \  \beta_{\bar G}^3\,\rho\  .
\end{equation*}
Actually, recalling from lemma {\rm\ref{generalisation1}} that, for given $G,T,K,(M,g)$, the total number of intersections of ${\rm rge}(G)$ with $\hat K_{n-2}$, {\em counting multiplicities}, is denoted by ${\mathcal N}\,$, one has $\ {\mathcal N}\leq \backslash \!\!\!{\mathfrak C}_6\,(\beta_{\bar G}/\rho)^2\ $ for any $G\in\bar G\,$.
\end{prop} 
\begin{proof} As $\bar G\in{\mathcal V}_c\,$, theorem \ref{text-1133} $(i), (ii), (iii)$ applies.
To each $n$-simplex, hit by some $G:=G_\vartheta\in\bar G$, correspond $\binom{n+1}{2}$ bones: {\it a priori}, all these bones can be hit by $G$. Searching a bound on the amount (including multiplicities) of intersected bones by a given $G\in\bar G$, one can divide $\binom{n+1}{2}$ by two, for a given bone is contained in at least two different $n$-simplices hit by $G$.
Thus, by lemma \ref{Riesqnb}, theorem \ref{text-1133} $(ii)$ (see remark \ref{bounded222}), we know the number (including multiplicities) ${\mathcal N}_\rho$ of bones hit by $G$ to be less than $\backslash \!\!\!{\mathfrak C}_6\,(\beta_{\bar G}/\rho)^2$. {\em This proves the last claim.} And $\backslash \!\!\!{\mathfrak C}_6\,(\beta_{\bar G}/\rho)^2$ is also an upper bound (for $G\in\bar G$) on ${\mathcal N}'_\rho\,$, the number of temperate points in ${\rm rge}(G)\cap\hat K_{n-2}$ (each {\it temperate} point may be a {\it neat} point). Now, choose the smallest integer $N'$ such that ${\mathcal N}'_\rho\leq N'^2$, so  
$N'\leq \sqrt{\backslash \!\!\!{\mathfrak C}_6}\,(\beta_{\bar G}/\rho)+1\,$ and ${\mathfrak C}_7:=
\sqrt{\backslash \!\!\!{\mathfrak C}_6}+1$ gives the bound on $N'\rho$ (requiring $\rho\leq\beta_{\bar G}$).
\par We know by theorem \ref{text-1133} $(iii)$ that, for almost any $G\!\in\!\bar G\,$, 
the image ${\rm rge}(G)$ is transverse to $\hat K_{n-3}$ (so ${\rm rge}(G)\cap\hat K_{n-3}=\emptyset$). By lemma
\ref{psquare5},
there exists a decomposition of $G$ into $N'^2$ temperate parametrised squares ($\backslash\!\!\! G_{i,j}$ is the restriction of $G$ to the sub-square $[s_{i,j},s_{i+1,j}]\times[t_{j},t_{j+1}]$, some of those $\backslash\!\!\! G_{i,j}$ may be trivial)  
$$\backslash\!\!\! G_{i,j}:[s_{i,j},s_{i+1,j}]\times[t_{j},t_{j+1}]\rightarrow M\hskip2mm\hbox{where}\hskip2mm i,j=0,\dots,N'-1 \  .
$$
Setting ${\mathfrak C}_8=2\,{\mathfrak C}_7\,\backslash\!\!\! {\mathfrak C}_5\,$, one has
$N'\ {\mathscr C}'\leq  {\mathfrak C}_7\,(\beta_{\bar G}/\rho)\,2\,\backslash\!\!\! {\mathfrak C}_5 \ \rho^2\leq  {\mathfrak C}_8\,\beta_{\bar G}\,\rho\, $.
Define ${\mathfrak C}_9=4\,{\mathfrak C}_7^3\,{\backslash\!\!\! {\mathfrak C}}_5^2\,$ to get the last inequality.
\end{proof}
\begin{cor} \label{goodboy} The above proposition {\rm \ref{Riesqnb2}} implies
\begin{equation*} N'\ \sup_{i,j}\ \Vert\,R_0^{\backslash\!\!\! G_{i,j}}\,\Vert_{0,0} \leq
N'\ {\mathscr C}'\leq  {\mathfrak C}_8\,\beta_{\bar G}\,\rho 
\rightarrow 0 \ \ \hbox{as}\ \  \rho \rightarrow 0\ .
\end{equation*} 
If $\beta_{\bar G}\,\rho\leq{\mathfrak C}_8^{-1}$, one has $N'\ {\mathscr C}'\in[0,1]$ and may apply proposition {\rm \ref{Regge6bis}}. 
\end{cor}
\begin{proof} Proposition {\rm \ref{diedre4}} ensures,
for $i,j=0,\dots,N'-1 $
\begin{equation*}\Vert\,R_0^{\backslash\!\!\! G_{i,j}}\,\Vert_{0,0} \leq{\mathscr C}':=2\, \backslash\!\!\!{\mathfrak C}_5 \  \rho^2\ .\qedhere
\end{equation*}
\end{proof}
\begin{rem} Given $p\in{\mathcal K}\,$, set $x=T(p)$. At $p$, we shall need to compare $\Vert\cdot\Vert_0$ (and $\Vert\cdot\Vert_{0,0}$) with $\Vert\cdot\Vert_x$ (with $\Vert\cdot\Vert_{x,x}$). Ask $\rho$ to satisfy $\rho \leq 1/\sqrt{2{\mathfrak C}_1}\,$, in order to apply
definition {\rm\ref{convpolyhedra}} $(ii) (b)$ and theorem {\rm \ref{theorem I}}.
\end{rem} 
\begin{defn}\label{Tconjug} If $L$ is in ${\rm End}(T_p{\mathcal K})$, define $[L]_T\in{\rm End}(T_xM)$
\begin{equation*}[L]_T:=dT(p)\circ L\circ dT^{-1}(x)\  .
\end{equation*}
\end{defn}
\begin{lem}\label{pff} Let $g_0$ and $T^\ast g_x$ verify as in definition {\rm \ref{convpolyhedra}} $(ii) (b)$ 
\begin{equation}\label{pfff}\vert T^\ast g_x-g_0\vert \leq {\mathfrak C}_1\rho^2\  T^\ast g_x\  .
\end{equation}
If $\rho\leq 1/\sqrt{2\,{\mathfrak C}_1}$, one gets
\begin{gather*}\frac{2}{3}\Vert\cdot\Vert_0\leq\Vert dT(T^{-1}(x))(\cdot)\Vert_x\leq 2\Vert\cdot\Vert_0\ ,\ \\ \  \   \frac{1}{3}\Vert\ast\Vert_{0,0}\leq\Vert [\,\ast\,]_T\Vert_{x,x}\leq 3\Vert\ast\Vert_{0,0}\  .
\end{gather*}
\end{lem}
\begin{proof} The claim is a direct consequence of  (\ref{pfff}) written as
\begin{equation*}(1-{\mathfrak C}_1\,\rho^2)\Vert dT(p)(\cdot)\Vert_x\leq \Vert \cdot\Vert_0\leq (1+{\mathfrak C}_1\,\rho^2)\Vert dT(p)(\cdot)\Vert_x\  .\qedhere
\end{equation*}
\end{proof}

\subsection{The main theorems}\label{main33}$ $
\begin{defn}\label{sqrestr} 
From now on, we shall consider a given square $G$ together with its restrictions $G_{\mid r}$, which are defined for any $r\in[0,1]$
$$G_{\mid r}:(s,t)\in[0,r]^2\mapsto G_{\mid r}(s,t)\in M\ .
$$
\end{defn}
\begin{rem} \label{help1} According to theorem {\rm\ref{text-1133}} $(iii)$, for {\em almost any} $G_\vartheta$ means for $G_\vartheta$ such that $\vartheta$ belongs to an open dense subset of ${\mathscr T}$.
\end{rem}
\begin{rem} \label{help2} Here, we exceptionally recall that section {\rm\ref{sssec:num3}} allows to define the parallel translation along $\Gamma$ if $G\!\in\!\bar G$ is such that $\Gamma\!=\!\partial{\rm rge}(G)$ meets $\hat K_{n-2}$.
\end{rem}
\begin{thm} \label{theorem F} Let $D$ be a compact domain in $(M,g)$ and, paired to $\rho\!>\!0$, bring in the polyhedron given by theorem {\rm\ref{theorem I}}, $T\!:\! K\!\rightarrow\! (M,g)\,$ (with $\bar D\!\subset \!T(K)$) and its metric $g_0$. Given any compact domain $\Theta\subset{\mathscr T}$, applying proposition {\rm \ref{Riesqnb2}} (thus lemma {\rm\ref{psquare5}}), divide almost any square $G_{\vartheta,{\mid r}}\!=\!G_{\mid r}\!:\![0,r]\times[0,r]\!\rightarrow \!D$, where $G_\vartheta=G$ is in $\bar G\in {\mathcal V}_c$ (theorem {\rm\ref{text-1133}}) and $\vartheta\in\Theta$, into $N_r'^2$ (with $N_r'\leq N'=N'_1$) temperate squares $\backslash\!\!\!G_{i,j}$ for $i,j\! =\!0,\dots, N_r'-1$. Set $x=G(0,0)$.
\par There exists ${\mathfrak C}_{10}>0$ such that, for $\rho$ small enough and any $r\in[0,1]$
\begin{multline*}\Vert\!\! \int_0^r\!\!\!\int_0^r\!\!{\mathcal P}_{B_{s,t}}^{-1}\!\circ\! R_G (s,t) \!\circ\! {\mathcal P}_{A_{s,t}} ds dt-[\sum_{i,j=0}^{N_r'-1} \!
(\!{\mathcal P}^0_{B_{i,j}}\!)^{\!-\!1}\!\!\circ\! R_0^{\backslash\!\!\!G_{i,j}}\!\circ\! {\mathcal P}^0_{A_{i,j}}]_T\Vert_{x,x} \!\leq\cr\leq {\mathfrak C}_{10}\,\beta_{\bar G}\max(1,\beta_{\bar G}^2)\,\rho\  ,
\end{multline*}
where $R_G(s,t)=R_{G(s,t)}(dG(\frac{\partial}{\partial s}),dG(\frac{\partial}{\partial t}))$ (notation {\rm\ref{notation007}}).
\end{thm}
\begin{proof} Recall that $[\sum_{i,j=0}^{N'-1} 
({\mathcal P}^0_{B_{i,j}})^{\!-\!1}\circ R_0^{\backslash\!\!\!G_{i,j}}\circ {\mathcal P}^0_{A_{i,j}}]_T$ stands for 
\begin{equation*}dT(T^{-1}(x))\circ (\sum_{i,j=0}^{N'-1} 
({\mathcal P}^0_{B_{i,j}})^{-1}\circ R_0^{\backslash\!\!\!G_{i,j}}\circ {\mathcal P}^0_{A_{i,j}}) \circ dT^{-1}(x)\,,
\end{equation*}
see definition \ref{Tconjug}.
Assume $\rho$ is in $]0,\rho_{\bar G}]\cap\,]0,1/\sqrt{2\,{\mathfrak C}_1}]$ ($\rho_{\bar G}$ refers to theorem {\rm\ref{text-1133}}, ${\mathfrak C}_1$ to theorem \ref{theorem I}).
We develop the case $r=1$.\par
The transversality of $\bar G$ to $\hat K_{n-3}\!=\!T(K_{n-3})$ makes us focus on the $G_\vartheta\!=\!G\!\in \!\bar G$ for which ${\rm rge}(G)\cap\hat K_{n-3}=\emptyset$ (this holds for almost any $G\in \bar G$, see theorem {\rm\ref{text-1133}}, definition \ref{testsquare1} $(iii)$ and remark \ref{help1}).\par
A family of parametrised squares $\bar G\in{\mathcal V}_c$ given through theorem {\rm\ref{text-1133}} can be viewed as the family $\bar\gamma:=\bar{\bf G}_1$ (or $\bar\delta:=\bar{\bf G}_2$) in $c+1$ parameters of curves $\gamma_t$ (or $\delta_s$); see also the ``squared squares'' in subsection \ref{adaptThom} and theorem \ref{text-1133bis}.  Theorem \ref{theorem I} and definition \ref{convpolyhedra} $(iii)$ imply, for those curves which avoid $\hat K_{n-2}=T(K_{n-2})$, that there exists ${\mathfrak C}_2:={\mathfrak C}_2(c+1)$ (depending on $D\subset(M,g)$) such that, for any curve $\gamma_t\in\bar\gamma$ or $\delta_s\in\bar\delta$ avoiding $\hat K_{n-2}$, one can write, bringing in the polyhedron $T,K,D\subset(M,g)$ associated with $\rho$ (${\mathcal P}_\gamma^0$ means ${\mathcal P}_{T^{-1}\gamma}^0$) 
\begin{gather*} \Vert{\mathcal P}_{\gamma_t}- dT \circ {\mathcal P}_{\gamma_t}^0\circ dT^{-1}\Vert_{\gamma_t(0),\gamma_t(1)}\leq {\mathfrak C}_2\,\rho\ \hbox{\rm length}_g(\gamma_t)\,;\\
\Vert{\mathcal P}_{\delta_s}- dT \circ {\mathcal P}_{\delta_s}^0\circ dT^{-1}\Vert_{\delta_s(0),\delta_s(1)}\leq {\mathfrak C}_2\,\rho\ \hbox{\rm length}_g(\delta_s)\,;
\end{gather*}
and thus, for $\Gamma = \Gamma_{0,0;1,1}\,$, if $\Gamma\cap\hat K_{n-2}=\emptyset$ (see (\ref{caraf3}) in proposition \ref{sqRiem}; $[{\mathcal P}_\Delta^0]_T$ is short for $dT \circ {\mathcal P}_\Delta^0 \circ dT^{-1}$ by definition \ref{Tconjug})
\begin{equation}\label{ahah}
\Vert{\mathcal P}_{\Gamma}- [{\mathcal P}_\Gamma^0]_T\Vert_{x,x}\leq {\mathfrak C}_2\,\rho\ \hbox{\rm length}_g(\Gamma)
\leq 4\,{\mathfrak C}_2\,\beta_{\bar G}\,\rho\  .
\end{equation}
If $\Gamma\cap\hat K_{n-2}\not=\emptyset$, the same inequality may be written, see remark \ref{help2}.
If $\beta_{\bar G}\,\rho\leq 1/{\mathfrak C}_8$, corollary \ref{goodboy},
propositions {\rm \ref{Riesqnb2}}, \ref{Regge6bis} give
\begin{equation*}\Vert({\mathcal P}^0_\Gamma)^{\!-\!1}\!-\!\hbox{\rm Id}\!-\!\sum_{i,j=0}^{N'-1}\!({\mathcal P}^0_{B_{i,j}})^{\!-\!1}\!\circ \!R_0^{\backslash\!\!\!G_{i,j}} \!\circ\! {\mathcal P}^0_{A_{i,j}}\,\Vert_{0,0}\leq N'^3{\mathscr C}'^2\leq {\mathfrak C}_9\,\beta_{\bar G}^3\, \rho\  .
\end{equation*}
As $\rho\leq 1/\sqrt{2\,{\mathfrak C}_1}$, it results from  lemma \ref{pff}
\begin{equation}\label{Regge330}\Vert\,[({\mathcal P}^0_\Gamma)^{-1}-\hbox{\rm Id}-\sum_{i,j=0}^{N'-1}({\mathcal P}^0_{B_{i,j}})^{-1}\circ R_0^{\backslash\!\!\!G_{i,j}} \circ {\mathcal P}^0_{A_{i,j}}]_T\,\Vert_{x,x}\leq 3\, {\mathfrak C}_9\,\beta_{\bar G}^3\, \rho\  .
\end{equation}
\par But one may bring back (\ref{formule8}) stated in theorem \ref{theorem A}, writing
\begin{equation}\label{buzz} ({\mathcal P}_{\Gamma})^{-1}-\hbox{\rm Id}=\int_0^1\!\!\int_0^1{\mathcal P}_{B_{s,t}}^{-1}\circ R_G(s,t) \circ {\mathcal P}_{A_{s,t}} \,ds\,dt\  .
\end{equation}
\par
Setting
${\mathfrak C}_{10}=4\,{\mathfrak C}_2+3 {\mathfrak C}_9$, (\ref{ahah}), (\ref{Regge330}), (\ref{buzz}) give the result.
\end{proof}
\par For coming use, we give a lemma generalising an observation already done if $\Delta=\Gamma$, while proving the last theorem \ref{theorem F}.
\begin{defn} \label{Bisounours1044} For $j=0,1,\dots,N'-1$ and $t\in[t_{j},t_{j+1}]$, define the loops $\Gamma_{i,j}^t$ and $\Gamma_{j}$ to be (see ({\rm\ref{paths0}}), ({\rm\ref{nota12}}), ({\rm\ref{paths2}}))
\begin{gather*} \Gamma_{i,j}^t=\Gamma_{s_{i,j},t_{j}; s_{i+1,j},t}=\delta_{s_{i,j}}^{t,t_{j}}\vee\gamma_{t}^{s_{i+1,j},s_{i,j}}\vee\delta_{s_{i+1,j}}^{t_{j},t}\vee\gamma_{t_{j}}^{s_{i,j},s_{i+1,j}}\  ,\\
 \Gamma_{t}=\Gamma_{0,0; 1,t}=\delta_0^{t,0}\vee\gamma_{t}^{1,0}\vee\delta_1^{0,t}\vee\gamma_0^{0,1}\  \  \hbox{and}\  \  
\Gamma_{j}=\Gamma_{0,0; 1,t_{j+1}}\  .
\end{gather*}
\end{defn}
\begin{defn} \label{Delta1044} Given $\sigma,\varrho,\tau,\upsilon\!\in\![0,1]$, while thinking to definition {\rm\ref{undeterm}}, {\em still denote by} $\Delta$ the loop $\gamma_0^{\sigma,0}\vee\delta_\sigma^{\tau,0}\vee\Gamma_{\sigma,\tau;\varrho,\upsilon}\vee\delta_\sigma^{0,\tau}\vee\gamma_0^{0,\sigma}\,$, extending ({\rm\ref{nota12}}), ({\rm\ref{blue7}}). Thus $\Delta$ now also denotes the boundary of the range by $G$ of a {\em sub-square $\Pi_\Delta$ parallel to $\Pi$}, as in definition {\rm\ref{undeterm}}, a loop based at $G(\sigma,\tau)$, {\em and} the above {\em extended loop} which is {\em based at} $x=G(0,0)$.
\end{defn}
\begin{lem}  \label{length} Given a family of squares $\bar G\in{\mathcal V}_c$ restricted to a compact domain $\Theta\!\subset\!{\mathscr T}$ (theorem {\rm\ref{text-1133}}), for any $\rho\!\in]0,\rho_{\bar G}]$, any $\Delta$ avoiding $\hat K_{n-2}$, with $\Delta\subset {\rm rge}(G)$ as in definition {\rm\ref{Delta1044}} ($\Delta$ is based at $x$), one has
\begin{equation} \label{parallprox} \Vert{\mathcal P}_{\Delta}-[{\mathcal P}_\Delta^0 ]_T\Vert_{x,x}\leq {\mathfrak C}_2\,\rho\ \hbox{\rm length}_g(\Gamma)\leq 4\,{\mathfrak C}_2\, \beta_{\bar G}\,\rho\  ,
\end{equation}
and the same inequality holds replacing together ${\mathcal P}_{\Delta}$ and ${\mathcal P}_\Delta^0$ by their inverses $({\mathcal P}_{\Delta})^{-1}$ and $({\mathcal P}_\Delta^0)^{-1}$.
\end{lem}
\begin{proof} Any $\Delta$ is a curve made of pieces of $\gamma_t$ and $\delta_s$ - some in the reverse direction - avoiding $\hat K_{n-2}$ which is the image by $G$ of a curve of length less than the length of $\partial \Pi\,$. 
Theorem \ref{theorem I} and definition \ref{convpolyhedra} $(iii)$ applies to each curve constituting $\Delta$, one has for any $s,\sigma,\varrho$ and $t,\tau,\upsilon$
\begin{gather*}  \Vert{\mathcal P}_{\gamma_t^{\sigma,\varrho}}- dT \circ {\mathcal P}_{\gamma_t^{\sigma,\varrho}}^0\circ dT^{-1}\Vert_{\gamma_t(\sigma),\gamma_t(\varrho)}\leq {\mathfrak C}_2\,\rho\ \hbox{\rm length}_g(\gamma_t^{\sigma,\varrho})\,;\\ 
 \Vert{\mathcal P}_{\delta_s^{\tau,\upsilon}}- dT \circ {\mathcal P}_{\delta_s^{\tau,\upsilon}}^0\circ dT^{-1}\Vert_{\delta_s(\tau),\delta_s(\upsilon)}\leq {\mathfrak C}_2\,\rho\ \hbox{\rm length}_g(\delta_s^{\tau,\upsilon})\,.
\end{gather*}
Using proposition \ref{sqRiem}, one gets (\ref{parallprox}) as expected.
Reversing the sense of $\Delta\,$, everything works the same. 
\end{proof}
Before stating and proving the next theorem in this section, it is worthwile to mention several things.
\begin{rem}
All $g$ or $g^+$-parallel translations considered in this section are elements in ${\bf G}_x={\bf GL}^+(T_xM)$ and their $\underline g_x$-distances $d_{\underline g_x}(\cdot,\ast)$ considered are their Riemannian distances in $({\bf G}_x,\underline g_x)$.
\end{rem}
\begin{nota}
Call ${\bf St}_2(T_xM)$ the Stiefel manifold of $g_x$-orthonormal $2$-frames and set 
\begin{equation}\label{BorneCourb}{\bf K}=\sup_{x\in D}\ \sup_{(u,v)\in{\bf St}_2(T_xM)} \Vert R(u,v) \Vert_{x,x}\  .
\end{equation}
\end{nota}
\begin{lem} \label{Pouf} Given $G$ with $x=G(0,0)$, one has for any $\Delta\subset{\rm rge}(G)$ as in definition {\rm\ref{Delta1044}} ($\Delta$ is based at $x$)
\begin{equation}\label{boundparallelsq}
\Vert {\mathcal P}_{\Delta}-\hbox{\rm Id}\Vert_{x,x}\leq {\bf K}\ \hbox{\rm area}_g(G(\Pi_\Delta))\leq {\bf K}\ \beta_G^2 \  ,
\end{equation}
implying $
\Vert {\mathcal P}_{\Delta}\Vert_{x,x}\!\leq \!
1+ {\bf K}\ \beta_G^2\,$.
The same holds with ${\mathcal P}_{\Delta}^{-1}$ replacing ${\mathcal P}_{\Delta}\,$.
\end{lem}
\begin{proof} Define $A_{\sigma,\tau}^{s,t}=\gamma_{t}^{\sigma,s}\vee\delta_\sigma^{\tau,t}$ and $B_{\sigma,\tau}^{s,t}=\delta_{s}^{\tau,t}\vee\gamma_{\tau}^{\sigma,s}$. Noticing that $\Vert {\mathcal P}_{\Delta}-\hbox{\rm Id}\Vert_{x,x}=\Vert {\mathcal P}_{\Gamma_{\sigma,\tau;\varrho,\upsilon}}-\hbox{\rm Id}\Vert_{y,y}$ with $y=G(\sigma,\tau)$ (if $A$ is a linear isometry $A:(E,\Vert\cdot\Vert_E)\rightarrow(F,\Vert\cdot\Vert_F)$ and $L:F\rightarrow F$ is linear, one has $\Vert L\Vert_{F,F}=\Vert L\circ A\Vert_{E,F}=\Vert A^{-1}\circ L\Vert_{F,E}\,$),
bringing back theorem \ref{theorem A}, using definition \ref{normcurv} as well as (\ref{formule10}), one can write
\begin{equation*} \Vert {\mathcal P}_{\Delta}-\hbox{\rm Id}\Vert_{x,x}=
\Vert\!\int\!\!\!\int_{[\tau,\upsilon]\times[\sigma,\varrho]} {\mathcal P}_{B_{\sigma,\tau}^{s,t}}^{-1} \circ R^G(s,t) \circ {\mathcal P}_{A_{\sigma,\tau}^{s,t}}\ \ G^\ast d\,{\rm area}\ \Vert_{y,y}\ .
\end{equation*}
Setting $z=G(s,t)$ (as ${\mathcal P}_{A_{\sigma,\tau}^{s,t}},{\mathcal P}_{B_{\sigma,\tau}^{s,t}}$ are linear isometries), one gets (\ref{boundparallelsq})
\begin{multline*} \sup_{{u \in T_xM}\atop{\Vert u\Vert_x=1}}\Vert({\mathcal P}_{\Delta}-\hbox{\rm Id})(u)\Vert_x\leq \!
\int\!\!\!\int_{[\tau,\upsilon]\times[\sigma,\varrho]}  \sup_{{u' \in T_zM}\atop{\Vert u'\Vert_z=1}} \Vert R^G(s,t) (u')\Vert_z\ \ G^\ast d\,{\rm area}
\leq \\\leq {\bf K}\ \hbox{\rm area}_g(G(\Pi_\Delta))\leq {\bf K}\ \beta_G^2\ ,
\end{multline*}
where the last inequality is from proposition \ref{sqRiem} (\ref{caraf2}).
Reversing $\Delta$, everything works the same for
${\mathcal P}_{\Delta}^{-1}$ in place of ${\mathcal P}_{\Delta}$.
\end{proof}

From lemmas \ref{length} and \ref{Pouf}, one deduces
\begin{lem}\label{club} Given a family $\bar G\in{\mathcal V}_c$ restricted to a compact $\Theta\!\subset\!{\mathscr T}$, for any $\rho\!\in]0,\rho_{\bar G}]$ such as $\beta_{\bar G}\rho\leq 1/{\mathfrak C}_8$ (corollary {\rm\ref{goodboy}}), any $\Delta$ avoiding $\hat K_{n-2}$, with $\Delta\subset {\rm rge}(G)$ as in definition {\rm\ref{Delta1044}} ($\Delta$ based at $x$), one has
\begin{equation} \label{parallprox2} \Vert [{\mathcal P}_\Delta^0 ]_T\Vert_{x,x}\leq 1+{\bf K}\ \beta_{\bar G}^2+ 4\,{\mathfrak C}_2\, \beta_{\bar G}\,\rho\leq1+{\bf K}\ \beta_{\bar G}^2+ 4\,{\mathfrak C}_2/{\mathfrak C}_8\  ,
\end{equation}
and the same inequality holds replacing ${\mathcal P}_\Delta^0$ by $({\mathcal P}_\Delta^0)^{-1}$.
\end{lem}

\begin{thm} \label{theorem G} Given a compact domain $D$ in $(M,g)$ and, paired to $\rho\!>\!0$, the polyhedron $T\!:\! K\!\rightarrow\! (M,g)$ produced by theorem {\rm\ref{theorem I}} ($\bar D\!\subset \!T(K)$) with its metric $g_0$, given a compact domain $\Theta\subset{\mathscr T}$ (theorem {\rm\ref{text-1133}}), apply proposition {\rm \ref{Riesqnb2}} to divide almost any square $G_{\vartheta,{\mid r}}\!=\!G_{\mid r}\!:\![0,r]\times[0,r]\!\rightarrow \!D$, where $G_\vartheta=G$ is in $\bar G\in {\mathcal V}_c$ and $\vartheta \in\Theta$, into $N_r'^2$ (with $N'_r\leq N'=N'_1$) temperate squares $\backslash\!\!\!G_{i,j}$ for $i,j\! =\!0,\dots, N_r'-1$. Set $x=G(0,0)$. 
There exists a polynomial ${\mathfrak C}_{11}(\beta_{\bar G}^2)$ of degree $3$ in $\beta_{\bar G}^2$ with nonnegative coefficients (independent of $\beta_{\bar G}$) such that, for $\rho$ small enough and any $r\!\in\![0,1]$, one has
\begin{multline*} \Vert\!\int_0^r\!\!\!\int_0^r\!{\mathcal P}_{A_{s,t}}^{-1}\circ R_G(s,t) \circ {\mathcal P}_{A_{s,t}} \,ds \,dt- [ (\sum_{i,j=0}^{N_r'-1} \!
({\mathcal P}^0_{A_{i,j}})^{-1}\circ R_0^{\backslash\!\!\!G_{i,j}}\circ {\mathcal P}^0_{A_{i,j}}) ]_T\Vert_{x,x}\!\leq \\\leq {\mathfrak C}_{11}(\beta_{\bar G}^2)\,\beta_{\bar G}\,\rho
\  .
\end{multline*}
\end{thm}
\begin{rem} In all proofs below, every $G\in\bar G$ is assumed to verify ${\rm rge}(G)\cap \hat K_{n-3}=\emptyset\,$ while the parallel translation along $\Gamma=\partial G$ may be extended according to section {\rm\ref{sssec:num3}} (see remarks {\rm\ref{help1}}, {\rm\ref{help2}}). 
\end{rem}
\begin{proof} For $\Gamma_{i,j}^t , \Gamma_{t} , \Gamma_{j}\,$, go back to definition \ref{Bisounours1044}. We assume $r=1\,$: the case  $r\in[0,1]$ is handled in the same way, quoting $N'_r\leq N'$.
\begin{assump}\label{clap}
We assume that all $\rho$ are in $]0,\rho_{\bar G}]\,\cap\,]0,1\!/\!\sqrt{2\,{\mathfrak C}_1}]$
(this points to theorems {\rm\ref{text-1133}} and \ref{theorem I}) and such as to verify $\beta_{\bar G}\rho\leq 1/{\mathfrak C}_8$
(this refers to proposition \ref{Riesqnb2} and corollary \ref{goodboy}).
\end{assump}
\begin{lem} \label{thG1} With the hypotheses of theorem {\rm \ref{theorem G}}, bring in the constant ${\mathfrak C}_9$ given in proposition {\rm \ref{Riesqnb2}}, one has 
\begin{equation*}\Vert [\sum_{j=0}^{N'-1}\! {\mathcal P}_{\Gamma_{j-1}}^0\circ ((\!{\mathcal P}_{\Gamma_j}^0\!)^{\!-\!1} \!-(\!{\mathcal P}_{\Gamma_{j-1}}^0)^{\!-\!1}\!) -\!\!\sum_{i,j=0}^{N'-1} \!
(\!{\mathcal P}^0_{A_{i,j}}\!)^{\!-1\!}\circ R_0^{\backslash\!\!\!G_{i,j}}\circ {\mathcal P}^0_{A_{i,j}}]_T
\Vert_{x,x}\!\leq \!3{\mathfrak C}_9\beta_{\bar G}^3\rho\,.
\end{equation*}
\end{lem}
\begin{proof} If $\rho\in]0,\rho_{\bar G}]$ is small enough, so that  proposition \ref{diedre4} applies, we know that there exists  
$\backslash\!\!\!{\mathfrak C}_5 =\mu\,{\mathfrak C}_5 >0$ depending only on $D\subset(M,g)$ such that,
for all $i,j=0,\dots,N'-1 $ one has
\begin{equation*}\Vert\,R_0^{\backslash\!\!\! G_{i,j}}\,\Vert_{0,0} \leq {\mathscr C}'= 2\,\backslash\!\!\!{\mathfrak C}_5 \,  \rho^2 \  .
\end{equation*}
From proposition \ref{ReggeGBbis} (extended \ref{ReggeGB}) we get 
\begin{equation*} \Vert\,\sum_{j=0}^{N'-1}{\mathcal P}^0_{\Gamma_{j-1}}\circ (({\mathcal P}^0_{\Gamma_{j}})^{-1}-({\mathcal P}^0_{\Gamma_{j-1}})^{-1})-\sum_{i,j=0}^{N'-1} R_0(\backslash\!\!\!\Delta_{i,j}) \,\Vert_{0,0}\leq N'^3{\mathscr C}'^2\  .
\end{equation*} 
The result follows from the inequalities $\rho\leq 1/\sqrt{2{\mathfrak C}_1}$ and $N'^3{\mathscr C}'^2\leq {\mathfrak C}_9\,\beta_{\bar G}^3\rho\,$ (proposition \ref{Riesqnb2}) joined to lemma \ref{pff}.
\end{proof}
\begin{rem}
Proving theorem {\rm\ref{theorem G}}, our strategy is based on applying proposition {\rm\ref{prox1}}. To control the distance between ${\mathcal P}_\Gamma^{-1}$ and $[({\mathcal P}_\Gamma^0)^{-1}]_T$ in $({\bf G}_x,\underline g_x)$, we construct a piecewise $\underline g_x$-geodesic (piecewise differentiable) joining $\hbox{\rm Id}$ to $[({\mathcal P}_\Gamma^0)^{-1}]_T$ through the finite sequence of vertices $[({\mathcal P}_{\Gamma_j}^0)^{-1}]_T$, as $j$ goes from $0$ to $N'-1$. The next lemmas {\rm\ref{forgottenfact}} and {\rm\ref{thG2}} are the central tools in order to establish the $\underline g_x$-proximity of this piecewise geodesic with the given curve $({\mathcal P}_{\Gamma_t})^{-1}$. 
\end{rem}
\begin{lem} \label{forgottenfact} With the hypotheses of theorem {\rm \ref{theorem G}}, one may find ${\mathfrak K}_1\!>\!0$ depending on $D\!\subset(M,g)$ and verifying for any $\rho\!>\!0$ (if small enough) and
for any $t\!\in[t_j,t_{j+1}]$ such that $\Gamma_t$ does not meet $\hat K_{n-2}$
\begin{equation} \label{forgetit1} d_{\underline g_x}([{\mathcal P}_{\Gamma_t}^0]_T,[{\mathcal P}_{\Gamma_j}^{0}]_T)
\leq {\mathfrak K}_1\,\beta_{\bar G}\,\rho  \  ,\  
\ d_{\underline g_x}([{\mathcal P}_{\Gamma_t}^0]_T,[{\mathcal P}_{\Gamma_{j-1}}^{0}]_T)
\leq {\mathfrak K}_1\,\beta_{\bar G}\,\rho  \  ,
\end{equation}
as well as
\begin{equation} \label{forgetit2} d_{\underline g_x}\!\!([\!(\!{\mathcal P}_{\Gamma_t}^0\!)^{\!-1\!}]_T,[\!(\!{\mathcal P}_{\Gamma_j}^{0}\!)^{\!-\!1}]_T\!)
\!\!
\leq\!\!  {\mathfrak K}_1\beta_{\bar G}\rho  ,\, d_{\underline g_x}\!\!([\!(\!{\mathcal P}_{\Gamma_t}^0\!)^{\!-1\!}]_T,[\!(\!{\mathcal P}_{\Gamma_{j\!-\!1}}^{0}\!)^{\!-\!1}]_T\!)
\!\!
\leq\!\!  {\mathfrak K}_1\beta_{\bar G}\rho  \, .
\end{equation}
\end{lem}
\begin{proof}
As above, we know from proposition \ref{diedre4} that, for 
$\backslash\!\!\!{\mathfrak C}_5\!=\!\mu\,{\mathfrak C}_5\! >\!0$, we have
for all $i,j=0,\dots,N'-1 $
\begin{equation*}\Vert\,R_0^{\backslash\!\!\! G_{i,j}}\,\Vert_{0,0} \leq {\mathscr C}'=2\,\backslash\!\!\! {\mathfrak C}_5 \,\rho^2 \  .
\end{equation*}
Then, we have (apply proposition \ref{Riesqnb2})
\begin{equation*} \Vert {\mathcal P}^0_{\Gamma_{j-1}}\circ({\mathcal P}^0_{\Gamma_j})^{-1}-\hbox{\rm Id}\Vert_{0,0}\leq  2\,N'\,{\mathscr C}'\leq
2\,{\mathfrak C}_8\,\beta_{\bar G}\,\rho\  .
\end{equation*}
Indeed, set $z\!=\!N'\,{\mathscr C}'\!\in\![0,1]$. Computing like in (\ref{binome1088}), one may write
\begin{multline*} \Vert {\mathcal P}^0_{\Gamma_{j-1}}\circ ({\mathcal P}^0_{\Gamma_j})^{-1}-\hbox{\rm Id}\Vert_{0,0}\leq \sum_{i=0}^{N'-1}\Vert\,R_0(\backslash\!\!\!\Delta_{i,j})\,\Vert_{0,0}+ \cr
+\!\!\!\!\sum_{0\leq i_1<i_2\leq N'-1}\!\!\!\!\!\!\!\!\!\!\!\Vert R_0(\!\backslash\!\!\!\Delta_{i_2,j}\!)\circ R_0(\!\backslash\!\!\!\Delta_{i_1,j}\!)\Vert_{0,0}+\cdots+
\!\Vert R_0(\!\backslash\!\!\!\Delta_{N'-1,j}\!)\!\circ\!\dots\!\circ\! R_0(\!\backslash\!\!\!\Delta_{0,j}\!)\Vert_{0,0}\!\leq
\\\leq z+\frac{z^2}{2}+\frac{z^3}{3!}+\dots+\frac{z^l}{l!}+\dots\leq\cr
\leq z\,(1+\frac{1}{2}+\frac{1}{2^2}+\dots+\frac{1}{2^{l-1}}+\dots)=2\,z\  .
\end{multline*}
As, for any $t\in [t_j,t_{j+1}]$, each of the loops $\Gamma_{i,j}^t$ (for a given $j$ and any $i=0,1,\dots,N'-1$) factoring out the decomposition of $\Gamma_{j-1}\vee\Gamma_t^{-1}$ is a curve included in one of the temperate squares $\backslash\!\!\!G_{i,j}$ coming out in the decomposition of $\Gamma_{j-1}\vee\Gamma_j^{-1}$ (the corresponding decomposition of $\Gamma_t\vee\Gamma_j^{-1}$ is handled in a similar way), one also gets (with $\Gamma_t\cap\hat K_{n-2}\!=\!\emptyset$)
\begin{gather} \label{Regge35} \Vert {\mathcal P}^0_{\Gamma_t}\circ({\mathcal P}^0_{\Gamma_j})^{-1}-\hbox{\rm Id}\Vert_{0,0}\leq  2\,N'\,{\mathscr C}'\leq 2\,{\mathfrak C}_8\,\beta_{\bar G}\,\rho
\  ,\\\label{Regge35bis}\  \Vert {\mathcal P}^0_{\Gamma_{j-1}}\circ({\mathcal P}^0_{\Gamma_t})^{-1}-\hbox{\rm Id}\Vert_{0,0}\leq 2\,N'\,{\mathscr C}'\leq 2\,{\mathfrak C}_8\,\beta_{\bar G}\,\rho\  .
\end{gather}
If $\rho\leq1/\sqrt{2\,{\mathfrak C}_1}$, applying lemma \ref{pff}, one derives
\begin{equation*}\Vert [{\mathcal P}^0_{\Gamma_t}\circ(\!{\mathcal P}^0_{\Gamma_j}\!)^{\!-\!1}\!-\!\hbox{\rm Id}]_T\Vert_{x,x} \!\leq\! 6\,{\mathfrak C}_8\,\beta_{\bar G}\,\rho
\  ,
 \Vert [{\mathcal P}^0_{\Gamma_{j\!-\!1}}\!\!\!\circ(\!{\mathcal P}^0_{\Gamma_t}\!)^{\!-\!1}\!-\!\hbox{\rm Id}]_T\Vert_{x,x}\!\leq \!6\,{\mathfrak C}_8\,\beta_{\bar G}\,\rho\,.
\end{equation*}
{\it Choose} $\beta_{\bar G}\,\rho \leq {\bf r}'_1/(6\,\sqrt{n}\,{\mathfrak C}_8)$, then {\it lemma {\rm\ref{nequiv2}} can be applied}, and set ${\mathfrak K}_1=6\,{\bf c}_2{\mathfrak C}_8$. It implies (\ref{forgetit2}) from the two inequalities just derived (left-invariance of $\underline g_x$ over ${\bf G}_x={\bf GL}^+(T_xM)$), while (\ref{forgetit1}) follows in the same way from (\ref{Regge35}) and (\ref{Regge35bis}), observing that ${\mathcal P}^0_{\Gamma_{j-1}}$ {\it and} ${\mathcal P}^0_{\Gamma_t}$ are $g_0$-isometries (apply twice the first equation in remark \ref{inv1}).
\end{proof}
\begin{lem} \label{thG2}  With the hypotheses of theorem {\rm \ref{theorem G}}, there exist uniform constants ${\bf A},{\bf B},{\mathfrak K}_2$ with $0\!<\!\max({\bf A},{\bf B})\!\leq\! {\mathfrak K}_2$ verifying for any $\rho>0$ (if small enough), in the division of $G\!\in\!\bar G\!\in\!{\mathcal V}_c$ into temperate squares of proposition {\rm \ref{Riesqnb2}}, 
the following inequalities hold true for any $j\!=\!0,1,\dots,N'-1$ and any $t\!\in\![t_j,t_{j+1}]$ such that $\Gamma_t$ does not hit a bone
\begin{gather} \label{Regge34}  d_{\underline g_x}\!({\mathcal P}_{\Gamma_t},[{\mathcal P}_{\Gamma_t}^{0}]_T)
\leq {\bf A}\,\beta_{\bar G}\rho\  ,\  \hbox{hence}\  \  d_{\underline g_x}\!({\mathcal P}_{\Gamma_j},[{\mathcal P}_{\Gamma_j}^{0}]_T)
\leq {\bf A}\,\beta_{\bar G}\rho\ ,\\ 
\label{Regge36}
d_{\underline g_x}\!([{\mathcal P}_{\Gamma_t}^0]_T,[{\mathcal P}_{\Gamma_j}^{0}]_T)
\!\leq\! {\bf B}\,\beta_{\bar G}\rho  \  \hbox{and}
\ d_{\underline g_x}\!([{\mathcal P}_{\Gamma_t}^0]_T,[{\mathcal P}_{\Gamma_{j-1}}^{0}]_T)
\!\leq \!{\bf B}\,\beta_{\bar G}\rho \ ,
\\ \label{Regge37}
d_{\underline g_x}\!([{\mathcal P}_{\Gamma_j}^0]_T,[{\mathcal P}_{\Gamma_{j-1}}^0]_T)\!\leq\! {\bf B}\,\beta_{\bar G}\rho
\  \  \ \hbox{and}
\ \ \  
d_{\underline g_x}\!({\mathcal P}_{\Gamma_j},{\mathcal P}_{\Gamma_{j-1}})\!\leq\! {\mathfrak K}_2\,\beta_{\bar B}'\rho \  ,\\ 
\label{Regge38}
d_{\underline g_x}\!({\mathcal P}_{\Gamma_t},[{\mathcal P}_{\Gamma_j}^0]_T)\leq {\mathfrak K}_2\,\beta_{\bar G}\rho\  \  \ \hbox{and}
\ \ \  
d_{\underline g_x}\!({\mathcal P}_{\Gamma_t},[{\mathcal P}_{\Gamma_{j-1}}^0]_T)\leq {\mathfrak K}_2\,\beta_{\bar G}\rho\  ,
\end{gather}
and the same list holds if {\em together} replacing all the maps ${\mathcal P}_{\Gamma_t}$,
${\mathcal P}_{\Gamma_j}$, ${\mathcal P}_{\Gamma_{j-1}}$, 
$[{\mathcal P}_{\Gamma_t}^0]_T$,
$[{\mathcal P}_{\Gamma_j}^0]_T$, $[{\mathcal P}_{\Gamma_{j-1}}^0]_T$ by their inverses.
\end{lem} 
\begin{proof} Set ${\mathfrak K}_2=3\,\max(4\,{\bf c}_2\,{\mathfrak C}_2\,,{\mathfrak K}_1))$.
\par
{\it Choosing} $\beta_{\bar G}\,\rho\leq  {\bf r}'_1/({4\,\sqrt{n}{\mathfrak C}_2})$, lemma \ref{length} implies
\begin{equation*}\Vert {\mathcal P}_{\Gamma_t}-[{\mathcal P}_{\Gamma_t}^{0}]_T\Vert_{x,x}
\leq 4\,{\mathfrak C}_2\,\beta_{\bar G}\,\rho\leq {\bf r}'_1/\sqrt{n}\ ,
\end{equation*}
thus 
$\Vert {\mathcal P}_{\Gamma_t}^{-1}\circ[{\mathcal P}_{\Gamma_t}^{0}]_T\!-\!{\rm Id}\Vert_{x,x}
\!\leq\! {\bf r}'_1/\sqrt{n}$ (since ${\mathcal P}_{\Gamma_t}$ is a $g_x$-isometry).
Apply lemma \ref{nequiv2}, we get for $\rho$ small enough and $t\!\in\![t_j,t_{j+1}]$ with $\Gamma_t\cap \hat K_{n-2}\!=\!\emptyset$, any $j\!=\!0,1,\dots,N'-1$, that (\ref{Regge34}) holds with ${\bf A}\!=\!4\,{\bf c}_2\,{\mathfrak C}_2$ since
$$ d_{\underline g_x}({\mathcal P}_{\Gamma_t},[{\mathcal P}_{\Gamma_t}^{0}]_T)
\leq 4\,{\bf c}_2\,{\mathfrak C}_2\,\beta_{\bar G}\,\rho\leq 
{\mathfrak K}_2\,\beta_{\bar G}\,\rho\  .
$$
\par Setting ${\bf B}\!=\!{\mathfrak K}_1$, lemma \ref{forgottenfact} implies (\ref{Regge36}).
\par  The first inequality in (\ref{Regge37}) reduces to (\ref{Regge36}).
\par The second inequality in (\ref{Regge37}) is derived from a combination of the first and from two applications of (\ref{Regge34}).
\par Finally (\ref{Regge38}) also follows from a combination of (\ref{Regge34}) and (\ref{Regge36}).
\par The case where ${\mathcal P}_{\Gamma_t}$,
${\mathcal P}_{\Gamma_j}$, ${\mathcal P}_{\Gamma_{j-1}}$, 
$[{\mathcal P}_{\Gamma_t}^0]_T$,
$[{\mathcal P}_{\Gamma_j}^0]_T$, $[{\mathcal P}_{\Gamma_{j-1}}^0]_T$ are replaced by their inverses is entirely parallel since
lemmas \ref{length} and \ref{forgottenfact} also hold true while doing this replacement.
\end{proof}
\begin{defn}\label{hardrain} 
One defines a piecewise $\underline g_x$-geodesic ${\mathfrak D}$ in ${\bf GL}^+(T_xM)$ which is built, for $j=0,1,\dots,N'-1$, by the shortest $\underline g_x$-geodesic ${\mathfrak D}_j$ joining $[({\mathcal P}_{\Gamma_{j-1}}^0)^{-1}]_T$ to $[({\mathcal P}_{\Gamma_j}^0)^{-1}]_T$ as $t$ runs in $[t_j,t_{j+1}]$ from $t_j$ to $t_{j+1}$.\par
Recall that ${\mathfrak G}$ is the path in ${\bf O}_g^+(T_xM)$ which consists in the curve $\{{\mathcal P}_{\Gamma_{\tau}}^{-1}={\mathcal P}_{\Gamma_{1,\tau}}^{-1}\mid \tau\in[0,1]\}$, joining $\hbox{\rm Id}$ to ${\mathcal P}_{\Gamma_{1,1}}^{-1}={\mathcal P}_\Gamma^{-1}\,$.
\end{defn}

\begin{lem} \label{thG3} With the hypotheses of theorem {\rm \ref{theorem G}}, given $\bar G\!\in\!{\mathcal V}_c$ there exist ${\mathfrak K}_3>0$ and a polynomial ${\mathfrak K}_4(\beta_{\bar G}^2)$ of degree $3$ in $\beta_{\bar G}^2$ with nonnegative coefficients (independent of $\beta_{\bar G}$) such that, for $\rho\!>\!0$ small enough, in the polyhedron paired with $\rho$ through theorem {\rm\ref{theorem I}}, one has
\begin{gather*} (i) \ \ \ d_{\underline g_x}({\mathfrak D},{\mathfrak G})\leq {\mathfrak K}_3\,\beta_{\bar G}\,\rho\ ;\\
(ii)\ \Vert[\sum_{j=0}^{N'\!-\!1}\! {\mathcal P}_{\Gamma_{j\!-\!1}}^0\!\!\circ\! (\!({\mathcal P}_{\Gamma_j}^0)^{\!-\!1}\!\!-\!({\mathcal P}_{\Gamma_{j-1}}^0)^{\!-\!1})]_T-\!\int_{{\mathfrak D}}\!A^{\!-\!1}\,dA\Vert_{x,x}\!\leq\!  {\mathfrak K}_4(\beta_{\bar G}^2)\,\beta_{\bar G}\rho\, .
\end{gather*}
\end{lem}
\begin{proof} As ${\mathfrak D}_j={\mathfrak D}_{\mid [t_j,t_{j+1}]}$ is the shortest $\underline g_x$-geodesic joining $[({\mathcal P}_{\Gamma_{j-1}}^0)^{-1}]_T$ to $[({\mathcal P}_{\Gamma_j}^0)^{-1}]_T$, one has, for any $t\in [t_j,t_{j+1}]$
\begin{multline} \label{thG33} \max(d_{\underline g_x}({\mathfrak D}_j(t),\![({\mathcal P}_{\Gamma_{j-1}}^0)^{\!-1}]_T),d_{\underline g_x}({\mathfrak D}_j(t),\![({\mathcal P}_{\Gamma_j}^0)^{\!-1}]_T)\!\leq\\ \leq \! d_{\underline g_x}([({\mathcal P}_{\Gamma_j}^0)^{-1}]_T,[({\mathcal P}_{\Gamma_{j-1}}^0)^{-1}]_T)\,.
\end{multline}
For any $j=0,1,\dots,N'-1\,,\,t\in[t_j,t_{j+1}]$, using the triangle inequality, (\ref{thG33}), then (\ref{Regge37}), (\ref{Regge38}) (lemma \ref{thG2}), setting ${\mathfrak K}_3=2\,{\mathfrak K}_2\,$, we get (i) 
\begin{equation*} d_{\underline g_x}\!(\!{\mathfrak D}_j(t),{\mathfrak G}(t)\!)
\!\leq \! d_{\underline g_x}\!(\!{\mathfrak D}_j(t),\![(\!{\mathcal P}_{\!\Gamma_{\!j\!-\!1\!}}^0)^{\!-\!1\!}]_T\!)\!+\!d_{\underline g_x}\!(\![(\!{\mathcal P}_{\!\Gamma_{\!j\!-\!1}\!}^0)^{\!-\!1\!}]_T,\!(\!{\mathcal P}_{\Gamma_t}\!)^{\!-\!1\!})\!\leq {\mathfrak K}_3\,\beta_{\bar G}\rho\,.
\end{equation*}  
\par
One has $\sum_{j=0}^{N'-1}[{\mathcal P}_{\Gamma_{j-1}}^0]_T\circ\int_{{\mathfrak D}_j}\,dA=\!\!
\sum_{j=0}^{N'-1} [{\mathcal P}_{\Gamma_{j-1}}^0\circ (({\mathcal P}_{\Gamma_j}^0)^{\!-1\!}-({\mathcal P}_{\Gamma_{j-1}}^0)^{\!-\!1})]_T\  .
$
Thus, one can write
\begin{multline} \label{thG32} \Vert\sum_{j=0}^{N'-1}\int_{{\mathfrak D}_j}([{\mathcal P}_{\Gamma_{j-1}}^0]_T-A^{-1})\circ dA
\Vert_{x,x}=\\=\Vert\sum_{j=0}^{N'-1}\int_{t_j}^{t_{j+1}}\big{(}([{\mathcal P}_{\Gamma_{j-1}}^0]_T-({\mathfrak D}_j(t))^{-1})\circ\frac{d{\mathfrak D}_j}{d t}(t)\big{)}\,d\,t
\Vert_{x,x}\leq\\\leq
\sum_{j=0}^{N'-1}\int_{t_j}^{t_{j+1}}\Vert[{\mathcal P}_{\Gamma_{j-1}}^0]_T-({\mathfrak D}_j(t))^{-1}\Vert_{x,x}\,\Vert \frac{d{\mathfrak D}_j}{d t}(t)\Vert_{x,x}\,d\,t
\leq
\\ \leq 2\, (1+{\bf K}\ \beta_{\bar G}^2+ 4\,{\mathfrak C}_2/{\mathfrak C}_8)\,{\bf c}_2\ {\mathfrak K}_2\,(\hbox{\rm length}_{\Vert \ \Vert_{x,x}}\!{\mathfrak D})\,\beta_{\bar G}\,\rho\, .
\end{multline}
{\em Indeed}, if $\rho>0$ small enough, thanks to (\ref{Regge37}), get
\begin{equation*} d_{\underline g_x}\!({\mathfrak D}_j(t),{\mathfrak D}(t_j)\!)\leq d_{\underline g_x}\!([(\!{\mathcal P}_{\Gamma_j}^0\!)^{\!-1\!}]_T,[(\!{\mathcal P}_{\Gamma_{j\!-\!1}}^0\!)^{\!-1\!}]_T\!)\!\leq\! {\mathfrak K}_2\,\beta_{\bar G}\,\rho\leq {\bf r}_1\,.
\end{equation*}
Use the left invariance of $d_{\underline g_x}$ and lemma \ref{nequiv2} to get
\begin{equation} \label{thG33bisss} \Vert ({\mathfrak D}(t_j)\!)^{\!-1}{\mathfrak D}_j(t)\!-\!{\rm Id}\Vert_{x,x}\!\leq\!  {\bf c}_2\,d_{\underline g_x}\!(({\mathfrak D}(t_j)\!)^{\!-1}{\mathfrak D}_j(t),{\rm Id})\!\leq\! {\bf c}_2\,{\mathfrak K}_2\,\beta_{\bar G}\,\rho\ .
\end{equation}
\begin{rem}
In general, $t\mapsto{\mathfrak D}^{-1}(t)$ is not a piecewise $\underline g_x$-geodesic. 
\end{rem}
\noindent  {\it If} $\beta_{\bar G}\rho\leq 1/2 {\bf c}_2{\mathfrak K}_2$ is true, (\ref{thG33bisss}) 
implies
$\,\Vert ({\mathfrak D}(t_j)\!)^{\!-1}{\mathfrak D}_j(t)\!-\!{\rm Id}\Vert_{x,x}\leq 1/2\,$. 
Thus 
lemma \ref{continv} gives
\begin{equation} \label{thG33bisss!}
\Vert({\mathfrak D}_j(t)\!)^{-1}{\mathfrak D}(t_j)-{\rm Id}\Vert_{x,x}
\leq 2\,\Vert ({\mathfrak D}(t_j)\!)^{\!-1}{\mathfrak D}_j(t)\!-\!{\rm Id}\Vert_{x,x}\  .
\end{equation}
By lemma \ref{club}, one has (using (\ref{parallprox2}))
\begin{equation}\label{thG33bisss!!}\Vert\!({\mathfrak D}(t_j)\!)^{-1}\!\Vert_{x,x}\!=\!\Vert[{\mathcal P}_{\Gamma_{j-1}}^0]_T\Vert_{x,x}\!\leq\! 1+{\bf K}\ \beta_{\bar G}^2+ 4\,{\mathfrak C}_2/{\mathfrak C}_8\,,
\end{equation} 
so that one derives from (\ref{thG33bisss}), (\ref{thG33bisss!}) and (\ref{thG33bisss!!})
\begin{multline}\label{thG1044}
\Vert{[\mathcal P}_{\Gamma_{j-1}}^0\!]_T\!-\!({\mathfrak D}_j(t))^{-1}\Vert_{x,x}\!=\!\Vert(\!({\mathfrak D}_j(t)\!)^{-1}{\mathfrak D}(t_j)\!-\!{\rm Id})(\!{\mathfrak D}(t_j)\!)^{-1}\Vert_{x,x}\!\leq\cr\leq
\!\Vert\!({\mathfrak D}(t_j)\!)^{-1}\!\Vert_{x,x}\,\Vert({\mathfrak D}_j(t)\!)^{-1}{\mathfrak D}(t_j)-{\rm Id}\Vert_{x,x}\leq\cr\leq
 2\, (1+{\bf K}\ \beta_{\bar G}^2+ 4\,{\mathfrak C}_2/{\mathfrak C}_8)\,{\bf c}_2\,{\mathfrak K}_2\,\beta_{\bar G}\,\rho
\  ,
\end{multline}
and (\ref{thG32}) {\em is finally seen to be true} from (\ref{thG1044}).
\par
Furthermore, one can apply lemma \ref{nequiv1} (use (\ref{thG33bisss})) to get (use (\ref{equiv222})) for any $j=0,\dots, N'-1$ and $t\in[t_j,t_{j+1}]$
\begin{equation}  \label{longueur0}
\Vert ({\mathfrak D}(t_j)\!)^{-1}\frac{d{\mathfrak D}_j}{d t}(t)\Vert_{x,x}\!\leq\! {\bf c}_1\Vert ({\mathfrak D}(t_j)\!)^{-1}\frac{d{\mathfrak D}_j}{d t}(t)\Vert_{\underline g_x}\!\!=\!{\bf c}_1\Vert \frac{d{\mathfrak D}_j}{d t}(t)\Vert_{\underline g_x}\,.
\end{equation}
Thus, using again (\ref{thG33bisss!!}), one obtains
\begin{equation}\label{club!}\hbox{\rm length}_{\Vert \ \Vert_{x,x}}\!{\mathfrak D} \!\leq\!
{\bf c}_1(1+{\bf K}\ \beta_{\bar G}^2+ 4\,{\mathfrak C}_2/{\mathfrak C}_8)\ \hbox{\rm length}_{\underline g_x}\!{\mathfrak D}
\end{equation}
Lemma \ref{forgottenfact} and the bound on $N'\rho$ given in proposition \ref{Riesqnb2} give
\begin{equation} \label{lengthgeod}
\hbox{\rm length}_{\underline g_x}\!\!{\mathfrak D}\!=\!\!\!\sum_{j=0}^{N'-1}\!\!d_{\underline g_x}\!(\![{\mathcal P}_{\Gamma_j}^0]_T,\![{\mathcal P}_{\Gamma_{j-1}}^0]_T\!)\!\leq\! N'\rho\,{\mathfrak K}_1\,\beta_{\bar G}\!\leq \!{\mathfrak C}_7{\mathfrak K}_1\ \beta_{\bar G}^2\, .
\end{equation}
Defining 
${\mathfrak K}_4(\beta_{\bar G}^2)=2\,{\bf c}_1\,{\bf c}_2\,{\mathfrak K}_1\,{\mathfrak K}_2\, {\mathfrak C}_7\,(1+{\bf K}\ \beta_{\bar G}^2+ 4\,{\mathfrak C}_2/{\mathfrak C}_8)^2\, \beta_{\bar G}^2\,$ establishes $(ii)$ thanks to (\ref{thG32}),(\ref{longueur0}) and (\ref{lengthgeod}). 
\end{proof}
\begin{lem} \label{thG4} With the hypotheses of theorem {\rm\ref{theorem G}}, one can find a polynomial ${\mathfrak K}_5(\beta_{\bar G}^2)$ of degree $1$ in $\beta_{\bar G}^2$ with nonnegative coefficients (independent of $\beta_{\bar G}$) such that,
for any $\rho>0$ small enough, one has in the polyhedron $T,K,(M,g)$ associated to $\rho$ through theorem {\rm\ref{theorem I}}
\begin{equation*}\Vert\int_{{\mathfrak D}}\!A^{-1}\,dA -\int_{{\mathfrak G}}A^{-1}\,dA\,\Vert_{x,x}\leq {\mathfrak K}_5(\beta_{\bar G}^2)\,\beta_{\bar G}\,\rho\  .
\end{equation*}
\end{lem}
\begin{proof}
This lemma is a consequence of lemma \ref{thG3} (i) and proposition \ref{prox1}. They tell us that there exists a constant ${\bf c}_3$ depending only on $({\bf GL}(T_xM),\underline g_x)$ such that, if $\beta_{\bar G}\rho\leq {\bf r}_1/2{\mathfrak K}_3$, one has
\begin{multline} \label{thG1099}\Vert\int_{{\mathfrak D}}\!\!A^{-1}\,dA-\int_{{\mathfrak G}}\!\!A^{-1}\,dA\Vert_{x,x}
\leq\\\leq ({\bf c}_3\,(\hbox{\rm length}_{\underline g_x}\!{\mathfrak G}+\hbox{\rm length}_{\underline g_x}\!{\mathfrak D})+1)\,{\mathfrak K}_3\,\beta_{\bar G}\rho\,.
\end{multline}
We also have the following
\begin{lem} For any ${\mathfrak G}$ as above, the length of ${\mathfrak G}$ verifies 
\begin{equation} \label{lengthinG} \hbox{\rm length}_{\underline g_x}\!{\mathfrak G}\leq {\bf K}\,\beta_{\bar G}^2\ ,
\end{equation}
where ${\bf K}$ defined in {\rm(\ref{BorneCourb})} (compare with lemma {\rm\ref{Pouf}}).
\end{lem}
\begin{proof} Set $y\!=\!G(\varrho,\tau)$. From (\ref{formulatralala}), as 
${\mathcal P}_{A_{\varrho,\tau}}, {\mathcal P}_{B_{\varrho,\tau}}$ are isometries between $(T_xM,g_x)$ and $(T_{y}M,g_{y})$ while ${\mathcal P}_{\Gamma^{\varrho,\tau}}$ is an isometry of $(T_{y}M,g_{y})$,  bringing in (\ref{extraloop1}) $\Gamma^{\varrho,\tau}=B_{\varrho,\tau}\vee A_{\varrho,\tau}^{-1}$, setting $R_y=R_{G(\varrho,\tau)}(\frac{\partial G}{\partial s},\frac{\partial G}{\partial t})$, get
\begin{multline*}  
\Vert\frac{\partial  {\mathcal P}_{\Gamma_{1,\upsilon}}^{-1}}{\partial \upsilon}(1,\tau)\Vert_{\underline g_x}\leq\int_0^1 \Vert{\mathcal P}_{B_{\varrho,\tau}}^{-1}\circ R_y \circ {\mathcal P}_{A_{\varrho,\tau}}\Vert_{\underline g_x}\ d\,\varrho=\cr=
\int_0^1 \Vert{\mathcal P}_{A_{\varrho,\tau}}^{-1}\circ {\mathcal P}_{\Gamma^{\varrho,\tau}}^{-1}\circ R_y \circ {\mathcal P}_{A_{\varrho,\tau}}\Vert_{\underline g_x}\ d\,\varrho=\cr=
\int_0^1 \Vert{\mathcal P}_{\Gamma^{\varrho,\tau}}^{-1}\circ R_y\Vert_{\underline g_{y}}\ d\,\varrho
=\int_0^1 \Vert R_y \Vert_{\underline g_{y}}\ d\,\varrho\  .
\end{multline*}
Thus, one gets, using definition \ref{normcurv} (\ref{Courbnorm})
\begin{equation*}  \hbox{\rm length}_{\underline g_x}{\mathfrak G} \leq \int\int_{\hbox{\tiny\rm rge}(G)}
\Vert R_\cdot^G \Vert_{\underline g_{\cdot}}\,d \,\hbox{\rm area}\leq {\bf K}\ \hbox{\rm area} (\hbox{\rm rge}(G))\, .
\end{equation*}
And, by proposition \ref{sqRiem} (\ref{caraf2}), the area of $\hbox{\rm rge}(G)$ is $\leq \beta_{\bar G}^2 \  .\qedhere $
\end{proof}
Define ${\mathfrak K}_5(\beta_{\bar G}^2)\!:=\!{\mathfrak K}_3\,({\bf c}_3\,({\mathfrak C}_7{\mathfrak K}_1\!+\!{\bf K})\,\beta_{\bar G}^2\!+\!1)$, lemma \ref{thG4} follows by plugging (\ref{lengthgeod}) and (\ref{lengthinG}) into (\ref{thG1099}).
\end{proof}
\noindent{\bf End of proof of theorem \ref{theorem G}}$ $

Setting ${\mathfrak C}_{11}(\beta_{\bar G}^2):=3{\mathfrak C}_9\,\beta_{\bar G}^2+{\mathfrak K}_4(\beta_{\bar G}^2)+{\mathfrak K}_5(\beta_{\bar G}^2)\,$, the proof uses lemmas \ref{thG1}, \ref{thG3} (ii), \ref{thG4} as well as equalities (\ref{formula15}) and (\ref{genangle}) of theorem \ref{theorem B}. The general case $r\in[0,1]$ is handled in the same way.
\end{proof}

\subsection{Pulling the nets}$ $\label{pulling33}

\noindent All parametrised squares $G$ are embedded, $x\!=\!G(0,0)$ is the base-point.
\begin{lem} \label{Reggeparallt} Given a piecewise flat polyhedron $(T,K,g_0)$ embedded in $M$ (as in theorem {\rm\ref{theorem I}}), a parametrised square $G\in\bar G\,$, where $\bar G\!\in\!{\mathcal V}_c$ (theorem {\rm\ref{text-1133}}), with $\hbox{\rm rge}(G)\!\subset\! T(K)$ transverse to $\hat K_{n\!-\!2}$ and $\hat K_{n\!-\!3}\subset M$, bring in the decomposition given by lemma {\rm\ref{psquare5}}. {\em Define} the {\em simply connected} set $\hbox{\rm rge}^\bullet(G)$ to be $\hbox{\rm rge}(G)$ deprived from the curves
$\gamma_{t_{j+1}}^{s_{i+1,j},1}\!\vee\!\delta_{s_{i+1,j}}^{\upsilon_{i,j},t_{j+1}}\!\vee\!\gamma_{\upsilon_{i,j}}^{\sigma_{i,j},s_{i+1,j}}$ where $(\sigma_{i,j},\upsilon_{i,j})\!\!\in]s_{i,j},s_{i+1,j}[\times]t_j,t_{j+1}[$ is such that $G_{i,j}(\sigma_{i,j},\upsilon_{i,j})\!=\!G(\sigma_{i,j},\upsilon_{i,j})$ belongs to a bone $\hat\xi\!\in\! \dagger\!\hat K_{n-2}\,$. 
\par There exists a canonical $g_0$-parallel translation on $T^{-1}(\hbox{\rm rge}^\bullet(G))$ {\rm identifying} $T_{T^{-1}(x)}{\mathcal K}$ with all tangent spaces
$T_{T^{-1}(y)}{\mathcal K}$ for $y\in\hbox{\rm rge}^\bullet(G)$. 
\par For any $y\!\in\! \hbox{\rm rge}(G)$, {\em define} a privileged path $\tilde A_y$ from $T^{-1}\!(x)$ to $T^{-1}\!(y)$.
If $y\!=\!G(s,t)$ with $s$ or $\,t\!=\!0\,$, set $\tilde A_y\!:=\!T^{-1}(\delta_0^{0,t})$ or $\tilde A_y\!:=\!T^{-1}(\gamma_0^{0,s})\,$. If $y\!=\!G_{i,j}(s,t)$ with $(s,t)\!\in]s_{i,j},s_{i+1,j}]\times]t_j,t_{j+1}]$ for some $i,j\,$, set 
\begin{equation}\tilde A_y\!:=\!T^{-1}(\gamma_t^{s_{i,j},s}\!\vee\!\delta_{s_{i,j}}^{t_j,t}\!\vee\!\gamma_{t_j}^{0,s_{i,j}}\!\vee\!\delta_0^{0,t_j})\,.
\end{equation}
For any $y\!\in\!{\rm rge}^\bullet(G)$ and $\epsilon>0$ small, there exists a well-defined path ${\tilde A}^\epsilon_y\!\subset\!T^{-1}\!(\hbox{\rm rge}^\bullet(G)\!)$ from $T^{-1}\!(x)$ to $T^{-1}\!(y)$ as close as wished to ${\tilde A}_y$, allowing us to {\em set} ${\mathcal P}^0_{{\tilde A}_y}\!\!:=\!{\mathcal P}^0_{{\tilde A}^\epsilon_y}$ (this is of use to apply theorem {\rm\ref{theorem I}} $(iii)$).
\par If $z\!=\!G(\sigma_{i,j},\upsilon_{i,j})\!\in\! \hat K_{n\!-\!2}$, {\em define} ${\mathcal P}^0_{{\tilde A}_z}\!\!:=\!lim_{\eta\rightarrow0}{\mathcal P}^0_{{\tilde A}_{y_{_{\eta}}}}$ where $y_\eta$ is chosen to be $y_\eta\!:=\!G(\sigma_{i,j}\!-\!\eta,\upsilon_{i,j})$ for $\eta\!>\!0$ small enough.
\end{lem}
\begin{proof} Indeed, ${\mathcal K}$ is an ambient flat manifold and 
$T^{-1}\!(\hbox{\rm rge}^\bullet(G)\!)\subset {\mathcal K}$ is simply connected. Any path joining $T^{-1}\!(x)$ to $T^{-1}\!(y)$ contained in $T^{-1}\!(\hbox{\rm rge}^\bullet(G)\!)$ produces the same $g_0$-parallel translation (apply theorem \ref{theorem A}), proving the first claim. 
\par Draw a scheme to prove the second claim. Given
$y\!=\!G_{i,j}(s,t)$ with $(s,t)\!\in]s_{i,j},s_{i+1,j}]\times]t_j,t_{j+1}]$ for some $i,j$ and $\epsilon\!>\!0$ small, define ${\tilde A}^\epsilon_y\!:=\!T^{-1}\!(\gamma_t^{s_{i,j}\!+\!\epsilon,s}\!\vee\!\delta_{s_{i,j}\!+\!\epsilon}^{t_j\!+\!\epsilon,t}\!\vee\!\gamma_{t_j\!+\!\epsilon}^{0,s_{i,j}\!+\!\epsilon}\!\vee\!\delta_0^{0,t_j\!+\!\epsilon})$ running from $T^{-1}\!(x)$ to $T^{-1}\!(y)$.
If $y$ is in some ${\rm rge}(G_{i,j})\setminus\partial {\rm rge}(G_{i,j})$ and $\epsilon\!>\!0$ is small enough, ${\tilde A}^\epsilon_y$ is included in $T^{-1}\!({\rm rge}^\bullet(G)\!)$ and $T({\tilde A}^\epsilon_y)$ is of $g$-length $\!\leq \!2\beta_{\bar G}\,$. Finally, given a point $y\!\in\!{\rm rge}^\bullet(G)$ and $\epsilon>0$ small enough, one also has  
${\tilde A}^\epsilon_y\!\subset\!T^{-1}\!(\hbox{\rm rge}^\bullet(G)\!)$.
\end{proof}
\begin{defn} \label{Reggemeas} The {\it Regge curving} ${\mathscr R}_0$ of the polyhedron $(T,K,g_0)$ in $(M,g)$ is, for almost any $G \!\in\!\bar G$ (with $\bar G\!\in\!{\mathcal V}_c\,$, see theorem {\rm\ref{text-1133}}) with $\hbox{\rm rge}(\bar G)\!\subset \!T(K)$ transverse to $\hat K_{n-2}$ and $\hat K_{n-3}$, at $x\!=\!G(0,0)\notin \hat K_{n-2}$, the $g_0$-skew-symmetric linear map ${\mathscr R}_0(G)\!\in\!\hbox{\bf End}(T_xM)$ which rephrases $\underline R_0^G$ of definitions {\rm\ref{Regge22}}, {\rm\ref{ReggeRegge!!}} and reads on $w\!\in\! T_xM$ (through $T$) with $\Vert\xi\Vert_0=1$ for all $\xi$ (use lemma {\rm\ref{Reggeparallt}})
\begin{equation*} {\mathscr R}_0(G)\,(w): = \sum_{\xi\in\dagger\! K_{n-2}}\!\!\!-\alpha_\xi\ \hbox{\rm index}_{G,\xi}\  dT\,(\ast_0(\xi\wedge dT^{-1} (w)))\  .
\end{equation*} 
\end{defn}

\begin{defn} \label{Riemannparallt}\label{Gaussmeas} Given a parametrised square $G$ in $(M,g)$, select the parallel translation ${\mathcal P}_{A_{\varrho,\tau}}$ on $\hbox{\rm rge}(G)$ which, for any $y=G(\rho,\tau)\in\hbox{\rm rge}(G)$, translates along $A_{\varrho,\tau}=\gamma_\tau^{0,\varrho}\vee\delta_0^{0,\tau}$ from $(T_xM,g_x)$ to $(T_yM,g_y)$.
The {\it Gauss curving} ${\mathscr R}$ is, for any $x\in M$, on any embedded square $G$ such that $x=G(0,0)$ and any $w\in T_xM$, the linear $g$-skew-symmetric map ${\mathscr R}(G)\in\hbox{\bf End}(T_xM)$ defined by (for $R_\cdot^G$, see definition {\rm\ref{normcurv}} ({\rm\ref{Courbnorm}}))
\begin{equation*} {\mathscr R}(G)\,(w) = \int\!\!\int_{{\rm rge}\,(G)}{\mathcal P}_{A_{\cdot}}^{-1}\circ R_\cdot^G \circ {\mathcal P}_{A_{\cdot}}(w) \,d\,{\rm area}_g\  .
\end{equation*} 
\end{defn}
Before stating our central convergence theorem, we quote a 
\begin{cor}\label{Regge-2} (of proposition {\rm\ref{Riesqnb2}}, lemmas {\rm\ref{Riesqnb}}, {\rm\ref{Regge-2033}}, definition {\rm\ref{ReggeRegge!!}}) Given any $G\in\bar G$ as above, $\rho>0$ small enough, one knows ${\mathcal N}\!\leq\! \backslash \!\!\!{\mathfrak C}_6(\beta_{\bar G}/\rho)^2$ ($\backslash \!\!\!{\mathfrak C}_6$ is defined in {\rm(\ref{red1})}). Thus,
setting ${\mathfrak C}\!=\!{\mathfrak C}_5^2\backslash \!\!\!{\mathfrak C}_6(\mu+1)$ 
\begin{equation*}\Vert \sum_{i,j=0}^{N'-1}\!\!R_0({\backslash \!\!\!\Delta_{i,j}})-\underline R_0^G\Vert\leq {\mathfrak C}\,\beta_{\bar G}^2\,\rho^2 \  .
\end{equation*}
\end{cor}
\begin{proof} Collecting and a direct computation establish the claim.
\end{proof}
\begin{thm} \label{convergence00} On any compact domain $D$ included in an open Riemannian manifold $(M,g)$, the Riemannian metric $g$ can be strongly approximated in the metric sense by a family of polyhedra $(T,K,(M,g))$ depending on a mesh $\rho>0$ (tending to $0$) and linked to a well-defined piecewise flat metric $g_0= g_0(\rho)$. In this approximation, given a family $\bar G\in{\mathcal V}_c$ (theorem {\rm\ref{text-1133}}), the Regge curving ${\mathscr R}_0(G)$ converges towards the Gauss curving ${\mathscr R}(G)$ uniformly in $G\in\bar G$.
\end{thm}
\begin{proof} The first part rephrases theorem \ref{theorem I}, the second (using (\ref{decloop3333}) of lemma \ref{generalisation1}) results from theorem \ref{theorem G} and above corollary \ref{Regge-2}.
\end{proof}
Before ending this section, we give another interpretation of the above defined curvings. First observe that, denoting as usual by $R$ the Riemannian curvature tensor in $(M,g)$, given two vectors $u,v\in T_xM$ one can interpret the skew-symmetric linear mapping $R(u,v)$ from $T_xM$ to $T_xM$ as an element of $\bigwedge^2(T_xM)$, denoted by $R(u\wedge v)$, through the following equality, true for any $w$ and $z\in T_xM$ (the natural scalar product induced by $g_x$ on $\bigwedge^2(T_xM)$ is still denoted by $\langle\cdot , \ast\rangle_x$)
\begin{equation}\label{curvoperator}\langle R(u\wedge v), z\wedge w\rangle_x=\langle R(u,v)w,z \rangle_x\  .
\end{equation}
Actually, $u\wedge v\mapsto R(u\wedge v)$ extended by linearity from $\bigwedge^2(T_xM)$ to itself is the {\it curvature operator}, a symmetric mapping, and the above equality, stated for any $u,v,w,z$, can be taken as its definition.
\begin{nota} \label{Regge10333} Given $G$ as in the above definitions {\rm\ref{Reggemeas}} and {\rm\ref{Gaussmeas}}, {\it one denotes by the same symbols the curvings} ${\mathscr R}_0(G)$ and ${\mathscr R}(G)$ which we understand this time {\it as elements of} $\bigwedge^2(T_xM)$.
\end{nota}
In accordance with (\ref{curvoperator}), one has for any $w,z\in T_xM$
\begin{equation*}\langle{\mathscr R}(G),z\wedge w\rangle=\langle{\mathscr R}(G)(w),z\rangle\  .
\end{equation*}
Viewed as an element of $\bigwedge^2(T_xM)$ and acting similarly to ${\mathscr R}(G)$, the Regge curving ${\mathscr R}_0(G)$ reads (definition \ref{Regge22} for the meaning of $\xi$ below)
\begin{equation} \label{Reggemeas2} {\mathscr R}_0(G)=\sum_{\xi\in \dagger\!K_{n-2}}\,\alpha_\xi\,\hbox{\rm index}_{G,\xi}\  dT(T^{-1}(x))(\ast_0 \,\xi)\  .
\end{equation}
{\em Indeed}, from the definition of the Hodge star operator $\ast_0$
 (see definition \ref{Regge-200}), if $\Omega_0$ is a canonical volume element of $(T_p{\mathcal K},g_0)$, one has
\begin{multline*}\langle\ast_0\,\xi,dT^{-1}(z\wedge w)\rangle_0=
\Omega_0(\xi\wedge dT^{-1}(z\wedge w))\!=\\=\!-\Omega_0((\xi\wedge dT^{-\!1}(w)\!)\wedge dT^{-\!1}(z)\!)\!=\!-\langle \ast_0(\xi \wedge dT^{-\!1}(w)\!),dT^{-1}(z)\rangle_0\  .
\end{multline*}
In view of lemma \ref{Reggeparallt}, the $g_0$-parallel translation $({\mathcal P}^0_{\tilde A_z})^{-1}={\mathcal P}^0_{\tilde A^{-1}_z}$ along the path $\tilde A^{-1}_z$ running from $T^{-1}(z), z=G(\varrho,\tau)\in\hat\xi\cap \hbox{\rm rge}(G)$ to $T^{-1}(x)$ with $x=G(0,0)$ {\em is not explicitly written} (${\mathcal P}^0_{\tilde A_y}$ in lemma \ref{Reggeparallt} is extended by letting $y\in{\rm rge}^\bullet(G)$ tend to $z$). As will soon be seen to be useful while introducing formula (\ref{ghost1}), we now present a detailed and apparently redundant version of (\ref{Reggemeas2}), based on the
\begin{defn}\label{restrindex} The restricted index of a parametrised square $G$ and a bone $\xi$ at $z$ is, for any family of open subsquares $G'\subset G$ sharing the local orientation of $G$ and decreasing to $z$
\begin{equation*}\hbox{\rm index}_{G,\xi,z}=\lim_{G'\searrow \{z\}}\hbox{\rm index}_{G',\xi}\  .
\end{equation*}
\end{defn}
\begin{lem} \label{Reggecomplet}  The Regge curving {\rm(\ref{Reggemeas2})} also
reads (this involves lemmas {\rm\ref{singular1}}, {\rm\ref{generalisation1}}, definition {\rm\ref{ReggeRegge!!}} and also lemma {\rm\ref{Reggeparallt}})
\begin{equation*} {\mathscr R}_0(G): = \sum_{{\xi\in \dagger\!K_{n-2}}\atop{z\in\hat\xi\cap\hbox{\tiny\rm rge}(G)}}\!\!\!\alpha_\xi\ \hbox{\rm index}_{G,\xi,z}\  dT(T^{-1}(x))\,(\ast_0\,{\mathcal P}^0_{\tilde A^{-1}_z}(\xi) )\  .
\end{equation*} 
\end{lem}

\section{The Regge theorem}\label{Reggetheorem}

\subsection{The theorem}$ $

We develop a proof of Regge's theorem inspired by Regge's original approach in his founding and intuitive paper \cite{R}.

\begin{convent} \label{bone1066} If $D\subset(M,g)$ is a compact domain in a Riemannian manifold, if $\rho>0$ is small, denoting by $T,K,(M,g)$ a polyhedral approximation on $D$ paired with $\rho>0$ as in theorem {\rm\ref{theorem I}} (we continue to skip $\rho$ in the writing of $T,K,(M,g)$, {\em which actually depend on} $\rho$), {\em $\xi$ bone} means all the time that $\xi\in \dagger\!K_{n-2}$ is such that $T(\xi)=\hat\xi\subset D\subset T(K\setminus \partial K)$. {\em Call} $\dagger\!K_{n-2}(D)$ the set of such bones.
\end{convent}

\begin{thm}\label{Regge1234} ({\rm See Regge \cite{R}}) As $\rho\!>\!0$ tends to $0$, if $(T,K,M)$ is a polyhedral approximation to $(M,g)$ on $D\subset T(K)$ according to convention {\rm\ref{bone1066}}
\begin{equation*} \sum_{\xi\in \dagger\!K_{n-2}(D)}\!\!\!\!\!\alpha_\xi\,\hbox{\rm Vol}_{n-2}^{g_0}(\xi) \ \ \rightarrow\ \  {\bf c}_1(n)\!\int_D\!\! \hbox{\rm scal}_g(p)\  dp\  ,
\end{equation*}
where ${\bf c}_1(n)$ is a constant depending only on $n$.
\end{thm}

The proof of this theorem fills this section \ref{Reggetheorem}.

\subsection{Riemannian and quasi-Riemannian squares}
\label{Riemannian1111}$ $

Among the parametrised squares, there are distinguished natural ones, connected to the given Riemannian manifold $(M,g)$, the {\em Riemannian squares} (definition \ref{Riemannian}). Such a square is embedded if its {\em side} is small enough. In general it does not share the features of theorem {\rm \ref{text-1133}} but can be approached, as close as desired, by elements of the residual ${\mathcal V}_c$ (theorem {\rm \ref{text-1133}}), elements to which the theorem applies.

\begin{assump} \label{tacitassump} On technical grounds, we enlarge $D$ (recall $D\!\subset\! T(K)\!\subset\!W\!\subset\!(M,g)$) to a close compact domain $D'\subset T(K)$ containing (the compact) $D$ in its interior (with $\partial D'$ close to $\partial D$).
 Choose ${\mathcal R}>0$ such that $B(x,2{\mathcal R})$ is convex (definition \ref{D.convex}) in $M$ for all $x\in D'\subset M$ and ensuring, for any $x\in D'$, that $\exp_x$ is a $1/2$-quasi-isometric map from $B(0_x,2{\mathcal R})$ to $B(x,2{\mathcal R})$ (put $C=1/2$ in proposition {\rm \ref{proposition A}}). And
 ${\mathcal R}$ is also assumed to be chosen according to lemma {\rm \ref{L.constante}}.
\end{assump}
\begin{rem} The choice of ${\mathcal R}>0$ relies on the geometry of $D$, for $D',W$ are as close to $D\!\subset\!(M,g)$ as wished. Thus, in the definition of some constants, {\it we often skip making an explicit reference to} ${\mathcal R}>0\,$.
\end{rem}

\begin{defn} \label{Riemannian} Given $x\!\!\in \!\!D'\!\subset \!(M,g)$, orthonormal vectors $u,v\!\in\! T_xM,\,r\!\in]0,{\mathcal R}/\sqrt{2}[$, call {\em Riemannian parametrised square} (over $D'$) of side $r$ the map 
\begin{equation*} G_r:(s,t)\in[0,1]\times[0,1]\longmapsto \exp_x (r(s\,u+t\,v)) \in B(x,{\mathcal R})\subset M\  ,
\end{equation*}
so it is the image under $\exp_x$ of a true planar square of side $r$
\begin{equation*}\Pi_r=\{r(s\,u+t\,v)\mid (s,t)\in [0,1]\times[0,1]\}\subset B(0_x,{\mathcal R})\  .
\end{equation*}
Choosing some open bounded interval $I$ containing $[0,1]$, we actually assume that the above definition of $G_r$ extends to all $(s,t)\in \bar I^2$.
\par
{\em Implicitly or explicitly}, a {\em Riemannian parametrised square} has a chosen {\em side} $r>0$.
Given ${\mathcal R}>0$ as above, a {\em Riemannian parametrised square} is surely well-defined and embedded if $r$ belongs to $]0,{\mathcal R}/\sqrt{2}[$. 
\end{defn}
\begin{defn} Given $r\!\in]0,{\mathcal R}/\sqrt{2}[$, call {\em exhaustive family of Riemannian squares $\underline G_r$ over $D'$} the family of Riemannian squares of side $r$ parametrised by the {\em whole} Stiefel manifold $\Theta':={\bf St}_2(TD')$ of orthonormal pairs of vectors in $(T_xM,g_x)$ as $x$ runs in $D'$; set also $\Theta:={\bf St}_2(TD)$ and ${\mathcal T}:={\bf St}_2(TM)$. Thus $c\!=\!\hbox{\rm dim}\,{\bf St}_2(TM)\!=\!n+ \hbox{\rm dim}\,{\bf O}(n)-\hbox{\rm dim}\,{\bf O}(n\!-\!2)=3\,(n-1)\,$.
\end{defn}
\begin{nota} \label{danslacase} We often write $\underline{G}$ thinking to any family $\underline{G}_r$, omitting $r\in ]0,{\mathcal R}/\sqrt{2}[$ implicit in the context. But also, {\em through extending}, $\underline{G}$ sometimes points at the whole collection
of those $\underline{G}_r$ for $r\in ]0,{\mathcal R}/\sqrt{2}[$.
\end{nota}
\begin{rem} \label{danslacasegolo} Given $f\!\in\!C^\infty(\bar I^2, M)$ and $\lambda\!\in]0,1]$, {\em set} $f_\lambda:=f(\lambda\cdot)\!\in\!C^\infty(\bar I^2, M)$. If $f$ and $h$ are close in $C^\infty(\bar I^2, M)$, so are $f_\lambda$ and $h_\lambda$ for any $\lambda\!\in]0,1]$: this follows by highlighting the dilations and restrictions in the definitions of $f_\lambda$ and $h_\lambda$. 
\end{rem}
\begin{lem} \label{familyRsq111} Choose a family $\bar G_{{\mathcal R}/\sqrt{2}}\!\in\!  C^\infty(\Theta'\!\times\! \bar I^2,M)$ which is ``Whitney-close to'' $\underline G_{{\mathcal R}/\sqrt{2}}$ (in the Whitney topology) and reads
$$\bar G_{{\mathcal R}/\sqrt{2}}\!:\!(\vartheta,s,t)\!\in\!\Theta'\!\times\! \bar I^2\!\mapsto \!\bar G_{{\mathcal R}/\sqrt{2}}(\vartheta,s,t)\!=\!G_\vartheta(({\mathcal R}/\sqrt{2})\,s,({\mathcal R}/\sqrt{2})\,t)\!\in\! M\ . 
$$
Define, for any $r\!\in ]0,{\mathcal R}/\sqrt{2}[$, the family of squares
$$\bar G_r(\vartheta,s,t):=\bar G_{{\mathcal R}/\sqrt{2}}\,(\vartheta,\lambda s,\lambda t)\ \ {\rm with}
\ \ \lambda=r\sqrt{2}/{\mathcal R}\ . 
$$ 
Then, for any $r\!\in ]0,{\mathcal R}/\sqrt{2}[$, $\underline G_r$ and $\bar G_r\!\in\!C^\infty(\Theta'\times\bar I^2, M)$ are close.
So, saying that {\em $\bar G\!\in\!{\mathcal V}_c$ is close to $\underline{G}$} of theorem {\rm\ref{text-1133}} {\em makes sense}.
\end{lem}
\begin{proof}
This expresses the theorem {\rm \ref{text-1133}} where $\Theta'\!:=\!{\bf St}_2(TD')$, the $c$-dimensional parameter space, builds the {\em exhaustive family of Riemannian squares} $\underline G$ {\em over} $D'$. Choosing any $r\!\in]0,{\mathcal R}/\sqrt{2}[$, remark \ref{danslacasegolo} establishes the claim (playing with $\lambda=r\sqrt{2}/{\mathcal R}$). 
\end{proof}
\begin{rem} We want to make clear that the symbols $\underline G$ and $\bar G$ carry two different meanings. First $\underline G$ (or $\bar G$) is a smooth mapping from $\Theta'={\bf St}_2(TD')$ to the Fr\' echet space $C^\infty(\bar I^2,M)$ (see \cite{G-G}, theorem {\rm1.11} page {\rm76}), second it is viewed as the image of this mapping, the family of parametrised squares embedded in $M$ and indexed by ${\bf St}_2(TD')$.  
\end{rem} 
\begin{lem}\label{CMS3}
\par
Any $\bar G$ close enough to $\underline G$ induces a one-to-one map
\begin{equation*}(x,u,v)\in{\bf St}_2(TD')\longmapsto \bar G(x,u,v)=G_{x,u,v}\in C^\infty(\bar I^2,M)\  .
\end{equation*}
\par\noindent
Further, one can find $\bar G\!\in\!{\mathcal V}_c$ (${\mathcal V}_c$ is defined according to theorem {\rm \ref{text-1133}} with ${\mathscr T}\!=\!{\bf St}_2(TD)$) as close as wished to $\underline G$ in the sense of lemma {\rm\ref{familyRsq111}}.
\end{lem}
\begin{proof} Given $G\in \bar G_r$, set 
${\mathfrak T}(G)\!:=\!(G,\frac{1}{r}\frac{\partial G}{\partial s},\frac{1}{r}\frac{\partial G}{\partial t})(0,0)\,$.
Introduce the map ${\mathfrak T}_{\bar G}$ for $\bar G=\bar G_r$ (and the corresponding ${\mathfrak T}_{\underline G}$ for $\underline G$)) 
\begin{equation*} {\mathfrak T}_{\bar G}:(x,u,v)\!\in\! {\bf St}_2(TD')\longmapsto {\mathfrak T}(G)\!:=\!(G,\frac{1}{r}\frac{\partial G}{\partial s},\frac{1}{r}\frac{\partial G}{\partial t})(0,0)
\in F_2(TM)\ .
\end{equation*}
where $G=G_{x,u,v}$ with $G\in\bar G$ (or with $G\in\underline G$), where $F_2(TM)$ is the bundle of free $2$-systems of vectors in $TM$ and $(\bar x,\bar u,\bar v):=(G, \frac{1}{r}\frac{\partial G}{\partial s}, \frac{1}{r}\frac{\partial G}{\partial t})(0,0)$ are $C^\infty$-functions of $x,u,v$. As ${\mathfrak T}_{\underline G}={\rm Id}$ and $\bar G$ near to $\underline G$ implies that ${\mathfrak T}_{\bar G}$ is close to ${\mathfrak T}_{\underline G}$ in the fine $C^1$-topology (extending the two maps nicely in an open neighborhood of $\,{\bf St}_2(TD')\,$), one infers (from theorem 3.10 of \cite{Mu} or
 \cite{Hi}, theorem 1.6 p.38) that ${\mathfrak T}_{\bar G}$ embeds $\,{\bf St}_2(TD')\,$ into $F_2(TM)$, which proves that $\bar G$ is one-to-one.
\par\vskip1mm  As for the second claim, notice that for $G\in\bar G$, the point $\bar x(x,u,v)=G_{x,u,v}(0,0)$ is in general distinct from $x\,$. Choose a domain $D''\supset D'$ such that $\bar x(x,u,v)\in D''$ for any $(x,u,v)\!\in\!{\bf St}_2(TD')\,$. If $D''$ is a diffeomorphic to a ball, the bundle $\,{\bf St}_2(TD'')\,$ is fibrewise diffeomorphic, through a fibre map $\theta$, to
$D''\times{\bf St}_2(\R^n)$ (localising, one may handle in this way trickier situations).
Set $\tilde x:=\bar x=G(0,0)$ for $G\in \bar G$. 
The Gram-Schmidt process gives from $(\bar u, \bar v)\in F_2(T_{\tilde x}D'')$ the new $(\tilde u=\bar u/\Vert \bar u\Vert_g, \tilde v)$ in $\,{\bf St}_2(T_{\tilde x}D'')\,$. Set $\tilde{\mathfrak T}_{\bar G}(x,u,v):=(\tilde x,\tilde u,\tilde v)$. Compose the map $\theta\circ \tilde{\mathfrak T}_{\bar G}$ from ${\bf St}_2(TD')$ to $D''\times{\bf St}_2(\R^n)$ by the natural projection on the factor ${\bf St}_2(\R^n)$. For any $x\!\in\!D'$, this composed map restricted to ${\bf St}_2(T_xD')$ is a diffeomorphism from the compact 
${\bf St}_2(T_xD')$ onto ${\bf St}_2(\R^n)$, which is connected if $n\!\geq\!3$ (if $n\!=\!2$, the surjectivity is also easy to handle). Since ${\mathfrak T}_{\bar G}$ and ${\mathfrak T}_{\underline G}$ are close enough on ${\bf St}_2(TD')$ in the fine $C^1$-topology, so are $\tilde{\mathfrak T}_{\bar G}$ and ${\mathfrak T}_{\underline G}={\rm Id}$ and one infers that
${\bf St}_2(TD)\subset\tilde{\mathfrak T}_{\bar G}({\bf St}_2(TD'))$.
\par {\em Given any $G\in\bar G$, define the new $\tilde G$} by setting for $(\tilde x,\tilde u,\tilde v)\!\in\!{\bf St}_2(TD)$
$$\tilde G_{\tilde x,\tilde u,\tilde v}(s,t)\!\!:=\!G_{x, u, v} (\!{\Lambda}(s,t)\!) \ \hbox{with}\ \Lambda\!:=\!(dG_{(x,u,v)}(0,0)\!)^{\!-\!1}\!\circ{\mathfrak L} \circ dG_{(x,u,v)}(0,0),
$$
and $\tilde{\mathfrak T}_{\bar G}(x,u,v)\!:=\!(\tilde x,\tilde u,\tilde v)$, where ${\mathfrak L}$ is the linear map in $T_{\tilde x}{\rm rge}(G)$ sending $(\bar u,\bar v)$ to $(\tilde u,\tilde v)$. Observe that ${\mathfrak T}_{\tilde G}=\hbox{\rm Id}$. If $\bar G$ is close enough to $\underline G$, then $(\bar x,\bar u, \bar v)\!\in\! F_2(TD'')$ is close to $(x,u,v)\!\in\! {\bf St}_2(TD')$, implying $(\tilde x,\tilde u,\tilde v)$ and $(x,u,v)$ are close in ${\bf St}_2(TM)$. Choose an open interval $J$ verifying $[0,1]\!\subset J\!\subset \bar J\!\subset I$. As $\bar G$ is a compact family over $D'$, if $\bar G$ is close enough to 
$\underline G$,
then $\Lambda$ is close to $\hbox{\rm Id}$
and sends the square $\{s\,e_1\!+\!t\,e_2\mid (s,t)\!\in \!\bar J^2\}$ into a parallelogram $\{s\,\Lambda (e_1)\!+\!t\,\Lambda (e_2)\mid (s,t)\!\in\!\bar J ^2\}$ {\it which is a subset of} $\bar I^2$ containing $[0,1]^2$. The $C^\infty$-mappings $\bar {\tilde G}$ and $\underline G$ both restricted to ${\bf St}_2(TD')\times \bar J ^2$ are close, the closer $\bar G$ to $\underline G$, the closer $\bar {\tilde G}$ to $\underline G$. All $(\tilde x,\tilde u,\tilde v)\!\in\!{\bf St}_2(TD)$ are described once in the reparametrised $\tilde G$. 
By definition of $\bar {\tilde G}$, the family of lines defined by $\tilde w\in{\bf S^1}$ (definition \ref{carrement}) inside some $\tilde G\in\bar {\tilde G}$ are reparametrised as lines defined by some $w\in{\bf S^1}$ inside $G\in\bar G$, the change of coordinates $\Lambda$ {\em being linear at the source}. Thus, the theorem {\rm \ref{text-1133}} also applies to $\bar {\tilde G}$.
Renaming $\bar G$ the derived
$\bar {\tilde G}$, the claim follows.
\end{proof}
\begin{defn} \label{CMS10333} A {\it generic exhaustive family of $r$-quasi-Riemannian squares} (with $r\in]0, {\mathcal R}/\sqrt{2}[$) in some domain $D\!\subset \!M$ is an element $\bar G_r$ of the residual set ${\mathcal V}_c\!\subset \!C^\infty({\bf St}_2(TD)\!\times \!\bar I^2,M)$ (defined in theorem {\rm \ref{text-1133}}). In view of lemma \ref{CMS3}, for any $r\in]0, {\mathcal R}/\sqrt{2}[$, $\bar G_r$ can be taken as $C^\infty$-close to $\underline G_r$ as desired as in lemma \ref{familyRsq111}.  
\end{defn}
\begin{rem} \label{shrink} Consider $G\in\bar G_r$ uniquely determined by $r\in]0,{\mathcal R}/\sqrt{2}[$ and $(x,u,v)\in {\bf St}_2(TD)$. 
Regarding questions like bounding the multiplicity of intersection points of a square $G\in\bar G_r\!\in\!{\mathcal V}_c$ with bones in some approximation $T,K,(M,g)$, what matters relies on the ``geometric support'' containing all squares $\hbox{\rm rge}(G_r)$.
\end{rem}

\subsection{Some integrals of curvatures}\label{someintegrals}$ $

\begin{defn} \label{elemvolume} Given an orientation on $D\subset(M,g)\,$, the Riemannian metric induced by $g$ on $\,{\bf St}_2(TD)\,$ defines a canonical volume element $d\,{\mathscr G}$ such that each
$({\bf St}_2(T_xD), d\,{\mathscr G}_x)$ is sent to $({\bf St}_2(\R^n),d\,{\mathpzc G})$ through an isometry sending $(T_xD,g_x)$ onto $(\R^n,{\rm std})$, where $d\,{\mathpzc G}$ is the probability ${\bf O}(n)$-invariant measure on ${\bf St}_2(\R^n)$.
Put this $d\,{\mathscr G}$ on $\bar G_r$ through its parametrisation by $\,{\bf St}_2(TD)\,$. If $f$ integrable satisfies ({\it the scaling invariance}) $f(G_r)\!=\!f(G_{r'})$ for any $r,r'\!\in\,]0,{\mathcal R}/\sqrt{2}]$ and the paired $G_r,G_{r'}\!\in\!\bar G$ (for $\bar G$ here, see notation \ref{danslacase}), then one has
\begin{equation*}\forall r,r'\in]0,{\mathcal R}/\sqrt{2}]\ \ \ \ \ \ \int_{G\in\bar G_r}f(G) \,d\,{\mathscr G}=\int_{G\in\bar G_{r'}}f(G) \,d\,{\mathscr G}\  .
\end{equation*}
\end{defn}
\begin{nota} Given a generic exhaustive family of quasi-Riemannian squares $\bar G$ on $D$ (definition {\rm\ref{CMS10333}}) and $G\!\in\!\bar G$, set $T_xG\!:=\!(\frac{\partial G}{\partial s}\!\wedge \!\frac{\partial G}{\partial t}){(0,0)}$, write $\Vert T_xG\Vert_g\!=\!1$, keeping similar {\rm(\ref{CMS10})}, {\rm(\ref{CMS9})}, {\rm(\ref{ghost2})} below, and define
\begin{equation}\label{CMS10} 
I^g(\bar G_r):=\frac{1}{r^2}\int_{G\in\bar G_r}\langle {\mathscr R}(G),\frac{T_xG}{\Vert T_xG\Vert_g}\rangle_g\  d {\mathscr G}\  .
\end{equation}
\end{nota}  
\begin{nota} In the sequel $g^0$ stands for $(T^{-1})^\ast g_0$ (since long ago, we write $T$ instead of $\hat T$; taking care of this abuse,
remark \ref{bibopalula} and definition \ref{D.metric1} ensure that $g^0$ is a local notation for the former $\hat g$).
\par
If a generic exhaustive family of quasi-Riemannian squares $\bar G$ is given on $D$, given a $\rho$-approximating polyhedron $K$ with its piecewise flat metric $g_0$, define the cousin 
\begin{equation}\label{CMS9}I^{g_0}(\bar G_r):=\frac{1}{r^2}\int_{G\in\bar G_r}\langle {\mathscr R}_0(G),\frac{T_xG}{\Vert T_xG\Vert_{g^0}}\rangle_{g^0}\  d {\mathscr G}\  .
\end{equation}
\end{nota} 
Looking at $\hat\xi\cap\hbox{\rm rge}(G)$ in $(M,g)$, observe (as $\bar G\!\in\!{\mathcal V}_c$) that it consists in finitely many temperate points $z$ of total index $\rm{index}_{G,\xi}$ the sum of subindices defined near each point $z\,$ by shrinking $G$, see lemma \ref{singular1}. 
In order to manage questions of convergence as $r\rightarrow0$, we introduce one more integral ${J}^{g_0}(\bar G_r)$, close to $I^{g_0}(\bar G_r)$, built on an analogue $\mathpzc{R}{}$ of the Regge curving $\mathscr{R}_0\,$, where the metric $g$ plays every role played by $g_0$ in the term $\ast_{0}  \,{\mathcal P}^{0}_{\tilde A_{z}^{-1}}(\xi)$  that appears in $\mathscr{R}_0\,$ (see lemmas \ref{Reggeparallt}, \ref{Reggecomplet}), thus turned into $\ast_g  \,{\mathcal P}_{\tilde A_{z}^{-1}}(\hat\xi_z)$, where ${\mathcal P}_{\tilde A_{z}^{-1}}$ stands for ${\mathcal P}_{T(\tilde A_{z}^{-1})}$ and $\hat\xi_z$ is the oriented $g$-normed $(n\!-\!2)$-vector supporting $T_z\hat\xi$, the orientation of $\hat\xi$ being induced by $\xi$ (see the corresponding equation (\ref{Reggemeas2}))
\begin{defn}  \label{Reggemeas2222}
Denoting by ${\mathcal P}_{\tilde A_{z}}$ the $g$-parallel translation from $x=G(0,0)$ to $z$, where $\tilde A_{z}$ is defined in lemma {\rm\ref{Reggeparallt}}), set
\begin{gather}\label{ghost1}\mathpzc{R}{}(G): = \sum_{{\xi\ {\rm bone}}\atop{z\in \hat\xi\cap{\rm rge}(G)}}\!\!\!\alpha_\xi\ \hbox{\rm index}_{G,\xi,z}\  \ast_g \!{\mathcal P}_{\tilde A_{z}}^{-1}(\hat\xi_z)\ ,\\\label{ghost2}
{J}{}^{g_0}(\bar G_r):=\frac{1}{r^2}\int_{G\in\bar G_r}\langle \mathpzc{R}{}(G),\frac{T_xG}{\Vert T_xG\Vert_g}\rangle_g\  d {\mathscr G}\  .
\end{gather}
Actually, ${J}{}^{g_0}(\bar G_r)$ as ${I}{}^{g_0}(\bar G_r)$ depend on $g$ and $g_0$ (through $\alpha_\xi$).
\end{defn}
\subsection{Proof of the Regge theorem... what remains to prove}
\label{sketchRegge}$ $
\begin{lem}\label{bornangles} Given a compact domain $D\subset (M,g)$ and a simplicial approximation resulting from theorem {\rm \ref{theorem I}}, there exists a constant ${\mathfrak A}$ depending on $D\subset(M,g)$ such that one has, uniformly in $\rho>0$ 
\begin{equation*}\sum_{\xi\in \dagger\!K_{n-2}(D)}\vert\alpha_\xi\vert\,
\hbox{\rm Vol}_{n-2}^{g}(\hat\xi)\leq {\mathfrak A}\,.
\end{equation*}
Moreover, theorem {\rm\ref{theorem I}} implies (as $g^0\!\!=\!(T^{-1})^*g_0\!\rightarrow \!g$ while $\rho$ tends to $0$)
\begin{equation*}\sum_{\xi\in\dagger\! K_{n-2}(D)}\vert\alpha_\xi\vert\,
\vert\hbox{\rm Vol}_{n-2}^{g}(\hat\xi)-\hbox{\rm Vol}_{n-2}^{g_0}(\xi)\vert\rightarrow0\  \  \hbox{as}\  \  \rho\rightarrow0\ .
\end{equation*}
\end{lem}
\begin{proof} The number of $n$-simplices $\hat\sigma\!\subset\! D, \sigma\!\in \!K$, is smaller than ${\rm Cst}\rho^{-n}$, where ${\rm Cst}>0$ depends on $D\subset(M,g)$: they are disjoint and their volumes have a uniform lower bound ${\mathfrak C}_4\rho^n$, see theorem \ref{theorem I} and remark \ref{remarquable1}.
Thus, the number of their $(n\!-\!2)$-faces is bounded by ${\rm Cst'}\rho^{-n}$, where ${\rm Cst'}>0$ depends on $D\subset(M,g)$.
\par  
A direct consequence of theorem \ref{theorem I} $(ii) (b)$ is that $T$ is a ${\mathfrak C}_1\rho^2$-quasi-isometry between the metric spaces $(\sigma,g_0)$ and $(\hat\sigma,g)$. So, there exists a uniform ${\mathfrak C}_1'>0$,  depending on $D\subset(M,g)$, such that
\begin{equation}\label{broutons}\vert\hbox{\rm Vol}_{n-2}^{g}(\hat\xi)-\hbox{\rm Vol}_{n-2}^{g_0}(\xi)\vert\leq {\mathfrak C}_1'\rho^2\,\hbox{\rm Vol}_{n-2}^{g}(\hat\xi)\  .
\end{equation}
\par But $\vert\alpha_\xi\vert$ has an upper bound ${\mathfrak C}_5\rho^2$, uniform in $\xi\!\in\! \dagger\!K_{n-2}(D)$, see convention \ref{bone1066} and proposition \ref{diedre4}.
As $\hbox{\rm Vol}_{n-2}^{g_0}(\xi)$ is $\leq\!{\rm Cst}_n\,\rho^{n-2}$ (here, $\xi$ is a $g_0$-Euclidean $(n-2)$-simplex of diameter $\leq \!{\rm Cst}'_n\,\rho$ and ${\rm Cst}_n,{\rm Cst}'_n\!>0$ depend on $n$), using (\ref{broutons}) one gets the first claim.
\par Using again (\ref{broutons}), the first claim implies the convergence.
\end{proof}
Now, an important technical definition sorting the quasi-Riemannian squares in terms of their situation with respect to the bones of a complex $K$ of mesh $\rho$.
\begin{defn}\label{CMS-7777} Given $\rho\!>\!0$ small, $T,K$ being a polyhedral approximation of $D\!\subset\!(M,g)$ of mesh $\rho$ (theorem {\rm\ref{theorem I}}), if $\bar G_r$ is an exhaustive family of quasi-Riemannian squares (see definition \ref{CMS10333}), {\em call} $\bar G_r^\rho$ the open sub-family of those $G\!\in\!\bar G_r$ such that for any bone $\xi$ verifying ${\rm rge}(G)\cap\hat\xi\!\not=\!\emptyset$
the distance from the origin $x_G$ to $\hat\xi$ is achieved at a point interior to $\hat\xi\,$.
\end{defn}
\begin{rem} \label{CMS-7777bise} The integrals ${I}{}^{g_0}(\bar G_r^\rho),{J}{}^{g_0}(\bar G_r^\rho)$ are well-defined.
\end{rem}
\begin{proof} Now the {\em proof of Regge's theorem:} it relies on theorem \ref{convergence00}, which implies the convergence of the above $I^{g_0}(\bar G_r)$ towards $I^g(\bar G_r)$ (see (\ref{CMS10}), (\ref{CMS9})) as $\rho$ tends to $0$, and on the two propositions below. 
\begin{prop}\label{CMS-7} Given a compact domain $D\subset(M,g)$ and an approximation $T,K,(M,g)$ on $D$ paired with $\rho\searrow 0$ through theorem {\rm \ref{theorem I}}, given an exhaustive family of quasi-Riemannian squares $\bar G$ over $D$, there exist constants ${\mathfrak C}$ (depending only on $D\subset(M,g)$), ${\bf c}_2(n)$ (lemma {\rm\ref{scalint2}}, ${\bf c}_2(n)$ depends only on $n$) such that (see definition {\rm\ref{CMS-7777}} for $\bar G_r^\rho$)
\par $(i)$ for $\rho$ small and any $r$ in some $]0,a]$, one has
\begin{equation*}
\vert\,{I}^{g_0}(\bar G_r^\rho)-{J}^{g_0}(\bar G_r^\rho)\vert\leq{\mathfrak C}\, \rho\  ;
\end{equation*}
\par $(ii)$  if $\rho>0$ is fixed, as $r>0$ tends to $0$, one has 
\begin{equation*} \vert\,{I}^{g_0}(\bar G_r\setminus\bar G_r^\rho)\vert
\rightarrow  0\ \ \ \hbox{and}\ \  \ 
\vert{J}^{g_0}(\bar G_r^\rho)-{\bf c}_2(n)\!\!\!\!\!\!\sum_{\xi\in \dagger\!K_{n-2}(D)}\!\!\alpha_\xi\,
\hbox{\rm Vol}_{n-2}^g(\hat\xi)\vert\rightarrow  0\  .
\end{equation*}
\end{prop}
\begin{proof} The proof fills all subsections following this one. 
\end{proof}
\begin{prop}\label{CMS11} Given a compact domain $D\subset(M,g)$ and $\bar G$ as in proposition {\rm\ref{CMS-7}}, let $dp$ be the $g$-volume-element of $(M,g)$ at a point $p\in M$. There exists a constant ${\bf c}_3(n)$, depending only on $n$, such that
\begin{equation*}
I^g(\bar G_r) \rightarrow {\bf c}_3(n)\int_D \hbox{\rm scal}_g(p)\  dp
\  \  \  \hbox{as} \  \  \  r \rightarrow 0\  .
\end{equation*}
\end{prop}
\begin{proof} For any $G\!\in\!\bar G_r\!\equiv\!{\bf St}_2(TM)$, if $r$ tends to $0$ then $\frac{1}{r^2}\,\langle {\mathscr R}(G),\frac{T_xG}{\Vert T_xG\Vert_g}\rangle_g$ tends to $\langle R(\frac{T_xG}{\Vert T_xG\Vert_g}),\frac{T_xG}{\Vert T_xG\Vert_g}\rangle_g$. In view of the definition of the scalar curvature $\hbox{\rm scal}_g(x)$, exploiting the transitive action of the orthogonal group ${\bf O}_{g_x}(T_xM)$ on ${\bf St}_2(T_xM)$ at $x\in D$, the result follows from the 
\begin{lem}\label{RCMS1033terrrrr} Given a Euclidean space $(V,\langle\cdot, \ast\rangle)$, an orthonormal basis $v_1,\dots,v_m$ of $V$, a subgroup $\Omega\subset {\bf O}(V)$ verifying $\varpi_i(v_1)=v_i$ for some $\varpi_i\!\in\!\Omega$ and any $i\!=\!1,\dots,m$, for any $L$ endomorphism of $V$, then
\begin{equation*}\hbox{\rm trace} \,L=\frac{m}{\hbox{\rm vol}_{d\varpi}(\Omega)}\int_{\{\varpi\in \Omega\}}\!
\langle L(\varpi(v_1)), \varpi(v_1)\rangle \, d \varpi\ ,
\end{equation*}
where $d \varpi$ is an invariant measure on $\Omega\,$.
\end{lem}
\begin{proof} For any $\varpi\in \Omega$, the system $\varpi(v_1),\dots,\varpi(v_m)$ is orthonormal too and one has $\hbox{\rm trace} \,L=\sum_{j=1}^m\langle L(\varpi(v_j)), \varpi(v_j)\rangle$. Using the invariance and, for $i\!=\!1,\dots,m$, the element $\varpi_i\!\in\!\Omega$ verifying $\varpi_i(v_1)\!=\!v_i\,$, compute
\begin{multline*}\hbox{\rm trace} \,L=\frac{1}{\hbox{\rm vol}_{d\varpi}(\Omega)}\int_{\{\varpi\in \Omega\}}\!
\sum_{j=1}^m\langle L(\varpi(v_j)), \varpi(v_j)\rangle \, d \varpi=\\=
\frac{m}{\hbox{\rm vol}_{d\varpi}(\Omega)}\int_{\{\varpi\in \Omega\}}\!
\langle L(\varpi(v_1)), \varpi(v_1)\rangle \, d \varpi\  .\qedhere
\end{multline*}
\end{proof}
Setting $\Omega\!=\!{\bf O}_{g_x}(T_xD), V\!=\!\bigwedge^2T_xD$, choosing $G\!\in\!{\bf St}_2(T_xD)$ and, in $V$, an orthonormal basis of decomposable vectors $v_1\!:=\!T_xG,v_2,\dots,v_m$, one may apply the above lemma and the equality
\begin{multline*}\int_{\{\varpi\in O_{g_x}(T_xD)\}}\!
\langle R(\varpi(T_xG)), \varpi(T_xG)\rangle \, \frac{d \varpi}{{\rm vol}_{d\varpi}(O_{g_x}(T_xD))}=\\=\int_{\{G\in {\bf St}_2(T_xD)\}}\!\langle R(T_xG),T_xG\rangle\,d\,{\mathscr G}_x\ 
\end{multline*}
(use definition \ref{elemvolume}) to get the result, setting ${\bf c}_3(n)\!:=\!{2}/{n(n-1)}\,$.
\end{proof}

\subsubsection{Taking proposition {\rm\ref{CMS-7}} for granted, proof of Regge's theorem.}$ $\label{Regge333333}

Postponing the proof of proposition \ref{CMS-7}, we establish Regge's theorem.
Choose $\epsilon>0$. 
Proposition \ref{CMS-7} (i),
theorems \ref{convergence00}, \ref{theorem I}, lemma \ref{bornangles} allow choosing $\rho\!>\!0$ such that for any $r\!\in]0,a]$ (for some chosen $a\!>\!0$)
\begin{gather}\label{CMS??}\vert {I}^{g_0}(\bar G_r^\rho)- {J}^{g_0}(\bar G_r^\rho)\vert\leq \epsilon/6\  ,
\\\label{Regge77} \vert\,I^g(\bar G_r)-I^{g_0}(\bar G_r)\,\vert \leq \epsilon/6\  ,
\\\label{CMS????}{\bf c}_2(n)\times\vert \sum_{\xi\in \dagger\!K_{n-2}(D)}\!\!\!\alpha_\xi\,(
\hbox{\rm Vol}_{n-2}^{g}(\hat\xi)-\hbox{\rm Vol}_{n-2}^{g_0}(\xi))\vert\leq \epsilon/6\  
\end{gather}
(one has $\hbox{\rm Vol}_{n-2}^{g_0}(\xi)=\hbox{\rm Vol}_{n-2}^{g^0}(\hat\xi)$). Once $\rho>0$ is chosen to fulfill the above inequalities, 
using proposition \ref{CMS-7} $(ii)$, one chooses $r>0$ small enough such that 
\begin{gather}\label{CMS?!}\vert\,{I}^{g_0}(\bar G_r\setminus\bar G_r^\rho)\vert\leq\epsilon/6
\\\label{CMS?}\vert {J}^{g_0}(\bar G_r^\rho)-{\bf c}_2(n)\!\!\sum_{\xi\in \dagger\!K_{n-2}(D)}\!\!\alpha_\xi\,
\hbox{\rm Vol}_{n-2}^{g}(\hat\xi)\vert\leq \epsilon/6\  .
\end{gather}
Thanks to proposition \ref{CMS11}, one can even choose $r>0$ small enough so that (\ref{CMS?!}), (\ref{CMS?}) hold together with the following inequality
\begin{equation}\label{CMS-11}\vert I^g(\bar G_r) -{\bf c}_3(n)\int_D \hbox{\rm scal}_g(p)\  dp\vert\leq \epsilon/6\  .
\end{equation}
Together, (\ref{CMS????}), (\ref{CMS?}), (\ref{CMS??}), (\ref{CMS?!}), (\ref{Regge77}),  (\ref{CMS-11}) prove the theorem, since one has  (the first inequality follows from (\ref{CMS?}) and (\ref{Regge77}))
\begin{multline}\vert{\bf c}_2(n)\!\!\sum_{\xi\in \dagger\!K_{n-2}(D)}\!\!\alpha_\xi\,
\hbox{\rm Vol}_{n-2}^{g}(\hat\xi)-{\bf c}_3(n)\int_D \hbox{\rm scal}_g(p)\  dp\vert\leq\\\leq\epsilon/3+\vert {J}^{g_0}(\bar G_r^\rho)-I^g(\bar G_r)\vert\leq 5\epsilon/6\  .
\end{multline}
One has ${\bf c}_1(n)={\bf c}_3(n)/{\bf c}_2(n)\,$.
\end{proof}
It remains to prove proposition \ref{CMS-7}. 

\subsection{An integral paired with an $(n\!-\!2)$-linear subspace of ${\mathbb R}^n$}
\label{integraleuclidean}$ $
\begin{lem} \label{CMS1000} Given an $(n-2)$-vector subspace ${\bf X}\subset ({\mathbb R}^n,{\rm std})$ (where ${\rm std}$ is the standard Euclidean structure), call ${\bf G}_r^{{\bf X},y}$ the sets of squares ${\bf G}\!:\!(\varrho,\varsigma)\!\in\![0,1]^2\!\mapsto \!x+\varrho ru+\varsigma rv\!\in\!\R^n$ uniquely defined by their side $r$, base point $x_{\bf G}\!:=\!x$ closest to $y\!\in\!{\bf X}$ and $(u,v)\in{\bf St}_2({\mathbb R}^n)$ (so $x_{\bf G}$ varies in ${\bf X}_y^\perp$, orthogonal plane to ${\bf X}$ at $y$). The ${\bf G}_r^{{\bf X},y}$ are invariant by rotations and  translations preserving ${\bf X}$. Such a translation sends ${\bf G}_r^{{\bf X},y_1}$ onto ${\bf G}_r^{{\bf X},y_2}$ if $y_1\in{\bf X}$ is sent to $y_2\in{\bf X}\,$. Given any isometry preserving orientation and ${\bf X}$, which maps ${\bf G}$ on ${\bf G}'$, one has ${\rm index}_{{\bf G},{\bf X}}= {\rm index}_{{\bf G}',{\bf X}}\,$. 
\par More, ${\bf G}_{r_2}^{{\bf X},y}$ is homothetic to ${\bf G}_{r_1}^{{\bf X},y}$ by the homothety of center $y\in{\bf X}\,$, dilation factor ${\bf r}\!=\!r_2/r_1$. Two {\rm homothetic} ${\bf G}_1\!\in\!{\bf G}_{r_1}^{{\bf X},y},\, {\bf G}_2\!\in\!{\bf G}_{r_2}^{{\bf X},y}$ {\rm together} intersect ${\bf X}$ or, {\rm together}, they do not, in fact ${\rm index}_{{\bf G}_1,{\bf X}}\!=\!{\rm index}_{{\bf G}_2,{\bf X}}$. Given $r$, ${\bf G}_{r}^{{\bf X},y}$ is a product ${\bf X}_y^\perp\times{\bf St}_2({\mathbb R}^n)$. In view of the last lines, we later may write ${\bf G}^{{\bf X},y}$ meaning ${\bf G}_r^{{\bf X},y}$.
\end{lem}

Assume orientations are given on ${\mathbb R}^n$ and ${\bf X}\,$. Select cylindrical coordinates $(t,\vartheta,y)\!\in\! {\mathbb R}_+^*\!\times\![0,2\pi[\times{\bf X}$ around ${\bf X}$. If $d{\mathpzc G}_x\!=\!d{\rm vol}_{{\bf St}_2(T_x{\mathbb R}^n)}$ is the usual probability measure induced by ${\rm std}$ on ${\bf St}_2(T_x{\mathbb R}^n)$, define on ${\bf G}_{r}^{{\bf X},y}$ the following product measure, invariant under changes of side $r$
\begin{equation}\label{meas1} d{\mathpzc G}_{\mid{\bf G}_{r}^{{\bf X},y}}\!=\!dx\wedge d{\mathpzc G}_x\!=\!t \,dt\!\wedge d\vartheta\!\wedge d{\mathpzc G}_x \ \ \ \hbox{with}\ \ \ x_{\bf G}=x\!\equiv\!(t,\vartheta,y)\,.
\end{equation}
\begin{nota} Given $r,t,\vartheta$, write $t=rs,\,x_{\bf G}\!=\!(rs,\vartheta,y)$ and define  
\begin{equation}\label{model111}{\bf G}^{\cap{\bf X},y}_{r}(s,\vartheta)\!:=\!\{{\bf G}\in {\bf G}^{{\bf X},y}_r\!\mid\! x_{\bf G}\!=\!(rs,\vartheta,y)\,\hbox{and} \ {\bf G}\cap{\bf X}\not=\emptyset\}\,.
\end{equation}
 Now, ${d {\mathpzc G}}_{\mid{\bf G}_{r}^{{\bf X},y}}$ restricted to 
${\bf G}^{\cap{\bf X},y}_{r}(s,\vartheta)$ is ${d {\mathpzc G}}\equiv{d {\mathpzc G}}_x$, which we also denote ${d {\mathpzc G}}_{s,\vartheta}$ (for ${d {\mathpzc G}}$, see definition \ref{elemvolume}). One has in those notations
\begin{equation}\label{CMS8bise33-1}
\frac{1}{r^2}\,{d {\mathpzc G}}_{\mid{\bf G}_{r}^{{\bf X},y}}=s \,ds\wedge d\vartheta\wedge {d {\mathpzc G}}_{s,\vartheta} 
\ .
\end{equation} 
\end{nota}
\begin{lem}\label{scalint2}  Here ${\bf G}$ also stands for $u\wedge v\,$. The Euclidean integral 
\begin{equation} \label{CMS8bise33-2} {\bf J}^{{\bf X},y}_{r}:=\frac{1}{r^2}\int_{{\bf G}^{{\bf X},y}_{r}}
\hbox{\rm index}_{{\bf G},{\bf X}}\,\cos\angle_{\rm std}(\ast_{\rm std}\, {\bf X}, {\bf G})\,d {\mathpzc G}_{\mid{\bf G}_{r}^{{\bf X},y}}\  .
\end{equation}
depends only on $n$, we call ${\bf c}_2(n)$ its value ${\bf c}_2(n):={\bf J}^{{\bf X},y}_{r}$.
\end{lem}
\begin{proof} One computes from (\ref{meas1}), (\ref{CMS8bise33-1}), (\ref{CMS8bise33-2})
\begin{equation*} {\bf J}^{{\bf X},y}_{r}=\int_0^{2\pi}\!\!\!\!\int_0^{\sqrt{2}}\!\!\!\!\!\!s\!\!\int_{{\bf G}^{{\bf X},y}_{r}(s,\vartheta)}
\!\!\!\!\!\!\!\!\!\!\!\!\hbox{\rm index}_{{\bf G},{\bf X}}\,\cos\angle_{\rm std}(\ast_{\rm std}\, {\bf X}, {\bf G})\,{d {\mathpzc G}}_{s,\vartheta}\,ds\,d\vartheta\,.
\end{equation*}
Set $f(\!{\bf G}\!)\!=\!\hbox{\rm index}_{{\bf G},{\bf X}}\,\cos\angle_{\rm std}(\ast_{\rm std}\, {\bf X}, {\bf G})$. As $r,\vartheta$ vary, because $f$ is invariant, all $\int_{{\bf G}^{\cap{\bf X},y}_{r}(s,\vartheta)}\!f(\!{\bf G}\!){d {\mathpzc G}}_{s,\vartheta}$ are so (use (\ref{model111}), apply lemma \ref{CMS1000}).
Taking into account the relative orientation of the pair $({\bf G},{\bf X})$, one has 
\begin{equation*}\label{scalint4} {\bf J}^{{\bf X},y}_{r}=2\pi\,\int_0^{\sqrt{2}}\!\!\!s\!\int_{{\bf G}^{\cap{\bf X},y}_{r}(s,\vartheta)}
\!\!\vert\,\cos\angle_{\rm std}(\ast_{\rm std}\, {\bf X}, {\bf G})\,\vert\,{d {\mathpzc G}}_{s,\vartheta}\,ds\  .
\end{equation*}
Independent of $r,y, {\bf X}$ and of the choice of the Euclidean structure ${\rm std}\,$, this integral ${\bf J}^{{\bf X},y}_{r}$ takes the value ${\bf c}_2(n):={\bf J}^{{\bf X},y}_{r}$.
\end{proof}

\subsection{Near to an $(n\!-\!2)$-face of an approximation}\label{bone1111}$ $

We first recall and adapt things at the heart of part \ref{M-T1}. Theorem \ref{theorem I}, gives us a polyhedral approximation $T,K,(M,g)$ on $D$ paired with a {\em thickness} $t_0\!>\!0$ and a {\em mesh} $\rho\!>\!0$. It is realised (see section \ref{rbs}) by {\it Riemannian barycentric simplices}, each $n$-simplex being defined by a {\em spread} set of $n+1$ vertices in a convex ball. Two Riemannian barycentric simplices are close if their vertices are.
\par A relatively compact open domain $W$ is chosen, $D$ is included in the interior of $T(K)\!\subset \!W\!\subset\! \bar W\!\subset\! M$. One chooses  $s'_{0,k}\!\in]0,s_{0,k}[$, where, for $k\!=\!0,1,\dots,n$, $s_{0,k}$ is small and related to $t_{0,k}$ as in definition \ref{texture33333}. We select open balls $B_i=B(x_i,\underline{r}), \tilde B_i=B(x_i,\rho_3)$ of center $x_i$ with $\underline{r}:=\min_{i}\rho_{x_i,k}$ (for $\rho_{x_i,k}$, see remark \ref{bon-bon}; $\rho_3$ appears in lemma \ref{text8}), with $\bar B_i\!\subset \!\tilde B_i$ and the $B_i$ building a finite covering of $\bar W$. Theorem \ref{texture333} holds for those choices of balls for any $k\!=\!0,1,\dots,n\,$: each $\tilde B_i$ is, for all $k\!=\!0,1,\dots,n\,$, the basis of a {\it texture} $X_i^k\!=\!X^k_{x_i,\rho_{x_i},{s'_{0,k}}}$ (definition \ref{texture3333}) having total space a finite-dimensional manifold.

\begin{fact} By lemma {\rm \ref{simpletext}}, we can describe precisely the set of small barycentric $k$-simplices $\hat\eta$ we are interested in (see remark {\rm\ref{compact1099}} below). Recall (details are given below) that any such $\hat\eta$ is included in some $\tilde B_i$ and extends as a $k$-submanifold $\tilde\eta$ such that $\tilde\eta_i\!:=\!\bar{\tilde B}_i\cap\tilde\eta$ is a differentiable embedded $k$-ball verifying $\partial\tilde\eta_i\!=\!\partial\bar{\tilde B}_i\cap\tilde\eta$ ($\tilde\eta$ is transverse to $\partial\bar{\tilde B}_i$).
\end{fact}

\par The construction of section \ref{RBT} not only comprises all $k$-faces of $n$-simplices contained in $W$ of the approximation $T,K,(M,g)$ on $D$, but also describes elements of a set, denoted by ${\mathfrak S}_i^k$, which are {\em all buildable $k$-Riemannian barycentric simplices $\hat\eta$ of thickness $\geq t_{0,k}\!>\!0$, diameter $\leq \!2\rho$ and barycentre $p$ belonging to $B_i$,}  as {\em submanifolds} of {\em leaves} in $X_i^k$ (if $k\!=\!n$, see definition \ref{texture33333}). We skip here $\rho,\,t_{0,k}$.

 In order to deal with {\em compact} (instead of open) sets below, {\it we now slightly shrink} each $B_i$ into a concentric ball, {\it still called} $B_i$, in such a way that all previous properties remain true
 (in particular, the shrinking keeps the new $B_i$ to still cover $\bar W$). The new ${\mathfrak S}_i^k$ {\it samples $k$-simplices having barycentre in the new $B_i$} while its closure $\overline{{\mathfrak S}_i^k}$ samples those with barycentre in the new $\bar B_i$.
Magnifying, one cylindrifies the cone manifold ${\mathbb F}^k_{\geq {s'_{0,k}}}(M)$ into $\Gamma^k_{\geq {{s'_{0,k}}}}(M)$ (definitions \ref{gloups1}, \ref{cylindrification}, proposition \ref{gloups2}, proof of lemma \ref{simpletext}). Thus, we have established the important
\begin{rem}\label{compact1099} Skipping $W,\rho,t_{0,k}$, the set $\,{\mathfrak S}^k$ of $k$-faces of all buildable Riemannian barycentric $n$-simplices with thickness $\geq {t_{0,k}}\!>\!0$, diameter $\!\leq\! 2\rho$, contained in $\bar W$ is compact and ${\mathfrak S}^k\!\subset\!\cup_i{\mathfrak S}_i^k$ is smoothly parametrised (lemma {\rm \ref{simpletext}}) by a compact subset (still called ${\mathfrak S}^k$) of the manifold below (for ${\mathbb F}_{\geq s'_{0,k},\underline{r}}$, see definition {\rm\ref{cylindrification}})
\begin{equation} \label{bone111133}{\mathpzc F}^k={\mathpzc F}^k(M):=
\{({\bf v},t)\mid{\bf v}\in 
{\mathbb F}_{\geq s'_{0,k},\underline r}^k(p),\,t\in[-1,1], \,p\in M\}\ .
\end{equation}
\end{rem}
For each given $i$, one has ${\mathfrak S}_i^k\subset \overline{{\mathfrak S}_i^k}\subset  {\mathpzc F}^k\,$. Here we  
focus on $k\!=\!n\!-\!2$.
\par
The {\em barycentric device} developped in sections {\rm \ref{rbs}}, {\rm \ref{6.3}}, {\rm \ref{RBT}} extends all $(n\!-\!2)$-simplices $\hat \xi\!\in\!{\mathfrak S}_i^{n\!-\!2}$ to local embedded $(n\!-\!2)$-submanifolds $\tilde\xi$.
\begin{lem} \label{compactsimpl} In sections {\rm \ref{6.3}}, {\rm \ref{RBT}}, all bones $\hat \xi\!\in\!{\mathfrak S}_i\!:=\!{\mathfrak S}_i^{n\!-\!2}$ are extended to compact differentiably embedded $(n\!-\!2)$-balls $\tilde\xi_i\!=\!\bar{\tilde {B}}_i\cap \tilde \xi$ of controlled size, with $\partial\tilde\xi_i\!=\!\partial\bar{\tilde {B}}_i\cap \tilde \xi$, all $\hat\xi$ being uniformly small with $\rho$ inside $\tilde\xi_i$ (this follows from remark {\rm \ref{bon-bon}}), where the $\tilde\xi_i$ are {\rm leaves of the texture} $X_i^{n\!-\!2}$ (theorem {\rm \ref{texture333}}). Moreover, $\bar{\mathfrak S}_i$ is smoothly parametrised by a compact subset of ${\mathpzc F}^{n\!-\!2}$, still called $\bar{\mathfrak S}_i\,$.
\end{lem}
\begin{assump} \label{unifsmall1} Working near $\hat\xi\in{\mathfrak S}_i$, we often skip $i$ in $\tilde\xi_i$.
\end{assump}
\begin{nota} \label{nortrunc}Given $\hat\xi$, call $\exp_{\tilde \xi}^\perp$ the normal exponential map to $\tilde \xi$. The normal $\tau$-tubular neighborhood around $\tilde \xi$ is the truncated ($\tau>0$)
$$[{\rm Tub}\,\tilde \xi]_\tau:=\{x\in M \mid x \in{\rm rge}(\exp_{\tilde \xi}^\perp)\ {\rm and}\  d_g(x,\tilde \xi)\leq \tau\}\  .
$$
Denote by ${\bf N}^1(\tilde\xi)$ the unit normal bundle to $\tilde\xi$ in $(M,g)$, i. e.
\begin{equation*}{\bf N}^1(\tilde \xi):=\cup_{y\in\tilde\xi}\{w\in T_{y}M\mid \langle w,T_y\tilde\xi\rangle_g=0\ \ \hbox{and}\ \ \Vert w\Vert_g=1\}\ .
\end{equation*}
\end{nota}
\begin{cor} \label{compactsimpl333} There exists $\tau\!>\!0$ small enough so that, for any $i$, any $\hat\xi\!\in\!\bar{\mathfrak S}_i$, the map $\exp_{\tilde \xi_i}^\perp$ from $\{v\!\in\!T_{\tilde \xi_i}M\mid v\perp \tilde\xi_i\ \hbox{and} \ \Vert v\Vert_g\!\leq\!\tau\}$ to $[{\rm Tub}\,\tilde \xi_i]_\tau$ is bijective
(by assumption {\rm\ref{unifsmall1}}, $\xi$ might stand for $\xi_i$).
\end{cor}
\begin{proof} The radius of injectivity of $\exp_{\tilde \xi_i}^\perp$ is positive and varies continuously with $\hat\xi\in\bar{\mathfrak S}_i$. But $\bar{\mathfrak S}_i$ is compact, see lemma \ref{compactsimpl}, and the indices $i$ are finitely many.
\end{proof}
We now make a development parallel to the previous section \ref{integraleuclidean}.
\par Choose orientations on both $\tilde \xi$ and an open neighborhood of $\tilde \xi$ in $M$ containing ${\rm [Tub\,\tilde \xi]}_\tau$. By the tubular neighborhood theorem, $\exp_{\tilde \xi}^\perp$ is a diffeomorphism from {\em the disk-bundle} 
\begin{equation}\label{diskbdle}{\bf D}^\tau(\tilde\xi)\!:=\!\{t\,w \mid w \in {\bf N}^1(\tilde\xi)\ \ \hbox{and}
\ \ t\in[0,\tau]\}
\end{equation} 
onto ${\rm [Tub\,\tilde \xi]}_\tau\,$. 
The {\em Fermi coordinates} (see \cite{Chav}, page 138 ff) send $x\in{\rm [Tub\,\tilde \xi]}_\tau\setminus \tilde \xi$  on  
$(t,w)\in[0,\tau]\times {\bf N}^1(\tilde\xi)\,$. Actually, $t\!=\!d_g(x,\tilde \xi)$ and $y\in\tilde \xi$ is the foot of the shortest $g$-geodesic $t\in[0,\tau]\rightarrow\exp_yt\,w\in M$ from $\tilde \xi$ to $x\!=\!\exp_ytw$ and $w\in{\bf N}^1(\tilde\xi)\,$. In the Fermi coordinates at $x\!=\!\exp_y(t\,w)$, one can decompose the volume element in $([{\rm Tub}\,\tilde \xi]_\tau,g)$ (the related Jacobian determinant is $\grave\psi\!\in\!C^\infty([0,\tau]\!\times\!{\bf N}^1(\tilde\xi),\R)$, $dw$ is the $(n-1)$-volume element induced on ${\bf N}^1(\tilde\xi)$ by $(M,g)$) as follows
\begin{multline}\label{CMS0} dx=\grave\psi(t,w)\ dt\wedge dw \  \  \hbox{with}\  \  \frac{\grave\psi(t,w)}{t}\rightarrow 1\  \  \hbox{as}\  \  t\rightarrow 0\ ,\\
\hbox{where}\  \  dw:=dy\wedge dw_y\  \  \hbox{and}\  \ dy\,,\,dw_y\  \hbox{are volume elements on}\  \tilde\xi\,,\ {\bf N}^1(T_y\tilde\xi)\ ,
\end{multline} 
since the bundle ${\bf N}^1(\tilde\xi)\rightarrow \tilde\xi$ is here a natural Riemannian submersion.
\par We use the next lemma to describe squares $G\in\bar G$ near to a bone $\xi$. Set ${\bf D}_y^\tau(\tilde\xi)={\bf D}^\tau(\tilde\xi)\cap T_yM$ (see ({\rm\ref{diskbdle}})).
\begin{lem} \label{CMS-5} Given a bone $\xi$, one has a {\it bundle-correspondence} 
\begin{equation*}
\hbox{\bf St}_2(T{\rm [Tub\,\tilde \xi]}_\tau)\simeq \cup_{y\in\tilde\xi}\ {\bf D}_y^\tau(\tilde\xi)\times \hbox{\bf St}_2(T_yM)\  .
\end{equation*}
Define $\Delta_\tau^{\xi,y}:=\exp_y({\bf D}_y^\tau(\tilde\xi))$.
\end{lem}
\begin{proof} This relies on the above Fermi coordinates: indeed, the disk $\Delta_\tau^{\xi,y}$ is described by points $x=\exp(tw)\equiv(t,w)$ (with $tw\!\in\!{\bf D}_y^\tau(\tilde\xi)$) having closest point $y\in\tilde \xi\setminus\partial\tilde\xi$ while each ${\bf St}_2(T_{(t,w)}M)$ is identified to ${\bf St}_2(T_yM)$ by parallel translation along $t\mapsto\exp(tw)$. 
\end{proof}
The volume element $d{\mathscr G}$ of definition \ref{elemvolume} reads on $\hbox{\bf St}_2(T[{\rm Tub}\,\tilde \xi]_\tau)$
\begin{equation}\label{CMS2} d{\mathscr G}=\check\psi(x)\,dx\wedge d{\mathscr G}_x\  ,
\end{equation}
where $dx$ is the volume element on $(M,g)$, $d{\mathscr G}_x\!=\!d{\mathscr G}_{\hbox{\tiny\bf St}_2(T_{x}M)}\,$, where 
$\check\psi\!\in\!C^\infty([{\rm Tub}\,\tilde \xi]_\tau,\R)$ and $\check\psi(x)\rightarrow1$ uniformly in $x\in {\rm [Tub\,\tilde \xi]}_{\tau'}$ as
$\tau'\rightarrow0\,$.
\begin{defn}\label{model1111} Relying on lemma {\rm\ref{CMS-5}}, given $\bar G$, a bone $\xi,\,y\!\in\!\tilde\xi$, we sort out the various $G\!\in\!\bar G$, defining (lemma {\rm\ref{CMS1000}}) with $S\!\subset\!\tilde\xi\setminus\partial\tilde\xi$
\begin{gather*}\bar G_r^{\xi,y}:=\{G\in \bar G_r\mid x=G(0,0)\!\in\! [{\rm Tub}\,\tilde \xi]_\tau\ \ \hbox{and}\ \ d(x,y)=d(x,\tilde\xi)\}\  ,
\\ \bar G_r^S:=\cup_{y\in S}\,\bar G_r^{\xi,y}
\ \ \hbox{and}\ \ \bar G_r^{\cap\xi}:=\{G\in \bar G_r\mid G\cap\hat\xi\not=\emptyset\}\ .
\end{gather*}
\end{defn}
\begin{rem} \label{notmissing} Observe $\bar G_r^{\cap\xi}\!\not\subset\!{\bar G}_r^{\hat\xi}\,$: all $G\!\in \!\bar G_r^{\cap\xi}$ are not given in Fermi coordinates {\it around}  $\hat\xi$ but {\it around $\tilde\xi$ if $x$ belongs to $[{\rm Tub}\,\tilde \xi]_\tau$} and $r$ small (lemma {\rm\ref{compactsimpl}}).
This ``lack'' later appears not to harm our argument.
\end{rem}
\begin{defn} \label{CMS-1000} Given $r_1\!\in]0,\tau/2]$, the {\em pseudo-homothety} ${\mathcal H}_{y,{\bf r}}$ of center $y\!\in\!\tilde\xi$ and ratio ${\bf r}\!\in]0,1]$ sending  $r_1$ into $r_2:={\bf r}r_1$-squares is 
\begin{equation*}{\mathcal H}_{y,{\bf r}}: G_1\!=\!(x_1,u_1,v_1)\in \bar G_{r_1}^{\xi,y}\mapsto G_2\!=\!(x_2,u_2,v_2)\in \bar G_{r_2}^{\xi,y}
\end{equation*} 
where one has 
$x_1\!=\!(d(x_1,y),w), x_2\!=\!(d(x_2,y),w)=({\bf r}d(x_1,y),w)$ and $(u_1,v_1),(u_2,v_2)$ are orthonormal $2$-systems of vectors matching under parallel translation along $\,t\mapsto(t,w)\equiv\exp(tw)$ from $x_1$ to $x_2\,$.
The measure space $(\bar G_{r_1,x_1}^{\xi,y},d{\mathscr G}_{x_1})$ is sent by this mapping onto $(\bar G_{r_2,x_2}^{\xi,y},d{\mathscr G}_{x_2})$, where $\bar G_{r_i,x_i}^{\xi,y}$ is the set of $G\in \bar G_{r_i}^{\xi,y}$ such that $G(0,0)=x_i$ for $i=1,2\,$.
\par Actually, one defines new coordinates $(s,w)$ on $\bar G_r^{\xi,y}$, setting
$(s,w)\!:=\!(d(x,y)/r,w)$, where $x_1,x_2$ also read
$x_1\!=\!(r_1s,w), x_2\!=\!(r_2s,w)$ with $s\!=\!d(x_1,y)/r_1\!=\!d(x_2,y)/r_2$.
\end{defn} 
Choose $\tau>0$ according to corollary \ref{compactsimpl333}.
\begin{lem} \label{newvariable} Given $G\!\in\!\bar G_r$ of side $r\!=\!r_G$, define $s_G\!:=\!t_G/r_G$, where $t\!=\!t_G\!=\!d(x,\tilde \xi)$ with $x\!=\!G(0,0)$, so $\,s_{{\mathcal H}_{y,{\bf r}}(G)}\!=\!s_G$ as ${\bf r}\!>\!0$ varies. If $\bar G$ is close to $\underline G$ and $G\!\in\!\bar G_r$ cuts $\hat\xi$, one has $s\!\leq\! 2,\,d(x,\tilde \xi)\!\leq\! d(x,\hat \xi)\!\leq\! 2r$ and $(x,u,v)\!\equiv\!(t,w,u,v)$ belongs to $\hbox{\bf St}_2(T[{\rm Tub}\,\tilde \xi]_\tau)$ if $2r\!\leq\!\tau$. 
\end{lem}
\begin{proof} If $r$ is small enough, $G\!\in\!\underline G_r$ cuts $\hat\xi$ only if $d(\hat \xi,x)\!\leq\!\sqrt{2}r$, since the furthest point from $x$ in this square ends the diagonal starting from $x$. For $\bar G$ (uniformly) close enough to $\underline G$, no $G\!\in\! \bar G_r$ meeting $\hat\xi$ is left out if $d(\tilde \xi,x)\leq2r$, since $\hat\xi$ is uniformly small in $\tilde\xi$ (lemma \ref{compactsimpl}).
\end{proof}
\begin{nota} If $G\!\in\!\bar G_r^{\xi,y}$ meets $\hat\xi$, it belongs to (definition \ref{model1111})
\begin{equation}\label{CMS-33}[\bar G^{\xi,y}]_r:=\{G\in \bar G_r^{\xi,y}\mid d(\tilde \xi,x)\leq2r\}\subset \bar G_r^{\xi,y}\,.
\end{equation}
\par By definition {\rm\ref{CMS10333}}, lemma \ref{CMS-5}, all $G\!\in\![\bar G^{\xi,y}]_r$ read
\begin{equation}\label{CMS-3}
G=G_{r,x,u,v}\equiv(r,x,u,v)\equiv(r,s,w,u,v)\,,
\end{equation}
with $(r,s)\in]0,\tau/2]\times[0,2]$, $w\in{\bf N}^1(\tilde\xi)$, where $y$ is the nearest point on $\tilde \xi$ to $x=x_G=G(0,0)\equiv(s,w)$ and $(u,v)$ is seen in ${\bf St}_2(T_yM)$ after a parallel translation along the $g$-ray from $x$ to $y\,$. 
\end{nota}

\begin{lem}\label{n'enjetezplus1} Given a bone $\xi$, given $r\!>\!0$, one has at $G\!\in\!\bar G_r^{\xi,y}\!\subset\!\bar G_r$ (see {\rm(\ref{CMS0})}, {\rm (\ref{CMS2})})
\begin{equation}
\label{n'enjetezplus11} (d{\mathscr G}_{\mid\bar G_r^{\xi,y}})_G=
\psi(G)\ dt\wedge dw_y\wedge d{\mathscr G}_x\   \  \hbox{with}\  \  \psi(G)=\grave\psi(x_G)\,\check \psi(x_G)\  .
\end{equation}
If $x=y\in\tilde \xi$, one has $d{\mathscr G}_y=d{\mathscr G}_{\mid{\bf St}_2(T_yM)}\equiv d{\mathpzc G}_{\mid{\bf St}_2(T_{0_y}T_yM)}$ (in $(T_yM,g_y)$).
\par
Sending diffeomorphically $[\bar G^{\xi,y}]_r$ onto $[\bar G^{\xi,y}]_{{\bf r}\,r}$, the {\em pseudo-homothety} ${\mathcal H}_{y,{\bf r}}$ of definition {\rm\ref{CMS-1000}} verifies
\begin{equation} \label{n'enjetezplus2}  
({\mathcal H}_{y,{\bf r}}^\ast\,{\frac{d{\mathscr G}_{\mid\bar G_{{\bf r}r}^{\xi,y}}}{({\bf r}r)^2}})_G=\theta_{\bf r}(G)\,\frac{(d{\mathscr G}_{\mid\bar G_r^{\xi,y}})_G}{r^2}\ \  \hbox{with}\  \ \theta_{\bf r}(G)=
\frac{\psi({\mathcal H}_{y,{\bf r}}(G))}{{\bf r}\psi(G)}\ ,
\end{equation}
and $\theta_{\bf r}(G)$ extends $C^\infty$ in $({\bf r},G)$ to $[0,1]\!\times\!\bar G_r^{\xi,y}$.
Set ${\bf X}\!=\!T_y\tilde\xi\,$. 
\par {\em Define} $\Omega^\xi\!:=\!sds\wedge dw_y\wedge d{\mathscr G}_{x}$ on $[{\rm Tub}\tilde \xi]_{\tau}$;
if $G\!\in\!\bar G_r^{\xi,y}$, then (see {\rm(\ref{CMS8bise33-1})})
\begin{equation}\label{n'enjetezplus3} {\mathcal H}_{y,{\bf r}}^\ast\,{\Omega^\xi}={\Omega^\xi}\ \ \hbox{and}\ \  
\lim_{{\bf r}\rightarrow0}{\mathcal H}_{y,{\bf r}}^\ast\,\frac{d{\mathscr G}_{\mid\bar G_r^{\xi,y}}}{({\bf r}r)^2}=\frac{d{\mathpzc G}_{\mid{\bf G}_r^{{\bf X},0_y}}}{r^2}={\Omega^\xi}\ .
\end{equation}
\end{lem}
\begin{proof} 
Along a geodesic ray orthogonal in $y$ to $\tilde \xi$, in the coordinates ({\rm\ref{CMS-3}}), the pseudo-homothety ${\mathcal H}_{y,{\bf r}}$ (definition \ref{CMS-1000}) sends $\{r\}\times[\bar G^{\xi,y}]_r$ onto
$\{{\bf r}\,r\}\times[\bar G^{\xi,y}]_{{\bf r}\,r}$ along the rule
\begin{equation} \label{magnifique1} {\mathcal H}_{y,{\bf r}} : (r,s,w,u,v)\mapsto ({\bf r}\,r,s,w,u,v)\  ,
\end{equation}
or else
\begin{equation} \label{magnifique2} 
{\mathcal H}_{y,{\bf r}} : G_{r,x,u,v}\!\in\![\bar G^{\xi,y}]_r\mapsto G_{{\bf r}r,x_{\bf r},u,v}\!\in\![\bar G^{\xi,y}]_{{\bf r}\,r}\ ,
\end{equation} 
where ${x_{\bf r}}$ is at distance ${\bf r}\,t$ from $y$ if $x$ is at distance $t\,$, and $d{\mathscr G}_{x}$ (see ({\rm \ref{CMS2}})) is sent to
$d{\mathscr G}_{x_{\bf r}}$. 
Compute from (\ref{n'enjetezplus11})
\begin{equation}\label{n'enjetezplus4}(\frac{d{\mathscr G}_{\mid\bar G_r^{\xi,y}}}{r^2})_{G}=\frac{\psi(G)}{t_G} \ \frac{t_G}{r}\,d(\!\frac{t_G}{r}\!)\wedge dw_y\wedge d{\mathscr G}_{x_G}=\frac{\psi(G)}{t_G} \,\Omega^\xi_G\  , 
\end{equation} 
so, by pull-back, for any $G\in[\bar G^{\xi,y}]_r\,,\,G'={\mathcal H}_{y,{\bf r}}(G)\in{\mathcal H}_{y,{\bf r}}([\bar G^{\xi,y}]_r)$, get
\begin{equation*} 
({\mathcal H}_{y,{\bf r}}^\ast\,\frac{d{\mathscr G}_{\mid\bar G_{{\bf r}r}^{\xi,y}}}{({\bf r}r)^2})_{G}=\frac{\psi(G')}{t_{G'}}\,\Omega^\xi_{G}=
\frac{\psi(G')}{{\bf r}\psi(G)}\,(\frac{d{\mathscr G}_{\mid\bar G_r^{\xi,y}}}{r^2})_{G}\  ,
\end{equation*}
proving (\ref{n'enjetezplus2}). 
As for (\ref{n'enjetezplus3}), observe, thanks to (\ref{CMS0}), (\ref{CMS2}), that
\begin{equation} \label{n'enjetezplus33}\frac{\psi({\mathcal H}_{y,{\bf r}}(G))}{{\bf r}t}=
\frac{\grave\psi({\bf r}\,t,w)}{{\bf r}\,t}\ \check \psi({\mathcal H}_{y,{\bf r}}(G))\rightarrow1\  \ \hbox{as} \ \ {\bf r}\rightarrow0\ ,
\end{equation} and use (\ref{n'enjetezplus4}). This implies that $\theta_{\bf r}(G)$ has a finite limit as ${\bf r}\rightarrow0$. More, it follows from (\ref{n'enjetezplus33}) that the $C^\infty$ function $\psi({\mathcal H}_{y,{\bf r}}(G))$ in $({\bf r},G)$ vanishes for ${\bf r}=0$, which proves that $\theta_{\bf r}(G)$ has a $C^\infty$ extension in $({\bf r},G)$ to $[0,1]\!\times\! \bar G_r^{\xi,y}$ (use the Taylor formula with integral remainder). 
\end{proof}

\subsection{Proving part $(i)$ of proposition \ref{CMS-7}}\label{part(i)} $ $
\begin{nota}
Given $\tau\!>\!0$ (corollary {\rm\ref{compactsimpl333}}) and $r\!\in]0,\tau/2]$ (lemma {\rm\ref{newvariable}}),  a bone $\xi$, $\bar G$ generic close to $\underline{G}$ on $D$, recall (lemmas {\rm\ref{Reggeparallt}}, {\rm\ref{Reggecomplet}}, definitions {\rm\ref{Riemannparallt}}, {\rm\ref{restrindex}}) the $g, g_0$-parallel translations $\mathcal{P}_{\tilde A_z^{-1}},\,\mathcal{P}^0_{\tilde A_z^{-1}}$ paired with $z\!\in\! {\rm rge}(G)\cap \hat\xi,\,x\!=\!G(0,0)$ and $\tilde A_z^{-1}$.  Define the bundle-morphisms $f^\xi,\ \mathfrak{f}{\,}^\xi$ from the fibre ${\bf St}_2(T_xD)\!\equiv\!\bar G_r$ to the fibre $\bigwedge^2(T^\ast_xD)$
\begin{gather}\label{CMS4} \mathfrak{f}{\,}^\xi: G\longmapsto  \mathfrak{f}{\,}^\xi(G)\!=\!\!\!\!\!\sum_{z\in{\rm rge}(G)\cap \hat\xi}\!\!\!\!\!\!
{\rm index}_{G,\xi,z}\,\langle{\ast_g \,\mathcal{P}_{\tilde A_z^{-1}}(\hat\xi_z)}, \cdot\rangle_g
\ ,\\
\label{CMS4bissee} f^\xi\!:\! G\mapsto f^\xi(G)\!=\!\!\!\!\!\!\!\sum_{z\in{\rm rge}(G)\cap \hat\xi}\!\!\!\!\!\!\!
{\rm index}_{G,\xi,z} \langle dT(T^{-1}(x))({\ast_{g_0} \,\mathcal{P}^0_{\tilde A_z^{\!-\!1}}(\xi)}), \cdot\rangle_{g^0}
\,,
\end{gather}
and $f^\xi(G)=\mathfrak{f}{\,}^\xi(G)=0$ if $G\!\in\!\bar G_r$ has empty intersection with $\hat\xi\,$.
\end{nota}
\begin{lem} \label{secti1} Given $r\!\in]0,\tau/2],\,\bar G_r$, the supports of $\mathfrak{f}{\,}^\xi,f^\xi$ defined on $\bar G_r$ are subsets of $\cup_{y\in\tilde \xi}[\bar G^{\xi,y}]_r\!:=\!\{G\in\bar G_r\mid x\!=\!G(0,0)\in[{\rm Tub}\,\tilde \xi]_{2r}\}\,$. Restricting $\mathfrak{f}{\,}^\xi,f^\xi$ to $\bar G_r^\rho\,$ (definition {\rm \ref{CMS-7777}}), those supports are $\subset\cup_{y\in\hat \xi}[\bar G^{\xi,y}]_r$.
\end{lem}
\begin{proof} This follows from lemma {\rm\ref{newvariable}} and ({\rm\ref{CMS-33}}).
\end{proof}
\begin{lem}\label{auxili1111} Given $G\!\in\!\bar G$, identifying $\bigwedge^2(T^\ast_xD)$ with $\bigwedge^2(T_xD)$ through $g_x$-duality, one has (see {\rm(\ref{Reggemeas2})}, definition {\rm\ref{Reggemeas2222}})
\begin{equation*}\mathscr{R}_0(G)=\sum_{\xi\in \dagger\!K_{n-2}(D)}  \!\!\!\alpha_\xi\,f^\xi(G)
\ \ \hbox{and}\ \ \mathpzc{R}{}(G)=\sum_{\xi\in \dagger\!K_{n-2}(D)} \!\!\!\alpha_\xi\,\mathfrak{f}{\,}^\xi(G)\ .
\end{equation*}
\end{lem}
\begin{proof} Indeed, (\ref{CMS4bissee}) is closely related to the ``Regge curving'' ${\mathscr R}_0(G)$, see definition \ref{Reggemeas}, (\ref{Reggemeas2}) and lemma \ref{Reggecomplet}, while (\ref{CMS4}) is related to the ``ghost Regge curving'' $\mathpzc{R}{}(G)$ defined in (\ref{ghost1}).
\end{proof}

\begin{nota} Setting as before $T_xG\!:=\!(\partial G/\partial s\wedge\partial G/\partial t)(0,0)$, define (recall $\Vert T_xG\Vert_g=1$ and $T^\ast g^0=g_0$)
\begin{gather}\label{courpleine}[[\mathfrak{f}{\,}^\xi(G)]]\!=\!\mathfrak{f}{\,}^\xi(G)({T_xG}),\ \ \ 
[[f^\xi(G)]]\!=\!f^\xi(G)(\frac{T_xG}{\Vert T_xG\Vert_{g^0}})\, .
\end{gather}
With these notations, on may rewrite (\ref{CMS9}) and (\ref{ghost2}), for instance
\begin{equation}\label{CMS9-3333}I^{g_0}(\bar G_r)=\sum_{\xi\in\dagger K_{n-2}(D)}\!\!\alpha_\xi\,\int_{G\in\bar G_r}\,[[f^\xi(G)]]\,\frac{d {\mathscr G}}{r^2}\  .
\end{equation}
\end{nota}
\begin{rem}\label{muborne}Recall from theorem {\rm\ref{text-1133}} (see remark {\rm\ref{bounded222}})
that, for any $\rho\in]0,\rho_{\bar G}]$ and $K$ of theorem {\rm\ref{theorem I}} paired with $\rho$, the total number ({\it counting multiplicities}) of intersections of any bone $\hat\xi$ for $\xi\in K$ with any $G\!\in\!\bar G$ (with $\bar G\!\in\!{\mathcal V}_c$) is bounded by a universal  $\mu=m_{n-2}(\!(n\!-\!1)n\!+\!c)$.
\end{rem}
\begin{lem}\label{CMS5} For $\tau\!>\!0, \,r\!\in]0,\tau/2],\,G\!\in\!\bar G_r$, all $f{\,}^\xi, \mathfrak{f}{\,}^\xi$ verify
\begin{equation*} \vert \,[[f^\xi(G)]] \,\vert\ \ \hbox{and}\ \ \vert \,[[\mathfrak{f}{\,}^\xi(G)]] \,\vert\leq\!\!\!\!\sum_{z\in{\rm rge}(G)\cap \hat\xi}\vert\,\hbox{\rm index}_{G,\xi,z}\,\vert\leq \mu\ 1_{[{\rm Tub}\,\tilde \xi]_{2r}}(x_G)\  ,
\end{equation*}
with $\mu$ as in remark {\rm\ref{muborne}}.
\end{lem}
\begin{proof} By definition and properties of the index (sections \ref{sssec:num1}, \ref{sssec:num2}), one has (using lemma \ref{secti1})
\begin{equation*} \vert \,[[\mathfrak{f}{\,}^\xi(G)]] \,\vert\!\leq\!\!\!\!\sum_{z\in{\rm rge}(G)\cap \hat\xi} \!\!\!\vert\hbox{index}_{G,\xi,z}\,\langle {\ast_g \mathcal{P}_{\tilde A_z^{-1}}\hat \xi_z},\frac{T_xG}{\Vert T_xG\Vert_g}\rangle_g\,\vert\!\leq\! \mu\,1_{[{\rm Tub}\,\tilde \xi]_{2r}}\!(x_G)\,,
\end{equation*}
and the corresponding holds for $f^\xi$.
\end{proof}
\begin{lem} \label{integrals1111}
Thanks to lemma {\rm\ref{CMS5}}, define the integrals (see {\rm(\ref{courpleine})})
\begin{gather}\label{CMS8bis1} I^{\xi,y}(\bar G_r)\!:=\!\!\int_{G\in[\bar G^{\xi,y}]_r} [[f^\xi(G)]]\,
\frac{d{\mathscr G}_{\mid \bar G^{\xi,y}_r}}{r^2}\, ,\\ \label{CMS8bis2}
 {J}^{\xi,y}(\bar G_r)\!:=\!\!\int_{G\in[\bar G^{\xi,y}]_r} [[\mathfrak{f}{\,}^\xi(G)]] \,\frac{d{\mathscr G}_{\mid \bar G^{\xi,y}_r}}{r^2}\,.
\end{gather}
One has (definition {\rm \ref{CMS-7777}} for $\bar G_r^\rho$, {\rm(\ref{CMS2})}, {\rm(\ref{CMS0})} for the splitting of $d{\mathscr G}$)
\begin{equation*}I^{g_0}(\bar G_r^\rho)=\!\!\!\!\sum_{\xi\in\dagger\! K_{n-2}(D)}\!\!\!\alpha_\xi\int_{\hat\xi}I^{\xi,y}(\bar G_r)\,dy\ ,\ {J}^{g_0}(\bar G_r^\rho)=\!\!\!\!\sum_{\xi\in \dagger\!K_{n-2}(D)}\!\!\!\alpha_\xi
\int_{\hat\xi}{J}^{\xi,y}(\bar G_r)\,dy\  .
\end{equation*}
\end{lem}
\begin{proof}
Use (\ref{CMS9}), (\ref{ghost2}), (\ref{CMS8bis1}), (\ref{CMS8bis2}) and lemma \ref{auxili1111}.
\end{proof}

Recall that, the mesh $\rho\!>\!0$ is related, through definition \ref{linkint}, to an integer $E(\rho)$ which defines an order of regular subdivision needed (see lemma \ref{2a} (\ref{2a-bis33})) in the establishment of theorem \ref{theorem I}.
\begin{lem}\label{paliprox} There exists ${\mathfrak C}_0$ depending on $D\!\subset\!(M,g)$ satisfying for any bone $\xi\!\in \!K_E$, $E$ large enough and $r\!\in\![0,\tau/2]$, for $\bar G_r$ close to $\underline G_r$
\begin{equation*} \label{uniformity!2}\vert I^{\xi,y}(\bar G_r)-{J}^{\xi,y}(\bar G_r)\vert\leq {\mathfrak C}_0\,\rho\  .
\end{equation*}
\end{lem}
\begin{proof} Making an explicit mention of $(\mathcal{P}^0_{\tilde  A_z})^{-1}=\mathcal{P}^0_{\tilde A_z^{-1}}$, one has
\begin{multline*} \label{paliprox1111}\cos\angle_{g_0}(\ast_{g_0} \,\mathcal{P}^0_{\tilde A_z^{-1}}({\xi}),dT^{-1}(x)(T_xG))=\\=
\langle dT(T^{-1}(x))({\ast_{g_0} \,\mathcal{P}^0_{\tilde A_z^{-1}}(\xi)}), \frac{T_xG}{\Vert T_xG\Vert_{g^0}}\rangle_{g^0}\ .
\end{multline*}
But, for $\omega\!\in\!\bigwedge^2(T_xD)$ with $x$ in the interior of an embedded $n$-simplex
\begin{equation}\label{paliprox3333}\ast_{g_0} \,dT^{-1}(x) (\omega)=dT^{-1}(x) \ast_{g^0} \,(\omega)\  .
\end{equation}
Link the unit $(n\!-\!2)$-vector $\xi$ at $T^{-1}(z)$ with $T_z\hat\xi$ to get (use (\ref{paliprox3333}))
\begin{multline*}\label{paliprox2222}\ast_{g_0} \,\mathcal{P}^0_{\tilde A_z^{-1}}(\xi)=\ast_{g_0} \,\mathcal{P}^0_{\tilde A_z^{-1}}(dT^{-1}(z)(T_z\hat\xi))=\\=\ast_{g_0} \, (dT^{-1}(x) \circ dT(T^{-1}(x))\circ\mathcal{P}^0_{\tilde A_z^{-1}}\circ dT^{-1}(z)(T_z\hat\xi))=\\=
dT^{-1}(x)\,\ast_{g^0} \,(dT(T^{-1}(x)) \circ\mathcal{P}^0_{\tilde A_z^{-1}}\circ dT^{-1}(z)(T_z\hat\xi))\  .
\end{multline*}
For $E$ chosen large enough (theorem \ref{theorem I}, lemma  \ref{Reggeparallt}, definition \ref{Riemannparallt}), 
\par\vskip2mm 
{\em in view of} what was said from the beginning of this proof,
\par {\em in view of} the $g$-proximity of the metrics $g$ and $g_0$ {\em via} $T$ (by $(ii)$ $(b)$ of definition \ref{convpolyhedra} defining a polyhedral approximation),
\par {\em in view of} the resulting $g$-proximity of $\ast_g$ and $\ast_{g^0}$, 
\par {\em in view of} the $g$-proximity of the $g$ and $g^0$-parallel translations along the curve $\tilde A_z^{-1}$ from $z\in\hbox{rge}(G)\cap\hat\xi$ to $x=G(0,0)$, in fact $\mathcal{P}_{\tilde A_z^{-1}}$ and $dT(T^{-1}(x))\circ \mathcal{P}^0_{\tilde A_z^{-1}}\circ 
dT^{-1}(z)$ (by definition \ref{convpolyhedra} $(iii)$ $(c)$, because all paths $\tilde A_z$ have uniformly bounded length in terms of the side $r$ of $G$ and of the geometry of $D\subset(M,g)$), 
\vskip2mm\par\noindent 
{\em one finds} 
${\mathfrak C}_0'\!>\!0$ depending on $D\!\subset(M,g)$ such that, for $z\!\in\! {\rm rge}(G)\cap \hat\xi$
\begin{equation*} \label{uniformity!}\vert\cos \angle_{g}\! (\ast_{g} \mathcal{P}_{\!\!\tilde A_z^{\!-\!1}}\! (T_z\hat\xi), T_xG)\!-\!\cos\angle_{g_0}\!(\ast_{g_0}\mathcal{P}^0_{\!\!\tilde A_z^{\!-\!1}}\!(\xi), dT^{\!-\!1}\!(\!x\!)(\!T_xG\!)\!)\vert\!\leq\! 
{\mathfrak C}_0'\rho\,.
\end{equation*}
For any fixed $i$ (notations of lemma \ref{compactsimpl}, remark \ref{compact1099}, (\ref{bone111133}), corollary {\rm\ref{compactsimpl333}}, notation \ref{nortrunc}), define
$\widetilde{{\mathfrak S}_i}\!\subset\!{\mathfrak S}_i\!\times\!M$ (fibred over ${\mathfrak S}_i$)
\begin{equation*}\widetilde{{\mathfrak S}_i}:= \coprod_{\hat\eta\in{\mathfrak S}_i} {}\{\hat\eta\}\!\times\![{\rm Tub}\,\tilde \eta]_{\tau}\ , \ \ \hbox{where}\ \ \tilde \eta \ \ \hbox{stands for}\ \ \tilde \eta_i \ \ \hbox{(assumption \ref{unifsmall1})}
\end{equation*}
and its compact closure $\widetilde{\overline{\mathfrak S}_i}=\overline{\widetilde{\mathfrak S}_i}\subset {\mathpzc F}^{n\!-\!2}\times M$ (where ${\mathpzc F}^{n\!-\!2}$ is a manifold, (\ref{bone111133})).
For fixed $i$ and any $\hat\eta\!\in\!{\mathfrak S}_i$, $(d {\mathscr G}_{\mid \bar G^{\eta,y}_r})_G/r^2$ differs from $\Omega^\eta_{G}\!=\!sds\wedge dw_y\wedge d{\mathscr G}_{x_G}$ by a {\it continuous density} in $(\hat\eta,x)\!\in\!\overline{\tilde{\mathfrak S}_i}$ (lemma \ref{n'enjetezplus1}, (\ref{n'enjetezplus4})), {\it uniformly bounded} on $\overline{\tilde{\mathfrak S}_i}\,$. The set of $i$ is finite: an integration over $\tilde\eta\!=\!\tilde\xi$ (use lemma \ref{CMS5}) gives the result.
\end{proof}
\subsubsection{Proving proposition {\rm\ref{CMS-7} (i)}}$ $$ $
\begin{proof} Lemma \ref{paliprox} gives (for $\rho>0$ small enough, for any $i$, uniformly in $r\in[0,\tau/2]$ and $(\hat\xi,y)\in\widetilde{\mathfrak S}_i$) the inequality
$ \vert I^{\xi,y}(\bar G_r)-{J}^{\xi,y}(\bar G_r)\vert\leq 
{\mathfrak C}_0\,\rho\  ,
$
and, thanks to lemmas \ref{integrals1111}, \ref{bornangles}, one gets by summing up over all $\xi\in K_{n-2}(D)$ (convention \ref{bone1066})
\begin{equation*} \label{uniformity!2bis}\vert I^{g_0}(\bar G_r^\rho)-{J}{}^{g_0}(\bar G_r^\rho)\vert\leq {\mathfrak A}\,{\mathfrak C}_0\,\rho\  .
\end{equation*}
Setting ${\mathfrak C}={\mathfrak A}\,{\mathfrak C}_0$ gives $(i)$.
\end{proof}

\subsection{Proving part $(ii)$ of proposition \ref{CMS-7}} \label{part(ii)}
$ $

With $\rho>0$, the approximation $T,K$ is chosen and fixed.
\subsubsection{Neglecting negligible contributions}\label{negligeons}$ $

Here, we come back to a problem mentioned earlier (remark \ref{notmissing}) and prove the first convergence stated in 
proposition \ref{CMS-7} $(ii)$. 
Yet we begin with a useful new observation which is a nice feature of polyhedral approximations (see section \ref{6.1} for notions and notations).
\begin{prop}\label{-12-CMS} There exist constants ${\mathcal C}_1,{\mathcal C}_2$ (one has ${\mathcal C}_2\!\geq\!1$) depending on $T,K,D\subset(M,g)\,$ such that, 
given $E$ (thus $\rho$), if $r\leq{\mathcal C}_1\,\rho$ and $x$ is at $g$-distance $\leq r$ from each of two simplices (of any dimensions) $\hat\eta_1$ and $\hat\eta_2$ ($\eta_1, \eta_2$ in $K_E$), then $\eta_1\cap\eta_2\not=\emptyset$ and $d_g(x,\hat\eta_1\cap\hat\eta_2)\leq {\mathcal C}_2\,r$.
\end{prop}
\begin{proof} See appendix \ref{nicefeature}.
\end{proof}
\begin{cor}\label{-22-CMS} There exist constants ${\bf C}_1,{\bf C}_2$ depending on $T,K,D\subset(M,g)\,$ such that, 
given any $E$ (thus $\rho$), if $r\leq{\bf C}_1\,\rho$ and $x$ is at $g$-distance $\leq r$ from $k$ simplices $\hat\eta_1,\dots,\hat\eta_k$ (of any dimensions) with $\eta_1,\dots,\eta_k$ in $K_E$, then $\hat\eta_1\cap\cdots\cap\hat\eta_k\!\not=\!\emptyset$ and 
$d_g(x,\hat\eta_1\cap\cdots\cap\hat\eta_k)$ is $\leq {\bf C}_2\,r$. 
\end{cor}
\begin{proof} This follows by induction on $k$ and proposition \ref{bds4}.
\end{proof}
\begin{rem} 
If in a piecewise flat $n$-dimensional polyhedron $K$ exists, for some $\epsilon\!>\!0$, an integer $2{\mathcal B}_0$ bounding, for any $p\!\in\!K$, the number of $n$-simplices $\sigma\!\in\!K$ such that
$B(p,\epsilon)\cap\sigma\not\!=\!\emptyset$, then proposition \ref{-12-CMS} and corollary \ref{-22-CMS} are true in $K$ with constants ${\mathcal C}_1,{\mathcal C}_2,{\bf C}_1,{\bf C}_2$ depending on $n$, on a lower bound $t_0>0$ of the thickness for all the $n$-simplices of $K$ and on ${\mathcal B}_0$ (in place of ${\mathcal B}$), see remark \ref{ohnon}.
\end{rem}

\par
By definition \ref{CMS-7777}, if a square $G$ belongs to $\bar G_r\!\setminus\!\bar G_r^\rho$, there exists a bone $\xi_0$ such that $\hat\xi_0\cap G\!\not=\!\emptyset$ and the closest point to $x_G\!=\!G(0,0)$ in $\hat\xi_0$ is in $\partial \hat\xi_0$. We want to put a hand on $I^{g_0}(\bar G_r\!\setminus\!\bar G_r^\rho)\!=\!\sum_{\xi\in \dagger K_{n-2}(D)}\alpha_\xi \int_{\bar G_r\!\setminus\!\bar G_r^\rho}[[f^\xi]]\,{d{\mathscr G}}/{r^2}$ (see (\ref{CMS9-3333})). For dimensional reasons analysed below, the $G_r$ that cut a given $\hat\xi$ with $x_G$ verifying $d(x_G,\partial\hat\xi)\!\leq \!{\rm Cst}\,r$ (with ${\rm Cst}\!>\!0$ fixed) build a set $\Omega_\xi$ for which $I^\xi(\Omega_\xi)\!=\!\int_{\Omega_\xi}[[f^\xi]]{d{\mathscr G}}/{r^2}$ tends with $r$ to $0$. We thus focus on the squares $G_r$ cutting bones $\hat\xi_0$ and $\hat\xi$, for which $x_G$ has closest point to $\hat\xi_0$ in $\partial \hat\xi_0$ and to $\hat\xi$ in
$\hat\xi\setminus\partial\hat\xi$, giving a definition and proving a lemma.
 \begin{defn}\label{Fermi-100} Given a bone $\xi\!\in\!\dagger\!K_{n-2}(D)$, define 
\begin{equation*}  \backslash \!\!\!\!\bar G_r^{\cap\xi}\!=\!\{G\in \bar G_r^{\cap\xi}\mid 
d(x_G,\partial\hat\xi)\!\leq\!2{\mathcal C}_2r\}\subset\bar G_r^{\cap\xi}\!=\!\{G\in \bar G_r\mid G\cap\hat\xi\not=\emptyset\}\,.
\end{equation*}
\end{defn}
\begin{lem}\label{blingbling}
 For $r\!>\!0$ small enough (one assumes in particular $2r\!\leq\!{\mathcal C}_1\rho$) and for any bone $\xi$,
the set $\bar G_r^{\cap \xi}\setminus \bar G_r^\rho$ is included in 
$\,\backslash \!\!\!\!\bar G_r^{\cap\xi}$.
Therefore, one has $\bar G_r\!\setminus \!\bar G_r^\rho=\!\cup_{\xi\in \dagger K_{n-2}(D)}(\bar G_r^{\cap\xi}\!\setminus \!\bar G_r^\rho)\subset\cup_{\xi\in \dagger K_{n-2}(D)}\,\backslash \!\!\!\!\bar G_r^{\cap\xi}$.
\end{lem}
\begin{proof} For any $G\!\in\!\bar G_r\setminus \bar G_r^\rho$, there exists (definition \ref{CMS-7777}) a bone $\xi_0$ such that ${\rm rge}(G)\cap\xi_0\not=\emptyset$ and $d(x_G,\hat\xi_0)\!=\!d(x_G,\partial\hat\xi_0)\!\leq\! 2r$. But for any bone $\xi$ such that $G\!\in\!\bar G_r^{\cap\xi}\setminus \bar G_r^\rho\,$, one also has $d(x_G,\hat\xi)\!\leq\! 2r$. Proposition \ref{-12-CMS} and definition \ref{Fermi-100} imply the claim. 
\end{proof}
\begin{lem} \label{bords} Given a compact domain $D$ in $(M,g)$, an approximation $K,T$ of given mesh $\rho\,$ (paired with an integer $E$), there exist constants ${\bf c}$ and ${\bf C}$, {\rm depending on ${\mathcal C}_2$ and $D\!\subset\!(M,g)$}, such that one has, for any bone $\xi\in \dagger\!K_{n-2}(D)$ and $r$ small ($2r\!\leq\!{\mathcal C}_1\rho$ holds) 
\begin{equation*} \hbox{\rm Vol}_{d{\mathscr G}} (\backslash \!\!\!\!\bar G_r^{\cap\xi})\!\leq\! {\bf c}\,\rho^{n-3}\,r^3 \ \hbox{\rm and} 
\!\!\!\!\sum_{\substack{\xi\in\dagger\!K_{n-2}(D)\\ z\in{\rm rge}(G)\cap \hat\xi}} \!\!\!\!\!\!\!\!\vert\alpha_{\xi}\vert\!\int_{G\in\backslash \!\!\!\!\bar G_r^{\cap\xi}}\! \vert\hbox{\rm index}_{G,\xi}(z)\vert\,d{\mathscr G}
\!\leq\! {\bf C}\,r^3/\rho\,.
\end{equation*}
\end{lem}
\begin{proof} The boundary $\partial\hat\xi$ is a stratified manifold, union of $\hat\zeta$ with $\zeta$ belonging to
$\dagger K_{n-l}(\!D\!)$ and $l\!\geq\!3$. Each $x\!\in\!{\rm Tub_{2{\mathcal C}_2r}\partial\hat \xi}$ is joined by a shortest geodesic to some $y_1$ in the relative interior of some $n\!-\!l$ dimensional $\hat\zeta$, so $x$ belongs to a Fermi system of coordinates associated to this $\zeta$. Call $T_{2{\mathcal C}_2r}^\zeta$ the set of $x$ at distance $\!\leq \!2{\mathcal C}_2r$ from $\hat\xi$ having closest point in the relative interior of $\hat\zeta$ and define $\bar G_r^\zeta\!:=\!\{G\!\in\! \bar G_r\mid x=G(0,0)\!\in\!T_{2{\mathcal C}_2r}^\zeta\}$. By definition \ref{Fermi-100}, one has $\ \backslash \!\!\!\!\bar G_r^{\cap\xi}\!\subset\!\cup_{\{\zeta \in K\mid \zeta \subset \partial\xi\}}\bar G_r^\zeta\,$. 
\par A description of $T_{2{\mathcal C}_2r}^\zeta$, as done for ${\rm Tub_{\tau}\hat \xi}$ in lemma \ref{CMS-5}, works, but instead of a normal {\em disk}, we get an $l$-dimensional normal ball: one can find a constant $c_{n-l}$ {\it depending on} $D\!\subset\!(M,g),{\mathcal C}_2,l$ such that, for each $\zeta$ of dimension $n\!-\!l$ (with $l\!\geq\!3$), 
one has ${\rm Vol}(T_{2{\mathcal C}_2r}^\zeta)\!\leq\!c_{n-l} \,\rho^{n-l}\,r^3$. This is clear in a Euclidean metric. Then, use the notion of quasi-isometry and the compactness of ${\mathfrak S}^k$, see remark \ref{compact1099} (first, all ${\rm Vol}_{n-l}^g(\hat\zeta)$ are $\leq\!{\rm cst}_{n-l}\,\rho^{n-l}$, where ${\rm cst}_{n-l}\!>\!0$ depends on $D\subset(M,g),l$). From the definition \ref{elemvolume} of $d{\mathscr G}$, one finds ${\bf c}_{n-l}$ depending on $D\!\subset\!(M,g),{\mathcal C}_2,l$ such that
$\hbox{\rm Vol}_{d{\mathscr G}} (\bar G_r^\zeta)\!\leq\! {\bf c}_{n-l}\,\rho^{n-l}\,r^3$, 
implying (since $r\!\leq\!\rho$)
\begin{equation*}  \hbox{\rm Vol}_{d{\mathscr G}} (\backslash \!\!\!\!\bar G_r^{\cap\xi})\leq 
\sum_{l=3}^{n}  \binom{n-1}{l-2}{\bf c}_{n-l}\,  \rho^{n-l}\,r^l  \leq {\bf c}\,\rho^{n-3}\,  
r^3\  .
\end{equation*}
\par The second inequality to prove follows from the first, from $\vert\hbox{\rm index}_{G,\xi}\vert\leq \mu\,$ (remark \ref{bounded222}, proof of lemma \ref{CMS5}), from the bounds on $\vert\alpha_\xi\vert$ (proposition \ref{diedre4},  (\ref{diedre4444})) and on the number of $n$-simplices $\!\subset\!D$,
which is $\leq\!{\rm Cst}\,\rho^{-n}$ (${\rm Cst}\!>\!0$ depends on $D\!\subset\!(M,g)$: proof of lemma \ref{bornangles}).
\end{proof}
\begin{rem}\label{neglige} If, in a measurable space, measures $\mu_1,\dots,\mu_N$ and measurable subsets $A_1,\dots,A_N$ 
verify, for any measurable $B$,  $\mu_i(B)=0$ if $B\cap A_i=\emptyset$, one has $\sum_{i=1}^N\mu_i(\cup_{j=1}^NA_j)=\sum_{i=1}^N\mu_i(A_i)$.
 \end{rem}
\begin{cor} \label{toutdoux} The first claim of proposition {\rm\ref{CMS-7} $(ii)$} holds, indeed
one has in the notations of the preceeding lemma {\rm\ref{bords}}
$$\vert I^{g_0}(\bar G_r\setminus\bar G_r^\rho)\vert \!\leq\!{\bf C}\,r/\rho\ ,\ \ \hbox{thus}\ \ 
\vert I^{g_0}(\bar G_r\setminus\bar G_r^\rho)\vert\rightarrow 0\ \ \hbox{as}\ \ r\rightarrow 0\  .
$$
\end{cor}
\begin{proof} 
Define $[\![I^{g_0}]\!] \!:=\sum_{\xi\in\dagger\!K_{n-2}(D)}\vert\alpha_\xi\vert \vert I^\xi\vert$, with $\vert I^\xi\vert$  defined on any measurable set $\Omega \subset \bar G$ as (see (\ref{CMS9-3333}))
$\vert I^\xi\vert(\Omega)\!:=\!\int_{G\in\Omega}\vert\,[[f^{\xi}(G)]]\,\vert \,{d{\mathscr G}}/{r^2}\ .
$ 
\par Set $A:=\cup_{\xi\in\dagger\!K_{n-2}(D)}\,\backslash \!\!\!\!\bar G_r^{\cap\xi}$, then define $\mu_\xi\!:=\! \vert\alpha_\xi\vert\vert I^\xi\vert_{\mid A}$ and $A_\xi\!:=\!\backslash \!\!\!\!\bar G_r^{\cap\xi}$. If, for 
$G\!\in\!A$, one has $f^{\xi'}(G)\not=0$, one infers from proposition \ref{-12-CMS} that $G$ belongs to
$A_{\xi'}$. Thus the support of $\mu_{\xi'}$ (defined as a measure on $A$) is included in $ A_{\xi'}$.
In view of lemma \ref{blingbling}, one may write 
\begin{equation*}\vert I^{g_0}(\bar G_r\setminus\bar G_r^\rho)\vert=\vert I^{g_0}(\cup_{\xi\in\dagger\!K_{n-2}(D)}(\bar G_r^{\cap\xi}\setminus\bar G_r^\rho))\vert\leq
[\![I^{g_0}]\!] (\cup_{\xi\in\dagger\!K_{n-2}(D)} \,\backslash \!\!\!\!\bar G_r^{\cap\xi})\, .
\end{equation*}
Use remark \ref{neglige} applied to the above defined $\mu_\xi$ and $A_\xi$ to get
\begin{multline*}[\![I^{g_0}]\!] (\cup_{\xi\in\dagger\!K_{n-2}(D)} \,\backslash \!\!\!\!\bar G_r^{\cap\xi})=\!\sum_{\xi\in\dagger\!K_{n-2}(D)}\!\!\vert\alpha_{\xi}\vert\!\int_{G\,\in\,\backslash \!\!\!\!\bar G_r^{\cap\xi}}\vert [[f{}^{\xi}(G)]]\vert\,\frac{d {\mathscr G}}{r^2}\leq\\ \leq
\sum_{\substack{\xi\in\dagger\!K_{n-2}(D)\\ z\in{\rm rge}(G)\cap \hat\xi}}\!\!\vert\alpha_{\xi}\vert\!\int_{G\in\,\backslash \!\!\!\!\bar G_r^{\cap\xi}}\! \vert\hbox{\rm index}_{G,\xi}(z)\vert\,\frac{d{\mathscr G}}{r^2}\leq 
{\bf C}\,r/\rho\  ,
\end{multline*}
inserting, thanks to lemma \ref{CMS5}, the second inequality of lemma \ref{bords} to conclude.
\end{proof}

\subsubsection{The Euclidean case, a model for the index}\label{model}$ $

Because of the scale invariance seen in section \ref{integraleuclidean} (lemmas \ref{CMS1000}, \ref{scalint2}) we only depict intersections of Euclidean squares with codimension $2$ linear subspaces when the side $r$ is $r\!=\!1$.
\par
Consider an $(n-2)$-vector subspace ${\bf X}$ in the Euclidean space $({\mathbb R}^n,{\rm std})$. 
Parallel to (\ref{CMS-33}), define the compact $[{\bf G}^{{\bf X},0}]:=[{\bf G}^{{\bf X},0}]_1$ (the invariance by translations preserving ${\bf X}$ observed in lemma \ref{CMS1000} allows to select the special point $y:=0\!\in\!{\bf X}$)
\begin{equation*} [{\bf G}^{{\bf X},0}]\!:=\!\{{\bf G}\in{\bf G}_1^{{\bf X},0}\mid d(x_{\bf G},{\bf X})\!=\!d(x_{\bf G},0)\!\leq\!2\}
\end{equation*}
which is the union of the open
\begin{equation*} {\bf Z}^{{\bf X},0}\!:=\{{\bf G}\in[{\bf G}^{{\bf X},0}]\mid \hbox{\rm index}_{{\bf G},{\bf X}}=0\}
\end{equation*}
and of its compact complementary 
\begin{equation*}{}^{c}{\bf Z}^{{\bf X},0}\!:=\{{\bf G}\in[{\bf G}^{{\bf X},0}]\mid \hbox{\rm index}_{{\bf G},{\bf X}}\not=0\}\  .
\end{equation*}
Single out the compact subset of {\it zero-measure} of unstable situations 
\begin{equation*} {\bf Un}^{{\bf X},0}\!:=\{{\bf G}\in[{\bf G}^{{\bf X},0}]\mid \partial {\bf G}\cap{\bf X}\not=\emptyset\}\!\subset\![{\bf G}^{{\bf X},0}]\,.
\end{equation*}
\begin{lem} \label{compcroiss1} The set ${}^c{\bf Un}^{{\bf X},0}$ is open and dense in $[{\bf G}^{{\bf X},0}]$. For any ${\bf G}\in{}^c{\bf Un}^{{\bf X},0}\,$, small smooth perturbations $X$ of ${\bf X}$ and $G$ of ${\bf G}$ verify $\hbox{\rm index}_{G,X}=\hbox{\rm index}_{{\bf G},{\bf X}}=0,1$ or $-1$.
\end{lem}

\begin{proof} We have $\hbox{\rm index}_{{\bf G},{\bf X}}\!=\!0,1\,$ or $-1$ {\it except} if ${\bf G}$ and ${\bf X}$ share a segment of line, then $\hbox{\rm index}_{{\bf G},{\bf X}}$ is not defined and $\partial {\bf G}\cap {\bf X}\!\not=\!\emptyset\,$.
Thus $\hbox{\rm index}_{{\bf G},{\bf X}}\!=\!0,1,-1$ is stable, {\it except} if ${\bf G}\!\in\!{\bf Un}^{{\bf X},0}$. For any ${\bf G}\!\in\!{}^c{\bf Un}^{{\bf X},0}$ and any small smooth perturbations $G$ of ${\bf G}$ and $X$ of ${\bf X}$ one has, {\it by transversality}, $\hbox{\rm index}_{G,X}\!=\!\hbox{\rm index}_{{\bf G},{\bf X}}$. 
\end{proof}

\subsubsection{Magnifying again}\label{magnificationagain1}$ $

We magnify the situation near interior points of $\hat\xi\,$, letting the side of squares $G$ tend to $0\,$. This will show up to give enough control to prove that $J^{\xi,y}(\bar G_r^\rho)$ converges to ${\bf J}^{T_y\hat\xi,0_y}$ (see lemma \ref{scalint2}), giving also the key to the second claim in proposition \ref{CMS-7} $(ii)$.
\vskip1mm
We first recall how a {\it magnification of center $y$ and ratio $1/{\bf r}$} acts in a ball $B\!:=\!B(y,\alpha)\!=\!\exp_y(b)$ such that $\exp_y$ is injective in $b\!:=\!B(0_y,\alpha)$ (see section \ref{gloups33} and \cite{GLP}).
Consider the {\em cylinder} $\tilde b\!:=\!I\!\times \!b$, where $I$ is an open interval containing $[0,1]$. Letting ${\bf r}$ tend to $0$, call $h_{y,{\bf r}}$ the diffeomorphism $x\!\mapsto \!x_{\bf r}\!=\!\exp_y({\bf r}\exp_y^{-1}x)$ from $B(y,\alpha)$ to $B(y,{\bf r}\alpha)$. Then, perform a homothety on the scale, i.e. on the metric $g$, of ratio $1/{\bf r}^2$, magnifying the ``small'' $B(y,{\bf r}\alpha)$ to the size of the ``original'' $B\,$. Pull everything on $T_yM$ (through $\exp_y^{-1}$) and observe that a homothety of center $0_y$ and ratio $1/{\bf r}$ produces the rescaling by $1/{\bf r}^2$ of the pulled back $g$, i. e. of $\exp_y^\ast g\,$: one still may view $x_{\bf r}\!=\!h_{y,{\bf r}}(x)$ as the point $u\!:=\!\exp_y^{-1}x\in b$. One catches this in a picture by defining on the cylinder 
$\tilde b$ the $C^\infty$-Riemannian metric equal to $dt^2+g_{\bf r}$ on each slice $\{{\bf r}\}\times b$, where $g_{\bf r}\!=\!\exp_y^\ast(h_{y,{\bf r}}^\ast g)/{\bf r}^2$ for ${\bf r}\!\not=\!0$ is extended by $g_y$ at ${\bf r}\!=\!0$. 
\par The cylinder 
$\tilde b$ replaces a {\em cone}, which naturally pictures the situation if no rescaling of the metric after mapping $B$ by $h_{y,{\bf r}}$ is done. 

\begin{lem} \label{Thetaprol} Define for ${\bf r}>0$
\begin{equation*}\label{magnif1}\Theta_{y,{\bf r}}:({\bf r},u)\!\in\!\{{\bf r}\}\!\times\!b\subset\tilde b\longmapsto \Theta_{y,{\bf r}}({\bf r},u):=u\in b:=B(0_y,\alpha)\subset T_yM\ .
\end{equation*}
Identifying as above $x_{\bf r}\!=\!h_{y,{\bf r}}(x)\!\in\!B(y,{\bf r}\alpha)$ with $({\bf r},u)\!\in\!\{{\bf r}\}\!\times\!b\subset\tilde b$, where $u\!=\!\exp_y^{-1}x$, one may also view the map $\Theta_{y,{\bf r}}$ for ${\bf r}\!>\!0$ as
\begin{equation*}\label{magnif11111}\Theta_{y,{\bf r}}:x_{\bf r}\!\in\!B(y,{\bf r}\alpha)\longmapsto \Theta_{y,{\bf r}}(x_{\bf r})=\frac{\exp_y^{-1}x_{\bf r}}{\bf r}\in B(0_y,\alpha)\subset T_yM\ .
\end{equation*}
In accordance with the first version of $\Theta_{y,{\bf r}}$, one may prolong its second version at ${\bf r}=0$ in a $C^\infty$-map, setting $\Theta_{y,0}\!:=\!{\rm Id}$ on $B(0_y,\alpha)$ (and ``blowing up'' $B(y,\alpha)$ at $y$).

For ${\bf r}\!\in\!I\!\setminus\!\{0\}$, the lines $(I\!\setminus\!\{0\})\!\times\!\{u\} \!\subset \! (I\!\setminus\!\{0\})\!\times\!B(0_y,\alpha) $ in $\tilde b$ are paired in the above identifications with the $g$-geodesics $\exp_y({\bf r}u)$ running through $y$. 
\par Recall the {\em pseudo-homothety} ${\mathcal H}_{y,{\bf r}}$ of definition {\rm\ref{CMS-1000}}.
The $2$-frame $(u,v)$ which is defined along the geodesic $x_{\bf r}$ by setting $(u_{x_{\bf r}},v_{x_{\bf r}})\!=\!\frac{1}{{\bf r}r}({\partial {\mathcal H}_{y,{\bf r}}(G_{r})}/{\partial s},{\partial {\mathcal H}_{y,{\bf r}}(G_{r})}/{\partial t})(0,0)$ is parallel and orthonormal, one has $(u_{x_0},v_{x_0})=(u_y,v_y)$; moreover
\begin{equation}\label{causetoujours3333} (u_y,v_y)=\lim_{{\bf r}\rightarrow0}(\frac{\partial \Theta_{y,{\bf r}}\!\circ\!{\mathcal H}_{y,{\bf r}}(G_{r})}{\partial s},\frac{\partial\Theta_{y,{\bf r}}\!\circ\!{\mathcal H}_{y,{\bf r}}(G_{r})}{\partial t})(0,0)\ .
\end{equation}
 \end{lem}
\begin{proof} See section \ref{gloups33}, definition {\rm\ref{CMS-1000}} and (\ref{madre}) (below).
\end{proof}
\begin{rem}\label{magnif11} Notice that, given any submanifold $N$ passing through $y$, the magnified $\Theta_{y,{\bf r}}(N\cap B(y,{\bf r}\alpha))$ evolves to $T_yN\cap B(0_y,\alpha)\subset T_yM$ as ${\bf r}$ tends to $0$. So does $\Theta_{y,{\bf r}}(\hat\xi\cap B(y,{\bf r}\alpha))$ to
$T_y\hat\xi\cap B(0_y,\alpha)\subset T_yM$.
\end{rem}
\begin{defn} \label{jigstraight} For $\bar G$ close to $\underline G$, given $i,\,\hat\xi\in{\mathfrak S}_i, \,y\in\tilde\xi,\,r>0$ (see lemma \ref{compactsimpl} for ${\mathfrak S}_i$), a {\it pseudo-straightening} ${\mathcal F}_{y,{\bf r}}$ is defined on $\bar G_r^{\xi,y}$ from the {\em pseudo-homothety} ${\mathcal H}_{y,{\bf r}}$ of definition {\rm\ref{CMS-1000}} by setting (${\mathcal F}_{y,{\bf r}}$ sends a square of origin $x$ to a square in $T_yM$)
\begin{equation*}\label{magnif2}{\mathcal F}_{y,{\bf r}}:G\in \bar G_r^{\xi,y}\longmapsto 
{\mathcal F}_{y,{\bf r}}(G):=\Theta_{y,{\bf r}}\circ{\mathcal H}_{y,{\bf r}}(G)\subset T_yM\ .
\end{equation*}
\end{defn}
\begin{lem} \label{limitsqu} First, ${\mathcal F}_{\cdot,\ast}(\diamond)$ is smooth in $(y,{\bf r},G)\!\in\!\hat\xi\times[0,\tau/2]\times\bar G_r^{\xi,y}\,$. Given $G_{x,u,v}\!\in\!\bar G_r^{\xi,y}$ (definition {\rm\ref{CMS10333}}), ${\mathcal F}_{y,{\bf r}}(G_{x,u,v})$ tends, as ${\bf r}\rightarrow0$, towards a Euclidean square ${\bf G}\!\in\!{\bf G}_r^{{\bf X},0_y}$ (with ${\bf X}\!=\!T_y\hat\xi$)
\begin{equation*}\label{Euclidsqu}{\mathcal F}_{y,0}(G_{x,u,v}):={\bf G}:(\varrho,\varsigma)\in\bar I^2\longmapsto \exp_y^{-1}x+\varrho  r u_y+ \varsigma  r v_y\in T_yM\ ,
\end{equation*}
where $u_y,v_y\in T_{\exp_y^{\!-\!1}\!x}T_yM\equiv T_yM$ are parallel translated from $u,v\in T_xM$ (with $G=G_{x,u,v}$) along the ray joining $x$ to $y\in\hat\xi$. 
One has 
\begin{equation*}\label{n'enjetezplus33-2}{\mathcal F}_{y,0}([\bar G^{\xi,y}]_r):=\lim_{{\bf r}\rightarrow0}{\mathcal F}_{y,{\bf r}}([\bar G^{\xi,y}]_r)=[{\bf G}^{{\bf X},0_y}]_r\ .
\end{equation*} 
For any $x\!\in\!\Delta_{2r}^{\xi,y}$ (lemma {\rm\ref{CMS-5}}), ${\mathcal F}_{y,0}$ {\em has an inverse} ${\mathcal F}_{x,y,0}^{-1}$ from $[{\bf G}^{{\bf X},0_y}]_r\cap\{{\bf G}\!\mid\! x_{\bf G}\!=\!\exp_y^{\!-1}x\}$ to $[\bar G^{\xi,y}]_r\cap\{G \!\mid\! x_G\!=\!x\}$ sending the Euclidean square ${\bf G}_{\exp_y^{\!-1}x,u,v}$ of side $r$ to $G_{r,x,u,v}$ of side $r$.
One has for $G\!\in\![\bar G^{\xi,y}]_r$
\begin{equation}\label{louloup(e)} (\!{\mathcal F}_{x,y,0}^{-1}\!)^\ast {{\frac{d{\mathscr G}_{\mid \bar G_r^{\xi,y}}}{r^2}}}_{\!\!\mid G}\!\!\!=\!\frac{\psi(\!G_r\!)d{\mathpzc G}_{\mid {\bf G}^{{\bf X},0_y}}}{t_{G_r}r^2}\!=\!\frac{\psi(\!G_r\!){\Omega^\xi}}{t_{G_r}}\,\hbox{where}\,\lim_{r\rightarrow0}\!\frac{\psi(\!G_r\!)}{t_{G_r}}\!=\!1.
\end{equation}
\end{lem}
\begin{proof} First, observe that through definition  \ref{CMS-1000} and (\ref{magnifique1}) 
${\mathcal H}_{y,{\bf r}}$ reads 
\begin{equation*} {\mathcal H}_{y,{\bf r}} : (r,x,u,v)\mapsto ({\bf r}\,r,x_{\bf r},u_{x_{\bf r}},v_{x_{\bf r}})\  ,
\end{equation*}
with $x_{\bf r}=h_{y,{\bf r}}(x)$, so one may understand $\Theta_{y,{\bf r}}\circ{\mathcal H}_{y,{\bf r}}$ as (lemma \ref{Thetaprol})
\begin{equation} \label{madre}\Theta_{y,{\bf r}}\circ{\mathcal H}_{y,{\bf r}} : (r,x,u,v)\mapsto (\exp_y^{-1}x,u_{x_{\bf r}},v_{x_{\bf r}})\  .
\end{equation}
Thus $\Theta_{y,{\bf r}}\circ{\mathcal H}_{y,{\bf r}}$ is $C^\infty$ on $[0,\tau/2]\!\times\!\tilde\xi\!\times\!\bar G_r^{\xi,y}$.
\par Write ${\mathcal H}_{y,{\bf r}}(G_{x,u,v})\!=\!G_{x_{\bf r},u,v}\!\in\!\bar G_{{\bf r}r}^{\xi,y}$ (\ref{magnifique2}).  Pull on $T_yM$ by $\exp_y^{-1}$
\begin{multline*} \lim_{{\bf r}\rightarrow0}{\mathcal F}_{y,{\bf r}}(G_{x,u,v})(\varrho ,\varsigma)
=\lim_{{\bf r}\rightarrow0}\Theta_{y,{\bf r}}\circ {\mathcal H}_{y,{\bf r}}(G_{x,u,v})(\varrho ,\varsigma)=\\=
\lim_{{\bf r}\rightarrow0}\frac{1}{{\bf r}}\exp_y^{-1}\circ\,{\mathcal H}_{y,{\bf r}}(G_{x,u,v})(\varrho ,\varsigma)=\lim_{{\bf r}\rightarrow0}\frac{1}{{\bf r}}\exp_y^{-1}\circ \,G_{x_{\bf r},u,v}(\varrho {\bf r},\varsigma {\bf r})=\\=\lim_{{\bf r}\rightarrow0}\frac{1}{{\bf r}}(\exp_y^{-1}\circ \,G_{x_{\bf r},u,v}(\varrho {\bf r},\varsigma {\bf r})-\exp_y^{-1}\circ \,G_{x_{\bf r},u,v}(0,0))+\\+\lim_{{\bf r}\rightarrow0}\frac{1}{{\bf r}}\exp_y^{-1}\circ \,G_{x_{\bf r},u,v}(0,0)= \varrho r u_y+\varsigma r v_y+\exp_y^{-1}x\  .
\end{multline*}
The last equality follows from $\exp_y^{-1}(x)\!=\!\frac{1}{{\bf r}}\exp_y^{-1}\circ \,G_{x_{\bf r},u,v}(0,0)$, true for any ${\bf r}\!\not=\!0$ under consideration, from the fact that the first term in the penultimate sum is a linear mapping in $(\varrho,\varsigma)$ and from (\ref{causetoujours3333}).
Concluding, lemma \ref{n'enjetezplus1} (\ref{n'enjetezplus2}), (\ref{n'enjetezplus4}), (\ref{n'enjetezplus3}) and (\ref{n'enjetezplus33}) imply (\ref{louloup(e)}).
\end{proof}

\subsubsection{Chasing $\hbox{\rm index}(G_r,\xi)$, with $G_r\in\bar G_r$ and $\xi$ fixed bone, as $r\!\rightarrow\!0$}\label{chasing1022}$ $

Until the end of the proof of Regge's theorem, we write $r, {\bf r}$ whenever they occur.
Near $y$ interior to $\hat\xi\,$, the smaller ${\bf r}\!>\!0$, the more  $\Theta_{y,{\bf r}}(\hat\xi\!\cap\! B(y,{\bf r}\alpha))$ looks like $T_y\hat\xi\cap B(0_y,\alpha)\!\subset\! (T_yM,g_y)$ (remark \ref{magnif11}). Pulling on $(T_yM,g_y)$, the smaller ${\bf r}\!>\!0$, the more ${\mathcal F}_{y,{\bf r}}([\bar G^{\xi,y}]_r)$ is smoothly $g_y$-close to the Euclidean $[{\bf G}^{T_y\hat\xi,0_y}]_r$ in $(T_yM,g_y)$ (lemma \ref{limitsqu}).
\par
With respect to intersections with $\Theta_{y,{\bf r}}(\hat\xi)$, the situation of the family ${\mathcal F}_{y,{\bf r}}(G_r)\!=\!
\Theta_{y,{\bf r}}\circ{\mathcal H}_{y,{\bf r}}(G_r)$, with $G_r\!\in\![\bar G^{\xi,y}]_r$, depends, for small enough ${\bf r}\!>\!0$, on the model-situation (for ${\bf r}\!=\!0$) of ${\mathcal F}_{y,0}(G_r)\!=\!{\bf G}_r\in[{\bf G}^{T_y\hat\xi,0_y}]_r$ with respect to
$T_y\,\hat\xi\!=\!\lim_{{\bf r}\rightarrow0}\Theta_{y,{\bf r}}(\hat\xi)$. Indeed, {\em since $\Theta_{y,{\bf r}}$ is a diffeomorphism} and the unstable ${\bf Un}_r^{T_y\hat\xi,0_y}$ is a zero-measure set in $[{\bf G}^{T_y\hat\xi,0_y}]_r$, by lemma \ref{compcroiss1} $\hbox{\rm index}_{{\mathcal H}_{y,{\bf r}}G_r,\xi}$ tends almost everywhere to $\hbox{\rm index}_{{\bf G}_r,T_y\hat\xi}$. 
\par Summarizing
 \begin{lem} \label{limite1} Given a bone $\xi$ and $y$ interior to $\hat\xi$, for almost any ${\bf G}_r\!\in\![{\bf G}^{T_y\hat\xi,0_y}]_r$ (all ${\bf G}_r\in {}^c{\bf Un}_r^{T_y\hat\xi,0_y}$) and for $G_r\in[\bar G^{\xi,y}]_r$ such that ${\mathcal F}_{y,0}(G_r)={\bf G}_r$
\begin{equation*}\label{scalint+1}\lim_{{\bf r}\rightarrow0}\hbox{\rm index}_{{\mathcal H}_{y,{\bf r}}(G_r),\xi}=\hbox{\rm index}_{{\bf G}_r,T_y\hat\xi}\  .
\end{equation*}
\end{lem}
\subsubsection{Lebesgue's theorem in approach}\label{lebesguons}$ $

\par The next goal is, for $y$ interior to $\hat\xi$, to prove $J^{\xi,y}(\bar G_r)\!\rightarrow\!{\bf J}^{T_y\hat\xi,0_y}={\bf c}_2(n)$ (lemma \ref{scalint2}) as $r\!\rightarrow\!0$. Starting from (\ref{CMS8bis2}), one may rewrite $J^{\xi,y}(\bar G_r)$ (see lemma \ref{CMS-5})
\begin{multline*}\label{Lebesguise}J^{\xi,y}(\bar G_r)\!:=\!\!\int_{G_r\in[\bar G^{\xi,y}]_r}\!\! [[\mathfrak{f}{\,}^\xi(G_r)]] \,\frac{d{\mathscr G}_{\mid \bar G^{\xi,y}_r}}{r^2}_{\mid G_r}
=\\=\int_0^{2r}\!\!\int_{{\bf N}^1(T_y\tilde\xi)}\int_{G_r\in[\bar G^{\xi,y}]_r\cap\{x_{G_r}\equiv(t,w)\}}\!\! [[\mathfrak{f}{\,}^\xi(G_r)]] \,\frac{d{\mathscr G}_{\mid [\bar G^{\xi,y,t,w}]_r}}{r^2}_{\mid G_r}\,dw\,dt\  ,
\end{multline*}
where $[\bar G^{\xi,y,t,w}]_r\!:=\![\bar G^{\xi,y}]_r\cap\{x_{G_r}\!\equiv\!(t,w)\}$.
Using the change of variables ${\mathcal F}_{x,y,0}^{-1}$ given in lemma \ref{limitsqu} and (\ref{louloup(e)}),  for any $\xi,y,x,r$, each $\bar G_r$ may be seen in $(T_yM,g_y)$. Indeed, one gets for $w\!\in\!{\bf N}^1(\tilde\xi)$ and $t\!\in\![0,2r]$
\begin{multline*}\!\!\!
\int_{G_r\in[\bar G^{\xi,y}]_r\cap\{x_{G_r}\equiv(t,w)\}}\!\! [[\mathfrak{f}{\,}^\xi(G_r)]] \,\frac{d{\mathscr G}_{\mid [\bar G^{\xi,y,t,w}]_r}}{r^2}_{\mid G_r}=\\=\int_{{\bf G}_r\in[{\bf G}^{T_y\hat\xi,0_y}]_r\cap\{x_{{\bf G}_r}\equiv(t,w)\}}\!\! [[\mathfrak{f}{}^\xi({\mathcal F}_{x,y,0}^{-1}({\bf G}_r)]]
\ ({\mathcal F}_{x,y,0}^{-1})^\ast \frac{d{\mathscr G}_{\mid [\bar G^{\xi,y,t,w}]_r}}{r^2}_{\!\mid G_r}=
\\=\int_{{\bf G}_r\in[{\bf G}^{T_y\hat\xi,0_y}]_r\cap\{x_{{\bf G}_r}\equiv(t,w)\}}\!\!\! [[\mathfrak{f}{}^\xi({\mathcal F}_{x,y,0}^{-1}({\bf G}_r)]]
\ \frac{\psi({\mathcal F}_{x,y,0}^{-1}({\bf G}_r))}{t_{{\mathcal F}_{x,y,0}^{-1}({\bf G}_r)}}\,\frac{d{\mathpzc G}_{\mid {\bf G}^{{\bf X},0_y,t,w}}}{r^2}\  .
\end{multline*} 
As ${\bf r}\rightarrow0$, the limit of $[[\mathfrak{f}{}^\xi({\mathcal F}_{x_{\bf r},y,0}^{-1}({\bf r}{\bf G}_r))]]\!=\![[\mathfrak{f}{}^\xi({\mathcal H}_{y,{\bf r}}(G_r)]]$ is given by 
\par 
\begin{lem} \label{limite1-2} Given a bone $\xi$, given $y$ interior to $\hat\xi$, one has for almost any ${\bf G}_r\in[{\bf G}^{T_y\hat\xi,0_y}]_r$ and $G_r\in[\bar G^{\xi,y}]_r$ such that ${\mathcal F}_{y,0}(G_r)={\bf G}_r$
\begin{equation*}\lim_{{\bf r}\rightarrow0}\,[[\mathfrak{f}{}^\xi({\mathcal H}_{y,{\bf r}}(G_r)]]= \hbox{\rm index}_{{\bf G}_r,T_y\hat\xi}\,\cos\angle_{g_y} (\ast_{g_y}\,T_y\hat\xi, {\bf G}_r)\  .
\end{equation*}
\end{lem}
\begin{proof} Thanks to lemma \ref{limite1}, for any $G_r\in[\bar G^{\xi,y}]_r$ verifying ${\mathcal F}_{y,0}(G_r)={\bf G}_r\in{}^c{\bf Un}^{T_y\hat\xi,0_y}$ and for any $r>0$ small enough, one has
\begin{equation*} \hbox{\rm index}_{{\mathcal H}_{y,{\bf r}}(G_r),\xi}=\hbox{\rm index}_{{\bf G}_r,T_y\hat\xi}=0,1\  \hbox{or}\  \  -1\  .
\end{equation*}
If $\hbox{\rm index}_{{\bf G}_r,T_y\hat\xi}\!=\!0$, then for ${\bf r}\!>\!0$ small enough, $[[\mathfrak{f}{}^\xi({\mathcal H}_{y,{\bf r}}(G_r))]]$ equals $0$. If 
$\hbox{\rm index}_{{\bf G}_r,T_y\hat\xi}\!=\!1$ or $-1$ (so ${\bf G}_r$ is in ${}^c{\bf Un}^{T_y\hat\xi,0_y}$), for $G_r$ such that ${\mathcal F}_{y,0}(G_r)\!=\!{\bf G}_r$ and for any ${\bf r}\!>\!0$ small enough, there is a unique
point $z_{{\bf r}} \!\in\! {\rm rge}({\mathcal H}_{y,{\bf r}}G_r)\cap \hat\xi$ (depending differentiably on ${\bf r}$, use the implicit function theorem).
Writing $x_{\bf r}\!=\!\exp_y({\bf r}\exp_y^{-1}x)$, as ${\bf r}$ tends to $0$
\begin{equation*}\langle \ast_g\, {\mathcal P}_{\tilde A_{z_{\bf r}}^{-1}} (T_{z_{\bf r}}\hat\xi/{\Vert T_{z_{\bf r}}\hat\xi\Vert_g}), T_{x_{\bf r}}{\mathcal H}_{y,{\bf r}}(G_r)/{\Vert T_{x_{\bf r}}{\mathcal H}_{y,{\bf r}}(G_r)\Vert_g}\rangle_g
\end{equation*} 
tends to the pending $\cos\angle_{g_y}(\ast_{g_y}\, T_y\hat\xi, {\bf G}_r)$, thus the conclusion.
\end{proof}
But ${\psi(G_{{\bf r}r})}/{t_{G_{{\bf r}r}}}\!\rightarrow\!1$ as ${\bf r}\!\rightarrow\!0$, see (\ref{louloup(e)}): if $y\!\in\!\hat\xi\!\setminus\!\partial\hat\xi$, 
the Lebesgue theorem thus implies $\lim_{r\!\rightarrow\!0}\!J^{\xi,y}(\bar G_r)\!=\!{\bf J}^{T_y\hat\xi,0_y}\!=\!{\bf c}_2(n)$ (see lemma \ref{scalint2}).
\par From lemma \ref{integrals1111}, recall that ${J}^{g_0}(\bar G_r^\rho)\!=\!\!\sum_{\xi\in \dagger\!K_{n-2}(D)}\!\alpha_\xi\!
\int_{\hat\xi}{J}^{\xi,y}(\bar G_r)\,dy$. Since the functions $[[\mathfrak{f}{}^\xi(\cdot)]]$ are bounded by lemma \ref{CMS5}, the Lebesgue theorem implies, in view of lemma \ref{bornangles}, the second claim in proposition \ref{CMS-7} $(ii)$, 
fulfilling the proof of Regge's theorem \ref{Regge1234}.

\vskip1cm

\section{Appendices}

\subsection{Proof of lemma \ref{squcurvparrtrslt}}\label{squcurvparrtrslt1}
\begin{proof} Let $w$ be a vector in $E_x$. Construct the following field along $G$. Set ${\mathcal W}_{0,\tau}$ to be the parallel translation of $w$ along $\delta_0^\tau$, then define ${\mathcal W}_{\varrho,\tau}$ as the parallel translation of ${\mathcal W}_{0,\tau}$ along $\gamma_\tau^\varrho$. This defines a differentiable vector field ${\mathcal W}$ along $G$ such that ${\mathcal W}(0,0)=w\,$, since lemma \ref{fait12} tells 
${\mathcal P}_{\delta;t} {\mathcal W}$ and ${\mathcal P}_{\gamma;s} {\mathcal P}_{\delta;t} {\mathcal W}$ are differentiable in $(s,t)$ and 
${\mathcal W}_{\varrho,\tau}={\mathcal P}_{\gamma;\varrho} {\mathcal P}_{\delta;\tau} w$. Lemma \ref{fait12} still gives ${\mathcal P}_{\Gamma_{\sigma,\upsilon}}^{-1}(w)={\mathcal P}_{\gamma;\sigma}^{-1} {\mathcal P}_{\delta;\upsilon}^{-1} {\mathcal W}(\sigma,\upsilon)$ is differentiable.
\par On the other hand, one has
\begin{multline}\label{blue1}(\frac{{\mathcal P}_{\Gamma_{\sigma,\upsilon}}^{-1}-\hbox{\rm Id}}{\sigma\upsilon}) \,w=
\frac{{\mathcal P}_{\gamma;\sigma}^{-1} {\mathcal P}_{\delta;\upsilon}^{-1} {\mathcal W}(\sigma,\upsilon)-w}{\sigma\upsilon}=\cr
=\frac{1}{\sigma}\ {\mathcal P}_{\gamma;\sigma}^{-1} (\frac{{\mathcal P}_{\delta;\upsilon}^{-1} {\mathcal W}(\sigma,\upsilon)-{\mathcal W}(\sigma,0)}{\upsilon})\  .
\end{multline}
By the definition of parallel translation, one also has
\begin{equation}\label{blue2}\nabla_{\!dG(\sigma,0)(\frac{\partial}{\partial \upsilon})}{\mathcal W}\!=\!\lim_{\upsilon\rightarrow0}
\frac{{\mathcal P}_{\delta;\upsilon}^{-1} {\mathcal W}(\sigma,\upsilon)\!-\!{\mathcal W}(\sigma,0)}{\upsilon}\hskip1.5mm
\hbox{and}\hskip1.5mm \nabla_{\!dG(0,0)(\frac{\partial}{\partial \upsilon})}{\mathcal W}\!=\!0\, .
\end{equation}
One derives from the lines (\ref{blue1}) and (\ref{blue2}) 
\begin{equation*}\lim_{\sigma\rightarrow0}\lim_{\upsilon\rightarrow0}(\frac{{\mathcal P}_{\Gamma_{\sigma,\upsilon}}^{-1}-\hbox{\rm Id}}{\sigma\upsilon}) \,w=\nabla_{dG(0,0)(\frac{\partial}{\partial \sigma})}\nabla_{dG(\sigma,0)(\frac{\partial}{\partial \upsilon})}{\mathcal W}\  ,
\end{equation*}
and this gives
\begin{equation*}\lim_{\sigma\rightarrow0}\lim_{\upsilon\rightarrow0}(\frac{{\mathcal P}_{\Gamma_{\sigma,\upsilon}}^{-1}-\hbox{\rm Id}}{\sigma\upsilon}) \,w=R_x(dG(\frac{\partial}{\partial \sigma}),dG(\frac{\partial}{\partial \upsilon}))w\  .
\end{equation*}
Indeed, by definition of ${\mathcal W}$, one has $\nabla_{dG(\sigma,\upsilon)(\frac{\partial}{\partial \sigma})}{\mathcal W}\equiv0$ and one also has $[\frac{\partial}{\partial \sigma},\frac{\partial}{\partial \upsilon}]=0$, so the images by $dG$ of $\frac{\partial}{\partial \sigma}$ and $\frac{\partial}{\partial \upsilon}$ commute too.
\end{proof}
\subsection{Proof of lemma \ref{movebase}}
\label{Taylor2} \begin{proof}
The loop $\Gamma_{\varrho,\tau;\sigma,\upsilon}$ is the boundary of the square $G_{\mid [\varrho,\sigma]\times[\tau,\upsilon]}$ and, in view of its definition (\ref{nota1233}), making the obvious change of variables, (\ref{blue3}), (\ref{blue4}), (\ref{formula12}) hold if $(\varrho,\tau),(\sigma-\varrho,\upsilon-\tau)$ play the former roles of $(0,0),(\sigma,\upsilon)$.  
\par If $(\varrho,\tau)=(0,0)$, in (\ref{formula12}) $o_x$ reads $o_x(\sigma^2+ \upsilon^2)=\omega_x(\sigma,\upsilon)\ 
((\sigma,\upsilon),(\sigma,\upsilon))$ where $\omega_x(\sigma,\upsilon)(p,p)$ is quadratic in $p\in I^2$ ($I$ open interval containing $[0,1]$) with values in ${\bf End}(E_x)$ and, comparing (\ref{formula12}) to a Taylor expansion with integral remainder of order $2$ (since ${\mathcal P}_{\Gamma_{\sigma,\upsilon}}^{-1}$ is $C^2$), one gets 
\begin{equation}\label{manquait1}
\omega_x(\sigma,\upsilon):=
\int_0^1 (1-t)\,d^{(2)}({\mathcal P}_{\Gamma_{\sigma,\upsilon}}^{-1})\,(t(\sigma,\upsilon))\,dt-\frac{1}{2}\,d^{(2)}({\mathcal P}_{\Gamma_{\sigma,\upsilon}}^{-1})(0,0)\ .
\end{equation}
So $\epsilon_x(\sigma,\upsilon):=\Vert \omega_x(\sigma,\upsilon)\Vert=o(1)$ verifies
\begin{equation}\label{manquait}\Vert o_x(\sigma^2+ \upsilon^2)\Vert\leq (\sigma^2+ \upsilon^2)\ \epsilon_x(\sigma,\upsilon)\ .
\end{equation}
Then, the last claim to prove follows from corresponding expressions (\ref{manquait1}) of $\omega_{G(\varrho,\tau))}(\sigma-\varrho,\upsilon-\tau)$ and corresponding (\ref{manquait}).
\end{proof}

\subsection{Proof of lemma \ref{MC3-2}}
\label{MC30}
\begin{proof} Starting from $A^{-1}\,A=\hbox{\rm Id}$, one gets $d A^{-1}\,A+A^{-1}\,dA=0$ and thus $d (A^{-1})=-A^{-1}\,d A\ A^{-1}$.
The following computation in coordinates establishes the lemma, using (\ref{MC2})
\begin{multline*} (d\,(A^{-1}\,dA))_i^j=(\sum_{k} (dA^{-1})_k^j\wedge d A_i^k)=\cr
=-\sum_{k,r,s}[(A^{-1})_r^j(d A)_s^r(A^{-1})_k^s]\wedge(d A)_i^k
=-\sum_{s}[A^{-1}(d A)]_s^j\wedge [A^{-1}(d A)]_i^s=\cr=
-[(A^{-1}\,d A)\wedge (A^{-1}\,d A)]_i^j\  .\qedhere
\end{multline*}

\end{proof}
\subsection{Proof of (\ref{diantre44}), needed in the proof of lemma \ref{prox6}}\label{diantre44bis}
\begin{proof}
Writing here
$L$ in place of ${\mathfrak G}{\mathfrak G}$, use the formalism of lemmas \ref{MC1}, \ref{MC3-2} and appendix \ref{MC30} to express the pull-back by $L:={\mathfrak G}{\mathfrak G}$ (on $\bar I^2$) of the $2$-form $(A^{-1}\,dA)\wedge(A^{-1}\,dA)$ (this $2$-form is vector valued) 
\begin{multline*}L^\ast\,((A^{-1}\,dA)\wedge(A^{-1}\,dA))=
L^\ast(\sum_{i,j}
[(A^{-1}\,dA)\wedge(A^{-1}\,dA)]_i^j(E_i^j)^*)=
\\=\sum_{i,j}
[L^\ast(A^{-1}\,dA)\wedge L^\ast(A^{-1}\,dA)]_i^j \  (E_i^j)^*=\\=\sum_{i,j}\sum_k
[L^{-1} dL]_k^j\wedge[L^{-1}dL]_i^k\, (E_i^j)^*\  .
\end{multline*}
So, one has (using the left-invariance of $\underline{g}_x$)
\begin{multline} \label{diantre11}\Vert \,L^\ast ((A^{-1}\,dA)\wedge(A^{-1}\,dA))\,\Vert_{x,x}^2=\Vert \,dL\,\wedge \,dL\,\Vert_{\underline{g}_x}^2=\\
=\sum_{i,j}\Vert\sum_k
[dL]_k^j\wedge[dL]_i^k\Vert^2\leq
n\sum_{i,j,k}\Vert
[dL]_k^j\wedge[dL]_i^k\Vert^2\leq\\\leq n\sum_{i,j,k,l}\Vert
[dL]_k^j\wedge[dL]_i^l\Vert^2\ .
\end{multline}
But $[dL]_k^j=[\frac{\partial L}{\partial s}ds+\frac{\partial L}{\partial t}dt]_k^j=[\frac{\partial L}{\partial s}]_k^j ds+[\frac{\partial L}{\partial t}]_k^j dt\,$, thus 
\begin{equation}\label{diantre22}[dL]_k^j\wedge[dL]_i^l=([\frac{\partial L}{\partial s}]_k^j[\frac{\partial L}{\partial t}]_i^l-
[\frac{\partial L}{\partial t}]_k^j[\frac{\partial L}{\partial s}]_i^l)\,ds\wedge dt\,.
\end{equation}
And $\frac{\partial L}{\partial s}=\sum_{j,k}[\frac{\partial L}{\partial s}]_k^j\,(E_k^j)^*$
(in ${\bf End}(E_x)$), so (use (\ref{diantre22}) for (\ref{diantre55}))
\begin{gather}\label{diantre33}\frac{\partial L}{\partial s}\wedge\frac{\partial L}{\partial t}=\frac{1}{2}\sum_{i,j,k,l}([\frac{\partial L}{\partial s}]_k^j [\frac{\partial L}{\partial t}]_i^l-[\frac{\partial L}{\partial t}]_k^j[\frac{\partial L}{\partial s}]_i^l)\,(E_k^j)^*\wedge(E_i^l)^*\  ,\\
\label{diantre55}4\Vert\,\frac{\partial L}{\partial s}\wedge\frac{\partial L}{\partial t}\Vert_{\underline g_x}^2=\sum_{i,j,k,l}\Vert\,[dL]_k^j\wedge[dL]_i^l\,\Vert^2\  .
\end{gather}
Thus (\ref{diantre11}) and (\ref{diantre55}) give (recalling $L:={\mathfrak G}{\mathfrak G}$)
\begin{multline*}  \Vert \,L^\ast ((A^{-1}\,dA)\wedge(A^{-1}\,dA))\,\Vert_{x,x}\leq \\\leq2\sqrt{n}\ 
\Vert \frac{\partial L}{\partial s}\wedge\frac{\partial L}{\partial t}\Vert_{\underline g_x}=2\sqrt{n}\ 
\Vert\frac{\partial {\mathfrak G}{\mathfrak G}}{\partial s}\wedge\frac{\partial {\mathfrak G}{\mathfrak G}}{\partial t}\Vert_{\underline g_x}\ .\qedhere
\end{multline*}
\end{proof}

\subsection{A remark on $\Delta$-complexes} \label{Deltacomplex}
\begin{rem} Thinking to a $\Delta$-complex (see {\rm\cite{Ha}}, p.{\rm102})- if $\eta$ is produced by an identification of {\it two} isometric distinct faces $\eta',\eta''\!\in\!\dagger\!K_{n\!-\!1}$ of the {\it same} $n$-simplex $\sigma\,$,  replace $\hbox{\rm st} (\eta)$ (see definition {\rm\ref{ast1}}) by the {\em reduced open star} $\check{\hbox{\rm st}} (\eta)$: it is the union of the interiors of a sub-simplex having vertices $\varpi_\sigma$ and the $p'_i\,$, of another having vertices  $\varpi_\sigma$ plus the $p''_i$ and of $\eta$ (and $p'_i,p''_i$ are identified in pairs, giving the vertices $p_i$ of  $\eta\subset\sigma$). There exists a local orientation around $\eta\,$ since we are led back to the situation of lemma {\rm\ref{aster3}}, i. e. $\check{\hbox{\rm st}} (\eta)$ is homeomorphic to a ball.
\end{rem}

\subsection{An elementary lemma in the Euclidean plane}\label{diedre-4bis}$ $

\begin{lem} \label{diedre-4} A linear map $L$ from ${\mathbb R}^2$ equipped with the usual Euclidean norm $\Vert\cdot\Vert$ to itself whose matrix in the canonical basis is
\begin{equation*} \begin{pmatrix} \cos \alpha -1&-\sin\alpha\\ \sin\alpha&\cos \alpha -1
\end{pmatrix}
\end{equation*}
has operator norm
\begin{equation*}\Vert\,L\,\Vert=\sup_{u\in{\mathbb R}^2\setminus\{0\}}\,\frac{\Vert\,L(u)\,\Vert}{\Vert\,u\,\Vert}=\sqrt{2(1-\cos\alpha)}\  .
\end{equation*}
\end{lem}
\begin{proof} In fact, one has
\begin{equation*}\Vert\,L\,\Vert=\sup_{u\in{\mathbb R}^2\setminus\{0\}}\,\frac{\Vert\,\begin{pmatrix} (\cos \alpha -1)u_1-\sin\alpha\, u_2\\
\sin\alpha\, u_1+(\cos \alpha -1)u_2\end{pmatrix}
\,\Vert}{\sqrt{u_1^2+u_2^2}}=\sqrt{2(1-\cos\alpha)}\  .\qedhere
\end{equation*}
\end{proof}

\subsection{Proof of proposition \ref{-12-CMS}} \label{nicefeature}
$ $

For the notations below, we refer to section \ref{6.1}. In lemma \ref{2a}(\ref{2a-bis33}), it is shown that in a regular refined simplicial complex $K_E$ of a given finite simplicial complex $K$ ($E$ is paired with $\rho$), each simplex has thickness $\geq\! {\bf t}_n\!>\!0$ and diameter bounded in terms of $\rho$ (above and below).
\par\vskip1mm 
Proposition {\rm\ref{-12-CMS}} holds for any $E$ and any $n$-simplex $\sigma\!\in\! K_E$ considered for itself, actually one has the
\begin{lem}\label{lem1033} There exist constants $C_1,C_3,C_2\!=\!1+2C_3$ depending on $T,K,D\!\subset\!(M,g)$ such that, 
given $E$ (thus $\rho$) and any $n$-simplex $\sigma\!\in\! K_E$, if $r/\rho\!\leq \!C_1$ and $x\!\in\!\hat\sigma\!=\!T(\sigma)$ is at $g$-distance $\!\leq\! r$ from each of two subsimplices (of any dimensions) $\hat\eta_1, \hat\eta_2$ of $\hat\sigma\,$, then 
\par $(i)$ $\hat\eta_1\cap\hat\eta_2\not\!=\!\emptyset$ (thus $\eta_1\cap\eta_2\not\!=\!\emptyset$), $(ii)$ $d_g(x,\hat\eta_1\cap\hat\eta_2)\!\leq\! C_2\,r$.
\par $(iii)$ If $p_1,p_2$ are closest to $x$ in $\hat\eta_1, \hat\eta_2$, some $q\!\in\!\hat\eta_1\cap\hat\eta_2$ verifies
\begin{equation*}d_g(x,q)\!\leq \!\min_{i=1,2}d_g(x,p_i)+d_g(p_1,q)+d_g(p_2,q)\!\!\leq\!  \min_{i=1,2}d_g(x,p_i)+C_3 \,d_g(p_1,p_2)\,.
\end{equation*}
\end{lem}
To prove $(ii),(iii)$, we shall need another
\begin{lem}\label{lem1044} There exists $C_3\!\geq\!1$ such that, for any intersecting simplices $\eta_1,\eta_2\!\subset\!\sigma\!\in \!K_E$, for any $p_1\!\in\!\hat\eta_1,p_2\!\in\!\hat\eta_2\,$, one can find $q\!\in\!\hat\eta\!:=\!\hat\eta_1\!\cap\!\hat\eta_2\!\not=\!\emptyset\,$ with 
\begin{equation*}d_{g\mid\hat\sigma}(p_1,q)+d_{g\mid\hat\sigma}(q,p_2)= d_{g\mid\hat\eta_1\cup \hat\eta_2}(p_1,p_2)\leq C_3\,d_{g\mid\hat\sigma}(p_1,p_2)\ .
\end{equation*}
\end{lem} 
Given ${\bf t}_n\!>\!0$, let ${\mathfrak K}$ be the set of {\em Euclidean} $n$-simplices of thickness $\!\geq\!{\bf t}_n$ and diameter $\!=\!1\,$. Up to isometries, ${\mathfrak K}$ is {\it compact}. Below  the distance $d$ is Euclidean or Riemannian ($d=d_g$) according to the context.

\begin{proof} We first prove lemma \ref{lem1044} when $\sigma\!\in\!{\mathfrak K}$.  
\begin{proof} The {\em Euclidean} sub-simplices $\eta_1, \eta_2$ having intersection $\eta\not\!=\!\emptyset$ generate a sub-simplex 
$\sigma'\!\subset\!\sigma\!\in\!{\mathfrak K}$, distinct from $\eta_1$ and $\eta_2$ (if this is not the case, the result holds). The union $\eta_1\cup \eta_2$ is a connected metric subspace of $\partial \sigma'$. There exists a constant $C_3^{\sigma}\!\geq\!1$ such that, for every sub-simplex $\sigma'\subset\sigma$ and intersecting pair $\eta_1, \eta_2$ of sub-simplices in $\sigma'$, one has, for $q\!\in\!\eta_1\cap\eta_2$ well-chosen, $d_{\eta_1\cup \eta_2}(p_1,p_2)\!=\!d(p_1,q)\!+\!d(q,p_2)\!\leq\! C_3^{\sigma}\,d_{\sigma'}(p_1,p_2)\!=\!C_3^{\sigma}\,d_{\sigma}(p_1,p_2)$. So, one gets a positive $C_3^{\mathfrak K}:=\max_{\sigma\in{\mathfrak K}}C_3^{\sigma}$ such that lemma \ref{lem1044} holds for all $\sigma$ in the {\it compact} ${\mathfrak K}$
with $C_3=C_3^{\mathfrak K}$. 
\end{proof} 
Below, when $\sigma\!\in\!{\mathfrak K}$ and $\rho\!=\!1$, the proof of lemma \ref{lem1033}.
\begin{proof}
Two non intersecting {\em Euclidean} subsimplices $\eta_1, \eta_2$ in a nondegenerate $n$-simplex $\sigma\!\in\!{\mathfrak K}$ are at distance $d(\eta_1, \eta_2)>0$, actually the minimum $m_\sigma$ of those finitely many distances, taken over all pairs of non intersecting sub-simplices in $\sigma$, is $>0$. Moreover, 
 $m_{\mathfrak K}=\inf_{\{\sigma\in{\mathfrak K}\}}\,m_\sigma\,$ must be $>0$. This shows that there exists $C_1^{\mathfrak K}$ such that, if $r\leq C_1^{\mathfrak K}$, the simplices $\eta_1, \eta_2$ must intersect, proving $(i)$ when $\sigma\!\in\!{\mathfrak K}$. 
\par To prove $(ii),(iii)$ of lemma \ref{lem1033} as $\sigma\!\in\!{\mathfrak K}$, take intersecting distinct $\eta_1,\,\eta_2\,$ (nothing to prove if $\eta_1\!\subset\!\eta_2$). If $p_1$ and $p_2$ are closest to $x$ in $\eta_1$ and $\eta_2$ and $q\in \eta=\eta_1\cap\eta_2\not=\emptyset$ is such that $d(p_1,q)+d(q,p_2)=d_{\eta_1\cup \eta_2}(p_1,p_2)\,$, one has (apply lemma \ref{lem1044})
\begin{multline}\label{lem1077}d(x,q)\!\leq \!\min_{i=1,2}(d(x,p_i)+d(p_i,q))\leq\\\!\leq \!\min_{i=1,2}d(x,p_i)+d(p_1,q)+d(p_2,q)\leq\!  \min_{i=1,2}d(x,p_i)+C_3^{\mathfrak K} d(p_1,p_2)\!\!\leq \!\!(1\!+\!2\,C_3^{\mathfrak K})r\, ,
\end{multline}
thus (if $\sigma\!\in\!{\mathfrak K}$) $(iii)$, then $(ii)$ (set $C_2^{\mathfrak K}=1+2\,C_3^{\mathfrak K}$).
\end{proof}
Lemmas {\rm\ref{lem1033}}, {\rm\ref{lem1044}} now hold for $\sigma\!\in\!{\mathfrak K}$. Viewed as stating relations among ratios of distances, if the hypothesis $r\!\leq\!C_1^{\mathfrak K}$ is replaced by $r\!\leq\!C_1^{\mathfrak K}\,{\rm diam}(\sigma)$, they hold for any Euclidean simplex rescaled from $\sigma$. For any Euclidean simplex of a regular subdivision of integer $E$, the diameter is bounded from below in terms of $E$ and thus of $\rho$ (see lemma \ref{2a-bis}(\ref{2a-bis33}) and (\ref{2l})), i. e. there exists a real $C>0$ such that $C\,\rho\!\leq\!{\rm diam}(\sigma)$, so these lemmas hold for any such Euclidean simplex and bounds $\tilde C_1^{\mathfrak K}\!=\!C\,C_1^{\mathfrak K},C_2^{\mathfrak K},C_3^{\mathfrak K}$, assuming this time $r\!\leq\!\tilde C_1^{\mathfrak K}\rho$.
\par\vskip1mm
{\em In the general case,} both lemmas follow from the previous study and, for $E$ large enough, from the uniform closeness of $g$ and $g^0=(T^{-1})^\ast g_0$ (use the $\rho^2$-quasi-isometry of definition \ref{convpolyhedra} $(ii)$): lemmas {\rm\ref{lem1033}}, {\rm\ref{lem1044}} hold in the $g$-distance as well as in the $g^0$-distance, with $C_1,C_3,C_2\,=\,1\,+\,2C_3$ depending on $T,K,D\!\subset\!(M,g)$ close to $\tilde C_1^{\mathfrak K},C_3^{\mathfrak K},C_2^{\mathfrak K}\,=\,1\,+\,2C_3^{\mathfrak K}$ and replacing them in the extended statements.
\end{proof}
\begin{rem}\label{chaud} Since $C_3$ is $\geq1$, one has $C_2\!=\!1+2C_3\!\geq\! 3\,$.
\end{rem}
{\bf Proof of the proposition \ref{-12-CMS}}$ $

First, we define constants ${\mathcal C}_1,{\mathcal C}_2$ in terms of the above $C_1,C_2,C_3$, for which proposition \ref{-12-CMS}
will be shown to hold.
\begin{nota} \label{lem1011} Recall the real ${\mathcal B}$ of proposition {\rm \ref{bds4}} (verifying  $2{\mathcal B}\!>\!1$) and call 
${\mathfrak b}$ the largest integer $<\!2{\mathcal B}$. Define ${\mathfrak c}_1:={C_1}/{2(\sum_{l=0}^{n-1} C_3^l)}$,
\begin{equation*}{\mathcal C}_1:=\frac{1}{2}\min(3,\frac{C_1}{2(1\!+\!C_3)^{{\mathfrak b}-1}\sum_{l=0}^{n-1} C_3^l}) \ ,\ \  {\mathcal C}_2:=2\,(1\!+\!C_3)^{{\mathfrak b}-1}(\sum_{l=0}^{\!\!n-1} C_3^l)\  .
\end{equation*}
\end{nota}

\begin{proof} Though we deal here with the general case, we shall drop $\hat\sigma,\hat\eta, \dots$ and write instead $\sigma,\eta,\dots$.
\par In general, $\eta_1$ and $\eta_2$ do not belong to the same $n$-simplex. Let $x$ be at distance $\leq r$ from $\eta_1$ and $\eta_2$ with $r/\rho\leq {\mathcal C}_1$, notation \ref{lem1011} implies $r\!< \!2\rho\,$. Thus {\it the $g$-shortest path} $\gamma$ joining $x_1$ to $x_2\,$, closest points to $x$ in $\eta_1,\eta_2$, being of length $\leq\!2 r\!\leq\! 2{\mathcal C}_1\rho\!<\!4\rho\,$, encounters at most $N$ {\it distinct} simplices of dimension $n$, where $N\!\leq\!{\mathfrak b}\!<\!2{\mathcal B}$ (see proposition {\rm \ref{bds4}}), some of which may be visited quite many times, in the Riemannian case. In the Euclidean case, one can build a $2$-dimensional polyhedron in which a minimising geodesic has two connected components in some $2$-simplex, hitting all its edges. So, we make a slightly deeper observation to 
handle the {\it succession of visits} of $n$-simplices by $\gamma$.
\begin{lem}\label{lem1066} If $r\!\leq\! {\mathfrak c}_1\,\rho$, in a given $n$-simplex $\sigma\!\in\! K_E\,$, all $(n\!-\!1)$-faces crossed by $\gamma$ (if any) together intersect along a simplex $\mu\not\!=\!\emptyset$ (thus, those faces are less than $n$). One can find a point $q\in\mu$ such that 
\begin{equation}\label{lem1088}\forall x'\in\gamma \ \ \ \ \ \ \ \ \ d(x',q)\leq C'_3\,r\ \ \ \hbox{with}\ \ C'_3=2(\sum_{l=0}^{n-1} C_3^l).
\end{equation}
\end{lem}
\begin{proof} The fact that $(n\!-\!1)$-faces having nonempty intersection are at most $n$ is true in any $n$-simplex $\sigma\,$. Any two $(n\!-\!1)$-faces $\lambda$ and $\lambda'$ crossed by $\gamma$ in $\sigma$ at points $p$ and $p'$ are at distance $\leq \!\hbox{\rm length}\,\gamma\!\leq\!2 r\!\leq\!2{\mathfrak c}_1\rho\!\leq\!C_1\rho\,$, and applying lemma \ref{lem1033}, one has $\lambda\cap\lambda'\!\not=\!\emptyset$ and finds $q\in\lambda\cap\lambda'$ such that  $d(p,q)$ and $d(p',q)\!\leq \!C_3\,d(p,p')\leq 2C_3r$. Thus, for any $x'\!\in\!\gamma$, $d(x',q)\!\leq\! 2(1 + C_3)r\!\leq\!C_1\rho$ (use $r\!\leq\!{\mathfrak c}_1\rho$ and \ref{lem1011}). If there is another $(n\!-\!1)$-face $\lambda''$ crossed by $\gamma$ at $p''\,$, again by lemma \ref{lem1033}, applied to $q\!\in\!\lambda\cap\lambda'$ and $p''\!\in\!\lambda''$ (as $p''$ is on $\gamma$, one has $d(p'',q)\!\leq\!C_1\rho$ by the previous analysis), one gets $\lambda\cap\lambda'\cap\lambda''\!\not=\!\emptyset$ and, applying lemma \ref{lem1033} $(iii)$ to choose $q'\!\in\!\lambda\cap\lambda'\cap\lambda''$ nearest to $p''$ and $q\,$, one also has
\begin{multline*}\forall x'\!\in\!\gamma\ \ \ d(x',q')\!\leq \!d(x',p'')\!+\!d(p'',q')\!\leq\! d(x',p'')\!+\!d(p'',q')\!+\!d(q',q)\!\leq\\\leq d(x',p'')+C_3\,d(p'',q)\leq 2(1+C_3(1 +C_3))\,r\leq C_1\,\rho\ ,
\end{multline*}
writing $d(x',p'')\leq 2\,r\,$, $d(p'',q)\!\leq\!2(1 +C_3)\,r$ and recalling \ref{lem1011} for ${\mathfrak c}_1$.
Exhausting all $(n\!-\!1)$-faces crossed by $\gamma$ in $\sigma\,$, one gets the claim.
\end{proof}
Enumerate all {\it distinct} $n$-simplices visited by their order of first appearance along $\gamma$ (keeping the index for any later visit), $\sigma_1,\sigma_2,\dots,\sigma_N$. Then, in each $\sigma_i$ consider  $\mu_i\!\subset\!\sigma_i$ defined according to the above lemma \ref{lem1066}: there exists $q_i\!\in\!\mu_i$ such that, for any $x'\!\in\!\gamma\,$, one has $d(q_i,x')\!\leq \!C'_3r\!\leq \!C'_3 {\mathcal C}_1\rho\!\leq \!C_1\rho$ (notation \ref{lem1011}, (\ref{lem1088}) and the choice $r\!\leq\!{\mathcal C}_1\rho$). The next step is to apply lemma \ref{lem1033} to the sequence $\mu_i$. 
\par The {\it distinct} $\mu_1$ and $\mu_2$ (if they are not distinct, proceed to the first distinct pair $\mu_{i-1},\mu_i$ for $i\!\geq\!3\,$; if there is no such pair, the proposition holds) are both contained in a common $(n\!-\!1)$-face $\sigma_1\cap\sigma_2\,$, so they are contained in $\sigma_2$ (in $\sigma_i$ if $\mu_{i-1},\mu_i$ is the first distinct pair)). Applying lemma \ref{lem1033} (with $C'_3r$ replacing $r$, which works since $C'_3r\!\leq\!C_1\rho$), they have a nonempty intersection $\mu^1_2\!\subset\!\sigma_2\,$ and one finds $q^1_2\in\mu^1_2$ such that, for any $x'\!\in\!\gamma\,$, $d(x',q^1_2)\!\leq\! d(x',q_2)\!+\!d(q_2,q^1_2)\!\leq\! (1\!+\!2C_3)C'_3r$ (by lemma \ref{lem1044}, $d(q_2,q^1_2)\!\leq\!C_3 d(q_1,q_2)\!\leq\!2C_3C'_3r$ since $ d(q_1,q_2)\!\leq\!2C'_3r$). But the new sequence $\mu^1_2,\mu^1_3\!\!=\!\!\mu_3,\mu^1_4\!\!=\!\!\mu_4\dots$ shares with $\mu_1,\mu_2,\mu_3,\mu_4\dots$ similar properties. For $ i\!=\!3,4\dots$, one has $q^1_i:=q_i\!\in\!\mu^1_i$ such that, for any $x'\!\in\!\gamma\,$, $d(q^1_i,x')\!\leq\!C'_3r$ and $d(x',q^1_2)\!\leq\! (1\!+\!2C_3)C'_3r$. {\it The win} is that the distinct $\mu^1_i$ for $i\!=\!2,3\dots$ are at most $N\!-\!1$. After at most $N\!-\!1$ similar operations (applying lemma \ref{lem1033} goes on since $r\!\leq\!{\mathcal C}_1\rho$, see \ref{lem1011} for ${\mathcal C}_1$), one gets a {\em nonempty intersection of all $\mu_i$, thus of $\eta_1\,$ and $\eta_2\,$}.
\vskip1mm 
\par For $i\!=\!1,\dots,N\!-\!1$, set $\mu_i^{i-1}\!\!:=\!\mu_{i-1}^{i-2}\cap\mu_i\!\!=\!(\mu_1\cap\dots\cap\mu_{i-1})\cap\mu_i$. As above, $\mu_i^{i-1}\!\subset\!\sigma_i$ holds. Get $q_i^{i-1}\!\in\!\mu_i^{i-1}$ verifying for any
$x'\!\in\!\gamma$ 
$d(x',q_i^{i-1})\!\leq\! \alpha_i\,C'_3r$ where $\alpha_i\!\leq\!(1\!+\!C_3)\alpha_{i-1}\!+\!C_3$ since (by lemma \ref{lem1044})
\begin{multline} d(x',q_i^{i-1})\!\leq\! d(x',q_{i-1}^{i-2})+d(q_{i-1}^{i-2},q_i^{i-1})\leq\\
\leq d(x',q_{i-1}^{i-2})+d(q_{i-1}^{i-2},q_i^{i-1})+d(q_{i}^{i-1},q_i)\leq\\
\leq d(x',q_{i-1}^{i-2})+C_3\,d(q_{i-1}^{i-2},q_i)\leq d(x',q_{i-1}^{i-2})+C_3\,(d(q_{i-1}^{i-2},x')+d(x',q_i))\ .
\end{multline}
For $i\!=\!1,2\dots$, set $\!A_1\!:=\!1,A_2\!:=\!1\!+\!2C_3, \dots, A_i(a,b)\!:=\!aA_{i-1}(a,b)\!+\!b,$ $\,\dots\ .$  Inserting
$\alpha_i\!=\!A_i(1+C_3,C_3)\!=\!2(1\!+\!C_3)^{i-1}\!-\!1$, one gets $d(x,q_i^{i-1})\!\leq\!2\,(1\!+\!C_3)^{N-1}C'_3r$ for $i\!=\!1,\dots,N\!-\!1$. 
Bringing in ${\mathcal C}_2$ (see notation \ref{lem1011}), the last claim of proposition \ref{-12-CMS} now follows. 
\end{proof}
\begin{rem}\label{ohnon} Proving lemmas \ref{lem1033}, \ref{lem1044}, we saw that $C_1^{\mathfrak K},C_2^{\mathfrak K},C_3^{\mathfrak K}$ depend only on $n$ and on a lower bound $t_0$ of the thickness for all the $n$-simplices of $K$. In notation \ref{lem1011}, the constants ${\mathcal C}_1,{\mathcal C}_2$, when defined directly from $C_1^{\mathfrak K},C_2^{\mathfrak K},C_3^{\mathfrak K}$ as they should be, dealing with a piecewise flat polyhedron $K$, depend on $n$, on a lower bound $t_0\!>\!0$ of the thickness for all the $n$-simplices of $K$ and on something which should play the role of ${\mathcal B}$. The role of $2{\mathcal B}$, in the ``Riemannian case'', in the argument of the proof of lemma \ref{lem1066} completing the proof of proposition  \ref{-12-CMS}, is only to bound from above the number of $n$-simplices laying sufficiently near any point of $K$ by $2{\mathcal B}$. So, if $K$ is a piecewise flat polyhedron such that, say, 
$2{\mathcal B}_0$ bounds the number of $n$-simplices which lay sufficiently near any point of $K$, proposition \ref{-12-CMS}
and corollary \ref{-22-CMS} are true for constants ${\mathcal C}_1,{\mathcal C}_2,{\bf C}_1,{\bf C}_2$ depending on $n$, on a lower bound $t_0\!>\!0$ of the thickness for all the $n$-simplices of $K$ and on this ${\mathcal B}_0$.
\end{rem}

\vfill\eject

\part{Some facts and results around singularity theory}$ $\label{annexe1033}
In the sequel, transversality and Thom's theorem as developped in \cite{C-M} come out to produce the needed generic local bounds.
\section{Textures, multiplicities and local degrees of intersection}
\label{multiple1088}$ $

\par
The material presented in this section \ref{multiple1088} is developped in \cite{C-M}. See \cite{B-L} and \cite{G-G} for an introduction to this field.
\begin{defn}\label{Texturedefn} (see \cite{C-M}, §{\rm 7}) Given an $n$-dimensional manifold $W$, an integer $k$ $(0\leq k\leq n)$, a $k$-{\em texture of corank $\kappa$} on $W$ is a pair $({\mathcal L},p)$ consisting of
\par\noindent -\! a smooth $k$-dimensional foliation ${\mathcal L}$ of a $(\kappa\!+k)$-manifold $X$;
\par\noindent -\! a proper smooth map $p:X\rightarrow W$ such that, for each leaf $L$ of ${\mathcal L}$, if $\iota_L$ is the inclusion map in $X$, which is an injective immersion (while any {\em embedding} is proper), the composed map $p\circ \iota_L$ is an immersion.
\par The manifold $X$ is the {\em total space} of the texture and $W$ its {\em basis}.
\par
{\em If the texture verifies that each $i_L\!:=\!p\circ\iota_L$ is an embedding of $L$} ($p_{\mid L}$ is one-to-one), we call it {\em field fresh} and write $L$ meaning $i_L(L)$.
\end{defn}
\par \begin{rem}\label{commeonfaitsonlit} This definition, taken from \cite{C-M}, is quite general; here {\em we only deal with field fresh textures}, built from a given local family of submanifolds $L$ of $W$ properly embedded through canonical embeddings $i_L\,$. The paradigm of a {\em field fresh} texture is the family of lines $L$ in the Euclidean plane $W$ which is ``unfolded'' by viewing each line as leaf $L$ of a foliation in the unit tangent bundle $X$ to $W\,$. If one is interested in the geodesics of a flat $2$-dimensional torus, one can still describe them as building a texture, which this time is no more field fresh. However, one who is interested in the local flow of the geodesics near a given point $p$ in the torus may consider the field fresh texture having basis a convex neighborhood $U$ around $p$, in which any connected geodesic arc through $U$ is thought for itself, despite of its recurrence as a global geodesic in the torus.
\par Actually, for given $k$ and $t_0\!>\!0\,$, the texture of interest to us (see definition \ref{texture3333}, theorem \ref{texture333}) describes the family of $k$-faces of all small Riemannian barycentric $n$-simplices of thickness $\!\geq\! t_0\,$. This family is {\em ``unfolded''} in a total space $X$, where the $k$-faces slip into $k$-disks $i_L(L)\equiv L$ embedded through $\iota_L:L\rightarrow X$ as {\em disjointed} leaves $L$ of a foliation of $X$, building {\em a field fresh texture ${\mathcal L}$.}
\end{rem}
\par
\begin{defn}\label{multiple1044}  (adapted from \cite{C-M} §2, §7 lemma 7.2.).
$ $
\begin{itemize}
\item
The {\em multiplicity} of $h\!\in\! C^\infty(L, {\mathbb R}^k)$ {\em at} $x\!\in \!L$, $L$ embedded $k$-manifold in $W$, is $m_x(h)\!=\!{\rm dim}\,C_x^\infty(L)/h^\ast {\mathcal M}_{h(x)}$, with ${\mathcal M}_{h(x)}$ maximal ideal in $C_{h(x)}^\infty({\mathbb R}^k)$. 
\item  Switching to $m$-jets in $J_{x}^m(L,{\mathbb R}^k)\!:=\!C_x^\infty(L,{\mathbb R}^k)/ {\mathcal M}_x^{m+1}$, the {\em multiplicity} of $j^mh(x)$ is $m_x(j^mh)\!=\!{\rm dim}\,C_x^\infty(L)/(h^\ast {\mathcal M}_{h(x)}\!+\!{\mathcal M}_x^{m+1})$.
\par Both multiplicities coincide if one is $\leq m$ (see §2, proposition 2.1 in \cite{C-M}), enabling a polynomial computation.
\item Given an $({n\!-\!k})$-submanifold $N\subset W\,$, {\em define}
$$A_N:=\{g\in C^\infty(W, {\mathbb R}^k)\mid N\subset g^{-1}\{0\} \ \hbox{and}\  0 \ \hbox{regular value of}\ g\}\ .
$$ 
\item Given a field fresh $k$-texture of basis $W$, given  $L\!\in\!{\mathcal L}$, $\,x\!\in \! L\!\subset W\,$, the {\em multiplicity} of $N$ versus $L$ at $x\!\in\! W$ is 
$$m_x(N,L):=\inf_{g\in A_N} m_x(g,L)\ \ \hbox{where}\ \ 
m_x(g,L):=m_x(g_{\mid L})\,.
$$
\item The {\em multiplicity} of $N$ with respect to ${\mathcal L}$ at $x$ is then
$$m_x(N,{\mathcal L}):=\inf_{g\in A_N}m_x(g,{\mathcal L})\ \ \hbox{where}\ \ 
m_x(g,{\mathcal L}):=\sup_{\{L\in{\mathcal L}\,\mid\, x\in L\}}m_x(g,L)\,.
$$
\item Given $x\in  W\,$, {\em define} ${\mathcal L}_x\!=\!\{L\in{\mathcal L}\mid x\in L\equiv i_L(L)\subset W\}$.
\end{itemize}
\end{defn}
\begin{rem}  \label{M-T-C-M} As the texture is {\em field fresh}, the above multiplicity $m_x(g,{\mathcal L})$ is the same as the integer called $\mu_{{\mathcal L},p,x}(g^{-1})$ in lemma {\rm 7.2} of \cite{C-M}. {\em Indeed}, set $h\!=\!g\circ p\,$; as $p\circ \iota_L=i_L$ is an embedding, if ${\bf x}\in \iota_L(L)$ verifies $p({\bf x})=x$, then  $C_{\bf x}^\infty(L)/h^\ast {\mathcal M}_{h({\bf x})}=C_x^\infty(L)/g^\ast {\mathcal M}_{g(x)}$.
\end{rem}
\begin{rem} To have $m_x(N,{\mathcal L})\!=\!m\in{\mathbb N}$ {\em means} that one can find $g\in A_N$ such that $m_x(g,L)\!\leq\! m$ for any $L\in{\mathcal L}_x\,$
and $m_x(g,L_m)\!=\!m$ for some $L_m\in{\mathcal L}_x\,$; as soon as exists $g'\in A_N$ for which $m_x(g',{\mathcal L})<\infty\,$, one may find $g\in A_N$ such that $m_x(N,{\mathcal L})\!=\!m_x(g,{\mathcal L})$.
\end{rem}
\begin{defn}\label{intdegree1044} (see \cite{C-M}, §{\rm7}) With $N$ in $W$ and ${\mathcal L}$ as above, let ${\mathcal N}_x$ be a system of neighborhoods $U$ of $x$ in $W\,$. Given $U\in{\mathcal N}_x\,,\,L\in{\mathcal L}\,$, call ${\bf n}_{\,U,L,N}$ the number of elements in 
$U\cap L\cap N\,$.
The {\em local degree of intersection} of $N$ with respect to  ${\mathcal L}$ at $x$ is 
$${\rm deg}_x(N,{\mathcal L})\!:=\!\inf_{U\in{\mathcal N}_x}\sup_{L\in{\mathcal L}_x}\ {\bf n}_{\,U,L,N}\ .
$$
\end{defn}
\begin{lem}\label{semicontmult} Given $\bar g\!\in\! C^\infty({\mathscr T}\times W,{\mathbb R}^k)$ (${\mathscr T}$ manifold), $m_x(g_{\vartheta},L)$, $L\!\in\!{\mathcal L}_x$ with $\bar g\!=\!\{g_\vartheta\!\mid \!\vartheta \!\in\!{\mathscr T}\}$ verify: if $m_{x_0}(g_{\vartheta_0},L)\!<\!\infty$, there exists an open set $\theta\!\times\! U$ around $(\vartheta_0,x_0)\!\in\!{\mathscr T}\!\times \!W$ such that, for any $b\!\in\!\R^k,\vartheta\!\in\!\theta\,$, the set $ g_\vartheta^{-1}\{b\}\!\cap\! L\!\cap\! U$ is a finite set $\{x_1,\dots,x_r\}$ and one also has
\begin{equation*}\label{polymult}{\bf n}_{\,U,L,g_\vartheta^{-1}\{b\}}\leq
\sum_{s=1}^rm_{x_s}(g_\vartheta,L)\leq m_{x_0}(g_{\vartheta_0},L)\ .
\end{equation*}
\end{lem}
\begin{proof} Paste together the proofs of proposition 2.4 on page 168 \cite{G-G} (this is proposition 2.2 in \cite{C-M}) and of corollary 2.3 in \cite{C-M}.
\end{proof}
\begin{lem}\label{degleqmult} (see  lemma {\rm 7.2} in \cite{C-M}) One has 
$$\hbox{\rm deg}_x(N,{\mathcal L})\leq m_x(N,{\mathcal L})\,.
$$
\end{lem}
\begin{proof} Skipping ${\mathscr T}$, apply lemma \ref{semicontmult} to any $g\in A_N\,$.
\end{proof}

\begin{rem}\label{codimension1066} In the set $J_{0,0}^m({\mathbb R}^k,{\mathbb R}^k)$ of jets $j^mf(0)$ of $C^\infty$-maps $f$ such that $f(0)\!=\!0$, define the set $\Sigma^m(k)$ of jets having multiplicity $>m$ (see definition \ref{multiple1044}). By Tougeron's theorem, this set is algebraic of codimension $c_k(m)\rightarrow\infty$ as $m$ tends to $\infty$ (see \cite{B-L} chapter {\rm 13}). So, $\Sigma^m(k)$ is stratified by a finite number of regular strata (which are submanifols of dimension lower than a maximal dimension).
\par{\em Define} the integer $\ m_k(l)\!:=\!\min \{m\mid c_k(m)>l\}$ 
(see \cite{C-M}, page {\rm 337}).
\end{rem}
Now, in the above notations, a byproduct of theorem 7.1, \cite{C-M}.
\begin{thm}\label{thm7.1C-Mbis} Given $n,k,\kappa\!\in\!{\mathbb N}^*$, a $c$-manifold ${\mathscr T}\!$, a {\em field fresh} $k$-texture of corank $\kappa$ on an $n$-manifold $W$, there exist open dense sets ${\mathcal W}_{0,c},{\mathcal W}_c$ of $\,\bar g\in C^\infty({\mathscr T}\!\times \!W,{\mathbb R}^k)$ such that (set $g_\vartheta(\cdot)\!:=\!\bar g(\vartheta,\cdot)$)
\par\noindent $(i)$ for any $\bar g\in {\mathcal W}_{0,c}\,$, $\vartheta\in{\mathscr T}$, at any $x\in g_\vartheta^{-1}\{0\}$
$${\rm deg}_x(g_\vartheta^{-1}\{0\},{\mathcal L})\leq m_x(g_\vartheta,{\mathcal L})\leq m_k(\kappa+c)\,;
$$
\par\noindent $(ii)$ for any $\bar g\in {\mathcal W}_c\,,\vartheta\in{\mathscr T}\,,b\in {\mathbb R}^k$ and $x\in g_\vartheta^{-1}\{b\}$
$${\rm deg}_x(g_\vartheta^{-1}\{b\},{\mathcal L})\leq m_x(g_\vartheta,{\mathcal L})\leq m_k(\kappa+k+c)\,.
$$
\end{thm}
This theorem is the door to theorems \ref{text-11}, \ref{text-1133}.
\begin{proof} We also refer the reader to the establishment of theorem 7.1 in \cite{C-M}. Here $B\!=\!{\mathbb R}^k$. {\em Below}, we shall focus on the collection of maps (see remark \ref{M-T-C-M} and the proof of theorem 7.1, bottom of p. 351 in \cite{C-M}) $\bar g\circ \tilde p\circ \iota_{{\mathscr T}\!\times\!L}\!\equiv \! \bar g_{\mid {\mathscr T}\!\times\!L}$ of restrictions of $\bar g$ to all embedded  {\em new} leaves ${\mathscr T}\!\times\!L\!\subset\! {\mathscr T}\!\times\!W$ (${\mathscr T}\!\times\!L_{\bf x}$ is the new leave through $(\vartheta,{\bf x})$) of the {\em new texture} $\tilde {\mathcal L}\,$: this texture $\tilde {\mathcal L}$ has total space ${\mathscr T}\!\times\! X$, basis ${\mathscr T}\!\times\! W$ and map $\tilde p$ defined by $\tilde p(\vartheta,{\bf x})\!:=\!(\vartheta,p({\bf x}))\!:=\!(\vartheta,x)$ (since $p$ is proper, so is $\tilde p$).
 
\begin{fact} \label{onsecouche} For each integer $m$, the set $J^m(\tilde {\mathcal L},{\mathbb R}^k)$ of jets $j^m \tilde f(\vartheta,{\bf x})$ of maps $ \tilde f\!:\!({\mathscr T}\!\times\! L_{\bf x},(\vartheta,{\bf x}))\!\rightarrow \!{\mathbb R}^k$ is a smooth fibre bundle over ${\mathscr T}\!\times\! X\!\times\!{\mathbb R}^k$. 
\end{fact}
\begin{proof}
\par See lemma 7.3 in \cite{C-M} and proof: $J^m({\mathcal L},{\mathbb R}^k)$ is a locally trivial bundle over $ X\!\times\!\R^k$ with fibre diffeomorphic to $J_{0,0}^m({\mathbb R}^k,{\mathbb R}^k)$.
And $J^m(\tilde{\mathcal L},{\mathbb R}^k)$ is the pull-back of $J^m({\mathcal L},{\mathbb R}^k)$ by $(\vartheta,{\bf x})\!\in\!{\mathscr T}\!\times\!X\mapsto {\bf x}\!\in\!X$. 
\end{proof}
Then, 
$\Sigma^m(\tilde {\mathcal L},{\mathbb R}^k)$ is the subset of jets which project onto $\Sigma^m({\mathcal L},{\mathbb R}^k)$ under the {\it submersion} from $J^m(\tilde{\mathcal L},{\mathbb R}^k)$ onto $J^m({\mathcal L},{\mathbb R}^k)$
\begin{equation}\label{jamais}\Pi : j^m\tilde f(\vartheta,{\bf x})\!\in\! J^m(\tilde{\mathcal L},{\mathbb R}^k)\rightarrow j^m \tilde f_\vartheta({\bf x})\!\in\! J^m({\mathcal L},{\mathbb R}^k)\ ,
\end{equation}
while $\Sigma^m({\mathcal L},{\mathbb R}^k)$ is the subbundle of jets $j^mf({\bf x})$ of $f\!:\!(L_{\bf x},{\bf x})\!\rightarrow\! {\mathbb R}^k$ showing $m_{\bf x} (f,L_{\bf x})\!>\!m$ (with $L_{\bf x}$ leaf through ${\bf x}$), of fibre diffeomorphic to $\Sigma^m(k)$ (remark \ref{codimension1066}). 
Moreover (and as $\Sigma^m({\mathcal L},{\mathbb R}^k)$ in $J^m({\mathcal L},{\mathbb R}^k)$), the set
$\Sigma^m(\tilde {\mathcal L},{\mathbb R}^k)\subset J^m(\tilde{\mathcal L},{\mathbb R}^k)$ is a subbundle in which a change of charts induces a polynomial automorphism of the algebraic fibre $\Sigma^m(k)$ (see remark \ref{codimension1066} and the proof of lemma 7.3 in \cite{C-M}).

\par A map $h\!\in\! C^\infty(X,\R^k)$  gives rise to a $C^\infty$-map $[j^m](h)$
\begin{equation}\label{cedre}[j^m](h):{\bf x}\!\in\!X\mapsto j^m(h\circ \iota_{L_{\bf x}})({\bf x})\!\in\!J^m({\mathcal L},{\mathbb R}^k)\ .
\end{equation} 

\begin{lem}\label{8.2} Let $({\mathcal L},p)$ be a texture (definition {\rm\ref{Texturedefn}}), {\em not} assumed to be field fresh. The set of $g\!\in\!C^\infty(W,\R^k)$
such that $[j^m](g\circ p):{\bf x}\mapsto j^m(g\circ p\circ \iota_{L_{\bf x}})(\bf x)$ is transversal to $\Sigma^m({\mathcal L},{\mathbb R}^k)$ is open and dense.
\end{lem}
\begin{rem} \label{8.2!} For completeness and clarity, we chose to deal here with general transversality and did not follow \cite{C-M}, where those considerations only appear at the very end in lemma 8.2.
\end{rem}
\begin{proof}
We first prove that $\Sigma^m({\mathcal L},{\mathbb R}^k)$
can be stratified according to definitions 3.4, 3.5 in  \cite{Fel}, i. e.  $\Sigma^m({\mathcal L},{\mathbb R}^k)=\bigcup S_j$ is a finite partition of relatively open smooth submanifolds $S_1,\dots,S_r$ of dimensions $s_1\!>\!s_2\!>\!\dots\!>\!s_r\!\geq\!0$ verifying a ``Whitney type condition $(a)$'', actually, in \cite{Fel} p.195, Feldman says that $S$ (here $S\!=\!(\Sigma^m({\mathcal L},{\mathbb R}^k)$) is {\em cohesive}:
\par {\em Given $j,k$ integers with $1\leq j<k\leq r$, if a sequence $y_i$ of points in $S_j$ converges to $y\in S_k$ and
$T_{y_i}S_j$ converges towards $T_\infty$, an $s_j$-dimensional subspace of $T_yJ^m({\mathcal L},\R^k)$, then $T_\infty\supset T_yS_k$.}

{\em Indeed}, first observe that such a stratification is a local feature and work (see lemma 7.3 in \cite{C-M} and its notations) in charts $(U_s\!\times\!\R^k,\varphi_s)$ covering $X\!\times\!\R^k$ and trivialising the bundle $\pi:J^m({\mathcal L},\R^k)\rightarrow X\!\times\!\R^k$. Each ``localised'' $\Sigma^m({\mathcal L},{\mathbb R}^k)\cap\,\pi^{-1}\!(U_s\!\times\!\R^k)\subset J^m({\mathcal L},\R^k)$ is {\em diffeomorphic} to the algebraic subset $\sigma:=\R^{\kappa+k}\!\times\!\R^k\!\times\! \Sigma^m(k)$ (see remark \ref{codimension1066}).
\par For any fixed $s$, one concludes from theorem 9.7.11 in \cite{B-C-R} that
\par
{\em The semi-algebraic set $\sigma\subset \R^d$ may be stratified into $N$ open regular disjoint submanifolds (which are semi-algebraic sets) $\sigma_l$ verifying, for any $l,l'$ integers between $1$ and $N$
\par\textbullet \hskip2mm$\sigma=\bigcup_l\sigma_l\,$;
\par\textbullet \hskip2mm if $l\not=l'$ and $\sigma_l\cap \bar \sigma_{l'}\not=\emptyset$, then $\sigma_l\subset \bar \sigma_{l'}$ and ${\rm dim}\,\sigma_l<{\rm dim}\,\sigma_{l'}\,$;
\par\textbullet\hskip2mm{\em (Whitney condition $(a)$)} if $z_i\!\in\!\sigma_{l'}$ is a sequence converging to a point $z\!\in\!\sigma_l$ and if $T_{z_i}\sigma_{l'}$ converges towards $T_\infty$, then $T_\infty\supset T_{z}\sigma_l$.
\par This stratification may be performed in a {\em canonical} way (this here important point is neatly exposed in \cite{G-W-P-L}, proposition {\rm 2.7 chapter {\rm I}}; in the neighboring context of the stratification of a map, see \cite{Ma2}, especially §{\rm 5}), implying that it is invariant by diffeomorphisms (see the proof of proposition {\rm 4.1} in \cite{C-M}), thus it piles up in a global stratification of $\Sigma^m({\mathcal L},{\mathbb R}^k)$ which fits with our first claimed goal in this proof.}
\vskip1mm
Indeed, for any fixed $s$, collect all $\sigma_l$ having the same dimension $s_k$ in a unique $S_{k,s}$, ordering  $S_{1,s},\dots,S_{r,s}$ from the top to the lowest existing dimension of $\sigma_l$. Then
{\em $\Sigma^m({\mathcal L},{\mathbb R}^k)$
satisfies the hypotheses of proposition {\rm 3.6} in \cite{Fel}} (we require here the case b) of this proposition), implying
\par
{\em The set of mappings
$\tau$ in $C^\infty(X,J({\mathcal L},\R^k))$ such that $\tau$ is transversal to $\Sigma^m({\mathcal L},{\mathbb R}^k)$ is open and dense}.
\par A jet $\sigma\!\in\!J^m(X,{\mathbb R}^k)$ is sent to the jet
$[\sigma]\!\in\!J^m({\mathcal L},{\mathbb R}^k)$ in the following well-defined way: choose a germ $f:(X,{\bf x})\rightarrow\R^k$ such that $\sigma=j^m(f)({\bf x})$ and {\em define} $[\cdot](\sigma)\!=\![\sigma]\!:=\!j^m(f\circ \iota_{L_{\bf x}})({\bf x})$, where $L_{\bf x}$ is the unique leaf in $X$ through ${\bf x}$. Doing a local computation in a foliated chart of $X$, one checks that $[\cdot]:\sigma\mapsto[\sigma]$ is a $C^\infty$-mapping.
\par So by propositions 3.4 and 3.5 page 46 in \cite{G-G}, the set of $h\!\in\!C^\infty(X,\R^k)$ such that $[j]^m(h)\!=\!([\cdot]\circ j^m)(h)$ is transversal to $\Sigma^m({\mathcal L},{\mathbb R}^k)$ is also open. As $p$ is proper, by proposition 3.9 p.49 in \cite{G-G}, the set ${\mathcal W}$ of $g$ in $C^\infty(W,\R^k)$ verifying that $[j^m](g\circ p)$ is transversal to $\Sigma^m({\mathcal L},{\mathbb R}^k)$ is open too.  
\vskip1mm
We sketch here a proof of the density in our context.

As usual, it results from lemma 4.6 in \cite{G-G} page 53. Choose a covering of $W$ by charts $(U_s,\varphi_s)$ and compact differentiable balls $\bar B_s\subset U_s$ verifying $W\!\subset\cup_s B_s$ and such that $\varphi_s(\bar B_s)\equiv \bar B_s$ can be thought as a compact ball in $\R^n$.
For any given $g\!\in\!C^\infty(\bar B_s,{\mathbb R}^k)$, the map $\Psi_m$
defined (for any polynomial map $b\!:\!\bar B_s\rightarrow {\mathbb R}^k$ of degree $\!\leq\! m$) by sending $(b,{\bf x})$ to $\Psi_m(b, {\bf x})\!:=\![j]^m((g\!+\!b)\circ p)({\bf x})\!\in\! J^m({\mathcal L},{\mathbb R}^k)\cap\,\pi^{-1}\bar B_s$ is {\em onto} ({\em adding a polynomial map} $b$ makes sense in $\bar B_s\subset\R^n$), thus transversal to $\Sigma^m({\mathcal L},{\mathbb R}^k)$. Indeed, $(b, {\bf x})\mapsto[j]^m(b\circ p)(b, {\bf x})$ is onto (do a local computation in a foliated chart of $X$) and one has $[j]^m((g+b)\circ p)\!=\![j]^m(g \circ p)+[j]^m(b\circ p)$. So, the set ${\mathcal W}_{s}$ of $g\!\in\!C^\infty(\bar B_s,{\mathbb R}^k)$ such that $[j^m](g\circ p)$ is transversal to $\Sigma^m({\mathcal L},{\mathbb R}^k)$ is dense. Applying lemma \ref{unlemme333} to the pair $\bar B_s,W$ to get from ${\mathcal W}_{s}$ a dense open ${\mathcal W}_{s}'\subset C^\infty(W,\R^k)$, one concludes that ${\mathcal W}\supset\cap_s{\mathcal W}_{s}'$ is dense in $C^\infty(W,\R^k)$.
\end{proof}

Built by a related procedure and stratified by a finite number of smooth submanifolds from the same algebraic set $\Sigma^m(k)$ of remark \ref{codimension1066}, $\Sigma^m({\mathcal L},{\mathbb R}^k)$ and $\Sigma^m(\tilde {\mathcal L},{\mathbb R}^k)$ are closed stratified submanifolds of the same codimension $c_k(m)$.
Set $m=m_k(\kappa+k+c)$ (see remark \ref{codimension1066}, $c_k(m)\rightarrow\infty$ as $m\rightarrow\infty$): the codimension of $\Sigma^m(\tilde{\mathcal L},{\mathbb R}^k)$ is $>\kappa+k+c$.

Since $\Sigma^m(\tilde {\mathcal L},{\mathbb R}^k)$ is the set of jets $j^m\tilde f(\vartheta,{\bf x})$ of $\tilde f$ verifying $m_{\bf x}(\tilde f_\vartheta,L_{\bf x})\!>\!m$ (with $L_{\bf x}$ leaf through ${\bf x}$), 
to say that there exists $x,\vartheta$ for which
$m_x(g_\vartheta,{\mathcal L})\!>\!m\,$ is the same as to say that
$\hbox{\rm rge}([j^m](\bar g\circ\tilde p))\cap\Sigma^m(\tilde{\mathcal L},{\mathbb R}^k)\not=\emptyset$
(see definition \ref{multiple1044}, remark \ref{M-T-C-M} and §2, §7 of \cite{C-M}): thus, by lemma \ref{8.2} applied to the texture $(\tilde{\mathcal L},\tilde p)$ (or by lemma 7.4 and the lines at the bottom of page 351 in \cite{C-M}), the set ${\mathcal W}_c$ of maps $\bar g\!\in\! C^\infty({\mathscr T}\!\times\! W,\R^k)$ which, for all $\vartheta\in{\mathscr T}$ and $x\!\in\! W$, satisfy $m_x(g_\vartheta,{\mathcal L})\!\leq\! m_k(\kappa+k+c)$  
 is {\em open and dense}.

\par Replace $\Sigma^m({\mathcal L},{\mathbb R}^k)$ by $\Sigma^m({\mathcal L},{\mathbb R}^k)_0\!:=\!\Sigma^m({\mathcal L},{\mathbb R}^k)\!\cap\!
\{j^m h({\bf x})\! \mid h({\bf x})\!=\!0\}$, which is by lemma 7.3 in \cite{C-M} a closed stratified manifold
of $J^m({\mathcal L},{\mathbb R}^k)$ of codimension $c_k(m)+k$. The same scheme as above applied to the set
$\Sigma^m(\tilde{\mathcal L},{\mathbb R}^k)_0$ of jets $j^m\tilde f(\vartheta,{\bf x})$ such that $j^mf_\vartheta({\bf x})\!\in\! \Sigma^m({\mathcal L},{\mathbb R}^k)_0$ provides with the expected open dense set ${\mathcal W}_{0,c}\!\subset \!C^\infty({\mathscr T}\times W,\R^k)$.
\par Then, in $(i)$ and $(ii)$, the inequalities on the left follow from definition \ref{intdegree1044} and lemma \ref{semicontmult} (see also lemma 7.2 in \cite{C-M}).
\end{proof}

\section{Adapting Thom's theorem in ``g\' eom\' etrie finie''}\label{adaptThom}$ $

To deal with intersections of embedded squares with faces of Riemannian barycentric $n$-simplices, beyond the classical theory, some results and consequences of a quite involved version of the Thom transversality theorem (see \cite{Th}) play a key role, see \cite{C-M}, this already came out in section \ref{RBT}.
\begin{rem}\label{Hirsch} If $X$ and $Y$ are manifolds, submersions, embeddings and proper maps constitute open sets in the Whitney $C^\infty$-topology of $C^\infty(X,Y)$, see \cite{G-G} chapter II §3, \cite{Hi}, chapter {\rm 2}, theorems {\rm 1.2}, {\rm 1.4}, {\rm 1.5}.
\end{rem}
For the convenience of the reader, we repeat definition \ref{carrement}.
\begin{defn} \label{carrementbis}Denote by $I$ a bounded open interval containing $[0,1]\,$. View $\bar I^2$ as $\Pi\!=\!\{\!(se_1\!+\!te_2\!)\!\in\! {\mathbb R}^2\mid (s,t)\!\in\!\bar I^2\}\,$, where $e_1,e_2$ is the canonical basis of ${\mathbb R}^2\,$. 
Given $w\in{\bf S}^1\subset{\mathbb R}^2\,$, the $2$-dimensional square
$\Pi$ may be thought as a $1$-{\it parameter} family of segments  $l$ of straight lines with $l\in\bar l_w\,$, where $\bar l_w$ is the set of segments $l\subset{\mathbb R}^2$ parallel to $w$ in $\Pi$, thus $\coprod_{l\in \bar l_w}l\!=\! \Pi\,$. Denote by $\bar l_1$ and $\bar l_2$ the sets $\bar l_{e_1}$ and $\bar l_{e_2}$ in $\Pi\,$.
\end{defn} 
\begin{defn} \label{squaredsq}We call {\em squared square} a data ${\bf G}^\boxplus\!=\!({\bf G}_1,{\bf G}_2,{\bf G})\!\,$, which, for given {\it embedded} parametrised square $G\in C^\infty(\bar I^2, W)$, reads
\begin{equation}\label{descriptsqsq}{\bf G}^\boxplus:=\!\{\big{(}G(l_1),\,G(l_2),\,\hbox{\rm rge}(G)\big{)} \mid l_1\in\bar l_1, \,l_2\in\bar l_2\}\ .
\end{equation} 
Observe that any ${\bf G}$ is embedded in some open ball $U$ of differentiable compact closure $\bar U\!\subset\! W\!\subset \!M$. 
\end{defn}
In order to apply below theorem \ref{thm7.1C-Mbis}, we present an embedded square ${\bf G}\!=\!\hbox{\rm rge}(G)$ as a subset of $g^{-1}\{0_{{\mathbb R}^{n-2}}\}$, i. e. $g\!\in\! A_{\bf G}\!\subset\!C^\infty(\bar U,{\mathbb R}^{n-2})$ (definition \ref{multiple1044}), doing a little more since we deal with {\em squared} squares.
\begin{defn}\label{dualsqu} 
Inverting $G$, set $(g_1,g_2)\!:=\!G^{-1}$. Extend it as an embedding ${\bf g}\!:=\!(g_1,g_2,g)\!\in\! C^\infty(\bar U, {\mathbb R}^n)$, with $g\!\in\!A_{\bf G}$ (this is done in a more precise setting in lemmas \ref{blue555} and \ref{blue555bister}). 
Set ${\mathpzc g}_1\!=\!(g_1,g),{\mathpzc g}_2\!=\!(g_2,g)$. Then ${\bf G}^\boxplus\!=\!({\bf G}_1,{\bf G}_2,{\bf G})$ in $U$ is described through ${\bf g}\!\in \!C^\infty(\bar U, {\mathbb R}^n)$
\begin{equation*}\label{Benalors}{\bf G}^\boxplus\subset
\{({\mathpzc g}_2^{-1}(\{y\}\times\{0_{{\mathbb R}^{n-2}}\}),{\mathpzc g}_1^{-1}(\{x\}\times\{0_{{\mathbb R}^{n-2}}\}),g^{-1}\{0_{{\mathbb R}^{n-2}}\})\mid x,y\in\bar I\}\ .
\end{equation*}
For given $w=(\alpha_1,\alpha_2)\in{\bf S}^1\subset{\mathbb R}^2$, set
$${\bf g}_w:=(\alpha_2\, g_1-\alpha_1\, g_2,g)\ \ \hbox{with}\ \ {\mathpzc g}_1={\bf g}_{e_2}\,,\,{\mathpzc g}_2={\bf g}_{-e_1}\,.
$$   
Any curve in the family (so ${\bf G}_1={\bf G}_{e_1}\,,\,{\bf G}_2={\bf G}_{e_2}$)
$${\bf G}_w:=\{{\bf G}_l:=G(l)\mid l\in\bar l_w\}
$$ 
verifies ${\bf G}_l\subset{\bf g}_w^{-1}\{(\alpha_2\,s_0-\alpha_1\,t_0,0_{{\mathbb R}^{n-2}})\}$ for some $(s_0,t_0)\in I^2\,$.
\par As usual, we write $\bar {\bf g}_w, \bar {\bf G}_w$ for $\{{\bf g}_{\vartheta,w}\mid\vartheta\!\in\!{\mathscr T}\}$, for $ \{{\bf G}_{\vartheta,w}\mid\vartheta\!\in\!{\mathscr T}\}\dots$ 
\end{defn}
\begin{defn}\label{topology1044} $ $ 
Define ${\mathscr O}_{\bar I^2}(W)\!:=\!\{G\!\in\!C^\infty(\bar I^2, W)\!\mid\! G \ \hbox{embedding}\}$. Define $\!{\mathscr O}^{\mathscr T}_{\bar I^2}(W)\!:=\!\{\bar G\!\in\!C^\infty({\mathscr T}\times\bar I^2, W)\!\mid\! \forall \vartheta\!\in\!{\mathscr T}\ \ G_\vartheta\ \,\hbox{is an embedding}\}$.
 These sets are open
(see remark {\rm\ref{Hirsch}}).
\par
Two {\em squared squares} (defined above) ${\bf G}_1^\boxplus,{\bf G}_2^\boxplus$ are {\em close} if they have a description through ({\rm\ref{descriptsqsq}}) by embeddings $G_1,G_2$ close in $C^\infty(\bar I^2, W)$. 
\end{defn}

\par In view of the localisations and globalisations to be done, a useful 
\begin{lem}\label{unlemme333} (also quoted as lemma {\rm\ref{unlemme}}) Let $M$ be an $n$-manifold, ${\mathscr T}$ a $c$-manifold, $B:=B^n$ the open unit ball in ${\mathbb R}^{n}$ and $\varphi$ a diffeomorphism from an open ${\mathcal O}\subset M$ onto an open of ${\mathbb R}^{n}$ that contains $\bar B$. Set $V=\varphi^{-1}B$. The restriction mapping
$$ {\mathcal R}_{\bar V}: \bar f\!\in\! C^\infty({\mathscr T}\times M,{\mathbb R}^l)\longmapsto 
{\mathcal R}_{\bar V}(\bar f)\!=\!\bar f_{\mid{\mathscr T}\times\bar V}\!\in\! C^\infty({\mathscr T}\times\bar V,{\mathbb R}^l)
$$
is onto, open and continuous in the Whitney topology. 
Thus, if ${\mathcal U}$ is open dense in $C^\infty({\mathscr T}\times M,{\mathbb R}^l)$, so is ${\mathcal R}_{\bar V}({\mathcal U})$ in $ C^\infty({\mathscr T}\times\bar V,{\mathbb R}^l)$. If ${\mathcal V}$ is open dense in $C^\infty({\mathscr T}\times\bar V,{\mathbb R}^l)$, so is $ {\mathcal R}_{\bar V}^{-1}({\mathcal V})$ in $C^\infty({\mathscr T}\times M,{\mathbb R}^l)$.
\end{lem} 
\begin{proof} See appendix \ref{unlemme333333}.
\end{proof}

\begin{lem} \label{blue555} Let $U$ be a smooth ball with $\bar U\!\subset\!W$. Define
$${\mathfrak O}_{\bar I^2}(U)\!:=\!\{{\bf g}\!\in\!C^\infty(\bar U,\R^n)\!\mid\!{\bf g}\ \hbox{embedding and} \ \bar I^2\!\times\!\{0_{{\mathbb R}^{n\!-\!2}}\}\!\subset\!{\bf g}(U)\}\ ,$$ which is an open set, 
and the corresponding ${\mathfrak O}^{\mathscr T}_{\bar I^2}(U)$.
\par
There exists an open, continuous map 
$\Phi$ sending $\bar{\bf g}\!\in\!{\mathfrak O}^{\mathscr T}_{\bar I^2}(U)$ to $\Phi(\bar{\bf g})\!=\!\bar G\!\in\!{\mathscr O}^{\mathscr T}_{\bar I^2}(U)$. 
\par If ${\mathscr T}$ is a {\em ball}, $\Phi$ is surjective onto ${\mathscr O}^{\mathscr T}_{\bar I^2}(U)$: thus, an open dense set in ${\mathfrak O}^{\mathscr T}_{\bar I^2}(U)$ produces an open dense set in ${\mathscr O}^{\mathscr T}_{\bar I^2}(U)$.
\end{lem}
\begin{proof} See appendix \ref{blue555555}.
\end{proof}

\begin{prop}\label{genericsqu} An open dense set ${\mathcal U}_c\!\subset \!C^\infty(\!{\mathscr T}\!\!\times\!\!W,\R^n\!)$ produces an open dense set of ${\mathscr T}\!\!$-families of squares $\bar G$ embedded in $W$ in ${\mathscr O}^{\mathscr T}_{\bar I^2}(W)$.
\end{prop}

\begin{proof} For any $\bar G\!\in\!{\mathscr O}^{\mathscr T}_{\bar I^2}(W)$ (see definition \ref{topology1044}), we build an open neighborhood ${\mathscr O}_{\bar G}\!\subset\!{\mathscr O}^{{\mathscr T}}_{\bar I^2}(W)$ of $\bar G$ and, induced by the given ${\mathcal U}_c\!\subset \!C^\infty(\!{\mathscr T}\!\!\times\!\!W,\R^n\!)$, a dense open subset ${\mathcal V}_{\bar G}\!\subset\!{\mathscr O}_{\bar G}$, see appendix \ref{genericsqu-333}.
\end{proof}

Now, theorem \ref{text-1133} and its proof (recall $\bar D\!\subset \!T(K)\!\subset\!W$).
\begin{thm}\label{text-1133bis}{\rm (companion to Thom's theorem {\rm\ref{text-11}}).} Denote by ${\mathscr T}$ a manifold, by $\Theta\!\subset\!{\mathscr T}$ a compact domain, both of dimension $c$. A polyhedral approximation $T,K$ of $D\!\subset\!(M,g)$ is given. In the open set ${\mathscr O}^{\mathscr T}_{\bar I^2}(W)$ (definition {\rm\ref{topology1044}}) exists a {\em residual} set ${\mathcal V}_c$ of ${\mathscr T}$-families $\bar G\!=\! \{G_\vartheta , \vartheta\!\in\!{\mathscr T}\}$ such that, for any $\bar G\!\in\!{\mathcal V}_c$,  one can find $\rho_{\bar G}>0$ (also depending on $\Theta$) verifying: 
\par for any $\rho\!\in]0,\rho_{\bar G}]$ exists a dense open set $\Theta_\rho$ of $\vartheta\!\in\!\Theta$ such that $G_\vartheta\!\in\!{\mathcal V}_c$ is testable for $\rho$ (see definition {\rm\ref{testsquare1}}). 
For $\bar G\!\in\! {\mathcal V}_c$, we say that {\em almost any $G_\vartheta\!\in\! \bar G$} (among $\vartheta$ in $\Theta$) {\em is testable for $\rho\,$}.
\end{thm}

\begin{rem} A residual set is dense in a Baire space and ${\mathscr O}^{\mathscr T}_{\bar I^2}(W)$ is Baire, see \cite{G-G} p.{\rm 44}, definition {\rm 3.2}, proposition {\rm 3.3}.
\end{rem}

\begin{rem} \label{bounded222bis} In view of definition {\rm\ref{testsquare1}}
\par 1. one even finds in ${\mathscr O}^{\mathscr T}_{\bar I^2}(W)$ an {\em open dense} of $\bar G$ satisfying only $(i)$ and $(ii)$: indeed, only properties involving squares $G$ {\em with respect to elements of} ${\mathfrak S}_{t_0}(W)$ (fact \ref{slip1033}) are then under consideration;
\par 2. lemma {\rm\ref{semicontmult}} together with the stated bounds $(i)$ and $(ii)$ on multiplicities directly imply bounds on {\em numbers of points of intersection}.
\end{rem}
\begin{proof} {\it We extend the proof of theorem {\rm  \ref{text-11}}}.
\par

We shall actually prove (see definitions \ref{squaredsq}, \ref{dualsqu}, \ref{topology1044}) that there exists in $C^\infty({\mathscr T}\!\times\!W,\R^n)$ an open dense set of $\bar{\bf g}$ such that ${\bf g}_\vartheta\!\in\!\bar{\bf g}$ verifies (for any $\vartheta\!\in\!\Theta$) {\em naturally paired} 
properties $(i)',(ii)'$ (stated below in this proof) with definition \ref{testsquare1}, $(i),(ii)$ ($\bar{\bf g}$ is linked to $\bar G$ through proposition \ref{genericsqu}, which produces open dense - residual - sets of $\bar G$ in ${\mathscr O}^{\mathscr T}_{\bar I^2}(W)$ from open dense - residual - sets of $\bar{\bf g}\!\in\!C^\infty({\mathscr T}\!\times\!W,\R^n)$). 
\vskip1mm
\par We use here some notations introduced at the beginning of section \ref{preparatorybds}:
the codimension $\kappa$ of the texture $X_i^{k}$ is $(k+1)n$. Towards $(i)$, we work with $X_i^{n-1}$ and $\kappa=n^2\,$, while, towards $(ii)$, with $X_i^{n-2}$ and $\kappa=(n-1)n$. 
From fact \ref{slip1033}, all small $\hat\sigma\!\in\!{\mathfrak S}_{t_0}(W)$ are subsets of $B(x_i,3\rho_3/4)$ for some $i$ and each $(n\!-\!1)$-face $\hat\eta$ (each $(n\!-\!2)$-face $\hat\xi$) of $\hat\sigma$ slips into a leaf of some $X_i^{n\!-\!1}$ (of some $X_i^{n\!-\!2}$). 
\par Fix the index $i$. 

\par\vskip1mm
We refer to theorem \ref{thm7.1C-Mbis} and proof for notions and notations below.

\par Consider the texture $\tilde{\mathcal L}$ with total space ${\mathscr T}\!\times\! X_i^{n\!-\!1}$, basis ${\mathscr T}\times B(x_i,\rho_3)$, whose leaves ${\mathscr T}\times L$ contain all ${\mathscr T}\times\hat\eta\,$.
Given $m$ integer, {\it define the closed stratified manifold $\tilde\Sigma_{0,n\!-\!2}^m\!:=\!\Sigma_{0,n\!-\!2}^m(\tilde{\mathcal L},{\mathbb R}^{n\!-\!1}\!)$ to be the subset of $m$-jets $j^m \!\tilde f\,(\vartheta,{\bf x})\!\in\!\Sigma^m(\tilde{\mathcal L},{\mathbb R}^{n\!-\!1})=\Pi^{-1}\Sigma^m({\mathcal L},{\mathbb R}^{n\!-\!1})$, at $(\vartheta,{\bf x})\!\in\!{\mathscr T}\!\times\! X_i^{n\!-\!1}$, such that $m_{\bf x}(\tilde f_\vartheta,L_{\bf x})\!>\!m$} (with $L_{\bf x}$ leaf through ${\bf x}$) {\em and} $j_{(\vartheta,{{\bf x}})}^0 \tilde f=\tilde f(\vartheta,{{\bf x}})\!\in\! {\mathbb R}\!\times\! \{0_{{\mathbb R}^{n\!-\!2}}\}\,$ (see (\ref{jamais}), fact \ref{onsecouche},
 remark \ref{codimension1066}). 
By its definition, 
$\tilde\Sigma_{0,n\!-\!2}^m$ has codimension $c_{n\!-\!1}(m)\!+\!n\!-\!2$ in $J^m(\tilde{\mathcal L},{\mathbb R}^{n\!-\!1}\!)$. \vskip1mm\par
\par Towards $(i)$, we state and prove the
\vspace{-2mm}
\begin{fact}\label{ifirstpart} There exists a dense open set ${\mathcal W}'_c\!\subset\!C^\infty({\mathscr T}\!\times\!M,\R^n)$ such that,
given $\bar{\bf g}\!\in\! {\mathcal W}'_c$, {\em one can find $\rho'_{\bar{\bf g}}\!>\!0$ for which} the property below is true {\em for all $\rho\!\in]0,\rho'_{\bar{\bf g}}]$, $\hat\sigma\!\in\!{\mathfrak S}_{t_0}\!(W)$ of diameter $\!\leq\! 2\rho$}
\vskip1mm
\par\noindent $(i)'$ (see definition {\rm\ref{testsquare1}} $(i)$ first \textbullet)\  {\em for any $w\in{\bf S}^1$, any $(\!n\!-\!1\!)$-face $\hat\eta$ of $\hat\sigma$, any $b\!\in\!\bar I\!\times\! \{0_{{\mathbb R}^{n\!-\!2}}\}$ and any $\vartheta\!\in\!\Theta$, one has}
$$\sum_{x\in \hat\eta\,\cap\,{\bf g}_{\vartheta,w}^{-1}\!\{b\}} m_x({\bf g}_{\vartheta,w},\hat\eta)\leq m_{n-1}(n^2+\!2\!+\!c)\, .
$$ 
Recall that $m_x({\bf g}_{\vartheta,w},\hat\eta)\!\geq\!m_x({\bf G}_{\vartheta,w},\hat\eta)$, see definitions {\rm\ref{dualsqu}, \ref{multiple1044}}.
\end{fact}

\begin{proof}
Given a map $\bar{\bf h}\!\in\! C^\infty(\!{\mathscr T}\!\times\! X_i^{n-1},{\mathbb R}^{n})$, define 
(see (\ref{cedre}))
$${\bf j}^m\bar{\bf h}:\!(\vartheta,{\bf x},w)\!\in\!{\mathscr T}\!\times\! X_i^{n-1}\!\times\!{\bf S}^1\mapsto\big{(}[j^m](\bar{\bf h}_w)(\vartheta, {\bf x}),w\big{)}\!\in\! J^m\!(\tilde{\mathcal L},{\mathbb R}^{n-1})\!\times\!{\bf S}^1,
$$
where $\tilde{\mathcal L}$ is the above texture with total space ${\mathscr T}\!\times\!  X_i^{n-1}$ and basis ${\mathscr T}\!\times\! B(x_i,\rho_3)$. Define the closed stratified manifold $\Upsilon^m:=\tilde\Sigma_{0,n\!-\!2}^m\!\times\!{\bf S}^1$. The codimension of $\Upsilon^m$ in 
$J^m\!(\tilde{\mathcal L},{\mathbb R}^{n-1})\!\times\!{\bf S}^1$ (equalling the one of $\tilde\Sigma_{0,n\!-\!2}^m$ in $J^m\!(\tilde{\mathcal L},{\mathbb R}^{n-1})$) is $c_{n\!-\!1}(m)\!+\!n\!-\!2$.

\par
To say that 
${\rm rge}({\bf j}^m\bar{\bf h})$ is transversal to $\tilde\Upsilon^m$ is equivalent to say that, for any $w\!\in\!{\bf S}^1$, the map $[j^m] (\bar{\bf h}_w)$ is transversal to $\tilde\Sigma_{0,n\!-\!2}^m\,$. 
\par Given $\bar{\bf g}\!\in\! C^\infty(\!{\mathscr T}\!\times\! B(x_i,\rho_3\!),{\mathbb R}^{n}\!)$ and setting $\bar {\bf h}:=\bar {\bf g}\circ\tilde p$, our concern is now focused on the mapping ${\bf j}^m\circ\tilde p^\ast$ below
$$\bar {\bf g}\!\in\!C^\infty(\!{\mathscr T}\!\times\! B(x_i,\rho_3\!),{\mathbb R}^n\!)\!\mapsto {\bf j}^m(\bar {\bf g}\circ\tilde p)\!\in\!C^\infty({\mathscr T}\!\times\!  X_i^{n-1}\!\times\!{\bf S}^1\!,J^m\!(\tilde{\mathcal L},{\mathbb R}^{n-1}\!)\!\times\!{\bf S}^1\!).
$$
 \par Leaving the developped proof of theorem \ref{thm7.1C-Mbis} (see remark \ref{8.2!}), we follow a direct route, reducing transversality to empty intersection. 
\par \noindent{\em Set} $m\!:=\!m_{n\!-\!1}(n^2\!+\!2\!+\!c)$ (see remark \ref{codimension1066}). As ${\rm dim}({\mathscr T}\!\times\! X_i^{n\!-\!1}\!\times\!{\bf S^1})$ equals $n(n\!+\!1)\!+\!c\,$, saying that $\bar {\bf g}\!\in\!C^\infty(\!{\mathscr T}\!\times\! B(x_i,\rho_3),{\mathbb R}^n\!)$ verifies that ${\bf j}^m\bar{\bf h}={\bf j}^m(\bar{\bf g}\circ\tilde p)$ is transversal to $\tilde\Upsilon^m$ {\em on a set} $S\!\subset\!{\mathscr T}\!\times\! B(x_i,\rho_3\!)\times\! {\bf S}^1$ is precisely to say that, {\em on this subset} $S$, $\bar {\bf g}$ verifies ${\rm rge}({\bf j}^m\bar {\bf h})\cap\tilde\Upsilon^m\!=\!\emptyset\,$. 
{\em Indeed}, as  {\em stratified manifold} (\cite{Ar2} page 217 or theorem \ref{thm7.1C-Mbis} and its proof), $\tilde\Upsilon^m$ is a {\em finite} union of submanifolds of codimensions $\!\geq\! c_{n\!-\!1}(m)\!+\!n\!-\!2$,  which are all $> n^2\!+\!n\!+\!c$ {\em in view of the definition of} $m_{n-1}$.

\par
For a $\bar {\bf g}\!\in\!C^\infty(\!{\mathscr T}\!\times\! B(x_i,\rho_3),{\mathbb R}^n\!)$ such that ${\bf j}^m\bar{\bf h}={\bf j}^m(\bar{\bf g}\circ\tilde p)$ is transversal to $\tilde\Upsilon^m$ {\em on a set} $S$, or, {\em equivalently}, such that for any $w\!\in\! {\bf S}^1$, the range of the map $[j^m]((\bar {\bf g}_{\mid S})_w\circ\tilde p)$ avoids $\tilde\Sigma_{0,n\!-\!2}^m$, {\em one has in symbols}
\begin{multline}\label{P'0}
({\mathcal P}_{{\bf S}^1})\ \ \forall b\!\in\!\bar I\!\times\! \{0_{{\mathbb R}^{n\!-\!2}}\},\ \ \forall(\vartheta,x,w)\!\in\!S\ \ \hbox{with} \ \ x\!\in\!{\bf g}_{\vartheta,w}^{-1}\!\{b\}\\ 
 m_x({\bf g}_{\vartheta,w}\,,X_i^{n-1})\leq m_{n-1}(n^2+\!2\!+\!c)\ .
\end{multline}
\par Set $B_i\!:=\!B(x_i,3\rho_3/4)$, choose a finite covering of $\Theta\!\times \!\bar B_i$ by compact balls $\bar \beta_j$
included in domains ${\mathcal O}_j\!\subset \!{\mathscr T}\!\times\!B(x_i,\rho_3)$ of charts. 
\vskip1mm  
\begin{fact}
The set
$${\mathcal W}_{c,i,j}'\!:=\!\{\bar{\bf g}\!\in\!C^\infty(\bar\beta_j,{\mathbb R}^n)\mid {\bf j}^m(\bar{\bf g}\circ\tilde p) \ \hbox{is transversal to}\  \Upsilon^m  \ \hbox{on}\ \bar \beta_j\}\ 
$$  is open and dense
in $C^\infty(\bar\beta_j,{\mathbb R}^n)$.
\end{fact}
\begin{proof}
As in the argument of the proof of theorem \ref{thm7.1C-Mbis},  proposition 3.6 of \cite{Fel} implies that the subset ${\mathfrak U}^m$ in $C^\infty(\bar\beta_j\!\times\!{\bf S}^1,J^m(\tilde{\mathcal L},{\mathbb R}^{n-1})\!\times\!{\bf S}^1)$ of maps transversal to ({\em meaning here}: avoiding) $\Upsilon^m$ is open.
\par {\it We now show} that the mapping ${\bf j}^m$
\begin{equation*}\bar{\bf h}\!\in\!C^\infty({\mathscr T}\!\times\!X_i^{n-1},\R^n)\mapsto {\bf j}^m \bar{\bf h}\!\in\!C^\infty({\mathscr T}\!\times\!X_i^{n-1}\!\times\!{\bf S}^1,J^m(\tilde{\mathcal L},\R^{n-1})\!\times\!{\bf S}^1)
\end{equation*}
is continuous in the $C^\infty$-topology. 
\par {\em Indeed}, {\em first} the map $\bar{\bf h}\mapsto j^m \bar{\bf h}$ is continuous (proposition 3.4 p.46 in \cite{G-G}). Next, $[\cdot]$ which maps an $m$-jet $\sigma\!\in\!J^m({\mathscr T}\!\times\!X_i^{n-1},\R^{n}\!)$ represented by a $C^\infty$ germ $\bar h\!:\!({\mathscr T}\!\times\!X_i^{n-1}\!,(\vartheta,{\bf x}))\!\rightarrow\!\R^n$ to $[\sigma]\!=\!j^m(\bar h\circ\iota_{{\mathscr T}\!\times\!L_{\bf x}})(\vartheta,{\bf x})$ in $J^m(\tilde{\mathcal L},\R^{n}\!)$, is $C^\infty$ (from local computation in a foliated chart of $X_i^{n-1}$). Thus the composed map $[j^m]:=[\cdot]\circ j^m$ is continuous in the $C^\infty$-topology by proposition 3.5, page 46 in \cite{G-G}. Then, the map $e_{{\bf S}^1}$, extending $f\!\in\!C^\infty({\mathscr T}\!\times\!X_i^{n-1},Z)$ to $e_{{\bf S}^1}\!(f)\!\in\!C^\infty({\mathscr T}\!\times\!X_i^{n-1}\!\times\!{\bf S}^1,Z\!\times\!{\bf S}^1)$ (with $Z$ manifold) defined for any $w\!\in\!{\bf S}^1$ by $e_{{\bf S}^1}\!(f)(\cdot,\ast,w):=
(f(\cdot,\ast),w),$ is continuous in the $C^\infty$-topology (from its definition and from proposition 3.6 p.47 in \cite{G-G}): so is $e_{{\bf S}^1}\circ [j^m]$ too. Given $\sigma\!\in\!J^m(\tilde {\mathcal L},\R^n)$,
select a germ $\bar h\!:\!({\mathscr T}\!\times\!X_i^{n-1}, (\vartheta,{\bf x}))\!\rightarrow \!\R^n$ such that
$[j^m]\bar h\, (\vartheta,{\bf x})\!=\!\sigma$; for any $w\!\in\!{\bf S}^1$, setting $\sigma_w\!:=\![j^m](\bar h_w)(\vartheta,{\bf x})$
defines unambiguously a {\em smooth} map $\varphi$ that sends $(\sigma,w)\!\in\!J^m(\tilde {\mathcal L},\R^n)\!\times\!{\bf S}^1$ to $(\sigma_w,w)\!\in\!J^m(\tilde {\mathcal L},\R^{n-1})\!\times\!{\bf S}^1$, both $\sigma$ and $\sigma_w$ having the same source $(\vartheta,{\bf x})$. One has ${\bf j}^m=\varphi\circ (e_{{\bf S}^1}\circ [j^m])$, implying, by proposition 3.5 p.46 in \cite{G-G}, that ${\bf j}^m$  is continuous. 
\par As $\tilde p$ is proper, the mapping $\tilde p^\ast:\bar {\bf g}\!\mapsto\!\bar {\bf g}\!\circ\!\tilde p$ is continuous (proposition 3.9 p.49 in \cite{G-G}).  Thus the composed map
\begin{multline*}{\bf j}^m\circ\tilde p^\ast:\bar {\bf g}\!\in\!C^\infty(\bar\beta_j,{\mathbb R}^n)\mapsto ({\bf j}^m\circ\tilde p^\ast)(\bar {\bf g}):=\\:={\bf j}^m(\bar {\bf g}\!\circ\!\tilde p)\!\in\!C^\infty((\tilde p^{-1}\bar\beta_j)\!\times\!{\bf S}^1,J^m(\tilde{\mathcal L},{\mathbb R}^{n-1}\!)\!\times\!{\bf S}^1)\cr
\end{multline*}
{\em is continuous} too and
${\mathcal W}'_{c,i,j}=({\bf j}^m\circ\tilde p^\ast)^{-1}{\mathfrak U}^m$ is open.

As usual, the density results from lemma 4.6 in \cite{G-G} page 53, observing that, for any given $\bar{\bf g}\!\in\!C^\infty(\bar\beta_j,{\mathbb R}^n)$, the map $\Psi_m$
defined for any polynomial map ${\bf b}\!:\!\bar \beta_j\rightarrow {\mathbb R}^n$ of degree $\!\leq\! m$ by sending $({\bf b},\vartheta, {\bf x},w)$ to $\Psi_m({\bf b},\vartheta, {\bf x},w)\!:=\!{\bf j}^m((\bar{\bf g}\!+\!{\bf b})\circ\tilde p)(\vartheta, {\bf x},w)\!\in\! J^m(\tilde{\mathcal L},{\mathbb R}^{n-1})\!\times\!{\bf S}^1$ ({\em adding a polynomial map} ${\bf b}$ makes sense in the domain of chart ${\mathcal O}_j$) is {\em onto}, thus transversal to $\Upsilon^m$. Indeed, $({\bf b},\vartheta, {\bf x},w)\mapsto{\bf j}^m({\bf b}\circ\tilde p)(\vartheta, {\bf x},w)$ is onto (from computation in a foliated chart of $X_i^{n-1}$) and one has ${\bf j}^m((\bar{\bf g}+{\bf b})\circ\tilde p)\!=\!{\bf j}^m(\bar{\bf g}\circ\tilde p)+{\bf j}^m({\bf b}\circ\tilde p)$. So, ${\mathcal W}_{c,i,j}'$ is dense. 
\end{proof}

For each $j$, apply lemma \ref{unlemme333} stated without ${\mathscr T}$, giving to $\beta_j$ the role of $V$: from the above ${\mathcal W}_{c,i,j}'$, one gets finitely many open dense subsets of $C^\infty({\mathscr T}\!\times\! M,{\mathbb R}^n)$, still named ${\mathcal W}_{c,i,j}'$, such that, for any $w\in {\bf S}^1$, any $\bar{\bf g}\in
{\mathcal W}_{c,i,j}'$ verifies a property $({\mathcal P}_{{\bf S}^1})$ (\ref{P'0}) on $S=\bar\beta_j\!\times\!{\bf S}^1$. 
\par Moreover, ${\mathcal W}_{c,i}':=\cap_j{\mathcal W}_{c,i,j}'\!\subset\!C^\infty({\mathscr T}\!\times\! M,{\mathbb R}^n)$ is an open dense set. 
\par We have established that
{\em for any $\bar {\bf g}\!\in \!{\mathcal W}'_{c,i}$ and $w\!\in \!{\bf S}^1$, one has}
\begin{multline}\label{P'1}
({\mathcal P}_w')\ \ \ \hbox{for any} \ b\!\in\!\bar I\!\times\! \{0_{{\mathbb R}^{n\!-\!2}}\},\,\vartheta\!\in\!\Theta\  \ \hbox{and} \ \ x\!\in\!{\bf g}_{\vartheta,w}^{-1}\!\{b\}\cap\bar B_i\\ 
 m_x({\bf g}_{\vartheta,w},X_i^{n-1})\leq m_{n-1}(n^2+\!2\!+\!c)\ .\hskip2cm
\end{multline}
\par

Define ${\mathcal W}'_c$ in $C^\infty({\mathscr T}\!\times\! M,{\mathbb R}^n)$ to be the open dense finite intersection ${\mathcal W}'_c\!:=\!\cap_i {\mathcal W}'_{c,i}$.
{\em Any $\bar{\bf g}\!\in\! {\mathcal W}'_c$ verifies $({\mathcal P}_w')$ {\rm(\ref{P'1})}  with $\bar W$ replacing $\bar B_i$.}

\vskip1mm
Now, we use property $({\mathcal P}'_w)$, the semi-continuity of lemma \ref{semicontmult} {\em and} fact \ref{slip1033} telling that all  $\hat\sigma\!\in\!{\mathfrak S}_{t_0}(W)$ we consider verify $\hat\sigma\!\subset\!B(x_i,3\rho_3/4)$: thanks to $({\mathcal P}_w')$ (\ref{P'1}), we cover the compact $\Theta\times{\bf S}^1\times \bar W$ by finitely many open sets ${\mathcal U}_j\,$, in each of which, for any $\hat\sigma\in{\mathfrak S}_{t_0}(W)$, any $(\!n\!-\!1\!)$-face $\hat\eta$ of $\hat\sigma$, one has, for all $(\vartheta ,w,x)\!\in\!{\mathcal U}_j$ and $b\!\in\!\bar I\!\times\!\{0_{\R^{n-2}}\}$
$$\sum_{x\in \hat\eta\,\cap\,{\bf g}_{\vartheta,w}^{-1}\!\{b\}} m_x({\bf g}_{\vartheta,w},\hat\eta)\leq m_{n-1}(n^2+\!2\!+\!c)\ .
$$  
\par
{\it If $(i)'$ was not true}, one could find a sequence $\rho_k\!>\!0$ tending to $0$, a sequence
$(\vartheta_k, w_k, a_{1,k},\dots,a_{n_k,k})\!\in\!\Theta\!\times\!{\bf S}^1\!\times \!\bar W$ converging to $(\vartheta,w, a)$, $b_k$ tending to $b\!\in\!\bar I\!\times\!\{0_{{\mathbb R}^{n\!-\!2}}\}$ and $(n\!-\!1)$-faces $\hat\eta_k $ of $\hat\sigma_k \!\in\! {\mathfrak S}_{t_0}(W)$ with ${\rm diam}(\hat\sigma_k)\,\leq 2\rho_k$ (use remark \ref{compact1099}) such that, for $q_k=1,\dots,n_k$
\begin{equation} \label{absurd-bis} 
 \sum_{a_{q_k,k}\in \hat\eta_k\,\cap\,{\bf g}_{\vartheta_k,w_k}^{-1}\!\{b_k\}}\!\!\! m_{a_{q_k,k}}({\bf g}_{\vartheta_k,w_k},\hat\eta_k)\geq 
m_{n-1}(n^2+\!2\!+\!c)+1\  .
\end{equation}
But $(\vartheta,w, a)$ belongs to some ${\mathcal U}_{j_0}$ and the choice of ${\mathcal U}_j$ expresses {\em the semi-continuity of lemma} \ref{semicontmult}, so the sum of all multiplicities at $x$ of ${\bf g}_{\vartheta,w}$ with $X_i^{n-1}$ for $(\vartheta ,w,x)\in {\mathcal U}_{j_0}$ is $\leq m_{n-1}(n^2+\!2\!+\!c)$, {\it a contradiction}: since $\hat\eta_k\!\subset\!\hat\sigma_k\!\subset\! B(a_{q_k,k},2\rho_k)$, one would have
$\{\vartheta_k\}\!\times\!\{w_k\}\!\times \!\hat\eta_k\!\subset\! {\mathcal U}_{j_0}$ for large $k\,$, together with 
(\ref{absurd-bis}), finishing the proof of fact \ref{ifirstpart}.
\end{proof}
\par
The second claim of $(i)$ (see definition {\rm\ref{testsquare1}} $(i)$) relies on the last argument in the proof of theorem \ref{theorem I}, replacing the bound $m_{n-1}(\kappa+c)$ by $m_{n-1}(n^2+\!2\!+\!c)$.

\vskip1mm 

\par
Recall again fact \ref{slip1033} and the notations of section \ref{preparatorybds}: all $\hat\sigma\!\in\!{\mathfrak S}_{t_0}(W)$ are subsets in $B(x_i,3\rho_3/4)$ for some $i$ and each $(n\!-\!2)$-face $\hat\xi$ slips into a leaf of some texture $X_i^{n\!-\!2}$. The codimension of the texture $X_i^{n-2}$ is $\kappa=(n-1)n\,$.  
\par As for $(ii)'$ below, one cares about $g_\vartheta^{-1}\{0\}$ with $\bar g\!\in\! C^\infty({\mathscr T}\!\times\! W,{\mathbb R}^{n-2})$. From theorem \ref{thm7.1C-Mbis} $(i)$ (in its notations $W=B(x_i,3\rho_3/4)$ is the basis of $X_i^{n-2}$, while $W$ is a distinct thing in this proof) and lemma \ref{unlemme333}, we get an open dense set ${\mathscr W}''_c$ of $\bar g\!\in\! C^\infty({\mathscr T}\!\times\!W,{\mathbb R}^{n\!-\!2})$ verifying for any $i$
$$({\mathcal P})\ \ \hbox{for any} \ \vartheta\!\in\!\Theta\,,\,x\!\in\! g_\vartheta^{-1}\{0\}\ \ \ m_x(g_\vartheta,X_i^{n-2})\!\leq\! m_{n-2}((n\!-\!1)n\!+\!c)\,.
$$
\par
Call $\pi_{1,2}:{\mathbb R}^n\rightarrow{\mathbb R}^{n-2}$ the projection forgetting $u_1,u_2$. Define $\Pi_{1,2}:\bar{\bf g}\!\in\! C^\infty({\mathscr T}\!\times\! W,{\mathbb R}^n) \mapsto \bar g\!:=\!\pi_{1,2}\circ\bar{\bf g}\!\in\! C^\infty({\mathscr T}\!\times\! W,{\mathbb R}^{n-2})$: this map is open, continuous, surjective. Then
$ {\mathcal W}''_c\!:=\!\Pi_{1,2}^{-1}({\mathscr W}''_c)$ is an open  dense set of $\bar{\bf g}\!=\!(\bar g_1,\bar g_2,\bar g)$ such that $({\mathcal P})$ is true over $\bar W\!\subset \!\cup_i\bar B_i\,$ for  $\bar g\!=\!\pi_{1,2}\circ\bar{\bf g}$. 
\par Given any fixed $\bar{\bf g}\!\in\!{\mathcal W}''_c$,
choose a finite covering ${\mathcal U}_j$ of the compact set $\Theta\times \bar W$ such that, for all $(\vartheta ,x)\in {\mathcal U}_j\,$, the multiplicities at intersection points $x$ of $g_\vartheta^{-1}\{0\}$ with any $(n\!-\!2)$-face $\hat \xi$ of any $\hat\sigma\in{\mathfrak S}_{t_0}(W)$ sum up to less than $m_{n-2}((n\!-\!1)n\!+\!c)$ (lemma \ref{semicontmult}). Arguing by contradiction as above ensures the existence of $\rho''_{\bar {\bf g}}>0$ such that for $\bar {\bf g}\in{\mathcal W}''_c$
$$(ii)'\ \ \hbox{for any} \ \vartheta\!\in\!\Theta\ \ \ \sum_{x\in g_\vartheta^{-1}\!\{0\}\cap\hat\xi}m_x(g_\vartheta,\hat\xi)\!\leq\! m_{n-2}((n\!-\!1)n\!+\!c)\,.
$$ 
{\it is} true {\it for any $\rho\in]0,\rho''_{\bar {\bf g}}]$, all $\hat\sigma\in{\mathfrak S}_{t_0}(W)$ of diameter $\leq 2\rho$ and any of its $(n-2)$-face $\hat\xi$}. 
\par
Set
$\rho_{\bar{\bf g}}\!:=\!\min(\rho'_{\bar {\bf g}},\rho''_{\bar {\bf g}})\!>\!0$. The open dense set
${\mathcal W}'''_c={\mathcal W}'_c\cap {\mathcal W}''_c\subset C^\infty({\mathscr T}\!\times\!W,\R^n)$ of $\bar {\bf g}$ is such that $(i)$ and $(ii)$ hold for $\bar G\!\in\!{\mathcal V}''_c$ defined by $\bar {\bf g}\!\in\!{\mathcal W}'''_c$, where ${\mathcal W}'''_c$ produces a dense open set ${\mathcal V}''_c\subset {\mathscr O}^{\mathscr T}_{\bar I^2}(W)$ through proposition \ref{genericsqu}.
\vskip1mm
We now skip to $(iii)$. {\em It might be valuable to notice here that proposition {\rm 4.5} on page {\rm 52} in \cite{G-G} is still true if one deals with transversality to a closed boundary manifold} (see also the last lines of page 54).
\par Set $\rho_{\bar G}:=\rho_{\bar{\bf g}}\,$, where $\bar G$ is defined through $\bar{\bf g}\,$.
For any $\rho\!\in]0,\rho_{\bar G}]$, consider the subdivision of integer $E$ paired with $\rho$ (see definition \ref{linkint}) giving the polyhedral approximation of theorem \ref{theorem I} and denote it, as usual, by $K$.
Apply the elementary transversality theorem (corollary 4.12 in \cite{G-G}), 
calling ${\mathcal V}_{1,E,c}\!\subset\! {\mathscr O}^{\mathscr T}_{\bar I^2}(W)$ the open dense set of $\bar G\!\in\!{\mathcal V}''_c$ transverse to all (finitely many) leaves in $X_i^{n\!-\!3}$ containing a simplex belonging to $\hat K_{n-3}$. Thom's tool in transversality (\cite{G-G}, lemma 4.6 page 53) implies that the set ${\mathcal V}_{1,E,c}$ of $\bar G$, verifies that $G_\vartheta\!\in\!\bar G$ is transverse to $\hat K_{n-3}$ (compact) for a dense open set of $\vartheta\in{\mathscr T}$, {\em thus} ${\rm rge}(G_\vartheta)\cap \hat K_{n-3}\!=\!\emptyset$. 
\par 
A similar procedure (applying corollary 4.12 in \cite{G-G}) can be followed, defining ${\mathcal V}_{2,E,c}\!\subset\! {\mathscr O}^{\mathscr T}_{\bar I^2}(W)$ to be the open dense set of $\bar G\!\in\!{\mathcal V}''_c$ such that $\partial \bar G$ is transverse to all (finitely many) leaves in $X_i^{n\!-\!2}$ containing a simplex belonging to $\hat K_{n-2}$, ensuring that $\partial ({\rm rge}(G_\vartheta))\cap \hat K_{n-2}\!=\!\emptyset$ for a dense open set of $\vartheta\in{\mathscr T}$. 
\par So, for the given $\rho$ (i.e. $K$ corresponding to the paired $E$ with $\rho$), any $\bar G\!\in\!{\mathcal V}_{E,c}\!:=\!{\mathcal V}_{1,E,c}\cap{\mathcal V}_{2,E,c}$ is such that, for an open dense set of $\vartheta\in{\mathscr T}$, one has $(iii)$ for $G_\vartheta$.
As the choice of ${\mathcal V}''_c$ ensures that properties $(i)$, $(ii)$ already hold for $\bar G\!\in\!{\mathcal V}''_c$, one gets the wanted residual set ${\mathcal V}_c$ by defining ${\mathcal V}_c:=\cap_E{\mathcal V}_{E,c}$.
\end{proof}

\section{Appendices} \label{Whitneyandal}
\subsection{Proof of lemma \ref{unlemme333}}\label{unlemme333333}$ $

As a useful preliminary to this proof and to the proof of lemma \ref{blue555}, we present Seeley's extension (see \cite{Se}) and its principal features.
Define $\R^n_+$ to be $\R^n_+\!:=\!\{x\!=\!(x_1\dots,x_n)\!\in\!\R^n\mid x_n\!\geq\!0\}$. A function is $C^\infty$ in $\R^n_+$ if there is an open $U\!\subset\!\R^n$ with $\R^n_+\!\subset\!U$ and a $C^\infty$-function $g$ in $U$ verifying $g_{\mid \R_+^n}\!=\!f$ (see \cite{D} p.18 ff.). Moreover, one has the
\begin{thm}\label{Seeley} (see {\rm\cite{Se}}) There is a continuous linear extension operator $E\!:\!C^\infty(\R_+^n,\R)\rightarrow C^\infty(\R^n,\R)$, where the topologies are those of uniform convergence of each derivative on compact subsets of $\R_+^n$ in the first place, of $\R^n$ in the second.
\end{thm}
Let us describe how is built such an operator $E$. One has the
\begin{lem}\label{Seeley33} (see {\rm\cite{Se}}) One can find sequences $a_k,b_k$ such that $(i)$ $b_k<0\,$, $(ii)$ $\sum_k \vert a_k\vert \vert b_k\vert^n\!<\!\infty$ for $n=0,1,2,\dots \,$, $(iii)$ $\sum_k a_k (b_k)^n=1$ for $n=0,1,2,\dots $ and $(iv)$ $b_k\searrow -\infty\ .$
\end{lem}
Then, choosing on $\R$ a $C^\infty$-function $\Phi$ with $1\!\geq\!\Phi\!\geq\!0$, $\Phi\!=\!1$ on $\{t\!\in\!\R\mid t\!\leq\!1\}$ and $\Phi\!=\!0$ on $\{t\!\in\!\R\mid t\!\geq\!2\}$, one defines $E\!:=\!E_{\Phi,a_k,b_k}$ satisfying theorem \ref{Seeley} by setting $E(f):=f$ if $x_n\!\geq\!0$ and, if $x_n\!\leq\!0\,$,
$$(E f)(x_1,\dots,x_{n-1},x_n):=\sum_k a_k\,\Phi(b_k\,x_n)\,f(x_1,\dots,x_{n-1},b_k\,x_n)\  .
$$
We now proceed to the proof of lemma \ref{unlemme333}.
\begin{proof}  Since there exists a $C^\infty$-function $\bar\phi\geq0$ equal to $1$ on ${\mathscr T}\!\times\!\bar V$, having a support in ${\mathscr T}\!\times\!{\mathcal O}$ as close as wished to ${\mathscr T}\!\times\!\bar V$, the definition of $C^\infty({\mathscr T}\!\times\!\bar V,{\mathbb R}^l)$ implies that the map ${\mathcal R}_{\bar V}$ is onto (use remark \ref{unlemme444}).
\vskip1mm
\par {\em For the openness, choose} a covering over ${\mathscr T}$ (over $M$) by open balls ${\mathscr T}_i$ (and ${\mathcal O_j}$) domains of charts $\psi_i$ (of charts $\varphi_j$, with $\varphi_1\!=\!\varphi,\,{\mathcal O}_1\!=\!{\mathcal O}\,$) with the $\overline{\mathscr T}_i$ compact building a locally finite covering, and choose open balls ${\mathscr T}'_i$, with $\bar {\mathscr T}'_i\!\subset\!{\mathscr T}_i$ for any $i$, satisfying ${\mathscr T}\!=\!\cup_i {\mathscr T}'_i$. Set $\omega_{i,j}\!:=\!(\psi_i,\varphi_j)$. Take $\bar f$
belonging to an open set ${\mathcal W}\!\subset C^\infty({\mathscr T}\!\times\! M,{\mathbb R}^l)$ in the {\em Whitney $C^s$-topology} and choose a continuous function $\delta=\delta_s>0$ on ${\mathscr T}\!\times\! M$ such that 
the open neighborhood $B_\delta(\bar f)$ of $\bar f$ in {\em this topology}
\begin{equation*} \{\bar g\in C^\infty({\mathscr T}\times M,{\mathbb R}^l)\mid \forall i,j\  \  \ \sum_{\vert r\vert=0}^s\vert d^r((\bar g-\bar f)\circ\omega_{i,j}^{-1})\vert<\delta\circ\omega_{i,j}^{-1}\}
\end{equation*}
($r$ is a multi-index) is contained in ${\mathcal W}$ (see \cite{G-G}, chapter II, §3 p. 42ff). 
\par Choose a small $\epsilon$-neighborhood $B_\epsilon$ of $B\!=\!B^n$ with $\bar B_\epsilon\!\subset\!\varphi({\mathcal O})$. Plugging Seeley's $C^\infty$-extension (theorem \ref{Seeley}) in the proof of Lemma 6.37 \cite{G-T} page 136 (or §5 in\cite{Ma}) gives that any $\bar h'\!\in\! C^\infty({\mathscr T}\times\bar{V},{\mathbb R}^l)$
{\it can be extended} to $\bar h''\!\in\! C^\infty({\mathscr T}\times \varphi^{\!-\!1}B_\epsilon,{\mathbb R}^l)$ so that (for any $i$ and $t\in\psi_i({\mathscr T}_i)$)
\begin{equation}\label{extension1066}\sup_{B_\epsilon}\,\sum_{\vert r\vert=0}^s\vert d^r(\bar h''\circ\omega_{i,1}^{-1})\vert\leq C_{s,\epsilon,\varphi,i}\ 
\sup_{\bar{B}^n}\,\sum_{\vert r\vert=0}^s\vert d^r(\bar h'\circ\omega_{i,1}^{-1})\vert\ .
\end{equation}
{\em Indeed}, choose $\lambda$ a partition of unity $\lambda_i$ subordinate to the covering ${\mathscr T}_i\!\times\!{\mathcal O}$ of ${\mathscr T}\!\times\!{\mathcal O}$. Let $\bar h_i''$ be the extension of $\bar h'_i\!=\!\bar h'\circ\omega_{i,1}^{-1}\,$ as in the cited lemma. One has, for any $i$ (among functions of $t\!\in\!\psi_i({\mathscr T}_i)$)
\begin{equation*}\label{extension1077}
\sup_{B_\epsilon}\,\sum_{\vert r\vert=0}^s\vert d^r(\bar h_i'')\vert\leq C'_{s,\epsilon,\varphi,i}\ 
\sup_{\bar{B}^n}\,\sum_{\vert r\vert=0}^s\vert d^r(\bar h'_i)\vert\ ,
\end{equation*}
where $C_i'\!:=\!C'_{s,\epsilon,\varphi,i}$ {\em does not depend} on $\bar h_i'$ (this is the essence of theorem \ref{Seeley}) {\em nor} on $t\!\in\!\psi_i({\mathscr T}_i)$, {\em since} the extension operator defined by Seeley involves only the variable of extension and {\em since} the straightening of $\partial B$ near one of its points involves a smooth map $\Psi$ independent of $t$.
\par
Set $\bar h''\!=\!\sum_j\!\lambda_j \,(\bar h_j''\circ\omega_{j,1})$, replace $C'_i$ by an appropriate $C_{s,\epsilon,\varphi,i}$ depending on $\lambda$ (only a finite number of ${\mathscr T}_j$ meet the compact $\overline{{\mathscr T}_i}$). This gives the wished extension $\bar h''$ and (\ref{extension1066}).
\vskip1mm\par Define $\!\bar f'\!:=\!{\mathcal R}_{\bar V}(\bar f)\!=\!\bar f_{\mid {\mathscr T}\!\times\!\bar V}$ and, for $\delta'\!\in\!C^0(\!{\mathscr T}\!\times\!\bar{V},\R_+\!\!\setminus\!\!\{0\}\!)$, $B_{\delta'}(\bar f')$ as
\begin{equation}\label{extension1088}\{\bar g'
\in C^\infty({\mathscr T}\times\bar{V},{\mathbb R}^l)\mid \forall i\ \ \ \sum_{\vert r\vert=0}^s\vert d^r((\bar g'-\bar f')\circ\omega_{i,1}^{-1})\vert<\delta'\circ\omega_{i,1}^{-1}\}\ .
\end{equation}
Thus $B_{\delta'}(\bar f')$ is a neighborhood of $\bar f'$ in the Whitney $C^s$-topology of
$C^\infty({\mathscr T}\times\bar{V},{\mathbb R}^l)$.
Now, for any $\bar g'\in B_{\delta'}(\bar f')$, set $\bar h':=\bar g'-\bar f'$ and apply the previous procedure to get $\bar h''$, then
set $\bar g :=\bar f+\phi_1\  \bar h''\,$, where $\phi_1$ is a function $\geq0$ equal to $1$ on $\bar{V}$ whose support is in $\varphi^{-1}B_\epsilon\subset {\mathcal O}={\mathcal O}_1$. 
\par In view of (\ref{extension1066}), (\ref{extension1088}), as $\bar g\!=\!\bar f$ off ${\mathscr T}\!\times\!\varphi^{\!-\!1}B_\epsilon$, one also has
\begin{equation*} \forall i\ \ \ \ \sum_{\vert r\vert=0}^s\vert d^r((\bar g-\bar f)\circ\omega_{i,1}^{-1})\vert < \delta'_i\ \ \hbox{where}\ \ \delta'_i\!:=\!C_i\ 
\sup_{\bar{B}^n}\,\delta'\circ\omega_{i,1}^{-1}\ ,
\end{equation*}
where $C_i$ depends on the above $C_{s,\epsilon,\varphi,i}$ and on the $\sup$ of derivatives up to order $s$ of the fixed function $\phi_1$ and where the continuous function $\delta_i'$ on $\psi_i({\mathscr T}_i)$ is $\!>\!0$.  We restrict $\delta_i'$ to $\psi_i(\bar{\mathscr T}'_i)$ and {\em view it as a function} on $\omega_{i,1}(\bar{\mathscr T}'_i\!\times\!{\mathcal O})$ which is constant in $u\in\varphi({\mathcal O})$. The cover $\bar{\mathscr T}'_i$ is locally finite, so choose $\delta'$ continuous verifying that, for any $i$, the above $\delta'_i$ verifies $\delta'_i\!\leq\!\delta\circ\omega_{i,1}^{-1}$ on $\omega_{i,1}(\bar{\mathscr T}'_i\!\times\!{\mathcal O})$. As $\bar g\!=\!\bar f$ off ${\mathscr T}\!\times\!{\mathcal O}$, one has $\bar g\!\in\!B_\delta(\bar f)$, $R_{\bar V}(\bar g)\!=\!\bar g'$, $\!B_{\delta'}(\bar f')\!\subset\! R_{\bar V}(\!B_\delta(\bar f)\!)\!\subset\! R_{\bar V}({\mathcal W})$, thus $R_{\bar V}$ is open.
\par\vskip1mm The continuity of ${\mathcal R}_{\bar V}$ is derived in \cite{G-G}, see proposition 3.9 and notes on page 49 (the inclusion ${\mathscr T}\!\times\!\bar V\!\hookrightarrow {\mathscr T}\!\times\!M$ is proper).
\par\vskip1mm Now, given a continuous surjective open function ${\mathcal F}:X\rightarrow Y$, if ${\mathcal U}$ is a dense open set in $X$, so is ${\mathcal F}({\mathcal U})$ in $Y$, and if ${\mathcal V}$ is a dense open set in $Y$, so is ${\mathcal F}^{\!-\!1}{\mathcal V}$ in $X$.
 {\em Indeed}, if ${\mathcal F}({\mathcal U})$ were not dense in $Y\,$, one would find an open set $V$ not meeting ${\mathcal F}({\mathcal U})$; but ${\mathcal U}\cap{\mathcal F}^{-1}(V)\not=\emptyset$, since ${\mathcal F}^{-1}(V)$ is open and ${\mathcal U}$ is dense in $X\,$, a contradiction. {\em Else}, if $x$ belongs to $X$ and ${\mathcal O}$ is any open containing $x\,$, the open set ${\mathcal F}({\mathcal O})$ meets the open dense set ${\mathcal V}$ in $Y$, say at $y'={\mathcal F}(x')$ for $x'\in{\mathcal O}$, and ${\mathcal F}^{-1}{\mathcal V}$ contains $x'$, so one deduces that ${\mathcal F}^{-1}{\mathcal V}$ is dense in $X\,$.
\end{proof}

\subsection{Proof of lemma \ref{blue555}}\label{blue555555}$ $

\begin{proof} For ${\bf g}\!\in\!{\mathfrak O}_{\bar I^2}$, set ${\bf h}\!:=\!{\bf g}^{-1}_{\mid\hbox{\tiny\rm rge}({\bf g})}$, $\ \hat I^2\!:=\!\bar I^2\!\times\! \{0_{{\mathbb R}^{n-2}}\}$ and $G(s,t)\!:=\!{\bf h}_{\mid \hat I^2}(s,t,0_{{\mathbb R}^{n-2}})$ for any $(s,t)\in \bar I^2\,$. We make the abuse $G\equiv{\bf h}_{\mid\hat I^2}$.
\par
\begin{lem} \label{blue555bister} Set ${\mathfrak O}_{\bar I^2}\!$ ($\!{\mathfrak O}_{\bar I^2}^{\mathscr T}\!$) short for ${\mathfrak O}_{\bar I^2}(U)\!$ (for ${\mathfrak O}_{\bar I^2}^{\mathscr T}(U)$). Then 
$$\Phi\!:\!{\bf g}\in {\mathfrak O}_{\bar I^2} \subset C^\infty(\bar U,\R^n)\!\longmapsto \!\Phi({\bf g}):=G\! \in \!{\mathscr O}_{\bar I^2}(U)\! \subset \!C^\infty(\bar I^2,U)\,,
$$
is a map which is open, continuous and surjective onto ${\mathscr O}_{\bar I^2}(U)$. 
\par For $\bar{\bf g}\!\in\! {\mathfrak O}_{\bar I^2}^{\mathscr T}$, define $\Phi(\bar{\bf g})\!\!=\!\bar G\! \in \!{\mathscr O}_{\bar I^2}^{\mathscr T}(U\!)$ through $\Phi({\bf g}_\vartheta)\!\!=\!G_\vartheta\!$ for any $\vartheta\!\in\!\!{\mathscr T}\!$: $\Phi$ is open, continuous;  if ${\mathscr T}$ is a ball, $\Phi$ is surjective onto ${\mathscr O}_{\bar I^2}^{\mathscr T}(U\!)$.
\end{lem} 
\begin{proof} First, the closure $\bar U\subset W$ is diffeomorphic to the closed ball $\bar B^n\!\subset\!{\mathbb R}^n$. Set $\bar B\!:=\!\bar B^n$. Viewing all mappings ${\bf g}\!\in\! {\mathfrak O}_{\bar I^2} \!\subset \!C^\infty(\bar U,\R^n)$ as mappings ${\bf g}\!\in\! {\mathfrak O}_{\bar I^2}\! \subset \!C^\infty(\bar B,\R^n)$ through this diffeomorphism, one shifts from $\Phi$ to the map (considered for $s=1,2,\dots,\infty$ and still called $\Phi$) 
$$\Phi\!:\!{\bf g}\in {\mathfrak O}_{\bar I^2} \subset C^s(\bar B,\R^n)\!\longmapsto \!G\! \in \!C^s(\bar I^2,B)\,
$$
and it is enough to prove the corresponding claims for the new $\Phi$.
The topologies on $\bar U$ (thus on $\bar B$) do not depend on the particular choice of a metric inducing those topologies, so we may work in this lemma with the ambient standard Euclidean metric on $\R^n$.

We prove each claim in two steps. First, we prove a claim without ${\mathscr T}$. It is worthwile to notice that the Whitney topology is then reduced to classical Banach space topology. Then, in presence of ${\mathscr T}$, we shift to deal with Whitney topology implicitly in the following way:
\begin{rem} \label{hitney} In presence of ${\mathscr T}$, choose locally finite coverings by relatively compact balls ${\mathscr T}_i',{\mathscr T}_i$ of ${\mathscr T}$ such that $\bar{\mathscr T}_i'\!\subset\!{\mathscr T}_i$ and ${\mathscr T}_i$ are domains of charts. The spaces
$C^s(\bar{\mathscr T}_i'\!\times\!\bar B,\R^n), C^s(\bar{\mathscr T}_i'\!\times\!\bar I^2,\R^n)$ are Banach spaces.  
\par In the Whitney topology, one may work with a sequence of inequalities $\Vert \cdot\Vert_{s,\mid \bar{\mathscr T}'_i\times \bar B}\!\leq\! \epsilon_i$ for given reals $\epsilon_i\!>\!0$ or, {\em equivalently}, with a single inequality $\Vert \cdot\Vert_{s,\mid {\mathscr T}\times \bar B}\!\leq \!\delta$ and a continuous function $\delta\!>\!0$.
\end{rem}
\par\vskip1mm
{\bf Proving the surjectivity.}
\par
Take $G\! \in \!{\mathscr O}_{\bar I^2}(B)$. 
Recall the normal {\em truncated} $\epsilon$-neighborhood of $G(\bar I^2)\equiv G(\hat I^2)$ (we denote it $[{\rm Tub} \ G(\hat I^2)]_\epsilon$, see also \ref{nortrunc}): $[{\rm Tub} \ G(\hat I^2)]_\epsilon$ is 
the {\em closure} of the set of points in $\R^n$ at distance $\!\leq \!\epsilon$ from $G(\hat I^2)$ such that the smallest distance is achieved at a point $G(I^2)$, where $I$ is the {\em open} interval. One can find $\epsilon \!>\!0$ such that $[{\rm Tub}\ G(\hat I^2)]_\epsilon$ is {\em injective}, meaning: each $x\!\in\![{\rm Tub}\ G(\hat I^2)]_\epsilon$ has a single closest point
in $G(\hat I^2)$.

\begin{lem}\label{blue555bisquater} Given $G\! \in \!{\mathscr O}_{\bar I^2}(B)$, if $\epsilon\!>\!0$ is small enough so that $[{\rm Tub}\ G(\hat I^2)]_\epsilon$ is injective, one can find an embedding $h\!:\!\bar{\mathcal O}\!\rightarrow \! B$ extending $G$ on $\bar{\mathcal O}$, where $\bar{\mathcal O}\!\subset\!\R^n$ is a differentiable ball such that ${\mathcal O}$ contains the truncated $\epsilon$-neighborhood $[{\rm Tub}\ \hat I^2]_\epsilon$ of $\hat I^2$.
\end{lem}
\begin{proof} 
First, $G$ extends as a smooth map on $J_1^2\!\times\!\{0_{{\mathbb R}^{n-2}}\}$, where $J_1$ is an open interval and $\bar I\!\subset\!J_1$, thus $G$ extends to $J^2\!\times\!\{0_{{\mathbb R}^{n-2}}\}$ with $\bar I\subset\!J\!\subset\!\bar J\!\subset\!J_1$ as a proper embedding of $\hat J^2:=\bar J^2\times \{0_{{\mathbb R}^{n-2}}\}$ into $B\!\subset\!\R^n$. Select a smooth orthonormal frame $\xi_1,\dots,\xi_n$ over ${\rm rge}(G)$ such that 
$\xi_i\!:=\!dG(e_i)/\Vert dG(e_i)\Vert$ for $i\!=\!1,2$. Using the uniform injectivity of the {\em normal truncated} exponential maps around
$\hat J^2$ and $G(\hat J^2)$ (normal to $T (\hat J^2)$ and $T (G(\hat J^2))$), one chooses injective {\em truncated} tubular neighborhoods $[{\rm Tub}\, \hat J^2]_{\epsilon_1}\!\equiv\!\hat J^2\!\times\!\bar B^{n\!-\!2}$ and $[{\rm Tub}\, G(\hat J^2)]_{\epsilon_1}$ of common radius $\epsilon_1>0$ (with $\epsilon_1>0$ small enough to guarantee $[{\rm Tub}\, G(\hat J^2)]_{\epsilon_1}\!\subset\!B$). Then, define $h$ to be such that, $h:=G$ on $\hat J^2$, $dh$ sends $e_i$ onto $\xi_i$ for any $i\!=\!3,\dots,n$ and, for any $u\!\in\!\hat J^2$, $h$ sends any point $w$ sitting on a normal line to $\hat J^2$ through $u$ at signed distance $t$ with $\vert t\vert\!<\!\epsilon_1$ to the point $h(w)$ at signed distance $t$ on the corresponding normal line to $G(\hat J^2)$ through $G(u)$. 
Completing the proof, choose $\epsilon\!\in]0,\epsilon_1[$ together with a differentiable ball ${\mathcal O}$ verifying $[{\rm Tub}\, \hat I^2]_{\epsilon}\!\subset\!{\mathcal O}$ and 
$\bar{\mathcal O}\!\subset\![{\rm Tub}\, \hat J^2]_{\epsilon_1}$.
\end{proof}
\begin{fact} \label{bluebird}There exists a diffeomorphism $\psi$ between $\overline{h({\mathcal O})}$ and $\bar B$ which keeps $G(\hat I^2)\!=\!h(\hat I^2)$ pointwise fixed.
\end{fact}
{\em Indeed}, as $\overline{h({\mathcal O})}\!\subset\!B\!\subset\!\R^n$ and $\!\bar B\setminus \!h({\mathcal O})$ is diffeomorphic to the complementary in $\bar B$ of an open round ball $B'$ with $\bar B'\!\subset B$, applying theorem 3.1, page 185 in \cite{Hi}, there exists a diffeomorphism $\Psi$ of $\bar B$ sending $\overline{B(0,1/2)}$ onto $\overline{h({\mathcal O})}$. If $\chi$ is the radial vector field on $\bar B$, if $\varphi$ is a $C^\infty$-function defined on $\bar B$, with $0\!\leq\!\varphi\!\leq\!1$, $\varphi\!=\!1$ on $\bar B\setminus h({\mathcal O})$, having support in
the complementary of the $\epsilon/2$-neighborhood of $G(\hat I^2)$ in $\bar B$, the field $\varphi\,d\Psi(\chi)$
gives by integration from $t\!=\!0$ to $t\!=\!\log 2$ {\em the} expected diffeomorphism $\psi$. 
This proves the fact.
\par\vskip1mm
So $G$ admits a smooth extension
 ${\bf h}\!:=\!\psi\circ h$ embedding a $C^\infty$-ball $\overline{{\mathcal O}}$ verifying $\hat I^2\!\subset\!{\mathcal O}\!\subset\!\R^n$ {\em onto} ${\rm rge}({\bf h})\!=\!\bar B\!\subset\!\R^n$, so ${\bf g}\!:=\!{\bf h}^{-1}$ is in ${\mathfrak O}_{\bar I^2}$ and $\Phi({\bf g})\!=\!G$. Thus $\Phi:{\mathfrak O}_{\bar I^2}\rightarrow{\mathscr O}_{\bar I^2}$ is {\em onto}. 
 \vskip2mm
{\em If ${\mathscr T}$ is a $c$-dimensional ball}, to prove that each $\bar G\! \in \!{\mathscr O}^{\mathscr T}_{\bar I^2}(B)$ may be written $\Phi(\bar{\bf g})$ with $\bar{\bf g}\! \in \!{\mathfrak O}^{\mathscr T}_{\bar I^2}(B)$, we rely on the construction above applied to a single arbitrary $\vartheta_0\!\in\!{\mathscr T}$. 
\par Given $I=]a,b[$ and $\alpha>0$, define $I_\alpha\!:=\,]a-\alpha,b+\alpha[$.
\par One can find a positive $C^\infty$-function $\epsilon_1\!:\!{\mathscr T}\!\rightarrow\! \R_+$ such that, for any $\vartheta\!\in\!{\mathscr T}$, each $G_\vartheta$ embeds $\bar I_{\epsilon_1}^2$ and each $\vartheta$-slice $\tau_{\vartheta,\epsilon_1}\!:=\!\{\vartheta\}\!\times\![{\rm Tub}\, \hat I_{\epsilon_1}^2]_{\epsilon_1}$ or $\tilde\tau_{\vartheta,\epsilon_1}\!=\!\{\vartheta\}\!\times\![{\rm Tub}\, G_\vartheta(\hat I_{\epsilon_1}^2)]_{\epsilon_1}\!\subset \!\!\{\vartheta\}\!\times\!B$ is an injective {\em truncated} tubular neighborhood. {\em Indeed}, given $\vartheta\!\in\!{\mathscr T}$, one can find $\eta(\vartheta)>0$ such that $(i)$ $G_\vartheta$ embeds $\hat I^2_{\eta(\vartheta)}$ and $(ii)$ $G_\vartheta(\hat I^2_{\eta(\vartheta)})$ has an $\eta(\vartheta)$-injective truncated tubular neighborhood $[{\rm Tub}\, \hat I_{\eta(\vartheta)}^2]_{\eta(\vartheta)}$. For any $\vartheta'$ near enough to $\vartheta$, properties $(i)$ and $(ii)$ still hold with $G_{\vartheta'}$ and $\hat I_{\eta(\vartheta)/2}^2$ in place of $G_\vartheta$ and $\hat I^2_{\eta(\vartheta)}$.
By definition of the Whitney $C^\infty$-topology, ${\mathscr T}$ being the union of compact sets, one may choose such a smooth function $\epsilon_1>0$. 
\par Using the scale $\epsilon_1$, one finds a $C^\infty$ $\vartheta$-family of diffeomorphisms $\bar h\!:\!\coprod_{\vartheta\in\mathscr T}\{\vartheta\}\! \times \!\tau_{\vartheta,\epsilon_1} \rightarrow\!\coprod_{\vartheta\in\mathscr T}\{\vartheta\}\! \times \!\tilde\tau_{\vartheta,\epsilon_1}$ by applying lemma \ref{blue555bisquater} coherently in $\vartheta$ to each pair $\tau_{\vartheta,\epsilon_1},\tilde\tau_{\vartheta,\epsilon_1}$: select a $C^\infty$ ${\mathscr T}$-family of orthonormal frames parallelising, for any $\vartheta\!\in\!{\mathscr T}$, the normal bundle to $G_\vartheta(\hat I^2_{\eta(\vartheta)})$ in $B$, which can be done {\em since ${\mathscr T}$ is a ball} (see \cite{Hu}, page 21, theorem 7.1, (H2), replacing $\cap$, misprint, by $\cup$).  One has a canonical $C^\infty$-family of diffeomorphisms $\bar E\!:=\!E_\vartheta:\tau_{{\vartheta_0,\epsilon_1}}\rightarrow\tau_{\vartheta,\epsilon_1}$, thus a $C^\infty$-family of diffeomorphisms $\bar {\tilde E}\!:=\!{\tilde E}_\vartheta:\tilde\tau_{\vartheta_0,\epsilon_1}\rightarrow\tilde\tau_{\vartheta,\epsilon_1}$, verifying $h_{\vartheta_0}\circ E_\vartheta^{-1}\!=\!
\tilde E^{-1}_\vartheta\circ h_\vartheta\,$.
\par
Consider a single $\vartheta_0\!\in\!{\mathscr T}$ and $\epsilon\!=\!\epsilon_1/2$. Apply fact \ref{bluebird} to $h_{\vartheta_0},\,\tau_{\vartheta_0,\epsilon_1},$ $\,\tilde\tau_{\vartheta_0,\epsilon_1}$. 
Let $\overline{{\mathcal O}}_{\vartheta_0}$ be a differentiable ball such that $\tau_{\vartheta_0,\epsilon}\!\subset\!{\mathcal O}_{\vartheta_0}\!\subset\!\overline{{\mathcal O}}_{\vartheta_0}\!\subset\!\tau_{\vartheta_0,\epsilon_1}$ and $\psi$ a diffeomorphism from $\overline{{\mathcal O}}_{\vartheta_0}$ onto $\bar B$.
Set $\overline{{\mathcal O}}_\vartheta\!:=\!E_\vartheta(\overline{{\mathcal O}}_{\vartheta_0})$. The $C^\infty$-family $\bar{\bf h}\!:=\!\psi\circ\tilde E^{-1}_\vartheta\circ (h_\vartheta)_{\mid \overline{{\mathcal O}}_\vartheta}\!=\! \psi\circ h_{\vartheta_0}\circ (E_\vartheta^{-1})_{\mid \overline{{\mathcal O}}_\vartheta}$ produces, inverting each diffeomorphism ${\bf h}_\vartheta$, the desired $\bar{\bf g}$, proving the surjectivity of $\Phi\,$.
\par\vskip1mm
{\bf Proving the openness.}
\par If ${\bf g}\!\in\!{\mathfrak O}_{\bar I^2}$, we prove, for any $s$, that $
G\!=\!\Phi({\bf g})$ has a neighborhood ${\mathscr V}_{G}$ of embeddings $G'\!\in\! C^\infty(\bar I^2,B)$ extendable to $C^\infty$-embeddings ${\bf h}'$ having source ${\bf g}(\bar B)$ and range ${\rm rge}({\bf h}')\!=\!\bar B$ and that the $C^s$-closeness of $G'$ to $G$ implies the $C^{s-1}$-closeness of ${\bf h}'$ to ${\bf h}\!=\!{\bf g}^{-1}$, thus the one of ${\bf g}'\!=\!({\bf h}')^{-1}$ to ${\bf g}$. As $G'\!=\!\Phi({\bf g}')$, this will allow us to derive that $\Phi$ is {\em open}. 
\par {\em Indeed}, in the Whitney topology, one may choose ${\mathscr V}_{G}$ so that, for any $G'\!\in\!{\mathscr V}_{G}$, all $G_t$ are embeddings into $B$, where $G_t\!:=\!G\!+\!t(G'\!-\!G)$ for $t\!\in\![0,1]$. For any $t\!\in\![0,1]$, any $v\!\in\!G_t(\hat I^2)$, define the $t$-dependent vector field on ${\rm rge}(G_t)$ 
\begin{equation}\label{blue555bisquinq}X_t(v)\!:=\!G'(G_t^{-1}(v))\!-\!G(G_t^{-1}(v))\ .
\end{equation} 
It is induced by the field $V\!:=\!G'-G:\hat I^2\rightarrow \R^n$ along $G_t$, verifying
\begin{equation}\label{blue555bisquinq1}\forall t\!\in\! [0,1],\ \ \forall u\!\in\!\hat I^2\ \ \  G_t(u)=\int_0^tV(u)\,d\tau+G(u)\ .
\end{equation}
The constructions we do are $C^\infty$ with $G$ and $G'$.
Extend $X$ defined on 
\begin{equation}\label{blue555bisquinq2} N=N^{G'}:=\coprod_{t\in[0,1]}\{t\}\!\times\!{\rm rge}(G_t)
\end{equation} 
as $X(t,v):=X_t(v)$ into a smooth field $Z:=E(X)$ defined on $[0,1]\!\times\!\bar B$ ($E$ is a linear continuous extension in the variable $v$, treating $t$ as a parameter, see \cite{Se}, \cite{Ma} and theorem \ref{Seeley}; use the proposition of §5 in \cite{Ma} and the proof of lemma \ref{unlemme333}). 
View $Z$ as a smooth $t$-dependent field $Z_t$ on $\bar B$ verifying $\Vert Z_t\Vert_{s}\!\leq\!C_{s,G'}\Vert X_t\Vert_{s}$ for any $t\!\in\![0,1]$. \par Actually, $C_{s,G'}$ depends {\em a priori} on $G'$ through the definition of $N\,$. At a small cost, one may select an $E$ controlled in terms of ${\mathscr V}_G$ only.
\begin{lem}\label{detendezvous} There exists $C_{s,{\mathscr V}_G}\!>\!0$, depending only on ${\mathscr V}_G\,$, such that, for any $G'\!\in\!{\mathscr V}_G$, the paired $C^\infty$-field $X$ defined through {\rm(\ref{blue555bisquinq})} may be extended to a $C^\infty$-field $Z=E(X)$ defined on $[0,1]\!\times\!\bar B$ verifying
$$\Vert Z_t\Vert_{s-1}\!\leq\!C_{s,{\mathscr V}_G}\Vert X_t\Vert_{s-1}\!\leq\!C_{s,{\mathscr V}_G}\Vert X_t\Vert_{s}\ .
$$
\end{lem}
\begin{proof}
Consider $N_0\!:=\![0,1]\!\times{\rm rge}(G)$. Doing the construction of the proof of lemma \ref{blue555bisquater}, picking up $\epsilon>0$ as there, extend ${\rm Id}_{[0,1]}\!\times\!G:[0,1]\!\times\!\hat I^2\rightarrow N_0$ in a bundle-isomorphism over $[0,1]$, denoted by $\varpi_0$, from $[0,1]\!\times\![\hbox{Tub}\,\hat I^2]_\epsilon$ onto $ N_{0,\epsilon}\!:=\![0,1]\!\times\![\hbox{Tub}\ G(\hat I^2)]_\epsilon\!\subset\![0,1]\!\times\!B$. 
\par Repeat the construction of the proof of lemma \ref{blue555bisquater} in the version with parameter (exposed above, taking here ${\mathscr T}\!=\![0,1]$), extending $(t,u)\!\in\![0,1]\!\times\!\hat I^2\mapsto (t,G_t(u))\!\in\! N$ (see (\ref{blue555bisquinq2})) into a diffeomorphism $\varpi_1$ from $[0,1]\!\times\![\hbox{Tub}\,\hat I^2]_\epsilon$ onto $N_{\epsilon}\!:=\!\prod_{t\in[0,1]}\{t\}\!\times\![\hbox{Tub}\ G_t(\hat I^2)]_\epsilon \!\subset\![0,1]\!\times\!B$. Indeed, making a good choice of ${\mathscr V}_G$, a {\em common} normal injectivity radius exists for all ranges of $G'\!\in\!{\mathscr V}_G$ and the selection of orthonormal fields  $\xi_1,\dots,\xi_n$ arising in the construction done for $G$ induces (using the proximity of the maps $G$ and $G'$) by projection and Gram-Schmidt orthonormalisation close well-determined $\xi'_1,\dots,\xi'_n$ sharing the corresponding wanted properties for $G'$.
\par Composing, $\varpi:=\varpi_1\circ\varpi_0^{-1}$ 
extends $(t,x)\!\in\!N_0\mapsto (t, G_t\circ G^{-1}(x))\!\in\!N$ as a diffeomorphism between $ N_{0,\epsilon}\!\subset\![0,1]\!\times\!B$ and $ N_\epsilon\!\subset\![0,1]\!\times\!B$. 
\par Then, choose a linear continuous extension operator $E_0$ to extend fields $\xi$ defined on $N_0$ to fields $\zeta$ defined on $[0,1]\!\times\!B$ having support in $N_{0,\epsilon}$: this is realised by making a good choice of $b_0$ in lemma \ref{Seeley33} and the control on $E_0$ is through a constant $C_{s,G}$. 
\par Now, use the previous $E_0,\varpi$ and $ \varpi^{-1}$ to extend fields $X$ defined on $N$ to fields $Z$ defined on $[0,1]\!\times\!\bar B$ having support in $N_\epsilon$: define $E$ to be
$$Z_{\cdot}=E(X)_{\cdot}:=d\varpi(\varpi(\cdot))\,(E_0(d\varpi^{-1}(\ast)(X_{\ast}))_{\varpi^{-1}(\ast)})\ \ .
$$
 One gets an inequality $C_{s-1,G'}\!\leq\!C_{s,\varpi}\,C_{s-1,G}$, where $C_{s,\varpi}$ depends only on $\Vert\varpi\Vert_s$ and $\Vert\varpi^{-1}\Vert_s$. But $\Vert\varpi\Vert_s$ and $\Vert\varpi^{-1}\Vert_s$ are controlled by $\Vert G'-G\Vert_s$ and $\Vert G\Vert_s$. Thus, for any $G'\!\in\!{\mathscr V}_G$, the constant $C_{s-1,G'}$ is controlled by $\Vert G'-G\Vert_s,\ \Vert G\Vert_s$ and $C_{s-1,G}$, thus by a constant $C_{s,{\mathscr V}_G}$.
\end{proof}
Choose an open ball $\beta$ concentric with $B$, satisfying $\bar\beta\!\subset B$ and such that any $G'\!\in\!{\mathscr V}_{G}$ verifies
$[{\rm Tub} \ G'(\hat I^2)]_\epsilon\!\subset\!\beta\,$: one may consider $G'$ to vary in a smaller neighborhood of $G$, still called ${\mathscr V}_{G}$. Choose a smooth function $\varphi$ having support in
$B$, with $0\!\leq\!\varphi\!\leq\!1$ and $\varphi\!=\!1$ on $\bar\beta$, and define $Y_t\!:=\!\varphi\,Z_t$, a compactly supported time-dependent field having integral lines over $[0,1]$ generating a family of diffeomorphisms $\{\Phi_{t_0,t},\,t_0,t\!\in\![0,1]\}$ of $\bar B$ (see \cite{Hi} chapter 8). For any $t\!\in\![0,1]$, one has
$\Vert Y_t\Vert_{s-1}\!\leq\!C_{s,{\mathscr V}_G,\varphi}\Vert X_t\Vert_{s}$
where $C_{s,{\mathscr V}_G,\varphi}$ is a positive real depending on ${\mathscr V}_G$ and $\Vert\varphi\Vert_{s-1}$.
\par
For any $u\in\hat I^2$, (\ref{blue555bisquinq}) and (\ref{blue555bisquinq1}) give 
$$G'(u)=\int_0^1X_t(G_t(u))\,dt+G(u)\ .
$$
Set ${\bf h}_0\!:=\!{\bf h}\!=\!{\bf g}^{-1}$. Define, for any $u\!\in\!{\bf g}(\bar B)$ and $ t\!\in\![0,1]$ 
$${\bf h}(t,u)\!=\!{\bf h}_u(t)\!=\!{\bf h}_t(u)\!:=\!\Phi_{0,t}({\bf h}(u))\ \ \hbox{i. e.}\ \ {\bf h}_t(u)\!=\!\int_0^t\!Y_\tau({\bf h}_\tau(u))\,d\tau\!+\!{\bf h}(u)\ .
$$
Set ${\bf h}'\!:=\!{\bf h}_1$. By definition of $Y_t$, all diffeomorphisms ${\bf h}_t$ send ${\bf g}(\bar B)$ onto $\bar B$,  so ${\bf g}'=({\bf h}')^{-1}$ is a diffeomorphism from $\bar B$ onto ${\bf g}(\bar B)$.
Both ${\bf h}'$ and ${\bf h}$ are $C^{s-1}$-close, due to $\Vert Y_t\Vert_{s-1}\!\leq\!C_{s,{\mathscr V}_G,\varphi}\Vert X_t\Vert_{s}\!\leq\!C'_{s,{\mathscr V}_G,\varphi}\Vert G'\!-\!G\Vert_{s}$ controlled by the $C^{s}$-closeness of $G,G'\!\in\!C^\infty(\bar I^2,B)$ (the last constant $C'_{s,{\mathscr V}_G,\varphi}$ takes care of $G_t^{-1}$ in (\ref{blue555bisquinq})).  

\par In $C^\infty(\bar B,{\bf g}(\bar B))$ (in $C^\infty({\bf g}(\bar B),\bar B)$) carrying the $C^{s-1}$-topology, call $E$ (and $F$) the subsets of diffeomorphisms ${\bf g}'$ from $\bar B$ onto ${\bf g}(\bar B)$ (diffeomorphisms ${\bf h}'$ from ${\bf g}(\bar B)$ onto $\bar B$).
Apply proposition 5 in §2 of \cite{Ma} (see also \cite{D}, prob.1, p.55):  a diffeomorphism ${\bf h}'\!\in\!F$ is sent to its inverse ${\bf g}'\!=\!({\bf h}')^{-1}\!\in\!E$ by a homeomorphism. Thus, if ${\mathscr V}_G$ is $C^s$-small enough in $C^\infty(\bar I^2,{\mathbb R}^n)$ around $G$, then ${\bf g}'$ belongs to a neighborhood ${\mathscr U}_{\bf g}\!\subset\!C^\infty(\bar B,{\bf g}(\bar B))$ of ${\bf g}$ which is correspondingly $C^{s-1}$-small: {\em any $C^s$-small enough ${\mathscr V}_G$ is contained in the image by $\Phi$ of ${\mathscr U}_{\bf g}$, whose $C^{s-1}$-smallness is controlled by the $C^s$-smallness of ${\mathscr V}_G$. Thus $\Phi$ is} open. 
\vskip2mm
\par To show that $\Phi$ is open in presence of ${\mathscr T}$, we follow the lines above using implicitly remark \ref{hitney} to deal with $C^s$-closeness, smallness... 

So ${\mathfrak O}_{\bar I^2}^{\mathscr T}$ replaces ${\mathfrak O}_{\bar I^2}$.
Choose $\bar{\bf g}\!\in\!{\mathfrak O}_{\bar I^2}^{\mathscr T}$ and set $\bar G\!:=\!\Phi(\bar{\bf g})$. Choose a neighborhood ${\mathscr V}_{\bar G}$ so that, for any $\bar G'\!\in\!{\mathscr V}_{\bar G}$, all $G_{t,\vartheta}$ are embeddings into $B$, with $\bar G_t\!:=\!\bar G\!+\!t(\bar G'\!-\!\bar G)$ for $t\!\in\![0,1]$. For any $t\!\in\![0,1],v\!\in\!G_t(\hat I^2)$, define, on the manifold ${\mathscr T}\!\times_{\bar G}\!\hat I^2\!:=\!\coprod_{\vartheta\!\in\!{\mathscr T}}\{\vartheta\}\!\times\!{\rm rge}(G_{t,\vartheta})$, the vector field  
\begin{equation}\label{blue555bisquinqseptimus}\bar X_t(v):=X_{t,\vartheta}(v)\!:=\!G_{\vartheta}'(G_{t,\vartheta}^{-1}(v))\!-\!G_{\vartheta}(G_{t,\vartheta}^{-1}(v))\ .
\end{equation} 
Again, we use the proposition of §5 in \cite{Ma} to extend the field $\bar X$ defined on $N:=\coprod_{t\in[0,1],\vartheta\in{\mathscr T}}\{(t,\vartheta)\}\!\times\!{\rm rge}(G_{t,\vartheta})$ as $\bar X({t,\vartheta},v):=X_{t,\vartheta}(v)$ into a smooth field $\bar Z:=E(\bar X)$ defined on $[0,1]\!\times\! {\mathscr T}\!\times\!\bar B$. As ${\mathscr T}$ is present, we proceed carefully. {\em Step by step, using remark {\rm\ref{hitney}}, produce a field} $\bar Z$ tangent to each slice $(t,\vartheta)\!=\!constant$, extending in the variable running $\bar I^2$, treating $t,\vartheta$ as parameters (use proposition of §5 in \cite{Ma}, proceed as in the proof of lemma \ref{detendezvous}, choosing a function $\epsilon(\vartheta)>0$): defining $N_i\!:=\!\coprod_{t\in[0,1],\vartheta\in\bar{\mathscr T}'_i}\{(t,\vartheta)\}\!\times\!{\rm rge}(G_{t,\vartheta})$, get first an extension $\bar Z_i$ on $[0,1]\!\times\! {\mathscr T}_i\!\times\!\bar B$ tangent to each slice $(t,\vartheta)\!=\!constant$ of each $\bar X_i\!:=\!\bar X_{\mid N_i}$. Denote by $\lambda_i$ a partition of unity subordinate to the covering ${\mathscr T}_i'$. Then, $\bar Z\!:=\!\sum_i\lambda_i\bar Z_i$ itself is tangent to each slice $(t,\vartheta)\!=\!constant$. View $\bar Z$ as a smooth $t$-dependent field $\bar Z_t$ on ${\mathscr T}\!\times\!\bar B$ whose $C^{s-1}$-smallness is governed ($\bar G$ and $s$ are fixed) by the $C^s$-smallness of $\bar X_t$ (and by $\Vert\lambda_i\Vert_{s-1}$), thus by the $C^s$-closeness of $\bar G$ and $\bar G'\,$, i. e. by the choice of ${\mathscr V}_{\bar G}$ and $\lambda_i\,$.
\par
Choose $b\!:=\!\coprod_{\vartheta\in{\mathscr T}}\{\vartheta\}\!\times\!\beta_\vartheta$ smooth, where $\beta_\vartheta, B$ are open concentric balls, with $\bar\beta_\vartheta\!\subset \!B$, so that any $\bar G'\!\in\!{\mathscr V}_{G}$ verifies
$[{\rm Tub}\,G_\vartheta'(\hat I^2)]_{\epsilon(\vartheta)}\!\subset\!\beta_\vartheta\,$ for any $\vartheta$: work successively in any $\bar{\mathscr T}_i'$ of remark \ref{hitney}, restricting ${\mathscr V}_{\bar G}$ at each step if necessary. Choose a smooth function $\bar\varphi$ having support in
${\mathscr T}\!\times\!B$, with $0\!\leq\!\bar\varphi\!\leq\!1$ and $\bar\varphi\!=\!1$ on $b$, and set $\bar Y_t\!:=\!\bar\varphi\,\bar Z_t$. For given $\bar G, s, \bar\varphi$, one sees that
$\bar Y_t$ is $C^{s-1}$-small if $\bar X_t$ is $C^{s}$-small. More, $\bar Y_t,\bar X_t$ are both {\em tangent to each slice $\vartheta=constant$}. Having compact support in any slice $\vartheta=constant$, the field $\bar Y_t$ integrates in global diffeomorphisms $\Phi_{t_0,t}$ of ${\mathscr T}\!\!\times\!\bar B$ for any $t_0,t\!\in\![0,1]$.
For any $u\!\in\!\hat I^2$, (\ref{blue555bisquinqseptimus}) gives 
$$\bar G'(u)=\int_0^1\bar X_t(\bar G_t(u))\,dt+\bar G(u)\ .
$$
Define $\bar{\bf h}':=\bar{\bf h}_1$ from $\bar{\bf h}_0=\bar{\bf h}$, where, for any $\vartheta\!\in\!{\mathscr T},u\!\in\!{\bf g}_\vartheta(\bar B), t\!\in\![0,1]$
$${\bf h}(t,\vartheta,u)={\bf h}_{t,\vartheta}(u):=\int_0^tY_{\tau,\vartheta}({\bf h}_{\tau,\vartheta}(u))\,d\tau+{\bf h}_\vartheta(u)=\Phi_{0,t}({\bf h}_\vartheta(u))\ .
$$
Both $\bar{\bf h}'$ and $\bar{\bf h}=\bar{\bf g}^{-1}$ are $C^{s-1}$-close, as $\bar Y_t$ is $C^{s-1}$-small in view of the $C^{s}$-closeness of $\bar G,\bar G'\!\in\!C^\infty({\mathscr T}\!\times\!\bar I^2,B)$ (choice of ${\mathscr V}_{\bar G}, \lambda_i$). Moreover, from its definition, ${\bf h}'_\vartheta$ is, for any $\vartheta$, with ${\bf h}_\vartheta$ a diffeomorphism from ${\bf g}_\vartheta(\bar B)$ onto $\bar B$, so ${\bf g}_\vartheta'=({\bf h}_\vartheta')^{-1}$ is a diffeomorphism from $\bar B$ onto ${\bf g}_\vartheta(\bar B)$. 

\par Call ${\mathscr T}\!\times_{\bar {\bf g}}\!\bar B$ the manifold $\coprod_{\vartheta\in{\mathscr T}}\{\vartheta\}\!\times\!{\bf g}_\vartheta(\bar B)$. In $C^\infty({\mathscr T}\!\times\!\bar B,\R^n)$ and $C^\infty({\mathscr T}\!\times_{\bar {\bf g}}\!\bar B,\bar B)$ equipped with the Whitney $C^{s-1}$-topology, call $E$ (and $F$) the subsets of $C^\infty$-families $\bar{\bf g}'$ ($C^\infty$-families $\bar{\bf h}'$) of diffeomorphisms ${\bf g}'_\vartheta$ from $\bar B$ onto ${\bf g}_\vartheta(\bar B)$ (of diffeomorphisms ${\bf h}'_\vartheta$ from ${\bf g}_\vartheta(\bar B)$ onto $\bar B$). 
Such a $\bar{\bf g}'$ gives rise to a diffeomorphism $\tilde{\bf g}'$
$$\tilde{\bf g}'\!:\!(\vartheta,v)\!\in\!{\mathscr T}\!\times\!\bar B\longmapsto\!(\vartheta,{\bf g}'_\vartheta(v))\!\in\!{\mathscr T}\!\times_{\bar {\bf g}}\!\bar B
$$
with inverse $\tilde{\bf h}'$
$$\tilde{\bf h}'\!:\!(\vartheta,u)\!\in\!{\mathscr T}\!\times_{\bar {\bf g}}\!\bar B\longmapsto\!(\vartheta,{\bf h}'_\vartheta(u))\!\in\!{\mathscr T}\!\times\!\bar B\  .
$$

One deduces by proposition 5 in §2 of \cite{Ma} that
the map sending $\bar{\bf h}'\!\in\!F$ to $\bar{\bf g}'\!\in\!E$ is a homeomorphism in the given topologies. Thus, if ${\mathscr V}_{\bar G}$ is $C^{s-1}$-small enough, $\bar{\bf g}'$ belongs to a correspondingly $C^{s-1}$-small neighborhood ${\mathscr U}_{\bar{\bf g}}\!\subset\!{\mathfrak O}_{\bar I^2}^{\mathscr T}\!\subset\!C^\infty({\mathscr T}\!\times\!\bar B,{\mathbb R}^n)$ of $\bar{\bf g}\,$: {\em so, any small enough $C^s$-neighborhood of $\bar G$ in $C^\infty({\mathscr T}\!\times\!\bar I^2,{\mathbb R}^n)$ is in the image by $\Phi$ of a neighborhood of $\bar{\bf g}$ in the $C^{s-1}$-topology on $C^\infty({\mathscr T}\!\times\!\bar B,\R^n)$: $\Phi$ is open}.
\par\vskip1mm
{\bf Proving the continuity.}
\par Select a ${\bf g}$ in the open ${\mathfrak O}_{\bar I^2}$. If $\beta$ is a ball concentric to $B$ with $\bar\beta\!\subset\!B$, such that $\hat I^2\!\subset\!{\bf g}(\beta)$, there exists an open $C^s$-neighborhood ${\mathscr U}_{\bf g}$ of embeddings ${\bf g}'$ verifying $\hat I^2\!\subset\!{\bf g}'(\beta)$ and ${\bf g}'(\bar\beta)\!\subset\!{\bf g}(B)$. Observe that the restriction ${\mathcal R}_{\bar\beta}\!:\!C^\infty(\bar B,\R^n)\!\rightarrow\! C^\infty(\bar\beta,\R^n)$ is $C^s$-continuous (lemma \ref{unlemme333}). Furthermore, if $\Vert{\bf g}_{\mid \bar\beta}'-{\bf g}_{\mid \bar\beta}\Vert_s<\eta$ with $\eta>0$ small enough, any ${\bf g}_{t\mid \bar\beta}={\bf g}_{\mid \bar\beta}+t({\bf g}_{\mid \bar\beta}'-{\bf g}_{\mid \bar\beta})$ is an embedding for any $t\in[0,1]$ with $\hat I^2\!\subset\!{\bf g}_{t\mid \bar\beta}(\beta)$ and ${\bf g}_{t\mid\bar\beta}(\bar\beta)\!\subset\!{\bf g}(B)$. 

For any $t\!\in\![0,1]$, define the $t$-dependent field at $v\!\in\!{\rm rge}({\bf g}_{t\mid \bar\beta})$ 
\begin{equation}\label{blue555bisquinqbisbis}{\bf X}_t(v)\!:=\!{\bf g}_{\mid \bar\beta}'({\bf g}_{t\mid \bar\beta}^{-1}(v))\!-\!{\bf g}_{\mid \bar\beta}({\bf g}_{t\mid \bar\beta}^{-1}(v))\ .
\end{equation} 
Again, extend the field ${\bf X}$ defined on ${\bf N}\!:=\!\coprod_{t\in[0,1]}\{t\}\!\times\!{\bf g}_{t\mid \bar\beta}(\bar\beta)$ as ${\bf X}(t,v)\!:=\!{\bf X}_t(v)$ into a smooth field ${\bf Z}\!:=\!E({\bf X})$ defined on $[0,1]\!\times\!{\bf g}(\bar B)$ (proposition of §5 in \cite{Ma}, following a similar path to lemma \ref{detendezvous}). View ${\bf Z}$ as a smooth $t$-dependent field ${\bf Z}_t$ on ${\bf g}(\bar B)$ which verifies (for any $ t\!\in\![0,1]$) $\Vert {\bf Z}_t\Vert_{s-1}\!\leq\!C_{s,{\mathscr U}_{\bf g}}\Vert {\bf X}_t\Vert_{s}$, where $C_{s,{\mathscr U}_{\bf g}}$ depends on ${\mathscr U}_{\bf g}\,$.
\vskip1mm\par
Choose a $C^\infty$-function $\varphi_1$ on $\R^n$ such that $0\!\leq\!\varphi_1\!\leq\!1$ and $\varphi_1\!=\!1$ on ${\bf g}(\bar\beta)$, having support in
${\bf g}(B)$, and define ${\bf Y}_t:=\varphi_1\,{\bf Z}_t$. One has
$\Vert {\bf Y}_t\Vert_{s-1}\!\leq\!C_{s,{\mathscr U}_{\bf g},\varphi_1}\Vert {\bf X}_t\Vert_{s}$
where $C_{s,{\mathscr U}_{\bf g},\varphi_1}$ is a positive constant depending on ${\mathscr U}_{\bf g}$ and $\Vert\varphi_1\Vert_{s-1}$. 
Thanks to the choice of ${\bf Y}_t$, for any $u\!\in\!\bar B$ and $t\!\in\![0,1]$, the integral curve $\check{\bf g}_{t}(u)$ of ${\bf Y}_t$ starting for $t\!=\!0$ from ${\bf g}(u)$ and ending at $\check{\bf g}\!:=\!\check{\bf g}_{1}$ builds a $C^\infty$-diffeomorphism $\check{\bf g}_{t}$ from $\bar B$ onto ${\bf g}(\bar B)$, whose $C^{s-1}$-closeness to ${\bf g}$ is controlled by $\Vert {\bf g}'-{\bf g}\Vert_s$. From (\ref{blue555bisquinqbisbis}), one gets for any $u\in\bar\beta$ (since, by construction, ${\bf X}_t={\bf Y}_t$ over $\bar\beta$)
\begin{equation}\label{blue555666}\check{\bf g}_{\mid \bar\beta}(u)={\bf g}_{\mid \bar\beta}'(u)=\int_0^1{\bf X}_t({\bf g}_{t\mid \bar\beta}(u))\,dt+{\bf g}_{\mid \bar\beta}(u)\ .
\end{equation}
\par
Applying again proposition 5 in §2 of \cite{Ma}, one gets inverses $\check{\bf h}$ to $\check{\bf g}$ whose $C^{s-1}$-closeness to ${\bf h}$ is also controlled by $\Vert {\bf g}'\!-\!{\bf g}\Vert_s$. From the choices of $\beta,\,{\bf Y}$ and thanks to  (\ref{blue555666}), any ${\bf g}'$ as above verifies ${\bf g}'\!=\!\check{\bf g}$ on $\bar\beta$, its inverse ${\bf h}'$ coincides with $\check{\bf h}$ on $\hat I^2$, implying the $C^s$ to $C^{s-1}$-continuity of $\Phi$ (proposition 3.9 p.49 in \cite{G-G}), which is true for any $s>0$ and thus implies its $C^\infty$-continuity. 

\vskip2mm
\par Finally, the continuity in presence of ${\mathscr T}$. Take a $\bar{\bf g}$ in the open ${\mathfrak O}_{\bar I^2}^{\mathscr T}$. Choose a smooth $b:=\coprod_{\vartheta\in{\mathscr T}}\{\vartheta\}\!\times\!\beta_\vartheta$ with $\beta_\vartheta, B$ concentric balls and $\bar\beta_\vartheta\!\subset\!B$, in such a way that, for any $\vartheta\in{\mathscr T}$, $\hat I^2\!\subset\!{\bf g}_\vartheta(\beta_\vartheta)$. 
\par There exists an open $C^s$-neighborhood ${\mathscr U}_{\bar{\bf g}}$ of families $\bar{\bf g}'$ verifying {\em first},  for any $\vartheta\in{\mathscr T}$, one has $\hat I^2\!\subset\!{\bf g}_\vartheta'(\beta_\vartheta)$ and ${\bf g}_\vartheta'(\bar\beta_\vartheta)\!\subset\!{\bf g}_\vartheta(B)$, {\em and also} that any $\bar{\bf g}'\!\in\!{\mathscr U}_{\bar{\bf g}}$ is so $C^s$-close to $\bar{\bf g}$ that, for any $\vartheta\!\in\!{\mathscr T}$, ${\bf g}_{t,\vartheta\,\mid \bar\beta_\vartheta}\!\!:=\!{\bf g}_{\vartheta\,\mid \bar\beta_\vartheta}\!+t({\bf g}_{\vartheta\,\mid \bar\beta_\vartheta}'\!-{\bf g}_{\vartheta\,\mid \bar\beta_\vartheta}\!)$ is an embedding for any $t\!\in\![0,1]$ with $\hat I^2\!\subset\!{\bf g}_{t,\vartheta\,\mid \bar\beta_\vartheta}(\beta_\vartheta)$ and ${\bf g}_{t,\vartheta\,\mid \bar\beta_\vartheta}(\bar\beta_\vartheta)\!\subset\!{\bf g}_\vartheta(\bar B)$ (use remark \ref{hitney} and a step by step argument, as in the proof of the openness in presence of parameter).

For any $t\!\in\![0,1],\vartheta\!\in\!{\mathscr T}$, define the $t$-dependent field at $v\!\in\!{\rm rge}({\bf g}_{t,\vartheta\mid\bar\beta_\vartheta})$ 
\begin{equation}\label{blue555bisquinqbisbisatrox}\bar{\bf X}_t(\vartheta,v)\!:=\!{\bf X}_{t,\vartheta}(v)\!:=\!{\bf g}_{\vartheta\,\mid \bar\beta_\vartheta}'({\bf g}_{t,\vartheta\,\mid \bar\beta_\vartheta}^{-1}(v))\!-\!{\bf g}_{\vartheta\,\mid \bar\beta_\vartheta}({\bf g}_{t,\vartheta\,\mid \bar\beta_\vartheta}^{-1}(v))\,.
\end{equation} 
Extend $\bar{\bf X}$ defined on ${\bf N}\!:=\!\coprod_{t\in[0,1],\vartheta\in{\mathscr T}}\{(t,\vartheta)\}\!\times\!{\bf g}_{t,\vartheta\,\mid \bar\beta_\vartheta}(\bar\beta_\vartheta)$ as $\bar{\bf X}(t,v)\!:=\!\bar{\bf X}_t(v)\!=\!{\bf X}_{t,\cdot}(v)$ into a smooth field $\bar{\bf Z}\!:=\!E(\bar{\bf X})$ on $[0,1]\!\times\!{\mathscr T}\!\times_{\bar g}\!\bar B$ {\em tangent to any slice $(t,\vartheta)\!=\!constant$}. To achieve this goal, use proposition of §5 in \cite{Ma}, remark \ref{hitney}, follow a path similar to lemma \ref{detendezvous} and a step by step argument as before. View $\bar{\bf Z}$ as a smooth $t$-dependent field $\bar{\bf Z}_t$ on ${\mathscr T}\!\times_{\bar{\bf g}}\bar B$ whose $C^{s-1}$-smallness is governed by the one of $\bar{\bf X}_t$, thus by the $C^s$-closeness of $\bar{\bf g}$ and $\bar{\bf g}'\,$.
Choose a smooth $\bar\varphi_1:{\mathscr T}\!\times \R^n\rightarrow \R_+$ verifying $0\!\leq\!\bar\varphi_1\!\leq\!1$ and $\bar\varphi_1\!=\!1$ on $b$, with support in
${\mathscr T}\!\times_{\bar{\bf g}}B$. Define $\bar{\bf Y}_t\!:=\!\bar\varphi_1\,\bar{\bf Z}_t\,$. The $C^{s-1}$-smallness of $\bar {\bf Y}_t$ is governed by the one of
$\bar {\bf X}_t \,$, and $\bar {\bf Y}_t, \bar{\bf Z}_t$ are {\em tangent to any slice $\vartheta\!=\!constant$}, hence the existence over $[0,1]$ holds for integral curves ($\bar{\bf Y}$ has compact support in any slice $(t,\vartheta)\!=\!constant$). In view of the choice of $\bar{\bf Y}_t$, for any $\vartheta\!\in\!{\mathscr T}, \ u\!\in\!\bar B$ and $t\!\in\![0,1]$, the integral curve $\check{\bf g}_{t,\vartheta}(u)$ of $\bar{\bf Y}_t$ starting for $t\!=\!0$ from ${\bf g}_{\vartheta}(u)$ and ending at ${\check{\bf g}}_\vartheta\!:=\!{\check{\bf g}}_{1,\vartheta}$ builds a $C^\infty$-family of diffeomorphisms ${\check{\bf g}}_{t,\vartheta}$ from $\bar B$ onto ${\bf g}_\vartheta(\bar B)$, whose $C^{s-1}$-closeness to $\bar{\bf g}$ is controlled by the $C^s$-closeness of $\bar{\bf g}'$ to $\bar{\bf g}\,$. But (\ref{blue555bisquinqbisbisatrox}) gives for any $\vartheta\!\in\!{\mathscr T}$ and any $u\in\bar\beta_\vartheta$ 
$$\check{\bf g}_{\vartheta\,\mid \bar\beta_\vartheta}(u)={\bf g}_{\vartheta\,\mid \bar\beta_\vartheta}'(u)=\int_0^1{\bf X}_{t,\vartheta}({\bf g}_{t,\vartheta\,\mid \bar\beta_\vartheta}(u))\,dt+\bar{\bf g}_{\vartheta\,\mid \bar\beta_\vartheta}(u)\,.
$$
The same argument developped while proving the openness in presence of ${\mathscr T}$ and proposition 5 in §2 of \cite{Ma} shows that the family $\bar{\check{\bf h}}$ of inverses ${\check{\bf h}}_\vartheta$ to ${\check{\bf g}}_\vartheta$ is $C^{s-1}$-close to $\bar{\bf h}$ since $\bar{\check{\bf g}}$ and $\bar{\bf g}$ are so, for $\bar{\bf g}'$ and $\bar{\bf g}$ are $C^s$-close. By construction, for any $\bar{\bf g}'$ as above, $\bar{\bf h}'$ coincides with $\bar{\check{\bf h}}$ on $\hat I^2$, so, using again proposition 3.9 of \cite{G-G}, one gets the $C^s$ to $C^{s-1}$-continuity of $\Phi$. And thus its $C^\infty$-continuity.
 \end{proof}

Lemma  \ref{blue555} follows from definition \ref{topology1044} and above lemma \ref{blue555bister} in the same way lemma \ref{unlemme333} was derived in the last lines of its proof.
\end{proof}

\subsection{Proof of proposition \ref{genericsqu}}\label{genericsqu-333}\begin{proof}
\begin{fact}\label{step0} Given any $\bar G\!\in\!{\mathscr O}^{\mathscr T}_{\bar I^2}(W),\,\vartheta_0\!\in\!{\mathscr T}$, select relatively compact open differentiable balls $U_0,{\mathscr T}_0$ with $\bar U_0\!\subset\!W\,,\,\vartheta_0\!\in\!{{\mathscr T}_0},\bar{{\mathscr T}_0}\!\subset\!{\mathscr T}$ such that $\hbox{\rm rge}(G_\vartheta)\!\subset\! U_0$ for any $\vartheta\!\in\!\bar{\mathscr T}_0$, i. e. $\bar G_{\mid\bar{\mathscr T}_0}\!\in\!{\mathscr O}^{\bar{\mathscr T}_0}_{\bar I^2}(U_0)$. By lemma {\rm\ref{blue555}}, one can find a family $\bar{\bf g}_0\!\in\!{\mathfrak O}^{\bar{\mathscr T}_0}_{\bar I^2}(U_0)$ defining $\bar G_{\mid \bar{\mathscr T}_0}=\{G_\vartheta\mid \vartheta\!\in\!\bar{\mathscr T}_0\}=\Phi(\bar{\bf g}_0)$. 
\end{fact}

Setting ${\mathcal U}_c^{\bar{{\mathscr T}_0}\times\bar U_0}\!:=\!{\mathcal U}_c\cap{\mathcal R}_{\bar{{\mathscr T}_0}\times\bar U_0}^{-1}\,{\mathfrak O}^{\bar{{\mathscr T}_0}}_{\bar I^2}(U_0)$, by lemma \ref{unlemme333} and its proof, as ${\mathcal R}_{\bar{{\mathscr T}_0}\times\bar U_0}:{\mathcal R}_{\bar{{\mathscr T}_0}\times\bar U_0}^{-1}\,{\mathfrak O}^{\bar{{\mathscr T}_0}}_{\bar I^2}(U_0)\mapsto{\mathfrak O}^{\bar{{\mathscr T}_0}}_{\bar I^2}(U_0)$ is continuous, open and surjective, ${\mathcal R}_{\bar{{\mathscr T}_0}\times\bar U_0}({\mathcal U}_c^{\bar{{\mathscr T}_0}\times\bar U_0})$ is a dense open set in the open set ${\mathfrak O}^{\bar{\mathscr T}_0}_{\bar I^2}(U_0)\!=\! {\mathcal R}_{\bar{{\mathscr T}_0}\times\bar U_0}{\mathcal R}_{\bar{{\mathscr T}_0}\times\bar U_0}^{-1}\,{\mathfrak O}^{\bar{{\mathscr T}_0}}_{\bar I^2}(U_0)$. Indeed, any $\bar {\bf g}_0\!\in\!{\mathfrak O}^{\bar{{\mathscr T}_0}}_{\bar I^2}(U_0)$ may be extended as a map of $C^\infty(\!{\mathscr T}\!\!\times\!\!W,\R^n\!)$, see the proof of lemma  \ref{unlemme333}.

\par By lemma \ref{blue555}, all $\bar {\bf g}'_0$ in ${\mathcal R}_{\bar{{\mathscr T}_0}\times\bar U_0}({\mathcal U}_c^{\bar{{\mathscr T}_0}\times\bar U_0})$ build, in the open ${\mathscr O}^{\bar{\mathscr T}_0}_{\bar I^2}(U_0)\!\subset\!C^\infty(\bar{\mathscr T}_0\!\times\!\bar I^2, U_0)$, an open dense set ${\mathcal V}'_0$ of families of embeddings $\bar G_0'$.

\par Now, relying on fact \ref{step0}, choose sequences of relatively compact balls ${\mathscr T}_i,U_i$ (with $i\!=\!1,2,\dots$) such that $\bar{\mathscr T}_i$ is a locally finite cover of ${\mathscr T}$ and choose $\bar{\bf g}_i\!\in\!{\mathfrak O}^{\bar{\mathscr T}_i}_{\bar I^2}(U_i)$ defining $G_\vartheta$ embedding into $U_i\!\subset\!W$ for any $\vartheta\in\bar{\mathscr T}_i$. Thus (for any $i$) $\bar G_{\mid\bar{\mathscr T}_i}$ belongs to ${\mathscr O}^{\bar{\mathscr T}_i}_{\bar I^2}(U_i)$.

\par Consider ${\mathcal R}_i\!:=\!{\mathcal R}_{\bar{\mathscr T}_i}$ restricting from $C^\infty(\!{\mathscr T}\!\times\!\bar I^2,W\!)$ to $C^\infty(\!\bar{\mathscr T}_i\!\times\!\bar I^2,W\!)$. The sets ${\mathcal R}_{\bar{{\mathscr T}_i}\times\bar U_i}({\mathcal U}_c^{\bar{{\mathscr T}_i}\times\bar U_i})$ build a dense open set ${\mathcal V}'_i$ in ${\mathscr O}^{\bar{\mathscr T}_i}_{\bar I^2}(U_i)$. Define the dense open set ${\mathcal V}_j\!:=\!\cap_{i=1}^j{\mathcal R}_i^{-1}{\mathcal V}_i'$ in the open ${\mathscr O}_j\!:=\!\cap_{i=1}^j{\mathcal R}_i^{-1}{\mathscr O}^{\bar{\mathscr T}_i}_{\bar I^2}(U_i)$. Define ${\mathscr O}_{\bar G}\!:=\!\lim_{j\rightarrow\infty}{\mathscr O}_j$: it contains $\bar G$. Set
${\mathcal V}_{\bar G}\!:=\!\lim_{j\rightarrow\infty}{\mathcal V}_j$. 
\par Call $j(i)$ the smallest $j$ such that $l\!\geq\!j$ implies $\bar{\mathscr T}_{l}\cap\bar {\mathscr T}_i\!=\!\emptyset$. Then, for $i$ given and  $j\!\geq\!j(i)$, one has ${\mathcal R}_i{\mathscr O}_{\bar G}\!=\!{\mathcal R}_i{\mathscr O}_j\!=\!{\mathcal R}_i{\mathscr O}_{j(i)}$ and ${\mathcal R}_i{\mathcal V}_{\bar G}\!=\!{\mathcal R}_i{\mathcal V}_j\!=\!{\mathcal R}_i{\mathcal V}_{j(i)}$.
In the Whitney $C^\infty$-topology a set ${\mathcal V}_{\bar G}\!\subset\!C^\infty({\mathscr T}\!\times\!\bar I^2,W)$ is {\em open} (is {\em dense} in a subset ${\mathcal S}$) {\em if and only} if each ${\mathcal R}_i({\mathcal V}_{\bar G})$ is {\em open} (is {\em dense} in ${\mathcal R}_i({\mathcal S})$) in the Whitney topological space $C^\infty(\bar{\mathscr T}_i\!\times\!\bar I^2,W)$. So, one gets the conclusion: indeed, ${\mathscr O}_{\bar G}$ is an open neighborhood of the initially given $\bar G$ in ${\mathscr O}^{{\mathscr T}}_{\bar I^2}(W)$ and ${\mathcal V}_{\bar G}$ is dense and open in ${\mathscr O}_{\bar G}$.
\end{proof}
\vfill\eject

\vskip2mm

daniel.meyer@imj-prg.fr 
\hskip2mm and \hskip2mm eric.toubiana@imj-prg.fr

\end{document}